\documentclass[a4paper,twoside, 10pt]{scrartcl}
\usepackage[bookmarks,bookmarksnumbered,plainpages=false,hyperfootnotes=false]{hyperref}
\usepackage[english]{babel}
\usepackage[fixlanguage]{babelbib}
\usepackage{cite, makeidx}
\usepackage{amssymb, amsmath, amstext, amsopn, amsthm, amscd, amsxtra, amsfonts}
\usepackage{yfonts, mathrsfs}
\usepackage{stackrel}
\usepackage{paralist}

\usepackage{multicol}
\usepackage[T1]{fontenc}
\usepackage[latin9]{inputenc}
\usepackage{sectsty}
\usepackage{float}
\usepackage{ifthen}
\usepackage{colonequals}

\usepackage[toc]{glossaries}

\usepackage[all]{xy}
\SelectTips{cm}{}
\usepackage{graphics}
\usepackage[pdftex]{graphicx}
\usepackage{tikz}

\usepackage{color}
\swapnumbers

\usepackage{fancyhdr}

{%

\fancyhead[LE, RO]{{\textbf{\thepage}}}
\fancyhead[RE]{{\textbf{ \small \leftmark}}}
\fancyhead[LO]{{\textbf{ \small \rightmark}}}

\cfoot[]{}
}

\sectionfont{\sffamily \bfseries}
\subsectionfont{\sffamily \bfseries}

\newtheoremstyle{standard}
 {16pt}  
 {16pt}  
 {}  
 {}  
 {\bfseries}
 {}  
 { } 
 {{{\thmnumber{#2.~}\thmname{#1}}}\thmnote{~(#3)}} 

\newtheoremstyle{kursiv}
 {16pt}  
 {16pt}  
 {\itshape}  
 {}  
 {\bfseries}
 {}  
 { } 
 {{{\thmnumber{#2.~}\thmname{#1}}}\thmnote{~(#3)}} 

\theoremstyle{standard}
\newtheorem{defn} [subsubsection]{Definition}
\newtheorem{ex} [subsubsection]{Example}
\newtheorem{rem}   [subsubsection]{Remark}
\newtheorem{nota}   [subsubsection]{Notation}
\newtheorem{setup} [subsubsection]{}

\newtheorem{con}   [subsubsection]{Construction} 
\newtheorem{conv}   [subsubsection]{Conventions}
\newtheorem{Asetup}[subsection]{}
\newtheorem{Acon}[subsection]{Construction}

\theoremstyle{kursiv}
\newtheorem{lemdef} [subsubsection]{Lemma and Definition}
 
\newtheorem{thm}[subsubsection]{Theorem}
\newtheorem{prop} [subsubsection]{Proposition}
\newtheorem{cor} [subsubsection]{Corollary}
\newtheorem{lem} [subsubsection]{Lemma}

\newtheorem{Alem}[subsection]{Lemma}

\newtheorem{Anota}[subsection]{Notation}

\newtheorem{Acor}[subsection]{Corollary}

\numberwithin{equation}{subsection}

\setdefaultenum{\normalfont (a)}{i.}{A.}{1.}

\newcommand{\NN}{\ensuremath{\mathbb{N}}}

\newcommand{\RR}{\ensuremath{\mathbb{R}}}
\newcommand{\SSS}{\ensuremath{\mathbb{S}}}

\newcommand{\ZZ}{\ensuremath{\mathbb{Z}}}

\newcommand{\la}{\ensuremath{\lambda}}

\newcommand{\vp}{\ensuremath{\varphi}}
\newcommand{\ve}{\ensuremath{\varepsilon}}

\newcommand{\aA}{\ensuremath{\mathcal{A}}}
\newcommand{\bB}{\ensuremath{\mathcal{B}}}
\newcommand{\cC}{\ensuremath{\mathcal{C}}}
\newcommand{\dD}{\ensuremath{\mathcal{D}}}
\newcommand{\eE}{\ensuremath{\mathcal{E}}}
\newcommand{\fF}{\ensuremath{\mathcal{F}}}
\newcommand{\gG}{\ensuremath{\mathcal{G}}}
\newcommand{\hH}{\ensuremath{\mathcal{H}}}
\newcommand{\iI}{\ensuremath{\mathcal{I}}}

\newcommand{\kK}{\ensuremath{\mathcal{K}}}
\newcommand{\lL}{\ensuremath{\mathcal{L}}}
\newcommand{\mM}{\ensuremath{\mathcal{M}}}
\newcommand{\nN}{\ensuremath{\mathcal{N}}}
\newcommand{\oO}{\ensuremath{\mathcal{O}}}
\newcommand{\pP}{\ensuremath{\mathcal{P}}}

\newcommand{\rR}{\ensuremath{\mathcal{R}}}
\newcommand{\sS}{\ensuremath{\mathcal{S}}}
\newcommand{\tT}{\ensuremath{\mathcal{T}}}
\newcommand{\uU}{\ensuremath{\mathcal{U}}}
\newcommand{\vV}{\ensuremath{\mathcal{V}}}
\newcommand{\wW}{\ensuremath{\mathcal{W}}}
\newcommand{\xX}{\ensuremath{\mathcal{X}}}

\newcommand{\zZ}{\ensuremath{\mathcal{Z}}}

\DeclareMathAlphabet{\mathpzc}{OT1}{pzc}{m}{it}

\newcommand{\PSI}{\ensuremath{{\bf \Psi}}}
\DeclareMathOperator{\Ch}{\mathcal{C}\it{h}}
\newcommand{\CH}[2]{\Ch_{#1,#2}}

\newcommand{\set}[1]{\left\{ #1 \right\}} 
\newcommand{\setm}[2]{\left\{ #1 \, \middle\vert \, #2\right\}} 

\newcommand{\tl}{\textquotedblleft}
\newcommand{\tr}{\textquotedblright \ }

\newcommand{\coloneq}{\colonequals}
\newcommand{\norm}[1]{\left\lVert #1 \right\rVert}
\newcommand{\opnorm}[1]{\norm{#1}_\text{op}}
\newcommand{\abs}[1]{\left\lvert #1 \right\rvert}
\newcommand{\BoundOp}[1]{\mathcal{L}\left(#1\right)}

\DeclareMathOperator{\im}{Im}
\renewcommand{\Im}{\im}

\DeclareMathOperator{\Ob}{Ob}
\DeclareMathOperator{\supp}{supp}

\DeclareMathOperator*{\Diff}{Diff}

\DeclareMathOperator{\id}{id}
\DeclareMathOperator{\germ}{germ}
\DeclareMathOperator{\comp}{comp}

\DeclareMathOperator{\inv}{inv}

\DeclareMathOperator{\res}{res}

\DeclareMathOperator{\dom}{dom}
\DeclareMathOperator{\cod}{cod}

\DeclareMathOperator{\Gl}{Gl}

\DeclareMathOperator{\ORB}{{\bf Orb}}
\DeclareMathOperator{\proj}{pr}

\DeclareMathOperator{\evol}{evol}
\DeclareMathOperator{\Evol}{Evol}

\newcommand{\LB}[1][\cdot \hspace{1pt} , \cdot]{\left[\hspace{1pt} #1 \hspace{1pt} \right]}

\newcommand{\lPI}[1]{\pP(#1)}

\newcommand{\vect}[1]{\mathfrak{X} \left( #1 \right)}
\newcommand{\Os}[1]{\mathfrak{X}_{\text{Orb}} \left( #1 \right)}
\newcommand{\Osc}[1]{\mathfrak{X}_{\text{\rm Orb}} \left( #1 \right)_c}
\newcommand{\OsK}[1]{\mathfrak{X}_{\text{Orb}} \left( #1 \right)_K}

\DeclareMathOperator{\Fl}{Fl}

\DeclareMathOperator{\Man}{Man}

\DeclareMathOperator{\ORBI}{Orb}
\newcommand{\Orb}[1]{\ORBI (#1)}
\newcommand{\ORBM}[1][(Q_1,\uU_1),(Q_2,\uU_2)]{\ORB \left(#1\right)}

\newcommand{\ido}[1]{\id_{(Q_{#1}, \uU_{#1})}}
\newcommand{\tpi}[1][Q,\uU]{\pi_{\tT(#1)}}

\newcommand{\Difforb}[1]{\textstyle{\Diff_{\text{\rm Orb}}}\left(  #1 \right)}
\newcommand{\Diffo}[1]{\Difforb{  #1 }_0}
\newcommand{\Diffc}[1]{\textstyle{\Diff_{\text{\rm Orb}}}\left(  #1 \right)_c}
\newcommand{\DiffK}[1]{\textstyle{\Diff_{\text{\rm Orb}}}\left(  #1 \right)_K}

\newcommand{\expo}{\exp_{\text{\rm Orb}}}
\newcommand{\Exp}{\text{Exp}}

\newcommand{\xra}{\xrightarrow}

\DeclareMathOperator{\zs}{{\bf 0}_\text{\rm Orb}} 		
\newcommand{\one}{{\bf 1}}			

\newcommand{\cmt}[1]{}

\newcommand{\no}[1]{\mbox{}}
\newcommand{\todo}[1]{}

\newcommand{\ind}[2]{\ifthenelse{\equal{#1}{}}{\emph{#2} \index{#2}}{\emph{#2}\index{#1}}}
\makeindex

 \newglossarystyle{stybob}{
\glossarystyle{listdotted} 
\renewcommand*{\glossaryentryfield}[5]{%
\item[]%
\glstarget{##1}{##2}%
\unskip\leaders\hbox to 2.9mm{\hss.}\hfill\strut ##3%
\space ##5}
}

\glsdisablehyper

\newglossaryentry{Difforb}{name={$\Difforb{Q,\uU}$},description={}, sort={Difforb}}
\newglossaryentry{Diffo}{name={$\Diffo{Q,\uU}$},description={},sort={Diffo}}
\newglossaryentry{Diffc}{name={$\Diffc{Q,\uU}$},description={},sort={Diffc}}
\newglossaryentry{DiffK}{name={$\DiffK{Q,\uU}$},description={},sort={DiffK}}
\newglossaryentry{DiffG}{name={$\Diff^G (M)$},description={},sort={DiffG}}
\newglossaryentry{DiffGc}{name={$\Diff^G_c (M)$},description={},sort={DiffGc}}

\newglossaryentry{vect}{name={$\vect{M}$},description={},sort={vect}}
\newglossaryentry{Os}{name={$\Os{Q}$},description={},sort={Os}}
\newglossaryentry{Osc}{name={$\Osc{Q}$},description={},sort={Osc}}
\newglossaryentry{OsK}{name={$\OsK{Q}$},description={},sort={OsK}}
\newglossaryentry{vfloc}{name={$X_\psi$},description={},sort={vfloc}}
\newglossaryentry{XdiamondY}{name={$X \diamond Y $},description={},sort={XdiamondY}}
\newglossaryentry{Xstar}{name={$X^*$},description={},sort={Xstar}}
\newglossaryentry{LambdavV}{name={$\Lambda_{\vV}$},description={},sort={LambdavV}}
\newglossaryentry{sigman}{name={$\sigma_{[n]}$},description={},sort={sigman}}

\newglossaryentry{PSI}{name={$\PSI(\vV)$},description={},sort={PSI}}
\newglossaryentry{CH}{name={$\CH{U}{V}$},description={},sort={CH}}
\newglossaryentry{ch}{name={$\Ch_\vV$},description={},sort={ch}}
\newglossaryentry{lambdabar}{name={$\overline{\lambda}$},description={},sort={lambdabar}}

\newglossaryentry{dom}{name={$\dom f$},description={},sort={domain}}
\newglossaryentry{cod}{name={$\cod f$},description={},sort={codomain}}

\newglossaryentry{ido}{name={$\ido{}$},description={},sort={id}}
\newglossaryentry{pito}{name={$\tpi{}$},description={},sort={pito}}
\newglossaryentry{expo}{name={$[\expo]$},description={},sort={expo}}
\newglossaryentry{rest}{name={$\res_U^W$},description={},sort={restriction}}
\newglossaryentry{evol}{name={$\evol$}, description={}, sort={evol}}
\newglossaryentry{Evol}{name={$\Evol$}, description={}, sort={Evol}}
\newglossaryentry{fwed}{name={$f^\wedge$}, description={}, sort={fwedge}}
\newglossaryentry{zs}{name={$\zs$},description={},sort={zero-orbisection}}
\newglossaryentry{theta}{name={$\theta_\psi$},description={},sort={theta}}
\newglossaryentry{fwedge}{name={$f^\wedge$},description={},sort={fwedge}}

\newglossaryentry{Orb}{name={$\Orb{\vV,\wW}$},description={},sort={Orb}}
\newglossaryentry{ORB}{name={$\ORB ((Q_1,\uU_1),(Q_2,\uU_2))$},description={},sort={ORB}}

\newglossaryentry{otsp}{name={$\tT_p$},description={},sort={tspofd}}
\newglossaryentry{tofd}{name={$\tT Q$},description={},sort={tofd}}
\newglossaryentry{TU}{name={$\tT\uU$},description={},sort={TU}}

\newglossaryentry{fixed points}{name={$\Sigma_G$},description={},sort={fixed point setexpo}}
\newglossaryentry{singular stratum}{name={$\Sigma_Q$},description={},sort={fixed stratum}}
\newglossaryentry{local group}{name={$\Gamma_z (Q)$},description={},sort={local group}}

\newglossaryentry{ballrm}{name={$B_\rho (0_x,\ve)$},description={},sort={ballrm}}
\newglossaryentry{omega1}{name={$\Omega_{r,K}$},description={},sort={omega1}}
\newglossaryentry{omegaK5}{name={$\Omega_{r,K_5}$},description={},sort={omegaK5}}

\newglossaryentry{F5}{name={$\fF_5(K)$},description={},sort={F5K}}

\newglossaryentry{KU}{name={$\lfloor K, U\rfloor$},description={},sort={KU}}
\newglossaryentry{KUr}{name={$\lfloor K, U\rfloor_r$},description={},sort={KUr}}
\newglossaryentry{Crco}{name={$C(M,E)_{\text{c.o}}$},description={},sort={Crco}}
\newglossaryentry{normKr}{name={$\norm{\cdot}_{K,r}$},description={},sort={normKr}}

\newglossaryentry{initr}{name={$\lvert r]$},description={},sort={|r]}}
\newglossaryentry{geoseg}{name={$[\hat{c}]_{\mid [a,b]}$},description={},sort={geodesic segment}}

\newglossaryentry{embtype1}{name={$\uU \Subset \vV$},description={},sort={embtype1}}
\newglossaryentry{embtype2}{name={$\uU \sqsubset \vV$},description={},sort={embtype2}}
\newglossaryentry{intervall}{name={$\iI$},description={},sort={intervall}}
\newglossaryentry{Ny}{name={$N_y$},description={},sort={Ny}}

\newglossaryentry{fderiv}{name={$f'(x)$},description={},sort={fderivative}}

\newglossaryentry{Fl}{name={$\Fl^f$},description={},sort={Flow}}
\newglossaryentry{Fldom}{name={$\mathfrak{D} (\xi)$},description={},sort={domain of the Flow}}

\newglossaryentry{partCr}{name={$C^{(r,s)}$},description={},sort={C partially r,s}}

\newglossaryentry{K5}{name={$K_5$},description={},sort={K5}}

\newglossaryentry{linmor}{name={$\lL (\RR^d)$},description={},sort={linear Endomorphisms}}

\newglossaryentry{lld}{name={$\delta^\ell$}, description={}, sort={left logarithmic derivative}}
\newglossaryentry{rld}{name={$\delta^r$}, description={}, sort={right logarithmic derivative}}
\newglossaryentry{lPI}{name={$\lPI{\gamma}$}, description={}, sort={left product integral}}

\title{The Diffeomorphism Group of a Non-Compact Orbifold}
\author{Alexander Schmeding\thanks{\small Mathematical Institute, University of Paderborn, Warburger Straße 100, 33098 Paderborn, Germany, \href{mailto:alsch@math.upb.de}{alsch@math.upb.de}}}
\date{}
\begin{document}

\maketitle
\thispagestyle{empty}
\begin{abstract}
We endow the diffeomorphism group $\Difforb{Q,\uU}$ of a paracompact (reduced) orbifold with the structure of an infinite dimensional Lie group modeled on the space of compactly supported sections of the tangent orbibundle. For a second countable orbifold, we prove that $\Difforb{Q,\uU}$ is $C^0$-regular and thus regular in the sense of Milnor. Furthermore an explicit characterization of the Lie algebra of $\Difforb{Q,\uU}$ is given. 
\end{abstract}

\tableofcontents 

\newpage
\section*{Introduction and Statement of Results}\addcontentsline{toc}{section}{Introduction and statement of results}

Diffeomorphism groups of compact manifolds and their subgroups are prime examples of infinite dimensional Lie groups. There are many well known results concerning the Lie group structure of these groups; e.g., a classical result states that the diffeomorphism group of a compact manifold is an infinite dimensional regular Lie group (see \cite{milnor1983}). For the algebraic structure of these groups, see \cite{ban1997}. More generally, Lie group structures on diffeomorphism groups of paracompact manifolds (even with corners) were constructed in \cite{michor1980} (also cf.\ \cite{hg2005} for the special case $\Diff (\RR^n)$). Furthermore, in \cite{conv1997} the diffeomorphism groups of manifolds were endowed with the structure of a regular Lie group in the \tl convenient setting of analysis\textquotedblright. We remark that the \tl convenient setting of analysis\tr (see \cite{conv1997}) is inequivalent to the setting of analysis adopted in this paper. Our studies are based on a concept of $C^r$-maps between locally convex spaces known as Keller's $C^r_c$-theory~\cite{keller1974} (see \cite{milnor1983}, \cite{hg2002a} and \cite{gn2007} for streamlined expositions, cf.\ also \cite{bgn2004}). The present paper generalizes the results on diffeomorphism groups of manifolds to diffeomorphism groups of reduced paracompact orbifolds. 

Orbifolds were first introduced by Satake in \cite{sat1956} as \textit{$V$-manifolds} to generalize the concept of a manifold. Later on they appear in the works of Thurston (cf.\ \cite{thur2002}), who popularized the term \tl orbifold\textquotedblright . One might think of an orbifold as a manifold with \tl mild singularities\textquotedblright. Objects with orbifold structure arise naturally, for example in symplectic geometry, physics and algebraic geometry (cf.\ the survey in \cite{alr2007}). It is well known that there are at least three different ways to define an orbifold: Orbifolds may be described by atlases of local charts akin to a manifold (see \cite{alr2007, follie2003,cgo2006}). Furthermore, orbifolds correspond to special classes of Lie groupoids (see \cite{follie2003} or the survey \cite{moer2002}). Finally one might think of them as Deligne-Mumford stacks (cf.\ \cite{ofdstack2008}). The author thinks that the first approach is suited best to apply methods from differential geometry to orbifolds. Hence in the present paper we define orbifolds in local charts. Unfortunately, this point of view makes it difficult to define morphisms of orbifolds. The literature proposes a variety of notions for these morphisms, e.g.\ the Chen-Ruan good map \cite{chen2001}, the Moerdijk-Pronk strong map \cite{osg1997}, or the maps in \cite{bb2008}. However, orbifolds in local charts are equivalent to certain Lie groupoids, whose morphisms are well understood objects. Thus orbifold morphisms should correspond to a class of Lie groupoid morphisms. The orbifold maps introduced by Pohl in \cite{pohl2010} satisfy these requirements, since they were modeled to be equivalent to groupoid morphisms.\footnote{Other concepts of orbifold maps are also widely believed to satisfy similar properties, cf. \cite[Section 2.4]{alr2007}. However in \cite{pohl2010} a counterexample to these claims may be found.} Furthermore these maps allow a characterization in local charts, which is amenable to methods of differential geometry and Lie theory. Therefore in the present paper, maps of orbifolds will be orbifold maps in the sense of Pohl \cite{pohl2010} (for a comprehensive introduction to these maps see Appendix \ref{sect: maps}). \\
To construct the Lie group structure on the diffeomorphism group of an orbifold we have to develop several tools from Riemannian geometry on orbifolds. These results are of interest in their own right and include the following:\\
We discuss geodesics on Riemannian orbifolds and prove that they are uniquely determined by their initial values. Then a detailed construction for a \emph{Riemannian orbifold exponential map} $[\expo]$ is provided. This map is an orbifold morphism in the sense of Pohl \cite{pohl2010}, which generalizes the concept of a Riemannian exponential map to Riemannian orbifolds (cf.\ \cite{cgo2006} and \cite{chen2001}, respectively for Riemannian exponential maps on geodesically complete orbifolds).

The Riemannian exponential map on a manifold may be used to construct the Lie group structure on the diffeomorphism group of the manifold (cf.  \cite{milnor1983}). The Riemannian orbifold exponential map allows us to follow this line of thought: We endow the diffeomorphism group of a paracompact reduced orbifold with the structure of an infinite dimensional locally convex Lie group in the sense of \cite{neeb2006}. More precisely the main results subsume the following theorem (cf. Theorem \ref{thm: diff:lgp}): 
\paragraph{Theorem A} \emph{The diffeomorphism group $\Difforb{Q,\uU}$ of a paracompact reduced orbifold $(Q,\uU)$ can be made into a Lie group in a unique way such that the following is satisfied:\\
 For some Riemannian orbifold metric $\rho$ on $(Q,\uU)$, let $[\expo]$ be the Riemannian orbifold exponential map. There exists an open zero-neighborhood $\hH_\rho$ in the space of compactly supported sections of the tangent orbibundle such that 
  \begin{displaymath}
   E \colon \hH_\rho \rightarrow \Difforb{Q,\uU}, [\hat{\sigma}] \mapsto [\expo] \circ [\hat{\sigma}]
  \end{displaymath}
 induces a well defined $C^\infty$-diffeomorphism onto an open submanifold of $\Difforb{Q,\uU}$. This condition is then satisfied for every Riemannian orbifold metric on $(Q,\uU)$. If $(Q,\uU)$ is a compact orbifold, then the Lie group $\Difforb{Q,\uU}$ is a Fr\'{e}chet-Lie group.}\\[2em]
 This result generalizes the classical construction of a Lie group structure on the diffeomorphism group $\Diff (M)$ of a paracompact manifold. For such a manifold, we may consider subgroups of $\Diff (M)$, whose elements coincide outside of a given compact set with the identity. It is known that these subgroups are Lie subgroups of $\Diff (M)$ (cf. \cite[Section 14]{hg2004}). Section \ref{sect: lgp:glob} contains a similar result for diffeomorphisms of orbifolds, which is a consequence of Theorem A: 
 \paragraph{Theorem B} \emph{Let $(Q,\uU)$ be a paracompact reduced orbifold. For each compact subset $K$ of $Q$ we define the group $\DiffK{Q,\uU}$ of all orbifold diffeomorphisms which coincide off $K$ with the identity morphism of the orbifold. Let $\Diffc{Q,\uU}$ be the group of all orbifold diffeomorphisms which coincide off some compact set with the identity morphism of the orbifold. Then the following holds: 
  \begin{compactenum} 
   \item The group $\Diffc{Q,\uU}$ is an open normal Lie subgroup of $\Difforb{Q,\uU}$.
   \item  For each compact subset $K$ of $Q$, there is a compact set $L \supseteq K$ such that $\Difforb{Q,\uU}_L$ is a closed Lie subgroup of $\Difforb{Q,\uU}$. The closed Lie subgroup  $\Difforb{Q,\uU}_L$ is modeled on the space of sections in the tangent orbibundle which vanish off $L$.
  \end{compactenum}
 If $(Q,\uU)$ is a trivial orbifold (i.e.\ a manifold), one may always choose $K=L$ in \textrm{(b)}.}\\[1em]
 We remark that Lie group structures for diffeomorphism groups of orbifolds were already considered by Borzellino and Brunsden. In \cite{bb2008} and the follow up \cite{bb2009}, the diffeomorphism group of a compact orbifold has been turned into a convenient Fr\'{e}chet-Lie group. The author does not know whether the orbifold morphisms introduced in \cite{bb2008} are equivalent to the class of orbifold maps considered in the present paper. If both notions were equivalent, the results of \cite{bb2008,bb2009} concerning the Lie group structure of the diffeomorphism group are subsumed in Theorem A. This follows from the fact that in the Fr\'{e}chet setting both notions of \tl smooth maps\tr coincide (cf. \cite{keller1974} and \cite[Theorem 4.11 (a)]{conv1997}). Hence Fr\'{e}chet Lie groups in the sense of  \cite{neeb2006} and \tl convenient Fr\'{e}chet Lie groups\tr coincide. However, we have to point out that the exposition in \cite{bb2008} contains several major errors (see Remark \ref{rem: bor} for further information on this topic). \\
 We also mention that in the groupoid setting, topologies for spaces of orbifold maps have been considered. Chen constructs in \cite{chen2006} a topology on the space of orbifold morphisms whose domain is a compact orbifold, turning the space into a Banach orbifold (also cf. similar results in \cite{haefliger2008}). The exposition of the present paper is independent of these results.
 
 After constructing the Lie group $\Difforb{Q,\uU}$, we determine the Lie algebra associated to this group. It is instructive to recall the special case of the diffeomorphism group $\Diff (M)$ of a compact manifold $M$. Milnor proves in \cite{milnor1983} that the Lie algebra associated to $\Diff (M)$ is the space of vector fields $\vect{M}$ on $M$, whose Lie bracket is the negative of the bracket product of vector fields. It turns out that an analogous result holds for the Lie algebra of the Lie group $\Difforb{Q,\uU}$. To understand the result we need the following facts:\\ 
 A map of orbifolds $[\hat{\sigma}]$, which is a section of the tangent orbibundle is called an orbisection. With respect to an orbifold chart of $Q$, each orbisection induces a unique vector field on the chart domain, called its canonical lift. In particular, each orbisection corresponds to a unique family of vector fields (cf. Section \ref{sect: tofd} for details). By construction, the local model for the Lie group $\Difforb{Q,\uU}$ is the space of compactly supported orbisections $\Osc{Q}$. We are now in a position to formulate the following result on the Lie algebra of the diffeomorphism group $\Difforb{Q,\uU}$ (Theorem \ref{thm: LALG}): 
\paragraph{Theorem C} \emph{The Lie algebra of $\Difforb{Q,\uU}$ is given by $(\Osc{Q}, \LB[\cdot,\cdot])$. Here the Lie bracket $\LB[\cdot,\cdot]$ is defined as follows:\\ For arbitrary $[\hat{\sigma}], [\hat{\tau}] \in \Osc{Q}$, their Lie bracket $\LB[[\hat{\sigma}], [\hat{\tau}]]$ is the unique compactly supported orbisection whose canonical lift on an orbifold chart $(U,G,\varphi)$ is the negative of the Lie bracket in $\vect{U}$ of their canonical lifts $\sigma_U$ and $\tau_U$.
}\\[2em]
Finally we discuss regularity properties of the Lie group $\Difforb{Q,\uU}$. To this end, recall the notion of regularity for Lie groups:\\
Let $G$ be a Lie group modeled on a locally convex space, with identity element $\one$, and $r\in \NN_0\cup\{\infty\}$. We use the tangent map of the right translation $\rho_g\colon G\to G$, $x\mapsto xg$ by $g\in G$ to define $v.g\coloneq T_{\one} \rho_g(v) \in T_g G$ for $v\in T_{\one} (G) =: L(G)$. Following \cite{Dahmen2012}, \cite{hg2015} and \cite{gn2007}, $G$ is called \emph{$C^r$-regular} if the initial value problem
\begin{displaymath}
\begin{cases}
\eta'(t)&= \gamma(t).\eta(t)\\
\eta(0) &= \one
\end{cases}
\end{displaymath}
has a (necessarily unique) $C^{r+1}$-solution $\Evol (\gamma):=\eta\colon [0,1]\rightarrow G$ for each $C^r$-curve $\gamma\colon [0,1]\rightarrow L(G)$, and the map
\begin{displaymath}
 \evol \colon C^r([0,1],L(G))\rightarrow G,\quad \gamma\mapsto \Evol (\gamma)(1)
\end{displaymath}
is smooth. If $G$ is $C^r$-regular and $r\leq s$, then $G$ is also $C^s$-regular. A $C^\infty$-regular Lie group $G$ is called \emph{regular} \emph{(in the sense of Milnor}) -- a property first defined in \cite{milnor1983}. Every finite dimensional Lie group is $C^0$-regular (cf. \cite{neeb2006}). Several important results in infinite-dimensional Lie theory are only available for regular Lie groups (see \cite{milnor1983}, \cite{neeb2006}, \cite{hg2015}, cf.\ also \cite{conv1997} and the references in these works).
We prove the following result (Theorem \ref{thm: diff:reg}):  
\paragraph{Theorem D}\emph{For a second countable orbifold, the Lie group $\Difforb{Q,\uU}$ is $C^k$-regular for each $k \in \NN_0 \cup \set{\infty}$. In particular this Lie group is regular in the sense of Milnor.}\\[1em]
Notice that in general the orbifolds in the present paper are not assumed to be second countable. However our methods require second countability of the orbifold to prove that the evolution map $\evol$ is smooth. It is known that the approach outlined in the present paper may not be adapted to orbifolds which are not second countable. Hence we pose the following question:\\[1em] 
\textbf{Open Problem:} Let $(Q,\uU)$ be a paracompact reduced orbifold which is not second countable. Is the Lie group $\Difforb{Q,\uU}$ a $C^r$-regular Lie group for some $r \in \NN_0 \cup \set{\infty}$? 

The present article commences with a brief introduction to infinite dimensional calculus, orbifolds and their properties (Section \ref{sect: prelim}). Our goal is to present a mostly self contained exposition of orbifolds and their morphisms. In particular, Appendix E contains all necessary information about orbifold maps in the sense of \cite{pohl2010}. However, the exposition avoids references to the groupoid morphisms after which these maps are modeled. The thesis is organized as follows:\\
In Sections \ref{sect: diffeo} and \ref{sect: tofd} classes of orbifold maps are discussed in the setting of \cite{pohl2010}. These include orbifold diffeomorphisms, partitions of unity and sections of the tangent orbibundle. Afterwards, we consider Riemannian geometry on orbifolds and develop important tools employed in the proof of the central results of this work. The main results of the thesis are contained in Section \ref{sect: lgp:str}. As an application we consider groups of equivariant diffeomorphisms $\Diff^G (\RR^n)$ associated to certain good orbifolds (i.e. orbifolds with a global chart). The Lie group structures obtained in these examples correspond to certain closed Lie subgroups of $\Diff (\RR^n)$ (considered as the Lie group constructed in \cite{hg2005}).\\
The less introductory material contained in the appendices should be taken on faith on a first reading. The presentation of this material in the text would have distracted from the main line of thought.

\bigskip
\textbf{Remark} (added in 2015) The open problem formulated in this introduction has been solved. 
A combination of \cite[Corollary 13.6]{hg2015} and the results contained in Section \ref{sect: regularity} shows that for a non-second countable orbifold the diffeomorphism group will be $C^1$-regular. 
\newpage
\pagestyle{fancy}
\thispagestyle{empty}
\section{Preliminaries and Notation}\label{sect: prelim}

\begin{conv}
 In this thesis, we work exclusively over the field $\RR$ of real numbers. All topological spaces will be assumed to be Hausdorff. We write $\NN \coloneq \set{1,2,\ldots}$ and $\NN_0 \coloneq \NN \cup \set{0}$.
\end{conv}

\subsection{Differential calculus in infinite dimensional spaces}

Basic references for differential calculus in locally convex spaces are \cite{hg2002a,keller1974,bgn2004, hg2002,hg2006}. Basic facts on infinite dimensional manifolds are compiled in Appendix \ref{app: mfd}. For the reader's convenience, we recall various definitions and results:

\begin{defn}\label{defn: deriv} \no{defn: deriv}
 Let $E, F$ be locally convex spaces, $U \subseteq E$ be an open subset, $f \colon U \rightarrow F$ a map and $r \in \NN_{0} \cup \set{\infty}$. If it exists, we define for $(x,h) \in U \times E$ the directional derivative $df(x,h) \coloneq D_h f(x) \coloneq \lim_{t\rightarrow 0} t^{-1} (f(x+th) -f(x))$. We say that $f$ is $C^r$ if the iterated directional derivatives
    \begin{displaymath}
     d^{(k)}f (x,y_1,\ldots , y_k) \coloneq (D_{y_k} D_{y_{k-1}} \cdots D_{y_1} f) (x)
    \end{displaymath}
 exist for all $k \in \NN_0$ such that $k \leq r$, $x \in U$ and $y_1,\ldots , y_k \in E$ and define continuous maps $d^{(k)} f \colon U \times E^k \rightarrow F$. If $f$ is $C^\infty$ it is also called smooth. We abbreviate $df \coloneq d^{(1)} f$.
\end{defn}

\begin{rem}
 If $E_1,E_2,F$ are locally convex spaces and $U \subseteq E_1,V \subseteq E_2$ open subsets together with a $C^1$-map $f \colon U \times V \rightarrow F$, then one may compute the \ind{}{partial derivative} $d_1f$ with respect to $E_1$. It is defined as $d_1 f \colon U \times V \times E_1 \rightarrow F, d_1f(x,y;z) \coloneq \lim_{t \rightarrow 0} t^{-1} (f(x+tz,y)-f(x,y))$. Analogously one defines the partial derivative $d_2 f$ with respect to $E_2$. The linearity of $df(x,y,\cdot)$ implies the so-called Rule on Partial Differentials for $(x,y) \in U\times V, (h_1,h_2) \in E_1\times E_2$:
  \begin{equation}\label{eq: rule:pdiff}
   df(x,y,h_1,h_2) = d_1 f(x,y;h_1) + d_2 f(x,y;h_2).
  \end{equation}
 By \cite[Lemma 1.10]{hg2002a}, $f \colon U \times V \rightarrow F$ is $C^1$ if and only if $d_{1}f$ and $d_{2}f$ exist and are continuous.
\end{rem}

\begin{defn}[Differentials on non-open sets]\mbox{}
 \begin{compactenum}
  \item The set $U \subseteq E$ is called \emph{locally convex} if every $x \in U$ has a convex neighborhood $V$ in $U$.
  \item Let $U\subseteq E$ be a locally convex subset with dense interior. A continuous mapping $f \colon U \rightarrow F$ is called $C^r$ if $f|_{U^\circ} \colon U^\circ \rightarrow F$ is $C^r$ and each of the maps $d^{(k)} (f|_{U^\circ}) \colon U^\circ \times E^k \rightarrow F$ admits a (unique) continuous extension $d^{(k)}f \colon U \times E^k \rightarrow F$. If $U \subseteq \RR$ and $f$ is $C^{1}$, we obtain a continuous map $f' \colon U \rightarrow E, f'(x) \coloneq df(x)(1)$. We shall write $\frac{\partial}{\partial x}f(x) \coloneq f' (x)$. In particular if $f$ is of class $C^r$, we define recursively $\frac{\partial^k}{\partial x^k}f (x) = (\frac{\partial^{k-1}}{\partial x^{k-1}}f)'(x)$ for $k \in \NN_0$ such that $k \leq r$, where $f^{(0)} \coloneq f$.   
 \end{compactenum}
\end{defn}

Using these definitions one may define infinite dimensional manifolds as usual. We refer to Appendix \ref{app: mfd} for definitions and comments on the notation used. To discuss regularity properties of Lie groups, the notion of $C^{r,s}$-mappings is useful. 
  \begin{defn}[$C^{r,s}$-mappings]
   Let $E_1$, $E_2$ and $F$ be locally convex spaces, $U$ and $V$ open subsets of $E_1$ and $E_2$, respectively, and $r,s \in \NN_0 \cup \set{\infty}$. A mapping $f\colon U \times V \rightarrow F$ is called a $C^{r,s}$-map if for all $i,j \in \NN_0$ such that $ i \leq r, j \leq s$, the iterated directional derivative 
   \begin{displaymath}
    d^{(i,j)}f(x,y,w_1,\dots,w_i,v_1,\dots,v_j) \coloneq (D_{(w_i,0)} \cdots D_{(w_1,0)}D_{(0,v_j)} \cdots D_{(0,v_1)}f ) (x,y)
   \end{displaymath}
exists for all $ x \in U, y \in V, w_1, \ldots , w_i \in E_1,  v_1, \ldots ,v_j \in E_2$ and yields continuous mappings
 \begin{align*} 
  d^{(i,j)}f\colon   & U \times V \times E^i_1 \times E^j_2 \rightarrow F,\\ 
                     &(x,y,w_1,\dots,w_i,v_1,\dots,v_j)\mapsto (D_{(w_i,0)} \cdots D_{(w_1,0)}D_{(0,v_j)} \cdots D_{(0,v_1)}f ) (x,y).
 \end{align*}
  \end{defn}

Again this concept may be extended to maps on non-open domains with dense interior: 
\begin{defn}\label{defn: crs} \no{defn: crs}
 Let $E_1$, $E_2$ and $F$ be locally convex spaces. Consider locally convex subsets with dense interior $U$ of $E_1$ and $V$ of $E_2,$ and $r,s \in \NN_0 \cup \set{\infty}$. We say that a continuous map $f\colon U \times V \rightarrow F $ is a $C^{r,s}$-map, if $f|_{U^\circ \times V^\circ } \colon {U^\circ \times V^\circ}\rightarrow F$ is a $C^{r,s}$-map and for all $i,j \in \NN_0$ such that $ i \leq r, j \leq s $, the map 
  \begin{displaymath}
   d^{(i,j)}(f\lvert_{U^\circ \times V^\circ}) \colon U^\circ \times V^\circ \times E^i_1 \times E^j_2 \rightarrow F 
  \end{displaymath}
 admits a continuous extension $d^{(i,j)}f\colon U \times V \times E^i_1 \times E^j_2 \rightarrow F$.
\end{defn}

For further results and details on the calculus of $C^{r,s}$-maps we refer to \cite{alas2012}. 

\begin{defn}
 Let $U,V$ be locally convex subsets with dense interior of locally convex spaces $E_1$ and $E_2$, respectively, and let $F$ be a locally convex space. For $r,s \in \NN_0 \cup \set{\infty}$, we define the spaces 
	\begin{align*}
	 C^r (U,F) &\coloneq \setm{f \colon U \rightarrow F}{f \text{ is a mapping of class } C^r}\\
         C^{r,s} (U\times V , F) &\coloneq \setm{f \colon U \times V \rightarrow F}{f \text{ is a mapping of class } C^{r,s}}.
	\end{align*}
 Furthermore, we define $C (U,F) \coloneq C^0 (U,F)$ and endow $C^r(U,F)$ with the compact-open $C^r$-topology (see Section \ref{sect: Crtop})
\end{defn}

In the following, we let $\Diff^r (M)$ be the group of $C^r$-diffeomorphisms from a $C^r$-manifold $M$ to itself for $r \in \NN_0 \cup \set{\infty}$. To shorten the notation, we write $\Diff (M) \coloneq \Diff^\infty (M)$ if $M$ is a smooth manifold.

\subsection{Orbifolds I: Moerdijk's definition}\label{sect: moerdijk}
In this section, we introduce orbifolds as in the works of Moerdijk et al. Our exposition follows \cite{follie2003}, but we slightly change the definition of orbifold charts (see Remark \ref{rem: ofd:charts}).

\begin{defn}[Orbifold charts] \label{defn: moerd:ofc} \no{defn: moerd:ofc}
Let $Q$ be a topological space. An \ind{}{orbifold chart} of dimension $n \geq 0$ is a triple $(U, G, \phi)$, where $U$ is a connected smooth paracompact $n$-dimensional manifold without boundary, $G$ is a finite subgroup of $\Diff (U)$ and $\phi \colon U \rightarrow Q$ is an open map which factors to a homeomorphism on the orbit space $U/G \rightarrow \phi (U)$.\\[1em]
If $(U,G,\phi )$ is an orbifold chart on $Q$ and $S$ an open $G$-stable subset of $U$, then $\set{g|_S \colon g \in G_S}$ is a group isomorphic to $G_S$ by Newman's Theorem \ref{thm: newman}. Thus by abuse of notation the triple $(S, G_S, \phi|_{S})$ is again an orbifold chart called the \ind{orbifold chart!restriction of}{restriction of} $(U,G,\phi )$ on $S$.\\ Let $(V,H, \psi )$ be another orbifold chart on $Q$. An \ind{orbifold chart!embedding}{embedding $\lambda \colon (V,H, \psi ) \rightarrow (U,G, \phi )$ of orbifold charts} is a topological embedding $\lambda \colon V \rightarrow U$ which is an \'{e}tale map\footnote{i.e.\ for each $p$ in the domain of $\lambda$, the tangent map $T_p \lambda$ is an isomorphism. On occasion these maps will also be called local diffeomorphisms.} that satisfies $\phi \circ \lambda = \psi$. \\[1em]
We say that two orbifold charts $(U,G,\phi )$ and $(V,H,\psi )$ of dimension $n$ on $Q$ are \ind{orbifold chart!compatibility}{compatible} if for any $z \in \phi (U) \cap \psi (V)$, there exist an orbifold chart $(W, K, \theta )$ on $Q$ with $z \in \theta (W)$ and embeddings between orbifold charts $\lambda \colon (W,K, \theta ) \rightarrow (U,G, \phi )$ and $\mu \colon (W,K, \theta ) \rightarrow (V,H, \psi )$.
\end{defn}


\begin{prop}[\hspace{-0.5pt}{\cite[Proposition 2.12]{follie2003}}]\label{prop: ch:prop} \no{prop: ch:prop}
Let $Q$ be a topological space.
\begin{compactenum}
 \item For any embedding $\lambda \colon (V, H, \psi ) \rightarrow (U, G, \phi )$ between orbifold charts on $Q$, the image $\lambda (V)$ is a $G$-stable open subset of $U$, and there is a unique isomorphism $\overline{\lambda} \colon H \rightarrow G_{\lambda (V)} \leq G$ for which $\lambda (hx) = \overline{\lambda} (h) \lambda (x)$.\glsadd{lambdabar}
 \item The composition of two embeddings between orbifold charts is an embedding between orbifold charts. 
 \item For any orbifold chart $(U, G, \phi )$, any diffeomorphism $g \in G$ is an embedding of $(U,G,\phi )$ into itself, and $\overline{g}(g') = gg'g^{-1}$.
 \item If $\la , \mu \colon (V,H,\phi ) \rightarrow (U,G, \phi )$ are two embeddings between the same orbifold charts, there exists a unique $g \in G$ with $\lambda = g \circ \mu$.
\end{compactenum}
\end{prop}
\begin{proof} The proof for \cite[Proposition 2.12]{follie2003} carries over verbatim to finite dimensional connected manifolds without boundary.
\end{proof}

\begin{defn}[Orbifolds I]\label{defn: moer:ofdI} \no{defn: moer:ofdI}
An \ind{}{orbifold atlas} of dimension $n$ for a topological space $Q$ is a set of pairwise compatible orbifold charts 
	\begin{displaymath}
		\uU \coloneq \setm{(U_i, G_i, \phi_i)}{i\in I}
	\end{displaymath}
of dimension $n$ on $Q$ such that $\bigcup_{i\in I} \phi_i (U_i ) = Q$. Two orbifold atlases of $Q$ are equivalent if their union is an orbifold atlas. An \ind{}{orbifold} of dimension $n$ is a pair $(Q, \uU )$, where $Q$ is a paracompact Hausdorff topological space and $\uU$ is an equivalence class of orbifold atlases of dimension $n$ on $Q$. 
\end{defn}

\begin{rem}\label{rem: ofd:charts}
 The definition of an orbifold does not exactly follow the exposition in \cite{follie2003}. We have to mention two changes: 
  \begin{compactenum}
   \item For an orbifold chart $(U,G,\pi)$ as defined in this section, the chart domain $U$ is a finite dimensional connected and paracompact manifold. In \cite{follie2003} one is only allowed to choose $U$ as an open subset of $\RR^n$. However, every orbifold in our sense uniquely determines one in the sense of \cite{follie2003}. This fact follows from Lemma \ref{lem: lin:ch}: Let $(U,G,\pi)$ be an orbifold chart as in Definition \ref{defn: moerd:ofc}. Then Lemma \ref{lem: lin:ch} allows the construction of an orbifold chart $(V_x,G_{V_x},\pi|_{V_x})$ for $x \in U$, where $V_x$ is diffeomorphic to an open subset of $\RR^n$. Hence the orbifolds defined in Definition \ref{defn: moer:ofdI} admit an orbifold atlas whose chart domains are open subsets of $\RR^n$.  
   \item Contrary to the treatment in \cite{follie2003}, we do not require the topological space $Q$ to be second countable. We do not need second countability of $Q$ for most of this work, whence we chose to omit it here (also compare Remark \ref{rem: psgp}).
  \end{compactenum}

\end{rem}

\subsection{Orbifolds II: Haefliger's definition}\label{sect: haef} \no{sect: haef}

We recall an equivalent definition of orbifolds as outlined in \cite{cgo2006}: 

\begin{defn}[Orbifolds II, {\cite{cgo2006}}]\label{defn: haef:ofdII} \no{defn: haef:ofdII}
 Let $Q$ be a paracompact Hausdorff topological space.
\begin{compactenum}
 \item Let $n$ be in $\NN_0$. A \ind{orbifold chart!reduced}{(reduced) orbifold chart} of dimension $n$ on $Q$ is a triple $(V,G,\varphi )$ where $V$ is a connected paracompact $n$-dimensional manifold without boundary, $G$ is a finite subgroup of $\Diff (V)$, and $\varphi \colon V \rightarrow Q$ is a map with open image $\varphi (V)$ that induces a homeomorphism from $V/G$ to $\varphi (V)$. In this case, $(V,G, \varphi )$ is said to \ind{set!uniformized}{uniformize} $\varphi (V)$.
 \item Two reduced orbifold charts $(V,G, \varphi )$, $(W,H,\psi )$ on $Q$ are called \emph{compatible} if for each pair $(x,y) \in V \times W$ with $\varphi (x) = \psi (y)$ there are open connected neighborhoods $V_x$ of $x$ and $W_y$ of $y$ and a $C^\infty$-diffeomorphism $h \colon V_x \rightarrow W_y$ such that $\psi \circ h = \varphi|_{V_x}$. The map $h$ is called a \ind{orbifold chart!change of charts}{change of charts}.
 \item A \ind{orbifold atlas!reduced}{reduced orbifold atlas} of dimension $n$ on $Q$ is a set of pairwise compatible reduced orbifold charts
	\begin{displaymath}
             \vV \coloneq \setm{(V_i, G_i , \varphi_i )}{ i \in I}
  	\end{displaymath}
 of dimension $n$ on $Q$ such that $\bigcup_{i\in I} \varphi_i (V_i) = Q$.
 \item Two reduced orbifold atlases are \ind{orbifold atlas!equivalence}{equivalent} if their union is a reduced orbifold atlas.
 \item A \ind{}{reduced orbifold structure} of dimension $n$ on $Q$ is an equivalence class of reduced orbifold atlases of dimension $n$ on $Q$.
 \item A \ind{orbifold!reduced}{reduced orbifold} of dimension $n$ is a pair $(Q,\uU )$ where $\uU$ is a reduced orbifold structure of dimension $n$ on $Q$.
\end{compactenum}
\end{defn}

The Definition \ref{defn: haef:ofdII} is equivalent to the Definition \ref{defn: moer:ofdI}, i.e. they yield the same equivalence classes of orbifold atlases. The compatibility conditions of both definitions coincide by \cite[Proposition 2.13]{follie2003}. The proof outlined in \cite{follie2003} carries over without any changes to our setting. 

\begin{rem}
\begin{compactenum}
 \item The term \tl reduced\tr refers to the requirement that for each reduced orbifold chart $(V, G, \varphi )$ in $\uU$ the group $G$ is a subgroup of $\Diff (V)$. Hence the action of $G$ on $V$ is effective. We will only consider reduced orbifolds (and maps between them). Thus to shorten our notation, we will drop the term \tl reduced\tr in the remainder of the paper. A \tl reduced\tr orbifold will thus simply be called an orbifold.  
 \item We will occasionally refer to the dimension of an orbifold as defined in \ref{defn: haef:ofdII} as the \ind{orbifold!dimension}{orbifold dimension}. We shall prove later that, as in the case of a manifold, the orbifold dimension is an invariant of the orbifold. More explicitly two orbifolds can only be diffeomorphic to each other if they have the same orbifold dimension. We postpone these considerations until we are ready to define morphisms of orbifolds.
 \item In general, maps of orbifolds (see Appendix \ref{sect: maps}) only admit local lifts in certain orbifold atlases contained in the equivalence class $\uU$ of the orbifold $(Q,\uU)$. Therefore we introduce the convention: An atlas $\vV$ contained in $\uU$ will be called a \ind{orbifold atlas!representative}{representative of $\uU$}.
 \item Notice that $\uU$ is only an equivalence class of orbifold atlases. We have not defined a maximal atlas, since the definition of orbifold charts would force the maximal atlas to be a proper class (and not a set). We avoid the set theoretic problems incurred by such a construction. However, by abuse of notation we will sometimes write $(U,G,\pi) \in \uU$ to denote an orbifold chart compatible with the given orbifold structure $\uU$.
\end{compactenum}
\end{rem}

For the rest of this paper we shall always assume that the orbifolds considered are defined as in Definition \ref{defn: haef:ofdII}. As we have already remarked, the definition of orbifolds given in the previous section is equivalent to our working definition of an orbifold. In particular the changes of orbifold charts restrict locally to open embeddings in the sense of Proposition \ref{prop: ch:prop}. On occasion it will turn out to be advantageous to work with embeddings of orbifold charts, as Proposition \ref{prop: ch:prop} is then available. 
\subsection{The topology of the base space of an orbifold}
\setcounter{subsubsection}{0}
In this section, we compile several facts about orbifolds which are well known in the literature (cf. \cite{alr2007,follie2003,bb2008,chen2001}).We give proofs for the reader's convenience.

\begin{lem}\label{lem: ofdch:bt} \no{lem: ofdch:bt}
  For any orbifold $(Q,\uU)$, the family of open subsets $\setm{\tilde{V} \coloneq \pi (V)}{(V,G,\pi) \in \uU}$ is a base for the topology on $Q$. 
\end{lem}

\begin{proof}
 Let $p \in Q$ and $U \subseteq Q$ an open neighborhood of $p$. Choose an orbifold chart $(V,G,\pi) \in \uU$ such that $p \in \tilde{V} = \pi (V)$. The map $\pi$ is given by the composition of the quotient map onto the orbit space with a homeomorphism onto an open set. Hence Lemma \ref{lem: orbitmap} shows that $\pi$ is continuous and open. The set $\pi^{-1} (U)$ is an open subset of $V$ containing some element $\hat{p} \in \pi^{-1} (p)$. By Lemma \ref{lem: st:nbhd} we can choose a $G_{\hat{p}}$-invariant open set $S$ such that $\hat{p} \in S \subseteq \pi^{-1} (U)$ and $(S,G_{\hat{p}}, \pi|_{S})$ is an orbifold chart. By construction, $p \in \pi (S) \subseteq U$, proving the lemma.
\end{proof}

To analyse the structure of the base space we need a well known fact from topology:
 
 \begin{prop} \label{prop: para:coco} \no{prop: para:coco}
 If $X$ is a Hausdorff space that is locally compact and paracompact, then each component of $X$ is $\sigma$-compact. If, in addition, $X$ is locally metrizable, then $X$ is metrizable and every component has a countable basis of the topology.
\end{prop}

\begin{proof}
 By  \cite[XI. Theorem 7.3]{dugun1966} each component is $\sigma$-compact. The space $X$ is paracompact, locally metrizable and Hausdorff, hence we may choose a locally finite closed cover consisting of metrizable subspaces. Then $X$ is metrizable by \cite[Theorem 4.4.19]{Engelking1989}). Each connected component $C$ is Lindelöf by \cite[XI. Theorem 7.2]{dugun1966}. We deduce from \cite[Corollary 4.1.16]{Engelking1989} that $C$ is second countable.
\end{proof}

\begin{prop}\label{prop: ofd:prop} \no{prop: ofd:prop}
 If $(Q,\uU)$ is an orbifold, then the topological space $Q$ has the following properties: 
\begin{compactenum}
                          \item $Q$ is a locally compact Hausdorff space.
			  \item $Q$ is connected if and only if $Q$ is path connected.
			  \item $Q$ is metrizable.
			  \item Every connected component $C$ of $Q$ is open, $\sigma$-compact and second countable.
\end{compactenum}
We remark that $Q$ is not necessarily second countable. 
\end{prop}

\begin{proof} \begin{compactenum}
               \item The space $Q$ is Hausdorff by definition of an orbifold. Clearly being a locally compact space is a local condition, i.e. may be checked within $\pi (U)$, where $(U, G, \pi) \in \uU$ is an arbitrary orbifold chart. Lemma \ref{lem: orbitmap} shows that $\pi (U)$ is a locally compact Hausdorff space, since every finite dimensional Hausdorff manifold $U$ is such a space.
	       \item The quotient map onto the orbit space is continuous and open (Lemma \ref{lem: orbitmap}) and manifolds are locally path-connected. Thus $Q$ is locally path connected, whence the assertion follows from general topology \cite[V. Theorem 5.5]{dugun1966}.
	       \item For every chart $(U,G,\pi) \in \uU$ the group $G \subseteq \Diff (U)$ is finite. The manifold $U$ is locally metrizable (since every chart is a homeomorphism) and a paracompact locally compact Hausdorff space. By Proposition \ref{prop: para:coco}, $U$ is metrizable. 
                The quotient map onto an orbit space is a closed-and-open map by Lemma \ref{lem: orbitmap}. Since metrizability is an invariant of closed-and-open maps by \cite[Theorem 4.2.13]{Engelking1989}, the space $Q$ is locally metrizable. Summing up, $Q$ is a locally metrizable, locally compact and paracompact Hausdorff space. Again by Proposition \ref{prop: para:coco} the metrizability of $Q$ follows.
		\item The space $Q$ is locally path-connected, which implies the openness of $C$ by \cite[V. 5.4]{dugun1966}. We already know that $Q$ is a Hausdorff space which is paracompact and locally compact. Every component of $Q$ is then $\sigma$-compact and second countable by Proposition \ref{prop: para:coco}.
              \end{compactenum}
 To prove the last remark, consider the following counterexample: Let $(Q,\uU)$ be an arbitrary orbifold modeled on a topological space $Q \neq \emptyset$ and $I$ be a set with cardinality at least $\aleph_1$. Construct the orbifold $(Q_I,\uU_I)$ by defining the topological space $Q_{I} \coloneq \coprod_{i \in I} Q$ as the disjoint union of copies of $Q$ and the orbifold charts on every copy of $Q$ as copies of charts in $\uU$. Then $(Q_I,\uU_I)$ is \emph{not second countable}, even if $Q$ is. 
\end{proof}

\subsection{Local groups and the singular locus}
Let $(Q,\uU)$ be an orbifold of dimension $n$, $(U,G,\pi) \in \uU$ an orbifold chart of $Q$ and $x \in U$. Let $z \coloneq \pi (x)$. We deduce from \cite[Lemma 2.10]{follie2003} that the differential at $x$ induces a faithful representation $G_x \rightarrow T_xU, g \mapsto T_xg$ and hence a faithful representation of $G_x$ in $\Gl (n,\RR)$ (cf.\ also Lemma \ref{lem: lin:ch}). The corresponding finite subgroup of $\Gl (n, \RR)$ is unique up to conjugation in $\Gl (n, \RR)$ (induced by the change of chart maps). This conjugacy class will be called $TG_x$. Since $G_{gx} = g G_xg^{-1}$ for any $g \in G$, we have $TG_x = TG_{gx}$. Let $\lambda \colon (V,H,\psi) \rightarrow (U,G,\pi)$ be an embedding of orbifold charts and $y \in V$ with $\lambda (y) =x$ and $\lambda \circ h = \overline{\lambda} (h) \circ \lambda$ for $h \in H$, entailing that $\overline{\lambda} (H_y) = G_x$ by Proposition \ref{prop: ch:prop} and  
\begin{displaymath}
  TG_x = T_y\lambda TH_y (T_y\lambda)^{-1}.
\end{displaymath}
 Thus the conjugacy class of $TG_x$ depends only on the point $z$ and not on the choice of the orbifold chart $(U,G,\pi)$ on $Q$ or on $x$. Hence the following definition is justified.

\begin{defn}[local group]\label{defn: locgp} \no{defn: locgp}
 Let $(Q, \uU)$ be an orbifold. For every $z\in Q$, by the above there is a group $\Gamma_z(Q) \subseteq \Gl (n,\RR)$ which is unique up to conjugation in $\Gl (n,\RR)$. We call $\Gamma_z (Q)$ the \ind{}{local group} of $z$. In the literature $\Gamma_z(Q)$\glsadd{local group} is also called the \ind{}{isotropy group} of $z$. We avoid this and reserve \tl isotropy group\tr for the subgroup of a group acting on a manifold, which fixes a given point.
\end{defn}

The singularities, i.e. points with non-trivial local group, generate a structure which distinguishes a non-trivial orbifold from a manifold. We claimed that orbifolds are manifolds with ``mild singularities''. To emphasize this point we shall investigate the singular locus (i.e. the set of all singularities). As a consequence of Newman's Theorem \ref{thm: newman}, the singular locus is a nowhere dense closed subset of the base space of an orbifold. In other words, the topological base space of an orbifold contains an open and dense manifold. A proof for this result is given in the rest of this section:    

\begin{defn}[Singular locus]\label{defn: sg:loc} \no{defn: sg:loc}
 Let $(Q,\uU)$ be an orbifold. The \ind{}{singular locus} of $Q$ is the subset \glsadd{singular stratum}
	\begin{displaymath} 
	 \Sigma_Q \coloneq \setm{z \in Q}{\Gamma_z (Q) \neq \set{1}}.
	\end{displaymath}
 In a chart $(U,G,\pi)$, one has $\Sigma_Q \cap \pi (U) = \pi (\Sigma_G)$, where $\Sigma_G$ is the set of points in $U$ with non trivial isotropy subgroup with respect to the action of $G$. An element $x \in Q$ is called a \ind{point!singular}{singular point} if $x \in \Sigma_Q$ and $x$ is called \ind{point!non-singular}{non-singular} if $x \not \in \Sigma_Q$. \\
 Since there are different orbifold structures on the same topological space, occasionally we have to indicate which one is meant. In these cases we shall write $\Gamma_z (Q,\uU)$ resp. $\Sigma_{(Q,\uU)}$, to avoid confusion. 
\end{defn}


\begin{prop}[Newman, Thurston]\label{prop: sl:ce} \no{prop: sl:ce}
 The singular locus $\Sigma_Q$ of an orbifold $(Q,\uU)$ is a closed set with empty interior.
\end{prop}

\begin{proof}
 Let $(U, G, \pi)$ be any chart at some point $p \in Q$. By definition $\Sigma_Q \cap \pi (U)$ is the image of $\Sigma_G$. As $G \subseteq \Diff (U)$ is finite, we deduce from Newman's Theorem \ref{thm: newman} that the set $\nN_{U}$ of non-singular points in $U$ is open and dense. Lemma \ref{lem: orbitmap} shows that the quotient map $\pi$ onto the orbit space is open, whence 
	\begin{displaymath}
	    \Sigma_Q = Q \setminus \bigcup_{(U, G , \pi) \in \uU} \pi (\nN_{U})
	\end{displaymath} 
 is a closed set. Since $\nN_U$ is dense in $U$, $\pi (\nN_U)$ is dense in $\pi(U)$. Then $(Q\setminus \Sigma_Q) \cap \pi (U)$ is dense in $\pi (U)$ and since the open sets $\pi (U)$ cover $Q$ (for some atlas), $Q\setminus \Sigma_Q$ is dense in $Q$. In particular $(\Sigma_Q)^\circ = \emptyset$ holds.
\end{proof}
\subsection{Orbifold atlases with special properties}

In this section, we construct special orbifold atlases. These atlases are needed later on, to construct charts for the diffeomorphism group of an orbifold.

\begin{defn}
 Let $(Q,\uU)$ be an orbifold and $\vV$ a representative of $\uU$. We say that another representative $\wW$ of $\uU$ \ind{orbifold atlas!refinement}{refines $\vV$} (or is a \emph{refinement of the atlas $\vV$}) if for every chart $(W, G, \psi) \in \wW$, there is a chart $(V,H,\pi) \in \vV$ and an open embedding of orbifold charts $\lambda_{W,V} \colon (W,G, \psi) \rightarrow (V,H,\pi)$. Given another representative $\vV'$ of $\uU$, we say that $\wW$ is a \emph{common refinement of $\vV$ and $\vV'$}, if $\wW$ refines $\vV$ and $\wW$ refines $\vV'$. 
\end{defn}


\begin{lem}\label{lem: ofa:cr} \no{lem: ofa:cr}
 For an orbifold $(Q,\uU)$ and two arbitrary representatives $\vV, \vV'$ of $\uU$, there exists a common refinement $\wW$ of $\vV$ and $\vV'$.
\end{lem}

\begin{proof}
 Since the union $\wW \coloneq \vV \cup \vV'$ is an orbifold atlas for $(Q,\uU)$, i.e. all charts are pairwise compatible, we may choose for each $x \in Q$ an orbifold chart whose image contains $x$ which and embedds into a chart in $\vV$ and a chart in $\vV'$ (cf.\ Definition \ref{defn: moerd:ofc}). The collection of all charts chosen this way is an atlas, which is a common refinement of $\vV$ and $\vV'$.
\end{proof}

\begin{lem}\label{lem: refn:bt} \no{lem: refn:bt}
  Let $(Q, \uU )$ be an orbifold. For any representative $\vV$ of $\uU$, consider the classes of orbifold charts \glsadd{embtype1}\glsadd{embtype2}
	\begin{align*}
	 \uU \Subset \vV 	&\coloneq \setm{(U,H,\phi) \in \uU}{\exists \lambda_{U,V} \colon (U,H,\phi) \rightarrow (V,G,\psi) \text{ embedding, for some } (V,G,\psi) \in \vV  }\\
	 \uU \sqsubset \vV	&\coloneq \setm{(U,H,\phi) \in \uU\Subset\vV}{\overline{\lambda_{U,V} (U)} \subseteq V \text{is compact}}.
	\end{align*}
 Then the sets $\setm{\phi (U)}{(U,H,\phi) \in \uU \Subset \vV}$ and $\setm{\phi (U)}{(U,H,\phi) \in \uU \sqsubset \vV}$ of open sets are bases for the topology on $Q$.\\
 Note that the compactness of $\overline{\lambda_{U,V} (U)}$ in $V$ implies that $\overline{\phi(U)} \subseteq \psi (V)$.
\end{lem}

\begin{proof}
 Consider an arbitrary open set $\Omega \subseteq Q$ and some point $x \in \Omega$. The set $\vV$ is an atlas and thus, there is some chart $(V,G,\psi) \in \vV$ with $x \in \im \psi$, say $x = \psi (y)$. Because $V$ is locally compact, $y$ has a compact neighborhood $K$ in $V$, contained in the open set $\psi^{-1} (\Omega)$. By Lemma \ref{lem: lin:ch}, $K$ contains an $G$-stable open neighborhood $W$ of $y$ in $V$. Then $(W,G_W,\psi|_W) \in U \sqsubset \vV$ (because $\lambda_{W,V}$ can be chosen as the inclusion map and $\im \psi|_W \subseteq \Omega$.
\end{proof}

\begin{defn}
 Let $(Q, \uU)$ be an orbifold. An orbifold atlas $\vV \coloneq \setm{(V_i,G_i,\pi_i )}{ i \in I}$ of $(Q,\uU)$ is called \ind{orbifold atlas!locally finite}{locally finite orbifold atlas} if the family $(\pi_i (V_i))_{i \in I}$ is a locally finite family of open sets.\footnote{We assume here that the atlas is ``indexed by $I$'' in the sense that the map $I \rightarrow \vV, i \mapsto (V_i,G_i,\psi_i)$ is injective.} 
\end{defn}

\begin{lem}\label{lem: locfin} \no{lem: locfin}
 Let $(Q, \uU)$ be an orbifold. Then the following holds:
\begin{compactenum}
\item There is a locally finite representative $\vV$ of $\uU$.
\item For each representative $\wW$ of $\uU$, there is a locally finite representative $\wW'$ which refines $\wW$.
\item The refinement $\wW'$ in (b) may be chosen with the following property: For each $(U,G,\psi) \in \wW'$, there are $(V,H, \varphi)\in \wW$ and an open embedding $\lambda_{U,V}$ of orbifold charts such that $\overline{\lambda_{U,V} (U)} \subseteq V$ is a compact set, whence $\overline{\tilde{U}} \subseteq \tilde{V}$. (using notation as in Lemma \ref{lem: ofdch:bt})
\end{compactenum}
Taking identifications, without loss of generality $\lambda_{U,V}$ is just the canonical inclusion (of sets) and $G$ is a subgroup of $H$. 
\end{lem}

\begin{proof}
 \begin{compactenum}
	\item The topological space $Q$ is a locally compact Hausdorff space. 
	For each $q \in Q$ pick a compact neighborhood $U_q$ of $q$. Then $(U_q^{\circ})_{q\in Q}$ is an open cover of $Q$. By paracompactness of $Q$, there is a locally finite open refinement $(W_j)_{j\in J}$ of $(U_q^{\circ})_{q\in Q}$. Note that every $\overline{W_j}$ is compact. By \cite[Lemma 5.1.6]{Engelking1989}, there exists a shrinking $(O_j)_{j \in J}$ of $(W_j)_{j \in J}$ that is an open cover of $Q$ such that $\overline{O_j} \subseteq W_j$ for each $j \in J$. The uniformized subsets of $Q$ form a basis of the topology by Lemma \ref{lem: ofdch:bt}. Thus for each $j \in J$, the compact set $\overline{O_j}$ is covered by finitely many uniformized sets which are contained in $W_j$, say $\overline{O_j} \subseteq \bigcup_{k=1}^{n_j} B_{j,k}$. Since the family $(W_j)_{j\in J}$ is locally finite,  
				\begin{displaymath}
                                          \setm{B_{j,k}}{j \in J, \ k =1,\ldots , n_j}
                            	\end{displaymath}
	is a locally finite open covering of $Q$ by uniformized subsets. The corresponding atlas $\vV$ is thus locally finite.  
   	\item[(b) and (c)] We may argue as in (a), but replace the set of all uniformized subsets of $Q$ by the set of all uniformized subsets, which are images of $\uU \Subset \wW$ (resp. images of $\uU \sqsubset \wW$ for (c)). Since Lemma \ref{lem: refn:bt} assures that these sets of images are bases of the topology, no further changes in the proof are needed. For the last statement identify $U$ and $\lambda_{U,V} (U)$ resp. $G$ with $\overline{\lambda} (G)$.\vspace{-2.25em}
 \end{compactenum}
 \end{proof}

\begin{lem}\label{lem: fin:al} \no{lem: fin:al}
 Let $(Q,\uU)$ be an orbifold and $\wW$ a locally finite orbifold atlas such that for each $(V,H,\varphi) \in \wW$ the uniformized subset $\varphi (V)$ is relatively compact. Consider a refinement $\wW'$ as in Lemma \ref{lem: locfin} (c) indexed by a set $I$. There exists a map $\alpha \colon I \rightarrow \wW$, which associates to each $i$ a chart $(V_{\alpha (i)}, H_{\alpha (i)}, \varphi_{\alpha (i)})$ into which $(U_i,G_i,\psi_i)$ embeds (as an orbifold chart) via an inclusion of sets $U_i \rightarrow V_{\alpha (i)}$. Furthermore, $I_V \coloneq \alpha^{-1} (V,H,\varphi) \subseteq I$ is finite for each $(V,H,\varphi) \in \wW$. 
\end{lem}

\begin{proof}
 Lemma \ref{lem: locfin} (c) ensures that for each $i \in I$, there is at least one chart in $\wW$ such that $(U_i,G_i,\psi_i)$ embeds into this chart via the inclusion of sets. Choose a chart $(V_{\alpha (i)}, H_{\alpha (i)}, \varphi_{\alpha (i)})$ such that $\overline{U_i} \subseteq V_{\alpha (i)}$ is compact, $G_i \subseteq H_{\alpha (i)}$ and $\psi_i = \varphi_{\alpha (i)}|_{U_i}$ holds. We obtain a map $\alpha \colon I \rightarrow \wW$ with the desired properties. For each $(V,H,\varphi) \in \wW$, the uniformized subset $\varphi (V)$ is relatively compact. Since $\wW'$ is locally finite, there is only a finite subset of $I$ such that $\psi_i (U_i) \cap \overline{\varphi (V)} \neq \emptyset$. Therefore $I_V \coloneq \alpha^{-1} (V,H,\varphi)$ is finite for each $(V,H,\varphi) \in \wW$.
\end{proof}

Later on an orbifold atlas will be needed which is adapted to a certain closed and discrete set. To construct such an atlas we need to deal with some technical difficulties in the following Lemma: 

\begin{lem}\label{lem: hilf:para}
 Let $X$ be a paracompact topological space, $D\subseteq X$ be a closed discrete subset (i.e. $X$ induces the discrete topology on $D$) Then there exist disjoint open neighborhoods $\Omega_x \subseteq X$ for $x \in D$ such that $(\Omega_x)_{x\in D}$ is locally finite.
\end{lem}
  
\begin{proof}
 For $x \in D$ let $V_x$ be an open neighborhood of $x$ such that $V_x \cap D = \set{x}$. Then $\vV \coloneq \setm{V_x}{x\in D} \cup \set{X \setminus D}$ is an open cover of $X$ and there is a locally finite open cover $(W_j)_{j \in J}$ subordinate to $\vV$. Let $J' \coloneq \setm{j \in J}{D\cap W_j \neq \emptyset}$. Then $(W_j)_{j \in J'}$ is an open cover of $D$ and for each $j \in J'$, there is $x_j \in D$ with $W_j \subseteq V_{x_j}$. Since $V_{x_j} \cap D = \set{x_j}$, $x_j$ is uniquely determined. Since $D \subseteq \bigcup_{j \in J'} W_j$, the map $J' \rightarrow D, j \mapsto x_j$ is surjective. For $x \in D$ choose $j(x) \in J'$ with $x_{j(x)} = x$. Then $(W_{j(x)})_{x\in D}$ is a locally finite open cover of $D$. Since every paracompact space is normal by \cite[Theorem 5.1.5.]{Engelking1989}, the space $X$ is a regular topological space. Hence there is a neighborhood $C_x \subseteq W_{j(x)}$ which is closed in $X$. The locally finite union $A_x \coloneq \bigcup_{y \in D\setminus \set{x}} C_y$ of closed sets is closed and $x \not \in A_x$ since $x \not \in V_y \supseteq C_y$. Define $\Omega_x \coloneq C_x^\circ \setminus A_x$. Then $(\Omega_x)_{x\in D}$ has the desired properties. 
\end{proof}

\begin{prop}\label{prop: atl:lfpts} \no{prop: atl:lfpts}
 Let $(Q,\uU)$ be an orbifold, $\vV \in \uU$ an orbifold atlas and $D$ a closed discrete subset of $Q$. There exist locally finite atlases $\aA = \setm{(U_i,G_i,\psi_i)}{i \in I}$ and $\bB = \setm{(W_j,H_j, \varphi_j)}{j \in J} \in \uU$ such that all of the following conditions are satisfied:
	\begin{compactenum}
	 \item the charts in $\aA, \bB$ are relatively compact, i.e.\ if $(U,G,\psi)$ is such a chart, then the set $\overline{\psi (U)}$ is a compact subset of $Q$, 
	 \item the atlas $\aA$ refines $\bB$ and $\bB$ refines $\vV$ as in Lemma \ref{lem: locfin} (c),
	 \item For $z \in D$, there are unique $i_z \in I$ and $j_z \in J$ with $z \in \psi_i (V_i)$ and $z \in \varphi_j (U_j)$, respectively, 
	 \item If $Q$ is $\sigma$-compact, then the sets $I$ and $J$ are countable.
	 \end{compactenum}
\end{prop}

\begin{proof}
 It suffices to construct $\bB$ with the asserted properties (to get $\aA$, we apply the same construction with $\bB$ instead of $\vV$). The space $Q$ is a metrizable locally compact space by Proposition \ref{prop: ofd:prop}. Using Lemma \ref{lem: hilf:para}, we may choose disjoint open neighborhoods $\Omega_z \subseteq Q$ for $z \in D$ such that $(\Omega_z)_{z\in D}$ is locally finite.
 As $Q$ is locally compact, we may choose for each $z \in D$ a compact neighborhood $L_{1,z} \subseteq \Omega_z$. By Lemma \ref{lem: refn:bt}, for each $z$ there is a relatively compact orbifold chart $(U_{z}, G_{z}, \varphi_{z}) \in \uU \sqsubset \vV$ such that $z \in \varphi_z (U_z) \subseteq \overline{\varphi_{z} (U_{z})} \subseteq L_{1,z}^\circ$. Furthermore, the inclusion of sets induces an embedding of orbifold charts. Again by local compactness, we may choose for each $z$ a compact neighborhoods $z \in L_{2,z} \subseteq \varphi_{z} (U_{z})$.\\
 The set $L_{2,z}$ is contained in $L_{1,z}$. Since each $L_{1,z}$ is contained in $\Omega_z$ and these sets form a locally finite family, the family $(L_{2,z})_{z \in D}$ is locally finite. The set $L \coloneq \bigcup_{z \in D} L_{2,z}$ is thus closed by \cite[Corollary 1.1.12]{Engelking1989} and we may consider the open subset $Q' \coloneq Q \setminus L$. Now $Q'$ is locally compact and as $Q$ is metrizable by Proposition \ref{prop: ofd:prop}, the subspace $Q'$ is paracompact. The images of the class $\rR \coloneq \setm{(V,H,\pi) \in \uU \sqsubset \vV}{\pi (V) \subseteq Q'}$ form a basis for the topology on $Q'$. Using an argument analogous to Lemma \ref{lem: locfin} (c), there is a locally finite orbifold atlas $\bB' = \setm{(W_j,H_j,\varphi_j)}{j\in J'} \subseteq \rR$ for $Q'$ such that each chart $(W,H,\varphi) \in \bB'$ is relatively compact and embeds into some member of $\vV$ as in Lemma \ref{lem: locfin} (c). 
Notice that by construction none of the charts in $\bB'$ contain elements of $D$.\\
 For each $z \in D$, the set $L_{z} \coloneq L_{1,z} \cap Q \setminus \varphi_{z} (U_{z}) \subseteq Q'$ is compact. The atlas $\bB'$ is locally finite and thus there are finite subsets $J'_z \subseteq J'$ such that $\varphi_j (W_j) \cap L_{z} \neq \emptyset$ iff $j \in J'_z$. Assume that $P$ is the image of an orbifold chart in $\bB'$ which is contained in 
   \begin{displaymath}
    O \coloneq Q \setminus \bigcup_{z \in D} L_{z} = \left(\bigcup_{z \in D} \varphi_z (U_z)\right) \cup \left( Q \setminus \bigcup_{z \in D} L_{1,z}\right).
   \end{displaymath}
 As each $L_{1,z}$ is a closed set and the family $(L_{1,z})_{z \in D}$ is locally finite, the union of the sets $L_{1,z}$ is closed by \cite[Corollary 1.1.12]{Engelking1989}. Therefore $O$ is an open set and by construction  
   \begin{displaymath}
     P = \left(\bigcup_{z \in D} \varphi_z (U_z) \cap P \right) \cup \left(P \cap \left(Q \setminus \bigcup_{z\in D} L_{1,z}\right)\right) 
   \end{displaymath}
 is a disjoint union of two open sets. As orbifold charts are connected, we deduce that their images are located as follows:\\ 
 Either the image is contained in $Q \setminus \bigcup_{z\in D} L_{1,z}$, or it intersects at least one of the $L_{z}, z\in D$, or it is contained in $\bigcup_{z \in D} \psi_{z} (U_z)$. Discarding the charts whose image is contained in $\bigcup_{z \in D} \varphi_{z} (U_z)$, we obtain the subset 
	\begin{align*}
	 J'' \coloneq \bigcup_{z \in D} J'_z \cup \setm{j\in J'}{\varphi_j (W_j) \cap \bigcup_{z\in D} L_{1,z} = \emptyset}
	\end{align*}
 of $J'$ such that the family $\bB'' \coloneq \setm{(W_j,H_j,\varphi_j)}{j \in J''}$ 
 covers $Q \setminus \bigcup_{z \in D} \psi_z (U_z)$.
 \\
 Set $J \coloneq J'' \coprod D$. The set indexes the atlas $\bB \coloneq \bB'' \cup \setm{(W_z,H_z,\varphi_z)}{z\in D}$. By construction, $\bB$ is a refinement of $\vV$ with the properties of Lemma \ref{lem: locfin} (c).\\
 It remains to prove that $\bB$ is locally finite: As $\bB'''$ is a locally finite atlas, it suffices to check the following condition: For each $z \in D$, only finitely many charts in $\bB''$ intersect the image of $(U_z, G_z,\varphi_z)$.
 For each $z \in D$, the charts indexed by $z$ are contained in $L_{1,z}$ and by construction only a finite number of charts in $\bB''$ intersect $L_{1,z}$. Thus at most finitely many images of charts in $\bB$ intersect a given $L_{1,z}^\circ, z\in D$, whence $\aA$ and $\bB$ are locally finite.\\[1em] 
 If $Q$ is $\sigma$-compact, then $Q$ is a countable union of compact sets, each of which meets $\im \varphi_j$ for only finitely many $j \in J$ (as $\bB$ is locally finite). Hence $J$ is countable. Likewise the index set $I$ of $\aA$ is countable. 
\end{proof}

The following lemma will allow us to control the local behavior of sections in the tangent orbifold.

\begin{lem}\label{lem: oemb:cov} \no{lem: oemb:cov}
 Let $(Q,\uU)$ be an orbifold and $\wW = \setm{(V_i,H_i,\varphi_i)}{i \in I}$ be a locally finite orbifold atlas. For each $i \in I$, let $K_i \subseteq V_i$ be a compact subset. Then, for each $i \in I$, there is an open cover $\set{Z_i^k}_{1\leq k \leq n_i}$ of $K_i \subseteq V_i$ such that 
	\begin{compactenum}
	 \item the sets $Z_i^k$ are $H_i$-stable for $i \in I$, $1\leq k \leq n_i$,
	 \item for each $j\in I$ with $Z_i^k \cap K_i \cap \varphi_i^{-1}\varphi_j (K_j) \neq \emptyset$ there is an embedding of orbifold charts $\lambda^k_{ij} \colon Z_i^k \rightarrow V_j$ 
	 \item The cover $\set{Z_i^k}_{1\leq k \leq n_i}$ may be chosen such that for each $i \in I, 1 \leq k \leq N_i$ there is a $H_i$-stable set open set $\hat{Z}_i^k$ such that $\overline{Z_i^k}$ is a compact set, contained in $\hat{Z}_i^k$ and each embedding $\lambda^k_{ij}$ is the restriction of an embedding on $\hat{Z}_i^k$.
	\end{compactenum}
\end{lem}

\begin{proof}
 The set $\tilde{K}_i \coloneq \varphi_i (K_i)$ is compact and since $\wW$ is locally finite, there is a finite subset $\fF_i$ of $\wW$ such that $\tilde{K}_i \cap \varphi (V) \neq \emptyset$ if and only if $(V,H,\varphi) \in \fF_i$. In particular, there is a finite set $J_i$ such that $\tilde{K}_{ij} \coloneq \tilde{K}_i \cap \varphi_j (K_j) \neq \emptyset$ if and only if $j \in J_i$. The compact sets $\tilde{K}_{ij}$ are contained in $\tilde{V}_i$. The set 
  \begin{align}
   K_{ij} \coloneq& K_i \cap \varphi_{i}^{-1} (\tilde{K}_{ij}) = K_i \cap \varphi_i^{-1} \varphi_j (K_j) = (\varphi|_{K_i})^{-1} (\varphi_j (K_j)). \label{ident KIJ}
  \end{align}
 is closed in $K_i$ and hence compact. For each $j \in J_i$, the set $K_{ij}$ is contained in $\varphi_{i}^{-1} \varphi_j (V_j)$. Thus each $K_{ij}$ may be covered with open $H_{i}$-stable subsets $\Lambda^r_{ij}$ of $V_{i}$ such that there is an open embedding of orbifold charts $\lambda^r_{ij} \colon \Lambda_r \rightarrow V_{j}$. Since $K_{ij}$ is compact, for each $j$ there is a finite family $(\Lambda_{ij}^r)_{1\leq r \leq m_j}$ which covers $K_{ij}$. As $J_i$ is finite, we obtain for each $x \in K_i$ an open neighborhood  
	\begin{displaymath}
	 N_x \coloneq \bigcap_{x \in\Lambda^r_{ij}} \Lambda_{ij}^r \cap \left( V_{i} \setminus \bigcup_{j \in J,x \not \in K_{ij}} K_{ij}\right).
	\end{displaymath}
 Choose an $H_{i}$-stable connected open neighborhood $x \in Z^x \subseteq N_x$. Each $y \in Z^x$ is contained in $K_{ij}$ only if $x$ is contained in $K_{ij}$ as well. For each $j\in J_i$ such that $x \in K_{ij}$, the open embeddings defined on $\Lambda_{ij}^r$ restrict to an open embedding of orbifold charts on $Z^x$. Since $K_i$ is compact we may select a finite open cover $\setm{Z^{x_k}}{x_k \in K_i , \ 1 \leq k \leq n}$ of $K_i$. Observe that $Z^{x_k} \cap K_i \cap \varphi_i^{-1} \varphi_j (K_j) = Z^{x_k} \cap K_{ij}$ holds by \eqref{ident KIJ}. If this intersection is non-empty, we derive $x_k \in K_{ij}$. By construction, there is an embedding of orbifold charts on $Z^{x_k}$ which satisfies (b). Hence the family $(Z^{x_k})_{1 \leq k \leq n}$ satisfies all properties of assertion (b).\\
 (c) follows directly from (b) and local compactness of each $V_i$: Before selecting a finite cover by some of the $Z^x$, we set $\hat{Z}^x \coloneq Z^x$ and choose for each $x$ a compact neighborhood $x \in C_x \subseteq \hat{Z}^x$. The $H_i$-stable sets form a base of the topology and we may select a new $H_i$-stable subset $x \in Z_x \subseteq C_x^\circ \subseteq \hat{Z}^x$. By compactness of $K_i$, we may select a finite covering from the family $(Z_x)_{x \in K_i}$ which satisfies (c).
\end{proof}
\newpage
\subsection{Examples of orbifolds}\label{ch: examples}
\setcounter{subsubsection}{0}
This section collects well known simple examples from the literature to illustrate the definition of an orbifold. We also fix some terminology  for later use.

\begin{ex}
 Every paracompact smooth finite dimensional manifold $M$ without boundary is an orbifold. An orbifold atlas for $M$ is given by the following set of charts: 
	\begin{displaymath}
	 \setm{(C, \set{\id_C} , \id_C)}{C \subseteq M \text{ a connected component}}
	\end{displaymath}
 where by abuse of notation $\id_C \colon C \rightarrow M$ is the inclusion map.
 We call this orbifold structure induced on the manifold $M$ the \index{trivial orbifold structure}\emph{trivial orbifold structure}.
\end{ex}

\begin{ex}[A mirror in $\RR^2$ {\cite[13.1.1]{thur2002}}]\label{ex: Z2hs} \no{ex: Z2hs}
 Consider $\RR^2$ together with the action of the linear diffeomorphism $\gamma \colon \RR^2 \rightarrow \RR^2, (x,y) \mapsto (-x, y)$. The map $\gamma$ fixes the points $(0,y), y \in \RR$. \\
 An orbifold structure is induced on the quotient $\RR^2/ \langle \gamma \rangle \sim H \coloneq \setm{ (x,y) \in \RR^2}{ x \geq 0 }$:
 
  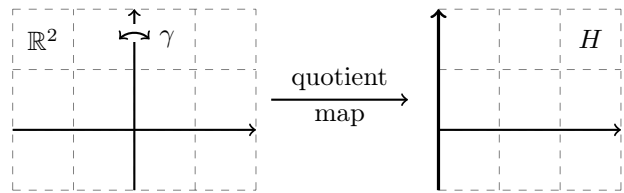
\begin{figure}[htb]\centering
   \begin{tikzpicture}[scale=0.8]
     \draw[style=help lines,dashed,gray] (-2,-1) grid (2,2);
    \coordinate (y) at (0,2);
    \coordinate (x) at (2,0);
    \draw[thick,->] (0,-1) -- (y); 
    \draw[thick,->] (-2,0) -- (x); 
    \draw(0,1.6) node[shape=rectangle,fill=white]{};
    \draw(-1.5,1.5) node[fill=white]{$\RR^2$};
    \draw[thick,<->] (-0.25,1.5) to [bend left =45](0.25,1.5) node[right,fill=white]{$\gamma$};
    \draw[thick,->] (2.25,0.5) -- node[midway,above]{quotient} node[midway,below]{\text{map}}(4.5,0.5);
    \draw[style=help lines,gray,dashed] (5,-1) grid (8,2);
    \draw(7.5,1.5) node[fill=white]{$H$};
    \draw[very thick,->] (5,-1) -- (5,2); 
    \draw[thick, ->] (5,0) -- (8,0);
  \end{tikzpicture}
 \caption[A mirror in $\RR^2$]{A mirror in $\RR^2$. The boundary of the half plane contains the singular points, while points outside the boundary are non-singular points.}
  \end{figure}
This example can be generalized in the following way to manifold with boundary (cf. \cite[13.2.2]{thur2002}):\\[1em]
Let $M$ be a (smooth) manifold  with boundary $\partial M$. Glue together two copies of $M$ along $\partial M$ to obtain the double $dM$ of $M$. Recall that by using a collar around the boundary (cf. \cite[Chapter 4, 6.]{hirsch1976}) the double may be endowed with the unique structure of a smooth manifold without boundary (see \cite[Definition 5.10 and Theorem 6.3]{munkres1966} for a full account on the construction). Again the diffeomorphism which interchanges both halves of the double generates a finite group $\Gamma$. By construction the orbifold $dM/\Gamma$ is isomorphic to $M$. Hence every manifold with boundary is in a natural way an orbifold, whose singular locus is the boundary of the manifold. 
\end{ex}
\begin{ex}[Good orbifolds]
 Let $M$ be a smooth finite dimensional paracompact manifold and $\Gamma \subseteq \Diff (M)$ be a subgroup. Assume that the canonical action of $\Gamma$ on $M$ is proper, i.e. there exists a metric $d$ on $M$ such that $\Gamma$ acts by diffeomorphisms and for each $x \in M$ there exists $r >0$ such that 
  \begin{displaymath}
   \setm{\gamma \in \Gamma}{\gamma . B_r^d (x) \cap B_r^d(x) \neq \emptyset}
  \end{displaymath}
 is finite. Then the orbit space $M/\Gamma$ may be endowed with an orbifold structure induced by the group action of $\Gamma$ on $M$ (cf. \cite[III.\gG 1.3]{msnpc1999} for details). An orbifold which arises in this way is called \ind{orbifold!developable}{developable} or \ind{orbifold!good}{good} orbifold.\\
 A particularly attractive situation arises if $M$ is a connected, paracompact manifold and $\Gamma$ is finite. Then the good orbifold obtained from these data possesses an atlas with one chart, i.e. $(M,\Gamma, \pi)$, where $\pi \colon M \rightarrow M/\Gamma$ is the canonical quotient map. In Example \ref{ex: sphere} we compute orbifold structures for $M = \SSS^2$. Several of these structures will be good orbifolds.  
\end{ex}

\begin{ex}[Symmetric products {\cite[Example 1.13]{alr2007}}]
 Suppose that $M$ is a smooth finite dimensional, paracompact manifold. Consider the \emph{symmetric product} $X_n \coloneq M^n /S_n$, where $M^n$ is the $n-$fold Cartesian product of $M$ and $S_n$ the symmetric group of $n$ letters which acts on $M^n$ by permutation of coordinates. Tuples of points have non trivial isotropy groups if they contain a number of repetitions in their coordinates. The diagonal of $M^n$ is fixed by each element of the finite group $S_n$.
\end{ex}

In the next example we consider two orbifold charts on the same topological space which induce non-diffeomorphic orbifolds.

\begin{ex}[\hspace{-0.5pt}{\cite[Example 2.2]{pohl2010}}]
 Let $Q \coloneq [0,1[$ be the topological space with the induced topology of $\RR$. The map $f \colon Q \rightarrow Q, x \mapsto x^2$ is a homeomorphism. Let $\rho \colon \RR \rightarrow \RR$ be the reflection in $0$. Consider the map $p \colon ]-1,1[ \rightarrow Q, x \mapsto \abs{x}$. Then $p$ induces a homeomorphism $]-1,1[ / \langle \rho\rangle$ and we derive orbifold charts $V_1 \coloneq (]-1,1[ , \langle \rho\rangle , p)$ and $V_2 \coloneq (]-1,1[,\langle \rho \rangle , f \circ p)$.\\
 However, these orbifold charts are not compatible. Assume to the contrary that they are compatible. Since $f\circ p (0) = 0 = p(0)$ there exist open connected neighborhoods $U_1$, $U_2$ of $0$ in $]-1,1[$ and a diffeomorphism $h \colon U_1 \rightarrow U_2$ such that $f \circ p = p \circ h$. This equation leads to $h(x) \in \set{\pm \sqrt{\abs{x}}}$. By continuity we have the following choices for $h$: 
 \begin{align*}
  h_1 (x) &\coloneq \sqrt{\abs{x}}				& h_2 (x) &\coloneq -\sqrt{\abs{x}} \\
  h_3 (x) &\coloneq \left\{\begin{array}{cl}
				-\sqrt{\abs{x}}, & x\leq 0 \\
				 \sqrt{\abs{x}}, & x\geq 0
	   		\end{array}\right. 	&h_4 (x) &\coloneq \left\{\begin{array}{cl}
										 \sqrt{\abs{x}}, & x\leq 0 \\
										-\sqrt{\abs{x}}, & x\geq 0
	   									\end{array}\right. \\
 \end{align*}
 Since none of the above is a differentiable, the two charts are not compatible. 
\end{ex}

\begin{ex}[Orbifold structures on the 2-sphere] \label{ex: sphere}
The following examples are all modeled on the $2$-sphere $\SSS^2$, i.e. the topological space of each of the orbifolds is the $2$-sphere with the topology turning it into a smooth manifold. Examples of this type first appeared in \cite{thur2002}. We give a detailed construction based on the exposition in \cite{cgo2006}:\\
Let $N$ be the north pole and $S$ the south pole of $\SSS^2$. Endow the sphere with the usual topology turning $\SSS^2$ into a smooth manifold. Define charts around $N$ and $S$, respectively, as follows:\\ 
Let $X_i \coloneq B^{\RR^2}_{\frac{3}{4}\pi} (0)$, $i=1,2$ be the open disc of radius $\frac{3}{4}\pi$ centered at $0$ in $\RR^2$. We describe points in polar coordinates $(r, \theta), 0 \leq r < \tfrac{3}{4}\pi,\ 0 \leq \theta < 2\pi$. Recall that the geodesics connecting $N$ and $S$ on $\SSS^2$ are the great circles connecting $N$ and $S$.
To construct the charts pick a great circle $C$ connecting $N$ and $S$. Every great circle connecting $N$ and $S$ can uniquely be identified by an angle of rotation $0 \leq \theta < 2\pi$. Furthermore, each $x$ on $\SSS^2 \setminus \set{S}$ is uniquely determined by a set of coordinates $(r, \theta ), \ 0 \leq r < \pi, 0 \leq \theta < 2\pi$. Here $r$ is the length of the geodesic segment between $x$ and $N$. Analogously we may identify each point $x$ in $\SSS^2\setminus \set{N}$ by a pair $(\pi - r,\theta), 0 \leq r < \pi, 0 \leq \theta < 2\pi$, where $\pi-r$ is the length of the geodesic segment between $x$ and $N$. We obtain (inverses of) charts
	\begin{align*}
        \psi_1 &\colon X_1 \rightarrow \SSS^2 , (r,\theta) \mapsto \begin{pmatrix} \cos \theta  & - \sin \theta & 0 \\ \sin \theta & \cos \theta & 0 \\ 0 & 0 & 1 \end{pmatrix} \begin{pmatrix} \sin r \\ 0 \\ \cos r  \end{pmatrix} ,\\ \psi_2& \colon X_2 \rightarrow \SSS^2 , (r,\theta) \mapsto \begin{pmatrix} \cos \theta  & - \sin \theta & 0 \\ \sin \theta & \cos \theta & 0 \\ 0 & 0 & 1 \end{pmatrix} \begin{pmatrix} \sin r \\ 0 \\ -\cos r  \end{pmatrix}
	\end{align*}
 for the manifold $\SSS^2$. These charts turn $\SSS^2$ into a smooth compact manifold in the usual way.\\
 We construct an orbifold structure on $\SSS^2$: Let $n_i \in \NN$ for $i=1,2$. Consider the subgroup $G_i \subseteq \Diff (X_i)$ which corresponds to a rotation $\sigma_i$ of order $n_i$ on $X_1$ and $X_2$. The map $p_i \colon X_i \rightarrow X_i , (r \cos \theta,r \sin \theta) \mapsto (r\cos (n_i \theta),r\sin (n_i\theta))$ identifies two points if and only if they are in the same $G_i$ orbit.\\
 Consider the quotient map $X_i \rightarrow X_i/G_i$ and canonically identify the orbit space with the \tl cone\tr 
    \begin{displaymath}
     C_i\coloneq \setm{(r,\theta) \in X_i}{ 0 \leq \theta < \frac{2\pi}{n_i}}
    \end{displaymath}
 endowed with the quotient topology with respect to $c_i \colon X_i \mapsto C_i , (r,\theta) \mapsto (r, \theta \mod \frac{2\pi}{n_i})$. A computation shows that $c_i \colon C_i \rightarrow X_i , (r,\theta) \mapsto (r,n_i\theta)$ is a homeomorphism of the topological spaces $C_i$ and $X_i$. Moreover, $p_i$ factors through the quotient $X_i/G_i \sim C_i$. We obtain orbifold charts $(X_i , G_i,q_i), \ i\in \set{1,2}$ with $q_i \coloneq \psi_i \circ p_i$. A computation shows that $A_{ij} \coloneq q_i^{-1} (\im q_j) = \setm{(r,\theta) \in X_i}{\tfrac{\pi}{4} \leq r \leq \frac{3\pi}{4}}$ is an open annulus. Furthermore, we obtain for each $(r, \theta) \in A_{ij}$ a neighborhood $\Omega_{r,\theta}$ such that the mapping 
 \begin{displaymath}
	\tau_{ij}|_{\Omega_{r,\theta}} \coloneq (q_i|^{q_j (\Omega_{r,\theta})})^{-1} \circ q_j\colon \Omega_{r,\theta} \rightarrow X_i , \ (r, \theta ) \mapsto \left(\pi -r , \frac{n_j}{n_i}\cdot \theta \right), \ i \neq j \in \set{1,2}
 \end{displaymath}
 makes sense. The maps $\tau_{ij}|_{\Omega_{r,\theta}}$ are local diffeomorphisms, which commute with the orbifold charts, i.e.\ $q_i \tau_{ij} = q_j|_{\dom \tau_{ij}},\ i \neq j \in \set{1,2}$ holds. Locally the restrictions of the maps $\tau_{ij}$ thus yield change of chart morphisms. Since we obtain change of charts for each $x \in A_{ij}$, the orbifold charts are compatible and induce the structure 
	\begin{displaymath}
	 \SSS^2_{(n_1, n_2)} \coloneq (\SSS^2 , \setm{(X_i,G_i,q_i)}{i=1,2})
	\end{displaymath}
 of a compact orbifold on $\SSS^2$. As a topological space, the base space of $\SSS^2_{(n_1, n_2)}$ coincides with $\SSS^2$ with the usual topology. We distinguish the following cases: 
\paragraph{$n_1, n_2 = 1$} In this case we have $q_i = \psi_i, \ i=1,2$ and thus $\SSS^2_{(1, 1)}$ is just the $C^\infty$-manifold $\SSS^2$. As a connected trivial orbifold, $\SSS^2$ is a good orbifold.
\paragraph{$n_1 > 1 , n_2 =1$:} We obtain a cone-shaped singularity of order $n_1$ in $N$, while $S$ is a regular point. The orbifold $\SSS^2_{(n_1, 1)}$ is called \emph{$\ZZ_{n_1}$-teardrop}. It is an example of a non developable orbifold. Indeed the orbifold $\SSS^2_{(n_1, n_2)}$ is developable (good) if and only if $n_1 = n_2$ holds (see \cite[Ch. III.$\gG$, Example 1.4 (1)]{msnpc1999}).
\paragraph{$n_1 \neq n_2,\ n_1,n_2 >1$:} We obtain an orbifold with two cone-shaped singularities of order $n_1$ respectively $n_2$. An orbifold of this kind is called $\ZZ_{n_1}$\hspace{-0.3em}-$\ZZ_{n_2}$-\emph{football}. As already mentioned, this orbifold is non-developable.
\paragraph{$n_1, n_2 = n > 1$:} Consider an action of a finite group of diffeomorphisms $\Gamma \subset \Diff (\SSS^2)$ generated by a rotation of order $n$ on $\SSS^2$ which fixes north and south pole. The group $\Gamma$ acts smoothly, effectively and almost free on $\SSS^2$. Hence the orbit space $\SSS^2/\Gamma$ is an orbifold using the global orbifold chart $\pi \colon \SSS^2 \rightarrow \SSS^2 /\Gamma$. By construction the orbifold structure of this orbifold agrees with $\SSS^2_{(n, n)}$. It is an example of a good orbifold.\\[1em]
On the level of topological spaces all of these orbifolds coincide. However the additional structure of cone-shaped singularities on the space is illustrated in the following picture: \begin{figure}[h]
  \centering
  \definecolor{cffffff}{RGB}{255,255,255}
\begin{tikzpicture}[y=0.80pt, x=0.8pt,yscale=-1, inner sep=0pt, outer sep=0pt]
\begin{scope}[shift={(0,-768.89763)}]
  \path[draw=black,line join=miter,line cap=butt,miter limit=4.00,line
    width=1.600pt] (290.0000,812.3622) .. controls (290.0000,812.3622) and
    (246.8972,841.6833) .. (242.6639,896.4495) .. controls (238.4305,951.2156) and
    (250.1003,990.7242) .. (290.7221,991.0337) .. controls (331.3439,991.3432) and
    (341.9012,952.6604) .. (337.2294,898.0354) .. controls (332.5577,843.4104) and
    (290.0000,812.3622) .. (290.0000,812.3622) -- cycle;
  \begin{scope}[shift={(39.75703,-25.59656)}]
    \path[draw=black,line join=miter,line cap=butt,miter limit=4.00,line
      width=0.800pt] (296.7338,932.5000) .. controls (296.7338,935.0163) and
      (275.8104,937.0561) .. (250.0000,937.0561) .. controls (224.1897,937.0561) and
      (203.2662,935.0163) .. (203.2662,932.5000);
    \path[draw=black,dash pattern=on 0.80pt off 0.80pt,line join=miter,line
      cap=butt,miter limit=4.00,line width=0.800pt] (203.2662,932.5000) .. controls
      (203.2662,929.9838) and (224.1897,927.9440) .. (250.0000,927.9440) .. controls
      (275.8104,927.9440) and (296.7338,929.9838) .. (296.7338,932.5000);
  \end{scope}
  \path[shift={(-70.0,758.89767)},draw=black,fill=cffffff,line join=miter,line
    cap=butt,miter limit=4.00,line width=1.600pt] (260.0000,143.4646) .. controls
    (260.0000,193.1702) and (219.7056,233.4646) .. (170.0000,233.4646) .. controls
    (120.2944,233.4646) and (80.0000,193.1702) .. (80.0000,143.4646) .. controls
    (80.0000,93.7589) and (120.2944,53.4646) .. (170.0000,53.4646) .. controls
    (219.7056,53.4646) and (260.0000,93.7589) .. (260.0000,143.4646) -- cycle;
  \begin{scope}[shift={(-10.0,-9.99997)}]
    \path[draw=black,fill=cffffff,line join=miter,line cap=butt,miter
      limit=4.00,line width=0.800pt] (199.2388,916.4786) .. controls
      (199.2388,922.8282) and (159.3522,927.9757) .. (110.1496,927.9757) .. controls
      (60.9470,927.9757) and (21.0604,922.8282) .. (21.0604,916.4786);
    \path[draw=black,fill=cffffff,dash pattern=on 0.80pt off 0.80pt,line
      join=miter,line cap=butt,miter limit=4.00,line width=0.800pt]
      (21.0604,916.4786) .. controls (21.0604,910.1289) and (60.9470,904.9815) ..
      (110.1496,904.9815) .. controls (159.3522,904.9815) and (199.2387,910.1289) ..
      (199.2387,916.4786);
  \end{scope}
  \path[fill=black] (0,850) node[above right] (text4311-8) {(a)};
  \path[fill=black] (220,850) node[above right] (text4311-8) {(b)};
  \path[fill=black] (380,850) node[above right] (text4311-8) {(c)};
  \path[draw=black,line join=miter,line cap=butt,miter limit=4.00,line
    width=1.600pt] (440.0000,812.3622) .. controls (440.0000,812.3622) and
    (396.8972,841.6833) .. (392.6639,896.4495) .. controls (388.4305,951.2156) and
    (409.8480,971.0338) .. (440.7221,991.0337) .. controls (470.0000,972.3622) and
    (491.9012,952.6604) .. (487.2294,898.0354) .. controls (482.5577,843.4104) and
    (440.0000,812.3622) .. (440.0000,812.3622) -- cycle;
    \path[draw=black,fill=cffffff,line join=miter,line cap=butt,miter
      limit=4.00,line width=0.800pt] (486.7719,907.9465) .. controls
      (486.7719,910.7352) and (465.8827,912.9960) .. (440.1145,912.9960) .. controls
      (414.3463,912.9960) and (393.4571,910.7352) .. (393.4571,907.9465);
    \path[draw=black,fill=cffffff,dash pattern=on 0.80pt off 0.80pt,line
      join=miter,line cap=butt,miter limit=4.00,line width=0.800pt]
      (393.4571,907.9465) .. controls (393.4571,905.1577) and (414.3463,902.8970) ..
      (440.1145,902.8970) .. controls (465.8827,902.8970) and (486.7719,905.1577) ..
      (486.7719,907.9465);
\end{scope}
\end{tikzpicture}
  \caption[Orbifold structures on the sphere]{Orbifold structures on $\SSS^2$: (a) the trivial orbifold $\SSS^2_{(1, 1)}$, i.e. the manifold $\SSS^2$,\\ (b) the teardrop $\SSS^2_{(n_1, 1)}$ and (c) the football $\SSS^2_{(n_1, n_2)}$.}
\end{figure}
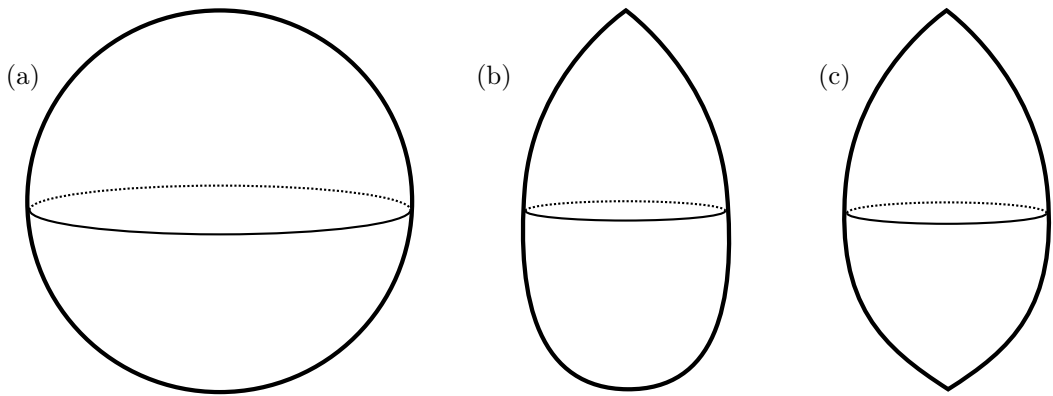
\end{ex}

\section{Maps of Orbifolds}\label{sect: diffeo}
\setcounter{subsubsection}{0}
In this thesis, we use maps of orbifolds as defined in \cite{pohl2010}. For the reader's convenience, we repeat the definitions and constructions of \cite{pohl2010} in Appendix \ref{sect: maps}. In this section, we obtain a characterization of orbifold diffeomorphisms. Then several tools and constructions for later chapters (such as open suborbifolds and orbifold partitions of unity) are provided.

\subsection{Orbifold diffeomorphisms}

Throughout this section, let $(Q_i,\uU_i), \ i\in \set{1,2}$ be arbitrary orbifolds. By definition, diffeomorphisms of orbifolds are the isomorphisms in the category of reduced orbifolds: 

\begin{defn}
 A morphism of orbifolds $[\hat{f}] \in \ORBM$ is called an \ind{orbifold map!diffeomorphism}{orbifold diffeomorphism} if there is $[\hat{g}] \in \ORB ((Q_2, \uU_2) , (Q_1, \uU_1))$ such that 
    \begin{displaymath}
     \ido{1} = [\hat{g}] \circ [\hat{f}] \text{ and } \ido{2} = [\hat{f}] \circ [\hat{g}].
    \end{displaymath}
 In this case, we also write $[\hat{f}]^{-1} \coloneq [\hat{g}]$. Let $\Difforb{(Q_1,\uU_1),(Q_2,\uU_2)}$ be the set of orbifold diffeomorphisms contained in $\ORB ((Q_1,\uU_1), (Q_2,\uU_2))$.\\ To shorten our notation, the \ind{}{orbifold diffeomorphism group} $\Difforb{(Q,\uU),(Q,\uU)}$ will be denoted by \gls{Difforb}.
\end{defn}

We will now characterize the lifts of orbifold diffeomorphisms. It will turn out that an orbifold diffeomorphism is completely determined by properties of its lifts. 

\begin{prop}\label{prop: ofd:diff} \no{prop: ofd:diff}
 Let $[\hat{f}] \in \ORB ((Q_1, \uU_1) , (Q_2, \uU_2))$ be a diffeomorphism of orbifolds. Each representative $\hat{f} = (f, \set{f_i}_{i\in I} , [P_f , \nu_f])$ satisfies the following properties: 
	\begin{compactenum}
	 \item the map $f$ is a homeomorphism and
	 \item every local lift $f_i$ of $\hat{f}$ is a local diffeomorphism. 
	\end{compactenum}
\end{prop}

\begin{proof}
 We first notice that since $[\hat{f}] \circ [\hat{f}]^{-1}$ and $[\hat{f}]^{-1} \circ [\hat{f}]$ are the respective identity morphisms, the maps $f \colon Q_1 \rightarrow Q_2$ and $f^{-1} \colon Q_2 \rightarrow Q_1$ (where $f^{-1}$ is the underlying continuous map of $[\hat{f}]^{-1}$) are homeomorphisms since composition yields the identity on $Q_2$ and $Q_1$, respectively. Hence (a) is true.\\ 
 Two representatives of the class $[\hat{f}]$ are related via lifts of the identity. Lifts of such mappings are local diffeomorphisms, whence locally lifts of different representatives of $[\hat{f}]$ are related via diffeomorphisms to each other. Thus the definition of $[\hat{f}]$ shows that it suffices to prove assertion (b) for any representative $\hat{f}$ of $[\hat{f}]$.\\ 
 Choose and fix representatives $\vV = \setm{(V_i,G_i,\pi_i)}{i \in I}$ of $\uU_1$, $\uU = \setm{(U_j, H_j , \psi_j)}{j \in J}$ of $\uU_2$ and $\wW = \setm{(W_k, L_k, \vp_k)}{k \in K} $ of $\uU_1$ such that the maps $[\hat{f}]$ and $[\hat{f}]^{-1}$  possess representatives $\hat{f} \in \Orb{\vV , \uU}$ and $\hat{g} \in \Orb{\uU , \wW}$, respectively. Let $\alpha \colon I \rightarrow J$ and $\beta \colon J \rightarrow K$ be the maps such that the mappings $f_i \colon V_i \rightarrow U_{\alpha (i)}$ and $g_j \colon U_j \rightarrow W_{\beta (j)}$ are local lifts of $\hat{f}$ and $\hat{g}$, respectively, with respect to orbifold charts $(V_i , G_i, \pi_i)$ and $(U_{\alpha (i)}, G_{\alpha (i)}, \psi_{\alpha (i)})$, $(U_j, G_j , \psi_j)$ and $(W_{\beta (j)} , G_{\beta (j)} , \vp_{\beta_{j}})$. To shorten the notation, set $\tilde{V}_i \coloneq \pi_i (V)$ and derive for every $i \in I$ a commutative diagram: 
	\begin{displaymath}
	 \begin{xy}
  \xymatrix{					
	V_i \ar[rr]^{f_i} \ar[d]^{\pi_i}		&	& U_{\alpha (i)} \ar[rr]^{g_{\alpha (i)}} \ar[d]^{\psi_{\alpha (i)}} 	&& W_{\beta (\alpha (i))} \ar[d]^{\vp_{\beta \alpha (i)}}\\
       \tilde{V}_i   \ar[rr]^{f|_{\tilde{V}_i}} 	&& \tilde{U}_{\alpha (i)} \ar[rr]^{f^{-1}|_{\tilde{U}_{\alpha (i)}}} && \tilde{W}_{\beta \alpha (i)}
  }
\end{xy}
	\end{displaymath}
Composition in the lower row induces the identity $\id_{Q_1}|_{\tilde{V}_i}$. We conclude that for each $i \in I$, the map $g_{\alpha (i)} \circ f_i$ is a local lift of the identity and thus a local diffeomorphism by Proposition \ref{prop: ll:id}. In particular, each $f_i$ is an immersion and hence $\dim Q_1 \leq \dim Q_2$. An analogous argument shows $\dim Q_2 \leq \dim Q_1$, whence $\dim Q_1 = \dim Q_2$. Since the orbifold dimensions coincide, $\dim V_i = \dim U_{\alpha (i)}$ holds. The inverse mapping theorem (see \cite[I, \S 5 Theorem 5.2]{langdgeo2001}) now implies that the immersion $f_i$ is a local diffeomorphism.
\end{proof}

\begin{cor}\label{cor: ofd;isdim} \no{cor: ofd:isdim}
 Two orbifolds $(Q_i, \uU_i),\ i \in \set{1,2}$ which are isomorphic have the same orbifold dimension.
\end{cor}

\begin{defn}
 Consider an orbifold map $[\hat{f}] \in \ORBM$ together with a corresponding representative of orbifold maps $\hat{f} = (f,\set{f_i} , [P_f,\nu_f])$. We say that $[\hat{f}]$ \ind{orbifold map!preserving local groups}{preserves local groups} if $f \colon Q_1 \rightarrow Q_2$ maps every element $p$ of $Q_1$ onto some element $f(p)$ of $Q_2$ such that $\Gamma_p (Q_1) \cong \Gamma_{f(p)} (Q_2)$.\\
 This property may be interpreted as preservation of the local structure of an orbifold. In particular, one would expect that this is a natural property of orbifold diffeomorphisms. Indeed this is true, as the following proposition shows:
\end{defn}

\begin{prop}\label{prop: ld:flgp} \no{prop: ld:flgp} 
 Let $[\hat{f}] \colon (Q_1,\uU_1) \rightarrow (Q_2, \uU_2)$ be a map of orbifolds, with a representative $\hat{f} = (f,\set{f_i}_{i \in I}, (P_f,\nu_f))$ such that $f$ is a homeomorphism and each $f_i$ is a local diffeomorphism. Then $[\hat{f}]$ preserves local groups. In particular, every orbifold diffeomorphism preserves local groups.
\end{prop}

\begin{proof}
 Let $p$ be in $Q_1$. There are orbifold charts $(V,G,\pi) \in \uU_1$ and $(U,H, \psi) \in \uU_2$ together with a local lift $f_V \colon V \rightarrow U$ of $\hat{f}$ such that $p \in \tilde{V}$, $q \coloneq f(p) \in \tilde{U}$ and $f_V$ is a local diffeomorphism. Fix some preimage $\hat{p} \in \pi^{-1}(p)$ and denote its image by $\hat{q} \coloneq f_V(\hat{p})$.\\ 
 Since $G_{\hat{p}}$ is finite, there is an open connected neighborhood $\Omega$ of $\hat{p}$ in $V$ such that for every $\gamma \in G_{\hat{p}}$, there is some $\mu_\gamma \in P_f$ with $\gamma|_{\Omega} = \mu_\gamma|_\Omega$. Thus one obtains 
	\begin{equation}\label{eq: comp:diffeo}
	 f_V (\gamma . x) = \nu_f (\mu_\gamma) f_V (x) \ \forall x \in \Omega , \gamma \in \Gamma_{\hat{p}}
	\end{equation}
 Shrinking $\Omega$ if necessary, we may assume that $\Omega$ to be a $G$-stable open connected subset with $G_\Omega = G_{\hat{p}}$ and $f_V|_\Omega$ is a diffeomorphism onto an open subset of $U$. By \eqref{eq: comp:diffeo}, $\psi \circ f_V$ factors over $\Omega/G_{\hat{p}}$ and it is an open map. Hence $(\Omega, G_{\hat{p}}, \psi \circ f_V)$ is an orbifold chart for $Q_2$. By construction $f_V$ is an embedding of orbifold charts from $(\Omega, G_{\hat{p}}, \psi \circ f_V)$ to $(U,H,\psi)$. Hence $(\Omega, G_{\hat{p}}, \psi \circ f_V) \in \uU_2$ and thus $\Gamma_p \cong G_{\hat{p}} \cong \Gamma_q$ (the groups are even conjugate in $\Gl (n,\RR)$). 
\end{proof}

\begin{rem}\label{rem: odi:tech} \no{rem: odi:tech}
 The proof of Proposition \ref{prop: ld:flgp} provides information about an orbifold map: Consider an orbifold map which satisfies the prerequisites of Proposition \ref{prop: ld:flgp}. Let $f_i \colon V_i \rightarrow W_i$ be its local lift with respect to the charts $(V_i,G_i,\pi_i)$ and $(W_i,H_i,\psi_i)$ and $x \in V_i$. Then there is an arbitrarily small open neighborhood $\Omega_x$ of $x$ in $V_i$ with the following properties: 
	\begin{compactenum}
	 \item $f_i|_{\Omega_x}$ is a diffeomorphism onto an open set $\Omega_{f_i(x)} \coloneq f_i(\Omega_x)$,
	 \item the set $\Omega_x$ is $G_i$-stable with $G_{i , \Omega_x} = G_{i, x}$,
	 \item for each $\gamma \in G_{i, x}$, the restriction $\gamma|_{\Omega_x}$ is an element of $P_f$,
	 \item the set $\Omega_{f_i(x)}$ is $H_i$-stable with $H_{i , \Omega_{f_i(x)}} = H_{i, f_i(x)}$. 
	\end{compactenum}
 In particular, $(\Omega_x ,G_{i,x}, \pi_i|_{\Omega_x})$ and $(\Omega_{f_i(x)}, H_{i, f_i(x)}, \psi_i|_{\Omega_{f_i(x)}})$ are orbifold charts contained in $\uU_1$ and in $\uU_2$,respectively. Locally, we may therefore always construct lifts which are diffeomorphisms.
\end{rem}

 It is possible to construct a charted orbifold map from a family of local lifts as in the last remark: 

\begin{prop}\label{prop: ld:qpgp} \no{prop: ld:qpgp}
 Let $(Q_i,\uU_i), i\in \set{1,2}$ be orbifolds, $f \colon Q_1 \rightarrow Q_2$ a homeomorphism and $\set{f_i}_{i\in I}$ be a family of local lifts of $f$ with respect to $\vV \in \uU_1$ and $\wW \in \uU_2$ such that each $f_i$ is a local diffeomorphism. Assume that $\vV$ satisfies (R2) from Definition \ref{defn: rep:ofdm}. Then there exists a pair $(P,\nu)$ such that $(f, \set{f_i}_{i \in I}, [P,\nu]) \in \Orb{\vV , \wW}$ is a representative of an orbifold map in $\ORBM$. The pair $(P,\nu)$ is unique up to equivalence.
\end{prop}

\begin{proof}
 Let $\vV = \setm{(V_i,G_i,\pi_i)}{i \in I}$ be the representative of $\uU_1$ such that every lift $f_i$ is a map $f_i \colon V_i \rightarrow W_i$ for some $(W_i,H_i, \psi_i) \in \uU_2$. As $f$ is a homeomorphism, $\wW \coloneq \setm{(W_i,H_i, \psi_i)}{i \in I}$ is an orbifold atlas. Define $F \coloneq \coprod_{i \in I} f_i $. Consider the set 
	\begin{displaymath}
	 P \coloneq \setm{h \in \PSI (\vV )}{h \text{ is a change of charts and }F|_{\dom h} , F|_{\cod h} \text{ are \'{e}tale embeddings}}.
	\end{displaymath}
 Clearly $P$ is a quasi-pseudogroup which generates $\PSI (\vV)$. Construct a map $\nu \colon P \rightarrow \PSI (\wW)$ as follows: For $\lambda \in P$ there are $i,j \in I$ such that $\dom \lambda \subseteq V_i$ and $\cod \lambda \subseteq V_j$. The map $F|_{\dom \lambda} = f_i|_{\dom \lambda}$ is a diffeomorphism onto an open set $U_\lambda \subseteq W_i$. We may now define  
	\begin{displaymath} 
	 \nu (\lambda) \coloneq f_j \lambda f_i|_{U_\lambda}^{-1} \colon U_\lambda \rightarrow f_j (\cod \lambda)
	\end{displaymath}
 The set $f_j (\cod \lambda )$ is open since $f_j$ is a local diffeomorphism. Following the definition of $P$, $\nu (\lambda)$ is a diffeomorphism. We compute $\psi_j \nu (\lambda ) = \psi_j f_j \lambda (f_i|_{\dom \lambda})^{-1} = f \pi_j \lambda (f_i|_{\dom \lambda})^{-1} = f \pi_i f_i|_{U_\lambda}^{-1} = f f^{-1} \psi_i|_{U_\lambda} = \psi_i|_{U_\lambda}$, which shows $\nu (\lambda ) \in \PSI (\wW)$. In addition, $F \circ \lambda = f_j \circ \lambda = \nu (\lambda ) \circ f_i|_{\dom \lambda} = \nu (\lambda) \circ F|_{\dom \lambda}$. Thus we have constructed a quasi-pseudogroup $P$ and a well-defined map $\nu \colon P \rightarrow \PSI (\wW)$ satisfying property (R4a) of Definition \ref{defn: rep:ofdm}. Reviewing (R4b)-(R4d) of the same definition, clearly these properties are satisfied by $\nu$. In conclusion, $(f, \set{f_i}_{i\in I} , P,\nu)$ is a representative of an orbifold map.\\
 To prove the uniqueness, assume that there is another pair $(P',\nu')$ turning $(f,\set{f_i}_{i \in I}, (P',\nu'))$ into a charted map. Consider $\lambda \in P$ and $\mu \in P'$ with $\germ_x \lambda = \germ_x \mu$ for some $x$ in their domeins. Then the mappings $f_j \circ \lambda = \nu(\lambda) \circ f_i|_{\dom \lambda}$ and $f_j \circ \mu = \nu' (\mu) f_i|_{\dom \mu}$ coincide in some neighborhood of $x$. Since $f_j$ is a local diffeomorphism, the mappings $\nu(\lambda)$ and $\nu' (\mu)$ coincide in some neighborhood of $F(x)$.
\end{proof}

Combining Remark \ref{rem: odi:tech} and Proposition \ref{prop: ld:qpgp}, we obtain the following corollary: 

\begin{cor}\label{cor: ll:diff} \no{cor: ll:diff}
 Let $f \colon Q_1 \rightarrow Q_2$ be a homeomorphism and $\set{g_i}_{i\in I}$ a family of local lifts of $f$ with respect to atlases $\vV'$ and $\wW'$ such that each $g_i$ is a local diffeomorphism. Assume that $\vV'$ satisfies (R2). Then there exist an orbifold atlas $\vV$ which refines $\vV'$ indexed by some $J$ and an orbifold atlas $\wW$ which refines $\wW'$ and a family of lifts $f_j$ with respect to $(V_j, G_j,\psi_j) \in \vV$, $(W_{\beta (j)}, H_{\beta (j)}, \varphi_{\beta (j)}) \in \wW$ such that each $f_j$ is a diffeomorphism. In addition there is a unique equivalence class $[P,\nu]$ with $P = \Ch_{\vV'}$ and $\nu (\lambda ) \coloneq f_k \lambda (f_j|_{\dom \lambda})^{-1}|_{f_{j} (\dom \lambda)}$ for $\lambda \in \CH{V_j}{V_k}$, $(V_r, G_r,\psi_r) \in \vV'$ for $r \in \set{ j,k}$ such that $\hat{f} \coloneq (f, \set{f_j}_{j \in J}, [P,\nu]) \in \Orb{\vV',\wW'}$.
\end{cor}

\begin{lem}\label{lem: add:cor} \no{lem: add:cor}
 Let $\vV = \setm{(V_i,G_i,\psi_i)}{i\in I}$ and $\wW = \setm{(W_j, H_j,\varphi_j)}{j \in J}$ be atlases for orbifolds $(Q_1,\uU_1)$ and $(Q_2,\uU_2)$, respectively. Consider a charted map of orbifolds $\hat{f} = (f, \set{f_i}_{i \in I}, [P,\nu]) \in \Orb{\vV,\wW}$ with the same properties as the map $\hat{f}$ in Corollary \ref{cor: ll:diff}. Then the following holds: 
	\begin{compactenum}
	 \item For each $G_i$-stable subset $\Omega \subseteq V_i$, the set $f_i (\Omega )$ is an $H_{\beta (i)}$-stable subset of $W_{\beta (i)}$ with isotropy group  $H_{\beta (i), f_i(\Omega)} \cong G_{i,\Omega}$.
	 \item After possibly shrinking $\vV$ and $\wW$, we may assume that the map\footnote{which assigns to each index $i$ an index $\beta (i)\in J$ such that $g_i \colon V_i \rightarrow W_{\beta (i)}$ holds.} $\beta \colon I\rightarrow J$ is bijective.
	 \item If $\beta$ is bijective, then $\nu \colon \Ch_\vV \rightarrow \Ch_\wW$ is a bijection. 
	\end{compactenum}
\end{lem}

\begin{proof}
 \begin{compactenum}
  \item Let $\Omega \subseteq V_i$ be a $G_i$-stable subset with isotropy subgroup $G_{i,\Omega}$ and $x \in \Omega$. Because $P = \Ch_\vV$ the proof of Proposition \ref{prop: ld:flgp} applies and we can take $\Omega_x = \Omega$ in Remark \ref{rem: odi:tech}.
  \item If there are $i,j \in I$ with $\beta (i) = \beta (j)$, we obtain a diffeomorphism $f_j^{-1}f_i \colon V_i \rightarrow V_j$. A quick computation shows that $\psi_j f_j^{-1} f_i = f^{-1} \varphi_{\beta (j)} f_i = \psi_i$ and thus $f_j^{-1} f_i$ is an embedding of orbifold charts. Reversing the roles of $i$ and $j$, also $f_i^{-1}f_j$ is an embedding of orbifold charts. Therefore we may omit one index of the pair $i,j$ with $\beta  (i) =\beta (j)$ and the set of orbifold charts indexed by the reduced set will again be an orbifold atlas. The axiom of choice allows us to shrink $\vV$ to obtain an orbifold atlas (which by abuse of notation will also be called $\vV$) such that $\beta$ is injective. Clearly since $f$ is a homeomorphism, the set of charts $\setm{(W_j,H_j,\varphi_j) \in \wW}{j = \beta (i) \text{ for some } i \in I}$ is an orbifold atlas. Thus by replacing $J$ with $\beta (I)$, we may assume that $\beta$ is surjective, hence bijective. 
  \item It is obvious that $\nu$ is injective. Let $\lambda \in \Ch_{W_k,W_l}$ be any change of charts morphism with $(W_r,H_r,\varphi_r) \in \wW, \ r = k,l$. There are unique $i,j \in I$ with $\beta (i ) = k$ and $\beta (j) = l$ and we obtain a diffeomorphism $\mu (\lambda) \coloneq f_j^{-1} \lambda f_i|_{f_i^{-1} (\dom \lambda)} \colon  f_i^{-1} (\dom \lambda) \rightarrow f_j^{-1} (\cod \lambda)$. A quick computation leads to $\psi_j \mu (\lambda) = f^{-1} \varphi_{l} \lambda f_i|_{f_i^{-1} (\dom \lambda)} = f^{-1} f \psi_i|_{\dom \lambda} = \psi_i|_{\dom \lambda}$ which proves that $\mu (\lambda) \in  \CH{V_i}{V_j}$. By construction $\nu (\mu (\lambda)) = \lambda$ holds and thus $\nu$ is a bijection. 
 \end{compactenum}

\end{proof}

The next proposition is the converse of Proposition \ref{prop: ofd:diff}, i.e. we shall prove that the properties of orbifold diffeomorphisms in Proposition \ref{prop: ofd:diff} actually characterizes those, and are equivalent to the categorical definition. The leading idea is to use the local properties of the lifts (i.e. that every lift may locally be inverted) to construct a family of lifts for $f^{-1}$. In general a given lift may not be inverted globally. Nevertheless it is possible to construct smaller charts and induced lifts, which may be inverted globally.

\begin{prop}\label{prop: ofd:iso} \no{prop: ofd:iso}
 Let $(Q_i,\uU_i), i\in\set{1,2}$ be orbifolds and $\vV \in \uU_1, \wW\in \uU_2$. Consider an charted map $\hat{f} \coloneq (f , \set{f_i}_{i\in I} , [P,\nu]) \in \Orb{\vV , \wW}$. If $f$ is a homeomorphism and $f_i \colon V_i \rightarrow W_{\alpha (i)}$ is a local diffeomorphism for each $i \in I$, then the orbifold map $[\hat{f}] \in \ORBM$ is a diffeomorphism of orbifolds.
\end{prop}

\begin{proof}
 Combining Corollary \ref{cor: ll:diff} and Lemma \ref{lem: add:cor}, there are orbifold atlases $\vV'$ indexed by $J$ and $\wW'$ indexed by $K$ together with a representative $\hat{g} \coloneq (f , \set{g_j}_{j\in J} , [P',\nu']) \in \Orb{\vV' , \wW'}$ of $[\hat{f}]$ such that each lift $g_j \colon V_{j} \rightarrow W_{\beta (j)}$ is a diffeomorphism and the map $\beta \colon J \rightarrow K$ is a bijection. We use the computation from the proof of Lemma \ref{lem: add:cor}: The inverse $g_j^{-1} \colon W_{\beta (j)} \rightarrow W_j$ of $g_j$ is a local lift of $f^{-1}$ with respect to $(W_{\beta (j)},H_{\beta (j)},\varphi_{\beta (j)})$ and $(V_j,G_j,\psi_j)$. Since $f$ is a homeomorphism, the family $\wW'$ is an atlas for $Q_2$ indexed by $K$. As each $g_j^{-1}$ is a diffeomorphism, by Proposition \ref{prop: ld:qpgp} there is a pair $Q \subseteq \PSI (\wW')$ and $\mu \colon P \rightarrow \PSI (\vV')$ such that $\hat{h} \coloneq (f^{-1}, \set{g_j^{-1}}_{j \in K}, [Q,\mu]) \in \Orb{\wW , \vV}$.\\
 Consider the compositions $\hat{h} \circ \hat{g}$ and $\hat{g} \circ \hat{h}$: 
The local lift for every $j \in J$ of $\hat{g}$ has been constructed as inverse maps of the local lift of $\hat{g}$ with respect to $(V_j , G_j,\psi_j)$ and $(W_{\beta (j)},H_{\beta (j)},\varphi_{\beta (j)})$. Thus the composition of both representatives gives a lift of the identity and we derive 
 \begin{displaymath}
 [\hat{f}] \circ [\hat{g}] = [\hat{h} \circ \hat{g}] = \ido{2} \quad \text{and} \quad [\hat{g}]\circ [\hat{f}] = [\hat{g} \circ \hat{h}] = \ido{1} .
 \end{displaymath}
\end{proof}

Observe that the proof of the last proposition yields the following fact: Assume that each member of the family of local lifts for an orbifold map is a diffeomorphism. Then this family uniquely determines the orbifold map. In particular, each orbifold diffeomorphism is uniquely determined by its family of local lifts:
 
\begin{cor}\label{cor: ofdiff:un} \no{cor: ofdiff:un}
 An orbifold diffeomorphism $[\hat{f}] \in \ORBM$ is uniquely determined by the family of local lifts $\set{f_i}_{i \in I}$ where $(f , \set{f_i}_{i \in I}, [P,\nu]) \in [\hat{f}]$ is an arbitrary representative.
\end{cor}

\begin{prop}\label{prop: diff:undmap}
 An orbifold diffeomorphism $[\hat{f}] \in \ORBM$ is uniquely determined by its underlying continuous map, i.e.\ for $\hat{f} = (f,\set{f_j}_{j \in J} , [P,\nu])$ the map $f$ uniquely determines $[\hat{f}]$.
\end{prop}

\begin{proof}
 Let $[\hat{g}] \in \ORBM$ be another orbifold diffeomorphism with underlying map $f$. Then the underlying map of $[\hat{g}]^{-1}$ is $f^{-1}$. Hence each representative $\hat{h}$ of $[\hat{g}]^{-1}\circ [\hat{f}]$ is given by $\hat{h} = (\id_{Q}, \set{h_i}_{i \in I} , [P',\nu'])$. Recall from Construction \ref{con: cp:chom} that the lifts $h_i, i \in I$ arise as composition of suitable lifts of representatives of $[\hat{f}]$ and $[\hat{g}]^{-1}$. Since all lifts of orbifold diffeomorphism are local diffeomorphisms by Proposition \ref{prop: ofd:diff}, we deduce that each $h_i$ is a local diffeomorphism. Now Proposition \ref{prop: ident:eq} implies $\ido = [\hat{h}] = [\hat{g}]^{-1} \circ [\hat{f}]$. Thus $[\hat{g}] = [\hat{f}]$ follows and proves the assertion. 
\end{proof}

Summarizing the preceding results, one obtains: 

\begin{cor}\label{cor: diff:char} \no{cor: diff:char}
 For an orbifold map $[\hat{f}] \in \ORBM$ the following are equivalent: 
	\begin{compactenum}
	 \item $[\hat{f}]$ is an orbifold diffeomorphism, 
	 \item each representative $(f , \set{f_i}_{i \in I}, [P,\nu]) \in [\hat{f}]$ satisfies: $f$ is a homeomorphism and each $f_i$ is a local diffeomorphism,
	 \item there is a representative $\hat{f} = (f , \set{f_i}_{i \in I}, [P,\nu])$ of $[\hat{f}]$ such that $f$ is a homeomorphism and each $f_i$ is a local diffeomorphism
	 \item there is a representative $\hat{f} = (f , \set{f_j}_{j \in J}, [P,\nu]) \in \Orb{\vV ,\wW}$ of $[\hat{f}]$ such that $f$ is a homeomorphism and each $f_j$ is a diffeomorphism. Furthermore, the assignment $\alpha \colon \vV \rightarrow \wW$ such that $f_j$ is a local lift with respect to the pair $(V_j,G_j,\varphi_j)$, $(W_{\alpha (j)}, G_{\alpha (j)}, \psi_{\alpha (j)})$ can be chosen bijective. 
	\end{compactenum}
 If $\hat{f}$ is as in (d), then a representative of $[\hat{f}]^{-1}$ is given by $(f^{-1}, \set{f_j^{-1}}, [\nu (P) , \theta]) \in \Orb{\wW , \vV}$. Here $\theta \colon \nu (P) \rightarrow \PSI (\vV)$ assigns to $\lambda \in \nu (P)$ with $\dom \lambda \subseteq W_{\alpha (i)}$ and $\cod \lambda \subseteq W_{\alpha (j)}$ the map $\theta (\lambda) \coloneq f_j^{-1} \lambda f_i|_{f_i^{-1} (\dom \lambda)}$.\\
 In particular, an orbifold diffeomorphism is uniquely determined by its underlying map and we obtain a natural inclusion of the orbifold diffeomorphisms into the set of homeomorphisms:
  \begin{align*}
   \Difforb{(Q_1,\uU_1),(Q_2,\uU_2)} &\rightarrow \text{\rm Homeo} ((Q_1,\uU_1),(Q_2,\uU_2)) \\
	  [(f, \set{f_i}_{i \in I}, [P,\nu])] &\mapsto f
  \end{align*}
 
\end{cor}

We remark that the characterization of orbifold diffeomorphisms via any family of lifts will be crucial for the rest of this work. It enables us avoid the technical details of the definition of orbifold maps. Instead we may think of an orbifold diffeomorphism as a family of compatible smooth lifts. In particular, these results enable an efficient investigation of orbifold diffeomorpismgroup.

\subsection{Open suborbifolds and restrictions of orbifold maps}
We define the notion of an open suborbifold to introduce the restriction of an orbifold map to an open subset. Any subset of a metrizable space with the induced topology is again a metrizable space. Every metrizable space is paracompact and Hausdorff by \cite[Theorem 5.1.3]{Engelking1989}. Since the base space $Q$ of the orbifold $(Q,\uU)$ is metrizable by Proposition \ref{prop: ofd:prop}, each of the subspaces in the following constructions will be a paracompact Hausdorff space.

\begin{defn}[open suborbifold]\label{defn: sofd:op}\no{defn: sofd:op}
 Let $(Q,\uU)$ be an orbifold. An orbifold $(X , \xX)$ is called an \ind{orbifold!open suborbifold}{open suborbifold} of $(Q,\uU)$ if there is a map $[\hat{\iota}] = [(\iota , \set{\iota_k}_{k\in I} , [P,\nu])] \in \ORBM[(X,\xX),(Q,\uU)]$ such that 
	\begin{compactenum}
	 \item $\iota$ is a topological embedding with open image,
	 \item every $\iota_k$ is a local diffeomorphism.
	\end{compactenum}
 A map $[\hat{\iota}]$ with the properties (a) and (b) is called an \ind{orbifold map!open embedding}{open embedding of orbifolds}.
\end{defn}
Since it will not be needed, we shall not define the general notion of a (possibly non-open) suborbifold. The reader is refered to \cite[Definition 2.3]{alr2007} for further information on this topic.

\begin{defn}[Restriction of an orbifold map to an open subset]\label{defn: osup}
Let $(Q,\uU)$ be an orbifold and $\Omega \subseteq Q$ be an open subset. Choose an atlas $\aA \in \uU$ such that the images of $(V,G,\psi) \in \aA$ which satisfy $\psi (V) \subseteq \Omega$ cover $\Omega$. Then $\aA|_\Omega \coloneq \setm{(V,G,\psi) \in \aA}{\psi (V) \subseteq \Omega}$ is an orbifold atlas for $\Omega$. Notice that the equivalence class $\uU_\Omega$ of $\aA|\Omega$ does not depend on the choice of $\aA$ and defines an unique orbifold structure on $\Omega$. The inclusion $\iota_\Omega \colon \Omega \hookrightarrow Q$ of sets induces an open embedding of orbifolds, which we denote by $[\hat{\iota}_\Omega] \colon (\Omega , \uU_\Omega) \rightarrow (Q,\uU)$. Define the \ind{orbifold map!restriction to open subset}{restriction $[\hat{f}]|_\Omega$ of $[\hat{f}] \in \ORBM[(Q,\uU),(Q_2,\uU_2)]$ to $\Omega$} via  
    \begin{displaymath}
      [\hat{f}]|_{\Omega} \coloneq [\hat{f}] \circ [\hat{\iota}_\Omega].
    \end{displaymath}
\end{defn}

\begin{defn}[Corestriction of an orbifold map]\label{defn: cores} \no{defn: cores}
 Let $(X,\xX)$ be an open suborbifold of $(Q,\uU)$ together with an open embedding of orbifolds $[\hat{\iota}]$. Consider another orbifold  $(Q',\vV)$ and a map $[\hat{f}] \in \ORBM[(Q',\vV),(Q,\uU)]$ with representative $\hat{f} = (f,\set{f_k}_{k \in I}, [P,\nu]) \in [\hat{f}]$ such that $\im f \subseteq \im \iota$.\\
 For $k \in I$, let the lifts be given as $f_k \colon V_k \rightarrow U_{\alpha (k)}$, where $(U_{\alpha (k)} , G_{\alpha (k)}, \psi_{\alpha (k)})$ is an orbifold chart. Then $\im f_k \subseteq \psi_{\alpha (k)}^{-1} (\im \iota)$ holds. As $\im f_k$ is connected, it is contained in a connected component of the invariant set $\psi_{\alpha (k)}^{-1} (\im \iota)$. The connected components of an invariant set are $G_{\alpha (k)}$-stable subsets of $U_{\alpha (k)}$. Hence these connected components can be made into orbifold charts for the subset $\im \iota$. Using these charts, Lemma \ref{lemdef: ind:ofdm} shows that there is a representative $\hat{g} \in \Orb{\vV',\uU'}$ of $[\hat{f}]$ such that each lift $g_k \colon V'_k \rightarrow U_{k}'$ of $\hat{g}$ satisfies $\varphi (U_k') \subseteq \im \iota$.
 Define the \ind{orbifold map!corestriction}{corestriction} of $[\hat{f}]$: 
	\begin{displaymath}
	 [\hat{f}]|^{\im \iota} \coloneq \left[(f|^{\im \iota} , \set{g_k}_{k} , [P', \nu'])\right] \in \ORBM[(Q',\vV), (\im \iota , \uU_{\im \iota})]
	\end{displaymath}
 Here $(P',\nu')$ is the pair obtained via Lemma \ref{lemdef: ind:ofdm} for $\hat{g}$. In particular, we obtain a unique map $([\hat{\iota}]|^{\im \iota})^{-1} \circ [\hat{f}]|^{\im \iota} \in \ORBM[(Q',\vV), (X , \xX)]$ into the open suborbifold. By definition of the equivalence relation (Definition \ref{defn: sim:ofd}), the class $[\hat{f}]|^{\im \iota}$ does not depend on any choices made in the construction.  
\end{defn}

\begin{rem}\label{rem: subofd} \no{rem: subofd}
\begin{compactenum}
 \item 
 An orbifold $(X,\xX)$ is an open suborbifold of $(Q,\uU)$ if and only if there is an orbifold diffeomorphism from $(X,\xX)$ to an orbifold which arises as the restriction of $\uU$ to an open subset.
 \item Consider an open subset $\Omega \subseteq Q$ and the representative $\hat{f}= (f, \set{f_k}_{k \in I}, [P,\nu])$ of $[\hat{f}] \in \ORBM[(Q,\uU), (Q',\wW')]$ such that there is $J \subseteq I$ with the following properties:\\ $\vV_\Omega \coloneq \set{(V_j , G_j, \pi_j)}_{j \in J} \subseteq \uU_\Omega$ and $\Omega = \bigcup_{j \in J} \pi_j (V_j)$ hold.\\
 Define $P_J \coloneq P \cap \Ch_{\vV_\Omega}$ and set $\nu_J \coloneq \nu|_{P_J}$. The composition in $\ORB$ is induced by composition of suitable representatives. A computation with the representative above yields $[\hat{f}]|_\Omega = [\hat{h}]$, where $\hat{h}\coloneq (f|_{\Omega}, \set{f_j}_{j \in J}, [P_J,\nu_J])$. 
 \item Let $(X,\xX)$ be an open suborbifold with open embedding of orbifolds $[\hat{\iota}]$. By construction $[\hat{f}]|_{\im \hat{\iota}} = [\hat{f}] \circ [\hat{\iota}] \circ [\hat{\iota}]|^{\im \iota})^{-1}$ holds.
 \item In Section \ref{sect: tofd} tangent spaces of orbifolds and the tangent orbifold are defined. As these objects are defined via an arbitrary orbifold chart, analogous to the manifold case, for each open suborbifold $(X,\xX)$ of $(Q,\uU)$ the tangent spaces $\tT_p^\xX X$ and $\tT_{\iota (p)}^\uU Q$ are canonically isomorphic.\footnote{Here the symbol $\tT_p^\uU Q$ denotes the tangent space of the orbifold $(Q,\uU)$. The notation was chosen to emphasize the dependence on the orbifold structures $\xX$ and $\uU$.} If the open suborbifold is an open subset, we shall identify the tangent spaces later on.
\end{compactenum}
\end{rem}
\newpage 

\subsection{Partitions of unity for orbifolds}

\begin{defn}
 Let $(Q,\uU)$ be an orbifold, $\vV = \setm{(V_i, G_i ,\pi_i)}{i \in I}$ a representative of $\uU$ and endow $\RR$ with the trivial orbifold structure (i.e. the one induced by its manifold structure). \\ A family $\set{(\chi_i , \set{\chi_{i,j}}_{j \in J} , [P_i , \nu_i])}_{i \in I}$ in $\Orb{\vV ,\set{\id_\RR}}$ is called a \ind{orbifold map!partition of unity}{smooth orbifold partition of unity} subordinate to $\vV$ if the family of continuous maps $\set{\chi_i}_{i\in I}$ is a partition of unity subordinate to the open covering $\set{\pi_i (V_i)}_{i\in I}$, i.e.
	\begin{compactenum}
	 \item $\supp \chi_i \subseteq \pi_i (V_i)$ for all $i \in I$,
	 \item the family $(\supp \chi_i)_{i \in I}$ is locally finite, 
	 \item $\chi_i \geq 0$, for all $i \in I$ and $\sum_{i \in I} \chi_i (x) = 1$ for each $x \in Q$.	
	\end{compactenum}
\end{defn}

\begin{prop}[Partition of Unity]\label{prop: pu:hofd} \no{prop: pu:hofd}
 Let $(Q, \uU)$ be an orbifold. For each representative $\vV$ of $\uU$ there exists a smooth orbifold partition of unity subordinate to $\vV$.
\end{prop}

\begin{proof}
 Each representative of $\uU$ allows a locally finite refinement by Lemma \ref{lem: locfin} (b), thus the assertion will be true if the existence of a smooth orbifold partition of unity for an arbitrary locally finite representative of $\uU$ can be verified.\\ Let $\vV \coloneq \setm{(U_\alpha, G_\alpha , \pi_\alpha )}{\alpha \in I}$ be a locally finite representative and  $\tilde{\vV} \coloneq \set{\pi_\alpha (U_\alpha)}_{\alpha \in I}$ be the family of open images of the charts in $\vV$. 
 Since $Q$ is a paracompact Hausdorff space, applying \cite[Lemma 5.1.6]{Engelking1989} twice, there are locally finite families of open sets $\tilde{W}_\alpha^1 \subseteq \overline{\tilde{W}_\alpha^1} \subseteq \tilde{W}_\alpha^2 \subseteq \overline{\tilde{W}_\alpha^2} \subseteq \pi_\alpha (U_\alpha)$ such that $\setm{\tilde{W}_\alpha^1}{\alpha \in I}$ covers $Q$ (here the closure means closure in $Q$). 
Let $W_\alpha^i \coloneq \pi_\alpha^{-1} (\tilde{W}_\alpha^i),\ i\in \set{1,2}$. Observe that since $\overline{\tilde{W}}_\alpha^i \subseteq \im \pi_\alpha$, it is closed in the subspace topology. On $\im \pi_\alpha^i$, we identify $\pi_\alpha$ with the quotient map onto the orbit space of the $G_\alpha$-action on $U_\alpha$. This map is surjective continuous, open and closed by Lemma \ref{lem: orbitmap}. Hence for $i=1,2$ \cite[III. Theorem 8.3 (5) and Theorem 11.4]{dugun1966} imply $\pi_\alpha (\overline{W_\alpha^i}) = \overline{\tilde{W}_{\alpha}^i}$ and $\overline{W_{\alpha}^i} \subseteq \pi_\alpha^{-1} (\overline{\tilde{W}_\alpha^i})$. Vice versa \cite[III. Theorem 11.2 (2)]{dugun1966} yields $\overline{W_\alpha^i}= \pi^{-1}_\alpha \left(\overline{\tilde{W}_\alpha^i}\right)$. By construction, every set $W^i_\alpha$ is $G_\alpha$-invariant.\\
 The manifold $U_\alpha$ is a smooth connected paracompact (hence second countable by Proposition \ref{prop: para:coco}) and finite dimensional manifold. By the smooth Urysohn Lemma (cf.\ \cite[Corollary 3.5.5]{cmfd01}) for manifolds, there is a smooth map $f^\alpha \colon U_\alpha \rightarrow [0,1]$ such that $f^\alpha|_{\overline{W_\alpha^1}} \equiv 1$ and $\supp f^\alpha \subseteq W_\alpha^2$. Define an equivariant smooth map $\theta_\alpha \colon U_\alpha \rightarrow \RR$ with values in $[0,1]$ by averaging over $G_\alpha$: 
 \begin{displaymath}
	 \theta_{\alpha} (y) \coloneq \frac{1}{\lvert G_\alpha \rvert} \sum_{\gamma \in G_\alpha} f^\alpha (\gamma . y).
	\end{displaymath}
 Notice that $W_\alpha^1 \subseteq \supp \theta_\alpha \subseteq W_\alpha^2$ still holds by $G_\alpha$-invariance of these sets. In particular, the map vanishes outside of $\overline{W^2_\alpha}$. For every $\beta \in I$, define a map as follows:

	\begin{displaymath}
	 \theta_{\alpha ,\beta} \colon U_\beta \rightarrow [0,1] , \ x\mapsto \begin{cases}
	                                                                       \theta_\alpha (y) 	& \pi_\beta (x) = \pi_\alpha (y) \text{ for some } y \in U_\alpha \\
										0						& \pi^{-1}_\alpha \pi_\beta (x) = \emptyset 
	                                                                      \end{cases}
	\end{displaymath}
 The $G_\alpha$-equivariance of $\theta_{\alpha}$ implies that $\theta_{\alpha , \beta}$ is well-defined, and it is $G_\beta$-equivariant. We claim that $\theta_{\alpha , \beta}$ is smooth: To see this, note that for each $x \in \pi_\beta^{-1} (\im \pi_\alpha)$, there is an open neighborhood $V_x\subseteq U_\beta$ of $x$ and a smooth change of charts $\lambda \colon V_x \rightarrow U_\alpha$. On the open set $V_x$, the map $\theta_{\alpha , \beta}$ is a composition of smooth maps: $\theta_{\alpha , \beta}|_{V_x} = \theta_\alpha \circ \lambda$. Hence on $\pi_\beta^{-1} (\im \pi_\alpha)$ the map $\theta_{\alpha , \beta}$ is smooth.\\
 By construction $\supp \theta_{\alpha}  \subseteq \overline{W^2_\alpha} \subseteq U_\alpha$ holds, i.e.\ we obtain 
 $\pi_\beta (\supp \theta_{\alpha , \beta}) \subseteq \overline{\tilde{W}^2_\alpha} \subseteq \im \pi_\alpha$. The above shows that $\theta_{\alpha , \beta}$ is a smooth map on the open neighborhood $\pi_\beta^{-1} (\im \pi_\alpha )$ of its support. On the open set $U_\beta \setminus \supp \theta_{\alpha,\beta}$ the map vanishes and in conclusion $\theta_{\alpha , \beta}$ is smooth.\\
 Notice that $\theta_{\alpha , \alpha} = \theta_\alpha$ holds by construction. Since the family $\tilde{\vV}$ is locally finite, for $x\in Q$ there are only finitely many $\alpha \in I$ such that $\pi_\alpha^{-1} (x) \neq \emptyset$. Define another $G_\beta$-equivariant smooth map on $U_\beta$: 
	\begin{displaymath}
	 \chi_{\alpha , \beta} \colon U_\beta \rightarrow [0,1] , \chi_{\alpha , \beta} \coloneq \frac{\theta_{\alpha , \beta}}{\sum_{\delta \in I} \theta_{\delta , \beta}}.
	\end{displaymath}
 The map $\chi_{\alpha , \alpha}$ satisfies $\chi_{\alpha , \alpha}|_{U_\alpha \setminus \overline{W^2_\alpha}} \equiv 0$. Since $\pi_\alpha$ is an open map and $\pi_\alpha (\overline{W_\alpha^2})$ is closed, the map $\chi_{\alpha,\alpha}$ descends to a continuous map on $Q$: 
 	\begin{displaymath}
 	 \chi_{\alpha} \colon Q \rightarrow [0,1] , x \mapsto \begin{cases}
 	                                                                        \chi_{\alpha , \alpha} (x) 	& x = \pi_\alpha (y) \text{ with } y\in U_\alpha \\
										0						& x \in Q \setminus U_\alpha			
 	                                                                           \end{cases}.
 	\end{displaymath}
 By construction $\supp \chi_\alpha \subseteq \pi_\alpha (U_\alpha )$. For every $\sigma \in I$, the smooth map $\chi_{\alpha , \sigma}$ is a lift of $\chi_{\alpha}$ in the chart $(U_\sigma , G_\sigma , \pi_\sigma) \in \vV$. The family $\tilde{\vV}$ covers $Q$ and we have constructed a family of continuous maps with smooth lifts in every orbifold chart of $\vV$. As $\RR$ is a trivial orbifold, the following data completes the construction of an orbifold map: 
 Choose the quasi-pseudogroup $P \coloneq \Ch_\vV$ which generates $\PSI (\vV)$ and take $\nu \colon \Ch_\vV \rightarrow \PSI (\set{(\RR , \set{\id_\RR} , \id_\RR)}) , f \mapsto \id_{\RR }$. These choices induce a map $(\chi_\alpha , \set{\chi_{\alpha , \sigma}} , [P , \nu])$ which clearly satisfies the requirements of Definition \ref{defn: eq:rofdm} (cf. Remark \ref{rem: psgp}) and $\left(\hat{\chi}_\alpha \coloneq (\chi_\alpha , \set{\chi_{\alpha , \sigma}} , [P , \nu])\right)_{\alpha \in I} \subseteq \Orb{\vV , \set{\id_\RR}}$ is a family of charted orbifold maps.\\
 The construction of $\chi_\alpha$ shows $\tilde{W}^\alpha_1 \subseteq \supp \chi_\alpha \subseteq \pi_\alpha (U_\alpha)$ and the sets $\tilde{W}^1_\alpha$ cover $Q$. Thus the family $\set{\supp \chi_\alpha}_{\alpha \in I}$ covers $Q$  and since $\tilde{\vV}$ is locally finite, this family is locally finite. A quick computation now shows for $x \in Q$: 
	\begin{align*}
	 \sum_{\alpha \in I} \chi_\alpha (x) 	&= \sum_{\alpha \in I, x \in \pi_\alpha (U_\alpha)} \chi_{\alpha, \alpha} \pi^{-1}_{\alpha} (x) = \sum_{\alpha \in I, x \in \pi_\alpha  (U_\alpha)} \frac{ \theta_{\alpha , \alpha}}{\sum_{\delta \in I} \theta_{\delta, \alpha}} (\pi^{-1}_\alpha (x)) \\
					 	&= \sum_{\alpha \in I , x \in \pi_\alpha (U_\alpha)} \frac{\theta_\alpha \pi^{-1}_\alpha (x)}{\sum_{\delta \in I , x \in \pi_\delta (U_\delta)} \theta_\delta \pi_\delta^{-1}(x)} = 1.
	\end{align*}
 The family $(\chi_\alpha)_{\alpha \in I}$ therefore is a partition of unity subordinate to $\vV$. In conclusion, $(\hat{\chi}_\alpha)_{\alpha \in I}$ is a smooth orbifold partition of unity subordinate to $\vV$.
\end{proof}

\begin{nota}\label{nota: opu} \no{nota: opu}
 Let $(Q,\uU)$ be an orbifold with a locally finite representative $\vV$ of $\uU$ indexed by $I$. Consider an orbifold partition of unity $\set{\hat{\chi}_\alpha}_{\alpha \in I}$ subordinate to $\vV$  as in Proposition \ref{prop: pu:hofd}. For any pair $(\alpha , \beta ) \in I \times I$, the lift of $\chi_\alpha$ on $U_\beta$ will be abbreviated as $\chi_{\alpha , \beta}$.
\end{nota}

\thispagestyle{empty}
\section{Tangent Orbibundles and their Sections}\label{sect: tofd}
\setcounter{subsubsection}{0}
In this section, we construct an analogue to tangent manifolds and tangent maps for an orbifold. Tangent orbifolds are well known objects (cf. \cite[Proposition 1.21]{alr2007}). We emphasize that the bundle map associated to a tangent orbifold is a map of orbifolds. This allows us to define orbisections, i.e.\ maps of orbifolds which are sections of the bundle map. In Chapter \ref{sect: lgp:ident}, suitable spaces of orbisections will serve as model space for the diffeomorphism group of an orbifold. Furthermore, it is possible to construct a tangent endofunctor for the category of reduced (smooth) orbifolds. Throughout this section, let $(Q,\uU)$ be an orbifold. We begin with the construction of tangent orbifolds: 

\subsection{The tangent orbifold and the tangent endofunctor}
\begin{con}[Tangent space of an orbifold]\label{con: ts:ofd} \no{con: ts:ofd} 
 Let $p \in Q$ and $(V_i,G_i,\pi_i) \in \uU,\ i\in \set{1,2}$ be arbitrary orbifold charts with $p \in \pi_i (V_i)$. Consider pairs $(\pi_i , v_i) ,\ i=1,2$ where $v_i \in T_{x_i} V_i$ with $x_i \in \pi_i^{-1}(p)$. Notice that by compatibility of orbifold charts, there exist open neighborhoods $x_i \in U_i \subseteq V_i$ and a change of charts $\lambda \colon U_1 \rightarrow U_2$ such that $\lambda (x_1)=x_2$. Identify the tangent spaces $T_{x_i} V_i$ with the corresponding tangent spaces of the open submanifolds $U_i \subseteq V_i$. Since every change of charts is a diffeomorphism, the tangent spaces $T_{x_1} V_1$ and $T_{x_2} V_2$ are isomorphic.\\
 Introduce an equivalence relation on the set of all possible pairs of this kind: We declare two pairs to be equivalent, $(\pi_1 , v_1) \sim (\pi_2 ,v_2)$, if there are open subsets $x_i \in U_i \subseteq V_i$ and a change of charts $\lambda \colon U_1 \rightarrow U_2$ such that $T\lambda (v_1) = v_2$. Here $T\lambda \colon T U_1 \rightarrow T U_2$ is the tangent map of $\lambda$. Since $T \colon \Man \rightarrow \Man$ is a functor ($\Man$ being the category of smooth manifolds), the relation $\sim$ is an equivalence relation. The equivalence class $[\pi , v]$ of $(\pi,v)$ is called a \ind{orbifold!formal tangent vector}{formal orbifold tangent vector} and define the \ind{set!of formal tangent vectors}{\rm set $\tT_p Q$ of all formal orbifold tangent vectors at $p$}.\\
 Consider $x_1 \in \pi^{-1}(p)$, $(U,G, \pi) \in \uU$. The isotropy subgroup $G_{x_1}$ acts on $T_{x_1} U$ via the linear diffeomorphisms $\gamma . v \coloneq T_{x_1} \gamma . v$. Every $\gamma \in G$ is a self-embedding of orbifold charts, whence 
	\begin{equation}\label{eq: sim:gamma}
	 (\pi , v) \sim (\pi , T\gamma . v)\quad \forall \gamma \in G .
	\end{equation} 
 Let $\tilde{v} \in T_{x_1} U /G_{x_1}$ be the equivalence class of $v \in T_{x_1} U$ for $x_1\in \pi^{-1}(p)$. We obtain a bijective map 
	\begin{displaymath}
             k^{x_1}_\pi \colon T_{x_1} U / G_{x_1} \rightarrow \tT_p Q , k^{x_1}_\pi (\tilde{v}) \coloneq T \pi (v) \coloneq [\pi,v]  .
        \end{displaymath}
 To see that this map is indeed injective, consider elements $k^{x_1}_\pi (\tilde{v}) = k^{x_1}_\pi(\tilde{w})$. Thus there is a change of charts $\lambda$ with $T\lambda (v) = w$. By \cite[Lemma 2.11]{follie2003} $\lambda|_O = g|_O$ holds for suitable $g \in G_{x_1}$ on an open neighborhood $O$ of $x_1$ . By definition of $T_{x_1} U /G_{x_1}$ this implies $\tilde{v} = \tilde{w}$.\\ Endow $\tT_p Q$ with the unique topology making the bijection $k_\pi^{x_1}$ a homeomorphism. The space $\tT_p Q$,\glsadd{otsp} is called \ind{orbifold!tangent space}{tangent space of $Q$ at $p$}. We claim that the topology on $\tT_p Q$ neither depends on the choice of charts nor on the preimage $x_1$ in a given chart. Choose some chart $(U,G,\pi)$. As a first step, we prove that the topology does not depend on the choice of the preimage in this chart: 
\\[1em]\textbf{Step 1:} Choose another $x_2 \in \pi^{-1}(p)$. There is some $\gamma \in G$ with $\gamma . x_1 = x_2$. The isotropy groups of $x_1$ and $x_2$ are thus conjugate $\gamma. G_{x_1} \gamma^{-1} = G_{x_2}$. The derived actions of $G_{x_i}$ on $T_{x_i} U, \ i\in \set{1,2}$ are conjugate via the linear isomorphism $T_{x_1} \gamma$, i.e.\ $g.v = T_{x_1}(\gamma^{-1} \circ g \circ \gamma) (v)$ for all $g \in G_{x_2}$. This induces a homeomorphism $\widetilde{T_{x_1}\gamma} \colon T_{x_1} U /G_{x_1}\rightarrow T_{x_2} U/G_{x_2} $. For $v \in T_{x_1} U$, let $\tilde{v}$ be its image in $T_{x_1} U/G_{x_1}$ and compute 
	\begin{displaymath}
	  ((k^{x_2}_\pi)^{-1} \circ k^{x_1}_\pi) (\tilde{v})  = (k^{x_2}_\pi)^{-1} [\pi , v] \stackrel{\eqref{eq: sim:gamma}}{=}(k^{x_2}_\pi)^{-1} [\pi , T_{x_1} \gamma .v]  = \widetilde{T_{x_1} \gamma} (\tilde{v}).
	\end{displaymath}
 Since $\widetilde{T_{x_1} \gamma}$ is a homeomorphism, so is $(k^{x_2}_\pi)^{-1} \circ k^{x_1}_\pi \colon T_{x_1} U/G_{x_1} \rightarrow T_{x_2}/G_{x_2}$. In conclusion the topology on $\tT_p Q$ does not depend on the choice of $x_i \in \pi^{-1}(p)$, whence the index $x_i$ of $k^{x_i}_\pi$ can now be omitted.
\\[1em]\textbf{Step 2:} Consider another chart $(W, H, \psi)$ with $ p\in \psi(W)$ and pick $y \in \psi^{-1} (p)$. By compatibility of charts, there are open subsets $x \in V_U \subseteq U$, $y \in V_W \subseteq W$ and a change of charts homomorphism $\lambda \colon V_U \rightarrow V_W$ with $\lambda (x) = y$. Shrinking the open sets $V_U, V_W$, we may assume that $(V_U, G_{x}, \pi|_{V_U})$ is an orbifold chart and $\lambda$ an open embedding of orbifold charts. This map conjugates (in the sense of Proposition \ref{prop: ch:prop} (a)) the $G_x$-action on $T_x U$ to the $H_y$-action on $T_yW$ again inducing a homeomorphism $\widetilde{T_x \lambda} \colon T_x V_U /G_x \rightarrow T_y V_W/H_y$. As in Step 1, a well-defined homeomorphism is given by 
		\begin{displaymath}
                  k_{\psi} \circ k^{-1}_\pi \colon T_x U/G_x \rightarrow T_y /H_y , \tilde{v} \mapsto \widetilde{T\lambda}(\tilde{v}).
		\end{displaymath}
Therefore the topology on $\tT_p Q$ is independent of the choice of charts. 
\end{con}

\begin{rem}\label{rem: loc:tb} \no{rem: loc:tb}
 Let $(U,G,\pi)$ be an orbifold chart with $p \in \im \pi$. The homeomorphism $\tT_p Q \cong T_x U/ G_x$ for $x \in \pi^{-1} (p)$ allows us to think of $\tT_p Q$ as an orbifold. In particular, the tangent space $\tT_p Q$ may be identified with a convex cone. In contrast to tangent spaces of manifolds, the tangent spaces of an orbifold will not be vector spaces. Nevertheless, each orbifold tangent space contains a zero element $0_p \coloneq [\pi ,0_x]$, where $(U,G,\pi)$ is a chart with $p = \pi (x)$ and $0_x \in T_x U$ the zero element.\\
 In the manifold case, our definition boils down to: The tangent space of a manifold (considered as a trivial orbifold) at $p$ is the tangent space of the manifold at $p$.
\end{rem}

\begin{defn}[Tangent orbifold]\label{defn: tofd} \no{defn: tofd}
 Consider the set $\tT Q \coloneq \bigcup_{p \in Q} \tT_p Q$. Since the tangent spaces are mutually disjoint, we derive a well-defined map \glsadd{pito}
	\begin{displaymath}
	 \pi_{\tT Q} \colon \tT Q \rightarrow Q , [\psi , v] \mapsto \psi (x)\text{, where }v \in T_x \dom \psi.
	\end{displaymath}
 If $(U, G ,\psi) \in \uU$ is an arbitrary chart, then $G$ acts on $TU$ via the derived action $\gamma . X \coloneq T\gamma (X)$. Define $\Pi \colon TU \rightarrow TU/G$ to be the quotient map to the orbit space with respect to this action. Using the notation of Construction \ref{con: ts:ofd}, we obtain a map for $(U,G,\psi) \in \uU$:
 \begin{displaymath}
   T \psi \colon TU \rightarrow \tT Q , v \mapsto [\psi ,v] 
 \end{displaymath} 
 In particular, each $v \in T_x U$ is mapped to some $[\psi , v] \in \tT_{\psi (x)} Q$. Choose an atlas $\aA \in \uU$. We equip $\tT Q$ with the final topology with respect to the family $(T\psi)_{(U,G,\psi) \in \aA}$.\\
 This topology induces a canonical orbifold structure on $\tT Q$. An atlas for this orbifold is given by the family $(TU,G,T\psi)$,\footnote{Notice that we should have written $\setm{Tg}{g \in G}$ instead of $G$ in the definition of $(TU,G,T\psi)$. Definition \ref{defn: moerd:ofc} requires the acting group to be a subgroup of $\Diff (TU)$ which is only satisfied by $\setm{Tg}{g \in G}$. However, we use the canonical identification $G \cong \setm{Tg}{g \in G}$ to justify the shorter (but in fact incorrect) notation $(TU, G,T\psi)$.} where $(U,G,\psi)$ runs through $\aA$. The $G$-action of the chart $(TU,G,T\psi)$ is the derived action of $G$, i.e. $\gamma . v \coloneq T\gamma (v)$. With respect to this structure $\pi_{\tT Q}$ induces an orbifold map. Its lifts are given by the bundle projections $TU \rightarrow U$, for ${(U,G,\pi) \in \uU}$.\\ 
 We define the \ind{}{tangent orbifold} $\tT (Q,\uU)$ of $(Q,\uU)$. It is the orbifold $(\tT Q , \tT \uU)$, where $\tT\uU$ is the orbifold structure induced by $\tT\aA$. A proof for the details of this construction will be given in the next lemma.  
\end{defn}

\begin{lem}\label{lem: pr:tofd} \no{lem: pr:tofd} Let $(Q,\uU)$ be an $n$-dimensional orbifold. Using the notation of Definition \ref{defn: tofd}, the following statements hold:
 \begin{compactenum}
   \item Let $(U,G,\psi), (V,H,\vp) \in \uU$ and $\lambda \colon U \supseteq W \rightarrow W' \subseteq V$ be a change of charts. Its tangent map $T\lambda \colon TW \rightarrow TW'$ is a diffeomorphism with $T\vp T\lambda = T\psi|_{TW}$.
  \item For any chart $(U,G,\psi) \in \uU$ we set $\tilde{U} \coloneq \psi (U)$ and $T \tilde{U} \coloneq \im T \psi $. Then $T\tilde{U}$ is an open set in $\tT Q$ and $T\psi$ is an open map. 
  \item The topology on $\tT Q$ does not depend on the choice of the atlas $\aA \in \uU$ in Definition \ref{defn: tofd}.
  \item For each $\aA \in \uU$, the set $\tT\aA \coloneq \setm{(TU,G, T\psi ) }{(U, G,\psi) \in \aA}$ is an orbifold atlas for $\tT Q$.\glsadd{TU} The orbifold charts in this atlas are compatible via the changes of charts computed in (a).
  \item The map $\pi_{\tT Q} \colon \tT Q \rightarrow Q , [\psi ,v] \mapsto \psi (x) , \ v \in T_x U$ is continuous and $\tT Q$ is a Hausdorff paracompact space. In conclusion, $\tT(Q,\uU)$ is an orbifold. 
  \item $\pi_{\tT Q}$ induces a morphism of orbifolds $\tpi \in \ORBM[\tT(Q,\uU) , (Q,\uU)]$.
  \item The topology on $\tT Q$\glsadd{tofd} induces on each $\tT_p Q$ the topology obtained in Construction \ref{con: ts:ofd}.
 \end{compactenum}
\end{lem}

\begin{proof}
 \begin{compactenum}
    	\item For the change of charts $\lambda$, the tangent map $T \lambda \colon TW \rightarrow TW'$ is a diffeomorphism. It suffices to prove the commutativity for each element of $T_r W$, where $r \in W$ is arbitrary. Since $\lambda$ is a change of charts, $\vp \lambda = \psi|_{\dom \lambda }$ holds. The definition of $\tT_{\psi (r)} Q$ yields $[\psi , v] = [\vp , T\lambda (v)]$. We obtain for $v \in T_r W$ the identity
	\begin{displaymath}
 	T\vp T \lambda (v) = [\vp , T \lambda (v)] = [\psi , v] = T\psi (v).
	\end{displaymath}
	\item The space $\tT Q$ is endowed with the final topology with respect to the mappings $T\pi$, where $(V,H,\pi)$ runs through $\aA$. To prove the assertion we need to show that $(T \pi)^{-1} (T \psi (V))$ is an open set for every $(W,H,\pi) \in \aA$ and open set $V \subseteq TU$. Define the \glsadd{CH}\ind{set!of changes of charts}{set of changes of charts from $U$ to $W$}: 
		\begin{displaymath}
	         \CH{U}{W} \coloneq \setm{\lambda \colon U \supseteq \dom \lambda \rightarrow \cod \lambda \subset W}{\lambda \text{ is a change of charts}}.
		\end{displaymath}
 	Then one computes its preimage as 
   	\begin{align*}
 	(T \pi)^{-1} (T\psi (V))  	&= \setm{w \in TW}{ [\pi,w] \in T \psi (V)} \\
					&= \setm{w \in TW}{\exists \lambda \in \CH{U}{W},  w = T\lambda (v) \text{ for some } v \in V} \\
					&= \bigcup_{\lambda \in \CH{U}{W}} T\lambda (\dom T \lambda \cap V) \subseteq TW.
        \end{align*}
	Each $T\lambda$ is a diffeomorphism onto its (open) image, whose domain is an open set. Thus every $T\lambda (\dom T\lambda \cap V)$ is an open subset in $TW$. This proves $(T \pi)^{-1} (T\psi (V))$ to be an open set, whence $T\psi$ is an open map with open image $T\tilde{U}$ in $\tT Q$. 
        \item To see that the topology does not depend on the choice of $\aA$, we consider the final topology $\oO'$ on $\tT Q$ with respect to the mappings $T\psi$, where $(U,G,\psi)$ runs through an atlas $\aA' \in \uU$. It suffices to prove that the topologies coincide if $\aA \subseteq \aA'$ holds. Thus without loss of generality the final topology $\oO$ with respect to $\aA$ is finer than the topology $\oO'$. However, the computation in (b) shows that $\oO'$ is also finer than $\oO$, whence $\oO = \oO'$ follows and the topology does not depend on the choice of $\aA$.
     	\item If $(U, G, \phi) \in \aA$ is an arbitrary chart, then $T\phi$ has an open image by (b). Consider the map $T\overline{\phi} \colon TU /G \rightarrow \im T\phi , v \mapsto [\phi , v]$. Combining Proposition \ref{prop: ch:prop} with the definition of the equivalence relation in Construction \ref{con: ts:ofd}, this map is a well-defined bijective map. We may factor $T\phi$ as $T\phi = T\overline{\phi} \circ \Pi$, where $\Pi$ is the quotient map to the orbit space associated to the $G$ action on $TU$. Since $\Pi$ is a quotient map and $T\phi$ is continuous, $T\overline{\phi}$ is continuous. If $V\subseteq TU/G$ is an open set, then $\Pi^{-1}(V)$ is an open set. Since $T\phi$ is open by (b) the set $T\overline{\phi}(V) = T\phi \Pi^{-1} (V)$ is an open set. Thus $T\overline{\phi}$ is an open map and in conclusion $T\phi$ may be factored as the quotient map to the orbit space associated to the group action composed with a homeomorphism. In particular, the set of orbifold charts
	\begin{displaymath}
		\tT \aA \coloneq \setm{(TU,G , T\pi )}{(U, G, \pi) \in \aA}
	\end{displaymath}
  	covers $\tT Q$. In (a), we have constructed a family of maps which are change of chart maps for $\tT\aA$. Using this family of changes of charts, the definition of the chart maps and tangent spaces $\tT_pQ$ shows that each pair of orbifold charts in $\tT\aA$ is compatible. Thus $\tT\aA$ is an orbifold atlas inducing a unique orbifold structure $\tT (Q,\uU)$ of dimension $2\cdot \dim (Q,\uU)$ on $\tT Q$.
  	\item The definitions of $\pi_{\tT Q}$ and $\tT Q$ together with the compatibility of orbifold charts yield $\pi_{\tT Q}^{-1} (\psi (U)) = T\psi (TU)$, for every $(U,G,\psi) \in \uU$. Hence the preimages of a basis of the topology under $\pi_{\tT Q}$ are open (cf. Lemma \ref{lem: ofdch:bt}) and thus $\pi_{\tT Q}$ is continuous by \cite[Proposition 1.4.1.]{Engelking1989}.\\[1em]
	\emph{The space $\tT Q$ is a Hausdorff space}: Let $x,y \in \tT Q$ be distinct points.\\ 
        First case: $\pi_{\tT Q} (x) \neq \pi_{\tT Q} (y)$. There are orbifold charts $(U_x,G_x,\psi_x), (U_y,G_y,\psi_y) \in \uU$ such that $\pi_{\tT Q} (x) \in \psi_x (U_x)$, $\pi_{\tT Q} (y) \in \psi_y (U_y)$ and $\psi_x (U_x) \cap \psi_y (U_y) = \emptyset$ hold. As the images of these charts do not intersect, the set $\Ch_{U_x,U_y}$ is empty. By construction of the equivalence relation, $T\psi_x (TU_x) \cap T\psi_y (TU_y) = \emptyset$. Hence $x \in \pi_{\tT Q}^{-1} (\psi_x (U_x))$ and $y \in \pi_{\tT Q}^{-1} (\psi_y (U_y))$ are contained in disjoint open sets.\\
        Second case: $\pi_{\tT Q} (x) = \pi_{\tT Q} (y)$. Choose any orbifold chart $(U,G,\psi)$ with $\pi_{\tT Q} (x) \in \psi (U)$. Then $x,y \in \pi_{\tT Q}^{-1} (\psi (U)) = T \psi (TU)$. Both $x$ and $y$ are contained in $T\psi (TU)$, which is homeomorphic to the orbit space $TU /G$. This space is Hausdorff by Lemma \ref{lem: orbitmap} and there are disjoint open subsets $x\in V_x, y \in V_y$ of $T\psi (TU)$. As $T\psi (TU)$ is open, both points are contained in disjoint open subsets of $\tT Q$. In conclusion the space $\tT Q$ is a Hausdorff space.\\[1em]
        \emph{The space $\tT Q$ is paracompact}: Connected components of $\tT Q$ are open and closed, therefore \cite[Theorem 5.1.35]{Engelking1989} implies that $Q$ will be paracompact if each connected component of $\tT Q$ is paracompact.
        We claim that each connected component $C$ of $\tT Q$ is second countable. If this is true, paracompactness of a component is assured by the following observations:
        The quotient map to an orbit space preserves locally compact spaces by Lemma \ref{lem: orbitmap}. Thus $\tT Q$ is locally compact, hence a regular space. Combining \cite[Theorem 3.8.1]{Engelking1989} and \cite[Theorem 5.1.2]{Engelking1989} second countability of a component implies paracompactness of that component.\\
	\textbf{Proof of the claim:} Every component $C' \subseteq Q$ is second countable (cf.\ Proposition \ref{prop: ofd:prop}). The continuous map $\pi_{\tT Q}$ maps $C$ into some component $C' \subseteq Q$. Since $C'$ is second countable, there is a countable base $\bB$ of the topology on $C'$. The images of orbifold charts in $C'$ also form a base of the topology by Lemma  \ref{lem: ofdch:bt}. Thus without loss of generality $\bB$ contains only (open) images of a set of orbifold charts $\rR = \setm{(U_i , G_i, \pi_i)}{i \in I}$ in $\uU$. By construction of $\pi_{\tT Q}$, the countable family of open sets $\tT \bB \coloneq (T\pi_i(TU_i))_{(U_i, G_i,\pi_i) \in \rR}$ covers $C$. Observe that $T\tilde{U}_i \cong TU_i/G_i$ and $TU_i$ is the tangent manifold of a connected paracompact manifold, thus connected paracompact and second countable by Proposition \ref{prop: para:coco}. The quotient map to the orbit space is continuous and open by Lemma \ref{lem: orbitmap} which implies that $T\tilde{U}_i$ is also second countable. As a countable union of open and second countable spaces, $C$ is second countable.
  	\item The map $\pi_{\tT Q}$ is continuous by (e) and we have to construct lifts for $\pi_{\tT Q}$: Consider an arbitrary orbifold chart $(TU , G, T\psi)  \in \tT \uU$. Let $\pi_{TU} \colon TU\rightarrow U$ be the bundle projection of the tangent bundle. This map is smooth, and we obtain a commutative diagram:
  		\begin{displaymath} 
   			\xymatrix{
 				TU \ar[rr]^{T\psi} \ar[d]^{\pi_{TU}} 	&	 & T\tilde{U} \ar[d]^{\pi_{\tT Q}} \\
				U \ar[rr]^{\psi}					&	 & \tilde{U}
 				 }
  		\end{displaymath}
  	Choose a representative $\aA \in \uU$ and define $P_{\pi_{\tT Q}} \coloneq \bigcup_{(U,W) \in \aA \times \aA} \setm{T\lambda}{\lambda \in \CH{U}{W}}$. We have to show that the quasi-pseudogroup  $P_{\pi_{\tT Q}}$ generates $\PSI (\tT \aA)$. Let $\vp \in \PSI (\tT \aA)$ and pick an arbitrary $v \in \dom \vp$. Then there are $(TU, G, T\pi), (TV,H,T\psi) \in \tT \aA$ and an open set $v \in \Omega \subseteq T U$ such that $\vp|_{\Omega}$ is a diffeomorphism onto an open set $\Omega' \subseteq TV$ which contains $w \coloneq \vp(v)$. Since $T\psi (w) = T\pi (v)$ holds, the equivalence relation shows that there are open sets $x \in W \subseteq U$, $y \in W' \subseteq V$ and a change of charts $\lambda \colon W \rightarrow W'$ such that $v \in T_x W,\ w \in T_y W'$ and  $T\lambda (v) = w$. Shrinking $W$ and $W'$ we may assume that $T\lambda \colon TW \rightarrow TW'$ is an embedding of orbifold charts. Thus on $TW$, the maps $T\lambda$ and $\vp|_{TW}$ are embeddings of orbifold charts. By Proposition \ref{prop: ch:prop}, there is an $h \in H_{w}$ such that $h . T\lambda = \vp|_{TW}$. The definition of the group action on charts in $\tT \aA$ yields $\vp|_{TW}= h . T\lambda = T (h\circ \lambda)$. Now $h\circ \lambda \in \PSI (\aA)$ implies $T(h.\lambda) \in P_{\tpi}$. In conclusion, $P_{\pi_{\tT Q}}$ generates $\PSI (\tT \aA)$. Define the map 
	\begin{displaymath}
	 \nu_{\pi_{\tT Q}} \colon P_{\pi_{\tT Q}} \rightarrow \PSI (\aA) , T\lambda \mapsto \lambda.
	\end{displaymath}
  	By construction, this map satisfies (R4a)-(R4d) of Definition \ref{defn: rep:ofdm} and therefore 
	\begin{displaymath}
	 \tpi \coloneq (\pi_{\tT Q} , \setm{\pi_{TU}}{(U,G, \pi)\in \aA} , [P_{\pi_{\tT Q}}, \nu_{\pi_{\tT Q}}]) \in \Orb{\tT \aA,\aA}
	\end{displaymath}
 	is a representative of an orbifold map. We call $\tpi$ the \ind{orbifold map!bundle projection}{bundle projection}. By abuse of notation, we let $\tpi$ also be the equivalence class of the charted map $\tpi$ in $\ORB (\tT(Q,\uU), (Q,\uU))$. Clearly any choice of $\aA$ in the above construction yields the same class $\tpi$.  In particular, for each chart $(TU,G,T\psi)$ in $\tT \uU$ there is a representative of $\tpi$ such that the bundle projection $\pi_{TU} \colon TU \rightarrow U$ is a local lift of $\tpi$. The triple $(\tT (Q,\uU), (Q,\uU),  \tpi)$ is an orbibundle, the \ind{}{tangent orbibundle} (cf.\ \cite[p.14]{alr2007}).  
  	\item Choose some orbifold chart $(U,G,\psi) \in \uU$ such that $p \in \psi (U)$. Shrinking the chart, we may assume $\set{z} = \psi^{-1} (p)$, i.e.\ $G \cong \Gamma_{p}$. By construction, $\tT_p Q \subseteq T \psi (TU)$ holds. Recall from (c) that $T\psi = T\overline{\psi} \circ \Pi$, where $\Pi$ is the quotient map to the orbit space with respect to the $G$-action on $TU$ and $T\overline{\psi}$ is a homeomorphism. Observe $(T\overline{\psi})^{-1} (\tT_p Q) = \Pi (T_z U)$. Notice that for manifolds the subspace topology of $T_z U \subseteq TU$ coincides with the usual topology of $T_z U$. As the quotient map to an orbit space is open, \cite[VI. Theorem 2.1]{dugun1966} proves that the subspace topology of $(T\overline{\psi})^{-1} (\tT_p Q)$ and the quotient topology on $\Pi (T_z U) = T_zU /G$ coincide. In Construction \ref{con: ts:ofd}, $\tT_p Q$ has been endowed with precisely the same topology. Hence the induced topology on $\tT_p Q$ coincides with the one from Construction \ref{con: ts:ofd}.
\end{compactenum}\vspace{-2em}
\end{proof}
\noindent
Notice that for any trivial orbifold (i.e. for a manifold), the tangent orbibundle coincides with the tangent bundle of the manifold. For a non-trivial orbifold, an explicit example of a tangent orbifold will be computed in Example \ref{ex: ad:nofd}. \\
Mappings to the tangent orbifold admit representatives which are charted maps whose range atlas is $\tT \aA$ for some $\aA \in \uU$. Thus orbifold maps to the tangent orbifold always posses representatives which may be computed in the canonical orbifold charts of the tangent orbifold.

\begin{lem}\label{lem: oto:nmap} \no{lem: oto:nmap}
Let $[\hat{f}] \in \ORBM[(Q,\uU), \tT(Q,\uU)]$ be an arbitrary orbifold map. There is a representative $\hat{f} \in [\hat{f}]$ such that the range atlas of $\hat{f}$ is contained in $\tT \wW$ for some $\wW \in \uU$. In other words, $\hat{f}$ is a charted orbifold map with $\hat{f} \in \Orb{\vV , \tT \wW}$, where $\vV$ and $\wW$ are some representatives of $\uU$.
\end{lem}

\begin{proof}
 Let $[\hat{f}]$ be as above. Consider the composition  $\tpi \circ [\hat{f}]$ of $[\hat{f}]$ with the bundle projection $\tpi$ (Lemma \ref{lem: pr:tofd}). Reviewing \cite[Lemma 5.17]{pohl2010} (cf.\ Section \ref{sect: ofd_cat}), the composition in $\ORB$ is induced by the composition of representatives of the equivalence classes. Fix a representative $\tpi \in \Orb{\tT \wW , \wW}$ for some $\wW \in \uU$. Then there are representatives $\vV , \vV''$ of $\uU$ respectively a representative $\vV'$ of $\tT \uU$ together with the following charted orbifold maps: $\hat{g} \in \Orb{\vV , \vV'}$ with $\hat{g} \in [\hat{f}]$ and $\hat{h} \in \Orb{\vV',\vV''}$ with $\hat{h} \in \tpi$ such that these maps induce the composition, i.e. $\tpi \circ [\hat{f}] = [\hat{h} \circ \hat{g}]$. Furthermore, these charted maps can be chosen such that the following diagram is commutative: \begin{displaymath}
	 \begin{xy}
  \xymatrix{					
				  &							& \tT \wW \ar[rr]^{\tpi} & & \wW &    \\
      \vV \ar[r]^{\hat{g}} 	  & \vV' \ar[ru]_{\ve_1} \ar[rrrr]^{\hat{h}}  		& & & & \vV'' \ar[lu]^{\ve_2}}
\end{xy}
	\end{displaymath} 
 Here the charted maps $\ve_1$ and $\ve_2$ are lifts of the identity (cf. Definition \ref{defn: orb:cat}) and composition in the diagram is composition of charted orbifold maps. By definition of the composition in $\ORB$ we obtain $[\hat{f}] = [\ve_1 \circ \hat{g}]$ with a representative $\ve_1 \circ \hat{g} \in \Orb{\vV, \tT \wW}$. 
\end{proof}

The rest of this section will be devoted to construct a tangent functor for the category $\ORB$. To achieve this goal, we have to construct tangent orbifold maps. We record several observations, which will allows us to introduce tangent orbifold maps. 
\begin{rem}\label{rem: ch:tofd} \no{rem: ch:tofd}
  \begin{compactenum}
   \item Let $\vV$ be a representative of $\uU$ for an orbifold $(Q,\uU)$. The $G$-action in a chart in $\vV$ acts on the tangent chart via the derived action. Since the tangent functor $T \colon \Man \rightarrow \Man$ (where $\Man$ is the category of smooth (not necessarily finite dimensional) manifolds) is functorial, Proposition \ref{prop: ch:prop} (e) and the definition of the tangent manifold imply that $\tT \PSI (\vV) \coloneq \setm{T\lambda}{\lambda \in \PSI (\vV)}$ is a quasi-pseudogroup which generates $\PSI (\tT \vV)$. Furthermore, if $P$ is some quasi-pseudogroup which generates $\PSI (\vV)$, then the quasi-pseudogroup $\tT P \coloneq \setm{T\lambda}{\lambda \in P}$ generates $\PSI (\tT \vV)$.
   \item Let $\lambda , \mu \in \CH{V}{W}$ be change of charts and $X \in \dom T\lambda \cap \dom T\mu$ such that $\germ_X T\lambda = \germ_X T\mu$ holds. Choose an open $X$-neighborhood $U_X \subseteq TV$ with $T\lambda|_{U_X} = T\mu|_{U_X}$. This implies $\lambda|_{\pi_{TV}(U_X)} = \mu|_{\pi_{TV}(U_X)}$. Since $\pi_{TV}$ is an open map, $\pi_{TV} (U_X)$ is open and contains $\pi_{TV} (X)$. Thus $\germ_{\pi_{TV} (X)} \lambda = \germ_{\pi_{TV} (X)} \mu$ holds. 
  \end{compactenum}
\end{rem} 

\begin{defn}\label{defn: tofd:map} \no{defn: tofd:map}
 Let $(Q_i,\uU_i)$, $i=1,2$ be orbifolds and $[\hat{f}] \in \ORBM[(Q_1,\uU_1),(Q_2,\uU_2)]$ be a morphism with representative $\hat{f} = (f, \set                                                                                                                                                                                                                         {f_i}_{i \in I}, [P,\nu]) \in \Orb{\vV,\wW}$.\\
 Furthermore, consider orbifold atlases $\vV = \setm{(V_i,G_i,\psi_i)}{i\in I}$ and $\wW = \setm{(W_j,H_j,\varphi_j)}{j \in J}$. For two changes of charts $T\lambda = T\mu$ is satisfied if and only if $\lambda = \mu$, whence $\tT \nu \colon \tT P \rightarrow \PSI (\tT \wW) , T\lambda \mapsto T\nu (\lambda)$ is a well defined map. Here $\tT P$ is the quasi-pseudogroup of some $(P,\nu)$ in the class $[P,\nu]$ as in Remark \ref{rem: ch:tofd} (a). The class $[\tT P ,\tT \nu]$ does not depend on the choice of $(P,\nu)$ in $[P,\nu]$ by the definition of equivalence (cf.\ Definition \ref{defn: eq:rofdm}).\\
 Combining Remark \ref{rem: ch:tofd} (b) and the properties (R4a)-(R4d) of Definition \ref{defn: rep:ofdm} for the map $\nu$ with respect to $F \coloneq \coprod_{i \in I} f_i$, we see that $\tT \nu$ satisfies properties (R4a)-(R4d) with respect to $F' \coloneq \coprod_{i\in I} Tf_i$. In particular, we derive $T\varphi_{\alpha (i)} Tf_i (T\lambda .x ) = T\varphi_{\alpha (j)} Tf_j (x)$ for each $\lambda \in \CH{V_j}{V_i}$. Thus there is a well-defined continuous map $\tT f \colon \tT Q_1 \rightarrow \tT Q_2, \tT f (x) \coloneq T\varphi_{\alpha (i)} Tf_i T\psi_i^{-1} (x), \ x \in \im T\psi_i$.\\
 In conclusion, a charted map of orbifolds is given by
  \begin{displaymath}
    \widehat{\tT f} \coloneq (\tT f, \set{Tf_i}_{i\in I}, [\tT P , \tT\nu]) \in \Orb{\tT \vV , \tT \wW}.
  \end{displaymath}
 The map $\widehat{\tT f}$ is a representative of the \ind{orbifold map!tangent map}{orbifold tangent map} $[\widehat{\tT f}]$ of $[\hat{f}]$. We have to check that the construction of this map is functorial.  
\end{defn}

\begin{lem}\label{lem: tmap:wd} \no{lem: tmap:wd} The assignment $\tT \colon \ORB \rightarrow \ORB , (Q,\uU) \mapsto \tT (Q,\uU), [\hat{f}] \mapsto [\widehat{\tT f}]$ is a functor, i.e. 
 \begin{compactenum}
  \item If $\hat{\ve} = (\id_Q ,\set{f_i}_{i\in I}, [P,\nu]) \in \Orb{\vV, \wW}$ is a lift of the identity $\ido{}$, then $\widehat{\tT \ve}$ is a lift of the identity $\id_{\tT (Q, \uU)}$. 
  \item Let $\hat{f} = (f,\set{f_i}_{i \in I}, [P_f,\nu_f]) \in \Orb{\vV, \wW}$ and $\hat{g} = (g,\set{g_j}_{j\in J}, [P_g,\nu_g]) \in \Orb{\wW, \aA}$. Then $\widehat{\tT g \circ f} = \widehat{\tT g} \circ \widehat{\tT f}$.
  \item Two representatives $\hat{f_1}, \hat{f_2}$ of $[\hat{f}] \in \ORBM{}$ induce equivalent charted orbifold maps, i.e. $[\widehat{\tT f_1} ] = [\widehat{\tT f_2}]$.
  \item $[\widehat{\tT g \circ f}] = [\widehat{\tT g}] \circ [\widehat{\tT f}]$ holds for $[\hat{f}] \in \ORBM[(Q_1,\uU_1),(Q_2,\uU_2)]$, $[\hat{g}] \in \ORBM[(Q_2,\uU_2) ,(Q_3,\uU_3)]$. 
 \end{compactenum}
\end{lem}

\begin{proof}
\begin{compactenum}
  \item For each $i \in I$ let the lifts $f_i \colon V_i \rightarrow W_{\alpha (i)}$ be given with respect to the charts $(V_i,G_i,\psi_i)$ and $(W_{\alpha (i)}, H_{\alpha (i)}, \varphi_{\alpha (i)})$. Here $\alpha \colon I \rightarrow J$ is the map which assigns to $f_i$ the chart $W_{\alpha (i)}$. Each $f_i$ is a local diffeomorphism by Definition \ref{defn: idmor}. Using functoriality of $T$, again $Tf_i$ is a local diffeomorphism. By Proposition \ref{prop: ident:eq} the assertion will be true if $\tT \id_Q = \id_{\tT Q}$ holds. Consider $x \in \tT Q$ with $x \in \im T\psi_i$ for some $i\in I$. Choose $z_x \in TV_i$ with $T\psi_i (z_x) = x$ and observe that by Proposition \ref{prop: ll:id}, we may choose orbifold charts $(S_x, G_x, \psi_x|_{S_x})$ and $(S_x', G_x', \psi_x|_{S_x'})$ with $\pi_{TV_i} (x) \in S_x$ such that $f_i$ induces the identity on $S_x$ with respect to $\id_{S_x}$ and $(f_i|_{S_x})^{-1}$. Hence $f_i|_{S_x}$ is a change of charts, which implies $\tT \id_Q (x) = \tT \id_Q (T\psi_i (z_x)) = T\varphi_{\alpha (i)} Tf_i (z_x) = T\varphi_{\alpha (i)} T(f_i|_{S_x}) (z_x) = x$. 
  \item Define $h_i \coloneq g_{\alpha (i)} \circ f_i$ and $h = g\circ f$. Then $\hat{g} \circ \hat{f}$ is given by $\hat{h} = (h, \setm{h_i}{i\in I}, [P_h,\nu_h])$. From Definition \ref{defn: tofd:map}, we infer $\widehat{\tT (g \circ f)} = (\tT h , \setm{Th_i}{i\in I} , [\tT P_h ,\tT \nu_h])$.\\ 
  By construction, one has $\widehat{\tT f} \in \Orb{\tT \vV , \tT \wW}$ and $\widehat{\tT g} \in \Orb{\tT \wW, \tT \aA}$. These charted orbifold maps may therefore be composed as in Construction \ref{con: cp:chom}: The charted orbifold map $\widehat{\tT f} \circ \widehat{\tT g}$ is given as $\hat{h}_\tT \coloneq (\tT g \circ \tT f , \set{Tg_{\alpha (i)} \circ Tf_i}_{i \in I}, [P_{h_\tT}, \nu_{h_\tT}])$. By functoriality of $T$, we have $h_i = T (g_{\alpha (i)} \circ f_i) = Tg_\alpha (i) Tf_i $ for $i\in I$. Hence the lifts of $\widehat{\tT (g\circ f)}$ and $\hat{h}_\tT$ coincide for each $i\in I$. We conclude $\tT h = \tT g \circ \tT f$. \\
  If $(\tT P_{h}, \tT \nu_{h}) \sim (P_{h_\tT},\nu_{h_\tT})$ holds, then both maps will be equivalent as charted orbifold maps. By construction of the quasi-pseudogroups this indeed follows directly from the functoriality of $T$ and property (R4b) of Definition \ref{defn: rep:ofdm}. However, since quasi-pseudogroups work with the germs of maps, the computation has to be carried out on the germ level. Here are the technical details:\\
  Let $\lambda , \mu \in \CH{TV_i}{TV_j}, i,j \in I$, $\lambda \in \tT P_h$, $\mu \in P_{h_\tT}$ and $X \in \dom \lambda \cap \dom \mu$ with $\germ_X \lambda = \germ_X \mu$. To establish the equivalence, we have to prove the identity 
  \begin{equation}\label{eq: goal}
   \germ_{Th_i (X)} \tT \nu_h (\lambda) = \germ_{Th_i (X)} \nu_{h_\tT} (\mu).
  \end{equation}
  Set $x \coloneq \pi_{TV_i} (X)$. By definition of the quasi-pseudogroups of $\hat{f}$ and $\hat{g}$ (combine Remark \ref{rem: ch:tofd} and Construction \ref{con: cp:chom}), we obtain the following data: 
  \begin{compactitem}
	\item[1.] $\eta, \rho \in P_f$, $x \in U_{\eta ,x}, U_{\rho, x}$ open and $\eta|_{U_{\eta , x}} , \rho|_{U_{\rho , x}}  \in P_h$ with $\lambda = T \eta|_{U_{\eta , x}}$ and $\germ_X \mu = \germ_X T\rho$,
	\item[2.] $\xi_{\eta ,x}, \xi_{\rho,x} \in P_g$ with $\nu_h (\eta|_{U_{\eta , x}}) = \nu_g (\xi_{\eta ,x})$ and $\germ_{f_i (x)} \xi_{\eta ,x} = \germ_{f_i (x)} \nu_f (\eta)$, respectively for $\nu_h (\rho|_{U_{\rho , x}}) = \nu_g (\xi_{\rho ,x})$ and $\germ_{f_i (x)} \xi_{\rho ,x} = \germ_{f_i (x)} \nu_f (\rho)$
	\item[3.] $\xi_{\mu , X} \in \tT P_g$ with $\nu_{h_\tT} (\mu) = \tT \nu_{g} (\xi_{\mu, X})$ and $\germ_{Tf_i (X)} \xi_{\mu ,X} = \germ_{Tf_i (X)} \tT \nu_f (T\rho)$.
 \end{compactitem}
 For $\phi, \psi \in P_f$ and $z \in \dom \phi \cap \dom \psi$ Remark \ref{rem: ch:tofd} (b) implies $\germ_z \phi = \germ_z \psi$ if and only if $\germ_X T\phi = \germ_X T\psi$ for some $X \in T_z V_i$. Exploiting property (R4b) for $\nu_f$ we obtain $\germ_{f_i (x)} \nu_f (\phi ) = \germ_{f_i (x)} \nu_f (\psi)$, whence $\germ_{Tf_i (X)} \tT \nu_f (T\phi ) = \germ_{Tf_i (X)} \tT\nu_f (T\psi)$ holds. Analogously the same holds for $\nu_g$ and $\nu_h$ by 1. and 2.: 
  \begin{displaymath}
   \germ_{Th_i (X)} \tT \nu_h (\lambda) = \germ_{Th_i (X)} T \nu_h (\eta|_{U_{\eta,x}}) = \germ_{Th_i (X)} T \nu_g (\xi_{\eta,x}).
  \end{displaymath}
 We already know $\germ_X T\eta = \germ_X \lambda = \germ_X \mu = \germ_X T\rho$ and by Remark \ref{rem: ch:tofd} (b) $\germ_x \eta = \germ_x \rho$ follows. Using property (R4b) for $\nu_f$ and 2. one obtains $\germ_{f_i (x)} \xi_{\eta,x} = \germ_{f_i (x)} \nu_f(\eta) = \germ_{f_i (x)} \nu_f(\rho)$.\\ 
 Together with 3. this yields $\germ_{Tf_i (X)} T\xi_{\eta,x} = \germ_{Tf_i (X)} \tT \nu_f (T \rho) = \germ_{Tf_i (X)} \xi_{\mu ,X}$. Again by 3. and property (R4b) for $\tT \nu_g$ we derive: 
  \begin{displaymath}
   \germ_{Th_i (X)} \tT \nu_h (\lambda) = \germ_{Th_i (X)}  \tT \nu_g ( T\xi_{\eta,x}) = \germ_{Th_i (X)}  \tT \nu_g (\xi_{\mu, X})= \germ_{Th_i(X)} \nu_{h_\tT} (\mu).
  \end{displaymath}
 Since $X$, $\lambda,\mu$ were arbitrary, $(\tT P_h , \tT \nu_h) \sim (P_{h_\tT} ,\nu_{h_\tT})$ holds and we conclude $\widehat{\tT (g\circ f)} = \widehat{\tT g} \circ\widehat{\tT f}$.  
  \item In view of of (a) and (b), we can apply $\tT$ to the diagram \eqref{eq: comp:ofm} which defines the equivalence of charted orbifold maps (cf. Definition \ref{defn: sim:ofd} and the assertion follows.
  \item This is just the combination of (b) and (c).
\end{compactenum}
\end{proof}

\begin{rem}
 Let $(Q_i,\uU_i), i\in \set{1,2}$ be orbifolds and $[\hat{f}] \in \ORBM[(Q_1,\uU_1),(Q_2,\uU_2)]$. The definition of the tangent orbifold map implies that the following diagramm is commutative 
  \begin{displaymath}
   \begin{xy}
  \xymatrix{
      \tT (Q_1,\uU_1) \ar[rr]^{\tT [\hat{f}]} \ar[d]_{\tpi[Q_1,\uU_1]}    &&   \tT (Q_2, \uU_2) \ar[d]^{\tpi[Q_2,\uU_2]}  \\
      (Q_1,\uU_1) \ar[rr]_{[\hat{f}]}             &&   (Q_2,\uU_2)     
   } 
  \end{xy}
  \end{displaymath}
In other words, the family $(\tpi)_{(Q,\uU) \in \ORB}$ defines a natural transformation relating the endofunctors $\tT$ and $\id_{\ORB}$. 
\end{rem}

\subsection{Orbisections}

We now study sections of an orbifold into its tangent orbibundle. These maps will be called \tl orbisections\tr and may be thought of as an analogue of vector fields on manifolds. In this section, $(Q,\uU)$ is an orbifold.

\begin{defn} \label{defn: orbsec} \no{defn: orbsec}
 A map of orbifolds $[\hat{\sigma}] \in \ORBM[(Q,\uU) , \tT (Q,\uU)]$ is called an \ind{}{orbisection} if it satisfies 
	\begin{displaymath}
	 \tpi \circ [\hat{\sigma}] = \ido{}.
	\end{displaymath}
 Its \ind{orbisection!support}{support} $\supp [\hat{\sigma}]$ is the closure of $\setm{x \in Q}{\sigma (x) \neq 0_x}$, where $0_x \in \tT_x Q$ is the zero-element. We define the \ind{space!of orbisections}{\rm set of all orbisections} \gls{Os} of the orbifold $(Q,\uU)$.\\
 An orbisection $[\hat{\sigma}] \in \Os{Q}$ with $\supp [\hat{\sigma}] \subseteq K$ for some compact subset $K \subseteq Q$ is called \ind{orbisection!compactly supported}{compactly supported (in $K$)}.\\
 For $K\subseteq Q$ compact define the set $\OsK{Q} \coloneq \setm{[\hat{\sigma}] \in \Os{Q}}{\supp [\hat{\sigma}] \subseteq K}$\glsadd{OsK} of \ind{space!of orbisections supported  in $K$}{\rm orbisections supported in $K$}. Let \gls{Osc} be the \ind{space!of compactly supported orbisections}{\rm set of all compactly supported orbisections in $\Os{Q}$}.
\end{defn}

If $M$ is a trivial orbifold (i.e.\ a manifold), then orbisections are vector fields on the manifold. It is well known that vector fields for a manifold form a vector space. In Section \ref{sect: sp:os} we will prove that the set $\Os{Q}$ (and the subspaces $\Osc{Q}, \OsK{Q}$ are topological vector spaces over $\RR$ for any orbifold. This fact is quite surprising for a non-trivial orbifold. Indeed, recall that at a singular point, the orbifold tangent space does not support a vector space structure. However, lifts of a special kind for orbisections, we may obtain a vector space structure: For vector fields, it is often advantageous to consider the representative of a vector field $X \colon M \rightarrow TM$ in charts. For a manifold chart $\Psi$, this representative is defined to be $X_\Psi \coloneq d\Psi \circ X \circ\Psi^{-1}$. It is possible to obtain lifts of a similar kind for orbisections on arbitrary orbifolds. 

\begin{defn}
 Consider $[\hat{\sigma}] \in \Os{Q}$ with $\hat{\sigma} = (\sigma, \set{\sigma_i}_{i \in I} , [P_\sigma , \nu_\sigma]) \in \Orb{\vV, \tT \vV}$. If for each $i \in I$, the lift is a vector field $\sigma_i \in \vect{V_i}$, then $(\sigma_i)_{i \in I}$ is called < \ind{orbisection!canonical lifts}{family of canonical lifts for the orbisection $[\hat{\sigma}]$ with respect to $\vV$}. If there is no risk of confusing which orbifold atlas is meant, we will also say that $\set{\sigma_i}_{i\in I}$ is a \ind{canonical family|see{orbisection, can. lifts}}{canonical family} for $[\hat{\sigma}]$.
\end{defn}

Representatives of orbisections with canonical lifts with respect to a given atlas are unique: 

\begin{lem} \label{lem: os:cucr} \no{lem: os:cucr}
 Let $[\hat{f}] \in \Os{Q}$ and $\vV \in \uU$ be an arbitrary orbifold atlas such that there exists a representative $\hat{h} = (f, \set{f_i}_{i \in I}, [P_h,\nu_h]) \in \Orb{\vV,\tT \vV}$ whose lifts form a canonical family for $[\hat{f}]$. Then $\hat{h}$ is unique, i.e.\ if there is another representative of $[\hat{f}]$ whose lifts form a canonical family with respect to  $\vV$, then the members of this family must coincide with $\set{f_i}_{i \in I}$.
\end{lem}

\begin{proof}
 Let $\hat{g} = (f, \set{g_i}_{i \in I}, [P_g , \nu_g]) \in \Orb{\vV, \tT \vV}$ be another representative of $[\hat{f}]$ whose lifts form a canonical family with respect to $\vV$. For each chart $(V_i,G_i,\psi_i) , \ i \in I$ we have $\pi_{TV_i} f_i = \id_{V_i} = \pi_{TV_i} g_i$. On the other hand, $g_i$ and $f_i$ are lifts of $f$, thus for every point $x \in V_i$, there is $\gamma_x \in G_i$ such that $T\gamma_x f_i(x) = \gamma_x . f_i (x) = g_i(x)$. Combining these observations, we obtain 
	\begin{equation}
	 x= \pi_{TV_i} f_i (x) = \pi_{TV_i} g_i (x) = \pi_{TV_i} T\gamma_x f_i(x) = \gamma_x . x.
	\end{equation}
 Thus for each $x \in V_i\setminus \Sigma_{G_i}$ (i.e.\ $x$ is non-singular), we derive $\gamma_x = \id_{V_i}$ and $f_i(x) = g_i (x)$. The continuous maps $f_i$ and $g_i$ coincide on the dense set $V_i \setminus \Sigma_{G_i}$, whence $f_i = g_i$.
\end{proof}

It turns out that analogous to vector fields on manifolds, one is able to construct a canonical family for each orbisection with respect to any given orbifold atlas. At first we have to assure that there is at least some representative with a family of canonical lifts for a given orbisection:

\begin{lem}\label{lem: os:msr} \no{lem: os:msr}
 For every orbisection $[\hat{f}] \in \Os{Q}$, there is a representative $\vV$ of $\uU$ indexed by some $I$ and a representative of an orbifold map $\hat{g} = (f, \set{f_i}_{i \in I}, [P_{\hat{g}} , \nu_{\hat{g}}]) \in \Orb{\vV , \tT \vV}$ such that 
	\begin{compactenum}
	 \item $\hat{g} \in [\hat{f}]$,
	 \item $\set{f_i}_{i \in I}$ is a canonical family for $[\hat{f}]$ with respect to $\vV$. 
	\end{compactenum}
\end{lem}

\begin{proof}
 Following Lemma \ref{lem: oto:nmap}, we choose orbifold atlases $\aA \in \uU$ and $\wW \in \uU$ indexed by $I$ such that there is a representative $\hat{h} = ( f, \set{h_i}_{i \in I} , [P_{\hat{h}} , \nu_{\hat{h}}]) \in \Orb{\wW , \tT \aA}$ of $[\hat{f}]$. For $i \in I$ let $h_i \colon V_i \rightarrow TU_{\alpha (i)}$ be the lift with respect to $(V_i, G_i, \psi_i) \in \wW$ and $(TU_{\alpha (i)}, G_{\alpha (i)}, \pi_{\alpha (i)}) \in \tT \aA $. By Lemma \ref{lem: oto:nmap}, the composition $h_i^1 \coloneq \pi_{TU_{\alpha (i)}} \circ h_i \colon V_i \rightarrow U_{\alpha (i)}$ is a local lift of $\id_Q$, since $\tpi \circ [\hat{h}] = \ido{}$. For each $v \in V_i$ there is an open $G_i$-stable set $V_i^v$ by Proposition \ref{prop: ll:id} such that $h_i^1|_{V_i^v}$ is an open embedding of orbifold charts. \\
 Thus $V_i$ can be covered by open $G_i$-stable subsets $\setm{V_i^j}{j \in J_i}$ such that $h_i^1|_{V_i^j}$ is an embedding of the orbifold chart $(V_i^j, G_{V_i^j}, \psi_i| _{V_i^j})$ into $W_{\alpha (i)}$. Define an orbifold atlas $\vV \in \uU$ via $\vV \coloneq \setm{(V_i^j, G_{V_i^j}, \psi_i|_{V_i^j})}{i \in I , j \in J_i}$. Since $h_i^1$ is invertible on each $V_i^j$, for $j\in J_i$, one can construct a family of lifts for $f$ as follows: Set 
	\begin{displaymath}
	 f_i^j \coloneq T(h_i^1|_{V_i^j})^{-1} \circ h_i|_{V_i^j} \colon V_i^j \rightarrow TV_i^j.
	\end{displaymath}
 A computation proves the identity $\pi_{TV_i^j} \circ f_i^j = \id_{V_i^j}$, i.e. $f_i^j \in \vect{V_i^j}$. Since $h_i^1|_{V_i^j}$ is an embedding of orbifold charts, the same holds for $T(h_i^1|_{V_i^j})^{-1} = (Th_i^1|_{TV_i^j})^{-1}$ (cf. Lemma \ref{lem: pr:tofd}). By construction the mooth maps $f_i^j$ are induced by the lifts $h_i$ of $\hat{h}$ with respect to the inclusion of $V_i^j$ and the open embedding $Th_i^1|_{TV_i^j}$. Hence Lemma \ref{lemdef: ind:ofdm} implies that there is a representative $\hat{g} \in [\hat{f}]$ whose local lifts are given by the family $(f_i^j)_{i \in  I, j\in J_i}$. Therefore, $\hat{g} \in \Orb{\vV,\tT \vV}$ is a representative of $[\hat{f}]$ whose lifts form a canonical family with respect to the atlas $\vV$.
\end{proof}

We now have canonical lifts for an orbisection at our disposal. With this tool, it is possible to deduce a surprising properties of orbisections:  

\begin{prop}\label{prop: os:plg} \no{prop: os:plg}
 Orbisections preserve local groups.
\end{prop}

\begin{proof}
 Consider $[\hat{f}] \in \Os{Q}$ together with a representative $\hat{f} = (f, \set{f_i}_{i \in I} , [P_f,\nu_f])$ such that $\set{f_i}_{i \in I}$ is a canonical family with respect to some orbifold atlas $\vV$. Consider $x \in Q$ together with an orbifold chart $(V_i,G_i,\psi_i)$ such that $x \in \psi_i (V_i)$. Abbreviate $G \coloneq G_i$. Recall $f_i \in \vect{V_i}$, i.e.\ it is a vector field on $V_i$. Choose $z \in V_i$ with $\psi_i (z) = x$. We have to prove that $G_z$ coincides with $G_{f_i(z)}$. To this end consider $\gamma \in G_{f_i(z)}$. By definition, $\gamma$ acts on $TV$; via the \ind{group action!derived}{\rm derived action} $\gamma . v \coloneq T\gamma (v)$. One computes  
	\begin{displaymath}
	 z = \pi_{TV_i} f_i (z) = \pi_{TV_i} (\gamma . f_i (z)) = \pi_{TV_i} T\gamma (f_i (z)) = \gamma . \pi_{TV_i} f_i (z) = \gamma . z.
	\end{displaymath}
 Thus every $\gamma \in G_{f_i (z)}$ is an element of $G_z$. Hence $\theta \colon G_{f_i(z)} \rightarrow G_z, \gamma \mapsto \gamma$ is an injective group homomorphism. 
 We claim that $\theta$ is surjective. To prove this, consider $\delta \in G_z$. Observe that every $\delta \in G_z$ is a change of charts (even an embedding of orbifold charts) and there is $g \in P_f$ together with an open (connected) neighborhood $\Omega_z \subseteq V_i$ of $z$ such that $\delta|_{\Omega_z} = g|_{\Omega_z}$ holds. The map $\nu_f (g)$ is a change of charts of $TV_i$ into itself. Restricting to the open connected component $C$ of $\dom \nu_f (g)$ which contains $f_i(z)$, \cite[Lemma 2.11]{follie2003} implies that there is a unique $\gamma \in G$ such that $\nu_f (g)|_{C} = \gamma|_{C}$. On the open set $\Omega_z \cap f_i^{-1} (C)$, the identity 
	\begin{equation}\label{eq: os:gsc}
	 f_i \circ \delta|_{\Omega_z \cap f_i^{-1}(C)} = \nu_f (g) f_i|_{\Omega_z \cap f_i^{-1}(C)} = \gamma . f_i|_{\Omega_z \cap f_i^{-1}(C)}
	\end{equation}
 holds. The set $\Omega_z \cap f_i^{-1}(C)$ is a non-empty open set and by Newman's Theorem \ref{thm: newman} there is a non-singular $y \in \Omega_z \cap f_i^{-1}(C)$. Specializing to $y$, equation \eqref{eq: os:gsc} yields: 
	\begin{displaymath}
	 	f_i(\delta . y) = \gamma . f_i(y) = T\gamma f_i (y) \quad \Rightarrow \	\delta .y = \pi_{TV_i} f_i (\delta .y) = \pi_{TV_i} T\gamma f_i(y) = \gamma . y.
	\end{displaymath}
 Then $\delta^{-1}\gamma . y = y$ and $y$ being non-singular forces $\gamma = \delta$. Applying this to \eqref{eq: os:gsc} we obtain: 
	\begin{displaymath}
	 f_i (z) = f_i (\delta . z) = T\delta . f_i (z) = \delta . f_i(z).
	\end{displaymath}
 In other words, $\delta$ fixes $f_i(z)$ and thus $\delta$ is an element of the isotropy subgroup $G_{f_i(z)}$. Thus $\theta$ is surjective. We conclude that $\theta \colon G_z \rightarrow G_{f_i(z)} , \gamma \mapsto \gamma$ is an isomorphism of groups and that the local groups $\Gamma_z$ and $\Gamma_{f(z)}$ are isomorphic.
\end{proof}

The property to preserve local groups limits the choice of images an orbisection may take on a given singular point. In particular, there are elements in the tangent space at a singular point which are not in the image of any orbisection. We refer to Example \ref{ex: ad:nofd} for such a case.  

\begin{prop}\label{prop: os:chch} \no{prop: os:chch}
 Let $[\hat{f}]$ be an orbisection and $\vV \in \uU$ be an orbifold atlas. Furthermore, let $\hat{f} = (f, \set{f_i}_{i \in I}, [P_f,\nu_f]) \in \Orb{\vV, \tT \vV}$ be a representative of $[\hat{f}]$ such that $\set{f_i}_{i \in I}$ is a family of canonical lifts. For each element $\phi$ of the set of changes of charts $\Ch_\vV$ of $\vV$  (cf. Notation \ref{nota: chch:abr}) with $\dom \phi \subseteq V_i$ and $\cod \phi \subseteq V_j$, $(V_\alpha , G_\alpha , \psi_\alpha) \in \vV, \ \alpha \in \set{i,j}$, the identity 
	\begin{equation}\label{eq: os:chch}
	 f_j \phi = T\phi f_i|_{\dom \phi} 
	\end{equation}
 holds. The pair $(\Ch_\vV , \nu)$ with
	\begin{displaymath}
	 \nu \colon \Ch_\vV \rightarrow \PSI (\tT\vV), \phi \mapsto T\phi,
	\end{displaymath}
 is a representative of $[P_f,\nu_f]$. Here $\Ch_\vV$ is the quasi-pseudogroup of all changes of charts for the atlas $\vV$ (cf.\ Notation \ref{nota: chch:abr}).
\end{prop}

\begin{proof}
 Pick an arbitrary change of charts $\phi$ as above and choose a representative $(P_f ,\nu_f)$ of $[P_f,\nu_f]$. It suffices to prove the identity \eqref{eq: os:chch} on small neighborhoods of arbitrary points in $\dom \phi$. Let $x_0 \in \dom \phi$ be such a point. Since $P_f$ generates $\PSI (\vV)$, there is an open $x_0$-neighborhood $U_{x_0} \subseteq \dom \phi \subseteq V_i$ together with $\gamma_{x_0}^\phi \in P_f$ such that $\gamma_{x_0}|_{U_{x_0}} = \phi|_{U_{x_0}}$ holds. By definition, we obtain a local lift of $f$: 
	\begin{equation}\label{eq: os:chch2}
	  f_j \phi|_{U_{x_0}} = f_j \gamma_{x_0}|_{U_{x_0}}^\phi = \nu_f( \gamma_{x_0}^\phi) f_i|_{U_{x_0}}.
	\end{equation}
 On the other hand, the composition $T\phi f_i|_{U_{x_0}}$ is defined, since $f_i|_{U_{x_0}} \in \vect{U_{x_0}}$. By Lemma \ref{lem: pr:tofd} (a), $T\phi$ is a change of charts of $\tT \vV$ and thus $T\phi f_i|_{U_{x_0}}$ is a local lift of $f$ on $U_{x_0}$. For every $y \in U_{x_0}$, we obtain 
	\begin{displaymath}
	 T\psi_j \nu_f (\gamma_{x_0}^\phi ) f_i (y) = T\psi_j T\phi f_i(y).
	\end{displaymath}
 Thus there is a unique group element $g_y \in G_j$ such that $g_y. \nu_f (\gamma_{x_0}^\phi ) f_i (y) = T\phi f_i(y)$ holds. In Proposition \ref{prop: os:plg} we have seen that orbisections preserve local groups, whence they preserve non-singular points. Therefore lifts of orbisections map non-singular points to non-singular points. The set $U_{x_0}$ is a non-empty open subset of $V_i$ and by Newman's Theorem \ref{thm: newman}, the non-singular points of the $G_i$-action on $V_i$ are dense in $U_{x_0}$. Using \eqref{eq: os:chch2} for non-singular $y \in U_{x_0}$ we obtain the identities 
	\begin{align*}
	  &T\phi f_i(y) = g_y . \nu_f (\gamma_{x_0}^\phi) f_i (y) = g_y . f_j \phi (y) = Tg_y (f_j \phi (y)), \text{ whence} \\
	  & \phi (y) = \pi_{TV_j}  T\phi f_i(y) = \pi_{TV_j}  Tg_y (f_j \phi (y)) = g_y . \phi (y) .
	\end{align*}
 As changes of charts preserve non-singular points and $y$ is non-singular, $g_y = \id_{V_j}$ follows. The maps $\nu_f (\gamma_{x_0}^\phi ) f_i$ and $T\phi f_i$ therefore coincide on the non-singular points of $U_{x_0}$. As these points form a dense subset in $U_{x_0}$, the continuous maps must coincide on $U_{x_0}$, whence $T\phi f_i|_{U_{x_0}} = \nu_f (\gamma_{x_0}^\phi) f_i|_{U_{x_0}}$ holds and indeed $T\phi f_i|_{U_{x_0}} = f_j \phi|_{U_{x_0}}$ follows.\\
 The quasi-pseudogroup $\Ch_\vV$ generates $\PSI (\vV)$ and our previous considerations show that $\nu$ (as defined above) satisfies property (R4a) of Definition \ref{defn: rep:ofdm}. The functoriality of $T$ implies properties (R4b)-(R4d) of Definition \ref{defn: rep:ofdm} for $(\Ch_\vV , \nu)$. Notice that we did not change the family of lifts $\set{f_i}_{i \in I}$. Thus $\hat{h} \coloneq (f, \set{f_i}_{i\in I}, [\Ch,\nu]) \in \Orb{\vV ,\tT \vV}$ is a charted map such that $[\hat{f}] = [\hat{h}]$.
\end{proof}

\begin{rem}
 Let $M,N$ be smooth manifolds and $f \colon M \rightarrow N$ be a smooth map. Recall that $\sigma \in \vect{M}$ and $\tau \in \vect{N}$ are called \ind{vector field!$f$-related}{$f$-related} if $Tf \circ \sigma = \tau \circ f$ holds. Hence Proposition \ref{prop: os:chch} shows that canonical families of an orbisection are families of pairs of $f$-related vector fields, where $f$ runs through the changes of charts of the domains of the pair. 
\end{rem}

\begin{lem}\label{lem: os:cfas} \no{lem: os:cfas}
 Let $[\hat{f}]$ be an orbisection and $\vV$ be an arbitrary representative of $\uU$. There is a refinement $\vV'$ of $\vV$ and a representative $\hat{h} = (f, \set{h_i}_{i \in I}, [P,\nu]) \in \Orb{\vV', \tT \vV'}$ of $[\hat{f}]$ such that $\set{h_i}_{i\in I}$ is a family of canonical lifts for $[\hat{f}]$. 
\end{lem}

\begin{proof}
 By Lemma \ref{lem: os:msr}, we may choose a representative $\wW$ of $\uU$ indexed by $I$ and a representative $\hat{g} = (f, \set{g_i}_{i \in I}, [P,\nu]) \in \Orb{\wW , \tT \wW}$ of $[\hat{f}]$ such that $\set{g_i}_{i \in I}$ is a canonical family. Choose a common refinement $\vV'$ of $\wW$ and $\vV$. The refinement $\vV'$ induces a common refinement $\tT \vV'$ of $\tT \vV$ and $\tT \wW$, since embeddings of orbifold charts are mapped to embeddings of orbifold charts by the tangential functor $T$. Let $\vV'$ be indexed by $J$ and $\alpha \colon J \rightarrow I$ be a map such that for $j \in J$ there is an embedding of orbifold charts $\lambda_j \colon (V_j',G_j',\pi_j') \rightarrow (W_{\alpha (j)}, H_{\alpha (j)} ,\psi_{\alpha (j)})$. The family $\set{g_i}_{i \in I}$ is a canonical family, therefore 
	\begin{displaymath}
	 g_{\alpha(j)} \lambda_j (V_j') = g_{\alpha (j)} (\im \lambda_j) \subseteq T\im \lambda_j.
	\end{displaymath}
 Define the maps $h_j \coloneq (T\lambda_j)^{-1} g_{\alpha (j)} \lambda_j \colon V_j' \rightarrow TV_j'$. Then Lemma \ref{lemdef: ind:ofdm} assures that there is a pair $(P ,\nu )$ such that $\hat{h} \coloneq (f, \set{h_j}_{j\in J}, [P,\nu])$ is a representative of $[\hat{f}]$. A computation yields 
	\begin{displaymath}
	 \pi_{TV_j} h_j = \pi_{TV_j} (T\lambda_j)^{-1} g_{\alpha (j)} \lambda_j = \lambda_j^{-1} \pi_{TW_{\alpha (j)}} g_{\alpha (j)} \lambda_j= \id_{V_j}
	\end{displaymath}
 for each $j \in J$. In conclusion, $\set{h_j}_{j\in J}$ is a canonical family and the domain atlas of $\hat{h}$ is a refinement of $\vV$.
\end{proof}

The results obtained so far show that each orbisection possesses representatives whose lifts form canonical families for suitable refinements of $\vV$. We will now prove a converse: For each orbisection and an arbitrary orbifold atlas, there is a representative whose lifts form a canonical family with respect to the given atlas. This result is quite surprising since in general maps of orbifolds need not have lifts on an orbifold chart chosen in advance.

\begin{prop}\label{prop: os:cl} \no{prop: os:cl}
 Let $[\hat{f}] \in \Os{Q}$ and $\wW$ be an arbitrary representative of $\uU$ indexed by $J$. There exists a representative $\hat{g}=(f, \set{g_j}_{j\in J}, [P,\nu]) \in \Orb{\wW, \tT \wW}$ such that $\set{g_j}_{j\in J}$ is a canonical family with respect to $\wW$ . 
\end{prop}

\begin{proof}
Lemma \ref{lem: os:cfas} allows us to choose a refinement $\vV$ of $\wW$ indexed by $I$ and a representative $\hat{h} \coloneq (f, \set{f_{V_i}}_{i \in I}, [P,\nu]) \in \Orb{\vV, \tT \vV}$ of $[\hat{f}]$ such that $\set{f_{V_i}}_{i\in I}$ is a family of canonical lifts for $[\hat{f}]$. Let $(W_j , G_j , \psi_j) \in \wW$ be an arbitrary orbifold chart. We have to construct a local lift of $f$ on $(W_j , G_j , \psi_j)$. To achieve this, consider $x \in W_j$. Since $\psi_j (x) \in Q$ and $\vV$ is an atlas, there is a chart $(V_i, G_i,\varphi_i) \in \vV$ together with a change of charts $\lambda_{x} \in \CH{V_i}{W_j}$ (cf.\ Notation \ref{nota: chch:abr}) such that $x \in \im \lambda_x$. Then we define 
    \begin{equation}\label{eq: assign}
     f_{W_j} (z) \coloneq T \lambda_x f_{V_i} \lambda_x^{-1} (z) \in T_z W_j
    \end{equation}
 for all $z \in \im \lambda$. The definition of $f_{W_j}$ neither depends on the choice of $\lambda_x$ nor on $(V_i,G_i,\varphi_i)$. To see this, consider another chart $(V_j,G_j,\varphi_j) \in \vV$ and a change of charts morphism $\mu_x \in \CH{V_r}{W_j}$ with $x \in \im \mu_x$. Denote the intersection $\im \lambda_x \cap \im \mu_x$ as $\Omega_x$. We will show that for each $z$ in the open $x$-neighborhood $\Omega_x$, equation \eqref{eq: assign} yields the same $f_{W_j} (z)$ if $\mu_x$ is used instead of $\lambda_x$. Observe that $h_x \coloneq \lambda_x^{-1} \mu_x|_{\mu_x^{-1} (\Omega_x)}$ is a change of charts in $\CH{V_r}{V_i}$. Using that the family $\set{f_{V_i}}_{i \in I}$ is a canonical family of lifts with respect to $\vV$, we compute for $z \in \Omega_x$ 
  \begin{displaymath}
   T\lambda _x f_{V_i} \lambda_x^{-1} (z) = T\lambda_x f_{V_i} h_x \mu_x^{-1} (z) = T\lambda_x Th_x f_{V_j} \mu_x^{-1} (z) = T \mu_x f_{V_j} \mu_x^{-1} (z). 
  \end{displaymath}
 Hence, on $\Omega_x$ the assignment \eqref{eq: assign} does not depend on any of the above choices. Thus it makes sense to define a map as follows 
  \begin{displaymath}
   f_{W_j} \colon W_j \rightarrow TW_j , x \mapsto T\lambda f_{V_i} \lambda^{-1} (x) \text{ if there is } (V_i,G_i,\varphi_i) \in \vV \text{ and } \lambda \in \CH{V_i}{W_j} \text{ with } x \in \im \lambda.
  \end{displaymath}
 For each $x \in W_j$ there is a change of charts such that the identity \eqref{eq: assign} holds in an open $x$-neighborhood. Hence, the map $f_{W_j}$ is a smooth and by construction, a smooth vector field. Repeating the construction for each chart in $\wW$, we obtain a family of vector fields $\set{f_{W_j}}_{j \in J}$ which lift $f$.\\
 We claim that the family of vector fields is a canonical family of lifts. It suffices to prove identity \eqref{eq: os:chch} for each $\phi \in \CH{W_j}{W_k}$ and $j,k \in J$. To this end fix $\phi \in \CH{W_j}{W_k}$ and consider $z \in \dom \phi$ together with a change of charts $\lambda_z \in \CH{V_i}{W_j}$ such that $z \in \im \lambda_z \subseteq \dom \phi$. Then $\phi \circ \lambda_z \in \CH{V_i}{W_k}$ implies 
   \begin{displaymath}
    f_{W_k} \circ \phi (z) \stackrel{\eqref{eq: assign}}{=} T(\phi \lambda_z) f_{V_i} (\phi \circ \lambda)^{-1} \phi (z) = T\phi T\lambda_z f_i \lambda_z^{-1} (z) \stackrel{\eqref{eq: assign}}{=}  T\phi f_{W_j} (z).
   \end{displaymath}
 Since $z \in \dom \phi$ was arbitrary, this proves identity \eqref{eq: os:chch}. Hence by Proposition \ref{prop: os:chch} we may choose $\nu$, such that the map $\hat{g} \coloneq (f,\set{f_{W_j}}_{j \in J}, [\Ch_\wW , \nu])$ is a representative of an orbisection with canonical lifts. 
 The atlas $\vV$ is a refinement of $\wW$, thus for every $i \in I$, there is an embedding of orbifold charts $\lambda_i \colon (V_i,G_i,\pi_i) \rightarrow (W_{\alpha (i)}, G_{\alpha (i)} , \psi_{\alpha (i)})$. By construction, we obtain $f_{V_i} = T\lambda_i^{-1} f_{W_{\alpha (i)}}\lambda_i$ and therefore every lift $f_{V_i}$ is induced by a suitable lift of $\hat{g}$. Following Definition \ref{defn: sim:ofd}, we have $\hat{g} \sim \hat{h}$ and the classes $[\hat{g}]$ and $[\hat{f}]$ coincide. Thus the lifts are a canonical family of $[\hat{f}]$ with respect to $\wW$.
\end{proof}

Proposition \ref{prop: os:cl} shows that every orbisection may be identified in every given atlas with a unique family of canonical representatives. In particular, orbisections satisfy analogous properties as $C^\infty$-sections in the tangent bundle in the sense of \cite[below Remark 4.1.8]{chen2001}.

\begin{rem}\label{rem: os:clpro} \no{rem: os:clpro}
\begin{compactenum}
 \item A family $\fF$ of vector fields on an orbifold atlas $\vV$ which satisfies Equation \eqref{eq: os:chch} induces a continuous map $F \colon Q \rightarrow \tT Q$ (cf.\ the proof of Proposition \ref{prop: os:vsp} for the explicit construction) such that  
	\begin{compactitem}
	 \item[-] $(F, \fF, [\Ch_\vV, \nu]) \in \Orb{\vV , \tT\vV}$ with $\nu \colon \Ch_\vV \rightarrow \PSI (\tT \vV), \lambda  \mapsto T\lambda$, 
	 \item[-] $\fF$ is a canonical family.
	\end{compactitem}
 Vice versa, if $(f, \set{f_i}_{i \in I} , [P_f,\nu_f])$ is a representative of an orbisection whose lifts form a canonical family with respect to an atlas $\vV$, then the above construction for $\set{f_i}_{i \in I}$ yields the map $f$. Lemma \ref{lem: os:cucr} implies that an orbisection is uniquely determined by its family of canonical lifts with respect to any atlas $\vV$. This induces a one to one correspondence between the set of orbisections and families of vector fields for some orbifold atlas $\vV$ which satisfy \eqref{eq: os:chch}. 
 \item Notice that (a) implies: For $[\hat{f}] \in \Os{Q}$ and $(U, G, \psi) \in \uU$, there is a unique vector field $\hat{f}_U \in \vect{U}$ such that for $\hat{f} = (f, \set{f_i}_{i \in I}, [P,\nu])$ the identity $T\psi \hat{f}_U = f\psi$ holds.
 \item The canonical lift of the \ind{orbisection!zero orbisection}{zero orbisection} $\zs$\glsadd{zs} with respect to some orbifold chart $(U,G,\psi)$ is the zero-section in $\vect{U}$. 
 If $[\hat{f}] \in \Os{Q}$ is an orbisection and $(U,G,\psi)\in \uU$ is some chart such that $\psi (U) \cap \supp [\hat{f}] = \emptyset$, then the canonical lift of $[\hat{f}]$ on $U$ is the zero-section in $\vect{U}$.
\item Proposition \ref{prop: os:plg} implies that orbisections in $\Os{Q}$ take their values in 
  \begin{displaymath}
   \tT Q^{\text{inv}} \coloneq \setm{[\pi,v]}{(U,G,\psi) \in \uU, v \in TU \text{ with } g.v = v \text{ for all } g \in G_{\pi_{TU} (v)} }.
  \end{displaymath}
Notice that $T_x U^{\text{inv}} \coloneq \setm{v \in T_x U}{g.v=v \text{ for all } g \in G_x}$ is a subvectorspace of $T_x U$ and $TU^{\text{inv}} \coloneq \bigcup_{x \in U} T_xU^{\text{inv}}$ is invariant with respect to the derived $G$-action on $TU$. Since the chart mapping $T\psi$ is an open map, \cite[VI.\ Theorem 2.1]{dugun1966} implies that the restriction $T\pi|_{TU^{\text{inv}}}^{\tT Q^{\text{inv}} \cap \im T\pi}$ is a quotient map. Furthermore, the map $T\pi|_{T_xU^{\text{inv}}}^{\tT Q^{\text{inv}} \cap \tT_{\pi(x)} Q}$ is bijective. Thus $\tT_{\pi (x)}Q^{\text{inv}} \coloneq T\pi(T_xU^{\text{inv}})$ is in a natural way a vector space, whence the fibres of $\tT Q^{\text{inv}}$ are vector spaces. 
Notice that this vector space structure induces a vector space structure on $\Os{Q}$ by pointwise operations on canonical lifts. The details are recorded in the next section.
\item The underlying continuous map $\sigma$ of an orbisection $[\hat{\sigma}] \in \Os{Q}$ uniquely determines the orbisection. To see this, we choose a family of canonical lifts $(\sigma_i)_{i \in I}$ with respect to some atlas $\set{(U_i,G_i,\psi_i)}_{i \in I} \in \uU$ for $[\hat{\sigma}]$. From part (d) we derive for $x \in U_i $ the identity 
  \begin{displaymath}
   \sigma_i (x) = (T\psi_i|_{T_xU^{\text{inv}}}^{\tT Q^{\text{inv}} \cap \tT_{\pi(x)} Q})^{-1} \circ \sigma \circ \psi_i (x). 
  \end{displaymath}
 Hence, the underlying map $\sigma$ uniquely determines the canonical lifts $\sigma_i$. By part (a), the canonical family $\set{\sigma_i}_{i \in I}$ uniquely determines $[\hat{\sigma}]$, whence the assertion follows. In particular, we obtain a canonical embedding 
  \begin{displaymath}
   \Os{Q} \rightarrow C ( Q, \tT Q) , [\hat{\sigma}] \mapsto \sigma
  \end{displaymath}
\end{compactenum}
\end{rem}
\newpage
\subsection{Spaces of orbisections}\label{sect: sp:os}

 We now study spaces of orbisections. For these spaces we will obtain the structure of a real topological vector space. The construction of the vector space structure is inspired by arguments first given in \cite{bb2008}.

\begin{prop}\label{prop: os:vsp} \no{prop: os:vsp}
 The set $\Os{Q}$ of orbisections is a real vector space with pointwise vector space operations on canonical lifts. The zero element $\zs \in \Os{Q}$ of $\Os{Q}$ is called the \ind{orbifold map!zero orbisection}{zero orbisection}. Endowing $\Os{Q}$ with this vector space structure, the sets $\OsK{Q} \subseteq \Osc{Q} \subseteq \Os{Q}$ become linear subspaces.
\end{prop}

\begin{proof}
 Let $[\hat{f}], [\hat{g}] \in \Os{Q}$ and choose an arbitrary representative $\vV$ of the orbifold structure $\uU$, indexed by some set $I$. By Proposition \ref{prop: os:cl} we may choose unque representatives of orbifold maps $\hat{f} = (f, \set{f_i}_{i\in I} , [P_f , \nu_f]) \in \Orb{\vV, \tT \vV}$ of $[\hat{f}]$ and $\hat{g} = (g, \set{g_i}_{i\in I}, [P_g, \nu_g]) \in \Orb{\vV, \tT \vV}$ of $[\hat{g}]$ such that the families of lifts are canonical families. Without loss of generality $P_f = P_g = \Ch_\vV$ and $\nu_f (\lambda ) = \nu_g (\lambda ) = T\lambda$ hold, by Proposition \ref{prop: os:chch}. By construction, for each $i \in I$ the lifts are vector fields $f_i,g_i \in \vect{V_i}$.  Recall from \cite[2.7]{cmfd01} that the vector space structure on $\vect{V_i}$ is induced by pointwise operations. We define the vector space operations on $\Os{Q}$ via the following construction:\\
 For $z \in \RR$ consider $f_i + zg_i\colon V_i \rightarrow TV_i \in \vect{V_i}$. Remember that tangent maps act as linear maps on each tangent space. For every change of charts $\lambda \in \PSI (\vV)$ with $\dom \lambda \subseteq V_i$ and $\cod \lambda \subseteq V_j$ we obtain: 
	\begin{align}\label{eq: os:comp}
	 (f_j + z g_j ) \lambda (p) 	&= f_j (\lambda (p)) + z g_j(\lambda (p)) = \nu_f (\lambda) f_i (p) + z \nu_g (\lambda) g_i (p) \notag\\
				 	&= T_p\lambda (f_i(p)) + z T_p \lambda (z g_i(p)) = T_p \lambda (f_i (p) +zg_i(p)) \\
					&\equalscolon \nu_{f+zg} (\lambda) (f_i (p) + zg_i (p)). \notag
	\end{align}
 Define the quasi-pseudogroup $P_{f+zg} \coloneq \Ch_\vV$ together with $\nu_{f+zg} \colon P_{f+zg} \rightarrow \PSI (\tT \vV) ,\ \lambda \mapsto T\lambda$.
 The pair $( P_{f+zg} , \nu_{f+zg})$ and the family $(f_i + z g_i)_{i \in I}$ satisfy properties (R4a)-(R4d) of Definition \ref{defn: rep:ofdm}. 
 Notice that by Identity \eqref{eq: os:comp} for a chart $(V_i,G_i,\psi_i) \in \vV$ the map $T\psi_i (f_i + zg_i)$ is constant on each fibre $\psi_i^{-1} (y)$. As $\psi_i$ is a quotient map, the map  
	\begin{displaymath}
	 f + zg|_{\psi_i (V_i)} \colon \psi (V_i) \rightarrow T\psi_i (TV_i) , x \mapsto T\psi_i \circ (f_i + z g_i)\psi^{-1}(x) 
	\end{displaymath}
 is continuous, by \cite[VI. Theorem 3.2]{dugun1966} 
 Furthermore, the map $f_i +z g_i$ is a smooth lift for $f + zg|_{\psi_i (V_i)}$. We claim that for every pair $(i,j) \in I \times I$, the maps $f + zg|_{\psi_i (V_i)}$ and $f + zg|_{\psi_j (V_j)}$ coincide on $\psi_j (V_j) \cap \psi_i (V_i)$. If this is true, then $f+zg \colon Q \rightarrow \tT Q , x \mapsto f+zg|_{\psi_i(V_i)}(x)$, for $x \in \psi_i (V_i)$ is a well-defined continuous map. We obtain a charted orbifold map 
	\begin{displaymath}
	 \widehat{f + zg} \coloneq (f+zg, \set{f_i +z g_i}_{i \in I}, [P_{f+zg} , \nu_{f+zg}]) \in \Orb{\vV , \tT \vV}
	\end{displaymath}
 such that each lift $f_i + z g_i$ is a vector field. Hence $\set{f_i + zg_i}_{i\in I}$ is a canonical family with respect to the atlas $\vV$ and $[\widehat{f+zg}] \in \Os{Q}$ holds. Proof of the claim: Consider $x \in \psi_i (V_i) \cap \psi_j (V_j)$. For every pair $y_\alpha \in \psi_\alpha^{-1} (x) , \ \alpha \in \set{i,j}$, there is a change of charts $\lambda \in \CH{V_i}{V_j}$ such that $\lambda (y_i) = y_j$. Again by \eqref{eq: os:comp}, the claim follows as 
	\begin{align*}
	 f+zg|_{\psi_j (V_j)} (x) &= T\psi_j (f_j +z g_j ) (y_j) = T\psi_j (f_j + t g_j ) (\lambda (y_i)) \\
					&= T\psi_j T\lambda (f_i + z g_i ) (y_i) = T\psi_i (f_i + z g_i)(p) = f+zg|_{\psi_i (V_i)} (x)
	\end{align*}
 It remains to show that the construction does not depend on the atlas $\vV$. Let $\vV'$ be another representative of $\uU$ and $\hat{f}'$ and $\hat{g}'$, respectively, be representatives of $[\hat{f}]$ resp. $[\hat{g}]$, whose families of lifts form canonical families with respect to $\vV'$. By Lemma \ref{lem: ofa:cr}, we may choose a common refinement of $\vV$ and $\vV'$. The definition of equivalence of orbifold maps implies that the classes will be equal if the induced lifts on this refinement coincide. Without loss of generality we may assume that $\vV'$ refines $\vV$: Let $\vV' = \setm{(W_k, H_k, \phi_k)}{k\in K}$ and $\alpha \colon K \rightarrow J$ be the map which assigns to $k \in K$ an element of $I$ such that there is an embedding of orbifold charts $\lambda_k \colon (W_k,H_k,\phi_k) \rightarrow (V_{\alpha (k)} , G_{\alpha (k)}, \psi_{\alpha (k)})$. The atlas $\tT \vV'$ for $\tT Q$ is a refinement of $\tT \vV$. In particular, $T\lambda_k$ is an embedding of $(TW_k, H_k, T\phi_k)$ into $(TV_{\alpha (k)}, H_{\alpha (k)} , T\psi_{\alpha (k)})$. Let $\hat{f}' = (f, \set{f_k'}_{k \in K}, [P_f' , \nu_f'])$ and $\hat{g}' = (g, \set{g_k'}_{k \in K}, [P_g', \nu_g'])$. 
 The families $\set{f_i}_{i \in I}$ and $\set{g_i}_{i \in I}$ are families of vector fields and we obtain induced vector fields on each chart $(W_k,H_k,\phi_k)$ since this chart embeds into a chart $(V_{\alpha (k)},G_{\alpha (k)},\psi_{\alpha (k)})$. 
 Combine Lemma \ref{lem: os:cfas} and the uniqueness assertion for canonical lifts (Lemma \ref{lem: os:cucr}) to obtain the following identity for the induced vector fields 
	\begin{displaymath} 
	 f_k' = T\lambda_k^{-1}f_{\alpha (k)}\lambda_k , \quad g_k' = T\lambda_k^{-1}g_{\alpha (k)}\lambda_k.
	\end{displaymath}
 Constructing $\widehat{f' + z g'} \in \Orb{\vV', \tT \vV'}$ as above, we deduce from the last identity that $\widehat{f + z g} \sim \widehat{f' +zg'}$. A vector space structure on $\Os{Q}$ is thus defined via the assignment: 
  \begin{displaymath}
   [\hat{f}] + z [\hat{g}] \coloneq [\widehat{f+zg}].
  \end{displaymath}
 Clearly $\zs \in \OsK{Q} \subseteq \Osc{Q}$ holds, whence these subsets are not empty. The last claim follows from the definitions: For $[\hat{f}] , [\hat{g}] \in \Osc{Q}$ with $\supp [\hat{f}] \subseteq K$ and $\supp [\hat{g}] \subseteq L$ with $K, L \subseteq Q$ compact, one obtains $\supp ([\hat{f}] + z[\hat{g}]) \subseteq \supp [\hat{f}] \cup \supp [\hat{g}] \subseteq K \cup L$. Therefore $\OsK{Q}$ and $\Osc{Q}$ are linear subspaces.
\end{proof}

Our goal for the remainder of this section is to topologize the vector spaces $\Os{Q}$ and $\Osc{Q}$. If $Q$ is a compact topological space, then $\Os{Q}$ will be a Fr\'{e}chet space.
 
\begin{lem}\label{lem: os:btop} \no{lem: os:btop}
 Let $(Q, \uU)$ be an orbifold and $\vV = \setm{(U_i,G_i,\psi_i)}{i \in I}$ an arbitrary representative of $\uU$ indexed by $I$. There is a bijection identifying each $[\hat{f}] \in \Os{Q}$ with a unique representative $\hat{f}_\vV$ whose lifts $\set{f_{U_i}}_{i \in I}$ form a canonical family for $[\hat{f}]$ with respect to $\vV$.
 \begin{compactenum}
  \item The map 
	\begin{displaymath}
	 \Lambda_\vV \colon \Os{Q} \rightarrow \prod_{i \in I} \vect{U_i} ,  \hat{f}_\vV \mapsto (f_{U_i} )_{i \in I}
	\end{displaymath}
 is a linear injection into a direct product of topological vector spaces (cf.\ Section \ref{sect: sections} for information on $\vect{U_i}$), whose image is the closed vector subspace \glsadd{LambdavV}
	\begin{displaymath}
	H \coloneq \setm{(f_i)_{i \in I} \in \prod_{i \in I} \vect{U_i}}{\forall \lambda \in \Ch_\vV, \dom \lambda \subseteq U_i, \cod \lambda \subseteq U_j,\ f_j \lambda = T\lambda f_i|_{\dom \lambda}}
	\end{displaymath}
  \item If $\vV$ is a locally finite atlas such that each chart in $\vV$ is relatively compact, then the map 
	\begin{displaymath}
	 \Lambda_\vV \colon \Osc{Q} \rightarrow \bigoplus_{i \in I} \vect{U_i} ,  \hat{f}_\vV \mapsto (f_{U_i} )_{i \in I}
	\end{displaymath}
 is a linear injection into the direct sum of topological vector spaces (cf. \cite[4.3]{jarchow1980}). 
 Making identifications, its image is the closed vector subspace $H \cap \bigoplus_{i \in I} \vect{U_i}$.
 \end{compactenum}
\end{lem}

\begin{proof}
\begin{compactenum}
  \item For $[\hat{f}] \in \Os{Q}$, we let $(\hat{f}_{U_i})_{i \in I}$ be the family of canonical lifts with respect to $\vV$. Proposition \ref{prop: os:chch} shows that $\im \Lambda_\vV$ is contained in $H$. Remark \ref{rem: os:clpro} (a) implies that $\Lambda_\vV$ is injective and $\Im \Lambda_\vV = H$ holds. The vector space operations of $\Os{Q}$ are defined via pointwise operations for families of vector fields. Hence by definition, $\Lambda_\vV$ is a linear map.\\ 
  We have to show that $H$ is a closed vector subspace. Consider $\lambda \in \CH{U_i}{U_j}$ and arbitrary $y \in \dom \lambda$. Each element in $H$ must satisfy $f_j (\lambda (y)) = T\lambda f_i (y)$, i.e. we observe $\text{ev}_{\lambda (y)} (f_j) = (T\lambda \circ \text{ev}_y) (f_i)$. Here $\text{ev}_y$ and $\text{ev}_{\lambda (y)}$ are point evaluation maps defined on $\vect{U_i}$ and $\vect{U_j}$, respectively. The choice of the topology on $\vect{U_i}$ (cf.\ Definition \ref{defn: top:vect}) implies that point evaluation maps are continuous mappings on these spaces. To see this, note that for a manifold chart $(\kappa , V_\kappa)$ the restriction map $\res_{V_\kappa}^{U_i}$ is continuous (cf.\ Notation \ref{nota: res}). By \cite[Proposition 3.20]{alas2012}, point evaluation maps are continuous for all spaces $C^\infty (V_{\kappa},T_{\lambda (y)} V_j)$, whence the claim. Since the projections $\text{pr}_k \colon \prod_{i \in I} \vect{U_i} \rightarrow \vect{U_k}, (f_i)_{i \in I} \mapsto f_k$ are continuous for all $k$, we derive a continuous mapping 
  \begin{displaymath}
   h_{\lambda, y} \colon \prod_{i \in I} \vect{U_i} \rightarrow T_{\lambda (y)} U_j , (f_i)_{i \in I} \mapsto (T\lambda \circ \text{ev}_y) (f_i) - \text{ev}_{\lambda (y)} (f_j).
  \end{displaymath}
We may now write the space $H$ as the intersection
 \begin{displaymath}
  H = \bigcap_{\lambda \in \Ch_\aA} \bigcap_{y \in \dom \lambda} h_{\lambda,y}^{-1} (0).
 \end{displaymath}
Since each $h_{\lambda,y}$ is continuous, the space $H$ is a closed subspace of $\prod_{i \in I} \vect{U_i}$ as an intersection of such spaces. 
 \item The atlas $\vV$ is locally finite and thus only finitely many charts intersect a given compact set. In particular, $\Lambda_\vV$ makes sense. The canonical injection $I \colon \bigoplus_{i \in I} \vect{U_i} \rightarrow \prod_{i\in I} \vect{U_i}$ is continuous by \cite[4.3.1]{jarchow1980} and thus $I^{-1} (H) = H \cap \bigoplus_{i \in I} \vect{U_i}$ is a closed subset of $\bigoplus_{i \in I} \vect{U_i}$. Again by Proposition \ref{prop: os:chch} $\im \Lambda_\vV$ is contained in $I^{-1} (H)$ and by Remark \ref{rem: os:clpro}, $\Lambda_\vV$ is injective and $\Im \Lambda_\vV = I^{-1} (H) = H \cap \bigoplus_{i \in I} \vect{U_i}$. 
\end{compactenum}
\end{proof} 

\begin{defn}\label{defn: os:top} \no{defn: os:top}
 \begin{compactenum}
  \item Let $\vV$ be a representative of $\uU$ for an orbifold $(Q,\uU)$. Endow $\Os{Q}$ with the locally convex vector topology making the linear map 
	\begin{displaymath}
	 \Lambda \colon \Os{Q} \rightarrow \prod_{(U,G,\psi) \in \vV} \vect{U}, \ [\hat{f}] \mapsto (f_U)_{(U,G,\psi) \in \vV}
	\end{displaymath}
 a topological embedding. Here we have used the unique lifts $f_U$ constructed in Remark \ref{rem: os:clpro}. We call this topology the \ind{orbisection!topology}{orbisection topology} and note that it is the initial topology with respect to the family of maps $\tau_U \colon \Os{Q} \rightarrow \vect{U}, [\hat{f}] \mapsto f_U,\ (U,G,\psi) \in \vV$.
 \item Let $\vV \coloneq \setm{(V_j, H_j, \psi_j)}{j \in J} \in \uU$ be a locally finite orbifold atlas such that each chart in $\vV$ is relatively compact. Endow $\Osc{Q}$ with the locally convex vector topology making the map
	\begin{displaymath}
	 \Lambda_{\vV} \colon \Osc{Q} \rightarrow \bigoplus_{j\in J} \vect{V_j}, \ [\hat{f}] \mapsto (f_{V_j})_{j\in J}
	\end{displaymath}
 from Lemma \ref{lem: os:btop} (b) a topological embedding. We call this topology the \ind{orbisection!topology!compactly supported}{compactly supported orbisection topology} (or  \emph{c.s.\ orbisection topology}).\glsadd{LambdavV} \\ With respect to this topology, the linear maps $\tau_{V_j} \colon \Osc{Q} \rightarrow \vect{V_j}, [\hat{f}] \mapsto f_{V_j}$ are continuous for each $(V_j,G_j,\psi_j) \in \vV$.
 \end{compactenum}
\end{defn}
A priori, the topologies defined on the spaces of orbisections might depend on the choice of orbifold atlas. However, as in the manifold case, we will see that neither the orbisection topology nor the c.s.\ orbisection topology depend on this choice. To prove the independence of the compactly supported orbisection topology of the choice of the orbifold atlas, relatively compact orbifold charts are needed. This explains the additional requirement in Definition \ref{defn: os:top}.

\begin{lem}\label{lem: ost:insat} \no{lem: ost:insat}
 Let $\wW = \setm{(W_i,G_i, \phi_i)}{i \in I} \in \uU$ be an arbitrary orbifold atlas for $Q$. 
	\begin{compactenum}
	 \item The orbisection topology with respect to $\vV$ is initial with respect to the family $(\tau_{W_i})_{(W_i,H_i,\phi_i) \in \wW}$.
	 \item Let $\wW$ be locally finite such that each chart in $\wW$ is relatively compact. The c.s.\ orbisection topology $\oO_\vV$ with respect to $\vV$ and the c.s.\ orbisection topology $\oO_\wW$ with respect to $\wW$ coincide. 
	\end{compactenum}
\end{lem}

\begin{proof}
	\begin{compactenum}
	 \item Consider the atlas $\wW \cup \vV$ obtained by joining the atlases $\vV$ and $\wW$. Clearly the orbisection topology induced by $\vV$ (respectively by $\wW$) is coarser than the orbisection topology induced by $\wW \cup \vV$. We claim that the orbisection topology induced by $\vV$ is finer than the one induced by $\vV \cup \wW$. If this is true then both orbisection topologies coincide. An analogous argument applies to the topology induced by $\wW$. Hence it suffices to prove that the orbisection topology induced by $\vV$ coincides with the one induced by $\wW \cup \vV$. Without loss of generality we may assume that $\vV$ is contained in $\wW$, i.e. $\wW = \wW \cup \vV$ holds.\\ Let $\tT$ be the initial topology on $\Os{Q}$ with respect to $(\tau_{W_i})_{(W_i,G_i,\phi_i) \in \vV}$. Fix $(U,H,\psi) \in \wW$, we have to show that $\tau_U \colon (\Os{Q} , \tT) \rightarrow \vect{U}$ is a continuous map. \\
 	 The open sets $\setm{\tilde{V}_i \coloneq \tilde{U} \cap \tilde{W}_i}{i \in I}$ form an open cover of $\tilde{U}$. Define $V_i \coloneq \psi^{-1} (\tilde{V}_i)$ to derive an open cover of $U$. By \cite[Lemma F.16]{hg2004}, the topology on $\vect{U}$ is initial with respect to the family $(\res^U_{V_i})_{i \in I}$. Since every $V_{i}$ satisfies $\psi (V_{i}) \subseteq \phi_i (W_i)$ by compatibility of orbifold charts, there is a family of changes of charts $(\lambda_{ik})_{k \in K_i}$ in $\CH{W_i}{U}$ such that $\bigcup_{k \in K_i} \cod \lambda_{ik} = V_{i}$.\\ 
 	 Another application of \cite[Lemma F.16]{hg2004} implies that the topology of $\vect{V_i}$ is initial with respect to $(\res^{V_i}_{\cod \lambda_{ik}})_{k \in K_i}$. Using transitivity of initial topologies, $\tau_U$ will be continuous with respect to $\tT$ if we can show that every 
	\begin{displaymath}
	 f_{ik} \coloneq \res^{U}_{\cod \lambda_{ik}} \circ \tau_U \colon \Os{Q} \rightarrow \vect{\cod \lambda_{ik}} 
	\end{displaymath}
 	is continuous for $i \in I ,k \in K_i$. But \cite[Lemma F.15 (a)]{hg2004} implies that the mapping $\res^{W_i}_{\dom \lambda} \colon \vect{W_i} \rightarrow \vect{\dom \lambda_{ik}}$ is continuous. Now we use that $$g_{\lambda_{ik}} \colon \vect{\dom \lambda_{ik}} \rightarrow \vect{\im \lambda_{ik}}, X \mapsto T\lambda \circ X \circ \lambda^{-1}$$ is continuous. To see this, observe that in charts (using Lemma \ref{lem: vect:top}), the mapping reduces to a pullback by a smooth map which is continuous, by \cite[Lemma 3.7]{hg2002}. We conclude from $f_{ik} = g_{\lambda} \res^{W_i}_{\dom \lambda} \tau_{W_i}$ that $\tau_U$ is continuous with respect to $\tT$ for every $(U,G,\psi) \in \vV$. Thus the orbisection topology with respect to $\vV$ is finer than $\tT$, whence both topologies coincide. 
	\item Consider $\vV = \setm{(V_j,H_j, \psi_j)}{j \in J}$. Notice that $\vV \cup \wW$ still is a locally finite atlas with relatively compact charts. After replacing $\wW = \setm{(W_i,G_i,\varphi_i)}{i\in I}$ with $\wW \cup \vV$, we may assume without loss of generality that $\vV \subseteq \wW$ holds. Let $\oO_\wW$ be the c.s.\ orbisection topology with respect to $\wW$ and $\oO_\vV$ be the c.s.\ orbisection topology with respect to $\vV$. Since $\vV$ is contained in $\wW$ the definition of the c.s.\ topology implies $\oO_\vV \subseteq \oO_wW$, i.e.\ the topology $\oO_\wW$ is finer than $\oO_\vV$.
	Conversely we have to prove that $\oO_\vV$ is finer than $\oO_\wW$. To see this, it suffices to prove that $\id_{\Osc{Q}} \colon (\Osc{Q}, \oO_\vV) \rightarrow (\Osc{Q}, \oO_\vV)$ is continuous, which follows from \cite[I.\ \S 1 6.\ Proposition 5]{bourbaki1987} if every zero-neighborhood in $\oO_\wW$ contains a zero-neighborhood in $\oO_\vV$. We proceed in three steps:
        \paragraph{Step 1: Zero-neighborhoods in $\vect{W_r}$ induce zero-neighborhoods in \rm $(\Osc{Q}, \oO_\vV)$}\mbox{}\\ 
        Consider an orbifold chart $(W_r,G_r,\varphi_r) \in \wW$. The projection $\text{pr}_r \colon \prod_{i \in I} \vect{W_i} \rightarrow \vect{W_r}$ and the canonical inclusion $I_\wW \colon \bigoplus_{i \in I} \vect{W_i} \rightarrow \prod_{i \in I} \vect{W_i}$ are continuous (cf.\ \cite[II.\ \S 4 5.\ Proposition 7]{bourbaki1987}). Furthermore, since $\vV \subseteq \wW$ holds, we identify each chart $(V_j,H_j,\psi_j)$ in $\vV$ with a chart $(W_{\alpha (j)}, G_{\alpha (j)}, \varphi_{\alpha (j)})$ in $\wW$. Then the canonical inclusion 
	\begin{displaymath}
	  I_{\vV,\wW} \colon \bigoplus_{j \in J} \vect{V_j} \rightarrow \bigoplus_{i \in I} \vect{W_i} , (f_j)_{j \in J} \mapsto (\tilde{f}_i)_{i \in I} \text{ with }\tilde{f}_i \coloneq  \begin{cases} 0 & \text{for } i \neq \alpha (j) \text{ for all } j \\ 
        f_j & \text{ if }  i = \alpha (j) \text{ for } j \in J                                                                                                                                                                                                                                                                                                                                                                                                                                                                                                                                                                                                         \end{cases}	 
	\end{displaymath}
        is continuous. Then $\Lambda_{W_r} \coloneq \text{pr}_r \circ I_{\wW} \circ I_{\vV,\wW} \colon \bigoplus_{i \in I} (V_i,H_i,\psi_i) \rightarrow \vect{W_r}$ is a continuous map. Now each zero-neighborhood $\Omega$ in $\vect{W_r}$ induces a zero-neighborhood $(\Lambda_{W_r}\circ \Lambda_\vV)^{-1} (\Omega)$ in $\oO_\vV$.\\[1em]
        Consider $[\hat{\sigma}] \in \Osc{Q}$ and denote its canonical lifts on $(W_i,G_i,\psi_i) \in \wW$ by $\sigma_{W_i}$. By Proposition \ref{prop: os:chch}, the canonical lift $\sigma_{W_r}$ is uniquely determined by the canonical lifts 
	\begin{displaymath}
	  \setm{\sigma_{V_j}}{(V_j,H_j,\psi_j) \in \vV \text{ with } \overline{\varphi_r (W_r)} \cap \im \psi_j \neq \emptyset}.
	\end{displaymath}
        Recall that all charts in $\wW$ are relatively compact and $\vV$ is a locally finite atlas. Thus for each $r \in I$, there is only a finite subset $J_r \subseteq J$ such that $\im \psi_i \cap \overline{\varphi_r (W_r)} \neq \emptyset$ holds if and only if $j \in I_r$.
        Denote the canonical inclusion $\bigoplus_{k \in J_r} \vect{V_k}\hookrightarrow \bigoplus_{j \in J} \vect{V_j}$ by $\iota_{J_r}$. By \cite[II. \S 4 5. Proposition 8 (i)]{bourbaki1987}, the map $\iota_{J_r}$ is continuous for each $J_r \subseteq J$. The maps $I_{\wW}$ and $I_{\vV,\wW}$ respectively, are (up to identification) just inclusions of subsets and $\text{pr}_r$ is a projection. Since the lift $\sigma_{W_r}$ of an orbisection $[\hat{\sigma}] \in \Osc{Q}$ is uniquely determined by the family of lifts indexed by $J_r$, we obtain for each open set $\Omega \in \vect{W_r}$ the following: 
	\begin{equation}\label{eq: umrechn} \begin{aligned}
	    \text{ The lift } \sigma_{W_r} \text{ is contained in } \Omega &\text{ if and only if } [\hat{\sigma}] \in (\Lambda_{W_r} \circ \Lambda_\vV)^{-1} (\Omega),\\ &\text{ if and only if } (\sigma_{V_k})_{k \in J_r} \in (\Lambda_{W_r} \circ \iota_{J_r})^{-1} (\Omega).
	  \end{aligned} 
	\end{equation} 
  	\paragraph{Step 2: The countable case} We shall assume for this step only that the atlases $\vV, \wW$ are indexed by countable sets $I,J$. \\ Consider the vector spaces $(\bigoplus_{i \in I} \vect{W_i})_{\text{box}}$ and $(\bigoplus_{j \in J} \vect{V_j})_\text{box}$ respectively, endowed with the box-topology. Since $I,J$ are countable, the box topology coincides with the locally convex direct sum topology by \cite[Proposition 4.1.4]{jarchow1980}. A typical zero-neighborhood in $\bigoplus_{i \in I} \vect{W_i}$ is given by $U \coloneq \oplus_{i \in I} U_i$, where $U_i \subseteq \vect{W_i}$ is an open set. For each $i \in I$ choose by Step 1 an open box neighborhoods $B^i \coloneq\bigoplus_{\alpha \in J_i} B^i_\alpha$ such that $B^i \subseteq (\Lambda_{W_i} \circ \iota_{J_i})^{-1} (U_i)$.
	Reformulating Condition \eqref{eq: umrechn} this yields: If $\sigma_{V_\alpha} \in B^i_\alpha$ holds for all $\alpha \in J_i$, then $\sigma_{W_i} \in U_i$ follows. 
        Using the boxes defined above, we construct sets $\Omega_j \coloneq \bigcap_{i \in I_j} B_j^i$. Recall that $\vV$ contains only relatively compact charts and $\wW$ is locally finite. Thus for fixed $j \in J$ the set $I_j \coloneq \setm{j \in J_i}{i \in I}$ is finite, whence the set $\Omega_j$ is an open zero-neighborhood of $\vect{V_j}$. Now $B \coloneq \bigoplus_{j \in J} \Omega_j$ is a box zero-neighborhood in $\bigoplus_{j \in J} \vect{V_j}$. The open box-neighborhood $B$ contains only elements of $\bigoplus_{j \in J} \vect{V_j}$ which are mapped by the projection $\bigoplus_{j \in J} \vect{V_j} \rightarrow \bigoplus_{k \in J_i} \vect{V_k}$ into $\oplus_{k \in J_i} B^i_k$ for each $i \in I$. We obtain the following condition for an orbisection $[\hat{\sigma}] \in \Osc{Q}$ with families of canonical lifts $(\sigma_{V_j})_{j \in J}$ with respect to $\vV$ and $(\sigma_{W_i})_{i \in I}$ with respect to $\vV$: 
	\begin{align*}
	 [\hat{\sigma}] \in \Lambda_\vV^{-1} (B) &\Leftrightarrow (\sigma_{V_j})_{j \in J} \in B \Rightarrow (\forall i \in I) (\sigma_{V_j})_{j \in J_i} \in \bigoplus_{k \in J_i} B^{i}_k \\
	  & \Rightarrow (\forall i \in I) \ \ \sigma_{W_i} \in U_i \Rightarrow [\hat{\sigma}] \in \Lambda_\wW^{-1} (U) 
%
	\end{align*}
 	In other words, the typical zero-neighborhood $\Lambda_\wW^{-1} (U)$ in $\oO_\wW$ contains the zero-neighborhood $\Lambda_\vV^{-1} (B) \in \oO_\vV$. As sets of the form $\Lambda_\wW^{-1} (\bigoplus_{i \in I} U_i)$ form a base of zero-neighborhoods in $\oO_\wW$, we deduce $\oO_\wW \subseteq \oO_\vV$ and thus $\oO_\wW = \oO_\vV$. Furthermore, the map $\rho \coloneq \Lambda_\wW|^{\im \Lambda_\wW} \circ (\Lambda_\vV|^{\im \Lambda_\vV})^{-1})$ is an isomorphism of topological vector spaces.
	\paragraph{Step 4: The general case} In general neither $\vV$ nor $\wW$ need to be countable (since the orbifolds we consider need not be $\sigma$-compact). Orbifold charts are connected, whence each chart is contained in exactly one connected component. Let $\cC$ be the family of connected components of $Q$ and for $C \in \cC$ and an atlas $\aA$ define $\aA_{C} \coloneq \setm{(V,H,\psi) \in \aA}{\psi (V) \subseteq C}$. The subset $\aA_C$ is an atlas of orbifold charts for the component $C$. We may split the atlases $\vV$, $\wW$ into disjoint unions $\vV = \bigsqcup_{C \in \cC} \vV_C$ resp. $\wW = \bigsqcup_{C \in \cC} \wW_C$. By construction, $\vV_C$ is still contained in $\wW_C$ 
        Decompose the direct sums 
	\begin{displaymath}
	 \bigoplus_{i\in I} \vect{W_i} = \bigoplus_{C \in \cC} \left(\bigoplus_{(W, G, \phi) \in \wW_C} \vect{W} \right) \quad \quad  \bigoplus_{j\in J} \vect{V_j} = \bigoplus_{C \in \cC} \left(\bigoplus_{(V, H, \psi) \in \vV_C} \vect{V} \right)
	\end{displaymath}
 	and observe that the maps $\Lambda_\vV$ and $\Lambda_\wW$ decompose as $\Lambda_\vV = (\Lambda_{\vV_C})_{C \in \cC}$ and $\Lambda_\wW = (\Lambda_{\wW_C})_{C \in \cC}$. Every connected component $C \subseteq Q$ is $\sigma$-compact by Proposition \ref{prop: ofd:prop} (d). Since $\wW_C$ and $\vV_C$ are locally finite, both atlases have to be countable. Step 3 yields for each connected component $\cC$ an isomorphism $\rho_C =\Lambda_{\wW_C}|^{\im \Lambda_{\wW_C}} \Lambda_{\vV_C}^{-1}|^{\im \Lambda_{\vV_C}} \colon \im \Lambda_{\vV_C} \rightarrow \im \Lambda_{\wW_C}$. Taking direct sums in the category of topological vector spaces is functorial. Therefore the map $\oplus_{C \in \cC} \rho_C \colon \bigoplus_{C \in \cC} \im \Lambda_{\vV_C} \rightarrow \bigoplus_{C \in \cC} \im \Lambda_{\wW_C}$ is an isomorphism of locally convex topological vector spaces. Observe that the families of canonical inclusions (of vector subspaces) $\iota_\cC \colon \im \Lambda_{\vV_C} \hookrightarrow \bigoplus_{(V, H, \psi) \in \vV_C} \vect{V}$ respectively $\iota_C' \colon \im \Lambda_{\wW_C} \hookrightarrow \bigoplus_{(W, G, \phi) \in \wW_C} \vect{W}$ induce continuous linear maps $\iota \coloneq \oplus_{C \in \cC} \iota_C$ and $\iota' \coloneq \oplus_{C \in \cC} \iota_C'$, respectively. By \cite[II.6, Proposition 8]{bourbaki1987}, the subspace topology on $\im \iota$ turns $\iota$ into an isomorphism of topological vector spaces and the same holds for $\iota'$ and the subspace topology on $\im \iota'$. We deduce that 
	\begin{displaymath}
	 \Lambda_\wW|^{\im \Lambda_\wW} \circ (\Lambda_\vV|^{\im \Lambda_\vV})^{-1} = \iota' \circ \bigoplus_{C \in \cC}\left( \Lambda_{\wW_C}|^{\im \Lambda_{\wW_C}} \circ (\Lambda_{\vV_C}|^{\im \Lambda_{\vV_C}})^{-1}\right) \circ \iota^{-1}
	\end{displaymath}
	is an isomorphism of topological vector spaces. Thus $\oO_\vV = \oO_\wW$ holds.
	\end{compactenum}
\end{proof}

To illustrate the construction of the orbisection topologies, we consider the special case of orbisections on an orbifold with a global chart. It turns out that we may then identify the orbisections with subspaces of vector fields on the global chart.

\begin{ex}\label{ex: os:good}  Let $d \in \NN$. Consider a finite subgroup $G \subseteq \Diff (\RR^d)$. We define an orbifold structure on $Q \coloneq \RR^n /G$ via the atlas $\vV \coloneq \set{(\RR^d , G ,\pi)}$, where $\pi \colon \RR^d \rightarrow \RR^d/G$ is the quotient mapping.
 \begin{compactenum}
  \item By Proposition \ref{prop: os:cl}, each orbisection $[\hat{\sigma}] \in \Os{Q}$ can be identified with a unique vector field in $\vect{\RR^d}$. Since the group elements are changes of charts, for the canonical lift of an orbisection on the global chart $g. X = Tg \circ X = X \circ g$ holds for each $g \in G$. Thus the canonical lifts are $G$-equivariant vector fields. Hence by Lemma \ref{lem: ost:insat}, the map $\Lambda_\vV \colon \Os{Q} \rightarrow \vect{\RR^d}$ (cf.\ Lemma \ref{lem: os:btop}) establishes an isomorphism of topological vector spaces between $\Os{Q}$ and the space of all $G$-equivariant vector fields $\mathfrak{X}^G (\RR^d)$.\\ 
  Observe that $\mathfrak{X}^G (\RR^d)$ is a closed subspace of $\vect{\RR^d}$. To prove this, recall that for each $p \in \RR^d$ the point evaluation $\text{ev}_p \colon C^\infty (\RR^d , \RR^d) \rightarrow \RR^d$ is continuous by \cite[Proposition 3.20]{alas2012}. Hence for each pair $(p,g) \in \RR^d\times G$, the map $E_{p,g} \colon C^\infty (\RR^d,\RR^d) \rightarrow \RR^d , f \mapsto d g (p ,\cdot) \circ \text{ev}_p (f) - \text{ev}_{g(p)} (f)$ is continuous. We may then identify $\mathfrak{X}^G (\RR^d)$ with the closed vector subspace $\bigcap_{p \in \RR^d} \bigcap_{g \in G} E_{p,g}^{-1} (0)$.   
  \item We identify the compactly supported orbisections $\Osc{Q}$ with the set of equivariant compactly supported vector fields of $\RR^d$. To this end, consider     
  \begin{displaymath} 
     \mathfrak{X}^G_c (\RR^d) \coloneq \setm{X \in \mathfrak{X}_c (\RR^d)}{\forall g \in G, Tg \circ X = X \circ g}
    \end{displaymath}
  as a subspace of $\mathfrak{X}_c (\RR^d)$ (cf.\ Definition \ref{defn: cs:vf}). We claim that $\Osc{Q}$ and $\mathfrak{X}^G_c (\RR^d)$ are isomorphic as topological vector spaces. To this end, choose a locally finite orbifold atlas $\wW = \setm{(U_i,G_i,\pi_i)}{i \in I}$ for $Q$ with $I$ countable. By Lemma \ref{lem: locfin}, we can choose $\wW$ such that for each $i \in I$ the set $U_i \subseteq \RR^d$ is a relatively compact open subset such that the inclusion of sets induces an embedding of orbifold charts. Then $\RR^d = \bigcup_{i \in I} G. U_i$ holds, as $\wW$ is an orbifold atlas for $Q = \RR^d /G$. Since $G$ is a finite group we may assume that for each $i \in I$ and $g \in G$ there is $j \in J$ with $U_j = g. U_i$ and $G_j = g.G_i.g^{-1}$. Thus $(U_i)_{i \in I}$ is a locally finite cover of $\RR^d$ by relatively compact open subsets such that the cover is countable. Recall from the definition of the topologies that the rows in the following commutative diagram are topological embeddings with closed image (cf.\ Lemma \ref{lem: os:btop} and Definition \ref{defn: cs:vf})
    \begin{displaymath}
	 \begin{xy}
  			\xymatrix{
                                                                    \Osc{Q} \ar[rrr]^{\Lambda_\wW} \ar@{.>}[d]&&& \bigoplus_{i \in I} \vect{U_i} \ar[d]^{\cong}_\theta\\
     			 	 \mathfrak{X}_c^G (\RR^d) \ar[r]^{\subseteq} & \mathfrak{X}_c (\RR^d) \ar[rr]^-{R_\wW} && \bigoplus_{i \in I} C^\infty (U_i,\RR^d) 
 				 }
			\end{xy}
    \end{displaymath}
 Here the isomorphism $\theta$ is defined via $(f_i)_{i \in I} \mapsto (\text{pr}_2 \circ f_i)_{i \in I}$. As canonical lifts of an orbisection are unique by Lemma \ref{lem: os:cucr}, Proposition \ref{prop: os:cl} and a trivial computation yield $\im \theta \Lambda_\wW \subseteq \im R_\wW$. Furthermore, the image $\im R_\wW^{-1} \theta \Lambda_\wW$ coincides with $ \mathfrak{X}_c^G (\RR^d)$. Denote by $\sigma_{\RR^d}$ the canonical lift of $[\hat{\sigma}] \in \Osc{Q}$ with respect to the global chart. Then $R_\wW^{-1} \theta \Lambda_\wW$ induces the isomorphism of topological vector spaces 
 \begin{displaymath}\Osc{Q} \rightarrow \textstyle\mathfrak{X}_c^G \displaystyle (\RR^d) , [\hat{\sigma}] \mapsto \sigma_{\RR^d}.\end{displaymath}
 Observe that $\mathfrak{X}_c^G (\RR^d)$ is a closed vector subspace of $\mathfrak{X}_c (\RR^d)$. This follows from part (a) and the following facts: The inclusion $\iota \colon \bigoplus_{i \in I} C^\infty (U_i,\RR^d) \rightarrow \prod_{i \in I} C^\infty (U_i,\RR^d)$ is continuous by \cite[II. \S 4, Proposition 7]{bourbaki1987}. By definition of the topology on $\vect{\RR^d}$, we may identify $\mathfrak{X}^G (\RR^d)$ with a closed vector subspace $A$ of $\prod_{i \in I} C^\infty (U_i,\RR^d)$ such that $\iota^{-1} (A) = R_\wW (\mathfrak{X}_c^G (\RR^d))$ holds. Hence the assertion follows by continuity, since $\im R_\wW$ is a closed subspace.
 \end{compactenum}
We conclude that for the orbifold $Q = \RR^d/G$, the space $\Os{Q}$ corresponds to $\mathfrak{X}^G (\RR^d)$. Also $\Osc{Q}$ corresponds to $\mathfrak{X}_c^G (\RR^d)$.\\ We remark that a similar result holds for arbitrary orbifolds with a global chart, by essentially the same argument.
\end{ex}

\begin{thm}\label{thm: dorb:top} \no{thm: dorb:top}
 Let $(Q,\uU)$ be a second countable orbifold, i.e.\ $Q$ is a second countable space (or equivalently $Q$ is a $\sigma$-compact space). The topological vector space $\Os{Q}$ is then a Fr\'{e}chet space.
\end{thm}

\begin{proof}
 As $Q$ is second countable, there is a countable orbifold atlas $\setm{(U_i,G_i,\psi_i)}{i \in \NN}$ for $Q$. By Lemma \ref{lem: ost:insat}, the orbisection topology is initial with respect to the maps 
	\begin{displaymath}
        \tau_{U_i} \colon \Os{Q} \rightarrow \vect{U_i}, [\hat{f}] \mapsto f_{U_i}.
	\end{displaymath}
 In particular, Lemma \ref{lem: os:btop} yields a linear topological embedding 
	\begin{displaymath}
	 \Lambda \colon \Os{Q} \rightarrow \prod_{i \in \NN} \vect{U_i}, [\hat{f}] \mapsto (f_{U_i})_{i\in I}
	\end{displaymath}
 onto a closed subspace. The manifolds $U_i$ are finite-dimensional, connected and paracompact manifolds. Thus by Proposition \ref{prop: para:coco}, every $U_i$ is $\sigma$-compact and second countable. The space $\RR^n$ is a Fr\'{e}chet space over the locally compact field $\RR$. Combining these observations with Lemma \ref{lem: vect:top} and \cite[Proposition 4.19]{hg2004}, $\vect{U_i}$ with the topology defined in Definition \ref{defn: top:vect} is a Fr\'{e}chet space for each $i \in I$. The countable product of Fr\'{e}chet spaces is a Fr\'{e}chet space (combine \cite[I.\ \S 3 2.]{bourbaki1987} with \cite[Proposition 3.3.6]{jarchow1980}) and thus $\prod_{i \in I} \vect{U_i}$ is a Fr\'{e}chet space. From Lemma \ref{lem: os:btop} and Lemma \ref{lem: ost:insat}, we deduce that $\Os{Q}$ is isomorphic to a closed vector subspace of the Fr\'{e}chet space $\prod_{i \in I} \vect{U_i}$. Thus $\Os{Q}$ is a Fr\'{e}chet space.
\end{proof} 

\begin{cor}\label{cor: ost:prop} \no{cor: ost:prop}
 \begin{compactenum}
  \item The spaces $\Os{Q}$ with the orbisection topology and $\Osc{Q}$ with the c.s.\ orbisection topology are Hausdorff and complete topological vector spaces.
  \item If $(Q,\uU)$ is a compact orbifold, then the locally convex vector spaces $\Os{Q}$ and $\Osc{Q}$ coincide. If $Q$ is compact, then both spaces are Fr\'{e}chet spaces.
  \item Let $\vV$ be a locally finite orbifold atlas for $Q$ which consists of relatively compact charts. The family $(\tau_V)_{(V,G,\psi)\in \vV}$ as in Definition \ref{defn: os:top} (b) forms a patchwork for $\Osc{Q}$, turning it into a patched locally convex space. The topological embedding is given by $\Lambda_\vV$ (cf.\ Definition \ref{defn: patloc}).
 \end{compactenum}
\end{cor}

\begin{proof}
 \begin{compactenum}
  \item We endow the space of vector fields on a finite-dimensional manifold with the topology introduced in Definition \ref{defn: top:vect}. Recall that direct products and direct sums of Hausdorff and complete locally convex vector spaces are again such spaces by \cite[Proposition 4.3.3, Proposition 4.3.6 and Proposition 4.4.3]{jarchow1980}. The assertion follows from \cite[Remark F.8]{hg2004}, since the spaces $\Os{Q}$ and $\Osc{Q}$ with the topology of Definition \ref{defn: os:top} are isomorphic to closed subspaces of complete and Hausdorff spaces. 
  \item For finite index sets products and direct sums are canonically isomorphic. As locally finite covers of compact spaces are finite, together with Theorem \ref{thm: dorb:top} this proves the claim.  
  \item Follows directly from the definition of the c.s.\ orbisection topology (Definition \ref{defn: os:top}). 
 \end{compactenum}
\end{proof}

\begin{lem}\label{lem: OsK} \no{lem: OsK}
 Let $K \subseteq Q$ be a compact subset and endow $\OsK{Q} \subseteq \Osc{Q}$ with the subspace topology. The space $\OsK{Q}$ is a closed subspace of $\Osc{Q}$.
\end{lem}

\begin{proof}
 Choose an arbitrary locally finite orbifold atlas $\vV \coloneq \setm{(V_i,G_i,\psi_i)}{i \in I}$ for $(Q,\uU)$. By Lemma \ref{lem: ost:insat} (b), there is a topological embedding  $\Lambda_\vV \colon \Osc{Q} \rightarrow \bigoplus_{i \in I} \vect{V_i}$ whose image is closed.
 For each $i\in I$, we obtain a (possibly empty) subset $U_i \coloneq \psi_i^{-1} (Q \setminus K)$. If $U_i = \emptyset$ holds, define $A_i \coloneq \vect{V_i}$. Otherwise, consider $x \in U_i$ and a manifold chart $(W_\psi , \psi)$ for $V_i$ such that $x \in W_\psi$. The evaluation map $\text{ev}^\psi_{x} \colon C^\infty (W_\psi, \RR^d) \rightarrow \RR^d , \xi \mapsto \xi (x )$ is continuous by \cite[Proposition 11.1]{hg2004}. As the topology on $\vect{V_i}$ is initial with respect to the maps $\theta_\psi \colon \vect{V_i} \rightarrow C^\infty (W_\psi, \RR^d), X \mapsto X_\psi$, the point evaluation $\text{ev}_x \colon \vect{V_i} \rightarrow \RR^d , \sigma \mapsto \text{ev}_x^\psi \circ \theta_\psi (\sigma)$ is continuous. Hence we obtain a closed set $A_i \coloneq \bigcap_{x \in U_i} \text{ev}_x^{-1} (0)$. From \cite[II.\ \S 4 5.\ Corollary 1]{bourbaki1987}, we conclude that $A \coloneq \oplus_{i \in I} A_i = \prod_{i \in I} A_i \cap \bigoplus_{i \in I} \vect{V_i}$ is closed. By construction, each orbisection in $\Lambda_\vV^{-1} (A)$ vanishes off $K$, whence its support must be contained in $K$. We deduce $\Lambda_\vV^{-1} (A) = \OsK{Q}$, whence $\OsK{Q}$ is a closed set.
\end{proof}


The results in this section suggest that orbisections behave in many ways as vector fields for finite dimensional manifolds. Before we end this section, we point out that in some ways orbisections do \emph{not} behave like vector fields. There may be formal orbifold tangent vectors which are \emph{not contained} in the image of any orbisection. In the manifold case, this may never occur. The following example was first considered by Borzellino et al.\ (see \cite[Example 43]{bb2008}) in the context of their notion of orbifold maps:

\begin{ex}\label{ex: ad:nofd} \no{ex: ad:nofd}
 Consider $\RR $, with an action induced by the linear diffeomorphism $\gamma \colon \RR \rightarrow \RR , x \mapsto -x$. Set $G \coloneq \langle \gamma \rangle$ and let $\psi \colon \RR \rightarrow \RR /G$ be the quotient map to the orbit space. The quotient is homeomorphic to $Q \coloneq [0,\infty[$ (as a subspace of $\RR$). By abuse of notation we obtain an orbifold atlas $\uU \coloneq \set{(\RR , G , \psi)}$ for $Q$. Now $(Q,\uU)$ is an orbifold and the local groups are trivial for every point except $0$ (where it is isomorphic to $G$). We may thus compute the tangent spaces of $Q$ at $x \in Q$ in the following way:\\
  For $x \neq 0$ we have $\tT_x Q \cong \RR$ and $\tT_0 Q \cong [0, \infty[$. An atlas for the tangent orbibundle is induced by the orbifold chart $(T\RR , G, T\psi)$, where $G$ acts on $T\RR$ via the derived action. Taking identifications we obtain $T\RR \cong \RR^2$. The group $G$ acts via elements of $O(1)$ on $\RR$. Hence its action on $T\RR$ is induced by the linear map $T\gamma \colon \RR^2 \rightarrow \RR^2, \ (x,y) \mapsto (-x,-y)$. The topological base space of the tangent orbibundle is thus $\tT Q = \RR^2 /G$. 
 The zero vector is the only fixed point of the derived action of $G$. Since orbisections preserve local groups by Proposition \ref{prop: os:plg}, every orbisection maps $0 \in Q$ to $0 \in \RR^2/G \cong \tT Q$. Thus all orbisections in $\Os{Q}$ must vanish in $0 \in Q$ and 
	\begin{displaymath}
	 Q' \coloneq \bigcup_{(f, \set{\hat{f}_{(\RR , \ZZ_2 , \psi)}}, P, \nu) \in \Os{Q}} \Im f \subsetneq \tT Q
	\end{displaymath}
 Is the topological subspace $Q'$ at least an orbifold? We shall prove that the answer to this question is negative. Indeed it will turn out that $Q'$ is not locally compact.\\
 Following Remark \ref{rem: os:clpro} (d), the set $Q'$ is homeomorphic to $T\RR^{\text{inv}}/G$, i.e.\ it is homeomorphic to $(T (\RR\setminus \set{0}) \cup \set{0 \in T_0\RR}) /G$. Since $T(\RR \setminus \set{0}) \cup \set{0 \in T_0\RR} \sim \RR \setminus \set{0} \times \RR \cup \set{(0,0)}$ is not locally compact, Lemma \ref{lem: orbitmap} (e) implies that $Q'$ is not locally compact.
\end{ex}
\thispagestyle{empty}
\thispagestyle{empty}
\section{Riemannian Geometry on Orbifolds}

\setcounter{subsubsection}{0}

In this section, the notion of a Riemannian orbifold metric is recalled. Our approach follows the construction of Riemannian metrics on manifolds (cf. \cite[Ch.\ 1.2, Proposition 2.10]{rg1992}). The corresponding construction of such an object for an orbifold is well known (see for example \cite[Proposition 2.20]{follie2003}; we also recommend the survey in \cite[Appendix 4.2]{chen2001}). Nevertheless, the results are repeated here for the readers convenience, and to fix some notation.

\begin{defn}[Riemannian orbifold metric]\label{defn: rieofdm} \no{defn: rieofdm}
 Let $(Q,\uU)$ be an orbifold and consider some orbifold atlas $\vV = \setm{(V_i,G_i,\psi_i)}{i\in I}$ for $(Q,\uU)$. A \ind{}{Riemannian orbifold metric} on $Q$ is a collection $\rho = (\rho_i)_{i \in I}$, where $\rho_i$ is a Riemannian metric on the manifold $V_i$ such that the following holds:
  \begin{compactitem}
   \item[\emph{(Compatibility)}]  For each $(i,j) \in I \times I$ and each open $G_i$-stable subset $S \subseteq V_i$, every embedding of orbifold charts $\lambda \colon (S, (G_i)_S, \psi_i|_{S}) \rightarrow (U_j,G_j,\psi_j)$ is a Riemannian embedding, i.e.
	\begin{displaymath}
	 \rho_j (T_x \lambda (v) , T_x \lambda (w)) = \rho_i (v,w)  \quad \forall v,w \in T_x V_i, \ x \in S.
	\end{displaymath}
  \end{compactitem}
 Let $(Q,\uU)$ be an orbifold endowed with a Riemannian orbifold metric $\rho$. The triple $(Q,\uU,\rho)$ is called a \ind{orbifold!Riemannian}{Riemannian orbifold}. 
\end{defn}

\begin{rem}\label{rem: isom} \no{rem: isom}
 Consider a Riemannian orbifold metric $\rho$ on some orbifold $(Q,\uU)$, associated tp an atlas $\vV$ as above. For a chart $(V,G, \psi) \in \vV$, the group $G$ acts by self-embeddings of orbifold charts. If $V$ is endowed with a member $\rho_i$ of $\rho$, each element of $G$ thus acts as a Riemannian isometry with respect to $\rho_i$.
\end{rem}

\begin{prop}[\hspace{-0.5pt}{\cite[Proposition 2.20]{follie2003}}] \label{prop: ex:rieofd} \no{prop: ex:rieofd}
 Any orbifold $(Q,\uU)$ admits a Riemannian orbifold metric $\rho$.
\end{prop}

\begin{proof}
 Let $\vV = \setm{(V_i, G_i, \psi_i)}{i \in I}$ be any representative of $\uU$, and $\set{\hat{\chi}_i}_{i \in I}$ be a smooth orbifold partition of unity subordinate to $\vV$, which exists due to Proposition \ref{prop: pu:hofd}. Recall from \ref{nota: opu} that for every pair $(i,j) \in I\times I$, there is a smooth lift $\chi_{i,j}$ of $\chi_i$ to $(V_j,G_j,\psi_j)$. For $i \in I$, choose some Riemannian metric $m^{(i)}$ on $V_i$ (cf.\ \cite[VII., \S 1, Proposition 1.1]{langdgeo2001}). As $G_i$ acts by diffeomorphisms, we obtain pullback metrics on $V_i$. Averaging over $G_i$, on every tangent space there is a positive definite bilinear form:
 	\begin{displaymath}
 	 \langle v , w \rangle^{(i)}_p \coloneq \frac{1}{\lvert G_i\rvert} \sum_{g \in G_i} m_{g.p}^{(i)} (T_pg . v , T_pg.w), \quad \forall v,w \in T_p V_i , \ p \in V_i
 	\end{displaymath}
 such that the family $\langle - , - \rangle^{(i)} \coloneq (\langle - , - \rangle^{(i)}_p)_{p \in V_i}$ defines a Riemannian metric on $V_i$. By construction, each element of $G_i$ is a Riemannian isometry with respect to $\langle - , - \rangle^{(i)}$.\\
 Define a Riemannian metric $\rho_i$ on $V_i$ as follows: Because $(\supp \chi_i)_{i \in I}$ is locally finite, $\psi_i(p)$ with $p \in V_i$ is contained in $\supp \chi_i$ for only finitely many $i \in I$. Therefore there is an open $G_i$-stable subset $p \in S_p \subseteq V_i$ such that for $y \in S_p$, $\psi_i (y) \in \supp \chi_k$ can hold only if $\psi_i (p) \in \supp \chi_k$ for $k \in I$. Shrinking $S_p$, without loss of generality for each $k \in I$ with $\psi_i (p) \in \supp \chi_k$ there is an embedding of orbifold charts $\lambda_k^p \colon (S_p, (G_i)_{S_p}, \psi_i|_{S_p}) \rightarrow V_k$. If $\psi_i (p) \not \in \supp \chi_k$ simply let $\lambda_k^p \colon S_p \rightarrow V_p$ be constant (whence $T_p \lambda_k^p = 0$) and define for $v,w \in T_p V_i$: 
	\begin{displaymath}
	 (\rho_i)_p (v,w) \coloneq \sum_{j \in I} \chi_{j,i} (p) \cdot \langle T_p \lambda_j^p (v), T_p \lambda_j^p(w) \rangle^{(j)}_{\lambda_j (p)}
	\end{displaymath}
 Since the $\chi_{j,i}$ are the lifts of an orbifold partition of unity, only finitely many terms are non-zero and $(\rho_i)_p$ is a positive definite bilinear map on $T_p V_i \times T_p V_i$. The definition of $(\rho_i)_p$ neither depends on $S_p$ nor on the choice of $\lambda_k^p$:\\ To prove this, consider another $G_i$-stable set $p \in S_p'$ with embeddings $\mu_k^p$. Since we are only interested in the tangent map at $p$ (which may be computed in an arbitrarily small open subset), we restrict $\mu_k^p$ and $\lambda_k^p$ to an open and $G_i$-stable subset $S \subseteq S_p \cap S_p'$ which contains $p$. If $\psi_i (p) \not \in  \supp \chi_k$, the contribution to $(\rho_i)_p (v,w)$ is zero. Otherwise, Proposition \ref{prop: ch:prop} (d) implies that there is a group element $g \in G_k$ such that $\mu_k^p|_S = g \circ \lambda_k^p|_S$. By construction, every $g \in G_k$ is a Riemannian isometry with respect to $\langle - , - \rangle^{(k)}$. Thus every choice induces the same map.\\
 The maps $\lambda_j^p$, $\chi_{k,i}$ are smooth and $\langle - ,-\rangle^{(k)}$ is a Riemannian metric for each $k \in I$, thus the family $\rho_i \coloneq ((\rho_i)_p)_{p \in V_i}$ defines a smooth map on each open set $TS_p \oplus TS_p \subseteq TV_i \oplus TV_i$. By construction the map does not depend on the set $S_p$ and thus $\rho_i$ is smooth on $TV_i \oplus TV_i$. Hence it is a Riemannian metric on $V_i$. \\
 We claim that the family $(\rho_i)_{i \in I}$ satisfies the compatibility condition of Definition \ref{defn: rieofdm}: Consider arbitrary $i,j \in I$ together with an open $G_i$-stable subset $S \subseteq V_i$ and an embedding of orbifold charts $\mu \colon (S, (G_i)_S, \psi_i|_S) \rightarrow (V_j , G_j , \psi_j)$. For $p \in S$ and $v,w \in T_p V_i$, we have to show that $(\rho_j)_{\mu (p)} (T_p \mu (v), T_p \mu (w))$ coincides with $(\rho_i)_p (v,w)$.\\
 Since $\mu$ is an embedding of orbifold charts and by construction one has $\chi_{k,j} = \chi_k \circ \psi_j$, we derive $\chi_{k,j} \circ \mu = \chi_{k,i}|_{\dom \mu}$. We compute:
	\begin{align*}
	 (\rho_j)_{\mu (p)} (T_p \mu (v), T_p \mu (w)) &= \sum_{k \in I} \chi_{k,j} (\mu (p)) \cdot \langle T_{ \mu (p)} \lambda_k^{\mu (p)} T_p\mu (v), T_{ \mu (p)} \lambda_k^{\mu (p)} T_p\mu (w) \rangle^{(k)}_{\lambda_k^{\mu (p)} \mu (p)} \\
		&= \sum_{k \in I} \chi_{k,i} (p) \cdot \langle T_p (\underbrace{\lambda_k^{\mu (p)} \mu}_{\theta_k^p \coloneq}) (v), T_p (\underbrace{\lambda_k^{\mu (p)} \mu}_{=\theta_k^p} )(w) \rangle^{(k)}_{ \lambda_k^{\mu (p)} \mu(p)}\\
		&= \sum_{k \in I} \chi_{k,i} (p) \cdot \langle T_p \theta_k^p (v), T_p \theta_k^p (w) \rangle^{(k)}_{\theta_k^p (p)}.
	\end{align*}
 Restrict every non-constant map $\theta_k^p$ to a small open $G_i$-stable neighborhood of $p$ such that the restriction of $\theta_k^p$ yields an embedding of orbifold charts (cf. \cite[Proposition 2.13]{follie2003}). As the definition of the metric does not depend on the choice of embedding, indeed we obtain
	\begin{displaymath}
	 (\rho_j)_{\mu (p)} (T_p \mu (v), T_p \mu (w)) = (\rho_i)_p (v,w).
\end{displaymath}
The family $\rho$ is compatible as in Definition \ref{defn: rieofdm}, whence it is a Riemannian orbifold metric.
\end{proof}

A Riemannian orbifold metric (uniquely) extends to each representative of the orbifold structure:

\begin{prop}\label{prop: ext:rom} \no{prop: ext:rom}
 Let $(Q,\uU)$ be an orbifold and $\vV = \setm{(V_i,G_i,\psi_i)}{i\in I}$ some representative of $\uU$ for which there is a Riemannian orbifold metric $\rho = (\rho_i)_{i \in I}$. For each representative $\aA$ of $\uU$, there exists a unique Riemannian orbifold metric $\hat{\rho}$ which extends $\rho$ to $\vV \cup \aA$.
\end{prop}

\begin{proof}
 We construct a Riemannian metric on $(U,H,\phi) \in \vV \cup \aA$ as follows: For $q \in U$ choose an $H$-stable subset $q \in S_q \subseteq U$ together with an embedding $\tau_i^q \colon (S_q, H_{S_q}, \phi|_{S_q}) \rightarrow (V_i,G_i,\psi_i)$ for some $i \in I$. Define for $v,w \in  T_q U$
    \begin{displaymath}
     (\hat{\rho}_U)_q (v,w) \coloneq \rho_i (T_q \tau_i^q (v) , T_q \tau_i^q (w))
    \end{displaymath}
 Repeating this construction for each $q \in U$, arguments as in the proof of Propostion \ref{prop: ex:rieofd} show that $\hat{\rho}_U$ is a well-defined Riemannian metric on $U$. In particular, the $\hat{\rho}_U$ does not depend on the choices involved in its construction. Since in the above construction, we may always choose the inclusion $S_q \subseteq U$ for a chart $(U,G,H) \in \vV$, one obtains $\hat{\rho}_U = \rho_U$ for $(U,H,\phi) \in \vV$.\\
 Finally, the family $(\hat{\rho}_U)_{(U,G,\phi) \in \vV \cup \aA}$ satisfies the compatibility condition of Definition \ref{defn: rieofdm}. To see this, consider a change of charts $\lambda \in \Ch_{\vV \cup \aA}$. It suffices to check the compatibility condition for each $q \in \dom \lambda \subseteq U$ separately. By construction, there are embeddings of orbifold charts $\tau_q^i \colon S_q \rightarrow V_i$ and $\tau_{\lambda (q)}^j \colon S_{\lambda (q)} \rightarrow V_j$ into charts $(V_i,G_i,\psi_i), (V_j,G_j,\psi_j)\in \vV$. Then we compute for $v,w \in T_q \dom \lambda$: 
  \begin{displaymath}
   (\hat{\rho}_{\cod \lambda})_{\lambda (q)} (T_q \lambda (v) , T_q \lambda (w)) = \rho_j (T_q \tau_{\lambda (q)}^j \lambda (v) , T_q \tau_{\lambda (q)}^j \lambda (w)) \stackrel{(\star)}{=} \rho_i (T_q \tau_q^i (v) ,T_q \tau_q^i (w)) = (\hat{\rho}_{U})_q (v,w) .
  \end{displaymath}
 Here the identity $(\star)$ follows from the compatibility of the Riemannian orbifold metric $(\rho_i)_{i \in I}$ and the fact that on a neighborhood $\Omega$ of $\tau_q^i (q)$ the mapping $(\tau_{\lambda(q)}^j \circ \lambda \circ \tau_q^i|^\Omega)^{-1}$ is a embedding of orbifold charts.
\end{proof}
\noindent
Instead of defining a Riemannian orbifold metric as in Definition \ref{defn: rieofdm}, Proposition \ref{prop: ext:rom} yields an equivalent definition of a Riemannian orbifold metric: It may be defined as a family of Riemannian metrics on the class of all compatible (with respect to the orbifold structure) orbifold charts, which satisfies the compatibility condition (cf. \cite[p.41]{follie2003}). From this point of view, a Riemannian orbifold metric on any representative of $\uU$ induces a \emph{uniquely determined} \ind{Riemannian orbifold metric!on the class $\uU$}{Riemannian orbifold metric on the equivalence class $\uU$}. We shall adopt this point of view in Lemma \ref{lem: rpbm:odi} below.\\
Either way, a Riemannian orbifold metric was defined using embeddings of orbifold charts. The reader may have noticed that our working definition of orbifolds (cf.\ Definition \ref{defn: haef:ofdII}) uses change of charts (but is equivalent to the approach using embeddings of orbifold charts). The definitions in this chapter are slightly easier to formulate using open embeddings of orbifold charts and therefore we chose this approach. Nevertheless, changes of orbifold charts are Riemannian isometries:

\begin{lem}\label{lem: ch:riso} \no{lem: ch:riso}
 Let $(Q,\uU, \rho)$ be a Riemannian orbifold and consider for some $(U,H,\phi), (V,G,\psi) \in \uU$ a change of charts $\lambda \colon U \supseteq \dom \lambda \rightarrow \cod \lambda \subseteq V$ . Furthermore, let $\rho_{\dom \lambda}$ be the pullback metric of $\rho_U$ with respect to the inclusion $\dom \lambda \subseteq U$. Then $\lambda \colon (\dom \lambda , \rho_{\dom \lambda}) \rightarrow (V,\rho_V)$ is a Riemannian embedding.
\end{lem}
 
\begin{proof}
 Let $p \in \dom \lambda$ be arbitrary and choose an open connected $H$-stable subset $p \in S \subseteq \dom \lambda$. Then $(S,H_S, \phi|_S)$ is an orbifold chart and $\lambda|_S$ is an embedding of orbifold charts. Since $\rho_{U}$ and $\rho_V$ are members of $\rho$, the map $\lambda|_{S}$ is a Riemannian embedding. In particular, $(\rho_{\dom \lambda})_p = (\lambda^* \rho_V)_p$ holds. Since $p \in \dom \lambda$ was arbitrary, $\lambda$ is a Riemannian embedding. 
\end{proof}

\begin{defn}\label{defn: isom}
 Let $(Q_i,\uU_i,\rho_i), i =1,2$ be Riemannian orbifolds and consider a map of orbifolds $[\hat{f}] \in \ORBM$. The map $[\hat{f}]$ is called \ind{orbifold map!isometric}{orbifold isometric}, if there is a representative $\hat{f} = (f,\set{f_i}_{i \in I}, P,\nu) \in \Orb{\vV,\wW}$ such that each lift $f_i \colon V_i \rightarrow W_{\alpha (i)}$ is an isometric immersion of the Riemannian manifold $(V_i,\rho_{1,i})$ to the Riemannian manifold $(W_{\alpha (i)}, \rho_{2,\alpha (i)})$.\\
 If $[\hat{f}]$ is a diffeomorphism of orbifolds which is orbifold isometric, $[\hat{f}]$ is called an \ind{orbifold map!orbifold isometry}{orbifold isometry}.
\end{defn}

\begin{rem}\label{rem: riofd:iso} \no{rem: riofd:iso}
 The condition to be an isometric immersion of Riemannian manifolds may be checked locally. Lemma \ref{lem: ch:riso} (i.e.\ the compatibility conditions of Riemannian orbifold metrics) combined with Proposition \ref{prop: ext:rom} that a map $[\hat{f}]$ will be orbifold isometric if and only if each representative $\hat{f} \coloneq (f,\set{f_j}_{j \in J}, [P,\nu])$ shares this property that the family of lifts $\set{f_j}_{j \in J}$ consists of isometric immersions.\\
 As an obvious first example, we mention that for a Riemannian orbifold $(Q,\uU,\rho)$ the identity morphism $\ido{}$ is an orbifold isometry. 
\end{rem}


\begin{lem}\label{lem: rpbm:odi} \no{lem: rpbm:odi}
 Let $(Q,\uU,\rho)$ be a Riemannian orbifold and $(Q_1,\uU_1)$ be an orbifold together with an orbifold diffeomorphism $[\hat{f}] \in \ORBM[(Q_1,\uU_1),(Q,\uU)]$. There exists a unique Riemannian orbifold metric $[\hat{f}]^* \rho$ on $(Q_1,\uU_1)$ such that $[\hat{f}]$ becomes an orbifold isometry with respect to $(Q_1,\uU_1, [\hat{f}]^* \rho)$ and $(Q,\uU,\rho)$. The Riemannian orbifold metric $[\hat{f}]^* \rho$ is called \ind{Riemannian orbifold metric!pullback}{pullback metric} induced by $[\hat{f}]$.
\end{lem}

\begin{proof}
 Following Corollary \ref{cor: diff:char} (d), we choose orbifold atlases $\vV = \setm{(V_i,G_i,\psi_i)}{i \in I} \in \uU_1$ and $\wW = \setm{(W_j,H_j,\varphi_j)}{j \in J}  \in \uU$ such that there is a representative $\hat{g} = (f, \set{f_i}_{i \in I}, [P,\nu])$ of $[\hat{f}]$ with the following properties: 
	\begin{compactenum}
	 \item $f_i \colon V_i \rightarrow W_{\beta (i)}$ is a diffeomorphism for each $i \in I$,
	 \item the map $\beta \colon I \rightarrow J$ is bijective,
	 \item $P = \Ch_\vV$ holds and for $\lambda \in \CH{V_i}{V_j}$, one has $\nu (\lambda ) = f_j \lambda f_i^{-1}|_{f_i (\dom \lambda)}$ (see Corollary \ref{cor: ll:diff}).
	\end{compactenum}
 Proposition \ref{prop: ext:rom} yields a unique family of compatible Riemannian metrics $(\rho_j)_{j\in J}$ induced by $\rho$ such that each chart $(W_j,H_j,\varphi_j)$ turns into a Riemannian manifold $(W_j, \rho_j)$. Endow each manifold $V_i$ with the pullback metric $f_i^* \rho_{\beta (i)}$, turning $f_i$ into a Riemannian isometry. \\
 \textbf{Claim:} The family $(f_i^* \rho_{\beta (i)})_{i \in I}$ turns each $\lambda \in \CH{V_i}{V_j}$, $i,j\in I$ into a Riemannian embedding.\\ An argument analogous to the proof of Lemma \ref{lem: add:cor} (c) shows that $\mu \coloneq f_j \lambda f^{-1}_{i}|_{f_i (\dom \lambda)} \in \CH{W_{\beta (i)}}{W_{\beta (j)}}$ and $f_j \lambda = \mu f_i|_{\dom \lambda}$ holds. Consider $p \in \dom \lambda$ and compute for $v,w\in T_pV_i$:
	\begin{align*}
	 (f_j^* \rho_{\beta (j)})_{\lambda (p)} (T_p \lambda (v), T_p \lambda (w)) 	&= (\rho_{\beta (j)})_{f_j \lambda (p)} (T_pf_j\lambda (v), T_pf_j\lambda (w)) \\
											&= (\rho_{\beta (j)})_{\mu f_i (p)} (T_p\mu f_i(v), T_p \mu f_i(w))\label{eq: compat}\\
											&= (f_i^* \rho_{\beta (i)})_p (v,w).
	\end{align*}
 The last identity is due to the compatibility condition of $\rho$, since $\mu$ is a change of orbifold charts (cf.\ Lemma \ref{lem: ch:riso}). In view of Proposition \ref{prop: ext:rom}, the compatible family $(f_i^*\rho_{\beta (i)})_{i \in I}$ yields a unique Riemannian orbifold metric $[\hat{f}]^*\rho$.\\ 
 We have to assure that  $[\hat{f}]^*\rho$ does not depend on the choice of $\hat{g}$. To this end, consider another representative $\hat{h} = (f, \set{h_k}_{k\in K}, [\Ch_{\vV'}, \nu']) \in \Orb{\vV',\wW'}$ of $[\hat{f}]$ with the same properties as $\hat{g}$. Write $([\hat{f}]^*\rho)'$ for the Riemannian orbifold metric induced by $\hat{h}$. Reviewing Proposition \ref{prop: ext:rom}, both metrics will coincide if the family $(f_i^*\rho_{\beta (i)})_{i \in I} \coprod (h_j^* \rho_{\beta' (j)})_{j\in J}$ of Riemannian metrics is compatible in the sense of Lemma \ref{lem: ch:riso}. To check this choose $i\in I$, $j \in J$ and some change of charts $\lambda \in \CH{V_i}{V_j'}$. Then $h_j \lambda f_i^{-1}|_{f_i (\dom \lambda)}$ is a change of charts. An analogous computation as above together with the compatibility of the metrics $\rho_{\beta (i)}$ and $\rho_{\beta' (j)}$ yields that $\lambda$ is a Riemannian embedding. Thus $[\hat{f}]^*\rho$ and $([\hat{f}]^*\rho)'$ coincide, proving the uniqueness of the pullback orbifold metric.
\end{proof}

\begin{rem}
 In Lemma \ref{lem: rpbm:odi} special representatives of an orbifold diffeomorphism were used in the construction. Their lifts were given by a family of diffeomorphisms. The proof of Lemma \ref{lem: rpbm:odi} may be adapted to work with an arbitrary family of lifts of the orbifold diffeomorphism. In general, these families will be families of local diffeomorphisms by Corollary \ref{cor: diff:char}. In this case, the identities computed in the proof will only hold locally. Hence the same arguments require cumbersome notation, which may be avoided in the construction if representatives are used whose lifts are diffeomorphisms. 
\end{rem}

Our goal in introducing Riemannian orbifold metrics on orbifolds is to obtain an analogue of the Riemannian exponential map on a manifold for a Riemannian orbifold. To this end, we need to introduce the notion of a geodesic on a Riemannian orbifold. 

\subsection{Geodesics on orbifolds}\label{Sect: Geod}

In this section let $(Q,\uU, \rho)$ be a Riemannian orbifold. Notice that by Proposition \ref{prop: ext:rom}, the Riemannian orbifold metric $\rho$ induces a family of compatible Riemannian metrics for each representative of $\uU$. As we introduced Riemannian orbifold metrics, the question arises how geodesics for a Riemannian orbifold may be defined. Furthermore, one would like these geodesics to share at least some properties of geodesics on a Riemannian manifold. Geodesics on Riemannian orbifolds have been considered in the literature (cf.\ Haefliger and collaborators \cite{cgo2006,msnpc1999}, Chen et al. \cite{chen2001}) in the context of different frameworks (i.e. \'{e}tale groupoids, respectively Chen-Ruan good maps). For the setting considered in this work, we shall give a definition of an orbifold godesic which shares the properties developed for geodesics on Orbifolds in the literature. In fact, the restriction of a geodesic to a compact interval corresponds to a unique $\gG$-geodesic in the sense of Haefliger. However, since geodesics should be maps of orbifolds, our proofs are independent of this equivalence. \\
Throughout this section, $\iI \coloneq ]a,b[ \subseteq \RR$ will always be an open interval with $a<b$. Endow $\iI$ with the canonical structure of an open submanifold of $\RR$ (i.e.\ a trivial orbifold structure) and denote its orbifold structure by $\uU_\iI$.
As a first step, we define smooth paths in orbifolds:

\begin{defn}
 An orbifold map $[\hat{c}] \in  \ORB (\iI ,(Q,\uU))$ is called a \ind{orbifold map!smooth path}{smooth orbifold path}. 
\end{defn}

\begin{ex}
	\begin{compactenum}
	 \item If $(Q,\uU)$ is a trivial orbifold (i.e.\ a manifold), a smooth orbifold path is just a smooth curve $\iI \rightarrow Q$. 
	 \item Reconsider Example \ref{ex: Z2hs}: The map $\gamma \colon \RR^2 \rightarrow \RR^2, (x,y) \mapsto (-x,y)$ is a reflection of $\RR^2$ and $H$ is the right half plane. Let $q \colon \RR^2 \rightarrow H$ be the quotient map to the orbit space with respect to the $\langle \gamma \rangle$-action. Then $H$ is an orbifold with global chart $(\RR^2,\langle \gamma \rangle , q)$. As the orbifold atlas contains only one chart, the changes of charts are generated by $\gamma$ and $\id_{\RR^3}$. Define $I_1 \coloneq ]0,\tfrac{3}{4}[$ and $I_2 \coloneq ]\tfrac{1}{4},1[$ which cover $]0,1[ = I_1 \cup I_2$. Let $\lambda \colon I_1 \supseteq I_1 \cap I_2 \rightarrow I_2$ be the inclusion. Then the quasi-pseudogroup $P\coloneq \set{\id_{I_1}, \id_{I_2}, \lambda , \lambda^{-1}}$ generates the change of charts of $\set{I_1, I_2}$.\\ Consider the smooth maps $c_1 \colon I_1 \rightarrow \RR^2, t \mapsto (1-2t,1-2t)$ and $c_2 \colon I_2 \rightarrow \RR^2, t \mapsto (2t-1,1-2t)$. We obtain a continuous map $c \colon ]0,1[ \rightarrow H, t \mapsto q \circ c_i (t)$, for $t \in I_i$. Set $\nu (\lambda ) \coloneq \gamma$, to uniquely determine $\nu \colon P \rightarrow \Psi (\uU)$, which satisfies (R4) of Definition \ref{defn: rep:ofdm}. Then $\hat{c} \coloneq (c, \set{c_1,c_2} , P , \nu)$ is a smooth path in $H$. We sketch the images of the lifts and the smooth path in $H$: 
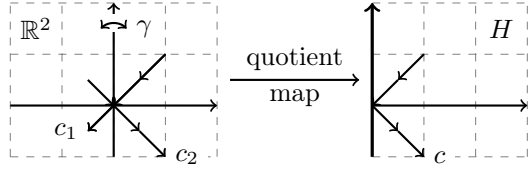
	  \begin{figure}[h!]\centering
	\begin{tikzpicture}[scale=0.68]
	  \draw[style=help lines,dashed,gray] (-2,-1) grid (2,2);
	\coordinate (y) at (0,2);
	\coordinate (x) at (2,0);
	\draw[thick,->] (0,-1) -- (y); 
	\draw[thick,->] (-2,0) -- (x); 
	\draw(0,1.6) node[shape=rectangle,fill=white]{};
	\draw(-1.5,1.5) node[fill=white]{$\RR^2$};
	\draw[thick,<->] (-0.25,1.5) to [bend left =45](0.25,1.5) node[right,fill=white]{$\gamma$};
	\draw[thick, ->] (1,1) -- (0,0); \draw[thick, ->] (0,0) -- (-0.5,-0.5) node[left]{$c_1$};
        \draw[thick, ->] (1,1) -- (0.5,0.5); \draw[thick, ->] (0.5,0.5) -- (0,0); \draw[thick, ->] (0,0) -- (-0.5,-0.5) node[left,fill=white]{$c_1$};
        \draw[thick, ->] (-0.5,0.5) -- (0,0); \draw[thick, ->] (0,0) -- (0.5,-0.5); \draw[thick, ->] (0.5,-0.5) -- (1,-1) node[right,fill=white]{$c_2$};
	\draw[thick,->] (2.25,0.5) -- node[midway,above]{\text{quotient}} node[midway,below]{\text{map}}(4.75,0.5);
	\draw[style=help lines,gray,dashed] (5,-1) grid (8,2);
	\draw(7.5,1.5) node[fill=white]{$H$};
	\draw[very thick,->] (5,-1) -- (5,2); 
	\draw[thick, ->] (5,0) -- (8,0);
	\draw[thick, ->] (6,1) -- (5.5,0.5); \draw[thick, ->] (5.5,0.5) -- (5,0); \draw[thick, ->] (5,0) -- (5.5,-0.5);  \draw[thick, ->] (5.5,-0.5)--(6,-1) node[right,fill=white]{$c$};
      \end{tikzpicture}
 \end{figure}   \\
 Notice that there is the weaker notion of a continuous path. It was introduced in \cite[Chapter III, 3]{msnpc1999} to obtain a fundamental group of an \'{e}tale groupoid. The map $\hat{c}$ induces a continuous path into $H$ in the sense of Haefliger (cf.\ \cite[III.\ Example  3.3 (2)]{msnpc1999}). Define a map $\nu' \colon P \rightarrow \Psi (\uU)$ via $\nu' (\lambda) = \id_{\RR^2}$. The tuple $(c,\set{c_1,c_2}, P ,\nu')$ does \emph{not} define a charted orbifold map, but it induces a continuous path in the sense of Haefliger (cf.\ \cite[III.\ Example 3.3 (2)]{msnpc1999}).
\end{compactenum} 
\end{ex}
\noindent
 In the last example, an orbifold path has been constructed with respect to a special orbifold atlas: Define the set of all orbifold charts $\aA_\iI = \setm{(V_\alpha, \set{\id_{V_\alpha}} , \pi_\alpha)}{\alpha \in A} \in \uU_\iI$ such that an orbifold chart $(V_\alpha , \set{\id_{V_\alpha}} , \pi_\alpha) \in \uU_\iI$ is contained in $\aA_\iI$ if and only if:
 $V_\alpha =\, ]l(\alpha), r(\alpha)[ \subseteq \iI$ is an open interval with $a \leq l(\alpha) < r(\alpha) \leq b$ and the map $\pi_\alpha \colon ]l(\alpha) , r(\alpha)[ \rightarrow \iI$ is the inclusion (of sets). By construction each change of orbifold charts in $\Ch_{V_\alpha, V_\beta}$ for two orbifold charts $(V_\alpha, \set{\id_{V_\alpha}} , \pi_\alpha), (V_\beta, \set{\id_{V_\beta}} , \pi_\beta) \in \aA_\iI$ is an inclusion of open sets.\\ 
 Consider a smooth orbifold path $[\hat{c}] \in \ORBM[\iI, (Q,\uU)]$ with representative $\hat{c} = (c,\set{c_k}_{k \in I}, [P,\nu])$ whose lifts are defined on charts $(\dom c_k , \set{\id_{\dom c_k}}, \pi_k)$. The chart maps of orbifold charts on $\iI$ are diffeomorphisms, since they are also manifold charts of the smooth manifold $\iI$. Define an orbifold atlas $\vV_{\hat{c}} \coloneq \setm{\pi_k (\dom c_k)}{k\in I}$ of $\iI$, where $\pi_k (\dom c_k) \subseteq \iI$ is a connected open interval. Hence $\vV_{\hat{c}} \subseteq \aA_\iI$ holds. Apply Lemma \ref{lemdef: ind:ofdm} together with this set of charts to obtain a representative $\hat{h} \in \Orb{\vV_{\hat{c}}, \wW}$ of $[\hat{c}]$, where $\wW$ is the range atlas of $\hat{c}$. In conclusion, for each smooth orbifold path, there is a representative whose domain atlas is contained in $\aA_\iI$. 

\begin{lem}\label{lem: osmpath:prop} \no{lem: osmpath:prop}
 Let $[\hat{c}] \in \ORBM[\iI,(Q,\uU)]$ be a smooth orbifold path and $P$ be some point in $\iI$. Identifying the tangent orbifold $\tT \iI$ with the tangent manifold $\iI \times \RR$, the element $\tT c (P,1) \in \tT_{c(P)} Q$ is called the \ind{orbifold map!smooth path!initial vector}{initial vector} of $[\hat{c}]$ at $P$. For each representative $\hat{c} = (c,\set{c_k}_{k \in I}, [P,\nu]) \in \Orb{\vV , \wW}$ of $[\hat{c}]$ with $\vV \subseteq \aA_\iI$ and $P \in \dom c_k$, the initial vector is induced by $T_P c_k (1) = c_k' (a)$.
\end{lem}

\begin{proof}
 Consider the lift $c_k \colon \dom c_k \rightarrow V_k$, where $(\dom c_k , \set{\id_{\dom c_k}}, \pi_k) \in \aA_\iI$ and $(V_k, G_k, \psi_k) \in \uU$. As $\iI$ is a trivial orbifold, the tangent manifold $T\iI \cong \iI \times \RR$ coincides with the tangent orbifold. We suppress the identification $T \id_\iI$ in the formulas: By Definition \ref{defn: tofd:map} of the tangent orbifold map, $\tT c (P,1) = T \psi_k Tc_k T(\pi_k)^{-1} (P,1)$ is well-defined. Hence it suffices to prove $T \pi_k^{-1} (P,1) = (P,1) \in \dom Tc_k \cong \dom c_k \times \RR$. As $(\dom c_k , \set{\id_{\dom c_k}}, \pi_k) \in \aA_\iI$, so $\pi_k$ is the inclusion of sets $\dom c_k \hookrightarrow \iI$, $\pi_k$ is the restriction of a linear continuous map. A computation in the identification proves $T\pi_k^{-1} (P,1) = (P,1)$, whence from $T_P c_k (1) = Tc_k (P,1)$ the assertion follows. 
\end{proof}


\begin{lem}\label{lem: nice:geod} \no{lem: nice:geod}
 Let $[\hat{c}] \in \ORBM[\iI , (Q,\uU)]$ be an orbifold path and $[a,b] \subseteq \iI$ some compact subset. There exists $\hat{g} = (c|_{]x,y[}, \set{g_k}_{1 \leq k \leq N}, [P_g,\nu_g])$ with $x<a<b<y$ and $N \in \NN$ such that 
    \begin{compactitem}
      \item[1.] $[\hat{c}]|_{]x,y[} = [\hat{g}]$,
      \item[2.] $\dom g_k = ]l(k), r(k)[$ for each $1\leq k\leq N$ such that $$x= l(1) < l(2) < r(1) < l(3) < r(2) < \cdots < l(N) < r(N-1) < r(N) = y$$ 
      \item[3.] $P_g = \set{\id_{]l(N),r(N)[}} \cup \setm{\id_{]l(k),r(k)[} , \iota_k^{k+1}, (\iota_k^{k+1})^{-1}}{1\leq k \leq N-1}$, where $\iota_k^{k+1}$ is the canonical inclusion $]l(k+1),r(k)[ \hookrightarrow ]l(k+1),r(k+1)[$. 
    \end{compactitem} 
\end{lem}

\begin{proof} Construct a refinement of the domain atlas of $\hat{c}$. A full proof is given in Appendix \ref{App: Haefgeodesics}.\end{proof}

In a neighborhood of a compact set, we may think of an orbifold path as a family of smooth paths, which are compatible in the following way: On each intersection of their domains, the inclusion of sets induces a change of orbifold charts in the range atlas which maps one lift to the other. The situation is sketched in the following figure for a smooth path in an orbifold $(Q,\uU)$: 
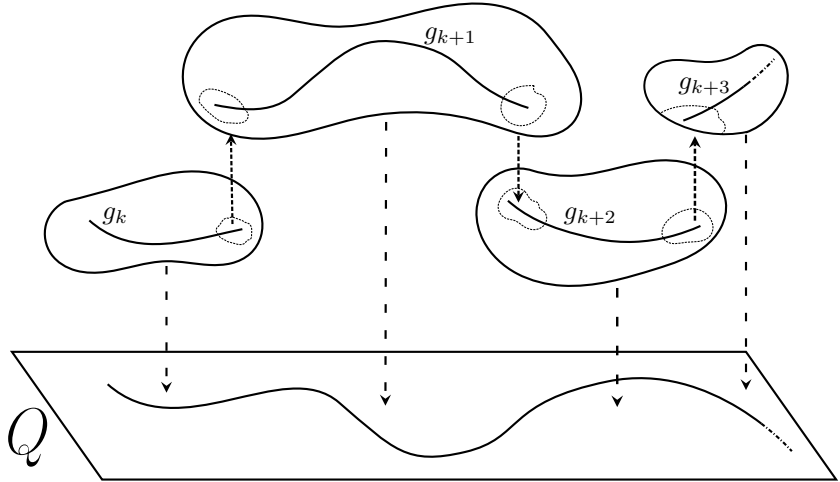
\begin{figure}[htb]
  \centering
  \begin{tikzpicture}[scale=0.85, y=0.80pt, x=0.8pt,yscale=-1, inner sep=0pt, outer sep=0pt, >=stealth]
\begin{scope}[shift={(0,-768.89762)}]
  \path[shift={(0,768.89762)},draw=black,line join=miter,line cap=butt,line
    width=0.800pt] (21.2598,198.4252) -- (70.8661,269.2913) -- (474.8031,269.2913)
    -- (425.1968,198.4252) -- (21.2598,198.4252) -- cycle;
  \path[fill=black] (14,1031.1023) node[above right] (text4334) {\Huge \textit{Q}};
  \path[draw=black,line join=miter,line cap=butt,line width=0.800pt]
    (50.6883,883.9211) .. controls (45.4913,886.3646) and (41.4684,891.1804) ..
    (39.9887,896.7292) .. controls (38.5090,902.2781) and (39.5999,908.4575) ..
    (42.8901,913.1643) .. controls (46.0866,917.7370) and (51.1556,920.7956) ..
    (56.5402,922.2559) .. controls (61.9249,923.7163) and (67.6217,923.6802) ..
    (73.1530,922.9509) .. controls (84.2155,921.4924) and (94.9002,917.3201) ..
    (106.0555,917.0634) .. controls (114.0580,916.8793) and (121.9681,918.7242) ..
    (129.9626,919.1243) .. controls (133.9599,919.3243) and (138.0041,919.1576) ..
    (141.8920,918.2074) .. controls (145.7798,917.2573) and (149.5178,915.4932) ..
    (152.4548,912.7744) .. controls (155.2487,910.1881) and (157.2524,906.7820) ..
    (158.2712,903.1137) .. controls (159.2901,899.4453) and (159.3345,895.5258) ..
    (158.5164,891.8076) .. controls (156.8802,884.3710) and (151.8221,877.9309) ..
    (145.4364,873.7834) .. controls (136.2964,867.8470) and (124.8169,866.4494) ..
    (113.9776,867.5855) .. controls (103.1383,868.7216) and (92.7035,872.1897) ..
    (82.2710,875.3430) .. controls (71.8272,878.4998) and (61.2940,881.3606) ..
    (50.6883,883.9211);
  \path[draw=black,line join=miter,line cap=butt,line width=0.800pt]
    (135.6887,784.1041) .. controls (127.3884,787.6777) and (120.5625,794.6384) ..
    (117.3543,803.0866) .. controls (114.1461,811.5348) and (114.6601,821.3652) ..
    (118.9226,829.3337) .. controls (122.9277,836.8211) and (129.9921,842.4211) ..
    (137.8829,845.5576) .. controls (145.7738,848.6941) and (154.4501,849.5002) ..
    (162.9202,848.9001) .. controls (179.8604,847.7000) and (195.8993,841.0992) ..
    (212.5010,837.5218) .. controls (236.1773,832.4198) and (261.1965,833.5883) ..
    (284.2444,841.0310) .. controls (292.0641,843.5561) and (299.7621,846.8109) ..
    (307.9392,847.6228) .. controls (312.0278,848.0288) and (316.2158,847.8024) ..
    (320.1372,846.5762) .. controls (324.0587,845.3501) and (327.7061,843.0846) ..
    (330.2538,839.8612) .. controls (332.7658,836.6830) and (334.1268,832.6751) ..
    (334.4131,828.6342) .. controls (334.6994,824.5932) and (333.9479,820.5199) ..
    (332.5646,816.7123) .. controls (329.7980,809.0970) and (324.6306,802.6406) ..
    (319.7262,796.1913) .. controls (317.2659,792.9561) and (314.8380,789.6710) ..
    (311.9280,786.8335) .. controls (304.8451,779.9268) and (295.1958,776.0379) ..
    (285.3944,774.6956) .. controls (275.5929,773.3534) and (265.6018,774.4085) ..
    (255.9015,776.3514) .. controls (236.5011,780.2371) and (217.4514,787.6951) ..
    (197.6844,786.8335) .. controls (187.2113,786.3770) and (176.9983,783.5846) ..
    (166.6417,781.9617) .. controls (156.2850,780.3389) and (145.3172,779.9587) ..
    (135.6887,784.1041);
  \path[draw=black,line join=miter,line cap=butt,line width=0.800pt]
    (307.6390,866.7650) .. controls (303.3765,865.3472) and (298.6963,865.2125) ..
    (294.3594,866.3828) .. controls (290.0224,867.5531) and (286.0453,870.0239) ..
    (283.0747,873.3935) .. controls (279.2738,877.7048) and (277.1869,883.4064) ..
    (276.9073,889.1470) .. controls (276.6277,894.8877) and (278.0972,900.6497) ..
    (280.7352,905.7560) .. controls (284.7438,913.5154) and (291.4033,919.7595) ..
    (299.1000,923.8870) .. controls (306.7967,928.0146) and (315.4938,930.0933) ..
    (324.2123,930.6063) .. controls (341.6495,931.6323) and (358.8144,926.5728) ..
    (375.4833,921.3524) .. controls (384.7665,918.4451) and (394.2400,915.3624) ..
    (401.9312,909.4063) .. controls (405.7768,906.4283) and (409.1222,902.7391) ..
    (411.3612,898.4212) .. controls (413.6002,894.1033) and (414.6940,889.1365) ..
    (414.0844,884.3110) .. controls (413.5944,880.4320) and (412.0160,876.7177) ..
    (409.6992,873.5682) .. controls (407.3825,870.4187) and (404.3415,867.8277) ..
    (400.9467,865.8881) .. controls (394.1572,862.0089) and (386.1071,860.7814) ..
    (378.2946,861.1160) .. controls (362.6697,861.7850) and (347.7997,868.3380) ..
    (332.2033,869.4944) .. controls (323.9405,870.1071) and (315.5658,869.1766) ..
    (307.6390,866.7650);
  \path[draw=black,line join=miter,line cap=butt,line width=0.744pt]
    (425.8479,799.1387) .. controls (419.8917,803.7074) and (412.4654,806.3279) ..
    (404.9609,806.5088) .. controls (398.9694,806.6532) and (393.0208,805.2913) ..
    (387.0284,805.3896) .. controls (384.0322,805.4387) and (381.0115,805.8632) ..
    (378.2333,806.9861) .. controls (375.4551,808.1091) and (372.9239,809.9674) ..
    (371.3038,812.4883) .. controls (369.7982,814.8312) and (369.1344,817.6722) ..
    (369.2655,820.4541) .. controls (369.3966,823.2360) and (370.2999,825.9557) ..
    (371.7023,828.3618) .. controls (374.5072,833.1740) and (379.1731,836.6301) ..
    (384.0752,839.2748) .. controls (396.7486,846.1119) and (411.8792,848.2720) ..
    (425.9595,845.2543);
  \path[draw=black,line join=miter,line cap=butt,miter limit=4.00,line
    width=0.812pt] (425.2019,845.5792) .. controls (434.9653,842.0389) and
    (442.9778,833.9667) .. (446.4453,824.1772) .. controls (448.2926,818.9621) and
    (448.8706,813.1370) .. (447.0790,807.9025) .. controls (446.1832,805.2853) and
    (444.6992,802.8550) .. (442.6863,800.9575) .. controls (440.6734,799.0600) and
    (438.1248,797.7102) .. (435.3987,797.2402) .. controls (432.0209,796.6579) and
    (428.4238,797.4729) .. (425.6268,799.4542);
  \path[draw=black,dash pattern=on 1.60pt off 0.80pt on 0.40pt off 0.80pt,line
    join=miter,line cap=butt,miter limit=4.00,line width=0.800pt]
    (433.5865,1007.1100) .. controls (437.1511,1009.8627) and (440.2092,1012.4869)
    .. (442.7024,1014.7623) .. controls (447.6888,1019.3132) and
    (450.4155,1022.4691) .. (450.4155,1022.4691);
  \path[draw=black,line join=miter,line cap=butt,line width=0.800pt]
    (73.9247,985.3048) .. controls (113.9167,1021.2573) and (170.8751,969.1464) ..
    (201.1721,995.4038) .. controls (231.4691,1021.6612) and (238.3365,1030.5483)
    .. (270.6533,1022.8731) .. controls (302.9701,1015.1979) and
    (292.0631,995.4038) .. (349.4255,983.6890) .. controls (378.1066,977.8315) and
    (403.3541,987.5266) .. (421.4314,998.6860) .. controls (425.9507,1001.4758)
    and (430.0218,1004.3572) .. (433.5865,1007.1100);
  \path[draw=black,dash pattern=on 1.60pt off 0.80pt on 0.40pt off 0.80pt,line
    join=miter,line cap=butt,miter limit=4.00,line width=0.800pt]
    (428.8512,816.1751) .. controls (430.3729,814.8906) and (431.7607,813.6657) ..
    (433.0073,812.5279) .. controls (435.5005,810.2525) and (437.4288,808.3258) ..
    (438.7338,806.9680) .. controls (440.0388,805.6101) and (440.7204,804.8211) ..
    (440.7204,804.8211);
  \path[draw=black,dash pattern=on 3.20pt off 6.40pt,line join=miter,line
    cap=butt,miter limit=4.00,line width=0.800pt] (423.8914,820.1803) .. controls
    (425.6737,818.8039) and (427.3294,817.4597) .. (428.8512,816.1751);
  \path[draw=black,line join=miter,line cap=butt,line width=0.800pt]
    (390.7042,839.0211) .. controls (391.4093,838.7558) and (392.1087,838.4844) ..
    (392.8023,838.2073) .. controls (394.1894,837.6530) and (395.5531,837.0758) ..
    (396.8924,836.4791) .. controls (402.2496,834.0925) and (407.2170,831.3941) ..
    (411.7363,828.6043) .. controls (416.2556,825.8144) and (420.3268,822.9331) ..
    (423.8914,820.1803);
  \path[draw=black,line join=miter,line cap=butt,line width=0.800pt]
    (132.6914,830.2738) .. controls (138.2112,831.4131) and (143.5699,832.3031) ..
    (148.6825,832.6187) .. controls (158.9078,833.2498) and (168.1484,831.5835) ..
    (175.7226,825.0191) .. controls (206.0196,798.7618) and (212.8870,789.8746) ..
    (245.2038,797.5499) .. controls (253.2830,799.4687) and (258.6607,802.1449) ..
    (263.0790,805.2630) .. controls (267.4973,808.3810) and (270.9562,811.9409) ..
    (275.1978,815.6271) .. controls (279.4394,819.3132) and (284.4636,823.1256) ..
    (292.0126,826.7486) .. controls (295.7871,828.5601) and (300.1928,830.3243) ..
    (305.4475,832.0017);
  \path[draw=black,line join=miter,line cap=butt,line width=0.800pt]
    (63.8257,894.4138) .. controls (68.8247,898.9078) and (74.0888,902.0259) ..
    (79.5328,904.0930) .. controls (84.9768,906.1602) and (90.6007,907.1764) ..
    (96.3192,907.4667) .. controls (107.7564,908.0474) and (119.5722,905.7247) ..
    (131.0850,903.0989) .. controls (136.8415,901.7861) and (142.5222,900.3974) ..
    (148.0419,899.2582);
  \path[draw=black,line join=miter,line cap=butt,line width=0.800pt]
    (294.1839,883.4924) .. controls (298.4255,887.1785) and (303.4498,890.9909) ..
    (310.9988,894.6139) .. controls (318.5478,898.2369) and (328.6215,901.6706) ..
    (342.9621,904.5993) .. controls (350.1324,906.0637) and (357.0881,906.5560) ..
    (363.7708,906.2964) .. controls (370.4535,906.0368) and (376.8632,905.0254) ..
    (382.9415,903.4821) .. controls (389.0199,901.9388) and (394.7668,899.8638) ..
    (400.1240,897.4771);
  \path[shift={(0,768.89762)},draw=black,dash pattern=on 0.80pt off 0.40pt,line
    join=miter,line cap=butt,miter limit=4.00,line width=0.400pt]
    (142.1940,123.4964) -- (140.5781,124.3043) .. controls (139.7603,124.1046) and
    (138.8632,124.2542) .. (138.1544,124.7082) .. controls (137.6705,125.0182) and
    (137.2812,125.4560) .. (136.9425,125.9201) .. controls (136.2584,126.8576) and
    (135.7478,127.9318) .. (134.9227,128.7479) .. controls (134.6516,129.0160) and
    (134.3482,129.2543) .. (134.1148,129.5558) .. controls (133.8043,129.9568) and
    (133.6314,130.4625) .. (133.6314,130.9696) .. controls (133.6314,131.4768) and
    (133.8043,131.9824) .. (134.1148,132.3835) -- (135.7306,134.8073) .. controls
    (136.1192,135.0987) and (136.5240,135.3686) .. (136.9425,135.6152) .. controls
    (137.4617,135.9211) and (138.0020,136.1913) .. (138.5583,136.4231) --
    (139.7702,137.2310) .. controls (141.2492,137.6429) and (142.7304,138.0468) ..
    (144.2138,138.4429) .. controls (145.1019,138.6800) and (146.0025,138.9160) ..
    (146.9218,138.9195) .. controls (147.3814,138.9215) and (147.8433,138.8639) ..
    (148.2806,138.7224) .. controls (148.7179,138.5810) and (149.1304,138.3538) ..
    (149.4652,138.0389) .. controls (149.7450,137.7758) and (149.9656,137.4575) ..
    (150.2095,137.1608) .. controls (150.4534,136.8641) and (150.7310,136.5810) ..
    (151.0811,136.4231) .. controls (151.2757,136.3353) and (151.4861,136.2893) ..
    (151.6934,136.2380) .. controls (151.9007,136.1867) and (152.1094,136.1283) ..
    (152.2929,136.0191) .. controls (152.4848,135.9050) and (152.6427,135.7382) ..
    (152.7601,135.5484) .. controls (152.8776,135.3586) and (152.9557,135.1461) ..
    (153.0055,134.9285) .. controls (153.1050,134.4934) and (153.0926,134.0417) ..
    (153.1009,133.5954) .. controls (153.1191,132.6207) and (153.2386,131.6470) ..
    (153.1890,130.6734) .. controls (153.1394,129.6998) and (152.9027,128.7005) ..
    (152.2930,127.9399) .. controls (151.7101,127.2128) and (150.8007,126.7581) ..
    (149.8692,126.7280) -- (147.8494,125.5162) -- (142.1940,123.4964);
  \path[shift={(0,768.89762)},draw=black,dash pattern=on 0.80pt off 0.40pt,line
    join=miter,line cap=butt,miter limit=4.00,line width=0.400pt]
    (140.5781,55.6311) .. controls (138.0814,54.2189) and (135.3313,53.2564) ..
    (132.4989,52.8034) .. controls (131.2940,52.6106) and (130.0480,52.5112) ..
    (128.8633,52.8034) .. controls (127.7388,53.0807) and (126.7079,53.7118) ..
    (125.9314,54.5710) .. controls (125.1548,55.4303) and (124.6329,56.5125) ..
    (124.4197,57.6509) .. controls (124.2702,58.4496) and (124.2702,59.2760) ..
    (124.4197,60.0746) -- (127.5591,63.7795) -- (129.2672,65.3766) --
    (132.4989,68.1538) -- (135.3266,70.1736) -- (138.9623,71.3855) .. controls
    (139.6351,71.3563) and (140.3092,71.3563) .. (140.9821,71.3855) .. controls
    (142.2752,71.4416) and (143.5631,71.6055) .. (144.8574,71.6220) .. controls
    (146.1517,71.6384) and (147.4716,71.5006) .. (148.6573,70.9816) .. controls
    (149.4582,70.6309) and (150.1967,70.0962) .. (150.6771,69.3657) .. controls
    (151.1182,68.6950) and (151.3242,67.8728) .. (151.2515,67.0734) .. controls
    (151.1788,66.2739) and (150.8280,65.5024) .. (150.2732,64.9222) --
    (149.0613,63.3063) .. controls (147.9031,60.8607) and (146.0659,58.7410) ..
    (143.8098,57.2469) .. controls (142.8032,56.5803) and (141.7154,56.0364) ..
    (140.5781,55.6311);
  \path[draw=black,dash pattern=on 0.80pt off 0.40pt,line join=miter,line
    cap=butt,miter limit=4.00,line width=0.400pt] (308.6255,819.6812) .. controls
    (307.2808,819.5954) and (305.9306,819.5954) .. (304.5859,819.6812) .. controls
    (302.9476,819.7857) and (301.3113,820.0190) .. (299.7384,820.4891) .. controls
    (297.0255,821.2999) and (294.4984,822.8496) .. (292.7433,825.0715) .. controls
    (290.9882,827.2934) and (290.0564,830.2075) .. (290.4473,833.0119) .. controls
    (290.7930,835.4915) and (292.1622,837.7933) .. (294.1095,839.3669) .. controls
    (296.0567,840.9404) and (298.5534,841.7913) .. (301.0564,841.8455) .. controls
    (303.5594,841.8997) and (306.0586,841.1739) .. (308.2036,839.8830) .. controls
    (310.3487,838.5921) and (312.1426,836.7485) .. (313.4730,834.6277) .. controls
    (314.5020,832.9874) and (315.2745,831.1080) .. (315.2000,829.1731) .. controls
    (315.1627,828.2057) and (314.9082,827.2396) .. (314.4084,826.4105) .. controls
    (313.9085,825.5813) and (313.1574,824.8950) .. (312.2612,824.5287) .. controls
    (311.4712,824.2058) and (310.5856,824.1311) .. (309.8374,823.7208) .. controls
    (309.1543,823.3461) and (308.6278,822.6955) .. (308.4039,821.9492) .. controls
    (308.1800,821.2030) and (308.2614,820.3700) .. (308.6255,819.6812);
  \path[shift={(0,768.89762)},draw=black,dash pattern=on 0.80pt off 0.40pt,line
    join=miter,line cap=butt,miter limit=4.00,line width=0.400pt]
    (303.7780,109.7617) .. controls (302.8406,108.7232) and (301.5353,108.0402) ..
    (300.1595,107.7861) .. controls (298.7837,107.5319) and (297.3438,107.6964) ..
    (296.0314,108.1811) .. controls (293.4067,109.1506) and (291.3753,111.3406) ..
    (290.0433,113.8013) .. controls (289.4309,114.9328) and (288.9400,116.1549) ..
    (288.8315,117.4370) .. controls (288.6949,119.0503) and (289.1903,120.7003) ..
    (290.1523,122.0028) .. controls (291.1142,123.3052) and (292.5280,124.2569) ..
    (294.0829,124.7083) .. controls (295.2762,125.0546) and (296.5424,125.1158) ..
    (297.7186,125.5162) .. controls (298.7834,125.8786) and (299.7368,126.5065) ..
    (300.6401,127.1767) .. controls (301.5434,127.8470) and (302.4126,128.5680) ..
    (303.3740,129.1518) .. controls (304.8829,130.0680) and (306.6243,130.6321) ..
    (308.3896,130.6370) .. controls (310.1548,130.6420) and (311.9385,130.0664) ..
    (313.2854,128.9255) .. controls (314.6324,127.7845) and (315.5065,126.0653) ..
    (315.5059,124.3000) .. controls (315.5057,123.4174) and (315.2912,122.5336) ..
    (314.8715,121.7571) .. controls (314.4518,120.9806) and (313.8263,120.3141) ..
    (313.0691,119.8607) .. controls (312.6505,119.6102) and (312.1831,119.4158) ..
    (311.8572,119.0528) .. controls (311.6626,118.8360) and (311.5298,118.5696) ..
    (311.4373,118.2934) .. controls (311.3448,118.0172) and (311.2907,117.7298) ..
    (311.2382,117.4433) .. controls (311.1856,117.1567) and (311.1340,116.8691) ..
    (311.0461,116.5914) .. controls (310.9582,116.3136) and (310.8325,116.0443) ..
    (310.6453,115.8211) .. controls (310.3110,115.4225) and (309.8105,115.2029) ..
    (309.3117,115.0550) .. controls (308.8130,114.9070) and (308.2949,114.8162) ..
    (307.8176,114.6093) .. controls (306.8391,114.1850) and (306.1271,113.3165) ..
    (305.5520,112.4184) .. controls (304.9770,111.5202) and (304.4925,110.5534) ..
    (303.7780,109.7617);
  \path[shift={(0,768.89762)},draw=black,dash pattern=on 0.80pt off 0.40pt,line
    join=miter,line cap=butt,miter limit=4.00,line width=0.400pt]
    (395.4769,119.0528) .. controls (390.6421,119.4985) and (385.9857,121.6812) ..
    (382.5502,125.1122) .. controls (381.0471,126.6134) and (379.7470,128.3988) ..
    (379.2249,130.4580) .. controls (378.9639,131.4876) and (378.9048,132.5746) ..
    (379.1122,133.6163) .. controls (379.3197,134.6580) and (379.7996,135.6523) ..
    (380.5304,136.4231) .. controls (381.2342,137.1654) and (382.1534,137.6854) ..
    (383.1288,137.9936) .. controls (384.1042,138.3019) and (385.1358,138.4052) ..
    (386.1586,138.3890) .. controls (388.2042,138.3565) and (390.2113,137.8554) ..
    (392.2452,137.6350) .. controls (394.5133,137.3892) and (396.8116,137.4927) ..
    (399.0721,137.1848) .. controls (400.2023,137.0308) and (401.3242,136.7723) ..
    (402.3747,136.3277) .. controls (403.4251,135.8831) and (404.4045,135.2472) ..
    (405.1720,134.4033) .. controls (406.4448,133.0036) and (407.0707,131.0732) ..
    (407.0175,129.1820) .. controls (406.9643,127.2908) and (406.2620,125.4467) ..
    (405.1720,123.9003) .. controls (403.0035,120.8239) and (399.2391,118.9418) ..
    (395.4769,119.0528);
  \path[draw=black,dash pattern=on 0.80pt off 0.40pt,line join=miter,line
    cap=butt,miter limit=4.00,line width=0.400pt] (378.5106,834.6277) .. controls
    (379.9320,833.8322) and (381.4201,833.1558) .. (382.9542,832.6079) .. controls
    (390.0181,830.0851) and (397.8080,830.3729) .. (405.1720,831.8000) .. controls
    (406.2451,832.0079) and (407.3369,832.2482) .. (408.2709,832.8161) .. controls
    (408.7379,833.1000) and (409.1599,833.4655) .. (409.4747,833.9123) .. controls
    (409.7894,834.3591) and (409.9937,834.8897) .. (410.0195,835.4356) .. controls
    (410.0386,835.8397) and (409.9614,836.2472) .. (410.0195,836.6475) .. controls
    (410.0651,836.9617) and (410.1936,837.2604) .. (410.3672,837.5263) .. controls
    (410.5409,837.7921) and (410.7591,838.0264) .. (410.9933,838.2409) .. controls
    (411.4617,838.6697) and (411.9971,839.0233) .. (412.4433,839.4752) .. controls
    (413.4772,840.5224) and (413.9530,842.0915) .. (413.6748,843.5365) .. controls
    (413.3967,844.9815) and (412.3724,846.2619) .. (411.0236,846.8504);
  \path[-<,draw=black,dash pattern=on 1.64pt off 0.82pt,line join=miter,line
    cap=butt,miter limit=4.00,line width=0.819pt] (142.6011,896.4220) .. controls
    (141.7903,849.6370) and (141.7903,849.6370) .. (141.7903,852.6370);
  \path[->,draw=black,dash pattern=on 1.72pt off 0.86pt,line join=miter,line
    cap=butt,miter limit=4.00,line width=0.861pt] (300,847.5927) --
    (300,884.5075);
  \path[->,draw=black,dash pattern=on 1.87pt off 0.94pt,line join=miter,line
    cap=butt,miter limit=4.00,line width=0.935pt] (397.1003,894.5494) .. controls
    (397.0963,845.9170) and (397.0933,863.2926) .. (397.0897,847.6642);
  \path[->,draw=black,dash pattern=on 3.20pt off 6.40pt,line join=miter,line
    cap=butt,miter limit=4.00,line width=0.800pt] (106.2992,919.6902) .. controls
    (106.2992,981.4960) and (106.2992,981.4960) .. (106.2992,990);
  \path[-<,draw=black,dash pattern=on 3.27pt off 6.55pt,line join=miter,line
    cap=butt,miter limit=4.00,line width=0.818pt] (227.1716,839.7815) .. controls
    (226.7718,995.6693) and (226.7718,995.6693) .. (226.7718,995.6693);
  \path[->,draw=black,dash pattern=on 3.63pt off 7.26pt,line join=miter,line
    cap=butt,miter limit=4.00,line width=0.907pt] (354.3307,931.8898) --
    (354.3307,998);
  \path[->,shift={(0,768.89762)},draw=black,dash pattern=on 3.20pt off 6.40pt,line
    join=miter,line cap=butt,miter limit=4.00,line width=0.800pt]
    (425.1968,77.9528) -- (425.1968,219.6850);
  \path[fill=black] (71,896.45667) node[above right] (text3859) {$g_k$};
  \path[fill=black] (248.03149,797) node[above right] (text3863)
    {$g_{k+1}$};
  \path[fill=black] (325,897) node[above right] (text3867)
    {$g_{k+2}$};
  \path[fill=black] (388,825) node[above right] (text3871)
    {$g_{k+3}$};
\end{scope}
\end{tikzpicture}
  \caption[A smooth orbifold path]{Image of a smooth orbifold path together with lifts on a special range atlas.}
\end{figure}

 We remark that representatives with the special properties discussed in Lemma \ref{lem: nice:geod} do not only exist around a given compact set (cf.\ Lemma \ref{lem: og:equiv}).  

\begin{defn}[Orbifold geodesic]\label{defn: ofd:geod} \no{defn: ofd:geod}
 Let $[\hat{c}] \in  \ORB (\iI ,(Q,\uU))$ be a smooth path in a Riemannian orbifold. The map $[\hat{c}]$ is an \ind{orbifold map!geodesic}{orbifold geodesic} if there is a representative $(c, \set{c_i}_{i \in I} , [P, \nu]) \in \Orb{\vV, \set{(V_j,G_j,\psi_j}_{j \in J}}$ with $\vV \subseteq \aA_\iI$ such that for each $i \in I$ the lift $c_i \colon ]l(i), r(i)[ \rightarrow V_{\alpha (i)}$ is a geodesic. Here $(V_{\alpha (i)}, \rho_{V_{\alpha (i)}})$ is the Riemannian manifold, where $\rho_{V_{\alpha (i)}}$ is the member of the Riemannian orbifold metric. If $[\hat{c}]$ is a geodesic, then the map $c \colon \iI \rightarrow Q$ is called a \ind{geodesic arc}{(geodesic) arc}. Sometimes we will by abuse of notation also call the image of $c$ a (geodesic) arc. 
\end{defn}

\begin{ex}\label{ex: geodesic} \no{ex: geodesic}
 Return to Example \ref{ex: Z2hs}: Consider $\gamma \colon \RR^2 \rightarrow \RR^2, (x,y) \mapsto (-x,y)$ and the orbifold $\RR^2/\langle \gamma\rangle \cong H$ (where $H$ is the right half plane in $\RR^2$). Endow the global chart $(\RR^2, \langle \gamma \rangle , \psi)$ with the flat Riemannian metric. 
As $\langle \gamma \rangle \subseteq O(2)$, this Riemannian metric is $\langle \gamma \rangle$-invariant. Non-trivial geodesics in this metric are straight lines, which induce geodesics of orbifolds. Geodesics contained either in the right or left half plane are mapped to straight lines in the quotient. Standard Riemannian geometry shows that a connected component of the set of points fixed jointly by a set of Riemannian isometries is a closed totally geodesic submanifold (cf.\ \cite[II. Theorem 5.1]{kobayashi1972}). Since $\langle \gamma \rangle$ acts by Riemannian isometries, geodesics which contain singular points either pass through the singular locus in one point or are contained in it. Furthermore, geodesics which pass through the singular locus, are reflected (as befits an example called mirror in $\RR^2$). The following figure depicts an arc of this type:\\
   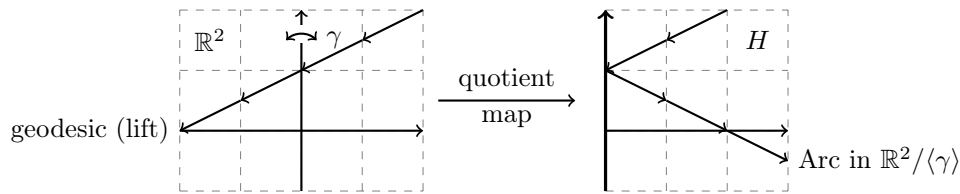
\begin{figure}[h]\centering
   \begin{tikzpicture}[scale=0.8]
     \draw[style=help lines,dashed,gray] (-2,-1) grid (2,2);
    \coordinate (y) at (0,2);
    \coordinate (x) at (2,0);
    \draw[thick,->] (0,-1) -- (y); 
    \draw[thick,->] (-2,0) -- (x); 
    \draw(0,1.6) node[shape=rectangle,fill=white]{};
    \draw(-1.5,1.5) node[fill=white]{$\RR^2$};
    \draw[thick,<->] (-0.25,1.5) to [bend left =45](0.25,1.5) node[right,fill=white]{$\gamma$};
    \draw[thick, ->] (2,2) -- (1,1.5); \draw[thick, ->] (1,1.5) -- (0,1); \draw[thick, ->] (0,1) -- (-1,0.5);  \draw[thick, ->] (-1,0.5)--(-2,0) node[left]{geodesic (lift)};
    \draw[thick,->] (2.25,0.5) -- node[midway,above]{\text{quotient}} node[midway,below]{\text{map}}(4.5,0.5);
    \draw[style=help lines,gray,dashed] (5,-1) grid (8,2);
    \draw(7.5,1.5) node[fill=white]{$H$};
    \draw[very thick,->] (5,-1) -- (5,2); 
    \draw[thick, ->] (5,0) -- (8,0);
    \draw[thick, ->] (7,2) -- (6,1.5); \draw[thick, ->] (6,1.5) -- (5,1); \draw[thick, ->] (5,1) -- (6,0.5);  \draw[thick, ->] (6,0.5)--(7,0);  \draw[thick, ->] (7,0) -- (8,-0.5) node[right]{Arc in $\RR^2/\langle \gamma\rangle$};
  \end{tikzpicture}
 \caption[Orbifold geodesics I: The reflected line]{Orbifold geodesic in $\RR^2/\langle \gamma\rangle$: Reflected line}
  \end{figure}  \\
In particular, orbifold geodesics behave differently from geodesics in Riemannian manifolds. It is well known that the arc of an orbifold geodesic may be not even locally length minimizing (cf. \cite[2.4.2]{cgo2006}). The following picture (which is slightly wrong to show the reflection) illustrates this behavior:\\                                                                                                                                                                                                                                                                                                                                                                                                                                                                                                                                                                                                                                                                                                                                                                                                                                                                                                                                                                                                                                                                                                                                                                                 
   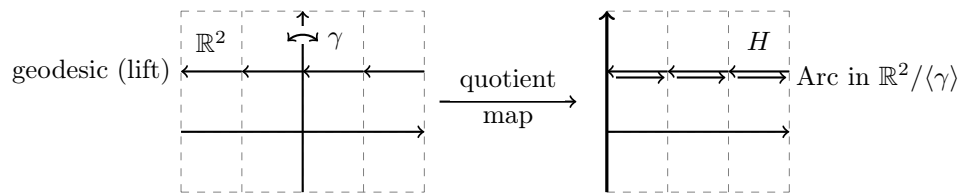
\begin{figure}[h] \centering
   \begin{tikzpicture}[scale=0.8]
     \draw[style=help lines,dashed,gray] (-2,-1) grid (2,2);
    \coordinate (y) at (0,2);
    \coordinate (x) at (2,0);
    \draw[thick,->] (0,-1) -- (y); 
    \draw[thick,->] (-2,0) -- (x); 
    \draw(0,1.6) node[shape=rectangle,fill=white]{};
    \draw(-1.5,1.5) node[fill=white]{$\RR^2$};
    \draw[thick,<->] (-0.25,1.5) to [bend left =45](0.25,1.5) node[right,fill=white]{$\gamma$};
    \draw[thick, ->] (2,1) -- (1,1); \draw[thick, ->] (1,1) -- (0,1); \draw[thick, ->] (0,1) -- (-1,1);  \draw[thick, ->] (-1,1)--(-2,1) node[left]{geodesic (lift)};
    \draw[thick,->] (2.25,0.5) -- node[midway,above,fill=white]{\text{quotient}} node[midway,below]{\text{map}}(4.5,0.5);
    \draw[style=help lines,gray,dashed] (5,-1) grid (8,2);
    \draw(7.5,1.5) node[fill=white]{$H$};
    \draw[very thick,->] (5,-1) -- (5,2); 
    \draw[thick, ->] (5,0) -- (8,0);
     \draw[thick, ->] (8,1) -- (7,1); \draw[thick, ->] (7,1) -- (6,1); \draw[thick, ->] (6,1) -- (5,1); \draw[thick, ->] (5.15,0.9) -- (5.95,0.9);  \draw[thick, ->] (6.15,0.9)--(6.95,0.9);  \draw[thick, ->] (7.15,0.9) -- (7.95,0.9) node[right]{Arc in $\RR^2/\langle \gamma\rangle$};
  \end{tikzpicture}
 \caption[Orbifold geodesics II: A non length-minimizing geodesic]{Orbifold geodesic in $\RR^2/\langle \gamma\rangle$: Not length minimizing in any neighborhood of the singularity.}
  \end{figure} \\
For further examples of orbifold geodesics (in particular, closed geodesics on orbifolds) we refer to \cite[2.4.5]{cgo2006}. 
\end{ex}

\begin{prop}\label{prop: geod:rep} \no{prop: geod:rep}
 Let $[\hat{c}] \in  \ORB (\iI,(Q,\uU))$ be an orbifold geodesic together with a representative $\hat{g} = (c, \set{g_j}_{j \in J}, [P , \nu])$ of $[\hat{c}]$. If the domain atlas of $\hat{g}$ is contained in $\aA_\iI$, then each lift $g_j$ is a geodesic.
\end{prop}

\begin{proof}
 As $[\hat{c}]$ is an orbifold geodesic, there is a representative $\hat{c}=(c, \set{c_i}_{i \in I}, [P',\nu']) \in \Orb{\vV , \vV'}$ such that every $c_i$ is a geodesic in $(V_i, \rho_i)$. Furthermore, the domain atlas of $\hat{c}$ is contained in $\aA_\iI$, $(V_i,G_i,\psi_i) \in \vV'$ and $\rho_i$ is the member of the Riemannian orbifold metric on this chart. Consider the lifts $g_j \colon \dom g_j \rightarrow W_j$ of $\hat{g}$ with respect to the charts $(W_i,H_i, \pi_i)$ in the range atlas of $\hat{g}$. Since $\hat{c} \sim \hat{g}$, the definition of equivalence for orbifold maps yields the following data: There are lifts $\ve$ and $\ve'$ of the identity on $\iI$, respectively $\ve''$  and $\ve'''$ on $(Q,\uU)$ together with a charted map of orbifolds $\hat{h}$ such that 
	\begin{align*}
	 \hat{c} \circ \ve = \ve' \circ \hat{h} & \text{ and } &  \hat{g} \circ \ve''= \ve''' \circ \hat{h}.
	\end{align*}
 We consider for $j \in J$ some $t \in \dom g_j$. As $t \in \iI$, there is an index $i \in I$ with $t \in \dom c_i$. Recall from Definition \ref{defn: idmor} that the lifts of $\ve , \ve', \ve''$, and $\ve'''$ are local diffeomorphism. In particular, they restrict to embeddings of orbifold charts on open sets by Proposition \ref{prop: ll:id}. Together with Lemma \ref{lem: lft:rel}, we obtain open neighborhoods $U \subseteq \dom c_i$ of $t$ and $V \subseteq V_j$ of $g_j (t)$ such that: There are changes of charts $\lambda \colon \dom c_i \ \supseteq U \rightarrow \dom g_j$ and $\mu \colon V_i \supseteq V \rightarrow W_j$ with 
	\begin{equation}\label{eq: geod:loc}
	 g_{j} \circ \lambda  = \mu \circ c_i|_U.
	\end{equation}
 The domain atlases are contained in $\aA_\iI$, whence $\dom c_i , \dom g_j \subseteq \iI$ and their chart maps are induced by the inclusions of sets. Hence the change of charts $\lambda \colon U \rightarrow \dom g_j$ is the inclusion of an open subset. Thus $g_{j}|_{U} = \mu \circ c_i|_{U}$. As $(Q,\uU, \rho)$ is a Riemannian orbifold, $\mu$ is a Riemannian isometry. Since isometries preserve geodesics (cf.\ \cite[IV. Proposition 2.6]{fdiffgeo1963}), the identity \eqref{eq: geod:loc} shows that in a neighborhood of $t$, the map $g_j$ is a geodesic in $(W_j,\rho_j)$. The construction did neither depend on $j\in J$ nor on $t$, whence $g_j$ is a geodesic for each $j \in J$. 
\end{proof}
\noindent
Two orbifold geodesics coincide on a joint interval $\iI$ if and only if their initial vectors coincide (cf.\ Lemma \ref{lem: og:equiv}). On a Riemannian manifold, geodesics are uniquely determined by their initial data in one point. The same holds for orbifold geodesics:

\begin{prop}\label{prop: unq:sdat} \no{prop: unq:sdat}
 Consider $p \in Q$, $\xi \in \tT_p Q$.
	\begin{compactenum}
	 \item There is an $\ve >0$ such that there exists an orbifold geodesic $\hat{c}_\xi \in \Orb{]-2\ve,2\ve[ , (Q,\uU)}$ with initial vector $\xi$ in $0$. 
	 \item Let $[\hat{c}] \in \ORBM[\iI, (Q,\uU)]$ and $[\hat{c}'] \in \ORBM[\iI',(Q,\uU)]$ be orbifold geodesics. If there exists $a \in \iI \cap \iI'$ such that the initial vectors of $\hat{c}$ and $\hat{c}$ in $a$ coincide, then the initial vectors of $[\hat{c}]$ and $[\hat{c}']$ coincide at each point in $\iI \cap \iI'$, whence $[\hat{c}']|_{\iI \cap \iI'} = [\hat{c}]|_{\iI \cap \iI'}$ holds.
	\end{compactenum}
\end{prop}
 
\begin{proof}
 \begin{compactenum}
 \item Choose some representative $(\pi , X) \in \xi$, where $(V,G,\pi) \in \uU$ and $X \in T_x V$ such that $T\pi (X) = \xi$. Set $x = \pi_{TV} (X)$. Let $\rho_V$ be the member of the  Riemannian orbifold metric on $V$, i.e.\ $(V,\rho_V)$ is a Riemannian manifold. Standard Riemannian geometry (cf.\ \cite[III. Theorem 6.4.]{fdiffgeo1963}) shows that there is an $\ve > 0$ and a geodesic $c_0 \colon ]-2\ve , 2\ve [  \rightarrow V$ with initial condition $(x,X)$, i.e. $c_0(0)=x$ and $T_0 c_0(1) = X$. Let $c \coloneq \pi \circ c_0$, $P \coloneq \set{\id_{]-2\ve ,2\ve[}}$ and $\nu \colon P \rightarrow \PSI (\uU)$ be the map which sends the element of $P$ to $\id_V$. We obtain an orbifold geodesic $\hat{c} \coloneq (c,\set{c_0}, P ,\nu)$. By construction the initial vector of $\hat{c}$ in $0$ is $\xi$.
 \item Since $\iI \cap \iI'$ is an open submanifold of $\iI , \iI'$ the orbifold maps restrict to orbifold maps in $\ORBM[\iI \cap \iI' , (Q,\uU)]$. To shorten our notation we may therefore assume that $\iI = \iI'$ and $a = 0$ holds. Choose representatives $\hat{c} = (c, \set{c_k}_{k \in I}, [P, \nu]) \in \Orb{\vV , \wW}$ and $\hat{c}'=(c',\set{c_r}_{r \in J}, [P', \nu']) \in \Orb{\vV',\wW'}$ whose domain atlases are subsets of $\aA_\iI$. We will check the condition of Lemma \ref{lem: og:equiv} (b), which is equivalent to the assertion:\\
 As a first step, show that there is $\ve >0$ such that for each $t \in ]-\ve , \ve[$ the condition of Lemma \ref{lem: og:equiv} (b) holds. Let $c_0$ be a lift of $\hat{c}$ and $c_0'$ be a lift of $\hat{c}'$ with $0 \in \dom c_0 \cap \dom c_0'$. Set $\cod c_0 = V_0$ and $\cod c_0' = V_0'$ for orbifold charts $(V_0,\set{\id_{V_0}},\pi_0)$ and $(V_0',\set{\id_{V_0'}},\pi_0')$, respectively.
 The geodesics pass through $c(0) = c'(0)$ with initial vector $\xi \in \tT_{c (0)} (Q,\uU)$. The construction of $\xi \in \tT_{c(0)} (Q,\uU)$ yields a change of charts $\lambda_0 \colon V_0 \supseteq U \rightarrow V \subseteq V_0'$ such that $T_0\lambda_0 c_0 (1) = T_0c_0' (1)$. The lifts $c_0$ and $c_0'$ are geodesics and $\lambda_0$ is an isometry. Uniqueness of geodesics on Riemannian manifolds now assures that there is an $\ve >0$ such that $T_t \lambda_0 c_0 (1) = T_t c_0' (1)$ for all $t \in ]-\ve ,\ve[$.\\
 We claim that the subset of $\iI$ where the condition of Lemma \ref{lem: og:equiv} (b) holds contains $\iI \cap [0,\infty[$.\\ Assume that this was not the case and consider
	\begin{displaymath}
		t_0 \coloneq \inf \setm{t \in \iI}{t > 0,\ \nexists \lambda \in \Ch_{\wW \cup \wW'}\ : t \in \dom c_k \cap \dom c_r' \text{ and } T_t\lambda c_k (1) = T_t c_r' (1)}.
	\end{displaymath}
 Let $c_k$ be the local lift of $c$ and $c_r'$ be the local lift of $c'$ such that $t_0 \in \dom c_k \cap \dom c_r'$. Their images are contained in $(V_k,G_k,\pi_k)$ and $(V_r,G_r,\pi_r)$, respectively. The first step assures that $t_0 >0$ and by construction, the condition of Lemma \ref{lem: og:equiv} (b) holds for all smaller $t$. This forces $c$ and $c'$ to coincide on $[0,t_0[$ and by continuity of these maps, we obtain $c(t_0) = c'(t_0)$. Thus there is a change of charts $\lambda \colon V_k \supseteq U \rightarrow V \subseteq V_r$ with $\lambda  c_k (t_0) = c_r'(t_0)$. Choose some $t < t_0$ with $c_k ([t,t_0]) \subseteq \dom \lambda$. Since $t < t_0$ holds, there is a change of charts $\mu$ with $T_t \mu c_k (1) = T_tc_r (1)$. Shrinking the domain of $\mu$, we may assume that $\mu$ is an embedding of orbifold charts and $\dom \mu \subseteq \dom \lambda$ is satisfied. Now $\lambda|_{\dom \mu}$ is an embedding of orbifold charts mapping $\dom \mu$ into $V_r$. By Proposition \ref{prop: ch:prop} (d) there is an element $h \in G_r$ such that $h \circ \lambda|_{\dom \mu} = \mu$. The change of charts $\lambda_{t_0} \coloneq h \circ \lambda$ is a Riemannian isometry which satisfies $T_t \lambda_{t_0} c_k (1) = T_t \mu c_k (1) = T_tc_r' (1)$. We deduce that on its domain, $\lambda_{t_0}$ maps the geodesic $c_k$ to $c_r$. There is some $\delta > 0$ such that $c_k(]t_0-\delta , t_0+\delta[) \subseteq \dom \lambda_{t_0}$ holds. Hence $T_{s} \lambda_{t_0} c_k (1) = T_s c_r'(1)$ holds for each $s \in ]t_0-\delta , t_0 + \delta[$. This contradicts our choice of $t_0$ and thus there may be no such point in $\iI \cap [0,\infty[$. An analogous argument for $t <0$ shows that the condition of Lemma \ref{lem: og:equiv} (b) holds for all of $\iI$, whence both orbifold geodesics coincide. 
 \end{compactenum}
\end{proof}

Since orbifold geodesics are uniquely determined by their initial vectors in some point, we may construct a join for two suitable geodesics: 

\begin{lem}\label{lem: join:ofg} \no{lem: join:ofg}
 Let $[\hat{c}] \in  \ORBM[\iI , (Q,\uU)]$ and $[\hat{c}'] \in  \ORBM[\iI' , (Q,\uU)]$ be orbifold geodesics such that for some $x_0 \in \iI \cap \iI'$, their initial vectors coincide. Then there is a unique orbifold geodesic $[\hat{c} \vee \hat{c}'] \in\ORBM[\iI \cup \iI' , (Q,\uU)]$ such that $[\hat{c} \vee \hat{c}']|_{\iI'} = [\hat{c}']$ and $[\hat{c} \vee \hat{c}']|_{\iI} = [\hat{c}]$.
\end{lem}

 \begin{proof} It is possible to \tl glue\tr two orbifold geodesics whose initial vector coincides in one point. This procedure, together with a full proof, can be found as Lemma \ref{GLUE!} in Appendix \ref{App: Haefgeodesics}.
  \end{proof}

Standard Riemannian geometry shows that the maximal domain $\iI$ has to be an open subset of $\RR$ (since the lifts of an orbifold geodesic are geodesics in suitable charts, whose maximal domain is always an open subset of $\RR$). Naturally we have to ask whether the orbifold geodesic constructed in Proposition \ref{prop: unq:sdat} (a) may uniquely (up to equivalence of orbifold morphisms) be extended to a maximal domain. In fact each geodesic with this initial vector in $0$ may then be derived as a restriction of the maximal geodesic. The next lemma is inspired by a lemma due to Chen and Ruan (cf.\ \cite[Lemma 4.2.6]{chen2001}): 

 \begin{lem} \label{lem: cpt:gcp} \no{lem: cpt:gcp}
  Let $p \in Q$ be any point and $\xi \in \tT_p Q$. 
    \begin{compactenum}
     \item There is a unique maximal interval $\iI_\xi$ such that an orbifold geodesic $[\hat{c}_\xi] \in \ORBM[\iI_\xi , (Q,\uU)]$ with initial vector $\xi$ in $0$ exists on $\iI_\xi$.
     \item If $Q$ is compact, then $\iI_\xi = \RR$ holds, for each $\xi \in \tT Q$. 
    \end{compactenum}
 \end{lem}

\begin{proof}
 \begin{compactenum}
  \item Let $S_\xi$ be the set of all orbifold geodesics whose initial vector at $0$ is $\xi$. Orbifold geodesics with initial vector $\xi$ at $0$ exist by Proposition \ref{prop: unq:sdat} (a), whence $S_\xi$ is non-empty. For two elements $[\hat{c}], [\hat{c}'] \in S_\xi$, there is a join $[\hat{c} \vee \hat{c}']$ by Lemma \ref{lem: join:ofg} which is again an element of $S_\xi$. Any finite number of elements in $S_\xi$ may be joined in this way. For $[\hat{c}] \in S_\xi$, we let $\iI_{\hat{c}}$ be the interval such that $[\hat{c}] \in \ORBM[\iI_{\hat{c}}, (Q,\uU)]$.\\
  Construct recursively an element $[\hat{c}_{\xi}] \in S_\xi$ on the open subset $\iI_{\xi} \coloneq \bigcup_{[\hat{c}] \in S_\xi} \iI_{\hat{c}}$. The set $\iI_{\xi}$ is an open connected subset of $\RR$ as a union of connected open subspaces with non-empty intersection (cf.\ \cite[Corollary 6.1.10]{Engelking1989}). Define $c \colon I_\xi \rightarrow Q$ via $c(t) \coloneq c'(t)$ if $t\in \dom c'$ with $\hat{c}' = (c',\set{c_i'}_{i\in I}, [P,\nu]) \in S_\xi$. This map is well-defined by Proposition \ref{prop: unq:sdat}. There exist numbers 
    \begin{displaymath}
     \cdots < a_{-2} < a_{-1} < a_0 = 0 < a_1 < a_2 < \cdots
    \end{displaymath}
 such that $I_\xi = \bigcup_{k \in \ZZ} [a_k,a_{k+1}]$ and such that, for each $k \in \ZZ$ a lift $c_k$ of some $\hat{c} \in S_\xi$ is defined on an open interval $I_k$ containing $[a_k,a_{k+1}]$, with image in $(V_k,G_k,\pi_k)$ and $c_k$ is a geodesic. Choose $l(k)$ so large and $r(k-1)$ so small that $a_k<r(k-1)<l(k+1)<a_{k+1}$ and $]l(k),r(k)[ \subseteq I_k$ hold and there exists a change of charts $\lambda_{k,k+1}$ with $c_k([l(k+1),r(k)]) \subseteq \dom \lambda_{k,k+1}$ with $\lambda_{k,k+1} \circ c_k|_{[l(k+1),r(k)]} = c_{k+1}|_{[l(k+1),r(k)]}$. Let $\wW$ be an atlas containing all $(V_k,G_k,\pi_k)$. Define $P \coloneq  \setm{\id_{]l(k),r(k)[}, \iota_{k}^{k+1}, \iota_{k+1}^k}{k \in \ZZ}$ where $\iota_k^{k+1} \colon ]l(k+1),r(k)[ \rightarrow ]l(k+1),r(k+1)[$ and $\iota_{k+1}^k \colon ]l(k+1),r(k)[ \rightarrow ]l(k),r(k)[$ are inclusions of sets. Now define $\nu \colon P \rightarrow \Psi (\wW)$ via $\nu (\id_{]l(k),r(k)[}) \coloneq \id_{V_k}, \nu (\iota_{k}^{k+1}) \coloneq \lambda_{k,k+1}$ and $\nu(\iota_{k+1}^{k}) \coloneq \lambda_{k,k+1}^{-1}$. Then $\hat{c} \coloneq (c,\set{c_k|_{]l(k),r(k)[}}_{k \in \ZZ} , [P,\nu])$ is a geodesic and $\hat{c} \in S_\xi$ because $T\pi_0 Tc_0 (0,1)$ also is the initial vector of some $\hat{c}' \in S_\xi$, and hence equal to $\xi$.
  \item Following (a), it is sufficient to prove that an orbifold geodesic $[\hat{c}] \in \ORBM[\iI , (Q,\uU)]$ with initial vector $\xi$ at $0$ and $\iI = ]a,b[$ may be extended in the following sense: If there is a sequence $(t_n)_{n \in \NN} \in ]a,b[$ such that $t_n \rightarrow b$ with $b < \infty$ and $\lim c (t_n)$ exists in $Q$, then there is an orbifold geodesic $[\hat{c}']$ defined on $]a,b'[,\ b' > b$, whose initial vector at $0$ is $\xi$. Set $q \coloneq \lim_{n \in \NN} c (t_n)$ and choose an orbifold chart $(V,G_x,\psi)$ with $q = \psi (x)$ for $x \in V$ and $G=G_x$. Notice that $\psi^{-1} (q) = \set{x}$ holds. Choose a compact neighborhood $U_x$ of $x$ and observe that $G_x . U_x$ is again a compact set. A compactness argument together with \cite[3.2 Proposition 2.5]{rg1992} proves that there are $\delta >0$ and $\ve >0$ such that for each $p \in U_x$ and $v \in B_{\rho_V} (0_q, \ve)$, there is a unique geodesic $\gamma_v \colon ]-\delta , \delta[ \rightarrow V$ with initial value $T_0 \gamma (1) = v$. Here $\rho_V$ is the member of the Riemannian orbifold metric on $V$. 
 For $N$ large enough one obtains $c(t_n) \in  \psi (U_x), \ \forall n \geq N$. The definition of an orbifold geodesic implies that for each $t_n$ there is some local lift $c_{n} \colon \dom c_{n} \rightarrow V_{n}$ of $c$ with $t_n \in \dom c_n$ and $(V_{n}, H_{n}, \varphi_{n}) \in \uU$. By compatibility of orbifold charts, $c(t_n) \in \im \varphi_{n} \cap \psi (U_x)$ for $n\geq N$ implies that there is some change of orbifold charts $\lambda_{n}$ with $\lambda_{n} c_{n} (t_n) \in G_x.U_x$. As each $\lambda_{n}$ is a Riemannian embedding, the definition of an orbifold geodesic yields $\norm{T_{t_n} \lambda_{n} c_{n} (1)}_{\rho_{V}} = K = \norm{T_{t_m} \lambda_{m} c_{m} (1)}_{\rho_{V}}$ for all $n,m \geq N$. Using homogeneity of geodesics on Riemannian manifolds (\hspace{-0.5pt}\cite[3.2 Lemma 2.6]{rg1992}),for each $q \in G_x.U_x$ there is some $\delta' >0$ such that for each $v \in B_{\rho_V} (0_q, K+1)$ the geodesic with initial value $v$ exists on $]-\delta',\delta'[$. Let $\gamma_{X}$ be the geodesic in $(V, \rho_V)$ with initial vector $X$.  Choose $n_0 > N$ so large that $b-t_{n_0} < \delta'$ holds. The geodesic $g_{n_0} \colon ]t_{n_0}- \delta , t_{n_0} + \delta'[ \rightarrow V, t \mapsto \gamma_{T_{t_{n_0}} \lambda_{n_0} c_{n_0} (1)} (t-t_{n_0})$ induces an orbifold geodesic $\hat{g} \coloneq (\psi \circ g_{n_0},\set{g_{n_0}}, \set{\id_{]t_{n_0}-\delta', t_{n_0}+\delta'[}}, \nu)$ where $\nu (\id_{]t_{n_0}-\delta', t_{n_0}+\delta'[}) \coloneq \id_V$. By construction, the initial vector of $\hat{g}$ in $t_{n_0}$ coincides with the initial vector of $\hat{c}$ in $t_{n_0}$. Thus Lemma \ref{lem: join:ofg} yields an orbifold geodesic $\hat{c} \vee \hat{g}$ which is defined on $]a,t_{n_0}+\delta'[$. The initial vector of $\hat{c} \vee \hat{g}$ in $0$ is $\xi$ and its domain strictly contains $]a,b[$.  
 \end{compactenum} 
\end{proof}
 \begin{rem}
  The maximal geodesics $[\hat{c}_\xi]$ on $\iI_{\xi}$ constructed in Lemma \ref{lem: cpt:gcp} (a) do not extend, i.e.\ if $[\hat{g}] \in \ORBM[\iI,(Q,\uU)]$ is a geodesic whose initial vector at $a \in \iI \cap \iI_{\xi}$ coincides with the initial vector of $[\hat{c}_\xi]$ in $a$, then $\iI \subseteq \iI_{\xi}$ and $[\hat{c}_\xi]|_{\iI} = [\hat{g}]$ hold.
 \end{rem}

\begin{thm}\label{thm: geod:cont} \no{thm: geod:cont}
 Let $(Q,\uU,\rho)$ be a Riemannian orbifold and $\xi \in \tT Q$.
	\begin{compactenum}
	 \item There exist $\delta , \delta'> 0$, an open neighborhood $O_\xi \subseteq \tT Q$ of $\xi$ and a continuous map $\alpha_\xi \colon ]-\delta , \delta'[\ \times O_\xi \rightarrow Q$ and for $\xi' \in O_\xi$ the path $\alpha_{\xi} (\cdot ,\xi') \colon ]-\delta , \delta'[\ \rightarrow Q,  t \mapsto \alpha (t,\xi')$ is the geodesic arc of an orbifold geodesic $[\hat{c}_{\xi'}]$ with initial vector $\xi'$ in $0$. We call $\alpha_\xi$ a \ind{orbifold geodesic flow}{orbifold geodesic flow,}
	 \item If $(\xi,\zeta) \in \tT Q \times \tT Q$ with $O_\xi \cap O_\zeta \neq \emptyset$, then $\alpha_\xi$ and $\alpha_\zeta$ coincide on the intersection of their respective domains. 
	 \item If the maximal orbifold geodesic $[\hat{c}_\xi]$ with initial vector $\xi$ in $0$ satisfies $[c,d] \subseteq \iI_\xi$, then the set $O_\xi$ in (a) may be constructed such that for $\zeta \in O_\xi$, the map $[\hat{c}_\zeta]$ is defined on $[c,d]$.
	\end{compactenum}
\end{thm}

\begin{proof}
 \begin{compactenum}
  \item By Proposition \ref{prop: unq:sdat} (a), there is some $\ve > 0$ together with the representative of an orbifold geodesic $\hat{c} = (c, \setm{g_i}{1\leq i \leq N}, P ,\nu)$ defined on $]-2\ve,2\ve[$ with initial vector $\xi$ in $0$. After shrinking the domain, without loss of generality $\hat{c}$ is defined on an open neighborhood of $[-\ve,\ve]$ with properties as in Lemma \ref{lem: nice:geod}. We show that there is an open neighborhood of $\xi$ such that each orbifold geodesic with initial vector in this set exists at least on $[0,\ve]$.\\ 
  To shorten the notation, relabel the charts as $\set{-t,-t+1,\ldots , 0, 1, \ldots ,s}$ for certain $s,t \in \NN_0$ such that $0 \in \dom g_0$. Let $g_i \colon ]l(i),r(i)[ \rightarrow U_i,\  -t \leq i \leq s$ be the lifts, where the $(U_i, G_i , \psi_i)$ are charts in $\uU$. By construction, for $-t \leq i < s$ there is a change of charts $\lambda_{i}^{i+1}$ satisfying $\lambda_{i}^{i+1} g_i|_{]l(i+1),r(i)[}= g_{i+1}|_{]l(i+1),r(i)[}$. Choose for $1 \leq i\leq s$ a point $z_i \in ]l(i),r(i-1)[$ with $z_0 \coloneq 0 < z_i < z_j$ for $i <j$. Define $X_i \coloneq T_{z_i} g_i (1)$ for $0 \leq i \leq s $ and observe that $g_i$ is uniquely determined by $X_i$. By construction $[\psi_0 ,X_0] = \xi$ holds. Finally choose $z_{s+1} \in \dom g_s$ with $z_{s+1} > \ve$ and $z_{s+1} > z_s$. \\
 Standard Riemannian geometry on manifolds shows that the geodesic flow depends smoothly on the initial data (cf.\ \cite[3.2 Proposition 2.5]{rg1992} and \cite[IV, \S 3 and VII, \S 7]{langdgeo2001}, respectively). On the Riemannian manifold $(U_i,\rho_i)$, there is a geodesic flow $\varphi_i \colon \dD_i \rightarrow T U_i$, defined on an open set $\dD_i \subseteq \RR \times TU_i$ (cf.\ \cite[IV, \S 4 Remark before Corollary 4.3]{langdgeo2001}). The map $\varphi_i$ is smooth by an application of \cite[IV, \S 2 Theorem 2.6]{langdgeo2001}. Since $g_s$ is a geodesic defined on $[z_s , z_{s+1}] \subseteq ]l(s),r(s)[$ with $T_{z_s} g_s (1)= X_s$, the compact set $[0 , z_{s+1}-z_s] \times \set{X_s}$ is contained in the open set $\dD_s$. An application of Wallace Theorem \cite[3.2.10]{Engelking1989} provides an open neighborhood $[0, z_{s+1}-z_s] \times \set{X_s} \subseteq ]- \delta_s , z_{s+1}-z_s +\delta_s[ \times V_s \subseteq \dD_s$. For each element $\zeta$ of this neighborhood in $TU_s$, the geodesic with initial data $\zeta$ exists on the interval $]z_s-\delta_s, z_s +\delta_s[$. \\
 Shrinking $V_s$ and $\delta_s$, we may assume that $V_s \subseteq \pi_{TU_s}^{-1} (\cod \lambda_{s-1}^s)$ and $z_s-\delta_s > r(s-2)$ hold. Identify $T\cod \lambda_{s-1}^s$ and $T\dom \lambda_{s-1}^s$ with open subsets of $TU_s$ respectively $TU_{s-1}$ and set $V_s' \coloneq (T\lambda_{s-1}^s)^{-1} (V_s) \subseteq TU_{s-1}$. The geodesic $g_{s-1}$ is determined by $X_{s-1}$ and its domain $]l(s-1),r(s-1)[$ contains $[z_{s-1}, z_s]$ with $T_{z_s} g_{s-1}(1) \in V_s'$. As the geodesic flow $\varphi_{s-1}$ is smooth, arguments as above applied to $\varphi_{s-1}$ yield an open set $V_{s-1} \subseteq TU_{s-1}$ with $V_{s-1} \subseteq T \cod \lambda_{s-2}^{s-1}$ and 
	\begin{compactitem}
	 \item[-] $[0, z_s-z_{s-1}] \times \set{X_{s-1}} \subseteq ]-\delta_{s-1}, z_s -z_{s-1} + \delta_{s-1}[ \times V_{s-1} \subseteq D_{s-1}$,
	 \item[-] $\varphi_{s-1} (z_s-z_{s-1},V_{s-1}) \subseteq V_s'$,
	 \item[-] $z_{s-1}-\delta_{s-1} > r(s-3)$.
	\end{compactitem}
 Again one obtains an open set $V_{s-1}' \coloneq (T\lambda_{s-2}^{s-1})^{-1} (V_{s-1}) \subseteq TU_{s-2}$. Repeating the argument for each $0 \leq i \leq s-2$, we arrive at an open neighborhood $V_0 \subseteq TU_0$ of $X_0$. For each $\zeta \in V_0$, there is a unique family of geodesics $(c_\zeta^i)_{0\leq i \leq s}$ such that $c_\zeta^i$ is defined at least on $]z_i-\delta_i, z_{i+1} + \delta_i[$. In addition these families satisfy $T_{z_i}\lambda_{i-1}^i c_{i-1} (1) = T_{z_i} c_i (1)$.\\ 
 Repeating the argument for $[-\ve , 0]$, we obtain an open set $V_0^-$ such that the geodesics are defined on $[-\ve,0]$. Set $V\coloneq V_0 \cap V_0^-$ and $\delta \coloneq z_{-t-1} - \delta_{-t}$ and $\delta' \coloneq z_{s+1} + \delta_s$. For each $\zeta \in V$ and $-t\leq i\leq s+1$, the geodesics $c_\zeta^i$ are defined on $[z_{i-1}-\delta_{i}, z_i + \delta_i]$. By construction one may restrict their domains such that $\lambda_{i}^{i+1} c_\zeta^i|_{]z_{i+1}-\delta_{i+1}, z_{i+1}+\delta_i[} = c_\zeta^{i+1}|_{]z_{i+1}-\delta_{i+1}, z_{i+1}+\delta_i[}$ holds. For each $\zeta \in V$, the family $(c_\zeta^i)_{-t\leq i \leq s}$ induces an orbifold geodesic. The continuity of the geodesic flows yields a well-defined continuous map
	\begin{displaymath}
	 \tilde{\alpha} \colon ]-\delta , \delta'[ \times V \rightarrow Q, (t,\zeta) \mapsto \psi_i (c^i_\zeta (t)) \text{ for each } t \in ]z_i - \delta_i , z_{i+1} +\delta_i[ .
	\end{displaymath}
 Consider the orbifold chart $(TU_0, G_0, T\psi_0) \in \tT \uU$ for the tangent orbifold $\tT (Q,\uU)$. Chart maps of orbifold charts are open maps and thus $O_\xi \coloneq T\psi_0 (V)$ is open in $\tT Q$. It contains $\xi = T\psi_0 (X_0)$ and the subspace topology on $O_\xi$ with respect to $Q$ coincides with the quotient topology induced on $O_\xi$ by $T\psi_0$ (since $T\psi_0$ factors via a homeomorphism with open image). The restriction $q\coloneq T\psi_0|_V^{O_\xi}$ is an open, continuous and surjective map. For each $\zeta \in O_\xi$, choose a preimage $\tilde{\zeta} \in q^{-1} (\set{\zeta}) \in V$. Notice that each choice of preimage for $\zeta$ induces an orbifold geodesic with initial vector $\zeta$ at $0$. Following Proposition \ref{prop: unq:sdat} (b) the geodesic arcs obtained from a choice of $q^{-1} (\zeta)$ coincide with the arc of $[\hat{c}_\zeta]$ on the intersection of their domains. Hence each choice defines the same continuous path in $Q$. As $\hat{c}_\zeta$ is defined at least on $]\delta, \delta'[$ the maximal geodesic with initial vector $\zeta$ is defined on this interval. We derive a well-defined map 
	\begin{displaymath}
	 \alpha \colon ]-\delta , \delta'[ \times O_p \rightarrow Q , \ (t,\zeta) \mapsto \tilde{\alpha} (t, \tilde{\zeta}) .
	\end{displaymath}
 The map $\id_{]-\delta , \delta'[} \times q$ is the product of open continuous surjective maps, whence it is itself open, continuous and surjective. In particular, this mapping is a quotient map such that $\tilde{\alpha} =  \alpha \circ (\id_{]-\delta , \delta'[} \times q)$ holds. As $\tilde{\alpha}$ is continuous, \cite[VI.\ Theorem 3.1]{dugun1966} implies that $\alpha$ is a continuous map.
 \item By Proposition \ref{prop: unq:sdat} (b), the arcs of two orbifold geodesics with the same initial data coincide. Hence for each $\omega \in O_\xi \cap O_\zeta$, the arcs of the geodesics coincide, therefore $\alpha_\xi (\cdot , \omega)$ and $\alpha_\zeta (\cdot , \omega)$ coincide on the intersection of their respective domains. This proves the assertion.
 \item Repeat the proof of (a) verbatim with $[c,d] \subseteq \iI_\xi$ instead of $[-\ve,\ve]$. 
 \end{compactenum}
\end{proof}

\begin{cor}\label{cor: geo:znb} \no{cor: geo:znb}
 For every $p \in Q$, there is an open neighborhood $W_p \subseteq \tT Q$ of $0 \in \tT_p Q$ and a continuous map $\alpha \colon ]-2,2[ \times W_p \rightarrow Q$ such that $]-2,2[\rightarrow Q, t \mapsto \alpha (t,\xi)$ is the unique geodesic arc with initial vector $\xi$ in $0$ defined on $]-2,2[$, for each $\xi \in W_p$.
\end{cor}

\begin{proof}
 Choose an arbitrary orbifold chart $(U,G,\psi)$ such that $p = \psi (x)$ for some $x \in U$. By definition $T\psi (0_x) =  0_p \in \tT_p Q$ holds, where $0_x \in T_xU$ is the zero element. Standard Riemannian geometry (see \cite[3.2 Proposition 2.7]{rg1992}) assures that there is a smooth mapping $\gamma \colon ]-2,2[ \times V \rightarrow U$, defined on some open set $V \subseteq TU$ such that each $x \in V$ induces a geodesic in $U$ defined at least on $]-2,2[$. Arguing as in the proof of Theorem \ref{thm: geod:cont}, we choose $W_p \coloneq T\psi (V)$ and $\alpha \colon ]-2,2[ \times W_p \rightarrow Q,  (t, \xi) \mapsto \psi(\gamma (t,x_\xi))$, where $x_\xi$ is an arbitrary preimage of $\xi$ under $T\psi$ in $V$.
\end{proof}

\begin{lem}\label{lem: geod:uarc}
 An orbifold geodesic $[\hat{c}]   \in \ORBM[\iI ,(Q,\uU)]$ is uniquely determined by its underlying map.
\end{lem}

\begin{proof}
 Let $[\hat{c}]$ and $[\hat{c}']$ be orbifold geodesics whose underlying map $c \colon \iI \rightarrow Q$ coincides. Shrinking the domains of definition of the lifts and composing with change of charts as necessary, we can achieve the following: There are representatives $\hat{c}$ of $[\hat{c}]$ and $\hat{c}'$ of $[\hat{c}']$, respectively such that their families of local lifts contain lifts $c_0 ,c_0' \colon ]-\ve , \ve[ \rightarrow V$ for some orbifold chart $(V,G,\psi)$. Since both $c_0$ and $c_0'$ lift $c$, we have $\gamma_x . c_0 (x) = c_0' (x)$ for every $x \in \, ]-\ve ,\ve[$ and $\gamma_x \in G$. Define for $\gamma \in G$ the set
 \begin{displaymath}
  U_\gamma \coloneq \setm{x \in ]-\ve,\ve[}{\gamma . c_0 (x) = c_0'(x)}.                                                                                                                                                                                                                                                                                                                                                                                                                                                                                                                                                               
 \end{displaymath}
 Notice that $U_\gamma$ is a closed set for $\gamma \in G$ and $]-\ve,\ve[\, = \bigcup_{\gamma \in G} U_\gamma$. Since $G$ is finite, Baires theorem asserts that at least some $U_\gamma$ must have non-empty interior. Hence the geodesics $c_0,c_0'$ coincide, up to composition by a group element in $G$, on an open subset of $]-\ve,\ve[$. Let $x$ be a point in the interior of $U_\gamma$. Since the geodesics $\gamma. c_0,c_0'$ coincide on an open neighborhood of $x$, their derivatives must coincide. By Lemma \ref{lem: osmpath:prop}, the initial vectors of both geodesics at $x$ coincide. Hence the assertion follows from Proposition \ref{prop: unq:sdat}.   
\end{proof}

 Albeit the quite similar behavior of orbifold geodesics to geodesics on Riemannian manifolds, not all properties of geodesics may be preserved in the orbifold case. For example, as is noted in \cite[2.4.2]{cgo2006}, orbifold geodesics may not even be locally length minimizing in the natural length metric on $Q$ (induced by piecewise differentiable paths). However, as we are only interested in geodesics as a tool to obtain an exponential map, we shall not investigate this behavior. 

\subsection{The Riemannian orbifold exponential map}
In this section, our main tool derived via Riemannian geometry on orbifolds, the Riemannian orbifold exponential map, is introduced. As before the triple $(Q,\uU , \rho)$ will be a Riemannian orbifold. 
By Lemma \ref{lem: cpt:gcp} (a), for each $\xi \in \tT Q$, there is a maximal orbifold geodesic $[\hat{c}_\xi]$ with initial vector $\xi$ in $0$. The geodesic arc of a maximal orbifold geodesic is unique by Proposition \ref{prop: unq:sdat}. Hence the continuous map of the base spaces $c_\xi \colon \iI_\xi \rightarrow Q$ is uniquely determined. 

\begin{defn}[Riemannian orbifold exponential map]
  Let $\Omega$ be the set of all $\xi \in \tT Q$ such that the orbifold geodesic $[\hat{c}_\xi]$ with underlying map $c_\xi\colon \iI_\xi \rightarrow Q$ satisfies $[0,1] \subseteq \iI_\xi$. The map 
	\begin{displaymath}
	 \expo \colon \Omega \rightarrow Q ,\ \xi \mapsto c_\xi (1) 
	\end{displaymath}
 is called \ind{Riemannian orbifold metric!exponential map}{Riemannian orbifold exponential map}. The set $\Omega$ is an open neighborhood of the zero section, by Theorem \ref{thm: geod:cont} (c) and Corollary \ref{cor: geo:znb}. We call $\Omega$ the \ind{Riemannian orbifold metric!exponential map!domain}{domain of the Riemannian orbifold exponential map}.
\end{defn}

\begin{lem}\label{lem: oem:cont} \no{lem: oem:cont}
 The Riemannian orbifold exponential map is continuous and for $0_p \in \tT_p Q$ the identity $\expo (0_p) = p$ holds.
\end{lem}

\begin{proof}
 Let $\xi \in \Omega$ be arbitrary. The geodesic $[\hat{c}_\xi]$ is defined on an open interval $\iI_\xi$ such that $[0,1] \subseteq \iI_\xi$ holds. By Theorem \ref{thm: geod:cont} (c), there is an open neighborhood $\xi \in O_\xi \subseteq \tT Q$ such that each orbifold geodesic $[\hat{c}_\omega]$ for $\omega \in O_\xi$ is defined on $[0,1] \subseteq ]-\delta , \delta'[$. Furthermore, $O_\xi \subseteq \Omega$ holds. There is a continuous map $\alpha_\zeta \colon ]-\delta , \delta'[ \ \times O_\xi \rightarrow Q, (t,\omega) \mapsto \hat{c}_\omega (t)$ such that by construction $\expo (\omega) = \alpha_\xi (1, \omega), \ \forall \omega \in O_\xi$ is satisfied. Hence $\expo$ restricts to a continuous map on the open set $O_\xi$. Theorem \ref{thm: geod:cont} (b) assures that for any $\zeta \in \Omega$ the maps $\alpha_\zeta (1, \cdot)$ and $\alpha_\xi (1, \cdot)$ coincide on $O_\xi \cap O_\zeta$. 
 From \cite[IV. Theorem 9.4]{dugun1966} we deduce that $\expo$ is continuous.\\
 Choose an arbitrary orbifold chart $(U, G, \psi) \in \uU$ such that $p \in \psi (x)$ for some $x \in U$. The chart $T\psi$ maps $0_x \in T_x U$ to $0_p \in \tT_p Q$. Standard Riemannian geometry assures that the geodesic $\gamma$ starting in $x$ with velocity $0$ is constant and hence defined on all of $\RR$. Setting $c \colon \RR \rightarrow Q, t \mapsto p$, we obtain a representative of an orbifold geodesic $\hat{c} \coloneq (c, \gamma , \set{\id_\RR} , \nu )$, where $\nu (\id_\RR ) \coloneq \id_U$. The orbifold geodesic $[\hat{c}]$ has initial vector $0_p$ in $0$ and its arc is uniquely determined by Proposition \ref{prop: unq:sdat}. This proves $\expo (0_p) = p$.  
\end{proof}

\begin{prop}\label{prop: exp:ofdm} \no{prop: exp:ofdm}
 Consider the open suborbifold $(\Omega , \uU_\Omega)$. The map $\expo$ induces a map of orbifolds $[\expo] \in \ORBM[(\Omega,\uU_\Omega),(Q,\uU)]$ also called  \ind{orbifold map!exponential map}{Riemannian orbifold exponential map}.\glsadd{expo}
\end{prop}

\begin{proof}
 The subset $\Omega \subseteq \tT Q$ is open. Hence the orbifold structure $\tT \uU$ induces a unique orbifold structure $(\Omega , \tT \uU_\Omega)$ (cf.\ Definition \ref{defn: osup}), turning this orbifold into an open suborbifold of $(\tT Q, \tT \uU)$. We claim that there is a representative $\vV$ of $\tT \uU_\Omega$ together with a family of lifts, turning $\expo$ into a charted orbifold map in $\Orb{\vV, \wW}$ for some $\wW \in \uU$. By Lemma \ref{lem: oem:cont}, the map $\expo$ is continuous. Construct smooth lifts of $\expo$: To this end, consider arbitrary $\xi \in \Omega$. By Theorem \ref{thm: geod:cont} and its proof, there is an open neighborhood $\xi \in O_\xi \subseteq \Omega$ together with the following data:
	\begin{compactitem}
	 \item[-]$(TU_1,G_1,T\psi_1) \in \tT \uU$, with $O_\xi = T\psi_1 (V) \subseteq T\psi_1 (TU_1)$ for some open $V \subseteq TU_1$,
	 \item[-]a family of orbifold charts $\set{(U_i, G_i, \psi_i)}_{1\leq i \leq N} \in \uU$,
	 \item[-]a continuous map  $\theta \colon V \rightarrow Q, X \mapsto \tilde{\alpha} (1,X)$ such that $\theta = \expo \circ T\psi_1|_{V}$ holds. The map $\theta$ is the composition of the geodesic flows $\varphi_{i}$ on $(U_i,\rho_i), 1\leq i\leq N$, changes of charts $\lambda_{ii+1}$ for $1\leq i < N$, the bundle projection of $TU_N$ and the orbifold chart $\psi_N$.
	\end{compactitem}
  Recall from the proof of Theorem \ref{thm: geod:cont} that there is a partition $0 = t_0 < t_1 < \cdots < t_N < 1$ such that a smooth map $\text{Exp}_\xi \colon  TU_1 \supseteq V \rightarrow U_N$ may be defined via
		\begin{align}
	\text{Exp}_\xi (X) \coloneq \pi_{TU_N} \varphi_{N} (1-t_N, \cdot) \circ T \lambda_{N-1N} \circ \varphi_{N-1} (t_N-t_{N-1} , \cdot) \circ \cdots \circ T\lambda_{12} \circ \varphi_1 (t_1, \cdot) (X). \label{eq: explift}
		\end{align} 
 Reviewing Theorem \ref{thm: geod:cont}, we see that $\theta = \psi_N \circ \text{Exp}_\xi$.\\ 
 Choose an open $G_1$-stable subset $W$ of $V$ which contains some preimage $x_\xi$ of $\xi$. Restricting $\text{Exp}_\xi$ to $W$, we obtain a smooth map $\text{Exp}_W$ on an orbifold chart $(W, G_W, T \psi_1|_{W})$. By construction, $\text{Exp}_W$ is a smooth lift of $\expo$ on $W$.\\
 We show that any local lift $\text{Exp}_W'$ of $\expo$ obtained via \eqref{eq: explift} with respect to $(W, G_W, T \psi_1|_{ W})$ and $(U_N, G_N, \psi_N)$ but taking other choices for the intermediary charts, geodesic flows and changes of charts, coincides with $\gamma . \text{Exp}_W$ for some $\gamma \in G_N$. \\
 The lifts $\text{Exp}_W$ and $\text{Exp}_W'$ are defined as restriction of a composition of geodesic flows $\varphi_i$, changes of charts $\lambda_{k k+1}$ and the bundle projection $\pi_{TU_N}$ (cf.\ \eqref{eq: explift}). Notice that the flows, changes of charts and the number $N$ may differ for $\text{Exp}_W'$. However, we fixed the chart $\varphi \coloneq \varphi_N = \varphi'_{N'}$. Each $\varphi_i (t_i - t_{i-1} , \cdot )$ is defined on an open subset of $TU_i$. It is a diffeomorphism from this subset onto its (open) image in $TU_i$ (this follows from \cite[IV, $\S$ 2, Theorem 2.9.]{langdgeo2001}). The change of chart $T\lambda_{kk+1}$ are \'{e}tale embeddings. In addition, the bundle projection $\pi_{TU_N}$ is an open map, whence $\text{Exp}_W$ is an open map as a composition of such maps. The same holds for $\text{Exp}_W'$ whose image is contained in $(U_N, G_N, \psi_N)$. The construction of the lifts $\text{Exp}_W$ and $\text{Exp}_W'$ shows that there are diffeomorphisms $\phi_W \colon W \rightarrow O$, $\phi_W' \colon W \rightarrow O'$ onto open sets $O, O'\subseteq TU_N$ with $\text{Exp}_W = \pi_{TU_N} \circ \varphi_{N} (1- t_N, \cdot ) \circ \phi_W$ and $\text{Exp}_W' = \pi_{TU_N} \circ \varphi_{N} (1- t_N', \cdot ) \circ \phi_W'$. Without loss of generality, taking the maximum of $t_N, t_N'$, we may assume $t_N = t_N'$. Observe that we obtain a diffeomorpism $\phi_W \circ \phi_W'^{-1} \colon O' \rightarrow O$. For each $X \in O'$, there are unique geodesics $\gamma_X' (t) \coloneq \pi_{TU_N} \varphi_N (t, X) \colon [0, 1 - t_N] \rightarrow U_N$ and $\gamma_X (t) \coloneq \pi_{TU_N} \varphi_N (t, \phi_W \circ \phi_W'^{-1} (X)) \colon [0, 1- t_N] \rightarrow U_N$. The geodesics $\gamma_X , \gamma_X'$ lift the same orbifold geodesic arc, since $\text{Exp}_W$ and $\text{Exp}_W'$ are restrictions of orbifold geodesic flows. By Lemma \ref{lem: og:equiv}, for $X\in O'$, there is some $g_X \in G_N$ with $T_{1-t_N} (g_X. \gamma_X) = T_{1-t_N}\gamma'_X$.\\
 The element $g_X$ acts as a Riemannian isometry, mapping geodesics to geodesics, which implies $g_X . \gamma_X (t) = \gamma_X' (t)$, for all $t \in [0, 1-t_N]$. For any non-singular $X \in O'$, the isometry $g_X$ is uniquely determined: To prove this, let $g_X' \in G_N$ be another isometry with $g_X'.\gamma_X = \gamma_X'$. Then,
  \begin{displaymath}
   Tg_X (X) = Tg_X . \varphi_N (0,\phi_W \circ \phi_W'^{-1} (X)) = \varphi_N (0, X) = Tg_X'.\varphi_N (0,\phi_W \circ \phi_W'^{-1} (X)) = Tg_X' (X).
  \end{displaymath}
 Since $X$ is non singular, $T_{\pi_{TU_N} (X)} g_X = T_{\pi_{TU_N}} g_X'$ and by \cite[Lemma 2.10]{follie2003}, $g_X = g_{X}'$ follows. The set $O' \subseteq TU_N$ is an open, connected set. Hence Lemma \ref{lem: ts:nspc} implies that $C \coloneq O' \setminus \Sigma_{TG_N}$ is connected. As we have seen, for each $X \in C$, there is a unique $g_X$ with $g_X . \gamma_X (0) = \gamma_X' (0)$. The set $H_{g_X} \coloneq \setm{c \in C}{g_X .\gamma_c (0) = \gamma_c' (0)} = \setm{c \in C}{g_X . \pi_{TU_N} \varphi_N' (1-t_N, c) = \pi_{TU_N} \varphi_N (1-t_N , \phi_W \circ \phi_W'^{-1}(c))}$ is a closed set by \cite[Theorem 1.5.4]{Engelking1989}. Uniqueness of $g_X$ proves that two such sets $H_g$ and $H_h$ are disjoint if and only if $g\neq h$ holds. Since $G_N$ is finite, the set $H_{g_X}$ is open and closed. By connectedness of $O'\setminus \Sigma_{TG_W}$, there is a unique $\gamma \in G_N$ with 
	\begin{equation}\label{eq: broken}
	 \gamma . \pi_{TU_N} \varphi_N (1-t_N , \cdot )|_{O'\setminus \Sigma_{TG_N}} = \pi_{TU_N} \varphi_N (1-t_N , \cdot ) \circ \phi_W\phi_W'^{-1}|_{O'\setminus \Sigma_{TG_N}}.
	\end{equation}
 The set $O' \setminus \Sigma_{G_N}$ is dense in $O'$ by Newman's Theorem \ref{thm: newman}. Hence, by continuity, \eqref{eq: broken} holds on all of $O'$. As $(\phi_W')^{-1} (O') = W$ by construction, we finally derive $\gamma . \Exp_{W} = \Exp_{W}'$.

 The construction of lifts did not depend on $\xi$, thus we may cover $\Omega$ with a set of orbifold charts $\vV \coloneq \setm{(W_i, G_i, \pi_i)}{i \in I}$ such that on each $(W_i, G_i, \pi_i)$ there exists a local lift $\text{Exp}_{W_i}$ of $\expo$ with respect to $(W_i, G_i , \pi_i)$ and a suitable chart $(U_i, G_i, \psi_i)$. Eliminating charts which occur severalfold, we may assume $(W_i, G_i, \pi_i) \neq (W_j,G_j,\pi_j)$ and $(U_i,G_i,\psi_i) \neq (U_j,G_j,\psi_j) \text{ for } i \neq j$ (by replacing charts $U_i$ with $U_i \times \set{i}$ if necessary). The charts in $\vV$ are compatible since they are contained in $\tT \uU$, their images cover $\Omega$ and we have $\vV \in \tT \uU_\Omega$. Define the atlas $\wW \coloneq \setm{(V,G,\psi) \in \uU}{V = \cod W_i \text{ for some } i \in I}$.\\
 We show that it is possible to construct a quasi-pseudogroup $P$ and a map $\nu$ such that the lifts commute with the changes of charts as in Definition \ref{defn: rep:ofdm}. To this end, consider arbitrary local lifts $\Exp_W$ and $\Exp_W'$ of $\expo$ with respect to the charts $(W,G,\pi)$, $(U, H, \psi)$ and $(W',G',\pi')$, $(U', H', \psi')$, respectively. Furthermore, let $h \in \Ch_\vV$ be a change of charts which induces a commutative diagram: 
	\begin{equation}\begin{gathered}\label{eq: comdiag:expo}
	\begin{xy}
  	\xymatrix{					
	\dom h \ar[r]^{\text{inc}}	\ar[dd]^h	&	W \ar[rrrr]^{\Exp_W} \ar[rd]^\pi 		& 	&		& 	&	U \ar[ld]^\psi	\\
							&				  			& \Omega \ar[rr]^\expo	& & Q	&			\\
	\cod h \ar[r]^{\text{inc}}			&  	W'\ar[rrrr]^{\Exp_{W'}}	\ar[ru]^{\pi'}		&	&		&	& 	U' \ar[lu]_{\psi'}.	 
  }
\end{xy}\end{gathered}
	\end{equation}
 Cover $\text{Exp}_W (\dom h)$ with the domains of suitable changes of charts. Our goal is to restrict $h$ to open subsets such that there are changes of charts which complement the right hand side of \eqref{eq: comdiag:expo} to a commuting triangle. By commutativity of \eqref{eq: comdiag:expo}, for each $X \in \dom h$ there is an embedding of orbifold charts $\lambda_X \in \Ch (U,U')$ such that $\lambda_X (\Exp_W (X)) = \Exp_{W'} (h(X))$. 
 Again let $\phi_W, \phi_{W'}$ denote the diffeomorphisms with $\Exp_W = \pi_{TU} \varphi_U (1-t_N, \cdot) \circ \phi_W$ and $\Exp_W' = \pi_{TU'} \varphi_U' (1-t_N', \cdot) \circ \phi_{W'}$. Since $\varphi_N (t, \phi_W (X))$ is defined for all $t \in [0,1-t_N]$, we deduce from the continuity of the flow that there is some $\ve > t_{N}, t_N'$ such that $\pi_{TU} \varphi_U (1-t , \varphi_U (\ve - t_N, \phi_W (X)) \in \dom \lambda_X$ holds for all $t \in [0 , 1-\ve]$. Define for $Y\in W$ the element $\tilde{Y} \coloneq \varphi_U (\ve -t_N, \phi_W (Y)) \in TU$. Now the open set $\varphi_U^{-1} (T\dom \lambda_X)$ contains $[0, 1- \ve] \times \set{\tilde{X}}$. The Wallace Theorem \cite[3.2.10]{Engelking1989} assures that there is an open neighborhood $\tilde{X} \in \tilde{V} \subseteq TU$ such that $[0,1-\ve] \times \tilde{V} \subseteq \varphi_U^{-1} (T\dom \lambda_X)$. By continuity of $\phi_W$, we can choose an open $G$-stable $X$-neighborhood $V \subseteq (\varphi_U (\ve -t_N, \cdot ) \circ \phi_W )^{-1} (\tilde{V}) \cap \dom h$ with $G_V = G_X$. For each $\tilde{Y}$ with $Y \in V$, the geodesic $\gamma_{\tilde{Y}} (t) \coloneq \pi_{TU} \varphi_U (t, \tilde{Y}), t \in [0, 1-\ve]$ is contained in $\dom \lambda_X$. We obtain two local lifts $\Exp_W'|_{h(V)}$ and $\lambda_X \circ \Exp_W \circ h^{-1}|_{h(V)}$ with respect to the charts $(h (V) , G'_{h(V)}, \pi'|_{h (V)})$ and $(U',H,\psi')$. 
The map $\lambda_X$ is a Riemannian embedding into $U'$ and thus commutes with parallel displacement (see \cite[IV.\ Proposition 2.6]{fdiffgeo1963}) of the open set $\dom \lambda_X$. Hence we derive $T\lambda_X \varphi_U (1- \ve ,\tilde{Y}) = \varphi_{U'} (1-\ve , T\lambda_X (\tilde{Y}))$ for $\tilde{Y} \in \tilde{V}$. In particular, the following holds: 
	\begin{align}
 	\lambda_X \circ \Exp_W \circ h^{-1}|_{h(V)} &= \pi_{TU'} T\lambda_X \varphi_{U}(1-\ve, \cdot ) \varphi_{U}(\ve-t_N, \cdot) \circ \phi_W \circ h^{-1}|_{h(V)} \notag\\
							&= \pi_{TU'} \varphi_{U'} (1-\ve, \cdot ) T\lambda_X \varphi_U (\ve - t_N , \cdot) \circ \phi_W \circ h^{-1}|_{h(V)}. \label{exp:lift}
	\end{align}
 The local lifts $\lambda_X \Exp_W h^{-1}|_{h(V)}$ and $\Exp_W'|_{h(V)}$ are therefore compositions of the bundle projection $\pi_{TU'}$, the geodesic flow on $U'$ and some diffeomorphism. As we have already seen, there is some $\gamma \in H'$ such that $\gamma . \lambda_X \Exp_W h^{-1}|_{h(V)} = \Exp_W'|_{h(V)}$ holds. Replacing $\lambda_X$ with the embedding of orbifold charts $\gamma . \lambda_X$, we derive
	\begin{equation}
	 \lambda_X \circ \Exp_W|_{V} = \Exp_{W'} \circ h|_{V}. \label{eq: diag:comm}
	\end{equation}
 We may thus cover $\dom h$ by open $G$-stable subsets $\setm{W_{X_i}}{i \in I_h}$ such that for each $h_i \coloneq h|_{W_{X_i}}$, there is a change of charts $\lambda_i^h$ which satisfies $ \lambda_i^h \circ \Exp_W|_{V} = \Exp_{W'} \circ h|_{V}$. Repeating this construction for every change of charts in $\Ch_\vV$, we obtain $P \coloneq\setm{h_i}{i \in I_h , \ h \in \Ch_\vV}$. By construction $P$ is a quasi-pseudogroup which generates $\PSI (\vV)$. For each element $f$ of $P$ choose and fix some $h \in \Ch_\vV$ with $f = h_i$ and define the map $\nu \colon P \rightarrow \PSI (\wW), f= h_i \mapsto \lambda_i^h$.\\ By construction, $\widehat{\expo}\coloneq(\expo , \setm{\Exp_W}{(W,G,\pi) \in \vV} , P , \nu)$ satisfies conditions (R1)-(R4a) of Definition \ref{defn: rep:ofdm}. We check condition (R4b), i.e.\ if $g,h \in P$ and $x \in \dom h \cap \dom g$ with $\dom g, \dom h \subseteq U$ and $\germ_x h = \germ_x g$, then $\germ_{\Exp_U (x)} \nu (h) = \germ_{\Exp_U (x)} \nu (g)$. \\
 Let $\dom \nu (h) \subseteq V$ and $\cod \nu (h) \subseteq V'$, where $(V,H, \psi), (V',H',\psi')$ are suitable orbifold charts. By construction we already know $\nu (h) (\Exp_U(x)) = \nu (g) (\Exp_U (x))$. Restricting to an open and $H_{\Exp_U (x)}$-stable subset $\Exp_U (x) \in S_x$ of $\dom \nu (g) \cap \dom \nu (h)$, the changes of charts $\nu (g)$ and $\nu (h)$ restrict to embeddings of orbifold charts. By Proposition \ref{prop: ch:prop}, there is a unique $\gamma \in H'$ such that $\gamma . \nu (g)|_{S_x} = \nu (h)|_{S_x}$. Now $\gamma . \nu (g) (\Exp_U (x)) = \nu(h) (\Exp_U (x)) = \nu (g) (\Exp_U (x))$ implies that $\gamma \in H'_{\nu (g) (S_x)}$ and from Proposition \ref{prop: ch:prop} we obtain some $\delta \in H$ with $\overline{\nu (g)}(\delta) = \gamma$.\\ 
 As $\Exp_U$ is an open map, the intersection $S_x \cap \im \Exp_U (\dom g \cap \dom h)$ is a non-empty open set. It contains at least one non-singular point $y$ by Newman's theorem \ref{thm: newman}. Both maps coincide on $\Exp_W (\dom g \cap \dom h)$, whence 
	\begin{displaymath}
	 \nu (g) (\delta . y) = \gamma . \nu (g) (y) = \nu (h) (y) = \nu (g) (y),
	\end{displaymath}
 which implies $\delta . y = y$. Since $y$ is non-singular, $\delta = \id_V$ follows. The mapping $\overline{\nu (g)}$ is a group homomorphism, from which we deduce $\gamma = \id_{V'}$. In conclusion, $\nu (g)|_{S_x} = \nu (h)|_{S_x}$ holds, whence their germs agree, proving property (R4b). The above shows that there is locally only one choice for $\nu (g)$. From this observation, one deduces that properties (R4c)-(R4d) are also valid for $\widehat{\expo}$.\\[1em]
 We have thus constructed a charted map 
    \begin{displaymath}
     \widehat{\expo} = (\expo , \set{\Exp_W}_{(W,G,\pi)\in \vV} , [P, \nu]) \in \Orb{\vV, \wW}
    \end{displaymath}
 for the range family $\wW \in \uU$ as defined above. To finish the proof, we need to check that every other choice of lifts yields a charted orbifold map which is equivalent to $\widehat{\expo}$.\\
 Let $\widetilde{\expo} = (\expo , \setm{E_{W'}}{(W',G',\psi') \in \vV'}, [P',\nu'])$ be another charted orbifold map whose lifts are constructed as above. Arguing as before, for each lift $\Exp_W$, we may cover $\im \Exp_W$ with the domains of embeddings $\mu_W^i, \ i \in I$ of orbifold charts such that:
 \begin{compactenum}
  \item $\dom \mu_W^i \neq \dom \mu_W^j$ for each $i \neq j$,
  \item for each $i$, there is a lift $E_{W'_i}$ of $\expo'$ and an embedding of orbifold charts $\lambda_W^i$ such that $\Exp_W (\dom \lambda_W^i)  \subseteq \dom \mu^i_W$ and $\mu_W^i \Exp_W|_{\dom \lambda_W^i} = E_{W'_i}\lambda_W^i$.
 \end{compactenum}
 Repeating this argument for each chart in $\vV$, we obtain an orbifold atlas $\aA$ of charts for $\Omega$ and a family $\fF$ of orbifold charts for $Q$. In particular, for each chart $A \in \aA$, there is a chart in $\fF$ together with two pairs of embeddings of orbifold charts: The first pair $(\iota_A^1 , \iota_A^2)$ being the canonical inclusion into $\dom \Exp_W $, respectively $\cod \Exp_W$ for a suitable lift of $\expo$, while the second pair is given by the embeddings $(\lambda_A, \mu_A)$ constructed above. It is now easy to check that the data $(\aA , \fF , (\iota_A^1, \iota_A^2)_{A \in \aA})$ and $(\aA , \fF , (\lambda_A, \mu_A)_{A \in \aA})$ satisfy the hypothesis of Lemma \ref{lemdef: ind:ofdm}. By construction, the induced lifts of $\widehat{\expo}$ and $\widetilde{\expo}$ coincide. In particular, the induced lifts satisfy an identity as in \eqref{exp:lift}, i.e.\ by construction they are given as the composition of geodesic flows, changes of charts and bundle projection of manifolds. An argument as above shows that locally there is just one choice for the change of charts in the image of $\nu$. Local uniqueness of the changes of charts relating the lifts thus forces $\widehat{\expo} \sim \widetilde{\expo}$ (cf. Definition \ref{defn: eq:rofdm}). Hence $[\widehat{\expo}] = [\widetilde{\expo}]$ follows and we abbreviate this unique map of orbifolds as $[\expo]$.  
\end{proof}

The above proof reveals several useful properties of the lifts for $\expo$, which we collect in the following 

\begin{rem}\label{rem: lift:expo} \no{rem: lift:expo}
 \begin{compactenum}
  \item The proof of Proposition \ref{prop: exp:ofdm} shows that arbitrary sets of lifts (which are given as lifts of orbifold geodesic flows evaluated at $1$) for $\expo$, where no two are defined on the same chart, may be complemented to a family of local lifts which satisfy (R2) of Definition \ref{defn: rep:ofdm}. Each of these families then induces a representative of $[\expo]$.
  \item The families of lifts we constructed in Proposition \ref{prop: exp:ofdm} have the additional property that for each $\Exp_W \colon (W,G_W,\pi) \rightarrow (U_W,G_{U_W} , \psi)$, there is an orbifold chart $(V,H,\varphi)$ such that $W \subseteq TV$ is an $H$-stable subset which is $G_W$-invariant.
 \end{compactenum}
\end{rem}

\thispagestyle{empty}
\section{Lie Group Structure on the Orbifold Diffeomorphism Group}\label{sect: lgp:str}\no{sect: lgp:str}
\setcounter{subsubsection}{0}

Throughout this section, we assume that $(Q, \uU , \rho)$ is a smooth Riemannian orbifold. We construct a Lie group structure on $\Difforb{Q,\uU}$ by an application of the construction principle outlined in Proposition \ref{prop: Lgp:locd}. To this end, the subgroup of all compactly supported orbifold diffeomorphisms will be turned into a Lie group.

\subsection{Lie group structure on \texorpdfstring{$\Diffo{Q,\uU}$}{an identity neighborhood}}\label{sect: lgp:ident} \no{sect: lgp:ident}

It turns out that our approach needs a framework, i.e. an orbifold atlas together with a collection of local data, which we fix now. Based on this preliminary work, we construct a locally convex Lie group structure modeled on $\Osc{Q}$ for the subgroup $\Difforb{Q,\uU}_0 \subseteq \Difforb{Q,\uU}$. This group is generated by elements in $\Difforb{Q,\uU}$ suitably close to the identity. In Section \ref{sect: lgp:glob}, this Lie group becomes the identity component for the Lie group $\Difforb{Q,\uU}$. 
 \begin{con}\label{con: sp:ofdat}\no{con: sp:ofdat}
 \begin{compactitem}
  \item[I.] Choose for each connected component $C \subseteq Q$ some $z_C \in C$. As $Q$ is locally path connected, each component of $Q$ is open. Hence $\setm{z_C}{C \subseteq Q, \text{ connected component}} $ is a discrete and closed subset. Combining Proposition \ref{prop: atl:lfpts} with Lemma \ref{lem: fin:al}, we may choose orbifold atlases $\aA, \bB \in \uU$ with the following properties: 
	\begin{compactenum}
	 \item the atlases $\aA = \setm{(U_i,G_i,\psi_i)}{i \in I}$ and $\bB = \setm{(W_j,H_j,\varphi_j)}{j \in J}$ are locally finite,
	 \item each chart in $\aA, \bB$ is relatively compact (i.e.\ its image in $Q$ is relatively compact),
	 \item For each connected component $C \subseteq Q$, there are unique $i_C \in I, j_C \in J$ with $z_C \in \psi_{i_C} (U_{i_C})$ (resp. $z_C \in \varphi_{j_C} (W_{j_C})$), 
	 \item $\aA$ is a refinement of $\bB$ and there is a map $\alpha \colon I \rightarrow J$ such that each $i \in I$ satisfies: 
		\begin{compactitem}
		 \item[i)] $\overline{U_i} \subseteq W_{\alpha (i)}$ and the canonical inclusion of sets is an embedding of orbifold charts, implying
		 $G_i \subseteq H_{\alpha (i)}$ and $\psi_i = \varphi_{\alpha (i)}|_{U_i}$,
 		 \item[ii)] $\alpha (i_C) = j_C$,	 
		 \item[iii)] $\alpha^{-1} (j)$ is finite for each $j \in J$.
		\end{compactitem}
	\end{compactenum}
 \item[II.] For each $i \in I$, the set $\overline{U_i} \subseteq W_{\alpha (i)}$ is compact and connected. By local compactness and local connectedness, there is a relatively compact connected open set $\overline{U}_i \subseteq O_i \subseteq W_{\alpha (i)}$. The set $H_{\alpha (i)} . O_i$ is open, $H_{\alpha (i)}$-invariant and $\overline{U_i}$ is a connected subset of $O_i \subseteq H_{\alpha (i)}.O_i$. Thus $\overline{U_i}$ is contained in a connected component of $H_{\alpha (i)} . O_i$. Replacing $O_i$ with this component, without loss of generality $O_i$ is an open, relatively compact, $H_{\alpha (i)}$-stable subset. Notice that $G_i \subseteq H_{\alpha (i), O_i}$ holds by construction.
 \item[III.] For each $j \in J$, define a compact, $H_j$-invariant subset $\kK_j \coloneq H_j . \bigcup_{i \in \alpha^{-1} (j)} \overline{O_i}$. Apply Lemma \ref{lem: oemb:cov} with respect to the family of compact sets $(\kK_j)_{j\in J}$ and the atlas $\bB$. There is a cover for each $\kK_j$ by a finite set $\zZ_j \coloneq \setm{Z_j^k}{1\leq k \leq N_j}$ of open $H_{j}$-stable sets such that: for each member of $\zZ_j$, there is a finite family of embeddings $(\lambda_{jh}^k \colon Z_j^k \rightarrow W_h)_{h \in Z(j,k)}$ with properties as in Lemma \ref{lem: oemb:cov}. By Part (c) of Lemma \ref{lem: oemb:cov}, each $Z_j^k$ is relatively compact and the embedding $\lambda_{jh}^k$ is the restriction of an embedding $\hat{\lambda}_{jh}^k$ whose domain contains $\overline{Z_j^k}$. 
 \item[IV.] Consider the open submanifold $\kK_j^\circ $, 
 which is $\sigma$-compact as an open subset of the second countable locally compact manifold $W_{j}$ (cf.\ the proof of Proposition \ref{prop: atl:lfpts} (d)). By Lemma \ref{lem: mfd:atl}, we may cover each $\kK_j^\circ, j \in J$ with a countable family $\set{(V^k_{5, j}, \kappa^j_k)}_{1 \leq k \leq l_j}$, $l_j \in \NN_0 \cup \set{\infty}$ of manifold charts such that the cover is locally finite and subordinate to the open cover $\setm{Z^j_{k} \cap \kK_j^\circ}{1\leq k\leq N_j}$ of $\kK_j^\circ$. Furthermore, these charts satisfy $\kappa_k^j (V_{5,j}^k) = B_5 (0)$ and the families $V_{r,j}^k \coloneq (\kappa_k^j)^{-1} (B_r (0))$, $1\leq k \leq j_l$ cover $\kK_j^\circ$ for each $r \in [1,5]$.\\
 Since $H_{\alpha (i)}$ is finite, the set $H_{\alpha (i)}. \overline{U}_i \subseteq \kK^\circ_{\alpha (i)}$ is compact. The atlas $\set{(V^k_{5, j}, \kappa^j_k)}_{1 \leq k \leq l_j}$ is locally finite, whence there is a finite subset $\fF_{5} (H_{\alpha (i)}.\overline{U_i})$\glsadd{F5} such that $V^k_{5,\alpha (i)} \cap H_{\alpha (i)}. \overline{U}_i \neq \emptyset$ if and only if the chart $(V^k_{5,\alpha (i)}, \kappa_k^j)$ belongs to $\fF_{5} (H_{\alpha (i)}. \overline{U_i})$.
 We define open sets \glsadd{omega1}
	\begin{displaymath}
	 \Omega_{r,i} \coloneq \bigcup_{(V_{5,\alpha (i)}^n, \kappa_{n}^{\alpha (i)}) \in \fF_5 (H_{\alpha (i)} .\overline{U_i})} V_{r,\alpha (i)}^n,\quad r \in [1,5]
	\end{displaymath}
 and compact sets $K_{5,i} \coloneq \overline{\Omega_{5,i}}$. There is a finite subset $\fF_5 (K_{5,i})$ such that a chart belongs to $\fF_{5} (K_{5,i})$ if and only if $V^k_{5,\alpha (i)} \cap H_{\alpha (i)}. K_{5,i} \neq \emptyset$ holds. Observe that $H_{\alpha (i)}.\overline{U_i} \subseteq \Omega_{1,i}$ is satisfied.  
\item[V.] Let $\rho_{j}$ be the Riemannian metric on $W_{j}$ and $\exp_{W_j} \colon D_{j} \rightarrow W_{j}$ the associated Riemannian exponential map. By compactness of $\kK_{j}$ and Lemma \ref{lem: rg:rmb}, there are constants $s_j > 0$ for $j \in J$ such that: The closure of $\hat{O}_j \coloneq \bigcup_{x \in \kK_j^\circ} B_{\rho_j} (0_x , s_j) \subseteq TW_{j}$
 is contained in $D_{j}$ and $\exp_{W_{j}}$ restricts to a diffeomorphism on $T_x W_{j} \cap \hat{O}_j$ for each $x \in \kK_j^\circ$. Moreover, $\overline{\Omega_{\frac{5}{4},K_{5,i}}}$ is compact for $i \in I$ and $\alpha^{-1}(j)$ is finite for $j \in J$. Shrinking the constants $s_j$, we can achieve that for each $i \in I$ and $x \in \Omega_{\frac{5}{4} ,K_{5,i}}$ the identity $\exp_{W_{\alpha (i)}} (B_{\rho_{\alpha (i)}} (0_x , s_{\alpha (i)})) \subseteq \Omega_{2,K_{5,i}}$ is satisfied. Since $\hat{\lambda}_{jh}^k (\overline{Z_j^k})$ is compact, Lemma \ref{lem: rg:rmb} yields a constant $0 < S_{jk} < \min \setm{s_h}{h \in Z(j,k)}$ such that $\exp_{W_h}$ restricts to a diffeomorphism on 
	\begin{displaymath}
	 T \hat{\lambda}_{jh}^k \left(B_{\rho_j} (0_x, S_{jk})\right) \subseteq T_{\hat{\lambda}_{jh}^k (x)} W_h ,\quad x \in \overline{Z_j^k}.
	\end{displaymath}
 Furthermore, since changes of charts are Riemannian embeddings, by choice of $S_{jk}$ 
    \begin{displaymath}
     T \lambda_{jh}^k \left(B_{\rho_j} (0_x, S_{jk})\right) \subseteq B_{\rho_h} (0_{\lambda_{jh}^k (x)}, s_h)
    \end{displaymath}
 holds for $x \in \dom \lambda_{jh}^k$. For each $j \in J$, we define $S_j \coloneq \min \setm{S_{jk}}{1\leq k \leq N_j}$.  
 The set $\fF_5 (K_{5,i})$ is finite and for each chart $(V_{5,\alpha (i)}^k, \kappa_k^{\alpha (i)}) \in \fF_5 (K_{5,i})$ the set $\bigcup_{x \in V_{3,\alpha (i)}^k} B_{\rho_{\alpha (i)}} (0_x, S_{\alpha (i)})$ is a neighborhood of the zero-section on the compact set $\overline{V_{2,\alpha (i)}^k}$. Hence the Wallace Lemma \cite[3.2.10]{Engelking1989} yields a constant $R_i >0$ with   
	\begin{displaymath}
	 B_2 (0) \times B_{R_i} (0) \subseteq T\kappa^{\alpha (i)}_k \left(\bigcup_{x \in V_{2,\alpha (i)}^k} B_{\rho_{\alpha (i)}} (0_x, S_{\alpha (i)})\right) \quad \forall (V_{5,\alpha (i)}^k, \kappa_k^{\alpha (i)}) \in \fF_5 (K_{5,i}).
	\end{displaymath}
 \end{compactitem}
 \end{con}
 For the rest of this section, we fix the data constructed in \ref{con: sp:ofdat} and use the symbols without further explanation. The next lemma is a rather technical statement. It is the first step in constructing orbifold diffeomorphisms using the Riemannian orbifold exponential map.

\begin{lem}\label{lem: eos:loc} \no{lem: eos:loc}
 Consider $(U_i , G_i ,\psi_i) \in \aA$ and for an orbisection $[\hat{\sigma} ] \in \Os{Q}$ denote by $\sigma_{\alpha (i)}$ its canonical lift on $W_{\alpha (i)}$ and by $\sigma$ its underlying continuous map. There exists an open neighborhood $\nN_i \subseteq \vect{W_{\alpha (i)}}$ of the form $\nN_i = (\res^{W_{\alpha (i)}}_{\Omega_{5,i}})^{-1} (\nN_i^{\Omega_{5,i}})$ of $0_{\alpha (i)}$ such that $\sigma_{\alpha (i)} \in \nN_i$ implies the following: 
\begin{compactitem}
        \item[i.] $\overline{\psi_i (U_i)} \subseteq \sigma^{-1}(\Omega)$, where $\Omega$ is the domain of $\expo$, 
	\item[ii.] $[\hat{E}^{\sigma}]|_{\psi_i (U_i)} \coloneq [\expo] \circ [\hat{\sigma}]|^\Omega_{\psi_i (U_i)}$ induces a diffeomorphism of orbifolds onto its image, 
	\item[iii.] $\sigma_{\alpha (i)} (\overline{\Omega_{2,i}}) \subseteq \hat{O}_{\alpha (i)}$ holds for $\hat{O}_{\alpha (i)}$ as in Construction \ref{con: sp:ofdat} V.
\end{compactitem}
 for some zero-neighborhood $\nN_i^{\Omega_{5,i}} \subseteq \vect{\Omega_{5,i}}$.
\end{lem}

\begin{proof}
 The set $O_i \subseteq \kK_{\alpha (i)}^\circ$ is open and $H_{\alpha (i)}$-stable, whence an $H_{\alpha (i)}$-stable open subset is given by $TO_i \cap \hat{O}_{\alpha (i)} \subseteq D_{\alpha (i)}$. We obtain an orbifold chart $(TO_i \cap \hat{O}_{\alpha (i)},H_{\alpha (i), TO_i \cap \hat{O}_{\alpha (i)}}, T\varphi_{\alpha (i)}|_{TO_i \cap \hat{O}_{\alpha (i)}})$ together with the lift $\Exp_{TO_i \cap \hat{O}_{\alpha (i)}} \coloneq \exp_{W_{\alpha (i)}}|_{TO_i \cap \hat{O}_{\alpha (i)}} \colon TO_i \cap \hat{O}_{\alpha (i)} \rightarrow W_{\alpha (i)}$ of $\expo$. By Remark \ref{rem: lift:expo} (a), there is a representative $\widehat{\expo} \in \Orb{\vV , \wW}$ of $\expo$ such that $\Exp_{TO_i \cap \hat{O}_{\alpha (i)}} $ is contained in the family of local lifts of $\widehat{\expo}$. 
 Notice that $\psi_i (U_i) \subseteq Q$ is an open subset,  whose inclusion $\iota_{\psi_i (U_i)}$ induces an open suborbifold structure (see Definition \ref{defn: sofd:op}). Consider an orbisection $[\hat{\sigma}]$ with $\im \sigma|_{\psi_i (U_i)} \subseteq \Omega$. Definitions \ref{defn: sofd:op} and \ref{defn: osup} together with Proposition \ref{prop: ocat:cp} imply that there is a well-defined map of orbifolds $\left[\hat{E}^\sigma|_{\psi_i (U_i)}\right] \coloneq [\expo] \circ [\hat{\sigma}]|_{\psi_i (U_i)}^{\Omega}$. 
 Now, we proceed in several steps:
	
 \paragraph{Step 1:} 
 Apply Lemma \ref{lem: vect:emb} to the family $\fF_5 (H_{\alpha (i)}.\overline{U_i})$ to obtain an open zero-neighborhood $N_i^{\Omega_{5,i}} \subseteq \vect{\Omega_{5,i}}$ (playing the role of $E_{5,K}$ in the lemma). Define $N_i \coloneq (\res^{W_{\alpha (i)}}_{\Omega_{5,i}})^{-1} (N_i^{\Omega_{5,i}}) \subseteq \vect{W_{\alpha (i)}}$ and observe that $0_{\alpha (i)} \in N_i$ and the following conditions hold: For each $X \in N_i$, the map $\exp_{W_{\alpha (i)}} \circ X|_{\Omega_{2,i}}$ is an \'{e}tale embedding into $W_{\alpha (i)}$. The set $\overline{\Omega_{2,i}} \subseteq \Omega_{5,i} \subseteq \kK_{\alpha (i)}^\circ$ is compact, which allows the construction of a $C^0$-neighborhood of the zero section $P_{1,i} \subseteq \vect{\Omega_{5,i}}$ such that $X \in P_{1,i}$ implies $X (\overline{\Omega_{2,i}}) \subseteq \hat{O}_{\alpha (i)}$. Set $\nN_i^{\Omega_{5,i}} \coloneq N_i^{\Omega_{5,i}} \cap P_{1,i}$ and $\nN_i \coloneq (\res^{W_{\alpha (i)}}_{\Omega_{5,i}})^{-1} (N_i^{\Omega_{5,i}} \cap P_{1,i})$. Each vector field in $\nN_i$ satisfies iii.\ and $\nN_i$ is a preimage as required. By construction, $\overline{\psi_i (U_i)} = \varphi_{\alpha (i)} (\overline{U_i}) \subseteq \varphi_{\alpha (i)} (O_i)$ holds and $\Exp_{TO_i \cap \hat{O}_{\alpha (i)}}$ is a lift of $\expo$, whence i.\ follows from property iii.
 In addition, if $\sigma_{\alpha (i)} \in \nN_i$ then the map $\exp_{W_{\alpha (i)}}\circ \sigma_{\alpha (i)}|_{H_{\alpha (i)}.U_i}$ is an \'{e}tale embedding. Specializing to $U_i$, the map $e^{\sigma_{i}} \coloneq \Exp_{TO_i \cap \hat{O}_{\alpha (i)}} \circ \sigma_{\alpha (i)}|_{U_i} = \Exp_{TO_i \cap \hat{O}_{\alpha (i)}} \circ \sigma_i$ is a \'{e}tale embedding, where $\sigma_i$ is the canonical lift of $[\hat{\sigma}]$ on $(U_i,G_i,\psi_i)$. From now on, consider $[\hat{\sigma}]\in \Os{Q}$ such that $ \sigma_{\alpha (i)} \in \nN_i$.

\paragraph{Step 2:} \emph{The map $e^{\sigma_{i}}$ is equivariant with respect to the inclusion $\nu \colon G_i \hookrightarrow H_{\alpha (i)}$}:
 Consider an $H_{\alpha (i)}$-invariant subset $R \subseteq \Omega_{2,i}$. We claim that $\exp_{W_{\alpha (i)}} \sigma_{\alpha (i)}|_{R}$ is equivariant with respect to $H_{\alpha (i)}$. If this is correct, then $e^{\sigma_{i}}$ commutes with any $\delta \in H_{\alpha (i), U_i} =G_i$, as $H_{\alpha (i)}. U_i \subseteq \Omega_{2,i}$ is invariant. To prove the claim, let $\delta \in H_{\alpha (i)}$ be arbitrary and $x \in R$. As $\delta . x \in R \subseteq \Omega_{2,i}$ holds, $\sigma_{\alpha (i)}$ is a canonical lift and $H_{\alpha (i)}$ acts by Riemannian isometries, we compute: $\exp_{W_{\alpha (i)}} \sigma_{\alpha (i)} (\delta . x) = \exp_{W_{\alpha (i)}} T\delta \sigma_{\alpha (i)} (x) = \delta . \exp_{W_{\alpha (i)}} \sigma_{\alpha (i)} (x)$, thus proving the claim.
 The map $e^{\sigma_{i}}$ is a local lift of $E^\sigma|_{\psi_i (U_i)} \coloneq (\expo \circ \sigma|_{\psi_i (U_i)})|^{\im \varphi_{\alpha (i)}}$.
\paragraph{Step 3:} \emph{The set $\im e^{\sigma_{i}}$ is $H_{\alpha (i)}$-stable with $H_{\alpha (i)}.\im e^{\sigma_{i}} \subseteq \Omega_{2,i}$}: Consider $\delta \in H_{\alpha (i)}$ such that $\delta . \im e^{\sigma_{i}} \cap \im e^{\sigma_{i}} \neq \emptyset$. For $x,y \in U_i$ with $e^{\sigma_{i}}(x) = \delta . e^{\sigma_{i}} (y)$, one obtains
	\begin{displaymath}
	 \exp_{W_\alpha (i)} \circ \sigma_{\alpha (i)} (x) = e^{\sigma_{i}} (x) = \delta . e^{\sigma_{i}} (y) = \exp_{W_{\alpha (i)}} \sigma_{\alpha (i)} (\delta . y).
	\end{displaymath}
 From Step 1, we conclude $x = \delta . y$, since on $H_{\alpha (i)} . U_i \subseteq \Omega_{2,i}$ the map $\exp_{W_{\alpha (i)}} \circ \sigma_{\alpha (i)}$ is a \'{e}tale embedding. By $H_{\alpha (i)}$-stability of $U_i$, $\delta \in G_i$ holds, whence $\delta . \im e^{\sigma_{i}} = \im e^{\sigma_{i}}$. This proves the $H_{\alpha (i)}$-stability of $\im e^{\sigma_{i}}$ and $G_{\im e^{\sigma_{i}}} = G_i$.\\ 
 The canonical lift $\sigma_{\alpha (i)}$ is contained in $\nN_i$. By construction of $\Omega_{1,i}$ (cf.\ Lemma \ref{lem: vect:emb}), the equivariance of this map implies: 
    \begin{displaymath}
     H_{\alpha (i)} . \im e^{\sigma_i} = \exp_{W_{\alpha (i)}} \sigma_i (H_{\alpha (i)}.U_i) \subseteq \exp_{W_{\alpha (i)}} \sigma_i (\Omega_{1,i}) \subseteq \Omega_{2,i}.
    \end{displaymath}
 \paragraph{Step 4:}  \emph{$E^\sigma|_{\psi_i (U_i)}$ is injective and a homeomorphism onto its open image}: Consider $x,y \in \psi_i (U_i)$ with $E^\sigma|_{\psi_i (U_i)} (x) = E^\sigma|_{\psi_i (U_i)} (y)$ and choose preimages $z_x \in \psi_{i}^{-1} (x) , z_y \in \psi_{i}^{-1}(y)$ of $x$ respectively $y$ in $U_i$. Since $e^{\sigma_i}$ is a lift of $E^\sigma|_{\psi_i (U_i)}$, there exists $\delta \in H_{\alpha (i)}$ such that $e^{\sigma_{i}} (z_x) = \delta . e^{\sigma_{i}} (z_y)$. By Step 3, we must have $\delta \in G_i$. Since $e^{\sigma_{i}}$ is an embedding, equivariance of this map yields $\delta . z_y = z_x$. Both points are in the same orbit, which forces $x$ and $y$ to coincide. Hence $E^\sigma|_{\tilde{U}_i}$ is injective.\\
 The local lift $e^{\sigma_{i}}$ is a \'{e}tale embedding and the maps of orbifold charts are continuous and open. For any open subset $S \subseteq \psi_i (U_i)$, $E^{\sigma}|_{\psi_i (U_i)} (S) = \varphi_{\alpha (i)} \circ e^{\sigma_{i}} \circ \psi_i^{-1} (S)$ is an open set. In conclusion, $E^\sigma|_{\psi_i (U_i)}$ is an open map, whose image is open in $Q$. In particular, $\im E^{\sigma}|_{\psi_i (U_i)}$ is an open suborbifold of $Q$. An atlas for $\im E^{\sigma}|_{\psi_i (U_i)}$ is given by $\set{(\im e^{\sigma_{i}}, G_i , \varphi_{\alpha (i)}|_{ \im e^{\sigma_{i}}})}$. \\
 Since composition in $\ORB$ is well-defined, a representative of $[\expo] \circ [\hat{\sigma}]|_{\psi _i (U_i)}^\Omega$, corestricted to $\im \varphi_{\alpha (i)}$ is given by $\hat{E}^\sigma|_{\psi_i (U_i)} = (E^\sigma|_{\psi_i (U_i)} , e^{\sigma_{i}}, G_i, \nu) \in \Orb{\set{(U_i,G_i,\psi_i)} , \set{(W_{\alpha (i)}, H_{\alpha (i)}, \varphi_{\alpha (i)})}}$. The map $E^\sigma|_{\psi_i (U_i)}$ is a homeomorphism mapping the open suborbifold $\psi_i (U_i)$ of $Q$ onto an open suborbifold such that the local lift of $E^\sigma|_{\psi_i (U_i)}$ is a diffeomorphism onto its (open) image. Proposition \ref{prop: ofd:iso} assures that $[\hat{E}^\sigma|_{\psi_i (U_i)}]$ is a diffeomorphism of orbifolds.
\end{proof}

\begin{setup}\label{setup: cC} \no{setup: cC}
 Later on, we shall apply patched mapping techniques (cf.\ Section \ref{sect: sections}) to prove the smoothness of several maps. To do so, we have to define an orbifold atlas, where charts may occur repeatedly: Let $\cC \coloneq \setm{(W_{\alpha (i)}, H_{\alpha (i)}, \varphi_{\alpha (i)})}{i \in I}$ be the orbifold atlas which arises from $\bB$ by collecting a different copy of $(W_j, H_j,\varphi_j)\in \bB$ for each $i \in \alpha^{-1}(j)$. Observe that this atlas is locally finite and each chart is relatively compact, as $\alpha^{-1} (j)$ is finite and $\bB$ is locally finite with relatively compact charts.
\end{setup}

\begin{prop}\label{prop: eos:imdiff} \no{prop: eos:imdiff}
 There are open zero-neighborhoods $\nN_i \subseteq \vect{W_{\alpha (i)}},\ i \in I$ which generate an open zero-neighborhood $\nN \subseteq \Osc{Q}$ such that each $[\hat{\sigma}] \in \nN$ induces an orbifold diffeomorphism $[\hat{E}^\sigma] \coloneq [\expo] \circ [\hat{\sigma}]|^\Omega \in \Difforb{Q,\uU}$. 
\end{prop}

\begin{proof}
 For each $i \in I$, construct via Lemma \ref{lem: eos:loc} a neighborhood $\nN_i \subseteq \vect{W_{\alpha (i)}}$. The construction shows that for each $[\hat{\sigma}]$ with $\sigma_{\alpha (i)} \in \nN_i$, the map $E^\sigma|_{\psi_i (U_i)}$ is an embedding of the open suborbifold $\psi_i (U_i)$. By definition of the direct sum topology, the box $\bigoplus_{i \in I} \nN_i \coloneq \left(\prod_{i \in I} \nN_i \right) \bigcap \bigoplus_{i \in I} \vect{W_{\alpha (i)}}$ is an open subset of $\bigoplus_{i \in I} \vect{W_{\alpha (i)}}$ (cf. \cite[4.3]{jarchow1980} respectively \cite[Proposition 7.1]{hg2003} for a proof).\\ Using the atlas $\cC$ introduced above, we define the set
	\begin{equation}\label{eq: eos:nb}
	 \nN \coloneq \Lambda_{\cC}^{-1} \left(\bigoplus_{i\in I} \nN_i \right),
	\end{equation}
 which is open in the c.s.\ orbisection topology by Lemma \ref{lem: ost:insat}. A combination of Definition \ref{defn: cores} and Remark \ref{rem: lift:expo} (a) shows that each $[\hat{\sigma}]$ contained in $\nN$ induces a well-defined map of orbifolds $[\hat{E}^\sigma] \coloneq [\expo] \circ [\hat{\sigma}]|^\Omega$ such that $E^{\sigma} \coloneq \expo \circ \sigma \colon Q \rightarrow Q$ is a local homeomorphism. In particular, $E^\sigma|_{\psi_i (U_i)}$ is an open embedding for each $i \in I$. Let $\widehat{\expo}$ be the representative of the Riemannian orbifold exponential map obtained from the family $(\text{Exp}_{TO_i \cap \hat{O}_{\alpha (i)}})_{i \in I}$ by Remark \ref{rem: lift:expo} (a). 
Then the domain atlas $\eE$ of $\widehat{\expo}$ contains the family $\set{(TO_i \cap \hat{O}_{\alpha (i)},H_{\alpha (i), TO_i \cap \hat{O}_{\alpha (i)}}, T\varphi_{\alpha (i)}|_{TO_i \cap \hat{O}_{\alpha (i)}})}_{i \in I}$  and for each $[\hat{\sigma}] \in\nN$, the canonical lifts $\sigma_i$ satisfy $\im \sigma_i \subseteq TO_i \cap \hat{O}_{\alpha (i)}$ for $i \in I$. Hence there is a representative $\hat{\sigma}|^\Omega \in \Orb{\aA, \eE}$ of $[\hat{\sigma}]|^\Omega$ whose lift on $(U_i,G_i,\psi_i), i \in I$ is just $\sigma_i|^{TO_i \cap \hat{O}_{\alpha(i)}}$. As composition in $\ORB$ is well-defined, we obtain $[\widehat{\expo} \circ \hat{\sigma}|^\Omega] = [\expo] \circ [\hat{\sigma}]|^\Omega$. Thus the lifts constructed in Lemma \ref{lem: eos:loc} yield a representative $\hat{E}^\sigma \coloneq \widehat{\expo} \circ \hat{\sigma} = (E^\sigma ,\set{e^{\sigma_{i}}}_{i\in I}, P, \nu) \in \Orb{\aA , \cC}$. Here each lift $e^{\sigma_{i}}$ is an \'{e}tale embedding and $(P,\nu)$ is obtained by an application of Construction \ref{con: cp:chom}. The image of such a lift is an orbifold chart $(\im e^{\sigma_{i}}, G_i, \varphi_{\alpha (i)}|_{\im e^{\sigma_{i}}})$.\\
 We have to check that $E^\sigma$ is surjective and injective for every $[\hat{\sigma}]\in \nN$ to prove the assertion. Reviewing the construction of $\nN_i$, the map $E^\sigma$ maps $\psi_i (U_i)$ into $\varphi_{\alpha (i)} (W_{\alpha (i)})$. Every orbifold chart is a connected set, whence its image is contained in a connected component of $Q$. Thus $E^\sigma$ maps every connected component of $Q$ into itself. In conclusion, it suffices to prove that the restriction of $E^\sigma$ to each component is bijective, whence we can assume that $Q$ is connected.\\ 
 As a first step, we show that for every orbisection $[\hat{\sigma}] \in \nN$ the map $E^\sigma$ is a proper map. To this end consider an arbitrary compact subset $L \subseteq Q$. The atlas $\bB$ is locally finite and thus $L$ meets only finitely many of the sets $\varphi_j (W_j), j \in J$, say $L \subseteq \bigcup_{r=1}^n \varphi_{j_r} (W_{j_r})$ and $L \cap \varphi_j (W_j) = \emptyset$ for all $j \in J \setminus \set{j_1, \ldots , j_n}$. For $[\hat{\sigma}] \in \nN$, we have $E^\sigma (\overline{\psi_i (U_i)}) \subseteq \varphi_{\alpha (i)} (W_{\alpha (i)})$. The closed set $(E^\sigma)^{-1} (L)$ is thus contained in 
	\begin{equation} \label{eq: proper}
	 (E^\sigma)^{-1}(L) \subseteq \bigcup_{r=1}^n \bigcup_{i \in \alpha^{-1} (j_r)} \overline{\psi_i (U_{i})}.
	\end{equation}
 By construction \ref{con: sp:ofdat}, each $\alpha^{-1} (j_r)$ is a finite set. Hence $(E^\sigma)^{-1} (L)$ is compact as a closed subset of a union of finitely many compact sets. Since $L$ was arbitrary, $E^\sigma$ is a proper map (cf.\ \cite[{Ch. I, \textsection 10, No. 3, Proposition 7}]{bourbaki1971}). 
 Combining the facts that $Q$ is locally compact by Proposition \ref{prop: ofd:prop} and $E^\sigma$ is a proper map, $E^\sigma$ is a closed map (cf.\ \cite[{Ch. I, \textsection 10, No.\ 2, Theorem 1}]{bourbaki1971}). The image of $E^\sigma$ is an open and closed set, since images of local homeomorphisms are open. But $Q$ is connected and thus $E^\sigma$ is surjective.\\ 
 The map $E^\sigma$ is a proper, surjective local homeomorphism of connected and path-connected locally compact spaces. Summing up, $E^\sigma$ is a covering of $Q$ onto $Q$ by \cite[Theorem 4.22]{forster1981}. Recall \ref{con: sp:ofdat} I.\ (c): There is some $z_Q \in Q$ such that $z_Q$ is contained in a unique pair of orbifold charts $\left( (U_{z_Q}, G_{z_Q}, \psi_{z_Q}) , (W_{z_Q}, H_{z_Q}, \varphi_{z_Q})\right) \in \aA \times \bB$. Since $E^\sigma (\psi_i (U_i)) \subseteq \varphi_{\alpha (i)} (W_{\alpha (i)})$ and $z_Q$ is not contained in any $\varphi_{j} (W_{j})$ except for $j = j_Q$ by \ref{con: sp:ofdat}, we derive from \eqref{eq: proper}: $\lvert(E^\sigma)^{-1} (z_Q) \rvert = 1$. The number of sheets of a covering is an invariant for the connected space $Q$ (cf.\ \cite[Theorem 4.16]{forster1981}), whence $E^\sigma$ is injective. \\
 In conclusion we have constructed a charted orbifold map $\hat{E^\sigma}$ such that $E^\sigma$ is a continuous, closed bijective map (i.e.\ a homeomorphism by \cite[III.\ Theorem 12.2]{dugun1966}) and each lift $e^{\sigma_i}, (V_i , G_i , \psi_i) \in \vV$ is a \'{e}tale embedding. Each lift is a local diffeomorphism, whence Proposition \ref{prop: ofd:iso} implies that $\hat{E^\sigma}$ is a representative of an orbifold diffeomorphism $[\hat{E^\sigma}] = [\expo]\circ [\hat{\sigma}]|^\Omega$.
\end{proof}

The mapping taking an orbisection from the zero-neighborhoods $\nN$ (see Proposition \ref{prop: eos:imdiff}) to an orbifold diffeomorphism will in general not be injective. However, on a sufficiently small zero-neighborhood one can always achieve this.
 
\begin{prop}\label{prop: diff:chart} \no{prop: diff:chart}
 Consider the family $(\nN_i)_{i \in I}$ as in Proposition \ref{prop: eos:imdiff}. For each $i \in I$, there is an open neighborhood $P_{2,i} \subseteq \vect{\Omega_{5,i}}$ of the zero-section and sets $\mM^{\Omega_{5,i}}_i \coloneq \nN_i^{\Omega_{5,i}} \cap P_{2,i}$, $\mM_i \coloneq (\res_{\Omega_{5,i}}^{W_{\alpha (i)}})^{-1} (\mM^{\Omega_{5,i}}_i)$ such that on the zero-neighborhood $\mM \coloneq \Lambda^{-1}_\cC \left(\bigoplus_{i \in I} \mM_i \right)$, the map 
	\begin{displaymath}
	 E \colon \mM \rightarrow \textstyle \Difforb{Q,\uU} \displaystyle ,\ E([\hat{\sigma}] ) \coloneq [\hat{E}^\sigma] = [\expo] \circ [\hat{\sigma}]|^\Omega, 
	\end{displaymath}
 is injective with $E(\zs )= \ido{}$.
\end{prop}

\begin{proof}
 Following Proposition \ref{prop: eos:imdiff}, each $[\hat{\sigma}] \in \nN = \Lambda_\cC^{-1} \left(\bigoplus_{i \in I} \nN_i \right)$ induces an orbifold diffeomorphism $[\hat{E}^\sigma]$.  Shrink $\nN_i$ to obtain an open $C^1$-neighborhood $\mM_i$ of the zero-section in $\vect{W_{\alpha (i)}}$: \\ 
 Choose for each $i \in I$ a non-singular point $z_i \in U_i$ (which exists due to Newman's Theorem \ref{thm: newman}, since $U_i$ is an open set) and an $H_{\alpha (i)}$-stable $z_i$-neighborhood $U_{z_i} \subseteq W_{\alpha (i)}$ with $H_{\alpha (i), U_{z_i}} = \set{\id_{W_{\alpha (i)}}}$. This is possible since $z_i$ is non singular. The family $\fF_5 (H_{\alpha (i)}.\overline{U_i})$ constructed in \ref{con: sp:ofdat} covers $\overline{U_i}$ and we may choose a chart $(V_{5,\alpha (i)}^k, \kappa^{\alpha (i)}_k)$ such that $z_i \in V_{3,\alpha (i)}^k$. Consider the open set $\hat{U}_{z_i} \coloneq TV_{5,\alpha (i)}^k \cap \hat{O}_{\alpha (i)} \cap \exp_{W_{\alpha (i)}}^{-1} (U_{z_i}) \subseteq TW_{\alpha (i)}$. The intersection $T_{z_i} W_{\alpha (i)} \cap \hat{U}_{z_i}$ is an open zero-neighborhood. 
 We obtain another open zero-neighborhood 
	\begin{displaymath}
	 \lfloor \kappa_k^{\alpha (i)} (z_i) , \text{pr}_2 (T\kappa_{k}^{\alpha (i)} (\hat{U}_{z_i} \cap (\kappa_k^{\alpha (i)} (z_i) \times \RR^d)))\rfloor \subseteq C^\infty (B_5(0), \RR^d)
	\end{displaymath}
 where $\text{pr}_2 \colon B_5 (0) \times \RR^d \rightarrow \RR^d$ is the projection. Define $P_{2,i}\subseteq \vect{\Omega_{5,i}}$ to be the open zero-neighborhood induced by $\lfloor \kappa_k^{\alpha (i)} (z_i) , \text{pr}_2 (T\kappa_{k}^{\alpha (i)} (\hat{U}_{z_i} \cap (\kappa_k^{\alpha (i)} (z_i) \times \RR^d)))\rfloor$. By construction, the map $\exp_{W_{\alpha (i)}} \circ \sigma_{\alpha (i)}$ maps $z_i$ into $U_{z_i}$ if $\sigma_{\alpha (i)}$ is contained in $P_{2,i}$. The intersection $\mM_i^{\Omega_{5,i}} \coloneq \nN^{\Omega_{5,i}}_i \cap P_{2,i} $ is a non-empty open zero-neighborhood in $\vect{\Omega_{5,i}}$. Define $\mM_i \coloneq (\res^{W_{\alpha (i)}}_{\Omega_{5,i}})^{-1} (\mM_i^{\Omega_{5,i}}) \subseteq \nN_i$. Then $\mM \coloneq \Lambda_\cC^{-1} \left(\bigoplus_{i\in I} \mM_i\right)$ contains $\zs$ and is an open subset of $\nN$ in $\Osc{Q}$.\\
 We show that the map $E$ (as in the statement of the proposition) is injective on $\mM$. Assume that there are $[\hat{\sigma}] , [\hat{\tau}] \in \mM$ such that $E([\hat{\sigma}]) = E([\hat{\tau}] )$. For $E ([\hat{\sigma}]) = [\hat{E}^\sigma]$, there is a representative $\hat{E}^\sigma$ in $\Orb{\aA , \cC}$, by Proposition \ref{prop: eos:imdiff}. By assumption, the orbifold maps induced by $\hat{E}^\sigma$ and $\hat{E}^\tau$ coincide, whence $E^\tau = E^\sigma$ follows. We will prove that for each $i \in I$, the lifts $e^{\sigma_{i}}$ and $e^{\tau_{i}}$ coincide. Fix $i \in I$ and observe that $E^\sigma = E^\tau$ implies that for each $z \in U_i$, there is some $\gamma_z \in H_{\alpha (i)}$ with $e^{\sigma_{i}} (z) = \gamma_z .e^{\tau_{i}} (z)$. Consider a component $C$ of $U_i \setminus \Sigma_{G_i}$. The set $\setm{c\in C}{\gamma . e^{\sigma_i} (c) = e^{\tau_i}(c)}$ is an open and closed subset of $C$. As $C$ is connected, there is a unique $\gamma_C \in  H_{\alpha (i)}$ with $e^{\sigma_{i}}|_{\overline{C}} = \gamma_C e^{\tau_{i}}|_{\overline{C}}$. For $x \in \overline{C} \cap \overline{C'}$, this yields the identity $T_x (\gamma_C e^{\tau_{i}}) = T_x e^{\sigma_{i}} = T_x (\gamma_{C'} e^{\tau_{i}})$. Since $e^{\tau_{i}}$ is a diffeomorphism, we derive $T_{e^{\tau_{i}}(x)} \gamma_{C'}^{-1} \gamma_C = T_{e^{\tau_{i}}(x)} \id_{W_{\alpha (i)}}$ and $\gamma_{C'}^{-1}\gamma_C \in H_{\alpha (i), e^{\tau_{i}} (x)}$. By \cite[Lemma 2.10]{follie2003}, $\gamma_C = \gamma_{C'}$ follows. Then $\gamma_C = \gamma_{C'}$ follows for each component such that there is a chain $C=C_1, C_2, \ldots , C_n= C'$ of components with $\overline{C}_k \cap \overline{C_{k+1}} \neq \emptyset$. Observe that by a combination of Lemma \ref{lem: lin:ch} and Lemma \ref{lem: RnGeom} each $x \in \Sigma_{G_i}$ is contained in some $\overline{C}$ and $\bigcup_{x \in \overline{C}} \overline{C}$ is a neighborhood of $x$. Hence there is a unique $\gamma$ with $\gamma . e^{\tau_{i}} = e^{\sigma_{i}}$. Specializing, we obtain $\gamma . e^{\tau_{i}} (z_i) = e^{\sigma_{i}} (z_i)$. The lifts $\sigma_{\alpha (i)},\tau_{\alpha (i)}$ are elements of $\mM_i$, whence by definition of $\mM_i$, $e^{\sigma_{i}} (z_i), e^{\tau_{i}}(z_i) \in U_{z_i}$ holds. The $H_{\alpha (i)}$-stability of $U_{z_i}$ forces $\gamma$ to be in the isotropy subgroup of $U_{z_i}$. Hence $\gamma = \id_{W_{\alpha (i)}}$ holds and we obtain $\exp_{W_{\alpha (i)}} \circ \sigma_{i} = \exp_{W_{\alpha (i)}} \circ \tau_{i}$. Lemma \ref{lem: eos:loc} iii.\ implies that $\im \sigma_i$ and $\im \tau_i$ are contained in $\hat{O}_{\alpha (i)}$. As $\exp_{W_{\alpha (i)}}$ is injective on $T_x W_{\alpha (i)} \cap \hat{O}_{\alpha (i)}$ for $x \in U_i$, we must have $\tau_{i} = \sigma_{i}$. Repeating the argument for $i \in I$, the families $\set{\tau_{i}}_{i \in I}$ and $\set{\sigma_{i}}_{i \in I}$ coincide. As those lifts are canonical lifts, Remark \ref{rem: os:clpro} (a) implies $[\hat{\sigma}] = [\hat{\tau}]$ and $E \colon \mM \rightarrow \Difforb{Q,\uU}$ is injective.
\end{proof}

We now apply the results of Section \ref{sect: barigeo} to construct a neighborhood $\hH$ of the zero-orbisection: 
\begin{con}\label{con: H} \no{con: H}
 Using the local data obtained in Construction \ref{con: sp:ofdat} IV., we define open sets \glsadd{omegaK5}
	\begin{displaymath}
	 \Omega_{r,K_{5,i}} \coloneq \bigcup_{(V_{5,\alpha (i)}^n , \kappa^{\alpha (i)}_n) \in \fF_5 (K_{5,i})} V_{r,\alpha (i)}^n,\quad r \in [1,5].
	\end{displaymath}
 By construction, $\Omega_{5,i} \subseteq \overline{\Omega_{5,i}} = K_{5,i} \subseteq \Omega_{r,K_{5,i}}$ holds for each $r \in [1,5]$.\\
 In Proposition \ref{prop: diff:chart} we have constructed sets $\mM_i^{\Omega_{5,i}}$ as intersections $\mM_i^{\Omega_{5,i}} = N_i^{\Omega_{5,i}} \cap P_{1,i} \cap P_{2,i}$, where $N_i^{\Omega_{5,i}}$ is an open zero-neighborhood as in Lemma \ref{lem: vect:emb}. Apply Construction \ref{con: comp:loc} with $R_i$ (see Construction \ref{con: sp:ofdat} V.) taking the role of $R$, $M \coloneq W_{\alpha (i)}$, $K \coloneq K_{5,i}$ and $P \coloneq P_{1,i} \cap P_{2,i}$ to construct an open zero-neighborhood $\hH_{R_i} \subseteq \mM_i \subseteq \vect{W_{\alpha (i)}}$. The set $E_{5,K}$ occuring in Lemma \ref{lem: vect:emb} is $N_i^{\Omega_{5,i}}$ from the proof of Lemma \ref{lem: eos:loc}. By construction, $\hH_{R_i} = \left(\res_{\Omega_{5, K_{5,i}}}^{W_{\alpha (i)}}\right)^{-1} (\hH_{R_i}^{K_{5,i}})$ holds for an open zero neighborhood $\hH_{R_i}^{K_{5,i}} \subseteq \vect{\Omega_{5, K_{5,i}}}$. Finally, for each $i \in I$ the construction yields a constant $0 <\tau_i, \nu_i < R_i$ with the following property:\\
 If $X \in \vect{W_{\alpha (i)}}$ such that for each $(V_{5,\alpha (i)}^k, \kappa_k^{\alpha (i)}) \in \fF_5 (K_{5,i})$, 
 the local representative $X_k$ satisfies $\norm{X_k}_{\overline{B_1 (0)},1} \leq \tau_i$, then $X$ is contained in $\mM_i$.\\
 Recall from Construction \ref{con: comp:loc} that for each pair $(X,Y) \in \hH_{R_i} \times \hH_{R_i}$, there are unique vector fields $X \diamond_i Y, X^{*_i}, Y^{*_i} \in \vect{\Omega_{\frac{5}{4},K_{5,i}}}$. Together with the definition of $R_i$ (\ref{con: sp:ofdat} V.), the estimates \eqref{eq: locdiamest} and \eqref{eq: locinvest} imply the following properties, which we note here for later use:
  \begin{equation}\label{eq: mix:work}
   X \diamond_i Y (x), X^{*_i} (x)  \in B_{\rho_{\alpha (i)}} (0_x, S_{\alpha (i)}) \subseteq \hat{O}_{\alpha (i)}, \quad \forall x \in \Omega_{\frac{5}{4}, K_{5,i}}.
  \end{equation}
  Moreover, for each chart $(V_{5,\alpha (i)}^n , \kappa_n^{\alpha (i)}) \in \fF_5 (K_{5,i})$, the vector field $X \diamond_i Y$ satisfies the estimate \eqref{eq: locdiamest}, i.e.\ $\norm{(X \diamond_i Y)_{[n]}}_{\overline{B_{\frac54}}(0)} < \nu_i$. Recall that $\nu_i$ in Construction \ref{con: comp:loc} is chosen exactly as in Lemma \ref{lem: eos:cov}. Hence Lemma \ref{lem: eos:cov} (b) yields for $X,Y \in \hH_{R_i}$ and $x \in V_{\frac54 , \alpha (i)}^n$ the identity 
  \begin{equation}\label{eq: local:tausch}
   \exp_n T\kappa_n^{\alpha (i)} (X \diamond_i Y) (x) = \kappa_n^{\alpha (i)} \exp_{W_{\alpha (i)}} (X \diamond_i Y) (x). 
  \end{equation}
Define the open subset $\hH \coloneq \Lambda_\cC^{-1} \left(\bigoplus_{i \in I} \hH_{R_i}\right)$ of $\Osc{Q}$. By construction, $\zs \in \hH \subseteq \mM$.
\end{con}

The vector fields $X \diamond_i Y$ and $X^{*_i}$ induced by orbisections in $\hH$ yield families whose members are $\lambda$-related for suitable changes of orbifold charts $\lambda$. The details are checked in the next lemma.
  
\begin{lem}\label{lem: vf:eq} \no{lem: vf:eq}
Consider orbisections $[\hat{\sigma}], [\hat{\tau}] \in \hH$ with families of canonical lifts $\set{\sigma_j}_{j \in J}, \set{\tau_j}_{j \in J}$ with respect to the atlas $\bB$. Let $\lambda \in \CH{W_k}{W_l}$ be a change of charts which satisfies $\dom \lambda \subseteq \Omega_{\frac{5}{4}, K_{5,i}}$ and $\im \lambda \subseteq \Omega_{\frac{5}{4},K_{5,j}}$ for $k = \alpha (i)$ and $l = \alpha (j)$. Then the following identities hold:
 \begin{align}\label{eq: emb:ch1}
    T\lambda  (\sigma_k \diamond_i \tau_k)|_{\dom \lambda} 	&= (\sigma_l \diamond_j \tau_l) \circ \lambda \\
    T\lambda  \sigma_k^{*_i}|_{\dom \lambda} 			&= \sigma_l^{*_j} \circ \lambda. \label{eq: emb:ch2}
 \end{align}
 Then the maps $\sigma_j \diamond_i \tau_j|_{U_i}$ and $\sigma_j^{*_i}|_{U_i}$ are equivariant with respect to the derived action of $G_{i}$. 
\end{lem}

\begin{proof}
 The identities \eqref{eq: emb:ch1} and \eqref{eq: emb:ch2} may be checked locally. Fix $x \in \dom \lambda \subseteq \Omega_{\frac{5}{4},K_{5,i}}$ together with a chart $(V_{5,k}^n,\kappa_n^k)\in \fF_5 (K_{5,i})$ such that $x \in V_{\frac{5}{4},k}^n$. The manifold atlas chosen for $\kK^\circ_k \subseteq W_k$ is subordinate to the cover $(Z_k^r \cap \kK^\circ_k)_{1\leq r \leq N_k}$. Hence there is some $Z^r_k$ with $V_{5,k}^n \subseteq Z_k^r$. As $x \in V_{5,k}^n \subseteq \kK_k$ and $\lambda (x) \in \Omega_{\frac{5}{4},K_{5,j}} \subseteq \kK_l$, by construction \ref{con: sp:ofdat} (cf.\ Lemma \ref{lem: oemb:cov}), there is an embedding of orbifold charts $\mu \colon Z_k^r \rightarrow W_l$ with $\mu (x) = \lambda (x)$. After possibly replacing $\mu$ with $\gamma \circ \mu$ for suitable $\gamma \in H_l$, there is an open neighborhood $U_x$ of $x$ in $\Omega_{\frac54 , K_{5,i}}$ with $\mu|_{U_x} = \lambda|_{U_x}$. By construction, we obtain $\mu (x) = \lambda (x) \in \Omega_{\frac{5}{4},K_{5,j}} \subseteq K_l^\circ$ and $T_x \mu = T_x \lambda$ holds. The definition of $S_k$ together with equation \eqref{eq: mix:work} implies $T\mu (\sigma_j \diamond_i \tau_j) (x) , T\mu \sigma_{j}^{*_i} (x) \in \hat{O}_l$ and $(\sigma_l \diamond_j \tau_l) \mu (x), \sigma_l^{*_j} \mu (x) \in \hat{O}_l$. Let $\exp_{n}$ be the Riemannian exponential map induced by the pullback metric on $B_5 (0)$ with respect to $(\kappa_n^k)^{-1}$. The map $\mu (\kappa_n^k)^{-1}$ is a Riemannian embedding of $B_5 (0)$ into $W_l$. From \cite[IV. Proposition 2.6]{fdiffgeo1963}, we deduce for $v \in \dom \exp_n$ that \begin{equation}\label{eq: riemb}
           \exp_{W_l} T\mu (\kappa_n^k)^{-1} (v) = \mu (\kappa_n^k)^{-1} \exp_{n} (v) .                                                                                                                                                                                                                                                                                                                                                                                                                                                                                                                                                                                                                                                                                                                                                                                                                                                                                                                                                                                                                                                                                                                                                                                                                                                                                                                                                                                                                                                                                                                                                                                                                             \end{equation}
 Recall from Construction \ref{con: H} that for $i\in I$, there is some open set $\hH_{R_i}$ with the same properties as in Lemma \ref{lem: vect:emb} such that $[\hat{\sigma}] \in  \hH$ implies $\sigma_{k} \in \hH_{R_i}$. For $X \in \hH_{R_i}$, we have: 
    \begin{compactenum}
     \item[i.] $\kappa_n^k \exp_{W_k} \circ X (z) = \exp_n T\kappa_n X(z)$ for each $z \in V_{3,k}^n$ (combine Lemma \ref{lem: eos:cov} (b) and (f)),
     \item[ii.] $\exp_{W_k} \circ X (V_{\frac{5}{4},k}^n) \subseteq V_{2,k}^n$ and $\exp_{W_k} \circ X (\overline{V_{2,k}^n}) \subseteq V_{3,k}^n$, (see Lemma \ref{lem: eos:cov} (d)),
     \item[iii.] $V_{\frac{5}{4},k}^n \subseteq \exp_{W_k} \circ X(V_{2,k}^n)$ (see Lemma \ref{lem: eos:cov} (d)).
    \end{compactenum}
 The families $\set{\sigma_k}$ and $\set{\tau_k}$ are canonical families, whence $\sigma_l \mu = T\mu \sigma_k$ holds. In addition, for the vector field $\sigma_k \diamond_i \tau_k$ on $V_{\frac{5}{4},k}^n$ the local identities \eqref{eq: glo:loc} and \eqref{eq: def2} are available. Combining these facts we compute: 
	\begin{align*}
	 \exp_{W_l} T_x \lambda (\sigma_k \diamond_i \tau_k) (x) &= \exp_{W_l} T_x\mu (\sigma_k \diamond_i \tau_k) (x)= \exp_{W_l} T(\mu (\kappa_n^k)^{-1}\kappa_{n}^k) (\sigma_k \diamond_i \tau_k) (x)\\ 
         &\hspace{-0.28cm}\overset{\eqref{eq: riemb}}{=} \mu (\kappa_n^k)^{-1} \exp_n T\kappa_n^k  (\sigma_k \diamond_i \tau_k) (x) \stackrel{\eqref{eq: local:tausch}}{=} \mu (\kappa_n^k)^{-1} \kappa_n^k \exp_{W_k} (\sigma_k \diamond_i \tau_k) (x) \\
         &\hspace{-0.3cm}\overset{\eqref{eq: glo:loc}}{=} \mu \exp_{W_k} (\exp_{W_k}|_{N_x})^{-1} \exp_{W_k} \sigma_{k} \exp_{W_k} \tau_k (x) \\
         &\hspace{-0.45cm}\stackrel{i. + \eqref{eq: riemb}}{=} \exp_{W_l} T\mu \sigma_{k} \exp_{W_k} \tau_k (x) = \exp_{W_l} \sigma_l \mu \exp_{W_k} \tau_k (x) \\
         &\hspace{-0.45cm}\stackrel{i.  + \eqref{eq: riemb}}{=} \exp_{W_l} \sigma_l \exp_{W_l} \tau_l \mu (x) =(\exp_{W_l} \sigma_l \exp_{W_l} \tau_l) \lambda (x)\\
         &\hspace{-0.3cm} \stackrel{\eqref{eq: def1}}{=} \exp_{W_l} (\sigma_l \diamond_j \tau_l) (\lambda (x)) .
	\end{align*}
 Since $\exp_{W_l}$ restricts to a diffeomorphism on $T_{\lambda (x)} W_l \cap \hat{O}_l$, the computation yields \eqref{eq: emb:ch1}.\\
 To obtain \eqref{eq: emb:ch2}, we use $x \in V_{\frac{5}{4},k}^n$ and compute with the facts from above: 
  \begin{displaymath}
   \exp_{W_l} T_x \lambda \sigma_k^{*_i} (x) = \exp_{W_l} T\mu \sigma_k^{*_i} (x) \stackrel{\eqref{eq: def2}}{=} \mu (\exp_{W_k} \circ \sigma_k|_{\Omega_{2,K_{5,i}}})^{-1} (x).
  \end{displaymath}
 As $x \in V_{\frac{5}{4},n}^k$, by iii.\ the image $(\exp_{W_k} \circ \sigma_k|_{\Omega_{2,K_{5,i}}})^{-1} (x)$ is contained in $V_{2,k}^n$. Since $T \kappa_n^k \sigma_k (V_{2,k}^n) \subseteq \dom \exp_n$, we conclude with \eqref{eq: riemb} that $\sigma_l \mu (V_{2,k}^n) = T\mu \sigma_k (V_{2,k}^n) \subseteq \dom \exp_{W_l}$. Thus we may consider: 
  \begin{align*}
   (\exp_{W_l} \sigma_l)\circ \exp_{W_l} T_x \lambda \sigma_k^{*_i} (x) &= \exp_{W_l} \sigma_l \mu (\exp_{W_k} \sigma_k|_{\Omega_{2,K_{5,i}}})^{-1} (x) \\ 
      &=\exp_{W_l} T \mu \sigma_k  (\exp_{W_k} \sigma_k|_{\Omega_{2,K_{5,i}}})^{-1} (x) \\
      &=\mu (\exp_{W_k} \sigma_k) (\exp_{W_k} \sigma_k|_{\Omega_{2,K_{5,i}}})^{-1} (x) =  \mu (x) = \lambda (x) \in \Omega_{\frac{5}{4},K_{5,j}}.
  \end{align*}
 Recall $\lambda (x) \in \Omega_{\frac{5}{4},K_{5,l}}$ and $T\mu \sigma_{k}^{*_i} (x) \in \hat{O}_l$. Now the definition of $\hat{O}_l$ in Construction \ref{con: sp:ofdat} V. yields $\exp_{W_l} T_x \mu \sigma_{k}^{*_i} (x) \in \Omega_{2, K_{5,j}}$. On $\Omega_{2, K_{5,j}}$ the map $\exp_{W_l} \circ \sigma_l$ is injective, by Step 1 in the proof of Lemma \ref{lem: eos:loc}. We deduce that $\exp_{W_l} T\lambda  \sigma_k^{*_i} (x) = \exp_{W_l} \sigma_l^{*_j} (\lambda (x))$ must hold. Since $\exp_{W_l}$ restricts to a diffeomorphism on $T_{\lambda (x)} \cap \hat{O}_j$, the computation yields \eqref{eq: emb:ch2}.
\end{proof}
The families $\set{\sigma_j \diamond_i \tau_j}_{i \in I}$ and $\set{\sigma_j^{*_i}}_{i \in I}$ obtained in this way induce orbisections:
 
\begin{prop}\label{prop: os:glue} \no{prop: os:glue}
 Consider orbisections $[\hat{\sigma}], [\hat{\tau}] \in \hH$, whose canonical families with respect to $\bB$ are given by $\set{\sigma_j}_{j \in J}$ and $\set{\tau_j}_{j \in J}$, respectively. Then 
	\begin{compactenum}
	 \item The family $\set{\sigma_{\alpha (i)} \diamond_i \tau_{\alpha (i)}}_{i\in I}$ induces an orbisection $[\widehat{\sigma \diamond \tau}] \in \mM$ whose family of canonical lifts with respect to the atlas $\aA$ is given by $(\sigma \diamond \tau)_i \coloneq \sigma_{\alpha (i)} \diamond_i \tau_{\alpha (i)}|_{U_i}$ for $i \in I$.
	 \item The family $\set{\sigma_{\alpha (i)}^{*_i}}_{i\in I}$ induces an orbisection $[\widehat{\sigma^*}] \in \mM$ whose canonical lifts with respect to the atlas $\aA$ are given by $(\sigma^*)_i \coloneq \sigma_{\alpha (i)}^{*_i}|_{U_i}$ for $i \in I$.
	\end{compactenum}
\end{prop}

\begin{proof}
 The families $\set{(\sigma \diamond \tau)_{i}}_{i \in I}$ and $\set{(\sigma^*)_i}_{i \in I}$ are compatible families of vector fields on the atlas $\aA$ by Lemma \ref{lem: vf:eq}. These families yield canonical families of lifts with respect to the atlas $\aA$. 
 In particular, the identities \eqref{eq: emb:ch1} and \eqref{eq: emb:ch2} allow the definition of continuous maps: 
	\begin{align*}
	 \sigma \diamond \tau &\colon Q \rightarrow \tT Q, x \mapsto T\psi_i (\sigma \diamond \tau)_i\psi_i^{-1}(x) \quad \text{ if } x \in \psi_i (U_i)\\
	 \sigma^* &\colon Q \rightarrow \tT Q, x \mapsto T\psi_i (\sigma^*)_i\psi_i^{-1}(x)     \quad \text{ if } x \in \psi_i (U_i).
	\end{align*}
 These data allow the definition of orbisections $[\widehat{\sigma \diamond \tau}]$ and $[\widehat{\sigma^*}]$ by Remark \ref{rem: os:clpro} (a). \\
 To complete the proof, we have to show that $[\widehat{\sigma \diamond \tau}] , [\widehat{\sigma^*}]$ are contained in $\mM$. To this end, we need to assure that $[\widehat{\sigma \diamond \tau}]$ and $[\widehat{\sigma^*}]$ are compactly supported. The orbisections $[\hat{\sigma}] , [\hat{\tau}] \in \hH$ are compactly supported, whence $\supp [\hat{\sigma}] \cup \supp [\hat{\tau}]$ is contained in a compact subset $K \subseteq Q$. Since $\bB$ is locally finite, there is a finite subset $\sS_{\sigma, \tau} \subseteq \bB$ such that $(W_j,H_j,\varphi_j) \in \sS_{\sigma , \tau}$ if and only if $\im \varphi_j \cap K \neq \emptyset$. Consider $(W_j , H_j, \varphi_j) \in \bB \setminus \sS_{\sigma , \tau}$. By Remark \ref{rem: os:clpro} (d) the canonical lifts of $[\hat{\sigma}], [\hat{\tau}]$ on $W_j$ are the zero-section in $\vect{W_j}$. The conclusion in Construction \ref{con: comp:loc} implies that $\sigma_j \diamond_i \tau_j \equiv 0$ and $\sigma_j^{*_i} \equiv 0$ for each $i \in \alpha^{-1} (j)$. Therefore the supports $\supp [\widehat{\sigma \diamond \tau}]$ and $\supp [\widehat{\sigma^*}]$ are contained in $K_{\sigma, \tau} \coloneq \bigcup_{(W_{\alpha (i)}, H_{\alpha (i)}, \varphi_{\alpha (i)}) \in \sS_{\sigma ,\tau}} \overline{\psi_i (U_i)}$. As $\sS_{\sigma, \tau}$ is finite and for $j \in J$ the set $\alpha^{-1} (j)$ is finite, $K_{\sigma,\tau}$ is a finite union of compact sets $\overline{\psi_i (U_i)}$. Hence the supports of $[\widehat{\sigma \diamond \tau}]$ and $[\widehat{\sigma^*}]$ are contained in a compact set, whence these orbisections are compactly supported.\\
 Following Proposition \ref{prop: os:cl}, we may consider the canonical lifts $(\sigma \diamond \tau)_k$ and $\sigma^*_k$ on each chart $(W_k, H_k , \varphi_k)\in \bB$. The orbisections $[\widehat{\sigma \diamond \tau}] , [\widehat{\sigma^*}]$ will be contained in $\mM$ if their respective canonical lifts are contained in $\mM_i$ for each $i \in \alpha^{-1} (k)$, $k \in J$.\\ Fix $i \in \alpha^{-1}(k)$ and define $(\sigma \diamond \tau)_k)_{[n]} \coloneq ((\sigma \diamond \tau)_k)_{\kappa_n} \circ \kappa_n^{-1}$\glsadd{sigman} and $(\sigma^*_k)_{[n]}\coloneq (\sigma^*_k)_{\kappa_n} \circ \kappa_n^{-1}$, respectively (cf.\ Definition \ref{defn: top:vect}) for $(V_{5,k}^n,\kappa_k^n) \in \fF_5 (K_{5,i})$. By construction \ref{con: H}, it suffices to prove for each chart $(V_{5,k}^n,\kappa_{n}^k)$ in $\fF_5 (K_{5,i})$, the condition $\norm{((\sigma \diamond \tau)_k)_{[n]}}_{\overline{B_1 (0)},1} < \tau_i$ respectively $\norm{(\sigma^*_k)_{[n]}}_{\overline{B_1 (0)},1} < \tau_i$ holds. Observe that the conditions may be checked on $\Omega_{\frac{5}{4},K_{5,i}}$. Uniqueness of canonical lifts together with \eqref{eq: emb:ch1} and \eqref{eq: emb:ch2} forces the canonical lifts $(\sigma \diamond \tau)_k$ respectively $(\sigma^*)_k$ to coincide with $\sigma_k \diamond_i \tau_k$ respectively $\sigma_k^{*_i}$ on $\Omega_{\frac{5}{4},K_{5,i}}$. Recall from the construction that the constant $\tau_i$ corresponds to the constant $\tau$ in Construction \ref{con: comp:loc}. Hence a combination of \eqref{eq: locinvest} with Corollary \ref{cor: comp:tau} yields $\norm{((\sigma \diamond \tau)_k)_{[n]}}_{\overline{B_1 (0)},1} =  \norm{(\sigma_k)_{[n]} \diamond (\tau_k))_{[n]}}_{\overline{B_1 (0)},1} < \tau_i$ and $\norm{(\sigma^*_k)_{[n]}}_{\overline{B_1 (0)},1} = \norm{((\sigma_k))_{[n]})^{*}}_{\overline{B_1 (0)},1}< \tau_i$. We conclude that each of the canonical lifts of $[\widehat{\sigma \diamond \tau}]$ and $[\widehat{\sigma^*}]$ on $(W_k,H_k,\varphi_k)$ is contained in $\mM_i$ with $i \in \alpha^{-1} (k)$. Summing up, $[\widehat{\sigma \diamond \tau}]$ and $[\widehat{\sigma^*}]$ are contained in $\mM$.
\end{proof}

\begin{rem}\label{rem: new:rem}
 \begin{compactenum}
  \item The last lemma implies that the map $E$ may be applied to $[\sigma \diamond \tau]$ and $[\sigma^*]$ for $[\sigma], [\tau] \in \hH$.
  \item Moreover, consider the canonical lifts $(\sigma \diamond \tau)_{W_{\alpha (i)}}$ and $\sigma^*_{W_{\alpha (i)}}$ of $[\sigma \diamond \tau]$ and $[\sigma^*]$, respectively, for $[\sigma], [\tau] \in \hH$ on a chart $(W_{\alpha (i)}, H_{\alpha (i)}, \varphi_{\alpha (i)}) \in \bB$ for $i\in I$. Let again $\sigma_{\alpha (i)}$ and $\tau_{\alpha (i)}$ be the canonical lifts of $[\hat{\sigma}]$ and $[\hat{\tau}]$, respectively, on $(W_{\alpha (i)}, H_{\alpha (i)}, \varphi_{\alpha (i)})$. Then uniqueness of canonical lifts shows that the restrictions of these vector fields to $\Omega_{\frac{5}{4},i}$ satisfy 
  \begin{displaymath}
   (\sigma \diamond \tau)_{W_{\alpha (i)}}|_{\Omega_{\frac{5}{4},i}} = \sigma_{\alpha (i)} \diamond_i \tau_{\alpha (i)} \quad \text{and} \quad \sigma^*_{W_{\alpha (i)}} = \sigma_{\alpha (i)}^{*_i},
  \end{displaymath}
by Lemma \ref{lem: vf:eq}.
 \end{compactenum}
\end{rem}

  In the rest of this section, these properties will be crucial for several key arguments. We shall now assure that the orbisections constructed satisfy the identities needed for composition and inversion in $E(\mM)$: 
 
\begin{lem}\label{lem: comp:cor} \no{lem: comp:cor}
 Consider $[\hat{\sigma}], [\hat{\tau}] \in \hH$. The following identities hold: 
	\begin{align}
	 E([\hat{\sigma}]) \circ E([\hat{\tau}]) &= E([\widehat{\sigma \diamond \tau}]) \label{eq: comp} \\
	 E([\hat{\sigma}])^{-1} &= E([\widehat{\sigma^*}]). \label{eq: inv}
	\end{align}
 \end{lem}

\begin{proof}
 Choose and fix arbitrary $[\hat{\sigma}], [\hat{\tau}] \in \hH$. The left hand and the right hand sides of the equations \eqref{eq: comp} resp.\ \eqref{eq: inv} are orbifold diffeomorphisms. As observed in Proposition \ref{prop: diff:undmap} and Corollary \ref{cor: ofdiff:un}, orbifold diffeomorphisms are uniquely determined by their underlying maps or their family of lifts. To prove the assertion it therefore suffices to show that their family of lifts or the underlying maps on both sides are equal.\\
 Consider the right hand sides of both equations: The orbisections $[\sigma \diamond \tau]$ and $[\sigma^*]$ have been constructed by a family of canonical lifts $\set{\sigma \diamond \tau)_i}_{i \in I}$ resp. $\set{(\sigma^*)_i}_{i \in I}$ with respect to the atlases $\aA$ and $\tT \aA$. Both orbisections are contained in $\mM$. Taking identifications $\im (\sigma \diamond \tau)_i , \im (\sigma^*)_i  \subseteq \hat{O}_{\alpha (i)}$ holds. Corestriction of each lift to $TU_i \cap \hat{O}_{\alpha (i)}$ yields representatives of $[\hat{\sigma}]|^\Omega$ and $[\hat{\tau}]|^\Omega$. 
 Thus representatives of $E([\widehat{\sigma \diamond \tau}])$ and $E([\widehat{\sigma^*}])$ are given by $(E^{\sigma \diamond \tau}, \set{e^{(\sigma \diamond \tau)_{i}}}_{i \in I} , P, \nu)$ respectively $(E^{\sigma^*}, \set{e^{\sigma^{*}_{i}}}_{i \in I} , P', \nu')$ in $\Orb{\aA, \cC}$. The lifts of these maps satisfy for each $i \in I$ by construction: 
	\begin{align}
	 \exp_{W_{\alpha (i)}} \circ \sigma_{\alpha (i)} \circ \exp_{W_{\alpha (i)}} \circ \tau_i &= \exp_{W_{\alpha (i)}} \circ (\sigma \diamond \tau)_i= e^{(\sigma \diamond \tau)_{i}} \label{eq:comploc}\\ 
	 (\exp_{W_{\alpha (i)}} \circ \sigma_{\alpha (i)}|_{\Omega_{2,i}})^{-1}|_{U_i} &= \exp_{W_{\alpha (i)}} \circ \sigma^*_i = e^{\sigma^*_{i}}. \label{eq: comp:inv}
	\end{align}
 We show that the lifts in \eqref{eq:comploc} coincide with the lifts of $E([\hat{\sigma}]) \circ E([\hat{\tau}])$. As $\ORB$ is a category, composition in $\ORB$ is associative. Hence lifts can be computed iteratively: $E([\sigma]) \circ E ([\tau]) = [\expo] \circ [\sigma]|^\Omega\circ [\expo]\circ [\tau]|^\Omega = [\expo] \circ ([\sigma]|^\Omega \circ [\expo] \circ [\tau]|^\Omega)$. As $\tau_{\alpha (i)}$ and $\sigma_{\alpha (i)}$ are contained in $\hH_{R_i}$, the composition of charted orbifold maps (cf. Construction \ref{con: cp:chom}) yields a lift of $E^\sigma \circ E^\tau$ on $U_i$ which coincides with the left hand side of \eqref{eq:comploc}. Therefore \eqref{eq: comp} follows from \eqref{eq:comploc} by an application of Corollary \ref{cor: ofdiff:un}.\\
 To prove the identity \eqref{eq: inv} we show that the underlying maps of both sides are equal. To this end, let $e^{\sigma^*}$ be the underlying map of $E([\hat{\sigma}^*])$. By Proposition \ref{prop: diff:undmap}, it suffices to check the identity 
  \begin{displaymath}
   \expo \circ \sigma \circ e^{\sigma^*} = \id_Q.
  \end{displaymath}
 If this identity holds, then assertion \eqref{eq: inv} follows. Clearly the identity can be chekced locally for each chart $(U_i,G_i,\psi_i) \in \aA$. By construction on $U_i$ we have $e^{\sigma^*} \psi_i = \varphi_{\alpha (i)} \circ e^{\sigma^*_i}$. Here $e^{\sigma^*_i}$ is the lift of $E([\widehat{\sigma^*}])$ in the chart $U_i$. Fix $x\in U_i$ and notice $\im e^{\sigma^*_i} \subseteq \Omega_{2,i}$ by \eqref{eq: comp:inv}. Choose a $H_{\alpha (i)}$-stable neighborhood $U_x \subseteq \Omega_{2,i}$ of $e^{\sigma^*_i} (x)$ in $W_{\alpha (i)}$. Restrict the canonical lift $\sigma_{\alpha (i)}$ of $[\hat{\sigma}]$ on $W_{\alpha (i)}$ to $U_x$. Then $\sigma_{U_x} \coloneq \sigma_{\alpha (i)}|_{U_x}^{TU_x}$ is a canonical lift of $[\hat{\sigma}]$ on the chart $(U_x,H_{\alpha (i), U_x}, \varphi_{\alpha (i)}|_{U_x})$. From $U_x \subseteq \Omega_{2,i}$ and $[\hat{\sigma}] \in \hH$, we deduce $\im \sigma_{U_x} = \sigma_{\alpha (i)} (U_x) \subseteq \hat{O}_{\alpha (i)}$, by Lemma \ref{lem: eos:loc} iii. Taking identifications, we may compose $\sigma_{U_x}$ and $\Exp_{TU_x \cap \hat{O}_{\alpha (i)}} \coloneq \exp_{W_{\alpha (i)}}|_{TU_x \cap \hat{O}_{\alpha (i)}}$. Recall from Lemma \ref{lem: eos:loc} that $\Exp_{TU_x \cap \hat{O}_{\alpha (i)}}$ is a lift of $\expo$. Moreover, Construction \ref{con: cp:chom} shows that $\Exp_{TU_x \cap \hat{O}_{\alpha (i)}} \circ \sigma_{U_x}$ is a lift of $\expo \circ \sigma$. Hence, we obtain the following identities 
    \begin{align*}
     \expo \circ \sigma \circ e^{\sigma^*} \psi_i (x) &= \expo \circ \sigma \circ \varphi_{\alpha (i)} \circ e^{\sigma^*_i} (x) = \varphi_{\alpha (i)} (\Exp_{TU_x \cap \hat{O}_{\alpha (i)}} \circ \sigma_{U_x} \circ e^{\sigma^*_i} (x)) \\
             &\overset{\eqref{eq: comp:inv}}{=} \varphi_{\alpha (i)} ((\exp_{W_{\alpha (i)}} \circ \sigma_{\alpha (i)}|_{U_x}) \circ (\exp_{W_{\alpha (i)}} \circ \sigma_{\alpha (i)}|_{\Omega_{2,i}})^{-1}|_{U_i} (x))\\
             &\overset{U_x \subseteq \Omega_{2,i}}{=} \varphi_{\alpha (i)} (x) = \psi_i (x)
    \end{align*}
  Since $x \in U_i$ has been chosen arbitrarily, we may repeat the construction for each $x \in U_i$, whence $\expo \circ \sigma \circ e^{\sigma^*} = \id_Q$ and thus \eqref{eq: inv} follow.
\end{proof}

We now turn our attention to the composition and inversion maps:

\begin{lem}\label{lem: osc:sm} \no{lem: osc:sm}
 The maps 
	\begin{align*}
	 \comp \colon & \hH \times \hH \rightarrow \mM \subseteq \Osc{Q},\hspace{4pt} ([\hat{\sigma}] ,[\hat{\tau}]) \mapsto [\widehat{\sigma \diamond \tau}] \\
	 \inv \colon &\hH \rightarrow \mM \subseteq \Osc{Q},\phantom{H \times H } \quad [\hat{\sigma}] \mapsto [\widehat{\sigma^*}]
	\end{align*}
 are smooth.
\end{lem}

\begin{proof}
The atlases $\aA$ and $\cC$ are indexed by $I$. Let $\sigma_i$ and $\sigma_{\alpha (i)}$ be the canonical lifts with respect to $(U_i, G_i, \psi_i) \in \aA$ and $(W_{\alpha (i)}, H_{\alpha (i)}, \varphi_{\alpha (i)}) \in \cC$, respectively. The continuous linear maps $\tau_i \colon \Osc{Q} \rightarrow \vect{U_i}, [\hat{\sigma}] \mapsto \sigma_i$ and $\lambda_i \colon \Osc{Q} \rightarrow \vect{W_{\alpha (i)}}, [\hat{\sigma}] \mapsto \sigma_{\alpha (i)}$ induce patchworks for $\Osc{Q}$, by Corollary \ref{cor: ost:prop}. 
 The product $\Osc{Q} \times \Osc{Q}$ is a locally convex vector space and we have the family of maps $\lambda_i \times \lambda_i \colon \Osc{Q} \times \Osc{Q} \rightarrow \vect{W_{\alpha (i)}} \times \vect{W_{\alpha (i)}}$ for $i\in I$. 
 Arguments as in the proof of Lemma \ref{lem: comp:sm} show that the family $(\lambda_i \times \lambda_i)_{i \in I}$ yields a patchwork for $\Osc{Q} \times \Osc{Q}$. Let $p$ be the corresponding topological embedding for this patched space (cf. Definition \ref{defn: patloc}).\\
 The patchwork on each of the spaces $(\Osc{Q} \times \Osc{Q}$, $(\lambda_i \times \lambda_i)_{i \in I})$ , $(\Osc{Q}, (\lambda_i)_{i\in I})$ and $(\Osc{Q}, (\tau_i)_{i\in I})$ is indexed by $I$. On the open set $\hH_{R_i}$ constructed in \ref{con: H} consider the maps 
	\begin{align*}
	 \comp_i \colon &\hH_{R_i} \times \hH_{R_i} \rightarrow \vect{U_i}, (X,Y) \mapsto X \diamond_i Y|_{U_i} \\
	 \inv_i \colon &\hH_{R_i} \rightarrow \vect{U_i}, X \mapsto X^{*_i}|_{U_i}.
	\end{align*}
 Since $\hH = \Lambda_\cC^{-1} (\oplus_{i \in I} \hH_{R_i})$, the identities for the patchwork established in the proof of Lemma \ref{lem: comp:sm} 
 yield $p (\hH \times \hH) \subseteq \oplus_{i \in I} (\hH_{R_i} \times \hH_{R_i})$ and $\Lambda_\cC (\hH) \subseteq \oplus_{i \in I} \hH_{R_i}$. By construction, we deduce from Proposition \ref{prop: os:glue}:
 \begin{displaymath}
  \left(\comp_i \right)_{i \in I} p|^{\oplus (\hH_{R_i} \times \hH_{R_i})}_{\hH \times \hH} = \Lambda_\aA \circ \comp \text{  and  } \left(\inv_i \right)_{i \in I} \Lambda_{\cC}|_{\hH}^{\oplus_{i \in I}\hH_{R_i}}= \Lambda_\aA \inv.
 \end{displaymath}
 These mappings make sense, since $\comp_i$ and $\inv_i$ vanish on the zero element. Hence $\comp$ and $\inv$ are patched mappings. By Proposition \ref{prop: pat:loc}, it is sufficient to prove that $\comp$ and $\inv$ are smooth on the patches, i.e.\ for each $i\in I$, the maps $\comp_i$ and $\inv_i$ are smooth. For the remainder of this proof we therefore fix $i \in I$ and prove the smoothness of $\comp_i$ and $\inv_i$:\\
 The open sets $\Omega_{r,K_{5,i}}$, $r \in [1,5]$ contain $U_i$. Consider the restriction maps $\res^{W_{\alpha (i)}}_{\Omega_{5,K_{5,i}}}$, $\res^{\Omega_{r,K_{5,i}}}_{U_i}$ which are linear and continuous, whence smooth by \cite[Lemma F.15 (a)]{hg2004}. Recall that the maps 
	\begin{align*}
	 c_i &\colon \hH_{R_i}^{\Omega_{5,K_{5,i}}} \times \hH_{R_i}^{\Omega_{5,K_{5,i}}} \rightarrow \vect{\Omega_{\frac{5}{4},K_{5,i}}}, X \mapsto X \diamond_i Y\\
	 \iota_i &\colon \hH_{R_i}^{\Omega_{5,K_{5,i}}} \rightarrow \vect{\Omega_{\frac{5}{4},K_{5,i}}} , X \mapsto X^{*_i}
	\end{align*}
 are smooth by Lemma \ref{lem: comp:sm}. By definition the maps $\comp_i$ and $\inv_i$ are given as compositions: 
  \begin{align*}
   \comp_i &=\res^{\Omega_{\frac{5}{4},K_{5,i}}}_{U_i} \circ c_i \circ (\res^{W_{\alpha (i)}}_{\Omega_{5,K_i}} \times \res^{W_{\alpha (i)}}_{\Omega_{5,K_i}}|_{\hH_{R_i} \times \hH_{R_i}})\\
   \inv_i  &= \res^{\Omega_{\frac{5}{4},K_{5,i}}}_{U_i} \circ \iota_i \circ \res^{W_{\alpha (i)}}_{\Omega_{5,K_i}}|_{\hH_{R_i}}.
  \end{align*}
 We conclude that $\comp_i$ and $\inv_i$ are smooth, whence $\comp$ and $\inv$ are smooth. 
\end{proof}

Endow $E(\mM)$ with the smooth manifold structure making $E\colon \mM \rightarrow E(\mM)$ a diffeomorphism.
We are now in a position to construct a Lie group structure on a subgroup of $\Difforb{Q,\uU}$:

\begin{prop}\label{prop: lgp:uncp} \no{prop: lgp:uncp}
 There is an open subset $\pP \subseteq E(\mM) \subseteq \Difforb{Q,\uU}$ which contains the identity such that the subgroup generated by $\pP$,
	\begin{displaymath}
	 \Diffo{Q,\uU} \coloneq \langle\pP \rangle,
	\end{displaymath}
 admits a unique smooth manifold structure turning $\Diffo{Q,\uU}$ into a connected Lie group modeled on $\Osc{Q}$ and $\pP$ into an open connected identity-neighborhood.
\end{prop}

\begin{proof}
 Endow $E(\mM)$ with the unique smooth manifold structure turning $E \colon \mM \rightarrow E(\mM)$ into a diffeomorphism. Consider $\pP_0 \coloneq E (\hH)$ as an open submanifold of $E(\mM)$. Combining Lemma \ref{lem: comp:cor} and Lemma \ref{lem: osc:sm} the composition and inversion 
	\begin{align*}
	m \colon& \pP_0 \times \pP_0 \rightarrow E(\mM) , ([\hat{f}], [\hat{g}]) \mapsto [\hat{f}] \circ [\hat{g}] = E (\comp (E^{-1} ([\hat{f}]) , E^{-1}([\hat{g}])))\\
	\iota \colon& \pP_0 \rightarrow E(\mM), [\hat{f}] \mapsto [\hat{f}]^{-1} = E(\inv (E^{-1}([\hat{f}]))
	\end{align*}
 are smooth maps. Observe that by Proposition \ref{prop: os:glue} and definition of $m$ and $\iota$ the images are contained in $E(\mM)$. The set $\pP_0$ is an open identity-neighborhood on which inversion and group multiplication of $\Difforb{Q,\uU}$ are smooth. Hence the preimage $\iota^{-1}(\pP_0) = \pP_0 \cap (\pP_0)^{-1}$ with $(\pP_0)^{-1} \coloneq \iota (\pP_0)$ is an open neighborhood of the identity in $\pP_0$. Thus $E^{-1} (\pP_0 \cap (\pP_0)^{-1})$ is an open zero-neighborhood in $\Osc{Q}$. Since this space is locally convex, we may choose a convex zero neighborhood $\hH_1 \subseteq E^{-1} (\pP_0 \cap (\pP_0)^{-1}) \subseteq \Osc{Q}$. Then $\pP_1 \coloneq E (\hH_1) \subseteq \pP_0 \cap (\pP_0)^{-1}$ is a connected, open identity neighborhood in $E(\mM)$. Since $\pP_1 \subseteq \pP_0 \cap (\pP_0)^{-1}$ holds, we have $\iota^{-1} (\pP_1) = \pP_0 \cap (\pP_1)^{-1} = (\pP_1)^{-1} = \iota (\pP_1)$. Being a preimage of an open set with respect to a continuous map, $(\pP_1)^{-1}$ is open. Furthermore it is connected as a continuous image of such a set. We obtain an open, connected identity-neighborhood $\pP \coloneq \pP_1 \cup (\pP_1)^{-1} \subseteq \pP_0$ in $E(\mM)$ by \cite[Corollary 6.1.10]{Engelking1989}. \\
 From the above, we deduce $m(\pP,\pP) \subseteq E(\mM)$ and the mapping $\pP \times \pP \rightarrow E(\mM) , ( [\hat{f}], [\hat{g}]) \mapsto [\hat{f}] \circ [\hat{g}] $ induced by $m$ is a smooth map. Furthermore, $\pP^{-1} = \pP \subseteq E(\mM)$ holds and the mapping $\pP \rightarrow E(\mM) ,  [\hat{f}] \mapsto [\hat{f}]^{-1} $ induced by $\iota$ is smooth. In conclusion all prerequisites of Proposition \ref{prop: Lgp:locd} (a) have been checked. Hence we derive a unique smooth manifold structure on \glsadd{Diffo}
	\begin{displaymath}
	 \Diffo{Q,\uU} \coloneq \langle \pP\rangle
	\end{displaymath}
 turning it into a Lie group such that $\pP$ is an open identity-neighborhood in $\Diffo{Q,\uU}$. In addition the manifold structure induced by $\Diffo{Q,\uU}$ coincides with the submanifold structure of $\pP \subseteq E(\mM)$. Therefore, $\pP \subseteq \Diffo{Q,\uU}$ is open and connected. As the group operations of $\Diffo{Q,\uU}$ are smooth, each of the sets $\pP^n$ (the elements of $\Diffo{Q,\uU}$, which are obtained by $n$-fold composition of elements in $\pP$, $n \in \NN$) is a connected identity-neighborhood. Since $\pP$ is a symmetric identity-neighborhood, we deduce from the proof of \cite[Theorem 5.7]{topgroup1979}:
  \begin{displaymath}
   \Diffo{Q,\uU} = \langle \pP\rangle = \displaystyle\bigcup_{n=1}^\infty \pP^n.
  \end{displaymath}
 Hence $\Diffo{Q,\uU}$ is a connected Lie group by \cite[Corollary 6.1.10]{Engelking1989}. 
\end{proof}

In the next section, we shall construct a Lie group structure on $\Difforb{Q,\uU}$. The Lie group structure on the subgroup $\Diffo{Q,\uU}$ of $\Difforb{Q,\uU}$ will turn this subgroup into the identity component of the Lie group $\Difforb{Q,\uU}$.

\subsection{Lie group structure on \texorpdfstring{$\Difforb{Q,\uU}$}{the orbifold diffeomorphism group}}\label{sect: lgp:glob} \no{sect: lgp:glob}

Unless stated otherwise, all symbols used in this section retain the same meaning as in Section \ref{sect: lgp:ident}. In particular, we shall always be working with a Riemannian orbifold $(Q,\uU,\rho)$. First, we will prove that the Lie group $\Diffo{Q,\uU}$ is independent of the choice of the atlases $\aA, \bB$ and the local data constructed in Section \ref{sect: lgp:ident}. Second, the construction does not depend on the choice of the Riemannian orbifold metric on $(Q,\uU)$. Having dealt with these preparations, an application of the Construction Principle \ref{prop: Lgp:locd} will yield a unique smooth Lie group structure on $\Difforb{Q,\uU}$. The strategy of the proof follows \cite{hg2006b} where a similar argument has been used to turn the diffeomorphism group of a manifold into a Lie group. 

\begin{lem}\label{lem: indep:ch} \no{lem: indep:ch}
 The Lie group $\Diffo{Q,\uU}$ constructed in Proposition \ref{prop: lgp:uncp} does neither depend on the choice of atlases $\aA$ and $\bB$, nor on the local data collected in Construction \ref{con: sp:ofdat}.
\end{lem}

\begin{proof}
 Let  $\aA^+$ and $\bB^+$ be orbifold atlases which satisfy the same properties as $\aA$ and $\bB$ in Construction \ref{con: sp:ofdat}. Replace $\aA$ and $\bB$ in the construction of Section \ref{sect: lgp:ident} with $\aA^+$ and $\bB^+$. Taking the Riemannian orbifold metric $\rho$ as before, we obtain another connected, smooth Lie group $\Diffo{Q,\uU}^+$ depending on the new set of data. As shown in Section \ref{sect: lgp:ident}, there is a $C^\infty$-diffeomorphism $E^+$, $E^+([\hat{\sigma}]) \coloneq [\expo] \circ [\hat{\sigma}]$ mapping the open convex zero-neighborhood $\hH_1^+$ (defined as in Proposition \ref{prop: lgp:uncp} with respect to $\aA^+$ and $\bB^+$, the open subset $\hH^+ \subseteq \Osc{Q}$ and the local data constructed for $\aA^+$, $\bB^+$) onto an open identity neighborhood in $\Diffo{Q,\uU}^+$. Then $O \coloneq \hH_1 \cap \hH_1^+$ is an open, convex (and hence connected) zero-neighborhood in $\Osc{Q,\uU}$. The map $E$ takes $O$ diffeomorphically onto an open identity neighborhood in $\Diffo{Q,\uU}$. As $\Diffo{Q,\uU}$ is a connected Lie group, $E(O)$ generates this group by \cite[Theorem 7.4]{topgroup1979}. Analogously, $E^+$ maps $O$ diffeomorphically onto an open identity neighborhood in $\Diffo{Q,\uU}^+$ which generates this group. Recall from Proposition \ref{prop: diff:chart} that $E([\hat{\sigma}]) = [\expo] \circ [\hat{\sigma}]|^{\Omega} = E^+([\hat{\sigma}])$ holds for each $[\hat{\sigma}] \in O$. Hence both maps coincide on $O$. We deduce that $\Diffo{Q,\uU} = \langle E(O)\rangle = \Diffo{Q,\uU}^+$ as an abstract group and also as a Lie group.
\end{proof}


\begin{lem}\label{lem: indep:rm} \no{lem: indep:rm}
 The Lie group $\Diffo{Q,\uU}$ constructed in Proposition \ref{prop: lgp:uncp} does not depend on the choice of the Riemannian orbifold metric $\rho$ on $(Q,\uU)$ (cf. Section \ref{sect: lgp:ident}).
\end{lem}

\begin{proof}
 Let $\rho^\#$ be another Riemannian orbifold metric on $(Q,\uU)$. By Lemma \ref{lem: indep:ch} we may use the same atlases $\aA =\setm{(U_i,G_i,\psi_i)}{i \in I}$ and $\bB = \setm{(W_j,H_j,\varphi_j)}{j \in J}$ as in Construction \ref{con: sp:ofdat}. Reviewing this, the local data constructed in Construction \ref{con: sp:ofdat} II.\ - IV.\ do \textbf{not depend} on the Riemannian orbifold metric. The constants $R_i,\ i \in I$ and $s_j, S_j,\ j \in J$ in Construction \ref{con: sp:ofdat} V.\ change for $\rho^\# = (\rho_j^\#)_{j \in J}$. The new constants depending on $\rho^\#$ will be denoted by $R_i^\#, \ i \in I$ and $s_j^\#, S_j^\#, j \in J$ (see Construction \ref{con: sp:ofdat} V.\ for their properties).\\
 Let $[\widehat{\expo}^\#]$ be the Riemannian orbifold exponential map with respect to $(Q,\uU, \rho^\#)$. As in Section \ref{sect: lgp:ident}, one constructs open zero-neighborhoods $\hH^\# \coloneq \Lambda_\cC^{-1} (\oplus_{i \in I} \hH_{R_i^\#})$ and $\hH^\# \subseteq \mM^\#$, which depend on the data in Construction \ref{con: sp:ofdat} I. - IV., the constants $R_i^\#, \ i \in I$ and $s_j^\#, S_j^\#, j \in J$, as well as on the Riemannian orbifold metric $\rho^\#$. Furthermore, we obtain an injective map $E^\# \colon \mM^\# \rightarrow \Diffo{Q,\uU}^\#$, a connected Lie group $\Diffo{Q,\uU}^\# = \langle \pP^\#\rangle$ and a convex zero-neighborhood $\hH_0^\# \subseteq \hH^\# \subseteq \Osc{Q}$ such that $E^\#|_{\hH_0^\#} \colon \hH_0^\# \rightarrow \pP^\# \subseteq \Diffo{Q,\uU}^\# , [\hat{\sigma}] \mapsto [\widehat{\expo}^\#] \circ [\hat{\sigma}]|^{\Omega^\#}$ is a diffeomorphism onto an open identity neighborhood.\\
 Fix some $i \in I$ and let $\fF_5 (K_{5,i})= \setm{(V_{5,\alpha (i)}^n, \kappa_n^{\alpha (i)})}{1\leq n \leq N_i}$ be the atlas of Construction \ref{con: sp:ofdat} IV.\footnote{To shorten our notation, we number all charts from $1$ to some $N_i \in \NN$, $i \in I$. It will always be clear from the context which charts are meant.}. For each $1 \leq n \leq N_i$ the Riemannian metrics induce pullback metrics with respect to the manifold charts $\kappa^{\alpha (i)}_n$. The charts $\kappa_n^{\alpha (i)}$ induce pullback metrics on $B_5(0)$ with respect to $\rho_{\alpha (i)}$ and $\rho_{\alpha (i)}^\#$. For $(V_{5,\alpha (i)}^n, \kappa_n^{\alpha (i)}), 1\leq n \leq N_i$ the associated Riemannian exponential maps will be denoted by $\exp_{W_{\alpha (i)},[n]}$ and $\exp_{W_{\alpha (i)},[n]}^\#$, respectively. Finally we define the local representatives of $X \in \vect{W_{\alpha (i)}}$ with respect to $\kappa_n^{\alpha (i)}$ via $X_{[n]} \coloneq X_{\kappa_{n}^{\alpha (i)}} \circ (\kappa_n^{\alpha (i)})^{-1} \in C^\infty( B_5(0) ,\RR^d)$.
 \\[2em] Observe that the open set $\hH_{R_i}$ in Construction \ref{con: H} was obtained by Construction \ref{con: comp:loc}. Reviewing Construction \ref{con: comp:loc}, for $1 \leq n \leq N_i$, real numbers $\ve_n , \delta_n >0$ have been chosen such that for each $x \in \overline{B_4 (0)}$, the map $\phi_{\alpha (i), [n], x} \colon B_{\ve_n} (0) \rightarrow \RR^d, y \mapsto \exp_{W_{\alpha (i)}, [n]} (x,y)$ is a diffeomorphism onto its open image which contains $B_{\delta_n}(0)$. Furthermore, by Lemma \ref{lem: exp:loc} the choice of $\ve_n$ yields the smooth map $b_{\alpha (i),[n]} \colon W_{\delta_n} \rightarrow B_{\ve_n} (0), b_{\alpha (i),[n]}(x,y) \coloneq \phi_{\alpha (i),[n],x}^{-1} (y)$. Recall that $\ve_n < \nu_i$ for $1\leq n \leq N_i$. Here $\nu_i$ is the constant constructed in Lemma \ref{lem: eos:cov} with respect to the finite family $\fF_5 (K_{5,i})$. Thus the assertions of Lemma \ref{lem: eos:cov} hold. For each $x \in V_{4,\alpha (i)}^n, 1\leq n \leq N_i$, there is an open set $N_x \subseteq T_x W_{\alpha (i)}$ with the following property: 
  \begin{equation}\label{eq:inc}
     B_{\delta_n} (\kappa_n^{\alpha (i)} (x)) \subseteq \exp_{W_{\alpha (i)}, [n]} (\kappa_n^{\alpha (i)} (x), B_{\ve_n} (0)) \subseteq \kappa_n^{\alpha (i)} \exp_{W_{\alpha (i)}} (N_x).                                                                                                                                                                   
  \end{equation}
Observe that the neighborhood $\hH_{R_i^\#}$ has been obtained by another application of Construction \ref{con: comp:loc} with respect to a family of constants $\ve_n^\#, \delta_n^\# >0$ for $1\leq n \leq N_i$. \\ 
 By Lemma \ref{lem: exp:loc} (c), we may choose constants $\ve_n^\# > \ve_{1,n}^\# >0$ for $1 \leq n \leq N_i$ so small that $\exp^\#_{\alpha (i), [n]} \left(\set{\kappa_n^{\alpha (i)} (x)} \times B_{\ve_{1,n}^\#} (0)\right)$ is contained in $B_{\delta_n} (\kappa_n^{\alpha (i)} (x))$ for $x \in \overline{V_{4,\alpha (i)}^n}$. For $1\leq n\leq N_i$, we choose for each $\ve_{1,n}^\#$ a constant $\delta_ n^\# > \delta_{1,n}^\# >0$ which satisfies the assertion of Lemma \ref{lem: exp:loc} (b), with $\ve$ replaced with $\ve^\#_{1,n}$. Apply Construction \ref{con: comp:loc} with $R \coloneq R_i^\#$ and $P\coloneq P_{1,i}^\# \cap P_{2,i}^\#$, but replace the pairs $(\ve_n^\#,\delta_n^\#)$ with $(\ve_{1,n}^\#, \delta_{1,n}^\#)$ to obtain an open zero-neighborhood $H_{R_i^\#} \subseteq \hH^{\Omega_{5,K_{5,i}}}_{R_i^\#}$. Thus the map
	\begin{equation}\label{eq:un}
	 u_n \colon B_4 (0) \times B_{\ve_{1,n}^\#} (0) \rightarrow B_{\ve_n} (0) , u_n (x,y) \coloneq b_{\alpha (i),[n]} (x, \exp^\#_{W_{\alpha (i)}, [n]} (x,y))
	\end{equation}
 makes sense and is smooth as a composition of smooth maps. By construction, $\ve_{1,n}^\# < \ve_n^\# < \nu^\#$, where $\nu^\#$ is the constant as in Lemma \ref{lem: eos:cov} with respect to the finite family $\fF_5 (K_{5,i})$. Hence we deduce with Lemma \ref{lem: eos:cov} (b) from equations \eqref{eq:un} and \eqref{eq:inc} that the map
	\begin{equation}\label{eq:nvf}
	 (E^{-1} E^\#)_i \colon H_{R_i^\#} \rightarrow \vect{\Omega_{1,K_{5,i}}}, (E^{-1} E^\#)_i (X)(x) \coloneq \exp_{W_{\alpha (i)}}|_{N_x}^{-1} \exp_{W_{\alpha (i)}}^\# \circ X (x)  
	\end{equation}
 makes sense. In addition, we show that $(E^{-1} E^\#)_i$ is a smooth map. To see this, let $1\leq n \leq N_i$ and recall that $H_{R_i^\#} \subseteq \vect{\Omega_{5,K_{5,i}}}$ is open and $\fF_5 (K_{5,i})$ covers $\Omega_{5,K_{5,i}}$. Hence for $1\leq n \leq N_i$, the maps $r_n \colon \vect{\Omega_{5,K_{5,i}}} \rightarrow C^\infty (B_5 (0), \RR^d) , X \mapsto X_{[n]}$ form a patchwork by Definition \ref{defn: top:vect}. Analogously, the maps $t_n \colon \vect{\Omega_{1,K_{5,i}}} \rightarrow C^\infty (B_1 (0), \RR^d), X \mapsto X_{[n]}|_{B_{1} (0)}$ yield a patchwork for $1\leq n\leq N_i$. Consider the open subset $\lfloor \overline{B_1 (0)},B_{\ve_{1,n}^\#} (0)\rfloor_\infty \subseteq C^\infty (B_5 (0), \RR^d)$. For $X\in H_{R_i^\#}$ we obtain $X_{[n]} (\overline{B_3 (0)}) \subseteq B_{\ve_{1,n}^\#} (0)$ (cf.\ Construction \ref{con: comp:loc} and Lemma \ref{lem: mb:exp}). Hence $r_n (H_{R_i^\#}) \subseteq \lfloor \overline{B_1 (0)}, B_{\ve_{1,n}^\#} (0)\rfloor_\infty$ holds. In addition, \cite[Proposition 4.23 (a)]{hg2004} with \eqref{eq:un} yields a smooth map 
  \begin{displaymath}
   U_n \colon \lfloor \overline{B_1 (0)}, B_{\ve_{1,n}^\#} (0)\rfloor_\infty \rightarrow C^\infty (B_1 (0), \RR^d), U_n (\sigma) \coloneq (u_n)_* (\sigma), 
  \end{displaymath}
 with $(u_n)_* (\sigma) (x) \coloneq u_n (x, \sigma (x))$ for $x \in B_1 (0)$. By \eqref{eq:un}, $U_n$ maps the zero-map to the zero-map. Evaluating \eqref{eq:un} pointwise for $(X,x) \in H_{R_i^\#} \times \Omega_{1,K_{5,i}}$, the local formula \eqref{eq:un} and Lemma \ref{lem: eos:cov} (b) yield the identity $t_n \circ (E^{-1} E^\#)_i = U_n \circ r_n$. Thus $(E^{-1}E^\#)_i$ is a patched mapping which is smooth on the patches, whence $(E^{-1} E^\#)_i$ is smooth by Proposition \ref{prop: pat:loc}.\\
 For each $j \in I$, construct in the same manner an open set $H_{R_j^\#} \subseteq \vect{\Omega_{5,K_{5,j}}}$ together with a smooth map $(E^{-1} E^\#)_j$. Define $H_i^\# \coloneq (\res^{W_{\alpha (i)}}_{\Omega_{5,i}})^{-1} (H_{R_i^\#}) \subseteq \hH_{R_i^\#} \subseteq \vect{W_{\alpha (i)}}$. By Construction \ref{con: H}, $H^\# \coloneq \Lambda_\cC^{-1} (\oplus_{i\in I} H_i^\#) \subseteq \hH^\#$ holds. For each $[\hat{\sigma}] \in H^\#$, the family $\set{(E^{-1} E^\#)_i (\sigma_{\alpha (i)}|_{\Omega_{5,K_{5,i}}})|_{U_i}}_{i \in I}$ is a family of vector fields. Since $[\hat{\sigma}]$ is compactly supported, only finitely many canonical lifts $\sigma_{\alpha (i)}$ are non-zero. By standard Riemannian geometry, the Riemannian exponential map composed with the zero section yields the identity. Hence \eqref{eq:nvf} shows that only finitely many of the vector fields $\set{(E^{-1} E^\#)_i( \sigma_{\alpha (i)}|_{\Omega_{5,K_{5,i}}})|_{U_i}}_{i \in I}$ will be non-zero. We claim that these vector fields form a canonical family of an orbisection. If this is true, then these vector fields define a compactly supported orbisection $E^{-1}E^\# ([\hat{\sigma}])$, whose lifts with respect to $\aA$ are given by $\set{(E^{-1} E^\#)_i (\sigma_{\alpha (i)}|_{\Omega_{5,K_{5,i}}})|_{U_i}}_{i \in I}$. On $U_i \subseteq \Omega_{1,i}$, these vector fields yield an orbisection if the following is satisfied:\\[0.5em]
 \emph{Let $[\hat{\sigma}] \in H^\#$ and $\lambda \in \CH{W_k}{W_l}$ be a change of charts which satisfies $\dom \lambda \subseteq \Omega_{1, i}$ and $\cod \lambda \subseteq \Omega_{1,j}$ for some  $k = \alpha (i)$ and $l = \alpha (j)$. Then the following identity holds:
 \begin{equation}\label{eq: os}
    T\lambda  (E^{-1} E^\#)_i (\sigma_k|_{\Omega_{5,K_{5,i}}})|_{\dom \lambda} = (E^{-1} E^\#)_j (\sigma_l|_{\Omega_{5,K_{5,j}}}) \circ \lambda .
 \end{equation}}
 The argument given in the proof of Lemma \ref{lem: vf:eq} may be repeated almost verbatim. We check the identity \eqref{eq: os} locally: Choose some $x \in \dom \lambda \subseteq \Omega_{1,i}$ and a chart $(V_{5,k}^n,\kappa_n^k)\in \fF_5 (K_{5,i})$ with $x \in V_{1,k}^n$. Again there is some $Z^r_k$ with $V_{5,k}^n \subseteq Z_k^r$. As $x \in V_{1,k}^n \subseteq \kK_k^\circ$ and $\lambda (x) \in \kK_l$, there is an embedding of orbifold charts $\mu \colon Z_k^r \rightarrow W_l$ with $\mu (x) = \lambda (x)$. After possibly composing $\mu$ with a suitable element of $H_l$, there is an open neighborhood $U_x$ of $x$ in $Z_k^r$ with $\mu|_{U_x} = \lambda|_{U_x}$ and thus $T_x \mu = T_x \lambda$.  \\
 Since $\rho$ and $\rho^\#$ are Riemannian orbifold metrics, each change of orbifold charts in $\CH{W_k}{W_l}$ is a Riemannian embedding of its domain endowed with the induced metrics into the Riemannian manifold $(W_l,\rho_l)$ respectively $(W_l, \rho_l^\#)$. By construction of $H_i^\#$, each $X \in H_i^\#$ satisfies 
  \begin{equation}\label{eq:esti:vf}
   \norm{\phi_{k , [n],x}^{-1} \exp_{W_k,[n]}^\# X_{[n]}}_{\overline{B_1 (0)},0} < \ve_n < R_i \text{ for each } 1 \leq n \leq N_i.
  \end{equation}
 Recall from Construction \ref{con: sp:ofdat} V.\ the properties of $R_i$ and $S_k$:\\
 The definitions imply that  $T\mu (E^{-1}E^\#)_i (\sigma_k|_{\Omega_{5,K_{5,i}}}) (V_{1,k}^n) \subseteq \hat{O}_l \subseteq \dom \exp_{W_l}$ for $[\hat{\sigma}] \in H^\#$. Computing locally on $V_{5,k}^n$, we use that $\mu (\kappa_n^k)^{-1}$ is a Riemannian embedding into $W_l$. Again by \cite[IV.\ Proposition 2.6]{fdiffgeo1963}, the identity $\exp_{W_l} T(\mu (\kappa_n^k)^{-1}) (v) = \mu (\kappa_n^k)^{-1} \exp_{W_k,[n]} (v)$ holds for each $v \in \dom \exp_{W_k,[n]}$. The family $\set{\sigma_k}_{k \in J}$ is a canonical family of lifts, whence $\sigma_l \mu = T\mu \sigma_k|_{\dom \mu}$. By definition of $H_i^\# \subseteq \hH_{R_i^\#}$, the identity $\kappa_n^k \exp_{W_k}^\# \circ X (z) = \exp_{W_k,[n]}^\# T\kappa_n X(z)$ holds for each $z \in V_{3,k}^n$ and $X \in H^\#_i$ (cf.\ the proof of Lemma \ref{lem: vf:eq}). Observe that $\lambda (x) \in  \Omega_{1,j}$ and $\sigma_l \in H_j^\#$. Combining these facts one computes: 
	\begin{align*}
	  \mbox{}&\exp_{W_l} T_x\lambda (E^{-1} E^\#)_i (\sigma_k|_{\Omega_{5,K_{5,i}}}) (x) =\exp_{W_l} T_x(\mu (\kappa_n^k)^{-1} \kappa_n^k) (\exp_{W_k}|_{N_x})^{-1} \exp_{W_k}^\# \sigma_k (x)\\
	  &=\mu (\kappa_n^k)^{-1} \exp_{W_k , [n]} T\kappa_n^k (\exp_{W_k}|_{N_{x}})^{-1} \exp_{W_k}^\# \sigma_k (x) \stackrel{\ref{lem: eos:cov} (b)}{=}\mu \exp_{W_k}^\# \sigma_k (x) \\
	  &=\mu (\kappa_n^k)^{-1} \exp_{W_k, [n]}^\# T \kappa_n^k \sigma_k (x) = \exp_{W_l}^\# \sigma_l (\mu (x)) = \exp_{W_l}^\# \sigma_l (\lambda (x))\\
          &\hspace{-0.28cm}\stackrel{\eqref{eq:nvf}}{=}\exp_{W_l} (E^{-1} E^\#)_j (\sigma_l|_{\Omega_{5,K_{5,j}}}) (\lambda (x)) .
   	\end{align*}
 As $x \in \kK_k^\circ$ and $\lambda (x) \in \Omega_{1,K_{5,j}}$, the definition of $R_i$ implies $T_x\lambda (E^{-1} E^\#)_k (\sigma_k|_{\Omega_{5,K_{5,i}}}) (x) \in \hat{O}_l$. By construction of $H_j^\#$, we deduce $(E^{-1} E^\#)_j (\sigma_l|_{\Omega_{5,K_{5,j}}}) \lambda (x) \in \hat{O}_l$. As $\exp_{W_l}$ is injective on $T_{\lambda (x)} W_l \cap \hat{O}_l$ and $x \in \dom \lambda$ was arbitrary, this proves \eqref{eq: os}. We conclude that the family of vector fields $\set{(E^{-1} E^\#)_i (\sigma_{\alpha (i)}|_{\Omega_{5,K_{5,i}}})|_{U_i}}_{i \in I}$ is a family of canonical lifts for a compactly supported orbisection $E^{-1}E^\# ([\hat{\sigma}])$. Define $E^{-1} E^\# \colon H^\# \rightarrow \Osc{Q}, [\hat{\sigma}] \mapsto E^{-1}E^\# ([\hat{\sigma}])$. Using the patchworks $(\lambda_i)_{i \in I}$ and $(\tau_i)_{i \in I}$ for $\Osc{Q}$ (see proof of Lemma \ref{lem: osc:sm}), a computation yields the identity 
	\begin{displaymath}
	 \res^{\Omega_{1, K_{5,i}}}_{U_i} (E^{-1}E^\#)_i \res^{W_{\alpha (i)}}_{\Omega_{5,K_{5,i}}}\lambda_i|^{H_i^\#}_{H^\#} = \tau_i E^{-1}E^\# ,\quad i \in I.
	\end{displaymath}
 We have already seen that $(E^{-1}E^\#)_i$ is smooth and $(E^{-1}E^\#)_i (0_{\alpha (i)}) = 0_i$ for each $i \in I$. By \cite[Lemma F.15 (a)]{hg2004}, the mappings $\res^{\Omega_{1, K_{5,i}}}_{U_i}$, $\res^{W_{\alpha (i)}}_{\Omega_{5,K_{5,i}}}$ are smooth, whence $E^{-1}E^\#$ is a patched mapping which is smooth on the patches. By Proposition \ref{prop: pat:loc}, $E^{-1}E^\#$ must be smooth and therefore it is continuous. Using continuity, there is an open, connected zero-neighborhood $\rR^\# \subseteq \hH_0^\# \cap H^\#$ such that $E^{-1} E^\# (\rR^\#) \subseteq E^{-1} (\pP)$. Uniqueness of canonical lifts proves that the canonical lifts of $E^{-1}E^\# ([\hat{\sigma}])$ on $W_{\alpha (i)}$ coincides on $U_i$ with $(E^{-1}E^\#)_i(\sigma_{\alpha (i)}|_{\Omega_{5,K_{5,i}}})|_{U_i}$. Recall the construction of the representative $\hat{E}^\sigma$ of $E([\hat{\sigma}])$ in Proposition \ref{prop: eos:imdiff}. Using \eqref{eq:nvf}, the construction yields for $E(E^{-1}E^\# ([\hat{\sigma}]))$ and $i \in I$ the lifts $\exp_{W_{\alpha (i)}}^\# \circ \sigma_{i}$. The same lifts are obtained, if this construction is carried out with respect to the Riemannian orbifold exponential map $[\expo^\#]$. As orbifold diffeomorphisms are uniquely determined by a family of lifts (cf.\ Corollary \ref{cor: ofdiff:un}), $E^\# ([\hat{\sigma}]) = E \circ (E^{-1}E^\# )([\hat{\sigma}]) \in E (E^{-1}(\pP)) = \pP$ holds for each $[\hat{\sigma}] \in \rR^\#$. The set $\rR^\#$ is an open and connected zero-neighborhood contained in $\hH_0^\#$. Since $\Diffo{Q,\uU}^\#$ is connected, $\langle E^\# (\rR^\#) \rangle = \Diffo{Q,\uU}^\#$ holds by \cite[Theorem 7.4]{topgroup1979}, which implies $\Diffo{Q,\uU}^\# \subseteq \Diffo{Q,\uU}$. In particular, the inclusion morphism $\Diffo{Q,\uU}^\# \rightarrow \Diffo{Q,\uU}$ is smooth on the open identity-neighborhood $E^\# (\rR^\#)$, hence smooth by \cite[Ch.\ III, \S 1, No.\ 2, Proposition 4]{bourbaki1989}. Reversing the roles of $\rho$ and $\rho^\#$, one deduces that also $\Diffo{Q,\uU} \subseteq \Diffo{Q,\uU}^\#$ and the inclusion morphism $\Diffo{Q,\uU} \rightarrow \Diffo{Q,\uU}^\#$ is smooth. In conclusion, $\Diffo{Q,\uU}$ and $\Diffo{Q,\uU}^\#$ coincide as Lie groups. 
\end{proof}

So far, we achieved that the Lie group structure on $\Diffo{Q,\uU}$ does neither depend on the local data (the atlases $\aA$, $\bB$ etc.) nor on the Riemannian orbifold metric. We exploit these facts to prove that the requirements of Proposition \ref{prop: Lgp:locd} (b) are satisfied:

\begin{prop}\label{prop: diffo:nor} \no{prop: diffo:nor}
 Let $[\hat{\phi}] \in \Difforb{Q,\uU}$ be an arbitrary orbifold diffeomorphism. Then for each $[\hat{f}] \in \Diffo{Q,\uU}$ we have $ [\hat{\phi}] \circ [\hat{f}] \circ [\hat{\phi}]^{-1} \in \Diffo{Q,\uU}$  and  
	\begin{displaymath}
	 c_{[\hat{\phi}]} \colon \Diffo{Q,\uU} \rightarrow \Diffo{Q,\uU}, [\hat{f}] \mapsto [\hat{\phi}] \circ [\hat{f}] \circ [\hat{\phi}]^{-1}
	\end{displaymath}
 is a smooth map. 
\end{prop}

\begin{proof}
 The proof will be quite simple, after some preparations:\\
 Following Corollary \ref{cor: diff:char} (d), we may choose orbifold atlases $\vV_i \coloneq \setm{(V^i_k,L^i_k, \pi^i_k) \in \uU}{k \in K} \subseteq \uU$, $i \in \set{1,2}$ together with a representative $\Phi = (\phi , \set{\phi_k}_{k \in K}, P, \nu) \in \Orb{\vV_1,\vV_2}$ of $[\hat{\phi}]$ such that each $\phi_k \colon V^1_k \rightarrow V_k^2$ is a diffeomorphism. Furthermore, Corollary \ref{cor: ll:diff} assures that we may choose $P = \Ch_{\vV_1}$ and $\nu (\lambda) = \phi_l \lambda \phi_k^{-1}|_{\phi_k (\dom \lambda)}$ for $\lambda \in \CH{V^1_k}{V^1_l}$.\\
 By Proposition \ref{prop: atl:lfpts} there are locally finite atlases $\aA$ and $\bB$ indexed by $I$ and $J$, respectively, which satisfy the properties of the atlases in Construction \ref{con: sp:ofdat} I. In addition, there is a map $\beta \colon J \rightarrow K$ such that $W_j$ is an open subset of $V^1_{\beta (j)}$, the inclusion of sets induces an embedding of orbifold charts and $\overline{W_j} \subseteq V_{\beta (j)}^1$ is compact for each $j \in J$. As a consequence of Lemma \ref{lem: indep:ch}, we may construct $\Diffo{Q,\uU}$ with respect to these atlases and the Riemannian orbifold metric $\rho$. Thus there are open sets $\hH_1 \subseteq \hH \coloneq \Lambda_\cC^{-1} (\bigoplus_{i \in I}\hH_{R_i})$ and a diffeomorphism $E|_{\hH_1}$ onto an identity neighborhood in $\Diffo{Q,\uU}$. \\
 By construction, the inclusions of sets $U_i \subseteq W_{\alpha (i)} \subseteq V_{\beta (\alpha (i))}^1$ and $\phi_{\beta (\alpha (i))}$ are changes of orbifold charts for each $i \in I$. For $i \in I$, the sets $W_{\alpha (i)}^+ \coloneq \phi_{\beta (\alpha (i))} (W_{\alpha (i)})$ and $U_i^+ \coloneq \phi_{\beta \alpha (i)} (U_i)$ are $L^2_{\beta (\alpha (i))}$-stable, open and relatively compact subsets of $V_{\beta (\alpha (i))}^2$ (cf.\ Lemma \ref{lem: add:cor} (a)). Define the following sets of orbifold charts for $Q$:
	\begin{displaymath}
	 \aA^+ \coloneq \setm{(U_i^+, G_i , \pi^2_{\beta \alpha (i)}|_{U_i^+})}{i \in I} \text{ and } \bB^+ \coloneq \setm{(W^+_j, H_j , \varphi^+_j \coloneq \pi^2_{\beta (j)}|_{W_j^+})}{j \in J}.
	\end{displaymath}
 The underlying map $\phi$ is a homeomorphism and each $\phi_k$ is a diffeomorphism. Hence $\aA^+$ and $\bB^+$ are orbifold atlases for $Q$ such that $\overline{U_i^+} \subseteq W_{\alpha (i)}^+$ for each $i \in I$ and the inclusions of sets induce embeddings of orbifold charts. Since $W_j^+$ is a relatively compact subset of $V_{\beta (j)}^2$ for each $j \in J$, we deduce from the continuity of $\pi_{\beta (j)}^2$ and \cite[Corollary 3.1.11]{Engelking1989} that the image of each chart in $\aA^+$ and $\bB^+$ is relatively compact. Exploiting that $\phi$ is a homeomorphism, $\aA^+$ and $\bB^+$ are locally finite atlases, since the same holds for $\aA$ and $\bB$. Furthermore, by construction of $\aA$ and $\bB$, for each connected component $C \subseteq Q$, there is a point $z_C$ which is only contained in the images of a unique pair of charts in $\aA \times \bB$. The homeomorphism $\phi$ permutes the connected components of $Q$, whence each $z_C$ is mapped into a separate component. Each element of $\setm{\phi (z_C)}{ C \subseteq Q \text{ a connected component}}$ is thus contained in the images of a unique pair in $\aA^+ \times \bB^+$ such that the images of different pairs are contained in different connected components. Summing up, the atlases $\aA^+$ and $\bB^+$ satisfy all properties required in Construction \ref{con: sp:ofdat} I. 
 \paragraph{} As $\bB$ is an atlas, a family of lifts for a representative of $[\hat{\phi}]$ is given by $\set{\Phi_j \coloneq \phi_{\beta (j)}|_{W_j}}_{j \in J}$. By construction, each of these lifts is a diffeomorphism and $\Phi_{\alpha (i)} (U_i) = U_i^+$ for each $i \in I$. Corollary \ref{cor: diff:char} assures that $\set{\Phi^{-1}_{j}}_{j \in J}$ is a family of lifts for a representative of $[\hat{\phi}]^{-1}$ in $\Orb{\bB^+,\bB}$. Observe that $\Phi^{-1}_j (U_i^+) = U_i$ holds for each $i \in \alpha^{-1} (j)$. 
 Before we prove the smoothness of $c_{[\hat{\phi}]}$, consider the following auxiliary maps:\\  
 Define $t_i \colon \hH_{R_i} \rightarrow \vect{U_i^+}, X \mapsto T\Phi_{\alpha (i)} X \Phi^{-1}_{\alpha (i)}|_{U_i^+}$ for $i\in I$. For $[\hat{\sigma}] \in \hH$, the family $\set{t_i(\sigma_i)}_{i \in I}$ defines a family of vector fields. We show that these vector fields are a family of canonical lifts of an orbisection: Let $\lambda \in \CH{U_i^+}{U_j^+}$ be any change of charts with arbitrary $i,j \in I$. As noted above, $\mu \coloneq \Phi_{\alpha (j)}^{-1} \lambda \Phi_{\alpha (i)}|_{\Phi_{\alpha (i)}^{-1} (\dom \lambda)}$ is a change of charts in $\CH{U_i}{U_j}$ and we compute 
	\begin{align*}
	 t_j (\sigma_j) \circ \lambda &= T\Phi_{\alpha (j)}\sigma_j \Phi_{\alpha (j)}^{-1}|_{U_j^+} \lambda = T\Phi_{\alpha (j)}\sigma_j \mu \Phi^{-1}_{\alpha (i)}|_{\dom \lambda}\\
				      & = T\Phi_{\alpha (j)} T\mu \sigma_i  \Phi^{-1}_{\alpha (i)}|_{\dom \lambda} = T\lambda t_i(\sigma_i)|_{\dom \lambda} .
	\end{align*}
 The family $\set{t_i(\sigma_i)}_{i \in I}$ is a family of canonical lifts with respect to $\aA^+$, whence it induces a unique orbisection $t([\hat{\sigma}])$. By construction, $t_i (\sigma_i)$ will be the zero-section if $\sigma_i$ is the zero-section. Hence $t([\hat{\sigma}])$ is compactly supported and we obtain a map $t \colon \hH \rightarrow \Osc{Q}, [\hat{\sigma}] \mapsto t([\hat{\sigma}])$. Consider the patchwork induced by the maps 
  \begin{displaymath}
   p_i \colon \Osc{Q} \rightarrow \vect{W_{\alpha (i)}}, p_i ([\hat{\sigma}]) = \sigma_{\alpha (i)} \text{ and } q_i \colon \Osc{Q} \rightarrow \vect{U_i^+}, q_i([\hat{\sigma}]) = \sigma_{U_i^+}  , \ i \in I,
  \end{displaymath}
 sending an orbisection to their canonical lifts. By construction of $\hH$ (cf. Construction \ref{con: H}), $p_i (\hH) \subseteq \hH_{R_i}$ holds. From $t_i \circ p_i|_{\hH}^{\hH_{R_i}} = q_i \circ t$ we deduce that $t$ is a patched mapping. 
 We claim that $t_i$ is smooth for each $i \in I$. If this were true, this implies the smoothness of $t$ by Proposition \ref{prop: pat:loc}.\\ To prove the claim, consider $t_i' \colon \hH_{R_i}^{\Omega_{5,K_{5,i}}} \rightarrow \vect{\Phi_{\alpha (i)} (\Omega_{5,K_{5,i}})}, X \mapsto T\Phi_{\alpha (i)} X \Phi^{-1}_{\alpha (i)}|_{\Phi_{\alpha (i)} (\Omega_{5,K_{5,i}})}$ and note the identity $t_i = \res^{\Phi_{\alpha (i)} (\Omega_{5,K_{5,i}})}_{U_i^+} t_i' \res_{\Omega_{5,K_{5,i}}}^{W_{\alpha (i)}}$. Since the restriction maps are smooth, it suffices to prove the smoothness of $t_i'$. By construction $\Omega_{5,K_{5,i}}$ is covered by the finite family of manifold charts $\fF_5 (K_{5,i}) = \setm{(V_{5,\alpha (i)}^n,\kappa_n^{\alpha (i)})}{1\leq n \leq N_i}$. Hence the sets $V_{5,\alpha (i)}^{n,+} \coloneq \Phi_{\alpha (i)} (V_{5,\alpha (i)}^n)$ cover $\Phi_{\alpha (i)} (\Omega_{5,K_{5,i}})$. Set $\gamma_n^{\alpha (i)} \coloneq \kappa_n^{\alpha (i)} \Phi^{-1}_{\alpha (i)}|_{V_{5,\alpha (i)}^{n,+}}$ to obtain a manifold atlas for $\Phi_{\alpha (i)} (\Omega_{5,K_{5,i}})$: $\fF^+_5(K_{5,i}) \coloneq \setm{(V_{5,\alpha (i)}^{n,+}, \gamma_n^{\alpha (i)})}{1\leq n\leq N_i}$. By Definition \ref{defn: top:vect} there are finite families of linear continuous mappings $\theta^n_{\alpha (i)} \colon \vect{\Omega_{5,K_{5,i}}} \rightarrow C^\infty (V_{5,\alpha (i)}^n, \RR^d), X \mapsto X_{\kappa_n}$ and $\theta^{n,+}_{\alpha (i)} \colon \vect{\Phi_{\alpha (i)}\Omega_{5,K_{5,i}}} \rightarrow C^\infty (V^{n,+}_{5,\alpha (i)}, \RR^d), Y \mapsto Y_{\gamma_n}$, with $1 \leq n \leq N_i$. The family $(\theta^n_{\alpha (i)})_{1\leq n \leq N_i}$ is a patchwork for $\vect{\Omega_{5,K_{5,i}}}$ and $(\theta^{n,+}_{\alpha (i)})_{1\leq n\leq N_i}$ is a patchwork for $\vect{\Phi_{\alpha (i)} (\Omega_{5,K_{5,i}})}$ by Lemma \cite[Lemma F.6]{hg2004}. As $\Phi_{\alpha (i)}^{-1}$ is smooth, the pullback $C^\infty (\Phi_{\alpha (i)}^{-1}|^{V_{5,\alpha (i)}^n}_{V_{5,\alpha (i)}^{n,+}}, \RR^d)$ is continuous linear and therefore smooth by \cite[Lemma 3.7]{hg2002}. A quick computation yields for $1\leq n \leq N_i$ the identity $\theta^{n,+}_{\alpha (i)} \circ t_i' = C^\infty (\Phi_{\alpha (i)}^{-1}|^{V_{5,\alpha (i)}^n}_{V_{5,\alpha (i)}^{n,+}}, \RR^d) \circ \theta^n_{\alpha (i)}$. We conclude that $t_i$ is a patched mapping, which is smooth on the patches, whence smooth by Proposition \ref{prop: pat:loc}.
 \paragraph{} The orbifold diffeomorphism $[\hat{\phi}]^{-1}$ induces a unique pullback metric $\rho^\# \coloneq ([\hat{\phi}]^{-1})^*\rho$ on $Q$ (cf. Lemma \ref{lem: rpbm:odi}). Denote by $\rho_j$ the members of $\rho$ on the orbifold charts $(W_j,H_j,\varphi_j),\ j \in J$. The Riemannian metric associated to $\rho^\#$ with respect to $(W^+_j,H_j,\varphi^+_j),\ j \in J$ are given by the pullback metric $\rho^\#_j \coloneq (\Phi_j^{-1})^*\rho_j$. For $j \in J$ let $\exp_j \colon D_j \rightarrow W_j$ be the Riemannian exponential maps with respect to $(W_j,\rho_j)$ and $\exp_j^\# \colon D_j^\# \rightarrow W_j^+$ be the exponential map with respect to $(W_j^+,\rho_j^\#)$. These pullback metrics turn $\Phi_j , \Phi_j^{-1}$ into Riemannian isometries and the map $[\hat{\phi}]$ into an orbifold isometry. In particular we derive $T\Phi_j (D_j) = D_j^\#$ and the exponential identity
	\begin{displaymath}
	 \exp_j^\# (T\Phi_j)|^{D_j^\#}_{D_j} = \phi_j \exp_j.
	\end{displaymath}
 Let $[\hat{\sigma}]$ be in $\hH$ and consider $e^{\sigma_i} $ as in Proposition \ref{prop: eos:imdiff}. From the last identity we deduce 
	\begin{align}\label{eq: conj:vf}
	\Phi_{\alpha (i)} \circ e^{\sigma_i} \circ \Phi_{\alpha (i)}|_{ U_i^+} = \Phi_{\alpha (i)} \exp_{\alpha (i)} \sigma_i \Phi^{-1}_{\alpha (i)}|_{U_i^+} = \exp_{\alpha (i)}^\# T\Phi_{\alpha (i)} \sigma_i \Phi^{-1}_{\alpha (i)}|_{U_i^+}.
	\end{align}
 Combining Lemma \ref{lem: indep:ch} with Lemma \ref{lem: eos:loc}, one may construct $\Diffo{Q,\uU}$ with respect to the atlases $\aA^+$, $\bB^+$ and the Riemannian orbifold metric $\rho^\#$. Hence there are an open connected zero-neighborhood $H^+_\# \subseteq \Osc{Q}$ and a map $E^+_\# \colon H^+_\# \rightarrow \Diffo{Q,\uU}, [\hat{\sigma}] \mapsto [\expo^\#]\circ [\hat{\sigma}]|^{\Omega^\#}$. Here $ [\expo^\#]$ is the Riemannian orbifold exponential map  associated to $\rho^\#$, whose domain is $\Omega^\#$. The map $E^+_\#$ is a diffeomorphism onto its image, which is an open identity-neighborhood in $\Diffo{Q,\uU}$. As $t$ is smooth and thus continuous, there is an open connected zero-neighborhood $A \subseteq \hH_1$ such that $t(A) \subseteq H^+_\#$. \\
 Recall from Corollary \ref{cor: ofdiff:un} that an orbifold diffeomorphism is uniquely determined by the lifts of any of its representatives. Hence for $[\hat{\sigma}] \in \hH_1 = E^{-1} (\pP)$ (cf.\ Proposition \ref{prop: lgp:uncp}), the orbifold diffeomorphism $[\hat{\phi}] \circ E([\hat{\sigma}]) \circ [\hat{\phi}]^{-1}$ is uniquely determined by $\set{\Phi_{\alpha (i)} \circ e^{\sigma_i} \circ \Phi_{\alpha (i)}^{-1}|_{U_i^+}}_{i \in I}$. In Proposition \ref{prop: eos:imdiff}, a representative of $E^+_\#([\hat{\sigma}])$ for $[\hat{\sigma}] \in H^+_\#$ in $\Orb{\aA^+,\bB^+}$ has been explicitly computed. Its lifts were given by $\set{\exp_{\alpha (i)}^\#  \circ \sigma_{U_i^+}}_{i\in I}$. Since the lifts uniquely determine the diffeomorphism, equation \eqref{eq: conj:vf} implies $c_{[\hat{\phi}]}E ([\hat{\sigma}]) = E^+_\# t([\hat{\sigma}]) \in \Diffo{Q,\uU}$ for every $[\hat{\sigma}] \in A$. In particular, $c_{[\hat{\phi}]} E(A) \subseteq \Diffo{Q,\uU}$. The set $E(A)$ is an open connected identity-neighborhood, whence it generates the connected Lie group $\Diffo{Q,\uU}$ by \cite[Theorem 7.4]{topgroup1979}. Therefore $c_{[\hat{\phi}]} (\Diffo{Q,\uU}) = c_{[\hat{\phi}]} (\langle E(A)\rangle) \subseteq \Diffo{Q,\uU}$. 
 We deduce from $c_{[\hat{\phi}]}|_{E(A)} = E^+_\# \circ t|_{A}^{H_\#^+} \circ (E|^{E(A)}_{A})^{-1}$ that the group automorphism $c_{[\hat{\phi}]}$ of $\Diffo{Q,\uU}$ is smooth on the open identity neighborhood $E(A)$, hence smooth by \cite[Ch.\ III, \S 1, No.\ 2, Proposition 4]{bourbaki1989}.
\end{proof}

The preceding proposition shows that for each $[\hat{\phi}]$, the conjugation map $c_{[\hat{\phi}]}$ is smooth and maps $\Diffo{Q,\uU}$ to itself. All requirements of Proposition \ref{prop: Lgp:locd} (b) have been checked. Applying this construction principle, we obtain a unique Lie group structure on $\Difforb{Q,\uU}$, turning $\Diffo{Q,\uU}$ into an open submanifold of $\Difforb{Q,\uU}$. Summarizing the results, we obtain:

\begin{thm}\label{thm: diff:lgp} \no{thm: diff:lgp}
 The group $\Difforb{Q,\uU}$ can be made into a Lie group in a unique way such that the following condition is satisfied:\\
 For some Riemannian orbifold metric $\rho$ on $(Q,\uU)$, let $[\expo]$ be the Riemannian orbifold exponential map with domain $\Omega$. There exists an open zero-neighborhood $\hH_\rho$ in $\Osc{Q}$ such that $$[\hat{\sigma}] \mapsto [\expo]\circ [\hat{\sigma}]|^\Omega$$ is a well-defined $C^\infty$-diffeomorphism of $\hH_\rho$ onto an open submanifold of $\Difforb{Q,\uU}$.\\
 The condition is then satisfied for every Riemannian orbifold metric on $(Q,\uU)$. The identity component of $\Difforb{Q,\uU}$ is the Lie group $\Diffo{Q,\uU}$ constructed in Section \ref{sect: lgp:ident}.
\end{thm}

\begin{cor}
 If $(Q,\uU)$ is a compact orbifold, then the Lie group $\Difforb{Q,\uU}$ is a Fr\'{e}chet-Lie group.
\end{cor}
\begin{proof}
 If $Q$ is compact, then $\Osc{Q} = \Os{Q}$ is a Fr\'{e}chet space, by Corollary \ref{cor: ost:prop}.
\end{proof}

We now consider subgroups of $\Difforb{Q,\uU}$ which turn out to be Lie subgroups of $\Difforb{Q,\uU}$. 
\begin{defn}
 Let $K \subseteq Q$ be a compact subset and denote for an orbifold map $[\hat{f}]$ its underlying map by $f$. Define the \ind{set!of orbifold diffeomorphisms supported in $K$}{set of all orbifold diffeomorphisms whose support is contained in $K$}: \glsadd{DiffK}
  \begin{displaymath}
   \DiffK{Q,\uU} \coloneq \setm{[\hat{f}] \in \Difforb{Q,\uU}}{f|_{Q\setminus K} \equiv \id_{Q\setminus K}}.
  \end{displaymath}
 We also say that the elements of $\DiffK{Q,\uU}$ coincide with the identity morphism of $Q$ off $K$.\\ Furthermore, we define the subset $\Diffc{Q,\uU} \subseteq \Difforb{Q,\uU}$\glsadd{Diffc} of all orbifold diffeomorphisms, whose underlying map coincides with $\id_Q$ outside some compact set in $Q$. Observe that the sets $\DiffK{Q,\uU}$ and $\Diffc{Q,\uU}$ are subgroups of $\Difforb{Q,\uU}$.
\end{defn}

\begin{rem}\label{rem: diff:ci} \no{rem: diff:ci}
 Notice that, by construction, $\Diffc{Q,\uU}$ contains $\Diffo{Q,\uU}$. The normal subgroup $\Diffc{Q,\uU}$ therefore is an open subgroup of $\Difforb{Q,\uU}$ by \cite[Theorem 5.5]{topgroup1979}. Hence it becomes a normal open Lie subgroup of $\Difforb{Q,\uU}$.  
\end{rem}

\begin{prop}\label{prop: DiffK:Lie}\no{prop: DiffK:Lie}
 Each compact subset $K$ of $Q$ is contained in a compact set $L$ such that the group $\Difforb{Q,\uU}_L$ is a closed Lie subgroup of $\Difforb{Q,\uU}$ modeled on $\Os{Q}_L$. 
\end{prop}

\begin{proof}
 We shall again use the notation of Section \ref{sect: lgp:ident}. The atlas $\aA$ is locally finite and the image of each chart in $\aA$ is relatively compact. Thus there are only finitely many charts $(U_i,G_i,\psi_i)$ in $\aA$ with $\psi_i (U_i) \cap K \neq \emptyset$. Let $I_K$ be the set indexing this family and consider the closed set $L \coloneq Q \setminus \left( \bigcup_{i \in I \setminus I_K} \psi_i (U_i)\right)$.
 By construction, $K \subseteq L \subseteq \overline{\bigcup_{i \in I_K} \psi_i (U_i)}$ holds, whence $L$ is a compact set. We claim that $\Difforb{Q,\uU}_L$ is a closed Lie subgroup of $\Difforb{Q,\uU}$ modeled on $\Os{Q}_L$.\\
 Choose for each $i \in I \setminus I_K$ a non singular point $x_i \in U_i$. By \cite[Theorem 1.9.5]{klingenberg1995}, we may choose $\ve_i >0$ with $\exp_{W_{\alpha (i)}} (B_{\rho_{\alpha (i)}} (0_{x_i}, \ve_i)) \cap H_{\alpha (i)} . x_i = \set{x_i}$. By definition of the topology on $\vect{W_{\alpha (i)}}$, there is an open neighborhood $\rR_i \subseteq \vect{W_{\alpha (i)}}$ of the zero-section such that $\sigma \in \rR_i$ implies $\sigma (x_i) \in B_{\rho_{\alpha (i)}} (0_{x_i}, \ve_i)$. Define the open neighborhood of the zero-orbisection 
  \begin{displaymath}
   \rR \coloneq \Lambda_\cC^{-1} \left(\bigoplus_{i \in I\setminus I_K} \rR_i \oplus \bigoplus_{j \in I_K} \vect{W_{\alpha (j)}}\right) \subseteq \Osc{Q}.
  \end{displaymath}
 Let $[\hat{\sigma}]$ be an element of $\hH_1 \cap \rR$, where $\hH_1$ is the open zero-neighborhood defined in Proposition \ref{prop: lgp:uncp}. Denote by $\set{\sigma_i}_{i \in I}$ the family of canonical lifts of $[\hat{\sigma}]$ with respect to $\aA$. Recall that $E([\hat{\sigma}])$ is a diffeomorphism, whose local lift with respect to $(U_i,G_i,\psi_i), i \in I\setminus I_K$ is the map $e^{\sigma_i} = \exp_{W_{\alpha (i)}}|_{\hat{O}_{\alpha (i)}} \circ \sigma_i$. Furthermore, $\exp_{W_{\alpha (i)}}|_{\hat{O}_{\alpha (i)} \cap T_x W_{\alpha (i)}}$ is a diffeomorphism for each $x\in U_i$, which maps $0_x$ to $x$. Since the canonical lift with respect to $(U_i,G_i,\psi_i)$ of the zero-orbisection is the zero-section, we deduce that $E(\hH_1 \cap \rR \cap \Os{Q}_L) \subseteq \Difforb{Q, \uU}_L$ holds.\\
 On the other hand, consider $[\hat{\sigma}] \in \hH_1 \cap \rR$ with $E([\hat{\sigma}]) \in \Difforb{Q,\uU}_L$. The underlying map of $E([\hat{\sigma}])$ coincides with $\id_Q$ on $Q \setminus L$. By construction, $\psi_i (U_i) \cap L \neq \emptyset$ holds for each $i \in I\setminus I_K$. Hence $\varphi_{\alpha (i)} \circ e^{\sigma_i} = \id_Q \circ \psi_i = \psi_i$. We deduce that $e^{\sigma_i} \colon U_i \rightarrow W_{\alpha (i)}$ must be an embedding of orbifold charts. Since the canonical inclusion $U_i \rightarrow W_{\alpha (i)}$ is an embedding of orbifold charts by Construction \ref{con: sp:ofdat} I.(d), Proposition \ref{prop: ch:prop} (d) yields $e^{\sigma_i} = h|_{U_i}$ for some $h \in H_{\alpha (i)}$. Specializing to the non-singular point $x_i\in U_i$, this yields $e^{\sigma_i} (x_i) = h(x_i) \in H_{\alpha (i)}.x_i$. Since $[\hat{\sigma}]$ is contained in $\rR$, $\sigma_i \in \rR_i$ and thus $e^{\sigma_i} (x_i) \cap H_{\alpha (i)} . x_i =\set{x_i}$. We obtain $h(x_i) = x_i$ and since $x_i$ is non-singular, $h= \id_{W_{\alpha (i)}}$ follows. Thus $e^{\sigma_i} = \id_{W_{\alpha (i)}}|_{U_i}$ and we deduce that $\sigma_i$ must be the zero-section in $\vect{U_i}$. Repeat the argument for each $i \in I \setminus I_K$. As $Q \setminus L = \bigcup_{i \in I\setminus I_K} \psi_i (U_i)$ holds by construction, $[\hat{\sigma}]$ is an element of $\Os{Q}_L$. Summarizing the preceding results, we obtain: 
  \begin{equation}\label{eq: osk}
   E(\hH_1 \cap \rR) \cap \Difforb{Q,\uU}_L = E (\hH_1 \cap \rR \cap \Os{Q}_L).
  \end{equation}
 Since $\pP = E(\hH_1)$ generates $\Diffo{Q,\uU}$, we deduce that $\Difforb{Q,\uU}_L$ is a Lie subgroup of $\Difforb{Q,\uU}$ modeled on $\Os{Q}_L$. The space $\Os{Q}_L$ is a closed vector subspace of $\Osc{Q}$ by Lemma \ref{lem: OsK}. Hence the identity \eqref{eq: osk} implies that $\Difforb{Q,\uU}_L$ is locally closed in the topological group $\Difforb{Q,\uU}$ and thus $\Difforb{Q,\uU}_L$ is a closed subgroup by \cite[Ch.\ III, \S 2, No.\ 1, Proposition 4]{bourbaki1971}.    
\end{proof}

For a trivial orbifold (i.e.\ a manifold) one need not refine the zero-neighborhood, i.e.\ we can always choose $K=L$ in Proposition \ref{prop: DiffK:Lie} for a trivial orbifold.  

\begin{rem}\label{rem: bor} \no{rem: bor}
  In \cite{bb2008} and the follow-up \cite{bb2009}, the diffeomorphism group of a compact orbifold was turned into a Fr\'{e}chet-Lie group in the sense of convenient differential calculus. The article \cite{bb2008} contains several errors, making it unclear whether the methods outlined in \cite{bb2008,bb2009} turn the orbifold diffeomorphism group into a convenient Lie group. Let us point out two serious problems: 
  \begin{compactitem}
   \item Lemma 23 in \cite{bb2008} states that the local lifts of an orbifold map are independent of local charts once the lifts are chosen. In particular, it is claimed that there is a unique extension of a lift defined on an open subset of a chart.  However, the lemma is false, as there may be several extensions to a lift. A counter-example can be obtained as follows: Let $\RR / \langle \gamma \rangle$ be the orbifold induced by the action of the reflection $\gamma$ at the origin. Consider a smooth map $f \colon ]-1,3[ \rightarrow \RR$ with $f(t) \neq 0$ if and only if $t \in ]0,1[ \cup ]2,3[$. If $q \colon \RR \rightarrow \RR/G$ is the global chart for this orbifold, $q\circ f$ is a continuous map, which induces a morphism of orbifolds in the sense of \cite{bb2008}. In fact, we may choose for example $f|_{]-1,1.5[}$ as a smooth lift at $0$. Clearly there are several possibilities to extend this lift smoothly to the pair of 
charts $]-1,3[, \RR$ thus contradicting the lemma.   
   \item In Definition 31 of \cite{bb2008}, the space of $C^r$-orbifold morphisms $C^r_{\text{Orb}} (\oO_1,\oO_2)$ is endowed with a topology following the classical construction of the $C^\infty$-Whitney topologies (cf.\ \cite{hirsch1976}). Unfortunately, the definition of the basic neighborhoods depends on a choice of lifts in the charts used. In \cite[Definition 20]{bb2008}, orbifold maps are equivalence classes of underlying maps with an assortment of local lifts. Two maps are equivalent if for each point $x \in Q$ their local lifts $f_x$ and $g_x$ at that point coincide as germs. Hence, if $x$ is a singular point, there are different lifts whose germ coincides in $x$ (simply compose with an element in the isotropy group of $x$). \\
   By definition of open sets in the topology, one has to compare local lifts at points on relatively compact open subsets of their domains. Thus at singular points there are maps with equivalent lifts which do not coincide outside of a small neighborhood. Hence, the distance depends on the choice of lifts, i.e.\ there are maps in a given open set, such that equivalent maps are not in the set. The topology and the distance defined in \cite[Definition 36]{bb2008} therefore depend on the representative of the equivalence class. This is problematic as the topology is defined on equivalence classes of orbifold maps, whence it should be independent of representatives.\\
   Unfortunately, these tools are needed in \cite{bb2008,bb2009} to turn the diffeomorphisms into a Lie group.   
  \end{compactitem}
\end{rem}
\subsection{The Lie algebra \texorpdfstring{of $\Difforb{Q,\uU}$}{}}\label{sect: LALG}

In this section, the Lie algebra $L(G)$ of the group $G \coloneq \Difforb{Q,\uU}$ constructed in Section \ref{sect: lgp:glob} will be determined. We stick to the notation introduced in Sections \ref{sect: lgp:ident} and \ref{sect: lgp:glob}. By construction, the map $E \colon \Osc{Q} \supseteq \hH_1 \rightarrow \pP \subseteq G, [\hat{\sigma}] \mapsto [\expo] \circ  [\hat{\sigma}]|^\Omega$ is a diffeomorphism of the open zero-neighborhood $\hH_1$ to an open identity-neighborhood $\pP$ in $G$. Furthermore, $E$ maps $\zs$ to $\ido{}$ by Proposition \ref{prop: diff:chart}. Use the natural isomorphism $T_{\zs} E$ to identify $T_{\ido{}} G$ with $\Osc{Q} \cong T_{\zs} \Osc{Q}$.\\
We modify the classical argument to compute the Lie algebra of the diffeomorphism group of a compact manifold via the adjoint action by Milnor (see \cite[pp.\ 1035-1036]{milnor1983}). To compute the Lie bracket, we have to understand the adjoint action of $T_{\ido{}} G$ on itself. Using the chart $E$, the product on $G$ pulls back to a smooth product operation 
    \begin{displaymath}
     [\hat{\sigma}] \ast [\hat{\tau}] \coloneq E^{-1} (E([\hat{\sigma}]) \circ E([\hat{\tau}]))
    \end{displaymath}
on the zero-neighborhood $\setm{([\hat{\sigma}],[\hat{\tau}])}{E([\hat{\sigma}])\circ E([\hat{\tau}]) \in \im E} \subseteq \hH_1 \times \hH_1 \subseteq \Osc{Q} \times \Osc{Q}$. By construction, $[\hat{\sigma}] \ast \zs = [\hat{\sigma}] = \zs \ast [\hat{\sigma}]$ holds. Hence the constant term of the Taylor series of $\ast$ in $(\zs ,\zs)$ (cf.\ \cite[Proposition 1.17]{hg2002a}) vanishes. Following \cite[Example II.1.8]{neeb2006}, the Taylor series is given as 
  \begin{displaymath}
   [\hat{\sigma}] \ast [\hat{\tau}] = ([\hat{\sigma}] + [\hat{\tau}]) + b([\hat{\sigma}], [\hat{\tau}]) + \cdots.
  \end{displaymath}
Here $b ([\hat{\sigma}], [\hat{\tau}]) = \left.\frac{\partial^2}{\partial s \partial t}\right|_{t,s=0} \hspace{-0.25cm} (t[\hat{\sigma}] \ast s [\hat{\tau}])$ is a continuous $\Osc{Q}$-valued bilinear map and the dots stand for terms of higher degree (cf.\ \cite{gn2007}). With arguments as in \cite[p.\ 1036]{milnor1983}, the adjoint action of $T_{\ido{}} G$ on itself is given by 
  \begin{displaymath}
   \text{ad} ([\hat{\sigma}]) [\hat{\tau}] = b([\hat{\sigma}],[\hat{\tau}]) - b([\hat{\tau}],[\hat{\sigma}]).
  \end{displaymath}
In other words, the skew-symmetric part of the bilinear map $b$ defines the adjoint action.\\ 
By \cite[Assertion 5.5]{milnor1983} (or \cite[Example II.3.9]{neeb2006}), the Lie algebra $L(G)$ of $G$ may be identified with $T_{\ido{}} G$ such that the Lie bracket coincides with the adjoint action: $\LB[x,y] = \text{ad}(x)y$. To compute the Lie bracket $\LB{}$, it is sufficient to compute the second derivative of the local product operation in $\Osc{Q}$. Consider the atlas $\aA$ as in Construction \ref{con: sp:ofdat} together with the linear topological embedding with closed image $\Lambda_\aA \colon \Osc{Q} \rightarrow \bigoplus_{i \in I} \vect{U_i}, [\hat{\sigma}] \mapsto (\sigma_i)_{i \in I}$. For fixed $[\hat{\sigma}], [\hat{\tau}] \in \Osc{Q}$, the map $(t,s) \mapsto t [\hat{\sigma}] \ast s [\hat{\tau}]$ factors through a finite subproduct of the direct sum. Hence the derivative of $s[\hat{\sigma}] \ast t [\hat{\tau}]$ may be computed from the derivatives of the canonical lifts $(t[\hat{\sigma}] \ast s[\hat{\tau}])_{i}$.\\
Recall from Lemma \ref{lem: comp:cor} that for each pair $[\hat{\sigma}], [\hat{\tau}] \in \hH_1$, there is an orbisection $[\widehat{\sigma \diamond \tau}] \in \Osc{Q}$ such that $E([\widehat{\sigma \diamond \tau}]) = E([\hat{\sigma}]) \circ E([\hat{\tau}])$ Returning for a moment to $E$ as a map on $\mM$ as in Proposition \ref{prop: diff:chart}. The mapping $E$ is bijective, whence we deduce for $i \in I$ the identity  
  \begin{displaymath}
   (t[\hat{\sigma}] \ast s[\hat{\tau}])_{i} = (t\sigma_{\alpha (i)} \diamond_i s\tau_{\alpha (i)})|_{U_i}.
  \end{displaymath}
 For the rest of the proof, fix $i \in I$ and compute $\left.\frac{\partial^2}{\partial s \partial t}\right|_{t,s=0}\hspace{-0.25cm} (t\sigma_{\alpha (i)} \diamond_i s\tau_{\alpha (i)})|_{U_i}$. By construction, the vector field $t\sigma_{\alpha (i)} \diamond_i s\tau_{\alpha (i)}$ is defined on $\Omega_{\frac54 , K_{5,i}}$. As the restriction map $\res_{U_i}^{\Omega_{\frac54 , K_{5,i}}}$ is continuous linear by \cite[Lemma F.15 (a)]{hg2004}, it commutes with the differential, i.e. $$\res_{U_i}^{\Omega_{\frac54 , K_{5,i}}} \left. \frac{\partial^2}{\partial s \partial t}\right|_{t,s=0} \hspace{-0.25cm} t\sigma_{\alpha (i)} \diamond_i s\tau_{\alpha (i)} = \left.\frac{\partial^2}{\partial s \partial t}\right|_{t,s=0} \hspace{-0.25cm} (t\sigma_{\alpha (i)} \diamond_i s\tau_{\alpha (i)})|_{U_i}.$$ Thus it suffices to compute the derivative in $\vect{\Omega_{\frac54,K_{5,i}}}$. \\ 
 The set $\setm{(V_{\frac54 ,\alpha (i)}^{n}, \kappa_n^{\alpha (i)})}{(V_{5,\alpha (i)}^n , \kappa_n^{\alpha (i)})) \in \fF_5 (K_{5,i}}$ is finite and covers $\Omega_{\frac54,K_{5,i}}$.  Hence the topology on the space $\vect{\Omega_{\frac54 , K_{5,i}}}$ is induced by the linear embedding with closed image 
    \begin{displaymath}
     \Gamma \colon \vect{\Omega_{\frac54 ,K_{5,i}}} \rightarrow \prod_{(V_{5,\alpha (i)}^n,\kappa_n^{\alpha (i)}) \in \fF_5 (K_{5,i})} C^\infty (V_{\frac54,\alpha (i)}^n , \RR^d) , X \mapsto (\text{pr}_2 T\kappa_n^{\alpha (i)} X|_{V_{\frac54,\alpha (i)}^n})_{\fF_{5,K_{5,i}}}.
    \end{displaymath}
 Here $\text{pr}_2$ is the linear projection onto the second component of $B_{\frac54} (0) \times \RR^d$. Since $(\kappa_n^{\alpha (i)})^{-1}|_{B_{\frac54} (0)}$ is a diffeomorphism onto $V_{\frac54 , \alpha (i)}^n$, the mapping 
  \begin{displaymath}
   C^\infty ((\kappa_n^{\alpha (i)})^{-1}|_{B_{\frac54} (0)}, \RR^d) \colon C^\infty (V_{\frac54 , \alpha (i)}^n ,\RR^d) \rightarrow C^\infty (B_{\frac54} (0),\RR^d), X \mapsto X \circ (\kappa_n^{\alpha (i)})^{-1}|_{B_{\frac54} (0)}
  \end{displaymath}
 is an isomorphism of topological vector spaces by \cite[Lemma A.1]{hg2004}. We derive an embedding of topological vector spaces with closed image $(C^\infty ((\kappa_n^{\alpha (i)})^{-1}|_{B_{\frac54} (0)} ,\RR^d))_{\fF_5 (K_{5,i}} \circ \Gamma $. Using this map, the derivative may be computed locally in $A \coloneq \prod_{\fF_5 (K_{5,i})} C^\infty (B_{\frac{5}{4} } (0) ,\RR^d)$. For $X \in \vect{W_{\alpha (i)}}$, define $X_{[n]} \coloneq \text{pr}_2 T\kappa_n^{\alpha (i)} X (\kappa_n^{\alpha (i)})^{-1}|_{B_{\frac54} (0)} \in C^\infty (B_{\frac54} (0), \RR^d)$. The map $\text{pr}_2$ is linear and each $T\kappa_n^{\alpha (i)}$ is linear in the vector space component. Hence the definition of the vector space operations of $\vect{W_{\alpha (i)}}$ shows that the identity $(tX)_{[n]} = tX_{[n]}$ holds for each $t \in \RR$ and $X \in \vect{W_{\alpha (i)}}$.\\
 To compute the derivative of $(t \sigma_{\alpha (i)} \diamond_i s \tau_{\alpha (i)})$ in $A$, more information on $(t \sigma_{\alpha (i)} \diamond_i s \tau_{\alpha (i)})_{[n]}$ is needed. Fortunately, by Construction \ref{con: H} a local formula is available. To write it down explicitly, we need to recall notation and facts from the construction: \\ 
 For each chart $(V_{5,\alpha (i)}^n , \kappa_{n}^{\alpha (i)})$, let $\exp_n$ be the Riemannian exponential map on $B_5 (0)$ associated to the pullback metric with respect to $\kappa_n^{\alpha (i)}$ and the member of the orbifold metric $\rho_{\alpha (i)}$ on $W_{\alpha (i)}$. Recall from Construction \ref{con: comp:loc} that for $x \in V_{\frac54,\alpha (i)}^n$, there is an open set $N_x \subseteq T_x W_{\alpha (i)}$ such that $T\kappa_n^{\alpha (i)} (N_x) \subseteq \dom \exp_n$ holds and $\exp_{n}$ restricts to an \'{e}tale embedding on this set (cf.\ Lemma \ref{lem: eos:cov}). By Construction \ref{con: H}, for $(t\sigma_{\alpha (i)} \diamond_i s\tau_{\alpha (i)})|_{\Omega_{\frac54 , K_{5,i}}}$ and each chart $V_{\frac54 , \alpha (i)}^n$ the local identity \eqref{eq: glo:loc} holds. We want to keep track of the local chart $(V_{5,k}^n , \kappa_n^k)$ in which we construct a new vector field via the operation $\diamond$ as in Construction \ref{con: comp:loc}. Hence we write $\diamond_{[n]}$ for $\diamond$ in the chart $(V_{5,k}^n , \kappa_n^k)$. Using the notation introduced, the identity \eqref{eq: glo:loc} yields the following formula for $x \in B_{\frac54} (0)$ 
  \begin{align*}
   &T\kappa_n^{\alpha (i)} t\sigma_{\alpha (i)} \diamond_i s\tau_{\alpha (i)} (\kappa_n^{\alpha (i)})^{-1} (x) = (x, (t\sigma_{\alpha (i)})_{[n]} \diamond_{[n]} (s\tau_{\alpha (i)})_{[n]} (x))\\ 
 = &(x, (\exp_n|_{T\kappa_n^{\alpha (i)} (N_x)})^{-1} \exp_n (\exp_n (x,(s\tau_{\alpha (i)})_{[n]} (x)) , (t\sigma_{\alpha (i)})_{[n]} (\exp_n (s\tau_{\alpha (i)})_{[n]} (x))))\\ 
 = &(x, (\exp_n|_{T\kappa_n^{\alpha (i)} (N_x)})^{-1} \exp_n (\exp_n (x,s(\tau_{\alpha (i)})_{[n]} (x)) , t(\sigma_{\alpha (i)})_{[n]} (\exp_n s(\tau_{\alpha (i)})_{[n]} (x)))). 
  \end{align*}
 Apply $\text{pr}_2$ to the formula above to obtain the desired identity for $(t \sigma_{\alpha (i)} \diamond_i s \tau_{\alpha (i)})_{[n]}$. To simplify the notation, we abbreviate $X \coloneq (\sigma_{\alpha (i)})_{[n]}$ and $Y \coloneq (\tau_{\alpha (i)})_{[n]}$. Recall the following properties of $\exp_n$ (cf.\ \cite[Theorem 1.6.12]{klingenberg1995}): $\exp_n (x,0) = x, d_2\exp_n (x,0) = \id_{\RR^d}$ for all $x \in B_{\frac{5}{4}} (0)$. Define $M_x \coloneq T\kappa_n^{\alpha (i)} (N_{(\kappa_n^{\alpha (i)})^{-1} (x)}) \subseteq T_x B_{\frac{5}{4}} (0)$ for $x \in B_{\frac{5}{4}} (0)$. Since $\exp_n$ is injective on $(x,0) \in T\kappa_n^{\alpha (i)} (N_x)$ with $\exp_n (x,0) = x$ and $d_2\exp_n (x,0) = \id_{\RR^d}$, we derive 
 \begin{displaymath}
  d(\exp_n|_{M_x}))^{-1} (x , \cdot) = \id_{\RR^d}.                                                                                                                                                                                                                                                                                                                                                                                                                                                                                                                                                                                                                                                                                                                                                                                                                                                                                                                                                                                                                                                                                                                                                                                                                         \end{displaymath}
 For $x \in B_{\frac54} (0)$, the facts collected above allow us to obtain 
    \begin{align}
     &\left. \frac{\partial^2}{\partial s \partial t}\right|_{t,s=0} \hspace{-0.25cm} (t\sigma_{\alpha (i)} \diamond_i s\tau_{\alpha (i)})_{[n]} (x) \notag \\
    =&\left. \frac{\partial^2}{\partial s \partial t}\right|_{t,s=0} (\exp_n|_{M_x})^{-1} \exp_n (\exp_n (x,sY (x)) , tX (\exp_n (x,sY (x)))) \notag\\
    =&\left.\frac{\partial}{\partial s}\right|_{s=0} d  (\exp_n|_{M_x})^{-1} \left.\frac{\partial}{\partial t}\right|_{t=0} \exp_n (\exp_n (x,sY (x)) , tX (\exp_n (x,sY (x)))) \notag\\
    =&\left.\frac{\partial}{\partial s}\right|_{s=0} d  (\exp_n|_{M_x})^{-1} (\exp_n (x,sY (x)) , X (\exp_n (x, sY (x)))). \label{eq: local:vf}
    \end{align}
  The map $d(\exp_n|_{M_x})^{-1}$ is linear in the second argument. Hence the rule on partial derivatives \eqref{eq: rule:pdiff} applied to \eqref{eq: local:vf} yields the following identity: 
  \begin{align*} 
   \left.\frac{\partial^2}{\partial s \partial t}\right|_{t,s=0}\hspace{-0.25cm} (t\sigma_{\alpha (i)} \diamond_i s\tau_{\alpha (i)})_{[n]} (x) =& d  (\exp_n|_{M_x})^{-1} \left(\exp_n (x,0),dX \left(\left.\frac{\partial}{\partial s}\right|_{s=0} \exp_n (x,sY(x))\right)\right)\\ 
       &+ d^{(2)} (\exp_n|_{M_x})^{-1} \left(\left.\frac{\partial}{\partial s}\right|_{s=0} \exp_n (x,sY(x)), X(\exp_n (x,0))\right)\\
    =& dX (x,Y(x)) + \underbrace{d^{(2)} (\exp_n|_{M_x})^{-1} (x,Y(x), X(x)))}_{S_{XY} \coloneq}.
  \end{align*}
 The derivative $d^{(2)} (\exp_n|_{T\kappa_n^{\alpha (i)} (N_x)})^{-1} (x,\cdot , \cdot)$ is a symmetric bilinear map by \cite[Proposition 1.13]{hg2002a}. Hence $S_{XY}$ is symmetric in $X$ and $Y$. An analogous computation yields: 
  \begin{displaymath}
  \left. \frac{\partial^2}{\partial s \partial t}\right|_{t,s=0}\hspace{-0.25cm} (t\sigma_{\alpha (i)} \diamond_i s\tau_{\alpha (i)})_{[n]} (x) = dY (x,X(x)) + S_{XY}.
  \end{displaymath}
 As $C^\infty (\kappa_n^{\alpha (i)}, \RR^d)$ is an isomorphism of topological vector spaces and evaluation at $x$ is continuous linear, $((\text{ad} ([\hat{\sigma}])[\hat{\tau}])_{\alpha (i)})_{[n]}$ is given by  
    \begin{align*}
     ((\text{ad} ([\hat{\sigma}])[\hat{\tau}])_{\alpha (i)})_{[n]} (x) &= \left.\frac{\partial^2}{\partial s \partial t}\right|_{t,s=0}\hspace{-0.25cm} (t\sigma_{\alpha (i)} \diamond_i s\tau_{\alpha (i)})_{[n]}(x) -\left. \frac{\partial^2}{\partial s \partial t}\right|_{t,s=0} \hspace{-0.25cm} (t\tau_{\alpha (i)} \diamond_i s\sigma_{\alpha (i)})_{[n]} (x)\\ 
      &= dX (Y(x)) - dY (X(x)) = (d  (\sigma_{\alpha (i)})_{[n]} (\tau_{\alpha (i)})_{[n]} - d(\tau_{\alpha (i)})_{[n]} (\sigma_{\alpha (i)})_{[n]})(x).
    \end{align*}
 Recall from \cite[Definition I.3.6]{neeb2006} that the Lie bracket of vector fields $V,W$ in $\vect{\Omega_{\frac54 , K_5,i}}$ is the unique vector field $\LB[V,W]_{i}$ such that for each chart $(V_{\frac54 , \alpha (i)}^n, \kappa_{n}^{\alpha (i)}) \in \fF_5 (K_{5,i})$ the identity 
    \begin{displaymath}
     (\LB[V,W]_{i})_{[n]} = dW_{[n]} V_{[n]} - dV_{[n]} W_{[n]}
    \end{displaymath}
 is satisfied. By the above computation, the negative of the Lie bracket of the vector fields $\sigma_{\alpha (i)}$ and $\tau_{\alpha (i)}$ coincides with $(\text{ad} ([\hat{\sigma}])[\hat{\tau}])_{\alpha (i)}$ on $\Omega_{\frac54,K_{5,i}}$. Since $U_i \subseteq \Omega_{\frac54 , K_{5,i}}$  holds, the canonical lift $(\text{ad} ([\hat{\sigma}])[\hat{\tau}])_{i}$ on $U_i$ coincides with the negative of the Lie bracket of the canonical lifts of $\sigma_i$ and $\tau_i$. By abuse of notation, let $\LB[\sigma_i ,\tau_i]$ be the Lie bracket of the lifts in $\vect{U_i}$. The families $\set{\sigma_{i}}_{i \in I}$ and $\set{\tau_{i}}_{i \in I}$ are families of canonical lifts of the orbisections $[\hat{\sigma}]$ and $[\hat{\tau}]$ with respect to the atlas $\aA$. Hence each pair of lifts $\sigma_i, \sigma_j$ (respectively $\tau_i$, $\tau_j$) for $i,j \in I$ is $\phi$-related for $\phi \in \Ch_{U_i,U_j}$ (i.e.\ \eqref{eq: os:chch} holds). By \cite[Assertion 4.6]{milnor1983}, $\LB[\sigma_i ,\tau_i]$ and $\LB[\sigma_j,\tau_j]$ are $\phi$-related for each $\phi \in \Ch_{U_i,U_j}$ and every pair $i,j \in I$. Hence $(\LB[\sigma_i,\tau_i])_{i \in I}$ is a family of canonical lifts for the compactly supported orbisection $\text{ad} ([\hat{\sigma}])[\hat{\tau}]$. The result of this section may now be summarized as follows:

 \begin{thm}[{Lie algebra of $\Difforb{Q,\uU}$}] \label{thm: LALG}
 Identify $T_{\ido{}} \Difforb{Q,\uU}$ via $T_{\zs} E$ with the space $\Osc{Q}$ and the Lie algebra of $\Difforb{Q,\uU}$ with $(\Osc{Q}, \LB[\cdot,\cdot])$. The Lie bracket $\LB[\cdot,\cdot]$ is defined as follows:\\ For arbitrary $[\hat{\sigma}], [\hat{\tau}] \in \Osc{Q}$, their Lie bracket $\LB[[\hat{\sigma}], [\hat{\tau}]]$ is the unique compactly supported orbisection whose canonical lift on an orbifold chart $(U,G,\varphi)$ is the negative of the Lie bracket in $\vect{U}$ of their canonical lifts $\sigma_U$ and $\tau_U$. 
 \end{thm}

If the orbifold is trivial (i.e.\ a manifold), Theorem \ref{thm: LALG} specializes to the well known description of the Lie algebra for the diffeomorphism group of a manifold (cf.\ \cite[Example II.3.14]{neeb2006}).
\subsection{Regularity properties of \texorpdfstring{$\Difforb{Q,\uU}$}{the Diffeomorphismgroup}}\label{sect: regularity}

In this section, we prove that $\Difforb{Q,\uU}$ is a regular Lie group in the sense of Milnor (cf.\ \cite[Definition 7.6]{milnor1983}). Unless stated otherwise the notation from Section \ref{sect: lgp:ident} and Section \ref{sect: lgp:glob} will be used. Another prerequisite is the definition of $C^k$-regularity as outlined in Appendix \ref{app: reg}. The philosophy in the proof of the Lie group properties for $\Difforb{Q,\uU}$ was to compute the relevant data locally in orbifold charts. Hence we investigate the situation on orbifold charts, where we study the flows of vector fields and their differentiability properties. Several facts from the calculus of $C^{r,s}$-mappings (see Definition \ref{defn: crs}, cf.\ \cite{alas2012}) are needed. We study the following differential equation:

 \begin{setup}\label{setup: regdiff} \no{setup: regdiff}
  Define $f \colon [0,1] \times B_5 (0) \times C^r ([0,1] , C^\infty (B_5 (0), \RR^d))  \rightarrow \RR^d$ via $f(t,x,\gamma ) \coloneq \gamma^\wedge (t,x) \coloneq \gamma (t)(x)$ for $r \in \NN_0 \cup \set{\infty}$. Consider the evaluation maps $\ve \colon C^\infty (B_5 (0),\RR^d) \times B_5 (0) \rightarrow \RR^d, \ve (\sigma ,x) \coloneq  \sigma (x)$ and $\ve_1 \colon C^r ([0,1] , C^\infty (B_5 (0), \RR^d)) \times [0,1] \rightarrow C^\infty (B_5 (0), \RR^d), (\gamma,t) \mapsto \gamma (t)$. By \cite[Proposition 3.20]{alas2012}, $\ve$ is smooth and $\ve_1$ is of class $C^{\infty, r}$. We may rewrite the map $f$ as $f(t,x,\gamma) = \ve (\ve_1 (\gamma , t) , x)$. Hence the chain rule \cite[Lemma 3.17]{alas2012} implies that $f$ is of class $C^{r,\infty}$ with respect to the product $[0,1] \times \left(B_5 (0) \times C^r([0,1],C^\infty (B_5 (0), \RR^d))\right)$. Thus the initial value problem 
    \begin{equation}\label{eq:diffeq}
     \begin{cases}
      x'(t)  &= f(t,x(t),\gamma) =\gamma^\wedge (t,x (t)),\\
      x(t_0) &= x_0, \quad \quad x_0 \in B_5 (0)
     \end{cases}
    \end{equation}
  admits a unique maximal solution $\varphi_{t_0,x_0,\gamma}$ by \cite[Theorem 5.6]{alas2012}. Fixing $t_0 = 0$, the flow\glsadd{Fl} of \eqref{eq:diffeq}, 
  \begin{displaymath}
   \Fl^f_0 \coloneq \Fl^f (0,\cdot) \colon [0,1] \times \left(B_5 (0) \times C^r([0,1],C^\infty (B_5 (0), \RR^d))\right) \supseteq \Omega_0 \rightarrow \RR^d, (t,(x_0,\gamma)) \mapsto \varphi_{0, x_0,\gamma} (t),
  \end{displaymath}
 is of class $C^{r+1,\infty}$ on the open subset $\Omega_{0}$ by \cite[Proposition 5.9]{alas2012}. 
 \end{setup}
  
 \begin{lem}\label{lem: loc:eb} \no{lem: loc:eb}
  Let $r \in \NN_0 \cup \set{\infty}$, $\gamma \in C^r ([0,1], C^\infty (B_5 (0), \RR^d))$ and consider $f$ as in \ref{setup: regdiff}.
    \begin{compactenum}
     \item  If $\gamma$ satisfies $\norm{\gamma(t)}_{\overline{ B_4 (0)},0} \leq 1$ for all $t \in [0,1]$, then the map $\Fl^f_0 ( \cdot , \gamma)$, is defined on $[0,1] \times B_3 (0)$ and $\Fl^f_0 ([0,1] \times B_3 (0) \times \set{\gamma}) \subseteq B_4 (0)$.
     \item Consider $\zeta > 0$ and a compact subset $K \subseteq B_3 (0)$. There exists $0 < \tau \leq 1$ such that for all $\gamma \in C^r ([0,1], C^\infty (B_5 (0),\RR^d))$ with $\sup_{t \in [0,1]} \norm{\gamma (t)}_{\overline{B_4 (0)},1} < \tau$ (cf.\ Definition \ref{defn: norm}), we have $\lVert\Fl_0^f (t, \cdot, \gamma)- \id_{B_3 (0)}\rVert_{K,1} < \zeta$ for all $t \in [0,1]$. 
     \item For $\tau$ as in {\rm (b)} and $B_{\tau} (0) \coloneq \setm{f \in  C^\infty (B_5 (0),\RR^d)}{\norm{f}_{\overline{B_4 (0)},0} < \tau}$, we obtain a smooth map 
      \begin{displaymath}
       F \colon C^{r} ([0,1], B_\tau (0) ) \rightarrow C^{r+1} ([0,1], C^\infty (B_3 (0),\RR^d)), \gamma \mapsto \Fl^f_0 (\cdot,\gamma)|_{[0,1] \times B_3 (0)}.
      \end{displaymath}
    \end{compactenum}
 \end{lem}
  
 \begin{proof}
  \begin{compactenum}
   \item For $x_0 \in B_3 (0)$, the maximal solution to the initial value problem \eqref{eq:diffeq} is the mapping $\Fl_0^f (\cdot ,x_0, \gamma)$. We claim that it is defined at least on $[0,1]$. Restricting $\Fl_0^f$, we obtain the maximal solution to the initial value problem \eqref{eq:diffeq} whose image remains inside of $B_4 (0)$: Denote this solution by $u \colon [0,t_0[ \rightarrow B_4 (0)$. Then $u$ is of class $C^1$. If $t_0 < 1$ holds, we deduce  $\norm{u(t)} \leq \norm{u(0)} + \lVert\int_0^{t} \gamma^\wedge (s,u(s)) ds\rVert \leq \norm{x_0} +1 \equalscolon \rho < 4$ from the Fundamental Theorem of Calculus \cite[Theorem 1.5]{hg2002a}. Therefore $u|_{[0,t_0[}$ does not leave the compact subset $\overline{B_\rho (0)} \subseteq B_4 (0)$. Close to $t_0$, the right hand side of the differential equation \eqref{eq:diffeq} is defined on an open subset of a finite-dimensional Banach space, whence by \cite[Lemma 3.11]{hg2007}, $C^k$ maps coincide with the $k$-times continuously Fr\'{e}chet differentiable maps considered in \cite{langdgeo2001}. One may therefore apply \cite[IV.\ Theorem 2.3]{langdgeo2001}: The maximal solution must be defined on an interval strictly larger than $[0,t_0[$, thus contradicting the choice of $t_0$. We conclude that $\Fl^f_0 (\cdot , \gamma)$ maps $[0,1] \times B_3 (0)$ into $B_4 (0)$.
   \item Observe that $\Fl_0^f (\cdot ,\gamma)$ is of class $C^{r+1,\infty}$ by \ref{setup: regdiff}. By \cite[Lemma 3.15]{alas2012} $\Fl_0^f (\cdot ,\gamma)$ is a $C^1$-mapping, whence the derivatives required for $\norm{\cdot}_{K,1}$ exist. The mapping $h \colon [0,1] \times B_3 (0) \rightarrow \RR^d, h(t,x) \coloneq \gamma^\wedge (t,\Fl_0^f (t,x,\gamma))$ is of class $C^{r,\infty}$ by the chain rule \cite[Lemma 3.19]{alas2012}. Fix $x \in B_3 (0)$ and consider the map $g \colon [0,1] \rightarrow \BoundOp{\RR^d}, g(t) \coloneq d_2 \Fl_0^f (t,x,\gamma; \cdot)$. Schwarz' Theorem \cite[Proposition 3.6 and Remark 3.7]{alas2012} implies that $g$ is a solution to   
	\begin{equation}\label{eq:diffhigh}
	  \begin{cases}
	              y'(t) &= d_2 \gamma^\wedge (t,\Fl_0^f (t,x,\gamma); \cdot) \circ y(t)\\
                      y (0) &= \id_{\RR^d} .     
	                 \end{cases}
	\end{equation}
 The domain of $\gamma^{\wedge} (t,\cdot)$ is an open subset of $\RR^d$. Hence the derivative $d_2 \gamma^\wedge (t,x;\cdot)$ is determined by the Jacobian matrix. As all norms on $\RR^d$ are equivalent, there is a constant $C>0$, depending only on $d$ and the choice of norm such that $\opnorm{d_2 \gamma^\wedge (t,x;\cdot)} \leq C \sup_{|\alpha|=1} \norm{\partial^\alpha \gamma^\wedge (t,x)}$ with partial derivatives in the $x$-variable. Furthermore, $\Fl_0^f (\cdot ,\gamma)$ maps $[0,1] \times B_3(0)$ into $B_4 (0)$ by (a) and $\norm{\cdot}_{\overline{B_4 (0)},1}$ controls the partial derivatives. Hence the above estimate yields 
  \begin{displaymath}
   \sup_{t \in [0,1]}  \opnorm{d_2 \gamma^\wedge (t,\Fl_0^f (t,x,\gamma); \cdot)} \leq \sup_{t \in [0,1]} C \norm{\gamma (t)}_{\overline{B_4 (0)},1}.
  \end{displaymath}
 Vice versa, there is a constant $c>0$, depending only on the norm and $d$ such that
   \begin{displaymath}
      \sup_{t\in [0,1]} \sup_{|\alpha|=1} \norm{\partial^\alpha (\Fl_0^f (t,\cdot,\gamma) - \id_{\RR^d})(x)} \leq \sup_{t\in [0,1]}c \opnorm{g(t)-g(0)}.
   \end{displaymath}
 Let $\theta>0$ be an upper bound for $\sup_{t \in [0,1]} C\norm{\gamma (t)}_{\overline{B_4 (0)},1}$. The mapping $g$ is of class $C^1$, whence the Fundamental Theorem of Calculus \cite[Theorem 1.5]{hg2002a} yields: 
  \begin{align*}
   \opnorm{g(t) - \id_{\RR^d}} &= \opnorm{g(t) -g(0)} = \opnorm{\int_{0}^t d_2\gamma^\wedge (s, \Fl_0^f (s,x,\gamma);g(s))ds} \\
                     &\leq \int_0^t \theta \opnorm{g(s)}ds =\int_0^t \theta (\opnorm{g(0)} + (\opnorm{g(s)}-\opnorm{g(0)}))ds \\
                     &\leq \int_0^t \theta \opnorm{\id_{\RR^d}} ds + \int_0^t \theta \opnorm{g(s)-\id_{\RR^d}} ds = \theta t + \int_0^t \theta \opnorm{g(s)-\id_{\RR^d}}ds.
  \end{align*}
 Apply Gronwall`s inequality \cite[6.1 Gronwall's Lemma]{amann1990} to choose $1> \tfrac{\tau_1}{C} >0$ such that $\sup_{t \in [0,1]} \norm{\gamma (t)}_{\overline{B_4 (0)},1} < \tfrac{\tau_1}{C}$ implies 
  \begin{equation}\label{eq: deriv:est}
    \sup_{t\in [0,1]} \sup_{|\alpha|=1} \norm{\partial^\alpha (\Fl_0^f (t,\cdot,\gamma) - \id_{\RR^d})(x)} \leq \sup_{t\in [0,1]}c \opnorm{g(t)-g(0)} < \zeta.
  \end{equation}
 Observe that the estimate \eqref{eq: deriv:est} holds for each $x\in B_3 (0)$, as the constants did not depend on $x$. We have to obtain an estimate for $\Fl_0^f$: The Fundamental Theorem of Calculus \cite[Theorem 1.5]{hg2002a} with equation \eqref{eq:diffeq} yields for $x\in B_3 (0)$: 
    \begin{align*}
     \norm{\Fl_{0}^f (t,x,\gamma)-\id_{B_3 (0)} (x)} &= \norm{\Fl_{0}^f (t,x,\gamma)-\Fl_0^f (0,x,\gamma)} = \norm{\int_0^t \gamma^\wedge (s , \Fl_{0}^f (s,x,\gamma)) ds} .
    \end{align*}
   Require $\sup_{t \in [0,1]} \norm{\gamma (t)}_{\overline{B_4 (0)},0} < \zeta$ to obtain $\sup_{t\in [0,1]} \norm{\Fl_{0}^f (t,x,\gamma)-\id_{B_3 (0)} (x)} <  \zeta$.
  The estimates show that $\tau \coloneq \min \set{\zeta , \tfrac{\tau_1}{C},1}$ is a constant with the desired properties.  
  \item Let $r \in \NN_0 \cup \set{\infty}$, $X$ be a Fr\'{e}chet space and $U \subseteq \RR^d$ an open subset. By Remark \ref{rem: crfs:prop}, each of the topological spaces $[0,1], C^r ([0,1], X)$ and $C^r (U, X)$ is metrizable. The set $C^r ([0,1], B_\tau (0))$ is an open subset of the Fr\'{e}chet space $C^r([0,1], C^\infty (B_5 (0), \RR^d))$ (cf.\ \cite[Lemma 3.6]{hg2002}), hence metrizable. Therefore each finite Cartesian product of these spaces is a $k$-space by \cite[XI.\ 9.3]{dugun1966} and we may use the Exponential Law for $C^{r,s}$-maps (cf.\ \cite[Theorem 3.28 (e)]{alas2012}):\\
  Since $\Fl_0^f (\cdot ,\gamma)$ is of class $C^{r+1,\infty}$, we deduce that $F(\gamma)$ is in $C^{r+1}([0,1],C^\infty (B_3 (0),\RR^d)$. Hence $F$ makes sense and we claim that $F$ is smooth. By \ref{setup: regdiff}, $\Fl_0^f$ is of class $C^{r+1,\infty}$ on the product $[0,1] \times (B_3 (0) \times C^{r} ([0,1],B_\tau (0)))$. The Exponential Law implies that $$(\Fl_0^f)^\vee \colon [0,1] \rightarrow C^\infty (B_3(0) \times C^r ([0,1],B_\tau (0))), t \mapsto ((x,\gamma) \mapsto \Fl_0^f (t,x,\gamma))$$ is a $C^{r+1}$-map. Now $(\Fl_0^f)^\vee$ coincides with the map 
  \begin{displaymath}
  (\Fl_0^f)^\dagger \colon [0,1] \rightarrow C^\infty (C^{r} ([0,1],B_\tau(0)) \times B_3 (0), \RR^d), t \mapsto ((\gamma,x) \mapsto \Fl_0^f (t,x,\gamma)),                                                                               
  \end{displaymath}
  except for the inessential order of $x$ and $\gamma$. Combine the Exponential Law  with \cite[Lemma 3.22]{alas2012} to establish the isomorphism
    \begin{displaymath}
     \Phi \colon C^\infty (C^r([0,1],B_\tau (0)) , C^\infty (B_3 (0),\RR^d)) \rightarrow C^\infty (C^{r} ([0,1],B_\tau(0)) \times B_3 (0), \RR^d), f \mapsto f^\wedge.
    \end{displaymath}
  Then $(F^\wedge)^\dagger \coloneq (\Phi^{-1} ((\Fl_0^f)^\dagger)) \colon [0,1] \rightarrow C^\infty (C^r([0,1],B_\tau (0)) , C^\infty (B_3 (0),\RR^d))$ is a mapping of class $C^{r+1}$ by the Exponential Law. Evaluating $(F^\wedge)^\dagger$ at $(t,\gamma) \in [0,1] \times C^r([0,1] , B_\tau (0))$ the definition yields $(F^\wedge)^\dagger (t)(\gamma) = F^\wedge (\gamma, t)$. Hence by \cite[Corollary 3.8]{alas2012} and the Exponential Law, the map $F^\wedge \colon C^r([0,1], B_\tau (0)) \times [0,1] \rightarrow C^\infty (B_3 (0),\RR^d), (\gamma, t) \mapsto F(\gamma)(t)$ is a $C^{\infty, r+1}$-map. By \cite[Theorem 3.28 (e)]{alas2012}, this proves $F$ to be a smooth map.
  \end{compactenum}
 \end{proof}
\noindent
 To prove the ($C^0$-)regularity of $\Difforb{Q,\uU}$, we have to construct a smooth evolution map $C^0 ([0,1], \Osc{Q}) \rightarrow \Difforb{Q,\uU}$. We will assure the smoothness of all relevant maps via patched mapping arguments. These are prepared by the following preliminary lemma. 

 \begin{lem}\label{lem: loc:comp} \no{lem: loc:comp}
  Consider $r \in \NN_0 \cup \set{\infty}$ and define for $\gamma \in C^r([0,1],\vect{W_{\alpha (i)}})$ and $(V_{5,\alpha (i)}^n , \kappa^{\alpha (i)}_n) \in \fF_5 (K_{5,i})$ the $C^r$-curves $\gamma_{\kappa_n^{\alpha (i)}} \coloneq \theta_{\kappa_n^{\alpha (i)}} \circ \gamma$ (cf.\ Definition \ref{defn: top:vect}) and $\gamma_{[n]} \coloneq C^\infty ((\kappa_{n}^{\alpha (i)})^{-1},\RR^d) \circ \gamma_{\kappa_n^{\alpha (i)}}$. For each $i \in I$, there is an open $C^1$-neighborhood $\eE^i \subseteq \vect{W_{\alpha (i)}}$ of the zero-section such that the following holds:
  \begin{compactenum}
    \item  For $\gamma \in C^r ([0,1], \eE^i)$, we obtain a map $e(\gamma) \in C^{r+1} ([0,1], \vect{\Omega_{2,K_{5,i}}})$, defined locally via 
    \begin{equation}\label{eq: def:vf}
     e(\gamma) (t)(x) = (\exp_{W_{\alpha (i)}}|_{N_x})^{-1} \circ (\kappa_n^{\alpha (i)})^{-1} \circ \Fl_0^f (t, \kappa_n^{\alpha (i)} (x),\gamma_{[n]}), \quad (t,x) \in [0,1] \times V_{2,\alpha (i)}^n                                                                                                                                                                                                                                                                                                                                      
    \end{equation}
   for $f$ as in \ref{setup: regdiff} and $N_x$ as in \ref{lem: eos:cov}. Furthermore, for $S_{\alpha (i)}$ as in Construction \ref{con: sp:ofdat} V. and $(t,x) \in [0,1] \times V_{2 , \alpha (i)}^{n}$, the following estimates hold:
	  \begin{align} \label{eq: est:vf}
	\exp_{W_{\alpha (i)}} \circ e(\gamma ) (t)(x) \in V_{3,\alpha (i)}^n \quad \text{ and } \quad e(\gamma) (t)(x) \in B_{\rho_{\alpha (i)}} (0_x, S_{\alpha (i)})  .
          \end{align} 
   \item For each $\gamma \in C^r ([0,1], \eE^i)$, the map $e(\gamma)(0)$ is the zero section in $\vect{\Omega_{2,K_{5,i}}}$. If $\gamma$ is the constant map $\gamma \equiv 0_{W_{\alpha (i)}}$, then $e (\gamma) (t)$ is the zero-section for each $t \in [0,1]$. 
   \item The following maps are smooth  
      \begin{align*}
       \omega_i &\colon C^r ([0,1], \eE^i) \rightarrow C^{r+1}([0,1], \vect{\Omega_{2,K_{5,i}}}), \gamma \mapsto e(\gamma)	\\
       \theta_i &\colon C^r([0,1],\eE^i) \rightarrow \vect{\Omega_{2,K_{5,i}}}, \gamma \mapsto e(\gamma)(1).
      \end{align*}
   \end{compactenum}
 \end{lem}

\begin{proof}
  The set $\fF_5 (K_{5,i})$ is finite, whence by Lemma \ref{lem: eos:cov} (a), we can choose and fix $\nu > 0$ with the following properties: For each $y \in \overline{\Omega_{4,K_{5,i}}}$, the map $\exp_{W_{\alpha (i)}}$ is injective on 
    \begin{displaymath}
     N_y = \bigcup_{(V_{5,\alpha (i)}^n,\kappa^{\alpha (i)}_n) \in I_y} (T\kappa^{\alpha (i)}_n)^{-1} \left(\set{\kappa^{\alpha (i)}_n (y)} \times B_{\nu}  (0)\right),
    \end{displaymath}
  where $I_y = \setm{(V_{5,\alpha i}^n, \kappa_n^{\alpha (i)}) \in \fF_5 (K_{5,i})}{x \in \overline{V_{4,\alpha (i)}^n}}$. Lemma \ref{lem: eos:cov} (b) holds for the exponential maps $\exp_n$ associated to the pullback metric on $B_{5} (0)$ with respect to $\rho_{\alpha (i)}$ and $\kappa_n^{\alpha (i)}$.\\
  Consider $(V_{5,\alpha (i)}^n, \kappa_n^{\alpha (i)}) \in \fF_5 (K_{5,i})$. By Lemma \ref{lem: exp:loc}, there are constants $\ve_n > 0$ and $1 > \delta_n >0$ such that $a_{n}^{\alpha (i)} \colon B_4 (0) \times B_{\delta_n} (0) \rightarrow B_{\ve_n} (0_x), a_{n}^{\alpha (i)}(x,y) \coloneq \exp_n|_{B_{\ve_n} (0_x)}^{-1} (x+y)$ is a smooth map. Shrinking $\ve_n, \delta_n$, without loss of generality $\ve_n < \min \set{R_i, \nu}$ holds for the constant $R_i$ from Construction \ref{con: sp:ofdat} V. Recall that $\kappa_n^{\alpha (i)} (V_{5,\alpha (i)}^n) = B_5 (0)$, whence by Lemma \ref{lem: loc:eb} (b) there is a constant $0 < \tau_n \leq 1$ such that for $\gamma \in C^r([0,1],C^\infty (B_5 (0) , \RR^d))$ with $\sup_{t \in [0,1]} \norm{\gamma (t)}_{\overline{B_4 (0)},1} \leq \tau_n$, one has 
    \begin{equation}\label{eq: delta}
         \sup_{t \in [0,1]}    \norm{\Fl_0^f (t, \cdot , \gamma) - \id_{B_3 (0)}}_{\overline{B_2 (0)},1} < \delta_n.                
    \end{equation}
 Observe that $\delta_n < 1$ together with \eqref{eq: delta} implies $\Fl_0^f (t, \cdot , \gamma) (\overline{B_2 (0)}) \subseteq B_3 (0)$. Consider the open zero-neighborhood $E_n \coloneq \setm{f \in C^\infty (B_5 (0), \RR^d)}{\norm{f}_{\overline{B_{4} (0)},1} < \tau_n}$ and let
   \begin{displaymath}
    \eE_n^i \coloneq \setm{\sigma \in \vect{\Omega_{5,K_{5,i}}}}{\sigma_{[n]} \coloneq  \text{pr}_2 \circ T\kappa_n^{\alpha (i)} \circ \sigma \circ (\kappa_n^{\alpha (i)})^{-1} \in E_n}
   \end{displaymath}
  be the open neighborhood of the zero-section in $\vect{\Omega_{5,K_{5,i}}}$ induced by $E_n$. Repeating this construction, we obtain open neighborhoods of the zero map (respectively the zero-section) for each chart in $\fF_5 (K_{5,i})$. Let $V_i \coloneq \bigcap_{\fF_5 (K_{5,i})} \eE_n^i \subseteq \vect{\Omega_{5,K_{5,i}}}$. We show that the open zero-neighborhood $\eE^i \coloneq (\res_{\Omega_{5,K_{5,i}}}^{W_{\alpha (i)}})^{-1} (V_i) \subseteq \vect{W_{\alpha (i)}}$ satisfies the assertion of the lemma.
  \begin{compactenum}
  \item Consider $\gamma \in C^r ([0,1], V_i)$ and $(V_{5,\alpha (i)}^n,\kappa^{\alpha (i)}_n) \in \fF_5 (K_{5,i})$. The map $h_n$ sending $\gamma (t)$ to $\gamma_{[n]} (t)$ for $t\in [0,1]$ is continuous linear by \cite[Lemma F.6 and Lemma 4.11]{hg2004}. We deduce from \cite[Lemma 1.2]{gn2012} that $(h_n)_* \colon C^r([0,1], \vect{\Omega_{5,K_{5,i}}}) \rightarrow C^r ([0,1],C^\infty (B_5 (0), \RR^d)), \gamma \mapsto \gamma_{[n]}$ is continuous linear. Since $\gamma \in V_i$, we have $\gamma_{\kappa_n^{\alpha (i)}} \in C^r ([0,1] , E_n)$. By construction, \eqref{eq: delta} holds, $a_n^{\alpha (i)}$ is smooth and $\Fl_0^f (\cdot, \cdot ,\gamma_{[n]})$ a $C^{r+1,\infty}$-mapping by \ref{setup: regdiff}. By the Exponential Law \cite[Theorem 3.28 (e)]{alas2012}, a map in $C^{r+1} ([0,1],C^\infty (B_2 (0), \RR^d))$ may be defined via  
   \begin{equation}\label{eq: local}
    e(\gamma)_n (t)\coloneq a_n^{\alpha (i)} \circ (\id_{B_2 (0)} , \Fl_0^f (t,\cdot,\gamma_{[n]}) - \id_{B_2 (0)}),\quad t \in [0,1].
   \end{equation}
  Observe that $e(\gamma)_n (t) (B_2 (0)) \subseteq B_{\ve_n} (0)$ for each $t \in [0,1]$. The construction may be repeated for each chart in $\fF_5 (K_{5,i})$.
  As $\ve_n < \min \set{\nu, R_i}$, we obtain by definition of $\nu$ and $R_i$ for $(t,x) \in [0,1] \times B_2 (0)$:
    \begin{equation}\label{eq: vflow:est}
     T(\kappa_n^{\alpha (i)})^{-1}(x, e(\gamma)_n (t) (x)) \in N_{(\kappa_n^{\alpha (i)})^{-1} (x)} \cap B_{\rho_{\alpha (i)}} (0_{(\kappa_n^{\alpha (i)})^{-1} (x)} , S_{\alpha (i)}).
    \end{equation}
 By Lemma \ref{lem: eos:cov} (b), the formula \eqref{eq: local} is equivalent to the right hand side of \eqref{eq: def:vf}. From the uniqueness of the flow $\Fl_0^f (\cdot, \gamma_{[n]})$, we deduce that the mappings $e(\gamma)_n$ coincide on the intersections of their domains, whence we obtain a map $e(\gamma) \in C^{r+1} ([0,1], \vect{\Omega_{2,K_{5,i}}})$. The local representative of this time dependent vector field on $(V_{5,\alpha (i)}^n,\kappa^{\alpha (i)}_n) \in \fF_5 (K_{5,i})$ is $e(\gamma)_n$. \\
 For $x \in V_{2,\alpha (i)}^n$, the formula of $e(\gamma)_n$ together with Lemma \ref{lem: eos:cov} (b) allows us to compute 
      \begin{align*}
       (\exp_{W_{\alpha (i)}}|_{N_x})\circ e(\gamma) (t)(x) &= (\kappa_n^{\alpha (i)})^{-1} \exp_n e(\gamma)_n (t)(\kappa_n^{\alpha (i)} (x)) \\&= \kappa_n^{\alpha (i)} \Fl_0^f (t, \kappa_n (x), \gamma_{[n]}) 
	\in V_{3,\alpha (i)}^n.
      \end{align*}
  Furthermore, \eqref{eq: vflow:est} shows that the estimate \eqref{eq: est:vf} holds. The map $\res_{\Omega_{5,K_{5,i}}}^{W_{\alpha (i)}}$ is continuous linear by \cite[Lemma F.15]{hg2004}, whence $(\res_{\Omega_{5,K_{5,i}}}^{W_{\alpha (i)}})_* \colon C^r([0,1], \vect{W_{\alpha (i)}}) \rightarrow C^r([0,1], \vect{\Omega_{5,K_{5,i}}})$ is continuous linear by \cite[Lemma 1.2]{gn2012}. Assign to a map $\gamma \in C^r ([0,1], \eE^i)$ the vector field $e(\res_{\Omega_{5,K_{5,i}}}^{W_{\alpha (i)}} (\gamma))$. By abuse of notation, we will omit $\res_{\Omega_{5,K_{5,i}}}^{W_{\alpha (i)}}$ from now on, i.e.\ for $\gamma \in C^r([0,1],\eE^i)$, $e(\gamma) \coloneq e(\res_{\Omega_{5,K_{5,i}}}^{W_{\alpha (i)}} (\gamma))$. 
  \item The map $\Fl^f_0 (\cdot, \kappa_n^{\alpha (i)} (x) , \gamma_{[n]})$ is a solution to the initial value problem \eqref{eq:diffeq} with initial value $\Fl^f_0 (0, \kappa_n^{\alpha (i)} (x) , \gamma_{[n]}) = \kappa_n^{\alpha (i)} (x)$. We obtain $e(\gamma)(0)(x) = (\exp_{W_{\alpha (i)}}|_{N_x})^{-1} (x) = 0_x$ from \eqref{eq: def:vf}, since $\exp_{W_{\alpha (i)}} (0_x) = x$ holds and on $N_x$ the map $\exp_{W_{\alpha (i)}}$ is injective.\\
  If $\gamma_{[n]} \equiv 0$, its flow is defined as $\Fl^f_0 (t, \kappa_n^{\alpha (i)} (x) , 0)= \kappa_n^{\alpha (i)} (x)$. Analogous to the previous argument, $e(\gamma) (t)$ is the zero-section for each $t\in [0,1]$. 
  \item We prove the smoothness of $\omega^i, \delta^i$ via a patched mapping argument. To this end, consider the continuous linear maps $p_n^s \colon \vect{\Omega_{s,K_{5,i}}} \rightarrow C^\infty (B_s (0) ,\RR^d) , \sigma \mapsto \sigma_{\kappa_n^{\alpha (i)}} \circ (\kappa_n^{\alpha (i)})^{-1}|_{B_s (0)}$ for $s \in [1,5]$. By Definition \ref{defn: top:vect}, $p^s \coloneq (p^s_n)_{(V_{5,\alpha (i)}^n , \kappa_n^{\alpha (i)}) \in \fF_5 (K_{5,i})}$ is a topological embedding with closed image. Thus  Lemma \ref{lem: pb:patch} yields a topological embedding with closed image  
    \begin{displaymath}
     p^s_* \colon C^r([0,1],\vect{\Omega_{s,K_{5,i}}}) \rightarrow \bigoplus_{\fF_5 (K_{5,i})} C^r ([0,1],  C^\infty (B_s (0), \RR^d)) , \gamma \mapsto (p^s_n \circ \gamma)_{\fF_{5, K_{5,i}}}.
    \end{displaymath}
 Consider the maps $h^i \colon C^r([0,1],V_i) \rightarrow C^{r+1}([0,1],\vect{\Omega_{5,K_{5,i}}}), \gamma \mapsto e (\gamma )$. We claim that there are smooth maps $D_n$ 
 such that the following diagram is commutative:   
  \begin{displaymath}
	 \begin{xy}
  	\xymatrix{
     	C^r ([0,1], \eE^i) \ar[r]^{\res_{\Omega_{5,K_{5,i}}}^{W_{\alpha (i)}}} & C^r([0,1],V_i) \ar[rr]^-{h^i} \ar[d]_{p^5_{*}}   & & C^{r+1} ([0,1],\vect{\Omega_{2,K_{5,i}}} \ar[d]^{p^2_*})  \\
					 &*+!!<0pt,\the\fontdimen22\textfont2>{\displaystyle \bigoplus_{\fF_5 (K_{5,i})} C^r ([0,1],  E_n)} \ar@{.>}[rr]^-{\bigoplus_{\fF_5 (K_{5,i})}D_n} & & *+!!<0pt,\the\fontdimen22\textfont2>{\displaystyle \bigoplus_{\fF_5 (K_{5,i})}C^{r+1} ([0,1], C^\infty (B_2 (0), \RR^d))} 
 				 }
	\end{xy}
	\end{displaymath}
 Observe that the vertical arrows are given by embeddings with closed image and composition in the upper row yields $\omega_i = h^i \circ \res$. Since $\res$ is a smooth map, $\omega_i$ will be smooth if $h^i$ is smooth. If the claim is true, then by Proposition \ref{prop: pat:loc} $h^i$ and thus $\omega_i$ will be smooth. Consider the open sets $ \lfloor \overline{B_2 (0)}, B_{\delta_n} (0) \rfloor_\infty \subseteq C^\infty (B_3 (0) ,\RR^d)$ and define 
  \begin{displaymath}
   (a_n^{\alpha (i)})_* \colon \lfloor \overline{B_2 (0)}, B_{\delta_n} (0) \rfloor_\infty \rightarrow C^\infty (B_2 (0) , \RR^d),\quad (a_n^{\alpha (i)})_* (g) (x) \coloneq a_n^{\alpha (i)} (x,g (x)).
  \end{displaymath}
 By \cite[Proposition 4.23 (a)]{hg2004}, $(a_n^{\alpha (i)})_*$ is smooth, since $a_n^{\alpha (i)}$ is smooth. From Lemma \ref{lem: loc:eb} and the definition of $E_n$, we deduce that $F_n \colon C^r ([0,1], E_n) \rightarrow C^{r+1} ([0,1], C^\infty (B_3 (0),\RR^d))$, $F_n(\gamma)(t) \coloneq \Fl_0^f (t,\cdot ,\gamma)|_{B_3 (0)} - \id_{B_3 (0)}, \ t\in [0,1]$ is smooth. The estimate \eqref{eq: delta} yields $F_n(\gamma)([0,1]) \subseteq \lfloor \overline{B_2 (0)},B_{\delta_n} (0)\rfloor_\infty$. Thus $(a_n^{\alpha (i)})_* \circ F_n^\wedge \colon C^r ([0,1], E_n) \times [0,1] \rightarrow C^\infty (B_2 (0) , \RR^d)$ is a $C^{\infty,r+1}$-map by the Exponential Law \cite[Theorem 3.28 (e)]{alas2012} and \cite[Lemma 3.18]{alas2012}. Apply \cite[Corollary 3.8 and Theorem 3.28 (e)]{alas2012} to obtain a smooth map: 
    \begin{displaymath}
     D_n \colon  C^r ([0,1],  E_n) \rightarrow C^{r+1}([0,1], C^\infty (B_2 (0), \RR^d), \gamma \mapsto ((a_n^{\alpha (i)})_* \circ F_n^\wedge)^\vee (\gamma) = (a_n^{\alpha (i)})_* \circ F_n (\gamma)
    \end{displaymath}
 with $D_n (0) = 0$. A computation with \eqref{eq: def:vf} and Lemma \ref{lem: eos:cov} (b) shows that $\oplus_{\fF_5 (K_{5,i})} D_n$ makes the above diagram commutative. By \cite[Proposition 2.20]{alas2012}, we consider the smooth evaluation $\ve_1 \colon C^{r+1} ([0,1], \vect{\Omega_{2,K_{5,i}}}) \rightarrow \vect{\Omega_{2,K_{5,i}}}, \gamma \mapsto \gamma (1)$. Since $\theta_i = \ve_1 \circ \omega_i$ holds, $\theta_i$ is smooth.
\end{compactenum}
\end{proof}
  
\begin{lem}\label{lem: Osgen} \no{lem: Osgen}
 In the setting of Lemma \ref{lem: loc:comp}, define the open set $\eE \coloneq \Lambda_\cC^{-1} (\bigoplus_{i \in I} \eE^i) \subseteq \Osc{Q}$, where $\cC$ is the orbifold atlas introduced in \ref{setup: cC}. Let $r \in \NN_0 \cup \set{\infty}$. For each $i \in I$ and $\gamma \in C^r([0,1],\Osc{Q})$, we define $\gamma_{\alpha (i)} \colon [0,1] \rightarrow \vect{W_{\alpha (i)}}, t \mapsto (\gamma (t))_{\alpha (i)}$, where $ (\gamma (t))_{\alpha (i)}$ is the canonical lift of $\gamma (t)$ with respect to the chart $(W_{\alpha (i)}, H_{\alpha (i)}, \varphi_{\alpha (i)})$. 
  \begin{compactenum}
    \item If $\gamma \in C^r([0,1], \Osc{Q})$, then the map $\gamma_{\alpha (i)}$ is of class $C^r$ and for $i\in I$, the map $p_i \colon C^r([0,1], \Osc{Q}) \rightarrow C^r([0,1], \vect{W_{\alpha (i)}}), \gamma \mapsto \gamma_{\alpha (i)}$ is continuous linear.
    \item For each $\gamma \in C^r([0,1],\eE)$, we obtain a path $e(\gamma) \in C^{r+1} ([0,1], \Osc{Q})$ whose canonical lifts with respect to $\aA$ are given by $e(p_i (\gamma))|_{U_i}$ for $i \in I$.
   \end{compactenum}
\end{lem}

\begin{proof}
 \begin{compactenum}
  \item Pick $\gamma \in C^r([0,1],\Osc{Q})$. By construction, $\Lambda_\cC \circ \gamma \in C^r([0,1], \bigoplus_{i \in I} \vect{W_{\alpha (i)}})$ has compact image. Arguing as in the proof of Lemma \ref{lem: pb:patch}, $\gamma$ induces a family of maps $(\gamma_{\alpha (i)})_{i \in I} \in \bigoplus_{i \in I} C^r([0,1], \vect{W_{\alpha (i)}})$. 
  Recall from the Definition \ref{defn: os:top} of the c.s.\ orbisection topology that each map $\tau_{W_{\alpha (i)}} \colon \Osc{Q} \rightarrow \vect{W_{\alpha (i)}}, [\hat{\sigma}] \mapsto \sigma_{W_{\alpha (i)}}$ is continuous linear. By \cite[Lemma 1.2]{gn2012}, $p_i$ is a continuous linear map, as $p_i = (\tau_{W_{\alpha (i)}})_*$ holds.   
  \item Consider the family of time-dependent vector fields $(s \mapsto e(\gamma_{\alpha (i)})(s)|_{U_i} )_{i \in I}$ constructed in Lemma \ref{lem: loc:comp} (a). We claim that for fixed $s \in [0,1]$, these vector fields are a canonical family of lifts of an orbisection. It is sufficient to check the following stronger condition: \\ 
  \emph{For all $i,j \in I$ and any change of charts $\mu \colon \Omega_{2,K_{5, i}} \supseteq \dom \mu \rightarrow \cod \mu \subseteq \Omega_{2, K_{5, j}}$, $e(\gamma_{\alpha (j)}) (s) \circ \mu = T \mu \circ  e(\gamma_{\alpha (i)}) (s)|_{\dom \mu}$ holds.}\\
  We check the condition locally: Pick $x \in \dom \mu$ together with charts $(V_{5,\alpha (i)}^n,\kappa_{n}^{\alpha (i)}) \in \fF_5 (K_{5,i})$, $(V_{5,\alpha (j)}^m, \kappa_m^{\alpha (j)}) \in \fF_5 (K_{5,j})$ such that $x \in V_{2,\alpha (i)}^n$ and $\mu (x) \in V_{2,\alpha (j)}^m \subseteq \Omega_{2,K_{5,j}}$. Since $\gamma_{\alpha (i)} \in \eE^i$, \eqref{eq: delta} yields maps 
    \begin{align*}
     \varphi_x &\colon [0,1] \rightarrow V_{3,\alpha (i)}^n, t \mapsto (\kappa_{n}^{\alpha (i)})^{-1} \Fl_0^f (t,  \kappa_{n}^{\alpha (i)}(x) , \gamma_{\alpha (i)[n]}) \\
     \varphi_{\mu (x)} &\colon [0,1] \rightarrow V_{3,\alpha (j)}^m, t \mapsto (\kappa_{m}^{\alpha (j)})^{-1} \Fl_0^f (t,  \kappa_{m}^{\alpha (j)}(x) , \gamma_{\alpha (j) [m]}).
    \end{align*}
 These maps are $C^1$-integral curves for the (time-dependent) vector field $\gamma_{\alpha (i)}$ with initial condition $\varphi_x (0) = x$, respectively for $\gamma_{\alpha (j)}$ with $\varphi_{\mu (x)} (0) = \mu (x)$ (using the terminology of \cite[IV, \S 2]{langdgeo2001}). The charts in $\fF_5 (K_{5,\alpha (i)})$ are contained in some $Z^{r}_{\alpha (i)}$, by Construction \ref{con: sp:ofdat}. Since $x \in \kK_{\alpha (i)}$ and $\mu (x) \in \kK_{\alpha (j)}$, there is a change of charts $\lambda \colon Z_{\alpha (i)}^r \rightarrow W_{\alpha (j)}$ with $\lambda (x) = \mu(x)$. Composing $\lambda$ with a suitable element of $H_{\alpha (j)}$, without loss of generality there is an open neighborhood $U_x$ of $x$ with $\mu|_{U_x} = \lambda|_{U_x}$. The set $V_{5,\alpha (i)}^n$ is contained in $\dom \lambda$, whence $\lambda \circ \varphi_x \colon [0,1] \rightarrow W_{\alpha (j)}$ defines a $C^1$-curve such that $\lambda \circ \varphi_x (0) = \lambda (x) = \mu (x) \in \Omega_{2,K_{5,i}}$. For fixed $t \in [0,1]$, the vector fields $\gamma_{\alpha (i)} (t)$ and $\gamma_{\alpha (j)} (t)$ are members of a canonical family of lifts of an orbisection, i.e.\ $\gamma_{\alpha (j)} (t) \circ \lambda = T\lambda \gamma_{\alpha (i)} (t)|_{\dom \lambda}$. We compute: 
      \begin{displaymath}
       \gamma_{\alpha (j)} (t) (\lambda \varphi_x (t)) = T\lambda \gamma_{\alpha (i)} (t)(\varphi_x (t)) = T\lambda (\tfrac{\partial}{\partial t}\varphi_x) (t) = \tfrac{\partial}{\partial t} (\lambda \circ \varphi_x) (t).
      \end{displaymath}
 Thus the $C^1$-curve $\lambda \circ \varphi_x$ is an integral curve for the time-dependent vector field $\gamma_{\alpha (j)}$ with initial condition $\lambda \circ \varphi_x (0) = \lambda (x) = \mu (x)$.  On the other hand, the same is true for the $C^1$ curve $\varphi_{\mu (x)}$. As integral curves for (time-dependent) vector fields are unique (cf.\ \cite[IV.\ Theorem 2.1]{langdgeo2001} with \cite[p.\ 71]{langdgeo2001}) we derive $\lambda \circ \varphi_x = \varphi_{\mu (x)}$. \\
 Computing locally, we exploit that $\lambda \circ (\kappa_n^{\alpha (i)})^{-1}$ is a Riemannian embedding of $B_{5} (0)$ into $W_{\alpha (j)}$. In particular, by \cite[IV.\ Proposition 2.6]{fdiffgeo1963} the identity 
  \begin{displaymath}
  \exp_{W_{\alpha (j)}} T(\lambda (\kappa_n^{\alpha (i)})^{-1}) (v) = \lambda (\kappa_n^{\alpha (i)})^{-1} \exp_n (v) \quad \forall v \in \dom \exp_n
  \end{displaymath}
 holds. Notably, the estimates \eqref{eq: est:vf} and \eqref{eq: vflow:est} hold. With Lemma \ref{lem: eos:cov} (b) and the identity \eqref{eq: def:vf} for $e(\gamma_{\alpha (i)})$ on $[0,1] \times V_{2,\alpha (i)}^n$, one deduces from the above identity 
      \begin{align*}
       \exp_{W_{\alpha (j)}} T \lambda e(\gamma_{\alpha (i)}) (s)(x) &= \exp_{W_{\alpha (j)}} T \lambda (T\kappa_n^{\alpha (i)})^{-1} T\kappa_n^{\alpha (i)}  e(\gamma_{\alpha (i)}) (s)(x) \\
      &= \lambda (\kappa_n^{\alpha (i)})^{-1} \exp_n T\kappa_n^{\alpha (i)} e(\gamma_{\alpha (i)}) (s) (x) \\ 
      &= \lambda (\kappa_n^{\alpha (i)})^{-1} \kappa_n^{\alpha (i)} \exp_{W_{\alpha (i)}}|_{N_x} e(\gamma_{\alpha (i)} )(s) (x) \\ 
      &= \lambda (\kappa_n^{\alpha (i)})^{-1} \Fl_0^f(s,\kappa_n^{\alpha (i)} (x), \gamma_{\alpha (i) [n]})   = \lambda \circ \varphi_x (s) = \varphi_{\mu (x)} (s) .
      \end{align*}
  On the other hand, the local formula \eqref{eq: def:vf} for $e(\gamma_{\alpha (j)})$ on $[0,1] \times V_{2,\alpha (j)}^m$ implies 
    \begin{displaymath}
     \exp_{W_{\alpha (j)}} \circ e(\gamma_{\alpha (j)}) (s)(\mu (x)) = \varphi_{\mu (x)} (s) = \exp_{W_{\alpha (j)}} T \lambda e(\gamma_{\alpha (i)}) (s)(x).
    \end{displaymath}
  By construction, $\lambda (x) = \mu (x) \in  \kK_j^\circ$. Moreover, the mappings $e(\gamma_{\alpha (j)}) (s)$ and $e(\gamma_{\alpha (i)}) (s)$ are vector fields which satisfy the estimate \eqref{eq: est:vf}. Together with these facts, the definition of the constants (cf.\ Construction \ref{con: sp:ofdat} V.) yields: 
  \begin{displaymath}
   e(\gamma_{\alpha (j)}) (s)(\mu (x)) , T \lambda e(\gamma_{\alpha (i)}) (s) (x) \in B_{\rho_{\alpha (j)}} (0_{\mu (x)}, s_{\alpha (j)}) \subseteq \hat{O}_{\alpha (j)}.
  \end{displaymath}
 The map $\exp_{W_{\alpha (j)}}$ is injective on the intersection $\hat{O}_{\alpha (j)} \cap T_{\mu (x)}W_{\alpha (j)}$. Hence from the above identity $e(\gamma_{\alpha (j)}) (s) \circ \mu (x) = T \mu \circ e(\gamma_{\alpha (i)}) (s) (x)$ follows, thus proving the claim. Since $U_i$ is contained in $\Omega_{2,K_{5,i}}$, we deduce that the family $(e(\gamma_{\alpha (i)})(s)|_{U_i})_{i \in I}$ is a canonical family for an orbisection. Thus Remark \ref{rem: os:clpro} (a) shows that this family induces an orbisection $e(\gamma) (s)$. Observe that $\Lambda_\cC \circ \gamma ([0,1])$ factors through a finite subset of $\cC$ by \cite[Ch.\ III, \S 1, No.\ 4, Proposition 5]{bourbaki1987}. We derive from Lemma \ref{lem: loc:comp} (b) that there are only finitely many members of $(e(\gamma_{\alpha (i)}) (s))_{i \in I}$ which are not the zero-section. Assume that the finite subset $F \subseteq I$ satisfies $e(\gamma_{\alpha (i)})(\cdot)|_{U_i} \not \equiv 0_{U_i}$ if and only if $i \in F$. Then $\supp [e(\gamma) (s)] \subseteq \bigcup_{i \in F} \varphi_{\alpha (i)} (W_{\alpha (i)})$. Since each $\varphi_{\alpha (i)} (W_{\alpha (i)})$ is a relatively compact subset of $Q$, the orbisection $[e(\gamma)(s)]$ is compactly supported.\\
 We are left to prove that the assignment $[0,1] \rightarrow \Osc{Q}, s \mapsto e(\gamma) (s)$ is of class $C^{r+1}$. Identify $\Osc{Q}$ via $\Lambda_\aA$ with a sequentially closed subspace of $\bigoplus_{i \in I} \vect{U_i}$. It suffices to prove that $\Lambda_\aA \circ e(\gamma)$ is contained in $C^{r+1}([0,1], \bigoplus_{i \in I} \vect{U_i})$. The path $\Lambda_\aA \circ e(\gamma)$ factors through the inclusion $\bigoplus_{i \in F} \vect{U_i} \hookrightarrow \bigoplus_{i \in I} \vect{U_i}$. Each component is given by the $C^{r+1}$-path $t \mapsto e(p_i (\gamma)) (t)|_{U_i}$, whence $\Lambda_\aA \circ e(\gamma)$ is a path of class $C^{r+1}$ as a map to $\bigoplus_{i \in I} \vect{U_i}$. 
 \end{compactenum}
\end{proof}

To assure the smoothness of the evolution map on the Lie group, we exploit the patched locally convex structure of $\Osc{Q}$. Unfortunately $C^r([0,1],\Osc{Q})$ will inherit this structure only if $\Osc{Q}$ is countably patched (cf.\ Lemma \ref{lem: pb:patch}). To assure this condition, we require: \\[1em] \textbf{Convention:} For the rest of this section, we let $Q$ be a $\sigma$-compact (or second countable) space. 

\begin{lem} \label{lem: OSC:pat} \no{lem: OSC:pat}
 Let $Q$ be a $\sigma$-compact space and $r \in \NN_0$. The maps 
  \begin{align*} \omega &\colon C^r([0,1], \eE) \rightarrow C^{r+1} ([0,1], \Osc{Q}), \gamma \mapsto e(\gamma)\\
   \evol &\colon C^r([0,1],\eE) \rightarrow \Osc{Q}, \gamma \mapsto e(\gamma) (1) 
  \end{align*}
are smooth and map the constant path $\gamma \equiv \zs$ to itself respectively to $\zs$.
\end{lem}

\begin{proof}
  The topological space $Q$ is $\sigma$-compact and $\aA, \cC$ are locally finite, whence $I$ is countable. Corollary \ref{cor: ost:prop} (c) shows that the mappings $\Lambda_\aA, \Lambda_\cC$ turn $\Osc{Q}$ into a patched locally convex space. As $r < \infty$ holds,  the spaces $ \bigoplus_{i \in I} C^{r+1} ([0,1], \vect{U_i})$ and $C^{r+1} ([0,1], \bigoplus_{i \in I} \vect{U_i})$ are isomorphic by the proof of Lemma \ref{lem: pb:patch}. The same is true if we replace each $U_i$ with $W_{\alpha (i)}$. For $\aA$ as in Construction \ref{con: sp:ofdat} and $\cC$ as in \ref{setup: cC}, we identify these spaces to consider the mappings 
      \begin{align*}
	P_\aA &\colon C^{r+1} ([0,1], \Osc{Q}) \rightarrow  \bigoplus_{i \in I} C^{r+1} ([0,1], \vect{U_i}), \gamma \mapsto \Lambda_\aA \circ \gamma = (\gamma_{U_i})_{i \in I}. \\
	P_\cC &\colon C^r ([0,1], \Osc{Q}) \rightarrow  \bigoplus_{i \in I} C^r ([0,1], \vect{W_{\alpha (i)}}), \gamma \mapsto \Lambda_\cC \circ \gamma = (\gamma_{\alpha (i)})_{i \in I}. 
      \end{align*}
 An application of Lemma \ref{lem: pb:patch} proves: $P_\aA , P_\cC$ are linear topological embeddings with closed image, whose components form patchworks, for $C^{r+1}([0,1], \Osc{Q})$ and $C^r([0,1], \Osc{Q})$, respectively.
 The maps $\omega$ and $\evol$ are well-defined by Lemma \ref{lem: Osgen} (b) and we claim that they are smooth.  
 For $i \in I$, let $\res_{U_i}^{\Omega_{2,K_{5,i}}} \colon \vect{\Omega_{2,K_{5,i}}} \rightarrow \vect{U_i}$ be the restriction map. These mappings are linear and continuous by \cite[Lemma F.15 (a)]{hg2004}. Thus $r_i \coloneq C^{r+1} ([0,1],\res_{U_i}^{\Omega_{2,K_{5,i}}})$ is continuous and linear by \cite[Lemma 1.2]{gn2012}, hence a smooth map. For $i \in I$, consider the smooth map $\omega_i$ defined in Lemma \ref{lem: loc:comp}. By Lemma \ref{lem: loc:comp} (b) the smooth map $r_i \circ \omega_i$ maps the constant path $\gamma \equiv O_{W_{\alpha (i)}}$ to the constant path whose image is the zero-section. From the definitions we obtain  
  \begin{align} \label{patchedevolution}
  \left( \bigoplus_{i \in I} r_i\circ \omega_i\right) \circ P_\cC|_{C^r([0,1],\eE)}^{\oplus_{i \in I} C^r([0,1],\eE^i)} = \Lambda_\aA \omega.
  \end{align} 
 Hence $\omega$ is smooth on the patches and we deduce from \eqref{patchedevolution} with Proposition \ref{prop: pat:loc} that $\omega$ is a smooth map. As the evaluation map $\text{ev}_1 \colon C^{r+1} ([0,1], \Osc{Q}) \rightarrow \Osc{Q} , \gamma \mapsto \gamma (1)$ is smooth (cf.\ \cite[Proposition 3.20]{alas2012}), the smoothness of $\evol$ follows from $\text{ev}_1 \circ \omega = \evol$ . The last assertion is a direct consequence of Lemma \ref{lem: loc:comp} (b).  
\end{proof}

\begin{lem}\label{lem: znhd:cst} \no{lem: znhd:cst}
 Let $\hH_\rho \subseteq \Osc{Q}$ be the open zero-neighborhood of Theorem \ref{thm: diff:lgp}. Consider an open identity-neighborhood $\sS \subseteq E(\hH_\rho)$ which is symmetric, i.e.\ $\sS = \sS^{-1}$. There is an open subset $\zs \in \rR \subseteq \eE \subseteq \Osc{Q}$ such that $\omega (C^r([0,1],\rR)) \subseteq C^{r+1}([0,1], E^{-1} (\sS))$.
\end{lem}

\begin{proof}
 Consider the $C^0$-neighborhood of the constant path $\gamma_{\zs} \equiv \zs$: 
  \begin{displaymath}
   C^{r+1} ([0,1], E^{-1} (\sS)) \coloneq C^0 ([0,1], E^{-1}(\sS)) \cap C^{r+1} ([0,1], \Osc{Q}).
  \end{displaymath}
 Specializing to $r=0$ in Lemma \ref{lem: OSC:pat} we see that $\omega \colon C^0 ([0,1], \eE) \rightarrow C^1([0,1], \Osc{Q})$ is smooth with $\omega (\gamma_{\zs}) = \gamma_{\zs}$. Then $\omega^{-1} (C^{1}([0,1], E^{-1} (\sS))) \subseteq C^0([0,1] ,\eE)$ is an open zero-neighborhood. The definition of the compact open topology yields an open set $\zs \in \rR \subseteq \Osc{Q}$ such that $\gamma_{\zs} \in C^0 ([0,1], \rR) \subseteq \omega^{-1} (C^{1}([0,1], E^{-1} (\sS)))$. The assertion follows.
\end{proof}
\noindent
Observe that, by construction, also $\evol (C^r([0,1], \rR)) \subseteq \hH_\rho$. We shall see presently that with the maps constructed in Lemma \ref{lem: OSC:pat}, a smooth evolution for the Lie group $\Difforb{Q,\uU}$ may be constructed. We would like to apply methods similar to the manifold case (cf.\ \cite[p.\ 1046]{milnor1983}) to prove the regularity of $\Difforb{Q,\uU}$.  However, if $(Q,\uU)$ is a non-trivial orbifold, it is more difficult to verify the existence of right logarithmic derivatives. We need representatives of the orbifold diffeomorphisms in $\sS$ tailored to this purpose:

\begin{lem}\label{lem: nlifts} \no{lem: nlifts}
 Consider $[\hat{f}] \in \sS$ with $[\hat{f}] = [\hat{E^\sigma}]$ for some $[\sigma] \in \hH_\rho$. For each $[\hat{g}] \in \sS$, there is a representative $E_{\hat{f}} (\hat{g})$ of $[\hat{g}]$ with lifts $\set{E_{\hat{f}} (\hat{g})_i}_{i\in I}$ such that the following properties are satisfied:
  \begin{compactenum}
   \item for each $i \in I$, the lift $E_{\hat{f}} (\hat{g})_i$ is an \'{e}tale embedding in $C^\infty (e^{\sigma_i} (U_i), W_{\alpha (i)})$ (cf.\ Lemma \ref{lem: eos:loc}),
   \item if $[\hat{g}] = [\hat{f}]^{-1}$ holds, then the lifts are given by $E_{\hat{f}} (\hat{f}^{-1})_i = (e^{\sigma_i})^{-1}$ for all $i \in I$.
\end{compactenum}
\end{lem}
 
\begin{proof}
 Let $[\hat{\tau}^g]$ be the unique preimage of $[\hat{g}]$ with respect to $E$. From $[\hat{g}] = E ([\hat{\tau}^g]) = [\expo ] \circ [\hat{\tau}^g]|^\Omega$ we deduce that the claim will hold if there are representatives of $[\expo]$ and $[\hat{\tau}^g]|^\Omega$ whose composition yields the desired representative. The map $[\hat{f}]$ is an orbifold diffeomorphism with representative $\widehat{E^\sigma} = (E^\sigma , \set{e^{\sigma_i}}_{i \in I}, [P,\nu])$. Hence the orbifold charts $\set{(e^{\sigma_i} (U_i) , G_i,\varphi_{\alpha (i)}|_{e^{\sigma_i} (U_i)})}_{i\in I}$ (cf.\ Lemma \ref{lem: eos:loc}) cover $Q$. Recall the following details from the proof of Lemma \ref{lem: eos:loc}:\\ 
 By Step 3, $H_{\alpha (i)} . \im e^{\sigma_i} \subseteq \Omega_{2,i}$ is an invariant subset such that $\im e^{\sigma_i}$ is $H_{\alpha (i)}$-stable. Using Lemma \ref{lem: eos:loc} iii., the canonical lifts $\tau^g_{\alpha (i)}$ map $\im e^{\sigma_i}$ into $\hat{O}_{\alpha (i)}$. Thus $(\tau^g_{\alpha (i)}|_{\im e^{\sigma_i}}^{\hat{O}_{\alpha (i)}})_{i\in I}$ is a family of lifts for a representative $\hat{\tau}'$ of $[\hat{\tau}^g]|^\Omega$. As $\Omega_{2,i} \subseteq \kK_{\alpha (i)}^\circ$, we obtain an open subset $T\im e^{\sigma_i} \cap \hat{O}_{\alpha (i)} \subseteq \dom \exp_{W_{\alpha (i)}}$ (cf.\ Construction \ref{con: sp:ofdat} IV.). This set is $G_i$-stable, whence $\exp_{W_{\alpha (i)}}|_{T \im e^{\sigma_i} \cap \hat{O}_{\alpha (i)}}$ is a lift of the orbifold exponential map $\expo$. By Remark \ref{rem: lift:expo} (a), there is a representative $\widetilde{\expo}$ of $[\expo]$ whose family of lifts contains $\set{\exp_{W_{\alpha (i)}}|_{T \im e^{\sigma_i} \cap \hat{O}_{\alpha (i)}}}_{i \in I}$. Composing $\widetilde{\expo}$ and $\hat{\tau}'$, we obtain a representative of $E([\hat{\tau}^g]) = [\hat{g}]$ whose lifts are the smooth mappings  
  \begin{equation}\label{eq: nmap}
   E (\hat{f}; \hat{E}^{\tau'})_i \coloneq (\exp_{W_{\alpha (i)}}|_{T\im e^{\sigma_i} \cap O_{\alpha (i)}}) \circ \tau^g_{\alpha (i)}|_{\im e^{\sigma_i}}^{O_{\alpha (i)}}.
  \end{equation}
 As a consequence of the proof of Lemma \ref{lem: eos:loc}, these maps are equivariant \'{e}tale embeddings.\\
 Since $e^{\sigma_i}$ is a lift for $[\hat{f}]$ for each $i \in I$, the map $E (\hat{f}; \hat{f}^{-1})_i \circ e^{\sigma_i}$ is a change of orbifold charts. Hence for each $i\in I$, there is a unique $\gamma_i^{f^{-1}} \in H_{\alpha (i)}$ such that $\gamma_i^{f^{-1}} \circ E (\hat{f} ; \hat{f}^{-1})_i = (e^{\sigma_i})^{-1}$. 
 The family $(\gamma_i^{f^{-1}})_{i \in I}$ induces a lift of the identity $\hat{\ve}$ by Proposition \ref{prop: lift:unit}. We obtain another representative $\hat{\ve} \circ \widetilde{\expo} \circ \tau'$ of $E([\hat{\tau}^g])$, whose lifts $E_{\hat{f}} (\hat{g})_i \coloneq \gamma_i^{f^{-1}} \circ E (\hat{f} ; \hat{g})_i$, $i \in I$ are  \'{e}tale embeddings. Furthermore, for $[\hat{g}] = [\hat{f}]^{-1}$, by construction assertion (b) holds.      
\end{proof}
\begin{rem}\label{rem: gp:ident}
 \begin{compactenum}
  \item The construction of $E(\hat{f}; \hat{g})$ in Lemma \ref{lem: nlifts} (combine $H_{\alpha (i)} . \im e^{\sigma_i} \subseteq \Omega_{2,i}$ (see step 3 of the proof of Lemma \ref{lem: eos:loc}) with Lemma \ref{lem: eos:loc} iii.) shows that we can define maps $E_{\hat{f}}^{\hat{g}} \coloneq \exp_{W_{\alpha (i)}} \circ \tau^g_{\alpha (i)}|_{H_{\alpha (i)}.(\im e^{\sigma_i})}$ with $E_{\hat{f}}^{\hat{g}}|_{\im e^{\sigma_i}} = E(\hat{f};\hat{f}^{-1})$. As each $\tau_{\alpha (i)}^g$ is a canonical lift of an orbisection, we deduce that $\eta \circ  E_{\hat{f}}^{\hat{g}} = E_{\hat{f}}^{\hat{g}} \circ \eta$ for each $\eta \in H_{\alpha (i)}$. 
  \item Let $[\hat{f}] = \ido{}$ and consider $\gamma_i^{f}$ as in the proof of Lemma \ref{lem: nlifts}. Then $\gamma_i^{\ido{}} = \id_{U_i}$, for each $i \in I$. To see this, observe the identities $\ido{} = \ido{}^{-1}$ and $E^{-1} (\ido{}) = \zs$. For $i \in I$, both lifts constructed in \eqref{eq: nmap} coincide as $\id_{U_i} = \exp_{W_{\alpha (i)}} \circ 0_{U_i}$. This forces the identity $\gamma_i^{\ido{}} = \id_{U_i}$. 
 \end{compactenum}
\end{rem}

\begin{defn}\label{defn: nspace} \no{defn: nspace}
 For $[\hat{\phi}]$ in $\sS$, let $[\widehat{\sigma^{\phi}}]$ be the unique orbisection in $\hH_\rho$ with $E ([\widehat{\sigma^{\phi}}]) = [\hat{\phi}]$. Apply Lemma \ref{lem: nlifts} to $[\hat{\phi}]^{-1} \in \sS$. By Part (b) of Lemma \ref{lem: nlifts}, we obtain a representative $\hat{\phi}$ of $[\hat{\phi}]$. For each $i \in I$ the lifts $g_i^\phi \coloneq E_{\hat{\phi}^{-1}} (\hat{\phi})_i$ of $\hat{\phi}$  are embeddings of $U_{\phi_i} \coloneq \exp_{W_{\alpha (i)}} (\sigma_i^{\phi^{-1}} (U_i)) \subseteq \Omega_{2,i}$ with $\im g_i^\phi = U_i$. The pointwise operations make 
  \begin{displaymath}
   C_i^\phi \coloneq \setm{f \in C^\infty (U_{\phi_i}, TW_{\alpha (i)})}{\pi_{TW_{\alpha (i)}} \circ f =g_i^\phi} 
  \end{displaymath}
 a vector space. Endow $ C_i^\phi$ with the unique topology turning $(g_i^\phi)^* \colon \vect{U_i} \rightarrow C_i^\phi, \sigma_i \mapsto \sigma_i \circ g_i^\phi$ into an isomorphism of topological vector spaces. We define a linear map 
  \begin{displaymath}
   \Lambda_{[\hat{\phi}]} \colon C_{[\hat{\phi}]} \coloneq \setm{[\hat{\sigma}] \circ [\hat{\phi}]}{[\hat{\sigma}] \in \Osc{Q}} \rightarrow \bigoplus_{i \in I} C_i^\phi , [\hat{\sigma}] \circ [\hat{\phi}] \mapsto (\sigma_i \circ g_i^\phi)_{i \in I},
  \end{displaymath}
where $\sigma_i$ is the canonical lift of $[\hat{\sigma}]$ on $U_i$. As orbisections are uniquely determined by a family of canonical lifts, the map $\Lambda_{[\hat{\phi}]}$ is injective. Endow $C_{[\hat{\phi}]}$ with the unique locally convex topology turning $\Lambda_{[\hat{\phi}]}$ into a topological embedding.\\ The lifts $g^{\ido{}}_i$ are the identity on $U_i$ for each $i \in I$, by Remark \ref{rem: gp:ident} (b). Therefore $C_{\ido{}}$ and $\Osc{Q}$ coincide, and hence the mappings $\Lambda_{\ido{}}$ and $\Lambda_\aA$ are the same. 
\end{defn}

For the rest of this section, fix the notation of Definition \ref{defn: nspace}. We obtain a structural result for the tangent manifold of $\Difforb{Q,\uU}$:

\begin{lem}\label{lem: char:tdiff} \no{lem: char:tdiff}
 Let $[\hat{\phi}]$ be an element of $\sS$ with $\sS$ as in Lemma \ref{lem: znhd:cst}. There is an isomorphism of topological vector spaces 
  \begin{displaymath}
    \alpha_{[\hat{\phi}]} \colon T_{[\hat{\phi}]} \Difforb{Q,\uU} \rightarrow \Im \Lambda_{[\hat{\phi}]},
  \end{displaymath}
whence $T_{[\hat{\phi}]} \Difforb{Q,\uU}$ is isomorphic as a topological vector space to $C_{[\hat{\phi}]}$. 
\end{lem}

\begin{proof}
 Fix $[\hat{\phi}] \in \sS$. As $\sS$ is a symmetric set (i.e.\ $\sS = \sS^{-1}$), the inverse $[\hat{\phi}]^{-1}$ of $[\hat{\phi}]$ is contained in $\sS$. By construction of $\sS$, there is a representative of $[\hat{\phi}]^{-1}$ with lifts $\set{(g^\phi_i)^{-1} \colon U_i \rightarrow W_{\alpha (i)}}_{i \in I}$. To shorten our notation, we set $U_{\phi_i} \coloneq (g^\phi_i)^{-1} (U_i)$ and recall $U_{\phi_i} \subseteq \Omega_{2,i}$ from Definition \ref{defn: nspace}. The family of lifts $\set{g^{\phi}_i}_{i \in I}$ uniquely determines a representative of $[\hat{\phi}]$, by Corollary \ref{cor: diff:char}. We proceed in several steps: 
\paragraph{Step 1:} \emph{Construct the mapping $\alpha_{[\hat{\phi}]}$}. For each $[\hat{g}] \in \sS$, denote by $[\hat{\sigma}^g]$ the compactly supported orbisection with $E([\hat{\sigma}^g]) = [\hat{g}]$. By Lemma \ref{lem: nlifts} (a) each $[\hat{g}] \in \sS$ possesses a representative $E_{\hat{\phi}^{-1}} (\hat{g})$ with lifts $(E_{\hat{\phi}^{-1}}(\hat{g}))_i \coloneq \gamma_i^{\phi}\exp_{W_{\alpha (i)}} \circ \sigma^g_{\alpha (i)}|_{U_{\phi_i}}$. Fix $i \in I, p \in U_{\phi_i}$ and consider the map 
      \begin{displaymath}
       \ve^{\phi_i}_p \colon \sS \rightarrow W_{\alpha (i)}, [\hat{g}] \mapsto E_{\hat{\phi}^{-1}}(\hat{g})_i (p).
      \end{displaymath}
 We show that $\ve^{\phi_i}_p$ is smooth. To this end. let $\tau_{W_{\alpha (i)}} \colon \Osc{Q} \rightarrow \vect{W_{\alpha (i)}}$ be the map which sends an orbisection to its canonical lift on $W_{\alpha (i)}$. By Definition \ref{defn: os:top} (b), this map is continuous linear, hence smooth. Choose a manifold chart $(V_p, \psi_p)$ of the manifold $W_{\alpha (i)}$ with $p \in V_p$. The map $\tau_{V_p} \colon \vect{W_{\alpha (i)}} \rightarrow C^\infty (V_p, \RR^d), X \mapsto X_{\psi_p} \coloneq \text{pr}_2 T\psi_p X|_{V_p}$ is continuous linear by Definition \ref{defn: top:vect}. Let $\ve_p \colon C^\infty (V_p,\RR^d) \rightarrow \RR^d , f \mapsto f(p)$ be the evaluation map in $p$. This map is a linear map, which is smooth by \cite[Proposition 3.20]{alas2012}. Finally define $\text{ev}_p \colon \vect{W_{\alpha (i)}} \rightarrow T_p W_{\alpha (i)}, X \mapsto X(p)$. As $\text{ev}_p = (T_p \psi_p)^{-1} (\psi_p (p), \cdot) \circ \ve_p \circ \tau_{V_p}$ holds, $\text{ev}_p$ is continuous linear. By construction of $\hH_\rho$, it is contained in the open subset $\mM$ constructed in Proposition \ref{prop: diff:chart} (cf.\ Construction \ref{con: H}). Hence Lemma \ref{lem: eos:loc} ii.\ implies that  $\text{ev}_p$ maps $\tau_{W_{\alpha (i)}} \circ E^{-1} (\sS) \subseteq \mM_i$ into the set $\hat{O}_{\alpha (i)} \cap T_p W_{\alpha (i)}$. The image of the smooth map $\text{ev}_p \circ \tau_{W_{\alpha (i)}} \circ E^{-1}|_\sS$ is thus contained in $\dom \exp_{W_{\alpha (i)}} \cap T_p W_{\alpha (i)}$. By construction of the lifts $E_{\hat{\phi}^{-1}}(\hat{g})_i$ in Lemma \ref{lem: nlifts}, one may rewrite $\ve^{\phi_i}_p$ as composition of smooth maps, thus establishing the desired smoothness: 
    \begin{displaymath}
     \ve_p^{\phi_i} = \gamma_i^\phi \circ \exp_{W_{\alpha (i)}}|_{T_p W_{\alpha (i)}} \circ \text{ev}_p \circ \tau_{W_{\alpha (i)}} \circ E^{-1}|_\sS.
    \end{displaymath}
 Repeating the construction for each pair $p \in U_{\phi_i}$, where $i$ runs through $I$, we obtain a map 
  \begin{align*}
   \alpha_{[\hat{\phi}]} \colon T_{[\hat{\phi}]} \Difforb{Q,\uU}& \rightarrow \prod_{i \in I} (TW_{\alpha (i)})^{U_{\phi_i}}\\
					      V& \mapsto (T_{[\hat{\phi}]} \ve^{\phi_i}_p (V))_{i \in I,p \in U_{\phi_i}}
  \end{align*}
 and abbreviate its image as $V_{[\hat{\phi}]} \coloneq \im \alpha_{[\hat{\phi}]}$.
\paragraph{Step 2:} \emph{Endow $V_{[\hat{\phi}]}$ with a vector space structure which turns $\alpha_{[\hat{\phi}]}$ into a linear map.}\\ The tangent space $T_{[\hat{\phi}]} \Difforb{Q,\uU}$ is the set of equivalence classes of $C^1$-curves $\eta \colon ]-\ve , \ve[ \rightarrow \sS$ with $\eta (0) = [\hat{\phi}]$, where $\eta \sim \theta$ if and only if $(E^{-1} \circ \eta )' (0) = (E^{-1} \circ \theta )' (0)$. Abbreviate the equivalence classes with respect to this relation by $[t \mapsto \eta (t)]_\sim$ (and likewise in $TW_{\alpha (i)}$). Since each $\ve_p^{\phi_i}$ is smooth and $\eta$ is of class $C^1$, for each $i \in I$ and $p \in U_{\phi_i}$ the curve $\ve_p^{\phi_i} \circ \eta$ is of class $C^1$. Hence the definition of $\alpha_{[\hat{\phi}]}$ yields  
    \begin{equation}\label{eq: sim:vf}
     \alpha_{[\hat{\phi}]} ([\eta]_{\sim}) = ([t \mapsto E_{[\hat{\phi}]^{-1}}(\eta (t))_i (p)]_\sim)_{i \in I,p \in U_{\phi_i}}.
    \end{equation}
 The curve $\eta$ in \eqref{eq: sim:vf} passes through $[\hat{\phi}]$ for $t=0$, whence by Lemma \ref{lem: nlifts} (b) for $i \in I$, $E_{\hat{\phi}^{-1}} (\eta (0))_i = g_i^\phi$ holds. Therefore we infer from \eqref{eq: sim:vf} the identity  
    \begin{equation}\label{eq: gna}
     V_{[\hat{\phi}]} \subseteq \setm{(f_i)_{i \in I} \in \prod_{i \in I} (TW_{\alpha (i)})^{U_{\phi_i}}}{\forall i \in I, p \in U_{\phi_{i}},\ f_i (p) \in T_{g^\phi_i (p)} W_{\alpha (i)}}. 
    \end{equation}
In particular, \eqref{eq: gna} shows that the pointwise operations turn $V_{[\hat{\phi}]}$ into a vector space. Furthermore, by \eqref{eq: gna} $T_{[\hat{\phi}]} \ve^{\phi_i}_p \colon T_{[\hat{\phi}]} \Difforb{Q,\uU} \rightarrow T_{g^\phi_i (p)} W_{\alpha (i)}$ is linear. By definition, the map $\alpha_{[\hat{\phi}]}$ becomes linear if $V_{[\hat{\phi}]}$ is endowed with the vector space structure induced by pointwise operations.
\paragraph{Step 3:} \emph{A formula relating $\alpha_{[\hat{\phi}]}$ to $\alpha_{\ido{}}$}. Let $\rho_{[\hat{\phi}]} \colon \Difforb{Q,\uU} \rightarrow \Difforb{Q,\uU} , [\hat{\psi}] \mapsto [\hat{\psi}] \circ [\hat{\phi}]$ be the right translation and define 
    \begin{displaymath}
     G^\phi \coloneq ((g^\phi_i)_{i \in I})^* \colon \prod_{i \in I} (TW_{\alpha (i)})^{U_{i}} \rightarrow \prod_{i \in I} (TW_{\alpha (i)})^{U_{\phi_i}}, (f_i)_{i \in I} \mapsto (f_i \circ g^{\phi}_i)_{i\in I}
    \end{displaymath}
Consider $[\eta]_\sim \in T_{\ido{}} \Difforb{Q,\uU}$. The composition in $\Difforb{Q,\uU}$ is continuous, as the latter is a Lie group. Since $\eta (0) = \ido{}$ holds, we may thus assume $\eta (t) \circ [\hat{\phi}] \in \sS$ for all $t$. By Lemma \ref{lem: nlifts} (a), there is a representative of $\eta (t) \circ [\hat{\phi}]$ with lifts $E_{\hat{\phi}^{-1}} (\eta (t) \circ \hat{\phi})_i = \gamma_i^\phi \exp_{W_{\alpha (i)}} \sigma^{\eta (t) \circ \phi}_{\alpha (i)}|_{U_{\phi_i}}$. Here $\sigma^{\eta (t) \circ \phi}_{\alpha (i)}$ is the canonical lift on $W_{\alpha (i)}$ of the compactly supported orbisection $[\widehat{\sigma^{\eta (t) \circ \phi}}]$ with $E([\widehat{\sigma^{\eta (t) \circ \phi}}]) = \eta (t) \circ [\hat{\phi}]$. The set $U_{\phi_i}$ is contained in $\Omega_{2,i} \subseteq \Omega_{\frac{5}{4}, K_{5,i}}$ (cf.\ Construction \ref{con: H}). By Remark \ref{rem: new:rem}, we thus have $\sigma^{\eta (t) \circ \phi}_{\alpha (i)}|_{U_{\phi_i}} = \sigma_{\alpha (i)}^{\eta (t)} \diamond_i \sigma_{\alpha (i)}^{\phi}|_{U_{\phi_i}}$. Recall that by construction of $\sigma_{\alpha (i)}^{\eta (t)} \diamond_i \sigma_{\alpha (i)}^{\phi}|_{U_{\phi_i}}$ (see \eqref{eq: def1} in Construction \ref{con: comp:loc}) the identity 
  \begin{displaymath}
   \exp_{W_{\alpha (i)}} \circ \sigma_{\alpha (i)}^{\eta (t)} \diamond_i \sigma_{\alpha (i)}^{\phi}|_{U_{\phi_i}} = \exp_{W_{\alpha (i)}} \circ \sigma^{\eta (t)}_{\alpha (i)} \circ \exp_{W_{\alpha (i)}} \circ \sigma_{\alpha (i)}^{\phi}|_{U_{\phi_i}}
  \end{displaymath}
 holds. Furthermore, $g_i^\phi = E_{\hat{\phi}^{-1}}(\hat{\phi})_i = \gamma^\phi_i \circ \exp_{W_{\alpha (i)}} \circ \sigma_{\alpha (i)}^\phi|_{U_{\phi_i}}$ and $\im g_i^\phi = U_i$. Hence we deduce $\exp_{W_{\alpha (i)}} \circ \sigma_{\alpha (i)}^\phi (U_{\phi_i}) \subseteq H_{\alpha (i)}. U_i \subseteq \Omega_{2,i}$. Analogous to Step 2 in the proof of Lemma \ref{lem: eos:loc}, one shows that $\gamma_i^\phi \in H_{\alpha (i)}$ commutes with $\exp_{W_{\alpha (i)}} \circ \sigma_{\alpha (i)}^{\eta (t)}|_{H_{\alpha (i)}.U_i}$. Summing up, we obtain: 
  \begin{align*}
 \alpha_{[\hat{\phi}]}(T\rho_{[\hat{\phi}]}([\eta]_\sim)) &=\left([t\mapsto E_{\hat{\phi}^{-1}} (\eta (t) \circ \hat{\phi})_i (p)]_\sim\right)_{i \in I, p \in U_{\phi_i}} \\
							    &=\left([t\mapsto \gamma_i^\phi \circ \exp_{W_{\alpha (i)}} \sigma^{\eta (t)}_{\alpha (i)} \diamond_i \sigma^{\phi}_{\alpha (i)} (p)]_\sim\right)_{ i \in I , p \in U_{\phi_i}} \\
							    &=\left([t \mapsto \gamma_i^\phi \exp_{W_{\alpha (i)}} \sigma^{\eta (t)}_{\alpha (i)} \exp_{W_{\alpha (i)} } \sigma^\phi_{\alpha (i)} (p)]_\sim \right)_{i \in I, p \in U_{\phi_i}} \\
							    &=\left([t \mapsto \exp_{W_{\alpha (i)}} \sigma^{\eta (t)}_{\alpha (i)} \underbrace{\gamma_i^\phi \exp_{W_{\alpha (i)} } \sigma^\phi_{\alpha (i)}}_{=g_i^\phi} (p)]_\sim \right)_{i \in I, p \in U_{\phi_i}} \\
							    &= G^\phi \circ \alpha_{\ido{}} ([\eta]_\sim) 
  \end{align*}
We derive $\alpha_{[\hat{\phi}]} \circ T\rho_{[\hat{\phi}]}|_{\dom \alpha_{\ido{}}} = G^\phi \circ \alpha_{\ido{}}$. Now $G^\phi (V_{\ido{}}) = V_{[\hat{\phi}]}$ follows, as $T\rho_{[\hat{\phi}]}$ is a diffeomorphism.
\paragraph{Step 4:} \emph{$G^\phi|_{V_{\ido{}}}$ is linear}. To see this, let $v,w \in T_{\ido{}} \Difforb{Q}$ and $r \in \RR$. Since $T\rho_{[\hat{\phi}]}$, $\alpha_{[\hat{\phi]}}$ and $\alpha_{\ido{}}$ are linear, the formula in Step 3 yields: 
  \begin{align*}
  G^\phi (\alpha_{\ido{}} (v+rw)) 			&= \alpha_{[\hat{\phi}]} (T\rho_{[\hat{\phi}]} (v+rw)) \\
							&=  \alpha_{[\hat{\phi}]} (T\rho_{[\hat{\phi}]} (v)) + r\alpha_{[\hat{\phi}]} (T\rho_{[\hat{\phi}]} (w))\\
							&= G^\phi (\alpha_{\ido{}} (v)) + r G^\phi (\alpha_{\ido{}} (w)).
  \end{align*}
\paragraph{Step 5:} \emph{$\alpha_{\ido{}}$ is an isomorphism of topological vector spaces and $V_{\ido{}} = \im \Lambda_\aA$.}\\ 
Consider the map $h \colon \Osc{Q} \rightarrow T_{\ido{}} \Difforb{Q,\uU}, [\hat{\sigma}] \mapsto [t \mapsto E (t[\hat{\sigma}])]$. For $i \in I$, we denote by $\sigma_i$ the canonical lift on $U_i$ of the orbisection $[\hat{\sigma}]$. Then \eqref{eq: sim:vf} together with Remark \ref{rem: gp:ident} (b) and \eqref{eq: nmap} implies:
    \begin{equation}\label{eq: compo}
     \alpha_{\ido{}} \circ h ([\hat{\sigma}]) = ([t \mapsto \exp_{W_{\alpha (i)}} (t \sigma_i (p))])_{i \in I, p \in U_i}.
    \end{equation}
 As $\exp_{W_{\alpha (i)}}$ is the Riemannian exponential map on $W_{\alpha (i)}$, the map $c_{i,p} (t) \coloneq \exp_{W_{\alpha (i)}} (t\sigma_i (p))$ is a geodesic with $c_{i,p}'(0) = \sigma_i (p)$. Therefore \eqref{eq: compo} yields $\alpha_{\ido{}} \circ h ([\hat{\sigma}]) = (\sigma_i )_{i \in I} = \Lambda_{\aA} ([\hat{\sigma}])$. Since  $E$ is a diffeomorphism $h = T_0 E (0,\cdot)$ is an isomorphism of topological vector spaces. Now $\alpha_{\ido{}} \circ h = \Lambda_\aA$ shows that $V_{\ido{}} = \im \alpha_{\ido{}}$ and $\alpha_{\ido{}}$ is an isomorphism of topological vector spaces. In particular, the formula shows that $\alpha_{\ido{}}$ is a linear isomorphism onto the closed subspace $V_{\ido{}} = \im \Lambda_\aA \subseteq \bigoplus_{i \in I} \vect{U_i}$.
\paragraph{Step 6:} \emph{$G^\phi|^{V_{[\hat{\phi}]}}_{V_{\ido{}}}$ is an isomorphism of topological vector spaces and $V_{[\hat{\phi}]} = \im \Lambda_{[\hat{\phi}]}$.} By definition, $G^\phi$ is the map $(g_i^\phi)^*_{i \in I}$ and each $g_i^\phi \colon U_{\phi_i} \rightarrow U_i$ is a diffeomorphism. The map $(g_i^\phi)^* \colon \vect{U_i} \rightarrow C_i^\phi$ is an isomorphism of topological vector spaces by Definition \ref{defn: nspace}. From \cite[Ch.\ II, \S 4, No.\ 5, Proposition 8 (i)]{bourbaki1987}, we deduce that the mapping $G^\phi|_{\bigoplus_{i \in I} \vect{U_i}}^{\bigoplus_{i \in I} C_i^\phi}$ is an isomorphism of topological vector spaces.  
 By Step 5, $V_{\ido{}}$ is a subspace of $\bigoplus_{i \in I} \vect{U_i}$ and $V_{[\hat{\phi}]} = G^\phi (V_{\ido{}})$ holds by Step 3. Since $G^\phi$ maps $\bigoplus_{i \in I} \vect{U_i}$ into $\bigoplus_{i \in I} C_{\phi_i}$, the set $V_{[\hat{\phi}]}$ is contained in $\bigoplus_{i \in I} C_{\phi_i}$. Endow $V_{\ido{}}$ with the subspace topology of $\bigoplus_{i \in I} \vect{U_i}$ and $V_{[\hat{\phi}]}$ with the subspace topology of $\bigoplus_{i \in I} C_i^\phi$. The map $G^\phi|_{V_{\ido{}}}^{V_{[\hat{\phi}]}}$ becomes an isomorphism of topological vector spaces. By construction, for $(f_i)_{i \in I} \in V_{[\hat{\phi}]}$ there is a unique $[\widehat{\sigma^f}] \in \Osc{Q}$ such that $(f_i)_{i \in I} = G^\phi \Lambda_\aA ([\widehat{\sigma^f}]) = (\sigma_i^f \circ g_i^\phi)_{i \in I}$. Hence the elements in $V_{[\hat{\phi}]}$ are of the form $(\sigma_i \circ g_i^\phi)_{i \in I}$, where $\sigma_i$ is the canonical representative on $U_i$ of some $[\hat{\sigma}] \in \Osc{Q}$. As a consequence of the definition of $\Lambda_{[\hat{\phi}]}$, as a set $\im \Lambda_{[\hat{\phi}]}$ and $V_{[\hat{\phi}]}$ coincide. By definition of the topology, they also coincide as topological vector spaces.
\paragraph{Step 7:} \emph{$\alpha_{[\hat{\phi}]}$ is an isomorphism of topological spaces for each $[\hat{\phi}] \in \sS$.}
 Endow $V_{[\hat{\phi}]}$ with the topology as in Step 6 and obtain a commutative diagram for $[\hat{\phi}] \in \sS$:
\begin{displaymath}
	 \begin{xy}
  	\xymatrix{
     	T_{\ido{}} \Difforb{Q,\uU} \ar[rr]^-{\alpha_{\ido{}}} \ar[d]^{T\rho_{[\hat{\phi}]}} && V_{\ido{}} \ar[d]^{G^\phi |_{V_{\ido{}}}^{V_{[\hat{\phi}]}}}\\
	T_{[\hat{\phi}]} \Difforb{Q,\uU} \ar[rr]^-{\alpha_{[\hat{\phi}]}} && V_{[\hat{\phi}]}}			
	\end{xy}
	\end{displaymath}
 As all arrows with the exception of the lower row are isomorphisms of topological vector spaces, so is $\alpha_{[\hat{\phi}]}$. By Step 6, $\im \alpha_{[\hat{\phi}]} = V_{[\hat{\phi}]} = \im \Lambda_{[\hat{\phi}]}$ holds, thus proving the assertion.  
\end{proof}

We are now in the position to obtain regularity properties for the Lie group $\Difforb{Q,\uU}$.  
 \begin{thm}\label{thm: diff:reg} \no{thm: diff:reg}
  Let $(Q,\uU)$ be $\sigma$-compact. Then the Lie group $\Difforb{Q,\uU}$ is $C^k$-regular for each $k \in \NN_0 \cup \set{\infty}$. In particular, this group is regular in the sense of Milnor.
 \end{thm}

\begin{proof}
 We claim that $\Difforb{Q,\uU}$ is a (strongly) $C^0$-regular Lie group. If this is true, then the assertion is a direct consequence of Definition \ref{defn: Ck:reg}. To prove the claim,  by Lemma \ref{lem: znbhd:reg} it suffices to obtain a smooth evolution and right product integrals for some zero-neighborhood $C^0([0,1],U)$. Let $E \colon \hH_\rho \rightarrow \Difforb{Q,\uU}, [\hat{\sigma}] \rightarrow [\expo]\circ [\hat{\sigma}]|^\Omega$ be the manifold chart at the identity introduced in Theorem \ref{thm: diff:lgp} (cf.\ Proposition \ref{prop: diff:chart}). Using the map $\evol$ introduced in Lemma \ref{lem: OSC:pat}, we define a map 
    \begin{displaymath}
     E_1 \coloneq E \circ \evol|_{C^0 ([0,1], \rR)} \colon C^0 ([0,1], \rR) \rightarrow \Difforb{Q,\uU},
    \end{displaymath}
 where $\rR$ is chosen as in Lemma \ref{lem: znhd:cst} with respect to the symmetric subset $\sS \subseteq \im E$. By Lemma \ref{lem: OSC:pat}, $\evol$ is a smooth map, whence $E_1$ is smooth as a composition of smooth maps. Identify $\Osc{Q}$ with $L(\Difforb{Q,\uU})$ via the isomorphism $T_0 E(0,\cdot) = \alpha_{\ido{}}^{-1} \circ \Lambda_\aA$ and recall from step 5 of Lemma \ref{lem: char:tdiff} $V_{\ido{}} = \Osc{Q}$. Following Lemma \ref{lem: znbhd:reg}, the Lie group $\Difforb{Q,\uU}$ will be (strongly) $C^0$-regular if we can show that each $\gamma \in C^0([0,1], \rR)$ has a right product integral $\pP (\gamma)$ with $\pP (\gamma ) (1) = E_1 (\gamma)$.\\
 \mbox{}\\ 
 We first need to understand the derivative of a $C^1$-curve $\eta \colon [0,1] \rightarrow \sS \subseteq \Difforb{Q,\uU}$. For $s\in [0,1]$, we let $[\hat{\sigma}^{\eta (s)}]$ be the preimage $E^{-1}(\eta(s))$ (cf.\ Definition \ref{defn: nspace}). Recall from Lemma \ref{lem: nlifts} that for all $s,t \in [0,1]$, there is a representative $E_{\eta (t)^{-1}} (\eta (s))$ of $\eta (s)$. Using the notation of Definition \ref{defn: nspace}, the lifts of this representative with respect to the atlas $\set{(U_{\eta (t)_i},H_{\alpha (i),U_{\eta (t)_i}}, \varphi_{\alpha (i)}|_{U_{\eta(t)_i}})}_{i \in I}$ are given as 
    \begin{displaymath}
	E_{\eta (t)^{-1}} (\eta (s))_i = \gamma_i^{\eta (t)} .\exp_{W_{\alpha (i)}} \circ \sigma_{\alpha (i)}^{\eta (s)}|_{U_{\eta (t)_i}} .
    \end{displaymath}
 The derivative of the lift with respect to $s$ may be computed locally in manifold-charts. To do so, we fix $p \in U_{\eta (t)_i}$ for some $t \in [0,1]$: Since $U_{\eta (t)_i} \subseteq \Omega_{2,i}$ by Definition \ref{defn: nspace}, we choose and fix a manifold chart $(V_{5,\alpha (i)}^{n_p}, \kappa_{n_p}^{\alpha (i)}) \in \fF_5 (K_{5,i})$ with $p \in V_{2,\alpha (i)}^{n_p}$. Observe that by \cite[Lemma F.6 and Lemma 4.11]{hg2004}, the map $K_{n_p}^{\alpha (i)} \colon \vect{V_{5,\alpha (i)}^{n_p}} \rightarrow C^\infty (B_{5} (0), \RR^d), X \mapsto X_{[n_p]}$ with $X_{[n_p]} = C^\infty ((\kappa_{n_p}^{\alpha (i)})^{-1}, \RR^d) (\theta_{\kappa_n^{\alpha (i)}}(X))$, is an isomorphism of topological vector spaces. As $\eta$ is of class $C^1$, the following composition yields a $C^1$-curve:  
  \begin{displaymath}
   \eta_{t,p,i} \coloneq K_{n_p}^{\alpha (i)} \circ \res^{W_{\alpha (i)}}_{V_{5,\alpha (i)}^{n_p}}\circ \tau_{W_{\alpha (i)}} \circ E^{-1}\circ \eta \colon [0,1] \rightarrow C^\infty (B_5 (0) , \RR^d).
  \end{displaymath}
 Let $\exp_{n_p}$ be the Riemannian exponential map induced on $B_{5} (0)$ by the pullback metric of the Riemannian metric on $W_{\alpha (i)}$ via $(\kappa_{n_p}^{\alpha (i)})^{-1}$. Since $E^{-1} (\sS) \subseteq \hH_\rho$ and $(V_{5,\alpha (i)}^{n_p}, \kappa_{n_p}^{\alpha (i)}) \in \fF_{5,K_{5,i}}$, the construction of $\hH_\rho$ (cf.\ Theorem \ref{thm: diff:lgp}, or more precisely Construction \ref{con: H} and Construction \ref{con: comp:loc}) shows $\eta_{t,p,i} ([0,1]) (\overline{B_3 (0)}) \subseteq B_{\ve_{n_p}} (0) \subseteq B_{\nu_{n_p}} (0)$, whence $\eta_{t,p,i} (s) \in \lfloor \overline{B_{2} (0)}, B_{\nu_{n_p}} (0)\rfloor_\infty \subseteq C^\infty (B_5 (0), \RR^d)$ holds for all $s \in [0,1]$. By choice of $\nu_{n_p}$, the set $B_{4} (0) \times B_{\nu_{n_p}} (0)$ is contained in $\dom \exp_{n_p}$ (cf.\ Lemma \ref{lem: eos:cov}). We deduce from \cite[Proposition 4.23]{hg2004} that 
  \begin{displaymath}
      (\exp_{n_p})_* \colon  \lfloor \overline{B_{2} (0)}, B_{\nu_{n_p}} (0)\rfloor_\infty \rightarrow C^\infty (B_2 (0), \RR^d), f \mapsto \exp_{n_p} (\id_{B_2 (0)} , f|_{B_2 (0)})   
  \end{displaymath}
 is smooth. We obtain a $C^1$-curve $(\exp_{n_p})_* \circ \eta_{t,p,i} \colon [0,1] \rightarrow C^\infty (B_2 (0), \RR^d)$. Furthermore, Lemma \ref{lem: eos:cov} (b) yields $\exp_{W_{\alpha (i)}} \circ T(\kappa^{\alpha (i)}_{n_p})^{-1}|_{B_2 (0) \times B_{\nu_{n_p}} (0)} = (\kappa_{n_p}^{\alpha (i)})^{-1} \circ \exp_{n_p}|_{B_2 (0) \times B_{\nu_{n_p}} (0)}$. The above considerations did not depend on $p \in U_{\eta (t)_i}$, whence they may be repeated for each $p\in U_{\eta (t)_i}, i \in I$. With Lemma \ref{lem: eos:cov} (b) and the Exponential law \cite[Theorem 3.28]{alas2012}, we may now compute the derivative as
    \begin{align}\notag
     \alpha_{\eta (t)} (\eta' (t)) &= \alpha_{\eta (t)} ([s \mapsto \eta (t+s)]_\sim) = ([s \mapsto E_{(\eta (t))^{-1}}(\eta (t+s))_i (p)]_\sim)_{i \in I,p \in U_{\gamma (t)_i}} \\
           &= ([s \mapsto \gamma_i^{\eta (t)} \exp_{W_{\alpha (i)}} \sigma_{\alpha (i)}^{\eta (t+s)} (p)]_\sim)_{i \in I,p \in U_{\gamma (t)_i}} \notag\\
           &= ([s \mapsto \gamma_i^{\eta (t)} \exp_{W_{\alpha (i)}} (T \kappa_{n_p}^{\alpha (i)})^{-1} (\kappa_{n_p}^{\alpha (i)} (p), \eta_{t,p,i} (t+s) (\kappa_{n_p}^{\alpha (i)} (p))])_{i \in I,p \in U_{\gamma (t)_i}} \notag \\
           &= ([s\mapsto \gamma_i^{\eta (t)} (\kappa_{n_p}^{\alpha (i)})^{-1} (\exp_{n_p})_* (\eta_{t,p,i}(t+s)) (\kappa_{n_p}^{\alpha (i)} (p))]_\sim)_{i \in I,p \in U_{\gamma (t)_i}} \notag\\   
	   &= ( T (\gamma_i^{\eta (t)} (\kappa_{n_p}^{\alpha (i)})^{-1}) \left.\frac{\partial}{\partial s}\right|_{s=t} ((\exp_{n_p})_*\circ \eta_{t,p,i})^\wedge (t, \kappa_{n_p}^{\alpha (i)} (p),1))_{i \in I,p \in U_{\gamma (t)_i}} .\label{locderiv}
    \end{align}
 Let $\xi \in C^0 ([0,1], \rR)$ be some continuous curve. By Lemma \ref{lem: znhd:cst}, we may consider the $C^1$-curve $\eta \coloneq E \circ \omega (\xi) \colon [0,1] \rightarrow \sS$. To compute the derivative $\eta'(t)$, we exploit the identity \eqref{locderiv}. The definition of the mappings implies 
    \begin{displaymath}
     \eta_{t,p,i} = K_{n_p}^{\alpha (i)} \circ \res^{W_{\alpha (i)}}_{V_{5,\alpha (i)}^{n_p}}\circ \tau_{W_{\alpha (i)}} \circ E^{-1}  \circ (E \circ \omega (\xi)) = (s \mapsto K_{n_p}^{\alpha (i)} (\omega (\xi)(s)_{\alpha (i)}|_{V_{5,\alpha (i)}^{n_p}})).
    \end{displaymath}
 The canonical lift $\omega (\xi)(s)_{\alpha (i)}$ is uniquely determined, whence $\omega (\xi)(s)_{\alpha (i)}$ coincides with $\omega_i (\xi_{\alpha (i)}(s))$ (cf.\ Lemma \ref{lem: loc:comp}) on $\Omega_{2,K_{5,i}}$ by the proof of Lemma \ref{lem: Osgen}. Since $(V_{5,\alpha (i)}^{n_p}, \kappa^{\alpha (i)}_{n_p}) \in \fF_5 (K_{5,i})$, we derive $V_{2,\alpha (i)}^{n_p} \subseteq \Omega_{2,K_{5,i}}$. Therefore the lift satisfies \eqref{eq: def:vf}. Summing up, for $(s,x) \in [0,1] \times V_{2,\alpha (i)}^{n_p}$:
  \begin{align*}
  \eta_{t,p,i} (s)(x) &= K_{n_p}^{\alpha (i)} (e(\xi)_{\alpha (i)}(s))(\kappa_{n_p}^{\alpha (i)} (x))\\
                     &= \proj_2 \circ T\kappa_{n_p}^{\alpha (i)} (\exp_{W_{\alpha (i)}}|_{N_x})^{-1} \circ (\kappa_{n_p}^{\alpha (i)})^{-1} \circ \Fl_0^f (s, \kappa_{n_p}^{\alpha (i)} (x),\xi_{\alpha (i)[n_p]}) .                                                                                                                                                                          
  \end{align*}
 Observe that $\exp_{n_p} T\kappa_{n_p}^{\alpha (i)} (\exp_{W_{\alpha (i)}}|_{N_{x}})^{-1}  = \exp_{n_p} T\kappa_{n_p}^{\alpha (i)} (\exp_{W_{\alpha (i)}}|_{N_{x}})^{-1}$. By construction of $N_x$ (see Lemma \ref{lem: eos:cov} (b)), we obtain: 
 \begin{displaymath}
  \exp_{n_p} T\kappa_{n_p}^{\alpha (i)} (\exp_{W_{\alpha (i)}}|_{N_{x}})^{-1} = \kappa_{n_p}^{\alpha (i)} \exp_{W_{\alpha (i)}} (\exp_{W_{\alpha (i)}}|_{N_x})^{-1} = \kappa_{n_p}^{\alpha (i)}.
 \end{displaymath}
 Insert this identity and the local formula for $\eta_{t,p,i}$ into \eqref{locderiv}: 
  \begin{align*}
   \alpha_{\eta (t)} (\eta' (t)) &= \left(T(\gamma_i^{\eta(t)} (\kappa_{n_p}^{\alpha (i)})^{-1})  \left.\frac{\partial}{\partial s}\right|_{s=t} ((\exp_{n_p})_* \circ \eta_{t,p,i})^\wedge (s, \kappa_{n_p}^{\alpha (i)} (p))\right)_{i \in I,p \in U_{\eta (t)_i}}\\
		 &= \left(T(\gamma_i^{\eta (t)} (\kappa_{n_p}^{\alpha (i)})^{-1})  \left.\frac{\partial}{\partial s}\right|_{s=t} \Fl_0^{f} (s,\kappa_{n_p}^{\alpha (i)} (p), \xi_{\alpha (i) [n_p]}))\right)_{i \in I,p \in U_{\eta (t)_i}} .
  \end{align*}
 Fixing $\kappa_{n_p}^{\alpha (i)} (p)$ and $\xi$, the flow $\Fl_0^{f} (\cdot,\kappa_{n_p}^{\alpha (i)} (p) , \xi_{\alpha (i) [n_p]})$ is a solution to the differential equation \eqref{eq:diffeq}. Thus
    \begin{align*}
      \left.\frac{\partial}{\partial s}\right|_{s=t} \Fl_0^{f} (s,\kappa_{n_p}^{\alpha (i)} (p) , \xi_{\alpha (i) [n_p]})  &= (\Fl_0^{f} (t,\kappa_{n_p}^{\alpha (i)} (p) , \xi_{\alpha (i) [n_p]}), \xi_{\alpha (i) [n_p]} (t)(\Fl_0^{f} (t,\kappa_{n_p}^{\alpha (i)} (p) , \xi_{\alpha (i) [n_p]})))\\
      &= T\kappa_{n_p}^{\alpha (i)} \xi (t)_{\alpha (i)}  \circ (\kappa_{n_p}^{\alpha (i)})^{-1} (\Fl_0^{f} (t,\kappa_{n_p}^{\alpha (i)} (p) , \xi_{\alpha (i) [n_p]})).
    \end{align*}
 Since  $\xi (t)_{\alpha (i)}$ is a canonical lift, it is equivariant with respect to $H_{\alpha (i)}$. Thus the last identity proves: 
 \begin{align*}
 \alpha_{\eta (t)} (\eta' (t)) 	 &=  (T(\gamma_i^{\eta (t)}) \xi (t)_{\alpha (i)} (\kappa_{n_p}^{\alpha (i)})^{-1} \Fl_0^{f} (t,\kappa_{n_p}^{\alpha (i)} (p), \xi_{\alpha (i) [n_p]}))_{i \in I,p \in U_{\eta (t)_i}}\\
				 &= (\xi (t)_{\alpha (i)} (\gamma_i^{\eta (t)}. (\kappa_{n_p}^{\alpha (i)})^{-1} \Fl_0^{f} (t,\kappa_{n_p}^{\alpha (i)} (p), \xi_{\alpha (i) [n_p]}))_{i \in I,p \in U_{\eta (t)_i}}.
  \end{align*}
 Moreover, $\omega (\xi) (t) = E^{-1} (\eta (t))$ holds by construction. Using the notation of Lemma \ref{lem: nlifts} and its proof, we obtain $\exp_{W_{\alpha (i)}} \circ \omega (\xi) (t)_{\alpha (i)} (p) = E(\eta (t)^{-1}, \eta(t))_i (p)$. On the other hand, \eqref{eq: def:vf} yields the identity
  \begin{displaymath}
   \exp_{W_{\alpha (i)}} \circ \omega (\xi) (t)_{\alpha (i)} (p) = (\kappa_{n_p}^{\alpha (i)})^{-1} \Fl_0^{f} (t,\kappa_{n_p}^{\alpha (i)} (p), \xi_{\alpha (i) [n_p]}).
  \end{displaymath}
 By choice of $\gamma_i^{\eta (t)}$ (see the proof of Lemma \ref{lem: nlifts}), we derive: 
  \begin{displaymath}
    \alpha_{\eta (t)} (\eta' (t)) = (\xi (t)_{\alpha (i)}  (g_i^{\eta (t)} (p))_{i \in I,p \in U_{\eta (t)_i}} =  (\xi (t)_{\alpha (i)} \circ g_i^{\eta (t)})_{i \in I} = \Lambda_{\eta (t)} (\xi (t) \circ \eta (t)).
  \end{displaymath}
 We may now use the structural results on the tangent space of $\Difforb{Q,\uU}$ at $\gamma (t)\in \sS$. To shorten the notation, abbreviate $\Psi \coloneq T_0 E(0,\cdot) = \alpha_{\ido{}}^{-1} \circ \Lambda_\aA$. From Lemma \ref{lem: char:tdiff} and its proof (in particular, the formula in Step 3), we infer 
  \begin{displaymath}
   \Lambda_{\eta (t)}^{-1} (\alpha_{\eta (t)} (\eta' (t))) = \xi (t) \circ \eta (t) = \Lambda_{\eta (t)}^{-1} (G^{\eta (t)} \Lambda_\aA (\xi (t))) =  \Lambda_{\eta (t)}^{-1} (\alpha_{\eta (t)} (T \rho_{\eta (t)} \Psi (\xi (t))).
  \end{displaymath}
 The map $\Lambda_{\eta (t)}^{-1} \circ \alpha_{\eta (t)}$ is an isomorphism of topological vector spaces, whence $\eta' (t) = T\rho_{\eta (t)} \Psi (\xi (t))$ follows. Recalling the definition of $\eta$ we have $\eta' (t) = \frac{d}{d t} E (\omega (\xi) (t)) = T\rho_{E (\omega (\xi)(t))} \Psi (\xi (t))$.\\
 The facts obtained so far allow the right logarithmic derivative of $\eta (t) = E(\omega (\eta) (t))$ to be computed: 
  \begin{equation}\label{eq: rlog}
   \delta^r (\eta) (t) = T\rho_{E(\omega (\xi) (t))^{-1}}  \frac{d}{d t}E(\omega (\xi) (t)) = T\rho_{E(\omega (\xi) (t))^{-1}} T\rho_{E(\omega (\xi) (t))} \Psi (\xi (t)) = \Psi (\xi (t)).
  \end{equation}
 By construction, $E_1 (\xi) = E(\omega (\xi) (1)) = \eta (1)$ and Lemma \ref{lem: OSC:pat} implies $\omega (\xi) (0) = \zs$. Thus $\eta (0) = E(\omega (\xi) (0)) = E(\zs) =\ido{}$ holds. Furthermore, the computation of the right logarithmic derivative \eqref{eq: rlog} shows that the curve $\xi$ possesses a right product integral $E(\omega (\xi)) =\eta$. We have already seen that the mapping $E_1$ is smooth, thus the proof is complete and $\Difforb{Q,\uU}$ is a (strongly) $C^0$-regular Lie group. 
\end{proof}

 The orbifolds in the present paper are not assumed to be second countable. We had to require second countability of the orbifold to assure that $\Osc{Q}$ is countably patched. In this case, we obtain an atlas indexed by the countable set $I$, whence the map 
  \begin{displaymath}
   \Lambda \colon \bigoplus_{i \in I} C^r([0,1], \vect{U_i}) \rightarrow C^r\left([0,1], \bigoplus_{i \in I} \vect{U_i}\right), (f_i) \mapsto \sum_{i \in I} (\iota_i)_* (f_i)
  \end{displaymath}
 is an isomorphism of topological vector spaces for $r \in \NN_0$ if the mapping spaces are endowed with the compact open $C^r$-topology (see Lemma \ref{lem: pb:patch}). This fact was crucial to prove the smoothness of the evolution map $\evol$. It is known that $\Lambda$ fails to be an isomorphism of locally convex spaces if $I$ is uncountable. We give a proof for this fact:\\[1em]
 Fix $r= 0$ and let $I$ be an uncountable set. Notice that arguments as in the proof of Lemma \ref{lem: pb:patch} assure that the map $\Lambda$ is an isomorphism of vector spaces which is continuous. We denote its inverse by $\Theta$  (see Lemma \ref{lem: pb:patch} for the construction). Hence we have to prove that $\Theta$ is discontinuous if $I$ is uncountable.\\
 For each $i \in I$, we choose and fix a one-dimensional subspace $E_i \subseteq \vect{U_i}$. The locally convex direct sum $\bigoplus_{i \in I} \RR \cong \bigoplus_{i\in I} E_i$ may be identified in a canonical way with a subspace of $\bigoplus_{i\in I} \vect{U_i}$ by \cite[Ch.\ II, \S 4 No.\ 5, Proposition 8]{bourbaki1987}. If we consider the subspaces $C([0,1],E_i) \subseteq C([0,1],\vect{U_i})$ for $i\in I$, we may analogously identify $\bigoplus_{i \in I} C([0,1], \RR) \cong \bigoplus_{i \in I} C ([0,1],E_i)$ with a subspace of $\bigoplus_{i \in I} C([0,1],\vect{U_i})$. A trivial computation yields the identity  
  \begin{displaymath}
   \Lambda \left(\bigoplus_{i \in I} C([0,1],E_i)\right) = C\left([0,1], \bigoplus_{i\in I} E_i\right).
  \end{displaymath}
 Hence the inverse $\Theta$ restricts to a map $T \coloneq \Theta|_{ C([0,1], \bigoplus_{i\in I} E_i)}^{\bigoplus_{i \in I} C([0,1],E_i)}$. We claim that $T$ is discontinuous, whence $\Theta$ must be discontinuous. To prove this claim, identify each of the spaces $E_i$ with $\RR$. The assertion then follows from the next lemma, whose proof was communicated to the author by D.\ Vogt and S.A.\ Wegner: 
 \begin{lem}
  The map $T \colon C([0,1], \bigoplus_{i \in I} \RR ) \rightarrow \bigoplus_{i \in I} C([0,1],\RR)$ is discontinuous for each uncountable set $I$.
 \end{lem}
 \begin{proof}
  Recall from \cite[\S 24]{meisevogt1997} and Remark \ref{rem: Cr:Norm} that the compact-open topology on the space $C([0,1],\bigoplus_{i\in I} \RR)$ is induced by the following system of seminorms: 
    \begin{displaymath}
     p_\delta (f) \coloneq \sup_{t \in [0,1]} \sum_{i \in I} \delta_i |(f(t))_i|, \text{ with } \delta = (\delta_i)_{i \in I} \text{ and } \delta_i > 0 \text{ for } i \in I .  
    \end{displaymath}
 Analogously, the topology on $\bigoplus_{i \in I} C([0,1],\RR)$ is induced by the following system of seminorms: 
    \begin{displaymath}
     q_\ve ((f_i)_{i \in I}) \coloneq \sum_{i \in I} \ve_i \sup_{t \in [0,1]} |f_i(t)|, \text{ with } \ve = (\ve_i)_{i \in I} \text{ and } \ve_i > 0 \text{ for } i \in I   .
    \end{displaymath}
 Arguing indirectly, we suppose that $T$ is a continuous map. Since $T$ is linear, it is continuous if and only if 
     \begin{align*}
      \forall \ve = (\ve_i)_{i \in I}\ &\exists \delta = (\delta_i)_{i \in I} , C \geq 0 \ \forall (f_i)_{i \in I} \in \bigoplus_{i \in I} C([0,1],\RR) ,\\  
                         & q_\ve ((f_i)_{i \in I}) \leq C p_\delta \left(\sum_{i \in I} (\iota_i)_* f_i\right) \\
               \Longleftrightarrow \forall \ve = (\ve_i)_{i \in I}\ &\exists \delta = (\delta_i)_{i \in I} \ \forall (f_i)_{i \in I} \in \bigoplus_{i \in I} C([0,1],\RR) ,\\
       & \sum_{i \in I} \ve_i \sup_{t\in [0,1]} |f_i (t)| \leq \sup_{t\in [0,1]} \sum_{i \in I} \delta_i |f_i (t)|.
     \end{align*}
 To obtain a contradiction, fix $\ve = (1)_{i\in I}$ and choose $\delta = (\delta_i)_{i \in I}$ as above. For $n \in \NN$, define the set $M_n \coloneq \setm{i \in I}{\delta_i \leq n}$. By construction, $I = \bigcup_{n \in \NN} M_n$ holds. Since $I$ is uncountable, there must be $N \in \NN$ with $|M_N|=\infty$. \\
 For $n \in \NN$, consider $E \subseteq M_N$ with $E = \set{i_1, \ldots ,i_n}$ and choose $f_{i_k} \in C([0,1],\RR)$ with $0 \leq f_{i_k} \leq 1$ such that $\supp f_{i_k} \cap \supp f_{i_j} = \emptyset$ if $k \neq j$. Furthermore, let there be $t_k \in [0,1]$ with $f_{i_k} (t_k) = 1$ for $1\leq k \leq n$. Define $(f_i)_{i \in I} \in \bigoplus_{i \in I} C([0,1],\RR)$ via $f_i \coloneq f_{i_k}$ if $i = i_k$ for $1\leq k \leq n$ and $f_i \coloneq 0$ otherwise. By choice of $\delta$, 
  \begin{equation}\label{eq: briso}
   \sum_{i \in I} \sup_{t \in [0,1]} |f_i(t)| \leq \sup_{t \in [0,1]} \sum_{i \in I} \delta_i |f_i (t)|.
  \end{equation}
 Compute both sides of the above inequality. For the left hand side of \eqref{eq: briso} the definition of the family $(f_i)_{i \in I}$ yields: 
  \begin{displaymath}
   \sum_{i \in I} \sup_{t \in [0,1]} |f_i (t)| = \sum_{1\leq k \leq n} \sup_{t \in [0,1]} |f_{i_k} (t)| = \sum_{1\leq k \leq n} 1 = n.
  \end{displaymath}
 On the other hand, since the supports of the maps $f_{i_k}$ are disjoint, the right hand side of \eqref{eq: briso}, evaluates as: 
   \begin{displaymath}
    \sup_{t \in [0,1]} \sum_{i \in I} \delta_i |f_i (t)| = \sup_{t \in [0,1]} \sum_{k=1}^n \delta_{i_k} |f_{i_k} (t)| = \sup_{1\leq k \leq n} \delta_{i_k} \leq \sup_{i \in M_N} \delta_i \leq N.
   \end{displaymath}
 Hence \eqref{eq: briso} yields $n \leq N$, where $N$ is fixed but $n$ may be chosen arbitrarily large. We derive a contradiction, whence $T$ may not be continuous. 
 \end{proof}

 Summing up, the inverse $\Theta$ of $\Lambda$ is discontinuous for uncountable index sets $I$. Hence $\Lambda$ fails to be an isomorphism of topological vector spaces if $I$ is uncountable. Thus our methods do not generalize to the setting of arbitrary paracompact orbifolds. As already stated in the introduction, with the methods in \cite{hg2015} one can prove that the diffeomorphism group of a non-second countable orbifold is $C^1$-regular. 
\thispagestyle{empty}
\section{Application to Equivariant Diffeomorphism Groups}
\setcounter{subsubsection}{0}
\setcounter{equation}{0}

In this section, we consider good orbifolds with an orbifold atlas of a single chart. If $(Q,\uU)$ is such an orbifold, we let $\set{(U,G,\pi)}$ be an orbifold atlas for $(Q,\uU)$ and call $(U,G,\pi)$ a \ind{orbifold chart!global chart}{global chart}. It turns out that for certain orbifolds with global chart, the group $\Difforb{Q,\uU}$ induces a Lie group structure on a subgroup of the diffeomorphism group of $U$. We begin our inquiry with several observations: 

\begin{setup}\label{setup: observ}\no{setup: observ}
 Let $(Q,\uU)$ be an orbifold with global chart $(U,G,\pi)$. Consider a diffeomorphism of $U$ which is a weak equivalence, i.e.\ a diffeomorphism $\tilde{h} \colon U \rightarrow U$ together with a group automorphism $\alpha \colon G \rightarrow G$ such that $\tilde{h} \circ g = \alpha (g) \circ \tilde{h}$ holds for all $g \in G$. Note that $\tilde{h}^{-1}$ is also a weak equivalence, with respect to the group automorphism $\alpha^{-1}$. In particular, $\tilde{h}$ and $\tilde{h}^{-1}$ induce mutually inverse continuous maps $h \colon Q \rightarrow Q$ and $h^{-1} \colon Q \rightarrow Q$, respectively. The pair $(h,\tilde{h})$ induces a representative of an orbifold map such that the corresponding orbifold map is a diffeomorphism of orbifolds by Proposition \ref{prop: ld:qpgp} and Proposition \ref{prop: ofd:iso}. Therefore, each diffeomorphism of the global chart which is a weak equivalence canonically induces a unique diffeomorphism of $(Q,\uU)$.\\[1em] Denote by $[\hat{h}]$ the diffeomorphism of orbifolds associated to  $\tilde{h} \in \Diff^G (U)$ by the above construction. We consider the map 
  \begin{displaymath}
   D \colon \textstyle \Diff^G \displaystyle (U) \rightarrow \Difforb{Q,\uU} , \tilde{f} \mapsto [\hat{f}] .
  \end{displaymath}
Each orbifold diffeomorphisms in the image of $D$ is induced by a lift in the global chart, i.e.\ by an element of $\Diff^G (M)$. Since orbifold diffeomorphisms are uniquely determined by their lifts (cf. Corollary \ref{cor: ofdiff:un}), the composition of the lifts in the global chart induces the composition of orbifold diffeomorphisms. The same argument shows that the image  $D(\tilde{h}^{-1})$ coincides with $D (\tilde{h})^{-1}$ (the inverse in $\Difforb{Q,\uU}$) by Corollary \ref{cor: diff:char}. Summing up, $D$ is a group homomorphism.\\
The map $D$ is not injective, as elements of $\Diff^G(U)$ which differ only by composition with an element of $G$ are mapped to the same diffeomorphism of orbifolds. From \cite[Lemma 2.11]{follie2003}, we deduce that the kernel of $D$ coincides with $G$. Hence $D$ induces an injective group homomorphism $\Delta$: 
  \begin{displaymath}
   \begin{xy}
  \xymatrix{
      \one \ar[r] & G\, \ar@{>->}[r] 	& \Diff^G (U) \ar@{->>}[r] \ar[rd]^D 	& \Diff^G (U)/G \ar[r] \ar@{>->}[d]^\Delta & \one  \\
		  & 		&				& \Difforb{Q,\uU} &	   
  }
\end{xy}
  \end{displaymath}
We now ask, whether all orbifold diffeomorphisms of $(Q,\uU)$ arise as quotients of elements on $\Diff^G(U)$. It will turn out that this is the case for certain orbifolds with a global chart, i.e.\ we prove that $\Delta$ is an isomorphism of groups in some cases. In this situation wndow $\Diff^G (U) /G$ via $\Delta$ with the unique Lie group structure turning the mapping into an isomorphism of Lie groups.
\end{setup}

\begin{setup} \label{setup: new}
We can also obtain a Lie group structure on a subgroup of $\Diff^G(U)$. Consider the subgroup of $\Diff^G (U)$ whose elements coincide with the identity off some compact subset:\glsadd{DiffGc}
  \begin{displaymath}
  \textstyle \Diff^G_c \displaystyle (U) \coloneq \setm{f \in \textstyle \Diff^G \displaystyle (U)}{\exists K \subseteq U \text{ compact}, f|_{U \setminus K} =\id_{U \setminus K}}.
  \end{displaymath}
By construction, $\Diff^G_c (U)$ is a subgroup of $\Diff^G (U)$. Then $D$ maps $\Diff^G_c (U)$ into the open Lie subgroup $\Diffc{Q,\uU}$ of $\Difforb{Q,\uU}$. If the intersection $G \cap \Diff^G_c (U)$ contains only $\id_{U}$, the mapping $D$ restricts to an injective group homomorphism $\Delta_c \colon \Diff^G_c (U) \rightarrow \Diffc{Q,\uU}$. By Newman`s theorem $G \cap \Diff_c (U) = \set{\id_U}$ holds, whenever $U$ is non-compact. Since $D$ is surjective for certain orbifolds, $\Delta_c$ becomes an isomorphism of groups for these orbifolds. Thus $\Diff^G_c (U)$ may be endowed with a Lie group structure induced by $\Diffc{Q,\uU}$. The construction principle outlined in Proposition \ref{prop: Lgp:locd} then allows the construction of a unique Lie group structure on $\Diff^G (U)$ which contains $\Diff^G_c (U)$ as an open subgroup.
\end{setup}

\begin{setup}\label{setup: good:ofd}\no{setup: good:ofd}
We introduce the class of orbifolds with global chart considered throughout this section: Let $d$ be in $\NN$ and $G$ be a finite subgroup of the orthogonal group $O(d)\subseteq \Diff (\RR^d)$ such that: \\[1em]
\mbox{}\quad \textbf{(IS)} The group $G$ satisfies $G_x = \set{\id_{\RR^d}}$ for all $x \in \RR^d \setminus \set{0}$.
\\[1em]
Recall that for odd dimension $d$, each element $g$ of $O(d)$ possesses at least one real eigenvalue $\lambda_g$. By orthogonality we must have $\lambda_g \in \set{-1,1}$. If an element $g \in G \setminus \set{\id_{\RR^d}}$, condition (IS) implies $\lambda_g = -1$. Then $g^2$ is an element of $G$ with real eigenvalue $1$. Again condition (IS) forces $g^2 = \id_{\RR^n}$ and thus all eigenvalues of $g$ must be $1$ or $-1$. Using condition (IS), all eigenvalues of $g$ are thus $-1$ and we obtain $g = -\id_{\RR^d}$. Hence for odd $d$ only $G= \set{\id_{\RR^d} , -\id_{\RR^d}}$ or $G = \set{\id_{\RR^d}}$ are possible. We are interested in the non trivial case, whence we assume for the rest of this section that $G \neq \set{\id_{\RR^d}}$ holds. We record the following observations  
  \begin{compactenum}
   \item If $d$ is odd, the group $G$ is generated by $-\id_{\RR^d}$. For $d=1$ this is a reflection, which will be denoted as $r \colon \RR \rightarrow \RR , x \mapsto -x$.
   \item If $d= 2$, the group $G$ may not contain reflections by condition (IS). In this case $G$ contains at least one (non-trivial) rotations of $\RR^2$ which fixes the origin.
  \end{compactenum}
 Let $\pi \colon \RR^d \rightarrow \RR^d /G$ be the quotient map onto the orbit space and $Q \coloneq \RR^d /G$. Then $\set{(\RR^d , G , \pi)}$ is an atlas for $Q$, turning the orbit space into a good orbifold with a global chart. We identify for $d \in \set{1,2}$ the orbit spaces with $[0,\infty[$ respectively a cone.
Each finite subgroup of $O(2)$ -- which is not a dihedral group -- is cyclic by \cite[Ch.\ 5, Theorem 3.4]{artin1991}. Hence the illustration above exhibits the general case for $d=2$.
 \end{setup}

\begin{prop}\label{prop: homeo:lift}
 Let $(Q,\uU)$ be an orbifold as in \ref{setup: good:ofd}. Consider $[\hat{h}] \in \Difforb{Q,\uU}$ with representative $(h, \set{h_i}_{i \in I}, [P,\nu]) \in [\hat{h}]$. The map $h$ lifts to a weak equivalence $\tilde{h} \colon \RR^d \rightarrow \RR^d$ with respect to the $G$-action.
\end{prop}

\begin{proof}
 Consider $[\hat{h}] \in \Difforb{Q,\uU}$ with representative $(h, \set{h_i}_{i \in I}, (P,\nu)) \in [\hat{h}]$. To construct a weak equivalence as required, we shall construct at first a lift on the set of non-singular points.\\
 For the orbifolds defined in \ref{setup: good:ofd}, there is only one singular point. The origin in $\RR^d$ is jointly fixed by all elements of $G$. Hence $\RR^d \setminus \set{0}$ corresponds to the set of non-singular points and we set $Q_{\text{reg}}\coloneq Q\setminus \set{0}$. Recall that the global chart $\pi \colon \RR^d \rightarrow Q$ is a branched covering in the sense of \cite[Section 10]{bcf2005}. Hence $q \coloneq \pi|_{\RR^d \setminus \set{0}}^{Q_{\text{reg}}}$ is a covering by \cite[Theorem 10.3]{bcf2005}.\\ Diffeomorphisms of orbifolds preserve singular points by Proposition \ref{prop: ld:flgp} and thus the homeomorphism $h \colon Q \rightarrow Q$ satisfies $f\pi (0) = \pi (0)$. The restriction $h|_{Q_{\text{reg}}}^{Q_{\text{reg}}}$ yields a homeomorphism. To construct a lift of $h$, we construct at first a lift on $\RR^d \setminus \set{0}$: \\[1em]
 \textbf{If $d=1$ holds}, then the space $\RR \setminus \set{0}$ is disconnected. Then the mapping $q|_{]0,\infty[} \colon ]0,\infty[ \rightarrow Q_{\text{reg}}$ is a homeomorphism and we obtain a well-defined homeomorphism $h^+ \coloneq (q|_{]0,\infty[})^{-1} h q|_{]0,\infty[}$, which maps $]0,\infty[$ to itself. This mapping extends to a homeomorphism via
  \begin{displaymath}
   h_{\text{reg}} \colon \RR \setminus \set{0} \rightarrow \RR \setminus \set{0}, x \mapsto \begin{cases}
												h^+ (x) & x > 0\\
												r \circ  h^+ \circ r (x) = -h^+(-x) & x < 0      .                                                                                    \end{cases}     
  \end{displaymath}
 By construction, $h_{\text{reg}}$ and also its inverse are equivariant maps with respect to $G = \langle r\rangle$. We deduce from \cite[II. Lemma 7.2]{bredon1972} that $h_{\text{reg}}$ extends to a continuous map $\tilde{h} \colon \RR \rightarrow \RR$ by $\tilde{h} (0) = 0$. Similarly we extend the inverse of $h_{\text{reg}}$, whence $\tilde{h}$ is an equivariant homeomorphism, i.e.\ an equivalence.\\[1em]
 \textbf{If $d \geq 2$ holds}, then the space $\RR^d \setminus \set{0}$ is (path-)connected. We have to construct a lift $f_{\text{reg}}$: 
  \begin{displaymath}
   \begin{xy}
  \xymatrix{
												      & &  &\RR^d \setminus \set{0} \ar[d]^q  \\
    \RR^d \setminus \set{0}	\ar@{.>}[rrru]^{f_{\text{reg}}} \ar[rrr]^{h|^{Q_{\text{reg}}}  \circ q}& &  & Q_{\text{reg}}  }
\end{xy}
  \end{displaymath}
 \textbf{For $d\geq 3$}, the space $\RR^d \setminus \set{0}$ is simply connected, path-connected and locally path-connected. Choose $x_0 \in \RR^d \setminus \set{0}$ and $y_0 \in q^{-1} h q (x_0)$. Then by \cite[Proposition 1.33]{hatcher2002}, there is a unique lift $h_{\text{reg}} \colon \RR^d \setminus \set{0} \rightarrow \RR^d \setminus \set{0}$ of $h|^{Q_{\text{reg}}} \circ q$ which maps $x_0$ to $y_0$.\\[1em]
 \textbf{For $d=2$}, the space $\RR^2 \setminus \set{0}$ is \emph{not} simply connected. However, it is path-connected and locally path-connected. We may still apply \cite[Proposition 1.33]{hatcher2002} if the fundamental group $\pi_1 (\RR^2 \setminus \set{0}, x_0)$ satisfies: 
  \begin{equation} \label{eq: homotopy}
   (h|^{Q_{\text{reg}}} \circ q)_* (\pi_1 (\RR^2 \setminus \set{0}, x_0)) \subseteq q_* (\pi_1 (\RR^2 \setminus \set{0}, y_0))
  \end{equation}
 Recall that the fundamental group $\pi_1 (\RR^2 \setminus \set{0}, x_0)$ can be identified with $\ZZ$ (cf.\ \cite[Example 1.15]{hatcher2002}). Moreover, as $G \subseteq \text{SO}(2)$ holds, the subgroup $G\subseteq O(2)$ must be a cyclic group, generated by a rotation $\gamma$ of order $m \in \NN$. As we have already seen, $Q$ is homeomorphic to a cone and $Q_{\text{reg}}$ may be identified with a cone whose tip has been removed. In particular, the space $Q_{\text{reg}}$ is homeomorphic to $\RR^2 \setminus \set{0}$ (see Example \ref{ex: sphere} for further details on this homeomorphism).\\
 Consider a generator $[e]$ of the fundamental group $\pi_1 (\RR^2 \setminus \set{0} , x_0)$, where $e$ is chosen as a circle around the origin passing through $x_0$. If $\gamma$ is a rotation of order $m$, then we have $q_*[e] = [q\circ e]$ is a loop in $Q_{\text{reg}}$, which passes $m$ times through $\pi (y_0)$.\\ Note that $\pi_1 (Q_{\text{reg}}, q (y_0))$ is isomorphic to $\ZZ$ and let $[f]$ be the generator of $\pi_1 (Q_{\text{reg}}, q (x_0))$. By abuse of notation we let $[f]$ be the generator of each fundamental group for points in $Q_{\text{reg}}$. From the arguments above, we deduce $q_* (\pi_1 (\RR^2 \setminus \set{0}, x_0)) = \langle m [f]\rangle$ and thus  
      \begin{displaymath}
       (h|^{Q_{\text{reg}}} \circ q)_* ([e])= (h|_{Q_{\text{reg}}}^{Q_{\text{reg}}})_* (m [f]) = m  ([h\circ f]) \in \langle m [f]\rangle = \im q_*      .                                                                                                                                                                                                                                                                                                                                                                              \end{displaymath}
 Therefore property \eqref{eq: homotopy} is satisfied and we obtain a unique lift $h_{\text{reg}} \colon \RR^2 \setminus \set{0} \rightarrow \RR^2 \setminus \set{0}$ of $h|_{Q_{\text{reg}}}^{Q_{\text{reg}}}$ mapping $x_0$ to $y_0$.
     
Analogous arguments allow the construction of a unique lift $(h^{-1})_{\text{reg}}$ for $h^{-1}|^{Q \setminus \set{0}} \circ q$ and $d \geq 2$, which maps $y_0$ to $x_0$.
We claim that $(h^{-1})_{\text{reg}}$ is the inverse of $h_{\text{reg}}$. If this is true, then $h_{\text{reg}}$ is a homeomorphism. To prove the claim, consider the map $f \coloneq h_{\text{reg}} \circ (h^{-1})_{\text{reg}}$ and compute
  \begin{displaymath}
   q \circ f = q \circ h_{\text{reg}} \circ (h^{-1})_{\text{reg}} = h\circ q \circ (h^{-1})_{\text{reg}} = q.
  \end{displaymath}
 Hence $f$ is a lift of $\id_{Q_{\text{reg}}}$ taking $y_0$ to $y_0$, and so is the map $\id_{\RR^d \setminus \set{0}}$. By the uniqueness of lifts between pointed spaces (see \cite[Proposition 1.34]{hatcher2002}), $h_{\text{reg}} \circ (h^{-1})_{\text{reg}} = f = \id_{\RR^d \setminus \set{0}}$. Likewise, $(h^{-1})_{\text{reg}} \circ h_{\text{reg}} = \id_{\RR^d \setminus \set{0}}$. Summing up, $h_{\text{reg}}$ is a homeomorphism with inverse $(h^{-1})_{\text{reg}}$.
  
 We now show that the homeomorphism $h_{\text{reg}}$ is a weak equivalence. To this end, let $g$ be in $G$ and $x$ in $\RR^d \setminus \set{0}$. We have 
    \begin{displaymath}
     q \circ h_{\text{reg}} \circ g  \circ h_{\text{reg}}^{-1} (x) = h h^{-1} q (x) = q(x) .
    \end{displaymath}
 Hence $h_{\text{reg}} \circ g \circ h^{-1}_{\text{reg}}$ is a lift of $\id_{\RR^d \setminus \set{0}}$ and so there is an unique element $\alpha (g) \in G$ such that $h_{\text{reg}} \circ g  \circ h_{\text{reg}}^{-1} (x_0)= \alpha (g)(x_0)$. By uniqueness of lifts, $h_{\text{reg}} \circ g  \circ h_{\text{reg}}^{-1}= \alpha (g)$ on $\RR^d \setminus \set{0}$. Repeat this construction for each $g \in G$ to obtain a map $\alpha \colon G \rightarrow G$ with $h_{\text{reg}} \circ g = \alpha (g) \circ h_{\text{reg}}$ on $\RR^d \setminus \set{0}$ for each $g \in G$. Since $\alpha (gk) \circ h_{\text{reg}} = h_{\text{reg}} \circ (g k) = \alpha (g) . h_{\text{reg}} \circ k = \alpha (g) . \alpha (k) . h_{\text{reg}}$ holds and $h_{\text{reg}}$ is a homeomorphism, the map $\alpha$ is an injective group homomorphism. As $G$ is finite, $\alpha$ is thus a group automorphism and $h_{\text{reg}}$ is a weak equivalence. \\
 We extend the weak equivalence $h_{\text{reg}}$ to a map $\tilde{h} \colon \RR^d \rightarrow \RR^d$ by defining $\tilde{h} (0) = 0$. This map is clearly bijective, equivariant with respect to $\alpha$ and lifts $h$. An analogous argument\footnote{The proof works exactly the same if we replace equivariant mappings with weak equivalences.} as in the proof of \cite[II.\ Lemma 7.2]{bredon1972} shows that this map and its inverse are both continuous. Hence $\tilde{h}$ is the desired weak equivalence of $\RR^d$ with respect to the $G$-action.   
\end{proof}

\begin{prop}\label{prop: sm:globlift}
 For an orbifold $(Q,\uU)$ as in \ref{setup: good:ofd}, the mapping $D$ introduced in \ref{setup: observ} is surjective. In particular, the induced map $\Delta \colon \Diff^G (\RR^d) /G \rightarrow \Difforb{Q,\uU}$ is a group isomorphism.
\end{prop}

\begin{proof}
 Consider an arbitrary orbifold diffeomorphism $[\hat{h}] \in \Difforb{Q,\uU}$ and fix a representative $\hat{h} = (h, \set{h_i}_{i \in I}, [P,\nu])$ of $[\hat{h}]$ with the following properties: Each $h_i \colon V_i \rightarrow W_i, i\in I$ is a diffeomorphism such that there are embeddings of orbifold charts $\lambda_i \colon V_i \rightarrow \RR^d$ and $\mu_i \colon W_i \rightarrow \RR^d$ into the global chart from above. A representative with these properties exists by compatibility of orbifold charts and a combination of Corollary \ref{cor: ll:diff} and Corollary \ref{cor: diff:char}. We have to prove that $[\hat{h}]$ is contained in the image of $D$.\\[1em]
 Construct a lift of the homeomorphism $h$ in $\Diff^G (\RR^d)$: Let $\tilde{h} \colon \RR^d \rightarrow \RR^d$ be the lift of $h$ constructed in Proposition \ref{prop: homeo:lift}. The lift $\tilde{h}$ is a weak equivalence and we denote by $\alpha \colon G \rightarrow G$ the associated group automorphism. We claim that $\tilde{h}$ is a smooth map with smooth inverse. If this is true, then $\tilde{h} \in \Diff^G (\RR^d)$ is a smooth lift of $h$ which is compatible with the family of lifts $\set{h_i}_{i \in I}$. Hence, Corollary \ref{cor: diff:char} implies $\Delta (\tilde{h}G) = [\hat{h}] \in \im D$. \\
 To prove the claim, recall that $\set{h_i}_{i \in I}$ is a family of smooth lifts for $h$. For each $i \in I$ the following diagram commutes with the exception of the outer square:
  \begin{displaymath}
   \begin{xy}
  \xymatrix{ V_i \ar[dd]^{\lambda_i} \ar[rrr]^{h_i} \ar[rd] & & & W_i \ar[ld] \ar[dd]_{\mu_i}	\\
		& Q \ar[r]^h & Q				\\
             \RR^d \ar[rrr]^{\tilde{h}} \ar[ru]^\pi & & & \RR^d \ar[lu]_\pi}
\end{xy}
  \end{displaymath} 
Notice that $h_i^{-1}$ is a lift for $h^{-1}$, whence $\pi \tilde{h} \lambda_i h_i^{-1} \mu_i^{-1} = \pi|_{\im \mu_i}$ holds. Consider the set of non-singular points in the image of $\mu_i$ which we abbreviate as $\text{Reg}_{\mu_i} \coloneq \im \mu_i \cap (\RR^d \setminus \set{0})$. Since each point in $\text{Reg}_{\mu_i}$ is non-singular, by the above we obtain a disjoint union $\text{Reg}_{\mu_i} = \sqcup P_g^i$ with $P_g^i \coloneq \setm{x \in \text{Reg}_{\mu_i}}{\tilde{h} \lambda_i h_i^{-1} \mu_i^{-1} (x) = g(x)}$. As $G$ is a finite group, the sets $P_g^i$ are open and closed. \\[1em]
\textbf{Case 1, $d \geq 2$}: Then $\text{Reg}_{\mu_i}$ is a connected set and we obtain $\text{Reg}_{\mu_i} = P_g^i$ for some $g \in G$. By continuity, we must have 
  \begin{displaymath}
   \tilde{h}|_{\im \lambda_i} = g \circ\mu_i \circ h_i \circ \lambda_i^{-1}.
  \end{displaymath}
 In other words, $\tilde{h}|_{\im \lambda_i}$ is a smooth mapping. Furthermore, for $x = g.y \in g.\im \lambda_i$ we have $\tilde{h} (x) = \alpha (g). \tilde{h} (y)$. Thus $\tilde{h}$ is smooth on all of $G.(\im \lambda_i)$ since each element of $G$ is smooth and $\tilde{h}$ is smooth on $\im \lambda_i$. The charts associated to the lifts $h_i \colon V_i \rightarrow W_i$ cover $Q$, whence $\RR^d = \bigcup_{i \in I} G. (\im \lambda_i)$ holds. Then $\tilde{h}$ is smooth and an analogous argument yields the same for $\tilde{h}^{-1}$. Summing up, for $d\geq 2$ the lift $\tilde{h}$ is contained in $\Diff^G (\RR^d)$.\\[1em]
 \textbf{Case 2, $d=1$}: The set $\text{Reg}_{\mu_i}$ is disconnected if $0 \in \im \mu_i$. However, analogous arguments as in the case $d \geq 2$ assure that $\tilde{h}$ is smooth on $\RR \setminus \set{0}$. It suffices to prove that $\tilde{h}$ is also smooth in $0$. To this end, consider $i \in I$ with $0 \in \im \mu_i$ (observe that this implies $0 \in \im  \lambda_i$ since $h_i$ is a lift of $h$, which fixes $\pi (0)$). Then $\lambda_i$ and $\mu_i$ are embeddings of orbifold charts whose images contain the point fixed by the reflection $r$ and are $G =\langle r\rangle$-stable sets. To shorten the notation, define $h_0 \coloneq \mu_i \circ h_i \circ \lambda_i^{-1}$. By construction, $\pi \circ h_0 \circ r \circ h_0^{-1} = \pi|_{\im \lambda_i}$ holds. Thus \cite[Lemma 2.11]{follie2003} yields a unique $\gamma \in G = \langle r \rangle$ with $h_0 \circ r = \gamma \circ h_0$. Since $h_0$ is a bijective map, we must have $\gamma =r$ and thus $h_0$ is equivariant. We claim that $\tilde{h}|_{\im \lambda_i} = g \circ h_0$ holds for a uniquely determined element $g \in \langle r\rangle$. If this is true, then $\tilde{h}|_{\im \lambda_i}$ is smooth. An analogous argument shows that $\tilde{h}^{-1}|_{\im \mu_i}$ is smooth. Summing up, both $\tilde{h}$ and $\tilde{h}^{-1}$ are smooth, whence $\tilde{h}$ is contained in $\Diff^G (\RR)$.\\[1em]
 \textbf{Proof of the claim:} Arguing as in the case $d\geq 2$, there are elements $\gamma^+ ,\gamma^- \in G$ such that the following is satisfied: 
 \begin{displaymath}
  \tilde{h}|_{\im \lambda_i \cap [0,\infty[} = \gamma^+ \circ h_0|_{\im \lambda_i \cap [0,\infty[}  \text{ and } \tilde{h}|_{\im \lambda_i \cap ]-\infty,0]} = \gamma^-\circ h_0|_{\im \lambda_i \cap ]-\infty,0]}.
 \end{displaymath}
 We have to prove that $\gamma^+$ and $\gamma^-$ coincide. Recall from the proof of Proposition \ref{prop: homeo:lift} that the map $\tilde{h}$ is equivariant. Since $G = \langle r\rangle$ is a commutative group, we obtain for $x \in \im \lambda_i \cap\, ]\hspace{-2pt}-\infty , 0[$: 
  \begin{displaymath}
    \gamma^- \circ h_0 (x) = \tilde{h}(x) = - \tilde{h} (-x) = - \gamma^+ \circ h_0 (-x) = r \circ \gamma^+ \circ h_0 (r(x)) = r \circ \gamma^+ \circ r \circ h_0(x) = \gamma^+ \circ h_0(x). 
  \end{displaymath}
 As $h_0 (x) \neq 0$ is a non-singular point, indeed $\gamma^- = \gamma^+$ follows.   
\end{proof}

\begin{cor}
  For an orbifold $(Q,\uU)$ as in \ref{setup: good:ofd}, the map $\Delta_c \colon \Diff_c^G (\RR^d) \rightarrow \Diffc{Q,\uU}$ introduced in \ref{setup: new} is an isomorphism of groups. 
\end{cor}

\begin{proof}
 If $(Q,\uU)$ is one of the orbifolds introduced in \ref{setup: good:ofd}, the only element in $G \cap \Diff_c^G (\RR^d)$ is the unit element $\id_{\RR^d}$. Hence $\Delta_c$ is an injective group homomorphism. We will prove that $\Delta_c$ is a surjective map.\\ 
 To this end consider $[\hat{h}] \in \Diffc{Q,\uU}$ with a representative $(h, \set{\tilde{h}}, [P,\nu])$. Here the lift $\tilde{h} \colon \RR^d \rightarrow \RR^d$ has been chosen with $\tilde{h} \in D^{-1} ([\hat{h}])$ (which is possible by Proposition \ref{prop: sm:globlift}). Let $K \subseteq Q$ be a compact set with $h|_{Q \setminus K} \equiv \id_{Q \setminus K}$. As $\pi \colon \RR^d \rightarrow Q$ is a proper map by Lemma \ref{lem: orbitmap}, the set $\pi^{-1} (K)$ is compact. Choose a compact set $L \subseteq \RR^d$ with $\pi^{-1} (K) \subseteq L$ and $\RR^d \setminus L$ being connected if $d\geq 2$. If $d=1$, we may assume that $0 \in L$ and $\RR \setminus L$ contains exactly two connected components. Recall from the proof of Proposition \ref{prop: homeo:lift} that the lift $\tilde{h}$ has been constructed with respect to an arbitrary pair $x_0 \in \RR^d \setminus \set{0}$ and $y_0 \in \pi^{-1} h \pi(x_0)$ such that $\tilde{h} (x_0) =y_0$ (if $d\geq 2$). Without loss of generality, choose $x_0 \in \RR^d \setminus L$. Since $h|_{Q \setminus \pi (L)} \equiv \id_{Q\setminus \pi (L)}$ holds, one can set $y_0 = x_0$. We claim that the lift $\tilde{h}$ with respect to these choices is contained in $\Diff^G_c (\RR^d)$. If this is true, then $\Delta_c (\tilde{h}) = D(\tilde{h}) = [\hat{h}]$ follows and $\Delta_c$ is a group isomorphism.\\[1em]
 To prove the claim, it suffices to prove that $\tilde{h}$ coincides with $\id_{\RR^d}$ outside the compact set $L$. We distinguish two cases: If $d \geq 2$, then $\tilde{h}$ is a lift of the identity on the connected set $\RR^d \setminus L$ which takes $x_0$ to $x_0$ ansd so is $\id_{\RR^d \setminus L}$. Hence, $\tilde{h}|_{\RR^d \setminus L} = \id_{\RR^d \setminus L}$ by uniqueness of lifts (cf.\ \cite[Proposition 1.34]{hatcher2002}). Hence $\tilde{h} \in \Diff_c^G (\RR^d)$ follows. \\
 If $d=1$, by choice of $L$ the space $\RR \setminus L$ contains two connected components $C_1,C_2$. Now \cite[Lemma 2.11]{follie2003} yields $\tilde{h}|_{C_i} = g_i|_{C_i}$ for some $g_i \in G$ and $i\in \set{1,2}$. By construction of $\tilde{h}$, we have $\tilde{h} (]0,\infty[) \subseteq \, ]0,\infty[$ and $\tilde{h} (]\hspace{-2pt}-\infty,0[) \subseteq \,]\hspace{-2pt}-\infty,0[$, whence $g_1 = g_2 = \id_\RR$ and thus $\tilde{h} \in \Diff_c^G (\RR)$. 
\end{proof}

For the rest of this section, if $(Q,\uU)$ is an orbifold as in \ref{setup: good:ofd}, we endow the group $\Diff_c^G (\RR^d)$ with the unique Lie group structure turning $\Delta_c$ into an isomorphism of Lie groups. We shall use this Lie group structure to construct a Lie group structure on $\Diff^G (\RR^d)$. 

\begin{prop}\label{prop: equi:lgp}
 Let $d$ be in $\NN$ and $G$ be a subgroup of $O(d)$ as defined in \ref{setup: good:ofd}. Then the group $\Diff^G (\RR^d)$ may be endowed with a unique Lie group structure such that $\Diff_c^G (\RR^d)$ becomes an open subgroup.
\end{prop}

\begin{proof}
 We use the Lie group structure on $\Diff_c^G (\RR^d)$ together with the construction principle in Proposition \ref{prop: Lgp:locd}. Clearly the subgroup $\Diff^G_c (\RR^d)$ of $\Diff^G (\RR^d)$ satisfies the requirements of Proposition \ref{prop: Lgp:locd} (a). To obtain a Lie group structure on $\Diff^G (\RR^d)$, we have to verify the following condition: For each $\tilde{g} \in \Diff^G (\RR^d)$, the mapping 
  \begin{displaymath}
   c_{\tilde{g}} \colon \textstyle\Diff^G_c \displaystyle (\RR^d) \rightarrow \textstyle \Diff^G_c \displaystyle (\RR^d) , \tilde{h} \mapsto \tilde{g} \circ \tilde{h} \circ \tilde{g}^{-1} 
  \end{displaymath}
 makes sense and is smooth. First, we notice that for each $\tilde{g} \in \Diff^G (\RR^d)$ the map $c_{\tilde{g}}$ makes sense because $\supp c_{\tilde{g}} (\tilde{h}) \subseteq \tilde{g}(\supp \tilde{h})$ is compact.\\
 Since $D$ is a group homomorphism, we have for $\tilde{g} \in \Diff^G (\RR^d)$ and $\tilde{h} \in \Diff^G_c (\RR^d)$ the identity 
    \begin{displaymath}
     D(\tilde{g}) \Delta_c (\tilde{h}) D(\tilde{g}) =  D(\tilde{g}) D (\tilde{h}) D(\tilde{g}) = D(\tilde{g} \tilde{h} \tilde{g}) = D (c_{\tilde{g}} (\tilde{h}) = \Delta_c  (c_{\tilde{g}} (\tilde{h}).
    \end{displaymath}
 Thus $c_{D(\tilde{g})} \circ \Delta_c = \Delta_c \circ c_{\tilde{g}}$ holds. 
Here $c_{D (\tilde{g})} \colon \Diffc{Q,\uU} \rightarrow \Diffc{Q,\uU}, [\hat{f}] \mapsto D(\tilde{g})\circ [\hat{f}] \circ D(\tilde{g})^{-1}$ is the conjugation map (cf.\ \ref{setup: observ}). Since $\Difforb{Q,\uU}$ is a Lie group which contains $\Diffc{Q,\uU}$ as an open subgroup, $c_{D (\tilde{g})}$ is a smooth map. Furthermore, $\Delta_c$ and $\Delta_c^{-1}$ are smooth, whence $c_{\tilde{g}}$ is a smooth map. Now Proposition \ref{prop: Lgp:locd} (b) proves the assertion. Furthermore, $G$ is a discrete normal subgroup of $\Diff^G (\RR^d)$ and $\Diff^G (\RR^d)/G$ is a Lie group such that $\Delta$ becomes an isomorphism of Lie groups.
\end{proof}

\begin{rem}
 The group $\Diff (\RR^d)$ has been turned into a Lie group modeled on the space $\mathfrak{X}_c (\RR^d)$ of compactly supported vector fields in \cite{hg2005}. The Lie group $\Diff (\RR^d)$ turns the subgroup $\Diff^G (\RR^d)$ into a closed Lie subgroup modeled on the space of equivariant compactly supported vector fields. The induced Lie group structure is precisely the Lie group structure on $\Diff^G (\RR^d)$ constructed in Proposition \ref{prop: equi:lgp}. We sketch a proof for these facts:\\[1em] 
 \textbf{$\Diff^G (\RR^d)$ is a closed Lie subgroup of $\Diff (\RR^d)$}: The Lie group structure on $\Diff (\RR^d)$ has been obtained by applying the construction principle (Proposition \ref{prop: Lgp:locd}) to the Lie group $\Diff_c (\RR^d)$. Hence it suffices to prove that $\Diff_c (\RR^d)$ contains $\Diff_c^G (\RR^d)$ as a closed Lie subgroup, whose induced Lie group structure coincides with the structure induced by $\Delta_c \colon \Diff_c^G (\RR^d) \rightarrow \Diffc{Q,\uU}$. For the following arguments we identify $\mathfrak{X}_c (\RR^d)$ with the space $C_c^\infty (\RR^d , \RR^d)$  (see Definition \ref{defn: cs:vf} for details). However, we suppress the identification in the notation below.\\
 In \cite[Theorem 6.5]{hg2005}, the group $\Diff_c (\RR^d)$ was turned into a Lie group using the following global chart at the identity element: 
  \begin{displaymath}
   \alpha \colon \mathfrak{X}_c (\RR^d) \supseteq \Omega \rightarrow \textstyle \Diff_c \displaystyle (\RR^d) , \sigma \mapsto \id_{\RR^d} + \sigma
  \end{displaymath}
 (for a suitable open zero-neighborhood $\Omega \subseteq \mathfrak{X}_c (\RR^d)$). If we endow $\RR^d$ with the flat Riemannian metric $\rho_f$ (i.e.\ the one associated to the euclidean metric), then the Riemannian exponential map with respect to this metric is given by $\exp (v) = v + p$ for $v \in T_p \RR^d$ (cf. \cite[Example after Definition 1.6.4]{klingenberg1995}). Hence we may rewrite $\alpha$ as $\alpha (\sigma) = \exp \circ \sigma$.\\
 By construction, $G$ is a subgroup of the orthogonal group, whence each element in $G$ is a Riemannian isometry with respect to $\rho_f$. In particular, $\exp$ commutes via $\exp \circ \gamma = \exp \circ (d \gamma) = \gamma . \exp$ with every $\gamma\in G$. Thus for $\sigma \in \mathfrak{X}_c (\RR^d)$, we derive $\exp \circ \sigma \circ \gamma = \exp \circ d\gamma \circ \sigma = \gamma . \exp \circ \sigma$ for each $\gamma \in G$. In other words, $\alpha (\Omega \cap \mathfrak{X}_c^G (\RR^d)) \subseteq \Diff^G_c (\RR^d)$.
 We remark that each element in $\Diff_c^G (\RR^d)$ is equivariant with respect to the $G$-action. Hence $\alpha (\sigma) \in \im \alpha \cap \Diff^G_c (\RR^d)$ implies $\exp \circ \sigma \circ \gamma = \gamma. \exp \circ \sigma = \exp \circ d\gamma . \sigma$. As the Riemannian exponential map $\exp|_{T_p \RR^d}$ is injective for each $p \in \RR^d$, the identity $d\gamma .\sigma =\sigma \circ \gamma$ follows for each $\gamma \in G$. We conclude  
  \begin{displaymath}
   \alpha (\Omega \cap  \mathfrak{X}_c^G (\RR^d)) = \textstyle \Diff_c^G \displaystyle (\RR^d).
  \end{displaymath}
 The restriction of $\alpha$ to $\Omega \cap  \mathfrak{X}_c^G (\RR^d)$ induces a submanifold chart $\tilde{\alpha}$ for $\Diff_c^G (\RR^d)$. Hence $\Diff_c^G (\RR^d)$ becomes a Lie subgroup of $\Diff_c (\RR^d)$ modeled on the closed vector subspace $\mathfrak{X}_c^G (\RR^d)$ of $\mathfrak{X}_c(\RR^d)$ (cf.\ Example \ref{ex: os:good}). This proves the first assertion. We denote by $\Diff_c^G (\RR^d)^*$ the group $\Diff_c^G (\RR^d)$ with the structure of a closed Lie subgroup of $\Diff_c (\RR^d)$. The symbol $\Diff_c^G (\RR^d)$ will denote the Lie group constructed in Proposition \ref{prop: equi:lgp}.\\[1em]
 \textbf{The Lie groups $\Diff_c^G (\RR^d)$ and $\Diff_c^G (\RR^d)^*$ coincide}: Observe that each element of $G$ is a Riemannian isometry with respect to the flat Riemannian metric $\rho_f$ on $\RR^d$. Since $\RR^d/G$ is an orbifold with global chart, the family $(\rho_f)$ induces a Riemannian orbifold metric $\rho$ on $(Q,\uU)$. Let $[\expo]$ be the Riemannian orbifold exponential map associated to $\rho$. By Lemma \ref{lem: indep:rm}, we may assume that the Lie group $\Difforb{Q,\uU}$ has been constructed with respect to the Riemannian orbifold metric $\rho$. Therefore a chart around the identity element for the open Lie subgroup $\Diffc{Q,\uU}$ is given by 
 \begin{displaymath}
  E \colon \Osc{Q} \supseteq \hH \rightarrow \Diffc{Q,\uU} , [\hat{\sigma}] \mapsto [\expo] \circ [\hat{\sigma}] 
 \end{displaymath}
 for a suitable open zero-neighborhood $\hH$. The Riemannian exponential map $\exp \colon T\RR^d \rightarrow \RR^d$ associated to $\rho_f$ is a lift of $[\expo]$. Since $(\RR^d,G,\pi)$ is a global chart for $\RR^d/G$, we obtain a representative $(\expo , \set{\exp} , [P,\nu])$ of $[\expo]$ by Remark \ref{rem: lift:expo} (a). Hence for each $[\hat{\sigma}]$, a representative of $E([\hat{\sigma}])$ is induced by the lift $\exp \circ \sigma_{\RR^d}$ (where $\sigma_{\RR^d}$ is the unique canonical lift of $[\hat{\sigma}]$ in the global chart). Recall from Example \ref{ex: os:good} that the mapping $H \colon \Osc{\RR^d/G} \rightarrow \mathfrak{X}_c^G (\RR^d), [\hat{\sigma}] \mapsto \sigma_{\RR^d}$ is an isomorphism of topological vector spaces. Shrinking the open zero neighborhood $\hH$, we may assume that $H (\hH) \subseteq \Omega \cap \mathfrak{X}_c^G (\RR^d)$ holds. Combining these facts, a trivial computation shows that the following diagram is commutative: 
  \begin{displaymath}
   \begin{xy}
  \xymatrix{ \Osc{Q} \ar[d]^H & \hH  \ar[l]_-\supseteq \ar[d]^{H|_\hH} \ar[rr]^-{E} && \Diffc{Q,\uU}  \ar[d]^{\Delta_c^{-1}} \\ 
             \textstyle \mathfrak{X}_c^G\displaystyle (\RR^d) &  **[r]{\Omega \cap \textstyle \mathfrak{X}_c^G\displaystyle (\RR^d)} \ar[l]_-\supseteq \ar[rr]^-{\tilde{\alpha}} && \Diff^G_c (\RR^d)^*}
    \end{xy}
  \end{displaymath}
 We deduce that the group homomorphism $\Delta_c^{-1}$ is smooth as a map into $\Diff_c^G (\RR^d)^*$ on a neighborhood of the identity in $\Difforb{Q,\uU}$. Hence $\Delta_c^{-1}$ is smooth as a map of the Lie group $\Diffc{Q,\uU}$ into $\Diff_c^G (\RR^d)^*$ by \cite[III.\ \S 1 2.\ Proposition 4]{bourbaki1989}. Vice versa the same holds for $\Delta_c \colon \Diff_c^G (\RR^d)^* \rightarrow \Diffc{Q,\uU}$. Therefore $\Delta_c$ is an isomorphism of Lie groups. However, $\Delta_c \colon \Diff_c^G (\RR^d) \rightarrow \Diffc{Q,\uU}$ is also an isomorphism of Lie groups. We obtain an isomorphism of Lie groups 
  \begin{displaymath}
    \textstyle \id_{\Diff^G_c (\RR^d)} \colon \Diff_c^G (\RR^d)^* \rightarrow \Diff_c^G (\RR^d),\ \id_{\Diff^G_c (\RR^d)} \coloneq \Delta_c^{-1} \circ \Delta_c,
    \end{displaymath} 
 whence the Lie groups $\Diff^G_c (\RR^d)$ and $\Diff^G_c (\RR^d)^*$ coincide. 
\end{rem}

\paragraph{Acknowledgement and remark:} The author would like to thank Helge Gl\"ockner for many helpful discussions on the subject of this work. This work is a version of the authors Ph.D-thesis.

\begin{appendix}\thispagestyle{empty}
\section{Hyperplanes and Paths in Euclidean Space}

The results in this appendix are part of the folklore. However, for the reader's convenience we provide full proofs for these known facts. As usual, a hyperplane $H$ in euclidean space $\RR^d$ is a linear subspace of codimension $1$ and a path is a continuous map from an interval to $\RR^d$.

\begin{Alem}\label{lem: avnhp}
 Let $d \in \NN$ and $X \subseteq \RR^d$ a linear subspace such that $\dim X \leq d-2$. Consider an open and path-connected subset $C \subseteq \RR^d$ and $x,y \in C \setminus X$. Then there exists a path $p \colon [0,1] \rightarrow C \setminus X$ connecting $x$ and $y$. In other words, $C\setminus X$ is path-connected.
\end{Alem}

\begin{proof}
 Without loss of generality, we may assume $X = \RR^{d-m}\times \set{0}$ and $m \geq 2$. The set $C$ is path-connected, whence there is a path $q\colon [0,1] \rightarrow C$ with $q(0)=x$ and $q(1)=y$. If the intersection $\im q \cap X$ is empty, there is nothing to prove. Otherwise we construct a path as follows:\\
 Consider the projections $\pi_X \colon \RR^d \rightarrow \RR^{d-m} \times \set{0}= X$ ans $\pi_2 \colon \RR^d \rightarrow \set{0} \times \RR^m$, respectively. The projections are continuous open maps, with $\pi_X + \pi_2 = \id_{\RR^d}$. Observe that $z \in X$ if and only if $\pi_2 (z) = 0$ holds. The set $\setm{q(t)}{t\in [0,1],\pi_2 (q(t))=0} = \im q \cap X$ is compact and does not contain $x$ and $y$. Therefore we can choose $x_i \in X, 1\leq i \leq N$ and $\ve >0$ with 
  \begin{displaymath}
   \im q \cap X \subseteq \bigcup_{1\leq i\leq N}B_{\ve } (x_i) \times B_{\ve} (0) \subseteq K \coloneq \bigcup_{1\leq i\leq N}\overline{B_{\ve } (x_i) \times B_{\ve} (0) } \subseteq C \setminus \set{x,y}.
  \end{displaymath}
 As each closed ball is path-connected, the sets $\overline{B_{\ve } (x_i) \times B_{\ve} (0)}$ are path-connected. Hence the set $K$ is a set with finitely many path-components $K_1 , \ldots , K_r$ (cf. \cite[p.\ 115]{dugun1966}). Each path-component is a union $K_i = \bigcup_{1\leq j \leq r_i} \overline{B_{\ve} (x_{i,j}) \times B_\ve (0)}$ and is thus compact. Furthermore, the boundary $\partial K$ satisfies $\partial K = \partial K_1 \cup \partial K_2 \cup \ldots \cup \partial K_r$, since the sets $K_i$ form a finite partition of closed and disjoint sets. As $\im q \cap X \subseteq K^\circ$ holds, we deduce that the boundary $\partial K_i$ does not contain elements of $\im q \cap X$. We construct the path by induction: The set $L_1 \coloneq q^{-1} (K_1)$ is a closed subset of $[0,1]$, which does not contain $0,1$ by construction. \textbf{Case 1}: If $L_1 = \emptyset$, set $q_1 \coloneq q$. \\ \textbf{Case 2}: If $L_1 \neq \emptyset$, the compactness of $L_1$ enables us to consider $s_1 \coloneq \min L_1$ and $t_1 \coloneq \max L_1$. For $t \in \set{s_1,t_1}$, we must have $q(t) \in \partial K_1$. As shown above, this implies $q(s_1),q (s_2) \not \in X$, i.e.\ $\pi_2 (q(s_1)), \pi_2 (q(t_1)) \in \overline{B_\ve (0)} \setminus \set{0}$ holds. Note that $\overline{B_\ve (0)} \setminus \set{0}$ is path-connected (by a variation of \cite[V.\ Theorem 2.2]{dugun1966}), since $m\geq 2$ is satisfied. Furthermore, $\pi_X (K_1)$ is path-connected, whence there is a path $\gamma_1 \colon [s_1,t_1] \rightarrow \pi_X (K_1) \times \left(\overline{B_\ve (0)} \setminus \set{0}\right) \subseteq K_1 \subseteq C$ with $\gamma_1 (s_1) = q(s_1)$ and $\gamma_1 (t_1) = q(t_1)$. Define a mapping
	\begin{displaymath}
	q_1 \colon [0,1] \rightarrow C,\ t \mapsto \begin{cases}
	                                               q (t) & t \in [0,1] \setminus ]s_1,t_1[ \\
						       \gamma_1 (t) & t \in [s_1,t_1].
	                                              \end{cases}
	\end{displaymath}
 By construction, $q_1$ is a path with $q_1(0) = x$ and $q_1 (1) =y$. Furthermore, $\im q_1 \cap K_1 = q_1 ([s_1,t_1])$ implies $\im q_1 \cap K_1 \cap X = \emptyset$. This also holds in Case 1. In either case, note that the definition of $q_1$ yields $\im q_1 \cap X \subseteq \bigcup_{2\leq i \leq r} K_i$.\\
 Assume that for all $i$ with $1\leq i < n \leq r$, we have already constructed a path $q_i$ connecting $x$ and $y$, whose image is contained in $C$ with $\im q_i \cap X \subseteq \bigcup_{i+1\leq j \leq r} K_j$. Consider the compact set $L_n \coloneq q_{n-1}^{-1} (K_n) \subseteq \, ]0,1[$. If $L_n$ is empty, simply set $q_n \coloneq q_{n-1}$ to obtain a path with the desired properties. Otherwise, we have to construct a path $q_n$ from $q_{n-1}$ such that the image does not intersect $(K_1 \cup \ldots \cup K_n) \cap X$. Apply the above construction verbatim with $L_n \neq \emptyset$ and $q_{n-1}$ instead of $L_1$ and $q$. Since $q_{n-1}$ does not intersect $K_i \cap X$ for each $1\leq i \leq n-1$, the construction yields a mapping $q_{n}$ with $\im q_n \cap X \subseteq \bigcup_{n+1 \leq i \leq r} K_i$, whose image is contained in $C$. Summing up, after finitely many steps the mapping $p\coloneq q_r$ satisfies: $\im p \subseteq C$, $p(0)=x, p(1) =y$ and $\im p \cap X \subseteq \bigcup_{r+1\leq i \leq r} K_i = \emptyset$. Hence $p$ is a path with the desired properties.   
\end{proof}
  
\begin{Alem}\label{lem: RnGeom} \no{lem: RnGeom}
 Let $d, m \in \NN$, $C$ be an open connected subset of $\RR^d$ and $(X_i)_{i=1 , \ldots , m}$ be a family of vector subspaces of $\RR^d$ such that $\dim X_i < d$ for all $1\leq i\leq m$ and $X_i \neq X_j$ for $i \neq j$.  
\begin{compactenum}
	\item For each pair $x,y \in C \setminus \bigcup_{i=1}^m X_i$, there is a path $p \colon [0,1] \rightarrow C$ such that 
		\begin{compactitem}
		 \item[1.] $p(0) = x,\ p(1)=y$,
		 \item[2.] $p([0,1]) \cap X_i = \emptyset$ for all $i$ such that $\dim X_i \leq d -2$,
		 \item[3.] $p([0,1])\cap X_i \cap X_j = \emptyset$ for all $i,j$ such that $i\neq j$.
		\end{compactitem}
	\item Assume there is $k \in \NN_0$ such that $\dim X_i = d -1$ if $1\leq i \leq k$ and $\dim X_i < d-1$ otherwise. Then the set $\RR^d \setminus \bigcup_{i=1}^m X_i$ with the subspace topology  has at most $2^k$ (path-)connected components.
        \item If $C \subseteq \RR^d$ is a convex open subset, then $C \setminus \bigcup_{i=1}^m X_i$ possesses at most $2^k$ connected components.
              \end{compactenum}
\end{Alem}

\begin{proof}
 \begin{compactenum}
  \item Since for $i \neq j$ we have $\dim X_i \cap X_j \leq d-2$, it suffices to construct a path $p$ which satisfies Properties 1.\ and 2.\ for an arbitrary finite number of subspaces $Y_i$ with $\dim Y_i \leq d-2$. Since $C$ is path-connected, $C \setminus Y_1$ is path-connected by Lemma \ref{lem: avnhp}. Iteratively, $C\setminus Y_1 \setminus Y_2 \setminus \cdots \setminus Y_m = C \setminus (Y_1 \cup \ldots \cup Y_m)$ is path connected by Lemma \ref{lem: avnhp}.
  \item The subspaces $X_i$ are closed in $\RR^d$, whence $\Omega \coloneq \RR^d \setminus \bigcup_{i=1}^m X_i$ is an open set. The components of $\Omega$ coincide with the path-components of $\Omega$ by \cite[V. 5.6]{dugun1966}. We claim that there are at most $2^k$ path-components. For a hyperplane $X_j$, we consider the two half spaces $H_j^+, H_j^-$ such that $\RR^d$ is the disjoint union $H_j^+ \cup X_j \cup H_j^-$. The half-spaces are the path-components of $\RR^d \setminus X_j$. Each half-space is a convex set. We observe that each intersection of half-spaces $H_{1}^{\sigma (1)} \cap \ldots \cap H_{k}^{\sigma (k)}$ with $\sigma \colon \set{1,2,\ldots , k} \rightarrow \set{+,-}$ is again a convex set. From (a) we deduce that these sets yield path-connected subsets of $\RR^d \setminus \bigcup_{1\leq j \leq m} X_j$ if we remove $\bigcup_{k+1 \leq j \leq m} X_j$. Hence $\RR^d \setminus \bigcup_{i=1}^m$ is realized as a union of no more than $2^k$ path-connected sets, from which the assertion follows.
  \item From the proof of (b), we deduce that the components are induced by intersections of $k$ half-spaces, which are convex sets. However the same holds for the subset $C \cap H_{j_1}^{\sigma (1)} \cap \ldots \cap H_{j_r}^{\sigma (k)}$. From part (a) we deduce with arguments as in (b) that all non-empty sets of this kind induce the connected components of $C \setminus \bigcup_{i=1}^m X_i$. As there are at most $2^k$ non-empty sets of this kind, the assertion follows. 
 \end{compactenum}
\end{proof}\thispagestyle{empty}
\section{Group Actions and Newman`s Theorem}

In this section, we recall several basic facts concerning group actions, orbit spaces and quotient mappings to orbit spaces. We are interested only in continuous group actions, whence each group action in this thesis will be required to be continuous. Several basic results will be repeated to fix some notation. For further information on group actions, we recommend \cite{bredon1972, dieck1987}.

\subsection{Group actions}

\begin{defn}[Group actions of topological groups]
Let $G$ be a topological group and $X$ a topological space. A \ind{group action}{$G$-action} on $X$ is a continuous map $\Theta \colon G \times X \rightarrow X$ such that: 
\begin{compactenum}
            \item $\Theta (\one , x) =x$ for all $x \in X$, where $\one$ is the identity element of $G$.
	    \item $\Theta (g_2,\Theta (g_1, x))= \Theta (g_2 g_1 , x)$ for all $g_1, g_2 \in G$ and $x \in X$.	
 \end{compactenum}
The pair $(X, \Theta)$ (or $(X,G)$ if the action is clear) is called a \ind{G-space@$G$-space}{$G$-space} and we denote it usually just by the underlying space $X$. We shall abbreviate $g.x \coloneq \Theta (g,x)$ if it is clear which action is meant.\\ 
For $x \in X$ the \ind{}{orbit} of $x$ is the set $G.x \coloneq \setm{g.x}{g\in G}$. Let $X/G \coloneq \setm{G.x}{x\in X}$ be the set of all orbits and endow it with the quotient topology induced by $p \colon X \rightarrow X/G, x\mapsto G.x$. The space $X/G$ is called the \ind{}{orbit space} of the $G$-space $X$. 
\end{defn}

\begin{defn}[Isotropy subgroups and fixed point sets]
 Let $X$ be a $G$-space. Define the \ind{group action!isotropy group}{isotropy group} $G_x\coloneq \setm{g\in G}{g.x=x}$ of $x\in X$. \\ For $g \in G$, the \ind{set!of fixed points}{set of fixed points of $g$} will be denoted by $\Sigma_g = \setm{x \in X}{g.x=x}$\glsadd{fixed points} and we write 
  \begin{displaymath}
      \Sigma_G \coloneq \setm{x \in X}{G_x \neq \set{\one}} = \bigcup_{g \in G\setminus \set{\one}} \Sigma_g.
  \end{displaymath}
For a subset $S \subseteq X$, we define $g.S \coloneq \setm{g.x}{x\in S}$ and let $G_S \coloneq \setm{g \in G}{g.S =S}$ be the \ind{group action!isotropy group of a set}{isotropy group of $S$}. A subset $S \subseteq X$ is called \ind{set!G-invariant@$G$-invariant}{$G$-invariant} if $G_S = G$ holds. Furthermore, a \ind{set!G-stable@$G$-stable}{$G$-stable} subset of $X$ is a connected set $S \subseteq X$ such that for $g \in G$ either $g.S =S$ or $g.S \cap S = \emptyset$ is satisfied.
\end{defn}

The elegant proof of the following lemma has been communicated to the author by A. Pohl:

\begin{lem}\label{lem: st:nbhd} \no{lem: st:nbhd}
 Let $X$ be a manifold, $G$ a finite topological group acting on $X$ via homeomorphisms, i.e. $\Theta (g,\cdot) \colon X \rightarrow X$ is a homeomorphism for each $g \in G$. Then, for each $x \in X$, there exist arbitrarily small open $G$-stable neighborhoods of $x$ whose isotropy groups coincide with $G_x$. In particular, the $G$-stable open sets form a base for the topology on $X$.
\end{lem}

\begin{proof} Let $U$ be any neighborhood of $x$ and $G.x = \set{x_1,x_2,\ldots ,x_n}$ be the distinct elements in the $G$-orbit of $x$, i.e. $x_i \neq x_j$ for $i \neq j$. Without loss of generality, $x=x_1$ holds. For $i=1,\ldots , n$, choose an open neighborhood $U_i$ of $x_i$ with the following property: For $i \neq j$, the sets $U_i$ and $U_j$ are disjoint and $U_1 \subseteq U$ holds. For $i=1,\ldots,n$, define $G_i^1 \coloneq \setm{g \in G}{g.x_i = x}$ and set 
    \begin{displaymath}
     S' \coloneq \bigcap_{1 \leq i \leq n} \bigcap_{g \in G_i^1} g.U_i.
    \end{displaymath}
 As $G$ acts by homeomorphisms, the set $S' \subseteq U_1 \subseteq U$ is an open neighborhood of $x$. Consider $h \in G$. If $h.x=x_i$ holds, this implies $h^{-1} \in G_i^1$. Therefore $S' \subseteq h^{-1}.U_i$ yields $h.S' \subseteq U_i$. For $i\neq 1$ we deduce from $U_i \cap U_1 = \emptyset$ and $S' \subseteq U_1$ for $h$ as above $h.S' \cap S' = \emptyset$. On the other hand, for $i=1$ we have $h \in G_x$, whence $hG_j^1 =G_j^1$ for all $j$ and thus 
    \begin{equation}\label{eq: perm}
     h.S' = \bigcap_{j=1}^n \bigcap_{g \in G_j^1} (hg).U_j = \bigcap_{j=1}^n \bigcap_{g \in G_j^1} g.U_j = S' .
    \end{equation}
Let $S$ be the connected component of $S'$ which contains $x$. As $X$ is locally path connected, $S$ is an open neighborhood of $x$ by \cite[V.\ Theorem 4.2]{dugun1966}. Since $G$ acts by homeomorphisms, by \eqref{eq: perm} $G_x$ permutes the connected components of $S'$ and fixes $x$. Combine \eqref{eq: perm} and the fact $h.S' \cap S'=\emptyset$ for $h \in G\setminus G_x$. We deduce that $G_S = G_x$ holds and $S$ is a $G$-stable open neighborhood of $x$ which is contained in $S' \subseteq U$.
\end{proof}

\begin{lem}[\hspace{-0.5pt}{\cite[Proposition 3.1 and Proposition 3.6]{dieck1987}}]\label{lem: orbitmap} \no{lem: orbitmap}
Let $X$ be a Hausdorff $G$-space and $G$ a compact topological group. Consider the quotient map $\pi \colon X \rightarrow X/G, x \mapsto G.x$ onto the orbit space. Then
 \begin{compactenum}
  \item $X/G$ is a Hausdorff space.
  \item $\pi$ is a continuous, open and closed map.
  \item $\pi$ is a proper map.
  \item $X$ is compact if and only if $X/G$ is compact.
  \item $X$ is locally compact if and only if $X/G$ is locally compact.
 \end{compactenum}
\end{lem}

\begin{rem}\label{rem: dgp:act}\no{rem: dgp:act}
 Let $M$ be a (possibly infinite-dimensional) manifold. The discrete topology is the unique Hausdorff topology turning a finite subgroup $G$ of $\Diff^r (M)$ into a topological group. The natural mapping $\Theta \colon G \times M \rightarrow M, (g,x) \mapsto g(x)$ is continuous since each element in $G$ is continuous and $G$ is endowed with the discrete topology. Hence each finite subgroup of $\Diff^r (M)$ induces a canonical action of a compact group on $M$ which satisfies the prerequisites of Lemma \ref{lem: orbitmap}.   
\end{rem}

\begin{defn}
 Let $f \colon X \rightarrow Y$ be a map from the $G$-space $X$ to the $H$-space $Y$.
\begin{compactenum}
        \item If there is a group homomorphism $\lambda \colon G \rightarrow H$ such that $f(g.x) = \lambda (g) . f(x)$ holds for all $x \in X, g \in G$, $f$ is called \ind{equivariant map!with respect to a morphism}{equivariant with respect to} $\lambda$.
        \item If $G$ and $H$ coincide and $f(g.x) = g.f(x)$ holds for all $x \in X$, $g \in G$, we call $f$ \ind{equivariant map}{equivariant}. An equivariant homeomorphism is called an \ind{equivariant map!equivalence}{equivalence}.
        \item Let $f$ be a homeomorphism and $G=H$. If there is a group automorphism $\alpha \colon G \rightarrow G$ with $f(g.x) = \alpha (g) . f(x)$ for all $x \in X, g \in G$, then the map $f$ is called a \ind{equivariant map!weak equivalence}{weak equivalence}. 
 \end{compactenum}
 Notice that the inverse of a (weak) equivalence is again a (weak) equivalence (cf. \cite[I 2.]{bredon1972}).
\end{defn}

\begin{defn}
 Let $M$ be a smooth manifold which is also a $G$-space. We define the set 
  \begin{displaymath}
   \textstyle \Diff^G (M) \displaystyle\coloneq \setm{f \in \Diff (M)}{f \text{ is a weak equivalence}}.
  \end{displaymath}\glsadd{DiffG}
\end{defn}

\begin{rem}
 It is easy to check the following facts about $\Diff^G (M)$: 
  \begin{compactenum}
   \item The set $\Diff^G (M)$ is a subgroup of $\Diff (M)$.
   \item If $G \subseteq \Diff (M)$ acts via the natural action on $M$, then $G \subseteq \Diff^G (M)$ follows. In this case, $G$ is a normal subgroup of $\Diff^G (M)$.
  \end{compactenum}
\end{rem}

\subsection{Newman`s Theorem}
The following theorem of M.H.A. Newman is an important tool to investigate the structure of orbifolds (for a proof see \cite{dress1969} also cf. \cite[III 9.]{bredon1972}): 

\begin{thm}[Newman 1931] \label{thm: newman} \no{thm: newman}
 Let $G$ be a finite group acting effectively by homeomorphisms on a connected paracompact finite dimensional manifold $M$. Then the set $M \setminus \Sigma_G$ of points with trivial isotropy group is  dense and open in $M$.
\end{thm}
\noindent
In the situation of Theorem \ref{thm: newman}, the elements of $\Sigma_G$ are called \ind{point!singular}{singular points} and the elements of $M \setminus \Sigma_G$ are called \ind{point!non-singular}{non-singular points}. If $G$ acts by $C^\infty$ diffeomorphisms on a paracompact smooth manifold, then Newman's Theorem is much easier to prove, see \cite[Lemma 2.10]{follie2003}.\\
We compile several interesting consequences of Newman's Theorem. For further information, we refer to \cite[Section 2.4]{follie2003}.

\begin{lem}[cf. {\cite[p. 36]{follie2003}}]\label{lem: lin:ch} \no{lem: lin:ch}
 Let $M$ be a smooth finite dimensional paracompact manifold, $G$ a finite subgroup of $\Diff(M)$ and $x \in M$. Then there exists arbitrarily small $G$-stable charts $(W, \kappa)$ with $x \in W$ such that $\kappa (x) = 0$ and $\kappa$ conjugates the isotropy group $G_x$ to a (finite) group of orthogonal transformations on $\kappa (W)$. Furthermore, $T_x g = \id_{T_x M}$ implies $g|_{W} = \id_W$ for each $g \in G_x$; if $M$ is connected it implies $g=\id_M$.
\end{lem}

\begin{proof}
 Since $G$ is finite, we may choose a $G$-invariant Riemannian metric on $M$ by \cite[Proposition 2.8]{follie2003}. The group $G$  thus acts via Riemannian isometries with respect to this metric. Let $\exp_M$ be the Riemannian exponential map with respect to this metric. By \cite[Theorem 1.6.12]{klingenberg1995}, we may choose $\ve >0$ such that $\exp_M$ induces a diffeomorphism from the open ball $B_\ve (0_x)$ centered at $0_x$ in $T_x(M)$ to an open neighbourhood $W$ of $x$, $\exp_{M ,x} \colon B_\ve (0_x) \rightarrow W \subseteq M$. As the metric is $G$-invariant, each $g \in G_x$ induces an orthogonal transformation $T_x g$ of $T_x M$. Since $\exp_M$ commutes with Riemannian isometries on its domain, we deduce $\exp_{M,x} \circ T_x g|_{\dom \exp_{M,x}} = g \circ \exp_{M,x}$. This formula shows that $T_xg = \id$ implies $g|_{W} = \id_W$, and also that $W$ is $G_x$-invariant. By continuity of $\exp_M$, we can shrink $\ve$ to ensure that $W$ is contained in $G$-stable neighborhood of $x$ (cf. Lemma \ref{lem: st:nbhd}). Hence there is $\ve >0$ such that $\exp_{M,x} (B_\ve (0_x)) = W$ is a $G$-stable subset with $G_W = G_x$. For such a $W$, define $\kappa \coloneq (\exp_{M,x}|_{B_\ve (0_x)})^{-1}$. The pair $(W,\kappa)$ satisfies the assertion. In particular, $W$ may be taken arbitrarily small.\\ For the final assertion, note that $g\neq \id_M$ implies $g|_W \neq \id_W$, by Newman's Theorem.
\end{proof}

\begin{lem}\label{lem: ts:nspc} \no{lem: ts:nspc}
 Let $M$ be a connected paracompact smooth manifold and $G$ be a finite subgroup of $\Diff (M)$. Denote by $\Sigma_{TG}$ the set of singular points with respect to the derived action $G\times TM \rightarrow TM, (g,X) \mapsto g. X \coloneq Tg (X)$ of $G$ on $TM$. For each open connected set $U \subseteq TM$, the set of non-singular points $U\setminus \Sigma_{TG}$ is (path-)connected. 
\end{lem}

\begin{proof}
 Without loss of generality we may assume $U \neq \emptyset$. Let $C$ be a component of $U \setminus \Sigma_{TG}$ and $\overline{C}$ be its closure in $U$. We will show that $\overline{C}$ is open. The connectedness of $U$ then entails $\overline{C} = U$. If there was another component $D \neq C$, then $\overline{C} \cap D =\emptyset$, because $D$ is open and $C\cap D =\emptyset$. But $D \subseteq U = \overline{C}$ yields a contradiction, whence $U \setminus \Sigma_{TG}$ is connected.\\
 To see that $\overline{C}$ is open, let $X \in \partial C$ (the boundary with respect to $U$). Then $X \in \Sigma_{TG}$ as $C$ is open and closed in the open subset $U \setminus \Sigma_{TG}$ of $U$.  
 By definition of the derived action for $g \in G$ we have $\pi_{TM} (X) = \pi_{TM} (g.X) = g.\pi_{TM} (X)$ if $g.X =X$. This implies $G_X \subseteq G_{\pi_{TM} (X)}$. By Lemma \ref{lem: lin:ch}, there is a $G$-stable manifold-chart $(W,\kappa)$ such that $\pi_{TM} (X) \in W$, $G_W = G_{\pi_{TM} (X)}$ and $\kappa$ conjugates $G_W$ to a finite group of orthogonal transformations on $\kappa (W) = B_\ve (0) \subseteq \RR^{d}$ for $d = \dim M$ and some $\ve >0$. For $g \in G_{\pi_{TM} (X)}$, let $\tilde{g}$ be the orthogonal transformation conjugate to $g$, i.e.\ $\tilde{g}$ is a linear map with $\tilde{g} \circ \kappa = \kappa \circ g$. The functoriality of the tangent functor implies $T \tilde{g} T\kappa = T\kappa Tg$. Taking suitable identifications, $T\tilde{g} = (\tilde{g}|_{B_\ve (0)} \circ \text{pr}_1, d\tilde{g}) = (\tilde{g}|_{B_\ve (0)} \times \tilde{g}$ is the restriction of a linear map. Thus $T\kappa$ conjugates the action of $G_W = G_{\pi_{TM} (X)}$ on $TW$ to a linear action on $T\kappa (TW) =B_\ve (0) \times \RR^d$. Since $W$ is $G$-stable with $G_W = G_{\pi_{TM} (X)}$, the set $TW$ is $G$-stable with $G_{TW} = G_{\pi_{TM} (X)}$ by definition of the derived action. Hence $TW \cap \Sigma_{TG} = TW \cap \Sigma_{TG_{\pi_{TW} (X)}}$ holds. Choose an open connected neighborhood $\Omega$ of $X$ in $TW \cap U$. If $\Omega \setminus \Sigma_{TG}$ is a connected set, then $(\Omega \setminus \Sigma_{TG}) \cap C = \Omega \cap C \neq \emptyset$ follows as $X \in \overline{C}$ and thus $\Omega \setminus \Sigma_{TG} \subseteq C$. As $\Omega \setminus \Sigma_{TG}$ is dense in $\Omega$ by Newman's Theorem, we deduce that $\Omega \subseteq \overline{C}$. Thus $\overline{C}$ will be open as required.\\ To verify this, observe that $\Omega \subseteq TW$ entails $\Omega \cap \Sigma_{TG} = \Omega \cap \Sigma_{G_{\pi_{TM} (X)}}$. Consider the open sets $\tilde{\Omega} \coloneq T\kappa (\Omega)$ and $\tilde{V} \coloneq T\kappa (\Omega \setminus \Sigma_{TG}) = \tilde{\Omega} \setminus T\kappa (\Sigma_{TG_{\pi_{TM} (X)}})$. We claim that $\tilde{V}$ is connected. If this is true, the same holds for $\Omega \setminus \Sigma_{TG}$, whence the proof is complete.\\
 \textbf{Proof of the claim:} As $T\kappa$ conjugates the group action to a linear action, the set $\tilde{\Omega} \cap T\kappa (TW \cap \Sigma_{TG})$ is the intersection of the open (path-)connected set $\tilde{\Omega}$ with a finite union of linear subspaces of $\RR^{2d}$. By Lemma \ref{lem: RnGeom}, the set $\tilde{V}$ will be connected if for each $g \in G_{\pi_{TM} (X)}$ the fixed point set of the associated linear map $T\tilde{g}$ is not a hyperplane in $\RR^{2d}$. For each $g \in G_{\pi_{TM} (X)} \setminus \set{\id_M}$, Lemma \ref{lem: lin:ch} implies that $\tilde{g}$ is not the identity map. From \cite[I.\ Proposition 2.18 (1)]{msnpc1999}, we deduce that the fixed points of $\tilde{g}$ are contained in a hyperplane $H \subsetneq \RR^d$. Each linear subspace fixed by $T\tilde{g}$ is thus contained in $H \times H$ and $\dim (H \times H) = 2d-2$. Hence $T\tilde{g}$ does not fix any hyperplane, whence $\tilde{V}$ is connected. 
\end{proof}
\thispagestyle{empty}
\section{Infinite Dimensional Manifolds and Lie Groups}

In this section, we briefly recall the notions of infinite dimensional manifolds and infinite dimensional Lie groups. Manifolds and Lie groups modeled on infinite dimensional spaces may be defined almost exactly as in the finite dimensional case.

\subsection{Manifolds modeled on locally convex spaces}\label{app: mfd}

\begin{defn}
We recall from \cite{gn2007} that a \ind{manifold!with rough boundary}{manifold with rough boundary} modeled on a locally convex space~$E$ is a Hausdorff topological space $M$ with an atlas of smoothly compatible homeomorphisms $\phi\colon V_\phi\to U_\phi$ from open subsets $V_\phi$ of $M$ onto locally convex subsets $U_\phi\subseteq E$ with dense interior.
If each $U_\phi$ is open, $M$ is an ordinary manifold (without boundary). 
In a similar fashion $C^r$-manifolds may be defined for $r \in \NN_0$. Unless stated otherwise, every manifold will be assumed to be without boundary. 
Direct products of locally convex $C^k$-manifolds, tangent spaces and tangent bundles may be defined as in the finite dimensional setting. We refer to \cite{neeb2006} for details. 
\end{defn}

\begin{nota}\label{nota: conno} \no{nota: conno}
 Let $M, N$ be $C^r$-manifolds (where $1\leq r \leq \infty$) and $f \colon M \rightarrow N$ a mapping of class $C^r$. We denote by $Tf \colon TM \rightarrow TN$ the tangent map. Abbreviate by $T_x f \colon T_x M \rightarrow T_{f(x)}N$ the restriction of $Tf$ to the tangent space $T_x M$ of $M$ at $x\in M$. If $N$ is an open subset of a locally convex space $F$, the tangent map $Tf\colon TM \rightarrow TN \cong N \times F$ is given by $(x,v) \mapsto (f(x) , df(x,v))$ for $x\in M$, $v \in T_x M$ and a map $df \colon TM \rightarrow F$. 
 If  $f \colon U \rightarrow V$ is a $C^r$-map, where $U,V$ are open subsets of locally convex spaces $E$ and $F$, it is convenient to think of $df(x,\cdot)$ as a differential. Hence we canonically identify $T_x U \cong E$ and $T_y V \cong F$ to obtain $df(x,v) = T_x f(v)$. \\
 We let $\pi_{TM} \colon TM \rightarrow M$ be the bundle projection. For $r=\infty$ we denote by $\vect{M}$ the space of smooth vector fields, i.e. smooth mappings $X \colon M \rightarrow TM$ with $\pi_{TM} \circ X = \id_M$.  
\end{nota}

\subsection{Function spaces and their topologies} \label{sect: Crtop}

Our exposition of the $C^r$-topology follows \cite{hg2004}, but we allow locally convex subsets. Albeit the definition of differentiability differs from the one used in \cite{hg2004}, on open subsets of locally convex spaces over the field $\RR$ are equivalent by \cite[Proposition 7.4]{bgn2004}.  
 
 \begin{defn}[Compact-open topology]
  Let $X$, $Y$ be Hausdorff topological spaces, $K\subseteq X$ compact and $U \subseteq Y$ open. We define the set \glsadd{KU} 
	\begin{displaymath}
	 \lfloor K, U \rfloor \coloneq \setm{f \in C(X,Y)}{f(K) \subseteq U}.
	\end{displaymath}
 Then the sets 
	\begin{displaymath}
	 \lfloor K_1, U_1 \rfloor \cap \lfloor K_2, U_2 \rfloor \cap \ldots \cap \lfloor K_n, U_n \rfloor
	\end{displaymath}
 with $n \in \NN$, $K_i \subseteq X$ compact and $U_i \subseteq Y$ open for $1 \leq i \leq n$, are a base for a topology on $C(X,Y)$ (cf. \cite[Section 3.4]{Engelking1989}). It is called the \ind{}{compact-open topology} and we denote by $C(X,Y)_{c.o.}$\glsadd{Crco} the space $C(X,Y)$ with this topology.
 \end{defn}

 \begin{defn}\label{defn: CrTop}
 Let $E,F$ be locally convex topological vector spaces, $U \subseteq E$ a locally convex subset with dense interior and $r \in \NN_0 \cup \set{\infty}$. Endow $C^r(U,F)$ with the unique locally convex topology turning 
 \begin{displaymath} 
   (d^{(j)} (\cdot))_{\NN_0 \ni j \leq r} \colon C^{r} (U,F) \rightarrow \prod_{\NN_0 \ni j \leq r} C (U\times E^{j}, F) , f \mapsto (d^{(j)}f) 
  \end{displaymath}
 into a topological embedding. We call this topology the \ind{compact-open topology!C^r topology@$C^r$-topology}{compact-open $C^r$-topology}. Notice that it is the initial topology with respect to the family $(d^{(j)} (\cdot))_{\NN \ni j \leq r}$. 
\end{defn}

\begin{rem} \label{rem: crfs:prop} \no{rem: crfs:prop} 
   \begin{compactenum}
   \item By \cite[Lemma 1.14]{hg2002a}, Definition \ref{defn: CrTop} coincides on open sets with the definition in \cite[Definition 3.1]{hg2002}. Hence if $U$ is an open subset of finite-dimensional space $E$ and $F$ is a Fr\'{e}chet space, then $C^r(U,F)$ is a Fr\'{e}chet space by \cite[Remark 3.2]{hg2002}. 
   \item For each compact subset $K \subseteq U$ and open subset $V \subseteq F$, the set \glsadd{KUr}
		\begin{displaymath}
		 \lfloor K, V\rfloor_r \coloneq \setm{\gamma \in C^r (U,F)}{\gamma (K) \subseteq V}
		\end{displaymath}
is open in $C^r (U, F)$ by \cite[Lemma 4.22]{hg2004}.\\
If $s, r \in \NN_0 \cup \set{\infty}$ with $r \leq s$, then $C^s (U,F) \subseteq C^r (U,F)$ holds by definition and the topology on $C^s (U,F)$ is finer than the subspace topology induced by $C^r (U,F)$. Let $\Omega$ be an open set in $C^s (U,F)$ such that $\Omega = C^s (U,F) \cap A$ holds for some open $A \subseteq C^r (U,F)$. Then we call \emph{$\Omega$ a $C^r$-open set in $C^s(U,F)$} or a \emph{$C^r$-neighborhood of $f \in C^s (U,F)$}, for any $f \in \Omega$.                                                      
   \end{compactenum}
 \end{rem}

\begin{defn}\label{defn: n:top}
 Let $E$ be a locally convex space and $M$ a $C^r$-manifold. Then we let $C^r(M,E)$ be the space of all $C^r$-mappings $\gamma \colon M \rightarrow E$. The pointwise operations turn $C^r(M,E)$ into a vector space. Endow $C^r (M,E)$ with the initial topology with respect to the family 
    \begin{displaymath}
     \theta_\kappa \colon C^r(M,E) \rightarrow C^r(V_\kappa, E), \gamma \mapsto \gamma|_{U_\kappa} \circ \kappa^{-1}
    \end{displaymath}
 where $\kappa \colon U_\kappa \rightarrow V_\kappa$ ranges through an atlas of $M$. The topology is independent of the choice of atlas by \cite[Lemma 4.9]{hg2004}. If $M$ is an open subset of a locally convex space, \cite[Lemma 4.6]{hg2004} proves that this topology coincides with the compact open $C^r$-topology.
\end{defn}

 \begin{defn}\label{defn: norm} \no{defn: norm}
 \begin{compactenum}
  \item  Let $U \subseteq \RR^d$ be an open subset $d \in \NN_0$ and $K \subseteq U$ compact. For $\xi \in C^r (U,\RR^d)$, $r \in \NN_0 \cup \set{\infty}$, the maximum norm $\norm{\cdot}_\infty$ and $k \in \NN_0$ with $k \leq r$, we use standard multiindex notation to set \glsadd{normKr} 
	\begin{displaymath}
		\norm{\xi}_{K,k} \coloneq \max_{\lvert \alpha\rvert \leq k} \max_{x \in K} \norm{\partial^\alpha \xi (x)}_\infty.
	\end{displaymath}  
   
  \item Let $E$ be a locally convex space and $r \in \NN_0 \cup \set{\infty}$. Endow $C^r([0,1],E)$ with the locally convex vector topology induced by the family of seminorms $\norm{\cdot}_{C^k,p}$ defined via 
    \begin{displaymath}
     \norm{\gamma}_{C^k,p} \coloneq \max_{j=0,\ldots ,k} \max_{t \in [0,1]} p \left(\frac{\partial^k}{\partial t^k} \gamma (t)\right)
    \end{displaymath}
  where $p$ ranges through the continuous seminorms on $E$ and $k \in \NN_0$ with $k \leq r$.
\end{compactenum}
 \end{defn}

\begin{rem}\no{rem: Cr:Norm}\label{rem: Cr:Norm}
 \begin{compactenum}
  \item Let $U \subseteq \RR^d$ be some open subset, where $d \in \NN_0$. As $U$ is $\sigma$-compact, there is a sequence of compact sets $(K_n)_{n \in \NN}$ such that $U = \bigcup_{n \in \NN} K_n$. By a variant of \cite[Proposition 4.4]{hg2002}, the locally convex topology induced by the family of seminorms $\setm{\norm{\cdot}_{K_n,k}}{n\in \NN , 0 \leq k \leq r}$ on $C^r(U,\RR^d)$ coincides with the compact-open $C^r$-topology.
  \item A variant of \cite[Proposition 4.4]{hg2002} shows that the topology introduced in Definition \ref{defn: norm} (b) is initial with respect to the mappings $d^{(j)} \colon C^r([0,1], E) \rightarrow C([0,1]\times \RR^j,E)_{\text{c.o}}, \gamma \mapsto d^{(j)} \gamma$, $0 \leq k \leq r$, i.e.\ it coincides with the compact-open $C^r$-topology.\\ In particular, then $C^r([0,1],U) \coloneq \setm{\gamma \in C^r ([0,1],E)}{\gamma ([0,1]) \subseteq U} = \lfloor [0,1], U\rfloor_r$ is an open subset for each open $U \subseteq E$. If $E$ is metrizable (respectively complete), $C^r ([0,1],E)$ is metrizable by \cite[2.8 Theorem 1]{jarchow1980} (respectively complete by \cite[Lemma 1.4]{gn2012}).
 \end{compactenum}
\end{rem}

\begin{nota}
 Let $U \subseteq E$ and $V \subseteq F$ be locally convex subsets with dense interior of locally convex topological vector spaces $E$ and $F$, respectively. Furthermore, let $G$ be a topological vector space and $f \colon U \rightarrow C(V,G)$ be a map. We associate to $f$ the map \glsadd{fwedge}
    \begin{displaymath}
    f^\wedge \colon U\times V \rightarrow G \text{ defined via } f^\wedge (u,v) \coloneq f(u)(v).
    \end{displaymath}
 
\end{nota}

\subsection{Spaces of sections and patched spaces}\label{sect: sections} \no{sect: sections}

In this section we endow the \ind{space!of smooth vector fields}{space of smooth vector fields $\vect{M}$}\glsadd{vect} on a smooth manifold $M$ with a topology. Furthermore, we use the concept of a \tl patched locally convex space\tr (cf.\ \cite{hg2003, hg2004}) to obtain a criterion for the differentiability of maps between spaces of sections. We recall the following facts from \cite[Appendix F]{hg2004}:

\begin{defn}\label{defn: top:vect} \no{defn: top:vect}
 Let $M$ be a smooth manifold modeled on the locally convex space $E$ and $\pi_{TM} \colon TM \rightarrow M$ be the bundle projection. Consider a maximal atlas $\aA$ of $M$ and a chart $(V_\psi , \psi) \in \aA$ with $\psi \colon V_\psi \rightarrow U_\psi$. Let $\text{pr}_2 \colon V_\psi \times E \rightarrow E$ be the canonical projection.\\ 
 For a vector field $X \in \vect{M}$, we define a local representative $X_\psi \coloneq \text{pr}_2 \circ T\psi \circ X|_{V_\psi} \colon V_\psi \rightarrow E$\glsadd{vfloc}. In particular $T\psi \circ X (y) =  (\psi (y), X_\psi (y))$ holds for all $y \in V_\psi$.\\
 We endow $\vect{M}$\glsadd{vect} with the unique locally convex topology turning the linear map 
	\begin{displaymath}
	 \Gamma \colon \vect{M} \rightarrow \prod_{(V_\psi, \psi) \in \aA} C^\infty (V_\psi , E) , \ X \mapsto (X_\psi)_{(V_\psi, \psi) \in \aA}
	\end{displaymath}
 into a topological embedding. Then the topology on $\vect{M}$ is the initial topology with respect to the family of linear maps $\theta_\psi \colon \vect{M} \rightarrow C^\infty (V_\psi , E) , X \mapsto X_\psi$.\glsadd{theta}
\end{defn}

\begin{lem}[\hspace{-0.5pt}{\cite[Lemma F.9]{hg2004}}]\label{lem: vect:top} \no{lem: vect:top}
 The topology on $\vect{M}$ is initial with respect to the family $(\theta_\phi)_{(V_\phi ,\phi) \in \bB}$, where $\bB \subseteq \aA$ is some atlas for $M$. 
\end{lem}

\begin{proof}
 Combine \cite[Lemma F.9]{hg2004} with \cite[Proposition 4.19]{hg2004}, which guarantees that the topology defined in \cite{hg2004} coincides with our definition of the compact-open $C^r$-topology over the field $\RR$.
\end{proof}

 \begin{nota}\label{nota: res}
  Let $M$ be a smooth manifold and $U$ an open subset of $M$.  We define the restriction map $\res^M_{U} \colon \vect{M} \rightarrow \vect{U}, X \mapsto X|_{U}^{TU}$.\glsadd{rest} For each open subset $U$ this map is continuous linear by \cite[Lemma F.15]{hg2004}.\footnote{The article \cite{hg2004} uses another concept of differentiability in locally convex vector spaces which is adapted to non-discrete topological fields. However as \cite[Proposition 7.4]{bgn2004} asserts, this concept of differentiability coincides with the one from Definition \ref{defn: deriv} on open sets of locally convex vector spaces over the field $\RR$. As we are only interested in this case, we may use the results of \cite{hg2004} without restriction.}
 \end{nota}

\begin{defn}\label{defn: cs:vf}
 Let $d \in \NN$. We define the \ind{space!of compactly supported vector fields}{space of compactly supported vector fields} $\mathfrak{X}_c (\RR^d)$.\footnote{Since this space is only needed in Example \ref{ex: os:good}, we shall only consider vector fields on $\RR^d$ (cf.\ \cite[Appendix F]{hg2004} for a more general definition).} The assignment $\theta \colon \mathfrak{X}_c (\RR^d) \rightarrow C^\infty_c (\RR^d,\RR^d) , X \mapsto \text{pr}_2 \circ X$ is a bijective map, where $\text{pr}_2$ denotes the canonical projection $T\RR^d \cong \RR^d \times \RR^d \rightarrow \RR^d, (x,y) \mapsto y$. We define a topology on $C^\infty_c (\RR^d,\RR^d)$ (and thus also on $\mathfrak{X}_c (\RR^d)$) turning $\theta$ into an isomorphism of topological vector spaces. Choose a locally finite cover $\uU = (U_i)_{i \in I}$ of $\RR^d$ by relatively compact open subsets $U_i \subseteq \RR^d$ such that the cover is countable. Then consider the map  
    \begin{displaymath}
     R_\uU \colon \mathfrak{X}_c (\RR^d) \rightarrow \bigoplus_{i \in I} C^\infty (U_i, \RR^d) , R_\uU (\sigma) \coloneq (\text{pr}_2 \circ \sigma|_{U_i})_{i \in I}
    \end{displaymath}
  We endow $\mathfrak{X}_c (\RR^d)$ with the unique locally convex topology induced by the linear map $R_\uU$. Here the right hand side has been endowed with the locally convex direct sum topology. By \cite[Lemma 8.10]{hg2004}, the topology constructed does not depend on the choice of covering $\uU$ (recall from \cite[Proposition 4.19]{hg2004} that the topology defined in \cite{hg2004} coincides with our definition of the compact-open $C^r$-topology over the field $\RR$). Furthermore $\mathfrak{X}_c (\RR^d)$ is a Hausdorff space and $R_\uU$ is a topological embedding with closed image by \cite[Proposition 8.13]{hg2004}. 
\end{defn}

\begin{defn}\label{defn: patloc} \no{defn: patloc}
 A \ind{}{patched locally convex space} over $\RR$ is a pair $(E, (p_i)_{i \in I})$, where $E$ is a topological $\RR$-vector space and $(p_i)_{i \in I}$ is a family of continuous linear maps $p_i \colon E \rightarrow E_i$ to topological vector spaces $E_i$ such that 
	\begin{compactenum}
	 \item for each $x \in E$, the set $\setm{i \in I}{p_i (x)\neq 0}$ is finite,
	 \item the linear map 
		\begin{displaymath}
			p \colon E \rightarrow \bigoplus_{i \in I} E_i , \ x \mapsto (p_i (x))_{i \in I} = \sum_{i \in I} p_i (x)
		\end{displaymath}
 	from $E$ to the direct sum $\oplus_{i \in I} E_i$ (equipped with the direct sum topology cf. \cite[Ch. II, \S 2, No.\ 5, Definition 2]{bourbaki1987}) is a topological embedding,
	\item the image $p(E)$ is sequentially closed in $\bigoplus_{i \in I} E_i$. 
	\end{compactenum}
 The mappings $p_i \colon E \rightarrow E_i$ are called \ind{patches}{patches}, and the family $(p_i)_{i \in I}$ is called a \ind{patchwork}{patchwork}. If $I$ is a countable set, we also say that $E$ is \ind{patched locally convex space!countably patched}{countably patched}.
\end{defn}


\begin{lem}\label{lem: pb:patch} \no{lem: pb:patch}
 Let $(E, (p_i)_{i \in I})$ be a patched topological $\RR$-vector space, with $p_i \colon E \rightarrow E_i$ and $p$ as in Definition \ref{defn: patloc}. For each $r \in \NN_0 \cup \set{\infty}$, the map 
  \begin{displaymath}
   p_* \colon C^r ([0,1] , E) \rightarrow C^r([0,1], \bigoplus_{i \in I} E_i) , g \mapsto p \circ g
  \end{displaymath}
 is a linear topological embedding whose image is sequentially closed. If $\lvert I\rvert < \infty$ or $E$ is countably patched and $r < \infty$, then the family $C^r ([0,1] , p_i) \colon C^r([0,1],E) \rightarrow C^r([0,1],E_i), \gamma \mapsto p_i \circ \gamma ,\ i  \in I$, turns $C^r([0,1], E)$ into a patched locally convex space over $\RR$.
\end{lem}

\begin{proof}
 The maps  $C^r ([0,1] , p_i)$ are continuous linear for $i \in  I$ and $p_*$ is a topological embedding by \cite[Lemma 1.2]{gn2012}. Without loss of generality we identify $E$ with a subspace of $F \coloneq \bigoplus_{i \in I} E_i$. Let $(f_n)_{n \in \NN} \subseteq \im p_*$ be a sequence which converges to some $f \in C^r([0,1], F)$. Since $E$ is sequentially closed, due to the continuity of the point evaluation maps (cf. \cite[Proposition 3.20]{alas2012}) for $t \in [0,1]$ the sequence $(f_n (t))_{n \in \NN}$ converges in $E$. Hence the image of $f$ is contained in $E$. Recall that directional derivatives may be computed as limits of sequences. As each element $f(t)$ is contained in $E$ and $E$ is sequentially closed, the mappings $d^{(k)} f$, for $\NN_0 \ni k\leq r$, take their images in $E$. Hence $f \in C^r ([0,1], E)$ holds and $\im p_*$ is sequentially closed as a subspace of $C^r([0,1], F)$.
 \paragraph{Case 1: $\lvert I \rvert < \infty$.} Since $I$ is finite, the coproduct $F \coloneq \bigoplus_{i \in I} E_i$ in the category of locally convex topological vector spaces coincides with the product of the $E_i$. Hence the canonical projection $\pi_i \colon F \rightarrow E_i$ and the canonical inclusion $\iota_i \colon E_i \rightarrow F$ are continuous linear for $i \in I$. From \cite[Lemma 1.2]{gn2012} we deduce that the mappings
 \begin{align*}
     ((\pi_i)_*)_{i \in I} \colon C^r([0,1], \bigoplus_{i \in I} E_i) \rightarrow \bigoplus_{i \in I} C^r([0,1], E_i), f \mapsto (\pi_i \circ f)_{i \in I},\\
     \bigoplus_{i \in I} C^r([0,1], E_i) \rightarrow C^r([0,1], \bigoplus_{i \in I} E_i), (f_i) \mapsto \sum_{i \in I} (\iota_i)_* (f_i)
    \end{align*}
 are continuous linear and mutually inverse. Thus $C^r([0,1], \bigoplus_{i \in I} E_i)$ and $\bigoplus_{i \in I} C^r([0,1], E_i)$ are isomorphic as locally convex spaces, whence the maps $(p_i)_*, i \in I$ form a patchwork for $C^r([0,1],E)$.
 \paragraph{Case 2: $\lvert I \rvert = \infty$ and $r < \infty$.} The canonical inclusions yield a family of continuous linear maps $((\iota_i)_*)_{i \in I}$ by \cite[Lemma 1.2]{gn2012}. As in the first case we obtain a linear and continuous map $\Lambda \colon \bigoplus_{i \in I} C^r ([0,1],E_i) \rightarrow C^r([0,1],F), (\gamma_i)_{i \in I} \mapsto \sum_{i \in I} (\iota_i)_* (\gamma_i)$. For the rest of the proof, we suppress the inclusions $\iota_i$ in the notation. To prove our claim, we have to construct an inverse mapping for $\Lambda$. To do so, pick $\gamma \in C^r([0,1],F)$. The compact set $\gamma ([0,1]) \subseteq F$ is contained in a finite partial sum by \cite[Ch.\ III, \S 4, No.\ 1, Proposition 5]{bourbaki1987}. As the inclusion of a finite partial sum is a topological embedding with closed image, from \cite[Lemma 1.2]{gn2012} and the isomorphism established for the finite case, we deduce that there are unique $\gamma_i  \in C^r ([0,1], E_i)$, for $i \in I$ with $\gamma = \Lambda ((\gamma_i)_{i \in I})$. Hence we obtain a well-defined inverse of $\Lambda$ via $\Theta \colon C^r ([0,1],F) \rightarrow \bigoplus_{i \in I} C^r ([0,1], E_i), \gamma \mapsto (\gamma_i)_{i \in I}$.\\
 We claim that $\Lambda$ is an isomorphism of locally convex spaces. To prove the claim, let $\Gamma_i$ be the set of all continuous seminorms on $E_i$. Consider $q = (q_i)_{i \in I}  \in \Gamma \coloneq \prod_{i \in I} \Gamma_i$ and obtain a continuous seminorm $r_q \colon F \rightarrow [0,\infty[, r_q (\sum_{i \in I} x_i) \coloneq \sup \setm{q_i (x_i)}{i \in I}$ with $x_i \in E_i$. Since the space $E$ is countably patched, the topology on $F$ coincides with the box topology by \cite[Proposition 4.1.4]{jarchow1980}. Hence the family $(r_q)_{q\in \Gamma}$ determines the locally convex topology on $F$. By definition of the topology on $C^r ([0,1], F)$, the continuous seminorms $s_q \colon C^r ([0,1], F) \rightarrow [0,\infty[$,
    \begin{displaymath}
     s_q (\gamma ) \coloneq \sup_{0\leq k \leq r} \sup_{x \in [0,1]} r_q (\tfrac{\partial^k}{\partial x^k} \gamma (x)) = \sup_{0\leq k \leq r} \sup_{x \in [0,1]} \sup_{i \in I} q_i (\tfrac{\partial^k}{\partial x^k} \gamma_i (x)),
    \end{displaymath}
 determine the locally convex topology on $C^r ([0,1], F)$ for $q$ ranging through $\Gamma$. Likewise, the locally convex topology on $C^r([0,1],E_i)$ is determined by the continuous seminorms $t_{q_i} \colon C^r ([0,1], E_i) \rightarrow [0,\infty[ , t_{q_i} (\gamma_i ) \coloneq \sup_{0\leq k \leq r} \sup_{x\in [0,1]} q_i (\tfrac{\partial^k}{\partial x^k} \gamma_i (x))$, where $q_i$ ranges through $\Gamma_i$. The locally convex sum topology, i.e.\ the box topology on $\bigoplus_{i \in I} C^r ([0,1], E_i)$, is induced by the family of seminorms $u_q \colon \bigoplus_{i \in I} C^r ([0,1], E_i) \rightarrow [0,\infty[,$ 
  \begin{displaymath}
   u_q ((\gamma_i)_{i \in I} \coloneq \sup_{i \in I} t_{q_i} (\gamma_i) = \sup_{i \in I} \sup_{ 0 \leq k \leq r} \sup_{x \in [0,1]} q_i (\tfrac{\partial^k}{\partial x^k} \gamma_i (x))
  \end{displaymath}
 for $q = (q_i)_{i \in I} \in \Gamma$. Observe that for each $q \in \Gamma$, we have $s_q \circ \Lambda = u_q$. We deduce that $\Lambda^{-1}$ is continuous (cf. \cite[Ch. II, \S 2, No.\ 4, Proposition 4]{bourbaki1987}), whence $\Lambda$ is an isomorphism of locally convex spaces.
\end{proof}

If $r = \infty$ and $\lvert I \rvert = \infty$, the map $\Lambda$ introduced in the proof of Lemma \ref{lem: pb:patch} still is a continuous linear bijection, but its inverse fails to be continuous in general. 

\begin{defn}
 Let $I$ be a set and $(E,(p_i)_{i \in I})$ and $(F,(q_i)_{i \in I})$ patched locally convex $\RR$-vector spaces with canonical embeddings $p \colon E \rightarrow \bigoplus_{i \in I} E_i$ and $q\colon F \rightarrow \bigoplus_{i \in I} F_i$ as in Definition \ref{defn: patloc}. 
	\begin{compactenum}
	 \item A map $f \colon U \rightarrow F$ defined on an open subset $U \subseteq E$ is called a \ind{}{patched mapping} if there exists a family $(f_i)_{i \in I}$ of mappings $f_i \colon U_i \rightarrow F_i$ on certain open neighborhoods $U_i$ of $p_i (U)$ in $E_i$, which is \textit{compatible with $f$} in the following sense: We have $0 \in U_i$ and $f_i (0) = 0$ for all but finitely many $i$, and $q_i (f(x)) = f_i (p_i(x))$ for all $i \in I$, i.e.\ $q \circ f = (\bigoplus_{i \in I} f_i) \circ p|_U^{\oplus U_i}$.
	 \item For $k \in \NN_0 \cup \set{\infty}$, we say that a patched mapping $f \colon U \rightarrow F$ is \textit{of class $C^k$ on the patches} if all of the mappings $f_i$ in (a) can be chosen of class $C^k$.
	\end{compactenum}
\end{defn}

\begin{prop}\label{prop: pat:loc} \no{prop: pat:loc}
 Let $I$ be a set and $(E,(p_i)_{i \in I})$, $(F,(q_i)_{i \in I})$ be patched topological $\RR$-vector spaces. Assume that $f \colon U \rightarrow F$ is a patched mapping from an open subset $U \subseteq E$ to $F$. If $f$ is of class $C^{k+1}$ on the patches, then $f$ is of class $C^k$. If $E$ and $F$ are countably patched and $f$ is $C^k$ on the patches, then $f$ is of class $C^k$. 
\end{prop}

\begin{proof}
 For $i \in I$, let $f_i \colon U_i \rightarrow F_i$ be the mappings compatible with $f$. Consider the box neighborhood $\oplus_{i \in I} U_i \coloneq \left(\prod_{i \in I} U_i\right) \cap \left(\bigoplus_{i \in I} E_i \right)$ which is open in the locally convex sum (cf.\ \cite[4.3]{jarchow1980}). The compatibility condition yields $q \circ f = (\bigoplus_{i \in I} f_i) \circ p|_U^{\oplus U_i}$. As shown in \cite[Proposition 7.1]{hg2003}, the map $\bigoplus_{i \in I} f_i$ is a $C^k$-map if each $f_i$ is of class $C^{k+1}$ (respectively if each $f_i$ is a $C^k$-map and $I$ is countable). By definition, this is the case if and only if $f$ is $C^{k+1}$ (respectively $C^k$ in the countable case) on the patches. The map $(\bigoplus_{i \in I} f_i) \circ p|_U^{\oplus U_i}$ is of class $C^k$ as a composition of a $C^k$-map and a smooth map. Thus $q \circ f$ is a $C^k$-map. Since the subspace $\im q$ is sequentially closed, the corestriction $(q\circ f)|^{\im q}$ is a $C^k$ map. As $q|^{\im q}$ is an isomorphism of topological vector spaces, $f$ is a $C^k$-map. 
\end{proof}
\subsection{Lie groups}

\begin{defn}
 A \ind{Lie group}{(locally convex) Lie group} is a group $G$ equipped with a smooth manifold structure(modeled on a locally convex space) turning the group operations into smooth maps. Denote its neutral element by $\one$ and recall that $L(G) \coloneq T_\one G$ is its Lie algebra (cf. \cite{hg2002a, neeb2006} for details).
\end{defn}

\begin{defn}
 Let $G$ be a Lie group. We denote by $\rho_g \colon G \rightarrow G, h \mapsto hg$ the right translation by $g\in G$. This yields a natural right action of $G$ on the tangent Lie group $TG$ (cf. \cite[Ch. III, \S 2]{bourbaki1989}):
  \begin{displaymath}
   v\cdot g \coloneq (T_x \rho_g )(v) \in T_{xg} G \quad \text{for  } x \in G , v \in T_x G.
  \end{displaymath}
\end{defn}

The following construction principle for Lie groups will be our main tool to construct Lie group structures (cf.\ \cite[Ch.\ III, \S 1, No.\ 9, Proposition 18]{bourbaki1989}).

\begin{prop}\label{prop: Lgp:locd} \no{prop: Lgp:locd}
 Let $G$ be a group and $U, V$ subsets of $G$ such that $\one \in V = V^{-1}$ and $V\cdot V \subseteq U$. Suppose that $U$ is equipped with a smooth manifold structure modeled on a locally convex space such that $V$ is open in $U$ and which turns $\iota \colon V \rightarrow V \subseteq U$ and $\mu \colon V \times V \rightarrow U$ - the mappings induced by inversion and the group multiplication respectively - into smooth maps. Then the following holds:
 \begin{compactitem}
  \item[\rm (a)] There is a unique smooth manifold structure on the subgroup $G_0 \coloneq \langle V\rangle$ of $G$ generated by $V$ such that $G_0$ becomes a Lie group, $V$ is open in $G_0$, and such that $U$ and $G_0$ induce the same smooth manifold structure on the open subset $V$.
  \item[\rm (b)] Assume that for each $g$ in a generating set of $G$, there is an open identity neighborhood $W \subseteq U$ such that $gWg^{-1} \subseteq U$ and  $c_g \colon W \rightarrow U, h \mapsto ghg^{-1}$ is smooth. Then there is a unique smooth manifold structure on $G$ turning $G$ into a Lie group such that $V$ is open in $G$ and both $G$ and $U$ induce the same smooth manifold structure on the open subset $V$.
 \end{compactitem}
\end{prop}

\subsection{Regular Lie groups}\label{app: reg}

\begin{defn}\label{defn: pI:lld} \no{defn: pI:lld}
 Let $G$ be a Lie group with Lie algebra $L(G)$. Consider a $C^k$-curve $p \colon [0,1] \rightarrow G$ with $k \geq 1$ and recall that 
  \begin{displaymath}
   \delta^r p \in C^{k-1} ([0,1],L(G)), (\delta^r p)(t) \coloneq p'(t) \cdot p(t)^{-1}
  \end{displaymath}
is called the \ind{logarithmic derivative, right}{right logarithmic derivative}\glsadd{rld} of $p$. Furthermore we call $p$ a \ind{product integral, right}{right product integral} for $\delta^r p$.\\
If $q \colon [0,1] \rightarrow G$ is another $C^k$-curve such that $\delta^r p = \delta^r q$ (i.e. both $p$ and $q$ are right product integrals for $\delta^r q$), then $q =p \cdot g_0$ holds for some constant $g_0 \in G$ (cf. \cite[Lemma 7.4]{milnor1983}).
\end{defn}

\begin{defn}\label{defn: nota:PI} \no{defn: nota:PI}
 If $\gamma \in C^k ([0,1],L(G))$ with $k\in \NN_0 \cup \set{\infty}$ admits a right product integral $p$, we define $\pP (\gamma ) \coloneq p \cdot p(0)^{-1} $. Thus $\pP (\gamma)$ is a right product integral for $\gamma$ such that $\pP (\gamma) (0) = \one_G$ is the identity element of $G$. The product integral is uniquely determined by this property.\glsadd{lPI}
\end{defn}


\begin{defn}\label{defn: Ck:reg} \no{defn: Ck:reg}
 Let $k \in \NN_0 \cup \set{\infty}$. A Lie group $G$ with Lie algebra $L(G)$ is called \ind{Lie group!Ckregular@$C^k$-regular}{(strongly) $C^k$-regular}, if for each $\xi \in C^k ([0,1],L(G))$, the initial value problem 
  \begin{equation}
   \gamma (0) = \one_G ,\quad \quad \delta^r (\gamma) = \xi
  \end{equation}
 has a solution $\lPI{\xi}$, which is then contained in $C^{k+1} ([0,1],G)$, and the corresponding evolution map\glsadd{evol} 
  \begin{displaymath}
   \evol_G \colon C^k([0,1],L(G)) \rightarrow G, \xi \mapsto \lPI{\xi} (1)
  \end{displaymath}
 is smooth. If $G$ is $C^k$-regular, we write \glsadd{Evol}
  \begin{displaymath}
   \Evol_G \colon C^k ([0,1], L(G)) \rightarrow C^{k+1} ([0,1],G),  \xi \mapsto \lPI{\xi}
  \end{displaymath}
 for the map on the level of Lie group-valued curves. For more information on regularity see \cite{hg2015}.\\
 The group $G$ is called \ind{Lie group!regular (in the sense of Milnor)}{regular (in the sense of Milnor)} if it is $C^\infty$-regular. For $k \leq r$ the $C^r$-regularity follows from $C^k$-regularity. 
\end{defn}

Notice that we have defined regularity properties of Lie groups using the right logarithmic derivative. Alternatively one may define \emph{left logarithmic derivative}, \emph{left product integrals} and regularity properties using these notions. However, it is well known that this results in the same concepts of regularity as defined in \ref{defn: Ck:reg}. See \cite[Proposition 1.3.6]{Dahmen2012} for a proof. 

The following lemma will be our main tool to prove the regularity of the orbifold diffeomorphism group. Its proof carries over almost verbatim from \cite[Proposition 1.3.10]{Dahmen2012}:

\begin{lem}\label{lem: znbhd:reg} \no{lem: znbhd:reg}
 Let $G$ be a smooth Lie group with Lie Algebra $L(G)$. Assume that there is a zero-neighborhood $U\subseteq C^k([0,1], L(G))$ for $k \in \NN_0 \cup \set{\infty}$ such that every $\xi \in U$ has a right product integral. Furthermore assume that $E_1 \colon U \rightarrow G , \xi \mapsto \lPI{\xi}(1)$ is smooth. Then $G$ is $C^k$-regular.
\end{lem}\thispagestyle{empty}
\section[Riemannian Geometry: Supplementary Results]{Riemannian geometry: Supplementary Results}\label{sect: barigeo} \no{sect: barigeo}
\setcounter{subsubsection}{0}
In this thesis we assume some basic familiarity with Riemannian metrics and geodesics. Our approach also requires standard results from Riemannian geometry as outlined in \cite{rg1992,fdiffgeo1963,klingenberg1995}. The results obtained in this section are a variation of ideas first developed in \cite{hg2006b}. Our goal is to fix the necessary notation and to provide estimates needed in the proof of the main theorems.  

\begin{Anota}
 The pair $(M, \rho_M)$ will always denote a finite dimensional smooth Riemannian manifold $M$, with Riemannian metric $\rho_M$. Notice that for each $x \in M$ the Riemannian metric yields a positive definite inner product $\rho_{M,x} \colon T_xM \times T_xM \rightarrow \RR$. We usually abbreviate 
    \begin{displaymath}
     \rho_M (X,Y) \coloneq \rho_{M,x} (X,Y) \quad \forall X,Y\in T_x M.
    \end{displaymath}
 We define the \ind{Riemannian metric!metric ball}{$\ve$-balls with respect to the Riemannian metric in $T_xM$ around the origin $0_x$} as $B_{\rho_M} (0_x , \ve) \coloneq \setm{X \in T_xM}{\sqrt{\rho_M (X,X)}<\ve}$.\glsadd{ballrm}
Recall that on every Riemannian manifold there exists a Riemannian exponential map 
	\begin{displaymath}
	 \exp_M \colon TM \supseteq D_M \rightarrow M 
	\end{displaymath}
whose domain $D_M$ is an open neighborhood of the zero-section. Each Riemannian exponential map on a smooth Riemannian manifold is smooth.
\end{Anota}

 Recall the following standard result of Riemannian geometry:   

\begin{Alem}\label{lem: rg:rmb} \no{lem: rg:rmb}
 Let $(M,\rho)$ be a Riemannian manifold with exponential map $\exp_M \colon D_M \rightarrow M$ and $K \subseteq M$ be a compact subset. There is $\ve>0$ and an open set $V\subseteq M$ containing $K$ such that the following holds \begin{compactenum}
     \item for each $x \in V$, the map $\exp_M|_{B_\rho (0_x , \ve)}^{\exp_M (B_{\rho} (0_x ,\ve))}$ is a diffeomorphism with open image in $M$,
     \item $\bigcup_{x \in V} B_\rho (0_x , \ve) \subseteq D_M$ is an open neighborhood of the zero section on $K$.                                                                                                                                                                                    
 \end{compactenum}
\end{Alem}
 
\begin{proof}
 Apply \cite[Theorem 1.8.15]{klingenberg1995} to each point $x \in K$. Since $K$ is compact, this yields a finite family $x_1, x_2, \ldots,x_n \in K$ and constants $\ve_1,\ldots , \ve_n$ such that: 
  \begin{compactitem}
   \item for each $1\leq k \leq n$ and $y \in \exp_M (B_{\rho} (0_{x_k},\ve_k))$, the mapping $\exp_M|_{B_{\rho} (0_{y},\ve_k)}$ is an embedding with open image,
   \item $K \subseteq V \coloneq \bigcup_{1\leq k\leq n} \exp_{M} (B_{\rho} (0_{x_k},\ve_k))$ holds.
  \end{compactitem}
 Set $\ve \coloneq \min \set{\ve_1, \ldots , \ve_n}$. The pair $(\ve,V)$ satisfies the assertion of the lemma since $\bigcup_{x \in V} B_\rho (0_x , \ve)$ is an open neighborhood of the zero section by the proof of \cite[Theorem 1.8.15]{klingenberg1995}.
\end{proof}

For the rest of this section, we endow $\RR^d$ (for $d\in \NN$) with the maximum norm $\norm{\cdot}_\infty$. We denote by $B_r (x)$ the metric ball around $x\in \RR^d$ with respect to $\norm{\cdot}_\infty$ and radius $r>0$. As a first step we discuss Riemannian exponential maps on metric balls in euclidean space. To this end, fix the metric ball $B_5 (0) \subseteq \RR^d, d \in \NN$ and endow it with an arbitrary Riemannian metric. 

\begin{Alem}\label{lem: exp:loc} \no{lem: exp:loc}
 Consider $B_5 (0)$ as a Riemannian manifold with arbitrary Riemannian metric. Let $\exp \colon D \rightarrow B_5 (0)$ be the associated Riemannian exponential map. There exist $\ve > 0$ and $1 > \delta >0$ such that 
	\begin{compactenum}
	 \item $\overline{B_4 (0)} \times \overline{B_\ve (0)} \subseteq D$ and $\phi_x \coloneq \exp (x ,\cdot )|_{B_\ve (0)}^{\exp (x,\cdot) (B_\ve (0))}$ is a diffeomorphism for each $x \in \overline{B_4(0)}$.
	 \item $B_\delta (x) \subseteq \exp (x , B_\ve (0))$ for each $x \in \overline{B_4 (0)}$ and  $b \colon W_\delta \rightarrow B_\ve (0), b(x,y) \coloneq \phi_x^{-1} (y)$ is a smooth map on the subset $W_\delta \coloneq \bigcup_{x \in \overline{B_4 (0)}} \set{x} \times B_\delta (x)$ of $B_5 (0) \times \RR^d$.
	 \item For each $t >0$, there exists $\sigma_t \in \, ]0,\ve[$ such that $\phi_x (B_{\sigma_t} (0)) \subseteq B_t (x)$ for each $x \in \overline{B_4 (0)}$.
    	\end{compactenum}
If $t \leq\frac{\delta}{2}$ holds in (c), we obtain a smooth map 
	\begin{displaymath}
	 f \colon B_3 (0) \times B_{\sigma_t} (0) \times B_{\sigma_t} (0) \rightarrow B_\ve (0), f(x,y,z) \coloneq b(x, \phi_{\phi_x(y)} (z)).
	\end{displaymath}
\end{Alem}
  
\begin{proof}
 \begin{compactenum}
  \item The set $\overline{B_4 (0)} \times \set{0}$ is a compact subset of $D$. Lemma \ref{lem: rg:rmb} yields an open neighborhood $\overline{B_4 (0)} \times \set{0} \subseteq W \subseteq D$, such that $\exp (x, \cdot)$ restricts to is a diffeomorphism on $W\cap T_xM$ for each $x \in \pi_{TB_5(0)} (W)$. An application of Wallace Lemma \cite[3.2.10]{Engelking 1989} yields $\ve>0$ such that $\overline{B_4(0)} \times \overline{B_\ve (0)} \subseteq W$ holds.
  \item For fixed $x \in \overline{B_4 (0)}$, we have $d_2 \exp (x,0; \cdot) = \id_{\RR^d}$ (cf. \cite[Proof of Theorem 1.6.12]{klingenberg1995}). Apply the parameter dependent Inverse Function Theorem \cite[Theorem 5.13]{hg2007} to the exponential map on $\overline{B_4 (0)} \times B_\ve (0)$. By compactness of $\overline{B_4 (0)}$, this yields some $\delta >0$ which satisfies the assertion of (b). Note that $W_\delta$ is relatively open in $\overline{B_4 (0)} \times \RR^d$ and thus a locally convex subset of $\RR^d \times \RR^d$ with dense interior.
  \item By uniform continuity of $\exp$ on $\overline{B_4 (0)} \times \overline{B_\ve (0)}$, we may choose $\sigma_t$ with the desired properties. If $t\leq\frac{\delta}{2}$ holds, we obtain $\phi_{\phi_x (y)} (z) \in B_\delta (x)$ for each $(x,y,z) \in B_3 (0) \times B_{\sigma_{t}} (0) \times B_{\sigma_{t}} (0)$. The assertion now follows from (b).
 \end{compactenum}
\end{proof}

The mappings defined in the last lemma will be used to obtain estimates for the growth of metric balls if certain maps are applied to these balls. We are interested in the composition of suitable vector fields on $B_5 (0)$ with the Riemannian exponential map. Recall that canonical lifts of orbisections are vector fields and lifts of the Riemannian orbifold exponential map are typically Riemannian exponential maps of the charts. Hence the following estimates describe the local behavior of a composition of such lifts. Moreover, the computations will enable us to control the composition of orbisections and the Riemannian orbifold exponential map.\\
In the proof of the next Lemma we use the space $\lL (\RR^d)$\glsadd{linmor} of linear and continuous endomorphisms of $\RR^d$. For the rest of this section we endow the space $\lL (\RR^d)$ with the operator norm $\opnorm{\cdot}$ with respect to $\norm{\cdot}_\infty$.  

\begin{Alem}\label{lem: mb:exp} \no{lem: mb:exp}
 Consider $B_5 (0)$ as Riemannian manifold with arbitrary Riemannian metric and exponential map $\exp$. Let $\ve , \delta$, and $D$ be as in Lemma \ref{lem: exp:loc}, and $\rho >0$. There exists an open $C^1$-neighborhood $\nN$ of the zero map in $C^\infty (B_5(0) , \RR^d)$ such that each $\xi \in \nN$ satisfies 
	\begin{compactenum}
	 \item $(\id_{B_5 (0)} , \xi) (\overline{B_3 (0)}) \subseteq \overline{B_3 (0)} \times B_{\ve} (0) \subseteq D$ and the estimate $\norm{\exp(x,\xi(x)) -x}_\infty \leq \min \set{\frac{1}{8}, \frac{\delta}{2}}$ holds for each $x \in B_3 (0)$,
	 \item the map $F_\xi \coloneq \exp \circ (\id_{B_3 (0)}, \xi|_{B_3 (0)})$ is an \'{e}tale embedding,
	 \item for $y \in B_3 (0)$, the following estimates are available:
	\begin{align}\label{eq: sz:expim}
	 B_{\frac{3}{4} s} (F_\xi (y)) \subseteq F_\xi (B_s (y)) \subseteq B_{\frac{5}{4}s} (F_\xi (y)) ,& & s \in ]0,3-\norm{y}] \\
	 B_{\frac{6s -1}{8}} (0) \subseteq F_\xi (B_s (0)) \subseteq B_{\frac{10s+1}{8}}(0)		,& & s \in ]0,3]\label{eq: exp:zero}\\
	 B_{\frac{8r-1}{10}} (0) \subseteq F_\xi^{-1} (B_r (0)) \subseteq B_{\frac{8r+1}{6}} (0)	,& & r \in \left]0, 2 +\tfrac{1}{8}\right] \label{eq: exp:inv}
	\end{align}
	\item there is a map $\xi^* \in C^\infty (\im F_\xi , \RR^d)$ such that $(F_\xi)^{-1} = \exp \circ (\id_{\im F_\xi} , \xi^*)$ is satisfied,
	\item $\norm{\xi^*}_{\overline{B_2 (0)},1} < \rho$ holds for each $\xi \in \nN$ and if $\xi \equiv 0$, then $\xi^* \equiv 0$,
  	\item the map 
	\begin{displaymath}
	 I \colon \nN \rightarrow C^\infty  (B_2 (0) , \RR^d) , \xi \mapsto \xi^*|_{B_2(0)}
	\end{displaymath}
  	is smooth. 
 \end{compactenum}
\end{Alem}

\begin{proof}
 We need preparatory estimates to control the derivatives of all relevant maps. \\
 Since $\ve, \delta $ were chosen as in Lemma \ref{lem: exp:loc}, we may consider the smooth map 
	\begin{displaymath}
	 a \colon B_4 (0) \times B_\delta (0) \rightarrow B_\ve (0), a(x,y) \coloneq b(x,x+y) = \phi_x^{-1} (x+y).
	\end{displaymath}
 Since $\exp (x,0) = x$ holds, we derive $a(x,0) = 0$ for each $x \in B_4(0)$. Thus $d_1 a(x,0;\cdot) = 0$ holds for all $x \in B_4(0)$. The set $\overline{B_3 (0)} \times \set{0} \subseteq a^{-1} (B_\rho (0))$ is compact, whence the Wallace Lemma \cite[3.2.10]{Engelking1989} allows us to choose $0< t\leq \min \set{\frac{1}{8},\frac{\delta}{2}}$ with 
 \begin{align}\label{eq: a1}
  a(\overline{B_3(0)} \times B_t (0)) \subseteq B_\rho (0) \text{ and}\\ \label{eq: a2}
  \opnorm{d_1 a(x,y;\cdot)} < \frac{\rho}{2} \quad \text{for} (x,y) \in \overline{B_3(0) \times B_t(0)} .
 \end{align}
 Set $m \coloneq \sup \setm{\opnorm{d_2a(x,y;\cdot)}}{x \in \overline{B_3 (0)}, y \in \overline{B_t(0)}} < \infty$. It is well known that the invertible matrices form an open subset $\lL (\RR^d)^\times$ of $\lL (\RR^d)$ and inversion is continuous on this set (cf. \cite[Proposition 1.33]{hg2007}). Hence there is $0 < \gamma < \frac{1}{4}$ such that for $A \in \lL (\RR^d)$ with $\opnorm{A -\id_{\RR^d}} < \gamma$ and thus $A \in \lL (\RR^d)^\times$, we have $\opnorm{A^{-1} -\id_{\RR^d}} < \frac{\rho}{2\cdot (m+1)}$.\\ 
 By Lemma \ref{lem: exp:loc}, we may choose $\sigma_t >0$ with respect to $\ve$ and $\delta$ such that $\ve > \sigma_t$ and $\phi_x (B_{\sigma_t} (0)) \subseteq B_t(x) \subseteq B_{\min \set{\frac{1}{8} , \frac{\delta}{2}}} (x)$ for each $x \in \overline{B_4 (0)}$. We obtain an open neighborhood of the zero-map $\lfloor \overline{B_3 (0)} , B_{\sigma_t} (0) \rfloor \subseteq C(B_5 (0) ,\RR^d)_{\text{c.o}}$ and by construction each $\xi \in \lfloor \overline{B_3 (0)} , B_{\sigma_t} (0) \rfloor$ satisfies the assertions of (a). We shrink $\lfloor \overline{B_3 (0)} , B_{\sigma_t} (0) \rfloor$ to construct $\nN$:\\
 For $\xi \in \lfloor \overline{B_3 (0)} , B_{\sigma_t} (0) \rfloor \cap C^\infty (B_5(0), \RR^d)$, we define the smooth maps $F_\xi \coloneq \exp \circ (\id_{B_3 (0)} , \xi|_{B_3 (0)})$ and $g_\xi \coloneq F_\xi - \id_{B_3 (0)}$. Our goal is to apply a quantitative version of the Inverse Function Theorem for Lipschitz continuous maps (cf.\ \cite[Theorem 5.3]{hg2007}). From \cite[Lemma 1.9]{hg2002a}, we deduce that the assignment $B_3 (0) \rightarrow \BoundOp{\RR^d} , x \mapsto dg_\xi (x,\cdot)$ is well defined and continuous. Since the domain of $g_\xi$ is convex, an estimate for $\opnorm{dg_\xi(z, \cdot)}$ will yield a Lipschitz constant for $g_\xi$: 
	\begin{align*}
	 dg_\xi (z;\cdot) &= d(F_\xi - \id_{B_3 (0)})(z;\cdot) = dF_\xi (z;\cdot) - \id_{\RR^d} \\
		  &= \underbrace{d_1\exp (z , \xi (z);\cdot) - \id_{\RR^d}}_{T_I (z)} + \underbrace{d_2\exp (z, \xi (z); d\xi (z;\cdot))}_{T_{II} (z)}, \quad z \in B_3 (0).
	\end{align*}
 The map $F \colon B_4 (0) \times B_\ve (0)  \rightarrow \lL (\RR^d) , (z,w) \mapsto d_1\exp (z,w;\cdot) - \id_{\RR^d}$ is continuous by \cite[Lemma 3.13]{hg2007} with $F(x,0) = 0$ for $x \in \overline{B_3 (0)}$. Using the Wallace Lemma as above, we find $s \in ]0,\ve]$ such that $F(\overline{B_3(0)}\times B_s (0)) \subseteq B_{\frac{\gamma}{2}}^{\opnorm{\cdot}} (0)$. Then $W_1 \coloneq \lfloor \overline{B_3 (0)}, B_s (0)) \rfloor \subseteq C (B_5 (0), \RR^d)_{\text{c.o.}}$ is an open neighborhood of the zero-map. For each $\xi \in \lfloor \overline{B_3 (0)} , B_{\sigma_t} (0) \rfloor \cap W_1 \cap C^\infty (B_5(0), \RR^d)$ and $x \in \overline{B_3 (0)}$, we derive $\opnorm{T_I (x)} \leq \frac{\gamma}{2} \leq \frac{1}{8}$.\\
 Since $\overline{B_3(0)} \times \overline{B_\ve (0)}$ is compact, there is an upper bound $\opnorm{d_2\exp (x,y;\cdot)} \leq C < \infty$ independent of $(x,y) \in \overline{B_3 (0)} \times \overline{B_\ve (0)}$. For each $\xi \in \lfloor \overline{B_3 (0)} , B_{\sigma_t} (0) \rfloor \cap W_1$ and $x \in \overline{B_3 (0)}$ we obtain the estimate $\opnorm{T_{II}(x)} \leq C \opnorm{d\xi (x;\cdot)}$.\\
 The topology on $C^\infty (B_5 (0), \RR^d)$ is initial with respect to the family of mappings $(d^{(k)})_{k \in \NN_0}$ by Definition \ref{defn: CrTop}. Thus we obtain an open $C^1$-neighborhood of the zero-map in $C^\infty (B_5(0), \RR^d)$ via
  \begin{displaymath}
		W_2 \coloneq \setm{\xi \in C^\infty (B_5 (0), \RR^d)}{d^{(1)} \xi \in \lfloor \overline{B_3 (0)} \times \overline{B_1 (0)} , B_{\frac{\gamma}{2C}} (0)\rfloor }.
  \end{displaymath}
 Define the $C^1$-neighborhood $\nN$ as $\nN \coloneq \lfloor \overline{B_3 (0)} , B_{\sigma_t} (0) \rfloor \cap W_1 \cap W_2$.
 For each $\xi \in \nN$, the construction shows $\text{Lip} (g_\xi) = \sup_{\norm{z}_\infty \leq 3} \opnorm{dg_\xi(z;\cdot)} \leq \gamma \leq \frac{1}{4}$.\\ 
 Since $\text{Lip} (g_\xi) < 1 = \frac{1}{\opnorm{\id_{\RR^d}}}$, the Lipschitz Inverse Function Theorem \cite[Theorem 5.3]{hg2007} yields: For $\xi \in \nN$, the map $F_\xi$ is a homeomorphism onto its image and \eqref{eq: sz:expim} is satisfied. Specializing \eqref{eq: sz:expim} to $y = 0$ together with (a) yields \eqref{eq: exp:zero}. Apply $F_\xi^{-1}$ to \eqref{eq: exp:zero} to obtain \eqref{eq: exp:inv}. We claim that $F_\xi$ is an \'{e}tale embedding. If this is true, (b) holds. To prove the claim, note that for each $z \in B_3 (0)$, one has $\frac{1}{4} \geq \opnorm{dg_\xi(z;\cdot)} = \opnorm{dF_\xi (z;\cdot) - \id_{\RR^d} (\cdot)}$. Hence $dF_\xi (z;\cdot)$ is in $\BoundOp{\RR^d}^\times$ for each $z \in B_3 (0)$. The Inverse Function Theorem (see \cite[I,4 Theorem 5.2]{langdgeo2001}) implies that $F_\xi$ is a local diffeomorphism and since it is already a homeomorphism onto its image, $F_\xi$ is an \'{e}tale embedding.\\ 
 We now prove the assertions (d)-(f). To this end, observe that by (c), the image of $F_\xi$ satisfies $B_{2+\frac{1}{8}} (0) \subseteq \im F_\xi \subseteq B_4 (0)$. Choose $x \in \im F_\xi$ and set $y \coloneq F_\xi^{-1} (x) \in B_3 (0)$. By construction of $\nN$, we have $\xi (y) \in B_{\sigma_t} (0)$, whence 
	\begin{equation}\label{eq: rho1}
	 x = F_\xi (y) = \phi_y (\xi (y)) \in B_t (y) \subseteq B_{\frac{\delta}{2}} (y)
	\end{equation}
 and thus $y \in B_t (x)$ holds. We may thus define $\xi^* (x) \coloneq b(x, F_\xi^{-1} (x))$ and obtain a smooth map $\xi^* \colon \im F_\xi \rightarrow \RR^d$ with $\im \xi^* \subseteq B_\ve (0)$. From the above estimates, we deduce that $h_{\xi^*} \coloneq \exp \circ (\id_{\im F_\xi} , \xi^*)$ is defined. A computation with $z \in B_3 (0)$ then shows 
  \begin{displaymath}
   h_{\xi^{*}} \circ F_\xi (z) = \exp (F_\xi(z), \xi^* ( F_\xi(z))) = \phi_{F_\xi(z)}(\xi^* F_\xi (z)) = \phi_{F_\xi(z)} (\phi_{F_\xi(z)}^{-1} F_\xi^{-1} (F_\xi (z))) = z.  
  \end{displaymath}
 Hence (d) holds. Notice that by construction $\xi^* (x) = a(x,(F_\xi)^{-1} (x) -x)$ for $x \in \im F_\xi$. In particular, if $\xi \equiv 0$, then $F_\xi = \id_{B_3(0)}$, whence $\xi^* (x) = a (x, F_\xi^{-1} (x)-x) = a(x,0)=0$. To obtain the estimate for (e), we computes the derivative: 
	\begin{equation}\label{eq: aabs}
	 d\xi^* (x; \cdot) = d_1 a(x, (F_\xi)^{-1} (x) -x; \cdot ) + d_2 a(x, (F_\xi)^{-1} (x) -x; d(F_\xi^{-1}) (x;\cdot) -\id_{\RR^d} (\cdot)). 
	\end{equation}
 By construction, we have $d(F_\xi^{-1}) (x;\cdot) = (dF_\xi (y;\cdot))^{-1}$ with $y \coloneq F_\xi^{-1} (x)$. By definition of $\nN$, $\opnorm{dF_\xi (y,\cdot) -\id_{\RR^d}} \leq \gamma$ and we derive $\opnorm{(dF_\xi (y;\cdot))^{-1} -\id_{\RR^d}} < \frac{\rho}{2\cdot (m+1)}$.\\
 Let $x \in \overline{B_2 (0)}$. Since $F_\xi^{-1} (x) -x \in B_t (0)$ by \eqref{eq: rho1}, the operator norm of the second summand in \eqref{eq: aabs} is smaller than $m \cdot \frac{\rho}{2(m+1)} < \frac{\rho}{2}$. Likewise, a combination of \eqref{eq: rho1} and \eqref{eq: a2} yields that the operator norm of the first summand is less than $\frac{\rho}{2}$. Summing up, $\opnorm{d\xi^* (x;\cdot)} < \rho$ holds for each $x \in \overline{B_2 (0)}$. As the operator norms on the compact set $\overline{B_2 (0)}$ were constructed with respect to $\norm{\cdot}_\infty$, we derive $\sup_{\lvert \alpha \rvert =1} \norm{\partial^\alpha \xi^*}_{\overline{B_2 (0)},0} \leq \sup_{x \in \overline{B_2 (0)}} \opnorm{d\xi^* (x;\cdot)} < \rho$. Moreover, by \eqref{eq: rho1} and \eqref{eq: a1} the estimate $\norm{\xi^* (x)}_\infty = \lVert a(x, F_\xi^{-1} (x)-x)\rVert_\infty < \rho$ follows. In conclusion, $\norm{\xi^*}_{\overline{B_2 (0)} , 1} < \rho$ and thus (e) holds.\\
 Recall that $\xi^* (x) = a (x, (F_\xi^{-1}|_{B_{2+\frac{1}{8}} (0)}  - \id_{B_{2+\frac{1}{8}} (0)}) (x))$ for $x \in B_{2+\frac{1}{8}} (0) \subseteq \im F_\xi$ (cf.\ \eqref{eq: exp:inv}). By construction of $\nN$, we obtain $F_\xi^{-1}|_{B_{2+\frac{1}{8}} (0)} - \id_{B_{2+\frac{1}{8}} (0)} \in \lfloor \overline{B_2 (0)} , B_\delta (0)\rfloor_\infty \subseteq C^\infty (B_{2+\frac{1}{8}} (0) , \RR^d)$. Let $a_*$ be the map $a_* \colon \lfloor \overline{B_2 (0)} , B_\delta (0) \rfloor_\infty \rightarrow  C^\infty (B_{2} (0), \RR^d)$ defined via $a_* (\eta ) (x) \coloneq a (x, \eta (x))$. This map is smooth by \cite[Proposition 4.23 (a)]{hg2004} and since $C^\infty (B_{2+\frac{1}{8}} (0), \RR^d)$ is a topological vector space, $\alpha \colon C^\infty (B_{2+\frac{1}{8}} (0), \RR^d) \rightarrow C^\infty (B_{2+\frac{1}{8}} (0), \RR^d), f \mapsto f - \id_{B_{2+\frac{1}{8}} (0)}$ is smooth. We claim that 
	\begin{displaymath}
	 h \colon \nN \rightarrow C^\infty (B_{2+\frac{1}{8}} (0), \RR^d) ,  \xi \mapsto F_\xi^{-1}|_{B_{2+\frac{1}{8}} (0)}
	\end{displaymath}
 is smooth. If this is true, the assertion of (f) follows, since $I = a_* \circ \alpha \circ h$.
 Remark \ref{rem: crfs:prop} (a) implies that the space $C^\infty (B_5 (0), \RR^d)$ is metrizable. Hence by \cite[Proposition E.3]{hg2004}, $h$ is a smooth map if and only if the map $h \circ c$ is smooth for each smooth curve $c \colon \RR \rightarrow \nN$. By the Exponential Law (see, e.g.\ \cite[Proposition 12.2]{hg2004}), the map $h\circ c \colon \RR \rightarrow C^\infty (B_{2+\frac{1}{8}} (0) ,\RR^d)$ will be smooth if $(h\circ c)^\wedge \colon \RR \times B_{2+\frac{1}{8}} (0) \rightarrow \RR^d, (\tau,x) \mapsto h(c (\tau)) (x)$ is smooth. To verify this, we adapt an argument from \cite[p.\ 455]{conv1997}: 
 Consider the map \begin{displaymath}
                  H \colon \RR \times B_{2+\frac{1}{8}} (0) \times B_3 (0) \rightarrow \RR^d , (\tau,x,y) \mapsto \exp (y, c^\wedge (\tau,y)) -x = F_{c (\tau)} (y) -x
                 \end{displaymath}
 which makes sense by construction of $\nN$. Furthermore, $H$ is smooth, as $c^\wedge \colon \RR \times B_5 (0) \rightarrow \RR^d$ is smooth by \cite[Theorem 3.28]{alas2012}. Since $F_{c (\tau)} \circ h(c (\tau)) (x) = x$ holds for each $\tau \in \RR $ and $x \in B_{2+\frac{1}{8}}(0)$, we obtain the identity $H(\tau,x,(h\circ c)^\wedge (\tau,x)) = 0$. A computation yields the following estimate for the derivative of $H$:
	\begin{align*}
	 \opnorm{d_3 H(\tau,x,y;\cdot) - \id_{\RR^d}} &= \opnorm{d_1 \exp (y, c^\wedge (\tau,y); \cdot) + d_2 \exp (y, c^\wedge (\tau,y); d_2 c^\wedge (\tau,y;\cdot)) - \id_{\RR^d}}\\
						   &\leq \opnorm{d_1 \exp (y, c^\wedge (\tau,y);\cdot) -\id_{\RR^d}} + \opnorm{d_2 \exp (y, c^\wedge (\tau,y); d_2 c^\wedge (\tau, y ; \cdot))} \\
						   &\leq \frac{\gamma}{2} + \frac{\gamma}{2} \leq \frac{1}{8} + \frac{1}{8} < 1.
	\end{align*}
 Here we used the estimates for $T_I$ and $T_{II}$ obtained above, which apply because $c (\tau) \in \nN$ holds for each $\tau \in \RR$. We deduce that $d_3 H(\tau,x,y;\cdot)$ is invertible for each $(\tau,x,y) \in \RR \times B_{2 \frac{1}{8}}(0) \times B_3 (0)$. Furthermore, for fixed $(\tau,x) \in \RR \times B_{2 +\frac{1}{8}}(0)$, the map $H(\tau,x, \cdot) = F_{c(\tau)} (\cdot ) - x$ is injective on $B_3 (0)$. Using the injectivity, we deduce with the Implicit Function Theorem \cite[Theorem 5.2]{hg2007} that $(h\circ c)^\wedge$ is smooth. In conclusion, (f) holds. 
\end{proof}

\begin{Alem}[\hspace{-0.5pt}{\cite[II.3 Theorem 3.3]{langdgeo2001}}]\label{lem: mfd:atl} \no{lem: mfd:atl}
 Let $M$ be a finite dimensional paracompact manifold of dimension $d$. Given an open cover $\oO$ of $M$, there exists a locally finite manifold atlas $\vV (\oO) \coloneq \set{(V_{5,k} ,\kappa_k)}_{k \in I}$ with the following properties:
	\begin{compactenum}
	 \item the cover $\vV (\oO)$ is subordinate to $\oO$ and each chartdomain $V_{5,k}$ is precompact,	
	 \item for each $k \in I$, one has $\kappa_k (V_{5,k}) = B_5 (0) \subseteq \RR^d$,
	 \item for each $\tau \in [1,5]$, the open sets $V_{\tau , k} \coloneq \kappa_k^{-1} (B_\tau (0))$ cover $M$ for $k \in I$.
	\end{compactenum}
If $M$ is $\sigma$-compact, then every atlas with properties (a) - (c) is countable.
\end{Alem}

\begin{proof}
 The manifold $M$ is locally compact and paracompact. Apply \cite[Lemma 5.1.6]{Engelking1989} together with local compactness of $M$ to obtain a refinement $\oO'$ of the covering $\oO$, such that the closure of each of the open sets in $\oO'$ is compact and contained in some open set in $\oO$. By Proposition \ref{prop: para:coco} each component of $M$ is second countable and thus we may apply \cite[II.3 Theorem 3.3]{langdgeo2001} to obtain a (countable) locally finite manifold atlas subordinate to $\oO'$ for each component. Thus the closure of any chart domain in this atlas is compact as a closed subset of a compact set. Taking the union of the atlases for the components, we obtain an atlas $\vV (\oO)$ for $M$ with the desired properties. If $M$ is $\sigma$-compact, say $M = \bigcup_{n \in \NN} K_n$ with compact sets $K_n$, then each $K_n$ meets $V_{5,k}$ for only finitely many $k$. Hence $I = \bigcup_{k \in \NN} \setm{k\in I}{V_{5,k} \cap K_n \neq \emptyset}$ is countable.
\end{proof}

 We shall combine our considerations to construct special neighborhoods of the zero-section in $\vect{M}$ for a paracompact Riemannian manifold $(M,\rho_M)$. Consider some atlas $\setm{(V_{5,k},\kappa_k)}{k\in I}$ for $M$ as in Lemma \ref{lem: mfd:atl}. For each chart $(V_{5,k}, \kappa_k)$, we define the pullback Riemannian metric $\rho_k$ on $B_5(0)$ with respect to $\kappa_k^{-1}$. Then $\kappa_k^{-1}$ becomes a Riemannian embedding. In particular, \begin{equation}\label{eq: tausch}
	 T\kappa_k^{-1} (B_{\rho_k} (0_{\kappa_k (x)}, r) )= B_\rho (0_x,r) \quad r >0
	\end{equation} holds for $x \in V_{5,k}$. Moreover, the Riemannian exponential map $\exp_k$ associated to the Riemannian pullback metric $\rho_k$ satisfies $T\kappa_k^{-1} (\dom \exp_k) \subseteq \dom \exp_M$ and
	\begin{equation}\label{eq: exp:komm}
	 \exp_M T\kappa_k^{-1}|_{\dom \exp_k}  = \kappa_k^{-1} \exp_k .
	\end{equation} 
  For the remainder of this section, we endow the image of a manifold chart with the pullback Riemannian metric just described. Whenever the constructions require a Riemannian metric on a chartdomain, we use the induced metric without further mention. In the next lemma, we use notation as in Definition \ref{defn: top:vect}. 

\begin{Alem}\label{lem: eos:cov} \no{lem: eos:cov}
 Let $(M, \rho_M)$ be a $d$-dimensional paracompact Riemannian manifold with Riemannian exponential map $\exp_M$ and some open cover $\oO$ of $M$. Choose via Lemma \ref{lem: mfd:atl} an atlas $\vV (\oO) \coloneq \setm{(V_{5,k},\kappa_k)}{k\in I}$ with respect to $\oO$. There are $\nu_k >0$ for $k\in I$ such that 
		\begin{compactenum}
		 \item for each $y \in M$, the map $\exp_M$ is injective on $N_y \coloneq \bigcup_{n \in I_y} T\kappa_n^{-1} (\set{\kappa_n (y)} \times B_{\nu_n}  (0)) \subseteq T_y M$\glsadd{Ny}, where the index set is defined as $I_y \coloneq \setm{k\in I}{y \in \overline{V_{4,k}}}$. 
		 \item $T \kappa_n (N_y) \subseteq \dom \exp_n$, $\exp_n|_{T \kappa_n (N_y) }$ is an \'{e}tale embedding and $\exp_n T\kappa_n|_{N_y} = \kappa_n \exp_M|_{N_y}$ for each $n \in I_y$. 
		\end{compactenum}
 If $J \subseteq I$ is finite, we may choose $\nu >0$ such that (a), (b) hold for each $k \in I$ with respect to $\nu_k =\nu$. Moreover, in this case there exists open $C^1$-zero-neighborhoods $\nN_k \subseteq C^\infty (V_{5,k},\RR^d)$ for $k \in J$ such that for each $X \in \theta_{\kappa_k}^{-1} (\nN_k) \subseteq \vect{V_{5,k}}$
	\begin{compactitem}
	 	 \item[\rm (c)] the map $\exp_M \circ X|_{\overline{V_{3,k}}}$ is defined, with $\im \exp_M \circ X|_{\overline{V_{3,k}}} \subseteq V_{5,k}$
		 \item[\rm (d)] the following estimates are available: $\exp_M \circ X (\overline{V_{\frac{5}{4},k}}) \subseteq V_{2,k}$, $V_{\frac{5}{4},k} \subseteq \exp_M \circ X(V_{2,k}) \subseteq V_{3,k}$ and $B_4 (0) \times B_{\nu} (0) \subseteq \dom \exp_k$.
		 \item[\rm (e)] the map $F_X^k \coloneq \exp_M \circ X|_{V_{3,k}}$ is an \'{e}tale embedding,
		 \item[\rm (f)] for each $x \in \overline{V_{3,k}}$, we have $X_{\kappa_k}(x) \in B_{\nu} (0)$.
	\end{compactitem}

\end{Alem}

\begin{proof}
 For each $k \in I$, Lemma \ref{lem: exp:loc} allows us to choose $\nu_k' >0$ such that $\exp_k (x,\cdot)$ restricts to an \'{e}tale embedding of $B_{\nu_k'} (x)$ for each $x \in \overline{B_4 (0)}$. Since $\overline{V_{4,k}}$ is compact and the cover $\vV$ is locally finite, there is a finite subset $F_k \subseteq I$ such that $V_{5,i} \cap \overline{V_{4,k}} \neq \emptyset$ if and only if $i \in F_k$. By compactness of $\overline{V_{4,k}} \cap \overline{V_{4,j}}$ for $j \in F_k$, there is some $\nu_k >0$ such that for each $j \in F_k$, one has $T(\kappa_k \circ\kappa_k^{-1}) (\set{\kappa_k (x)} \times B_{\nu_j} (0)) \subseteq \set{\kappa_k (x)} \times B_{\nu_k'} (0)$ for all $x \in \overline{V_{4,k}} \cap \overline{V_{4,j}}$.
 The choice of $\nu_k'$ together with \eqref{eq: exp:komm} shows that the open sets $N_x$ induced by the family $(\nu_k)_{k \in I}$ satisfy the assertion of (a). Since $T\kappa_n (N_x) \subseteq \set{\kappa_n (x)} \times B_{\nu_n'} (0)$ holds for each $n \in I_x$ by construction, the set $T\kappa_n (N_x)$ is contained in the domain of $\exp_n$ for each $n \in I_x$. Hence \eqref{eq: exp:komm} yields $\exp_M|_{N_x} = \exp_M T\kappa_k^{-1}|_{\dom \exp_k} T\kappa_k|_{N_x}  = \kappa_k^{-1} \exp_k T\kappa_k|_{N_x} $. We deduce that (b) must hold.\\
 If $J \subseteq I$ is finite, choose $\nu \coloneq \min \setm{\nu_k}{k \in J}$. We are left to construct the open sets $\nN_k$. Fix $k \in J$ and consider the chart $(V_{5,k}, \kappa_k)$. Reviewing Lemma \ref{lem: mb:exp}, the construction of $\nN_k' \subseteq C^\infty (B_5(0) , \RR^d)$ may be carried out using arbitrarily small $\ve$, since by hypothesis $\ve$ must have the same properties as in Lemma \ref{lem: mb:exp}, where it may be chosen arbitrarily small. 
 The map $\kappa_k$ is a diffeomorphism, whence the pullback $C^\infty (\kappa_k, \RR^d) \colon C^\infty (B_{5} (0) , \RR^d) \rightarrow C^\infty (V_{5,k}, \RR^d) , f \mapsto f \circ \kappa_k$ is linear bijective and continuous by a combination of \cite[Lemma 4.11]{hg2004} and \cite[Proposition 7.4]{bgn2004}. Define the open $C^1$-neighborhood $\nN_k \coloneq C^\infty (\kappa_k,\RR^d)^{-1} (\nN_k') \subseteq C^\infty (V_{5,k},\RR^d)$. The Riemannian exponential map $\exp_k$ is related to $\exp_M$ via \eqref{eq: exp:komm} and the identity in (b). Hence the properties obtained via Lemma \ref{lem: mb:exp} for vector fields with $X_{\kappa_k} \in \nN_k$ imply (c) - (f).
\end{proof} 

 \begin{Asetup}
 In the setting of Lemma \ref{lem: eos:cov}, consider a compact subset $K \subseteq M$. As $\vV (\oO)$ is locally finite, there is a finite subset $\fF_{5} (K) \coloneq \setm{(V_{5,k_j} , \kappa_{k_j})}{1 \leq j \leq N}$\glsadd{F5} of $\vV(\oO)$ such that $V_{5,k} \cap K \neq \emptyset$ holds if and only if $(V_{5,k},\kappa_k) \in \fF_{5} (K)$. Notice that $\fF_5 (K)$ induces a family of open neighborhoods \glsadd{omega1} of $K$ via 
	\begin{displaymath}
	 K \subseteq \Omega_{r,K} \coloneq \bigcup_{l=1}^N V_{r,k_l} , \quad r \in [1,5]
	\end{displaymath}
 The set $\fF_5 (K)$ is finite, whence the set $K_5 \coloneq \bigcup_{l =1}^N \overline{V_{5,k_l}}$\glsadd{K5} is compact. Again, we define a finite subset $\fF_5 (K_5) \coloneq \setm{(V_{5,n}, \kappa_n)}{n\in I, V_{5,n} \cap K_5 \neq \emptyset}$ of $\vV (\oO)$ as the set of charts which intersect the compact set $K_5$. As above, one defines open neighborhoods $\Omega_{r,K_5}$\glsadd{omegaK5} of $K_5$ for $r \in [1,5]$.
\end{Asetup}

We will now construct a neighborhood of the zero section such that the composition of sections in this neighborhood with the Riemannian exponential map yields an \'{e}tale embedding. The arguments in the proof of the following lemma are inspired by \cite[2.\ Theorem 1.4]{hirsch1976}.

\begin{Alem}\label{lem: vect:emb} \no{lem: vect:emb}
 Let $K \subseteq M$ be a compact set and $\fF_5 (K) = \setm{(V_{5,k}, \kappa_{k})}{1\leq k \leq N}$ as above. Construct for each $1\leq k \leq N$ a $C^1$-zero-neighborhood $\nN_k \subseteq C^\infty (V_{5,k},\RR^d)$ as in Lemma \ref{lem: eos:cov} (c)-(f) applied with the finite set $J=\set{1,\ldots, N}$. Furthermore, consider the continuous maps $\theta^{\Omega_{5,K}}_{\kappa_k} \colon \vect{\Omega_{5,K}} \rightarrow C^\infty (V_{5,k}, \RR^d), X \mapsto X_{\kappa_k}$.
 There are open $C^1$-zero-neighborhoods $\mM_{k} \subseteq \nN_{k}$ such that, setting $E_{5,K} \coloneq \bigcap_{k=1}^N (\theta^{\Omega_{5,K} }_{k})^{-1} (\mM_k) \subseteq \vect{\Omega_{5,K}}$ and $E \coloneq (\res_{\Omega_{5,K}}^M)^{-1} (E_{5,K}) \subseteq \vect{M}$ (cf.\ Notation \ref{nota: res}), the map $F_X \coloneq \exp_M \circ X|_{\Omega_{2,K}}$ is an \'{e}tale embedding for each $X \in E$, and $F_X (\overline{\Omega_{1,K}}) \subseteq \Omega_{2,K}$ holds.
\end{Alem}

\begin{proof}
 By Lemma \ref{lem: eos:cov} for each $X \in \theta^{-1}_{\kappa_k} (\nN_k)$ the map $\exp_M \circ X|_{ \overline{V_{2, k}}}$ is defined and its image is contained in $V_{3,k}$ for each $(V_{5,k},\kappa_k) \in \fF_5 (K)$. The manifold $M$ is locally compact, hence a regular topological space. Thus by \cite[Theorem 3.1.6]{Engelking1989}, we may separate the compact set $\overline{V_{2, k}}$ from the closed set $M \setminus V_{3,k}$. We obtain disjoint open sets $A_k, B_k \subseteq M$ such that $\overline{V_{2, k}} \subseteq A_k$ and $M\setminus V_{3,k} \subseteq B_k$ hold for each $(V_{5,k},\kappa_k) \in \fF_5 (K)$. \\ \textbf{Claim:} There are open neighborhoods $\mM_k \subseteq \nN_k$ of the zero-map, $1\leq k \leq N$, such that for $X \in E_{5,K}$ the following holds: 
 $F_X (\overline{V_{2,k}} \cap \Omega_{2,K})  \subseteq A_k$ and $F_X (\Omega_{2,K} \setminus V_{3,k}) \subseteq B_k$ for each $1\leq k \leq N$. 
 If this is true, then the proof may be completed as follows:\\
 Let $X$ be contained in $E_{5,k}$. Observe that the construction of $E_{5,k}$ implies that for each $1\leq k \leq N$ the map $F_X|_{V_{3,k} \cap \Omega_{2,K}} =F_X^k|_{V_{3,k}\cap \Omega_{2,k}}$ is an \'{e}tale embedding by Lemma \ref{lem: eos:cov} (e). Consider distinct $x,y \in \Omega_{2,K}$ and choose $1\leq k \leq N$ with $x \in \overline{V_{2,k}}$. If $y \in V_{3,k}$ we must have $F_X (x) \neq F_X (y)$ since the map is an \'{e}tale embedding on $V_{3,k} \cap \Omega_{2,K}$. On the other hand, if $y \in \Omega_{2,K} \setminus V_{3,k} \subseteq M \setminus V_{3,k}$, by the above $F_X (x) \in  F_X (\overline{V_{2,k}}\cap \Omega_{2,K}) \subseteq A_k$ and $F_X (y) \in F_X(\Omega_{2,K} \setminus V_{3,k}) \subseteq B_k$. Since $A_k$ and $B_k$ are disjoint, again $F_X (x) \neq F_X (y)$ follows, whence $F_X$ must be injective. Thus each $X \in E$ yields an injective local diffeomorphism $\exp_M \circ X|_{\Omega_{2,K}}$, i.e.\ $\exp_M \circ X|_{\Omega_{2,K}}$ is an \'{e}tale embedding. Furthermore, $F_X$ maps $\overline{V_{1,k}}$ into $V_{2,k}$ by Lemma \ref{lem: eos:cov} (d). Hence the definition of $\Omega_{1,K}$ and $\Omega_{2,K}$ yield $F_X  (\overline{\Omega_{1,K}}) \subseteq \Omega_{2,K}$.
 \\[1em] \textbf{Proof of the claim:} For $k \neq j$, we obtain a sets  
  \begin{displaymath}
   K_{kj} \coloneq \kappa_k (\overline{V_{2,k}} \cap (M \setminus V_{3,j})) \subseteq \overline{B_2 (0)} \text{ and } B_{kj} \coloneq T\kappa_k (TV_{5,k} \cap \exp_M^{-1} (B_j \cap V_{3,k})).
  \end{displaymath}
  By construction each set $K_{kj} \subseteq B_5 (0)$ is compact and each set $B_{kj}$ is an open subset of $TB_5 (0)$. Define $A_{kk} \coloneq T\kappa_k (TV_{5,k} \cap \exp_M^{-1} (A_k))$ for $1\leq k \leq N$. Recall the identity $\exp_M \circ 0_M = \id_M$, where $0_M \in \vect{M}$ is the zero section. This yields the inclusions $K_{kj} \times \set{0} \subseteq B_{kj}$ for each pair $(k,j) \in \setm{1\leq k,j \leq N}{k \neq j}$ and $\overline{B_2 (0)} \times \set{0} \subseteq A_{kk}$. Hence, the Wallace Lemma \cite[3.2.10]{Engelking1989} yields constants $\ve_{kj} >0$ for $1\leq j \leq N$ which satisfy $K_{kj} \times B_{\ve_{kj}} (0) \subseteq B_{kj}$ and $\overline{B_{2} (0)} \times B_{\ve_{kk}} (0) \subseteq A_{kk}$ for each pair $(k,j) \in \setm{1\leq k,j \leq N}{k \neq j}$. Moreover, for $1\leq k\leq N$ we obtain an open neighborhood 
	\begin{displaymath}
	 M_k \coloneq \lfloor \overline{B_2 (0)} , B_{\ve_{kk}} (0) \rfloor \cap \bigcap_{\substack{j=1 \\ j \neq k}}^N \lfloor K_{kj}, B_{\ve_{kj}} (0)\rfloor \subseteq C^\infty (B_5 (0) , \RR^d) 
	\end{displaymath}
 of the zero-map. Define the $C^1$-open set $\mM_k \coloneq C^\infty (\kappa_k, \RR^d)^{-1} (M_k) \cap \nN_k \subseteq C^\infty (V_{5,k}, \RR^d)$. By construction, each vector field $X \in E_{5,K}$ (defined as in the statement of the lemma) may be composed on $\Omega_{3,K}$ with $\exp_M$. With the identities \eqref{eq: exp:komm} and Lemma \ref{lem: eos:cov} (b), the mapping $F_X$ may be evaluated locally on $\overline{V_{2,k}}$ in the chart $(V_{5,k},\kappa_k) \in \fF_{5} (K)$. For any $X \in E_{5,K}$, we note that $X_{\kappa_k} \in C^\infty (\kappa_k, \RR^d)^{-1} (\lfloor \overline{B_2 (0)} , B_{\ve_{kk}} (0) \rfloor)$ holds. Observe that $\overline{B_2 (0)} \times B_{\ve_{kk}} (0) \subseteq A_{kk}$ and the definition of $A_{kk}$ imply $F_X (\overline{V_{2,k}}) \subseteq A_k$. Furthermore, each element $y \in \Omega_{2,K}\setminus V_{3,k}$ is contained in $\overline{V_{2,n}}$ for some $1\leq n \leq N$. Thus $\kappa_n (y)$ is contained in $K_{nk}$ by construction. Furthermore, $X_{\kappa_n} \in C^\infty (\kappa_n, \RR^d)^{-1} (\lfloor K_{nk} , B_{\ve_{nk}} (0) \rfloor)$ and  $K_{nk} \times B_{\ve_{nk}} (0) \subseteq B_{nk}$ hold. By definition of $B_{nk}$, a computation in the chart $(V_{5,n},\kappa_n)$ yields $F_X (y) \in B_k$. As $y \in  \Omega_{2,K}\setminus V_{3,k}$ and $k$ were chosen arbitrarily, $F_X (\Omega_{2,K}\setminus V_{3,k}) \subseteq B_k$ holds for each $1 \leq k \leq N$. 
\end{proof}

 We are interested in vector fields which yield, after composition with the Riemannian exponential map, an inverse for $F_X$ (respectively, the composition $F_Y \circ F_X$). In the rest of this section, we construct $C^1$-neighborhoods of the zero-section, whose elements permit such vector fields. Furthermore, the mappings sending a vector fields to the vector field which induces $F_X \circ F_Y$ respectively $F_X^{-1}$ should be smooth on these neighborhoods. 
 The leading idea is to construct these fields locally in a cover of charts, which will enable us to obtain them as global objects from the local data. For reasons which are explained in Section \ref{sect: lgp:str}, we construct a neighborhood of the zero-section depending on an open $C^1$-neighborhood of the zero-section chosen in advance and on a positive constant $R$.
	
\begin{Acon}\label{con: comp:loc} \no{con: comp:loc}
 Consider the setting of Lemma \ref{lem: vect:emb}: Let $K \subseteq M$ be compact and $E_{5,K} \subseteq \vect{\Omega_{5,K}}$ an open neighborhood of the zero-section as in Lemma \ref{lem: vect:emb}. Fix $R > 0$ and an arbitrary open $C^1$-neighborhood of the zero-section $P \subseteq \vect{\Omega_{5,K}}$. By construction of the manifold atlas, $\Omega_{5,K} \subseteq \Omega_{1,K_5}$ holds by Lemma \ref{lem: mfd:atl} (c). As the family $\fF_5 (K_5)$ is a manifold atlas for $\Omega_{5,K_5}$, the topology on $\vect{\Omega_{1,K_5}}$ is initial with respect to the family $\setm{\theta_{\kappa_k|_{V_{1,k}}}}{(V_{5,k}, \kappa_{k}) \in \fF_5 (K_5)}$ by Definition \ref{defn: top:vect}. Thus there is a family of open $C^1$-neighborhoods of the zero-map $W_k \subseteq C^\infty (B_1 (0),\RR^d) \cong  C^\infty (V_{1,k}, \RR^d), (V_{5,k}, \kappa_k) \in \fF_5 (K_5)$ with
 \begin{displaymath}
        (\res_{\Omega_{5,K}}^{\Omega_{1,K_5}})^{-1} (E_{5,K} \cap P) \supseteq \bigcap_{\fF_5 (K_5)} (\theta_{\kappa_k|_{V_{1,k}}} \circ C^\infty (\kappa_k|_{V_{1,k}}, \RR^d))^{-1} (W_k).                                                                                                                                                                                                                                                                                                                                                                                                                                                                                                                                                                                                                                                                                                                                        
  \end{displaymath}
 Here $C^\infty (\kappa_k|_{V_{1,k}}, \RR^d) \colon C^\infty (B_1 (0), \RR^d) \rightarrow C^\infty(V_{1,k},\RR^d)$ denotes the pullback $f \mapsto f\circ \kappa|_{V_{1,k}}$, which is continuous by \cite[Lemma 4.4]{hg2004}. Since $B_1 (0) \subseteq B_5 (0) = \kappa_k (V_{5,k})$ holds, Remark \ref{rem: Cr:Norm} (a) implies that we may choose $\tau >0$ such that for $f \in B_\tau^{k} \coloneq \setm{f \in C^\infty (B_5 (0),\RR^d)}{\norm{f}_{\overline{B_1 (0)}, 1} < \tau}$ the condition $f|_{B_1 (0)} \in W_k$ is satisfied. Shrinking $\tau$ if necessary, we may assume $\tau  < R$. Define the open $C^1$-neighborhood of the zero-section 
  \begin{displaymath}
   E' \coloneq \bigcap_{\fF_{5} (K_5)} (\theta^{\Omega_{5,K_5}}_{\kappa_k} \circ C^\infty (\kappa_k, \RR^d))^{-1} (B_{\tau}^{k}) \subseteq \vect{\Omega_{5,K_5}}.
  \end{displaymath}
 Then, the inclusions $E' \subseteq (\res^{\Omega_{5,K_5}}_{\Omega_{5,K}})^{-1} (E_{5,K} \cap P)$ and $(\res_{\Omega_{5,K_5}}^M)^{-1} (E') \subseteq E \cap (\res_{\Omega_{5,K}^M})^{-1} (P)$ hold. 
\paragraph{Step 1: A vector field inducing the composition $\exp_M \circ X \circ F_Y$:}
 Since the family $\fF_5 (K_5)$ is finite, we may fix a constant $\nu >0$ with $\nu < R$ as in Lemma \ref{lem: eos:cov}. Consider arbitrary $(V_{5,n}, \kappa_n) \in \fF_{5} (K_5)$ and shrink the $C^1$-open set $B_\tau^{n}$:
 Choose $\ve_n > \sigma_{\frac{\delta_n}{2}} > 0$ and $1>\delta_n >0$ with properties as in Lemma \ref{lem: exp:loc}, such that $\ve_n < \min \set{\tau , \nu}$ holds. Set $\sigma_n \coloneq \sigma_{\frac{\delta_n}{2}}$ and $\rho_n  \coloneq \min \set{\nu ,\tau}$. Apply Lemma \ref{lem: mb:exp} with the constants $\ve_n,\delta_n, \rho_n$ taking the roles of $\ve, \delta, \rho$ to obtain a $C^1$-neighborhood $\nN_n$ of the zero-map in $C^\infty (B_5 (0), \RR^d)$. Then each $X \in C^\infty (\kappa_n, \RR^d) (\nN_n) \subseteq C^\infty (V_{5,n},\RR^d)$ satisfies the assertions of Lemma \ref{lem: eos:cov} (c)-(e) with respect to $\nu$. By choice of the constants (cf.\ Lemma \ref{lem: exp:loc}), there is a smooth map $f_n \colon B_3 (0) \times B_{\sigma_n} (0) \times B_{\sigma_n} (0) \rightarrow B_{\ve_n} (0)$ with 
  \begin{equation}\label{eq:fn:ident}
   f_n (x,0,0) = 0, f_n (x,y,0)= y \text{ and } f_n(x,0,z) = z, \quad (x,y,z) \in B_3 (0) \times B_{\sigma_n} (0) \times B_{\sigma_n} (0).
  \end{equation}
 Hence the partial derivative satisfies $d_1f_n (x,0,0;\cdot) = 0$, for all $x \in B_3 (0)$. The continuous map $B_3(0) \times B_{\sigma_n} (0) \times B_{\sigma_n} (0), (x,y,z) \mapsto \opnorm{df_n (x,y,z; \cdot)}$ is bounded on $\overline{B_2 (0)} \times \overline{B_{ \frac{\sigma_n}{2}}(0)} \times \overline{B_{ \frac{\sigma_n}{2}} (0)}$ by some $t_n \geq 1$. As the partial derivative with respect to $x$ vanishes in $\overline{B_2 (0)} \times \set{0} \times \set{0}$, a compactness argument yields $0 < \mu_n < \min \set{\nu, \frac{\sigma_n}{2}, \frac{\tau}{6dt_n}}$ such that for all $(x,y,z) \in \overline{B_2 (0)} \times \overline{B_{\mu_n} (0)} \times \overline{B_{\mu_n} (0)}$ the estimate $\opnorm{d_1f_n (x,y,z; \cdot)} < \frac{\tau}{3}$ holds. Define the open $C^1$-zero-neighborhood
  \begin{displaymath}
   \hH_n' \coloneq \nN_n \cap \setm{f \in C^\infty (B_5 (0) , \RR^d)}{\norm{f}_{\overline{B_3 (0)}, 1} < \mu_n} \subseteq C^\infty (B_5(0), \RR^d)
  \end{displaymath}
 Since $\mu_n  < \tau$ holds, we deduce $\hH_n' \subseteq B_{\tau}^n$. Set $\hH' \coloneq \bigcap_{\fF_{5,K_5}} (\theta^{\Omega_{5,K_5}}_{\kappa_n})^{-1} C^\infty (\kappa_n, \RR^d) (\hH_n') \subseteq \vect{\Omega_{5,K_5}}$ to obtain a $C^1$-neighborhood of the zero-section contained in $E'$.
 \\
 Let $\xi, \eta$ be elements of $\hH_n'$. By Lemma \ref{lem: mb:exp}, $F_\xi (B_2 (0)) \subseteq B_3 (0)$ holds, whence the composition $F_\eta \circ F_{\xi}|_{B_2 (0)}$ is defined. Since $\mu_n < \sigma_n$, we have $F_\eta F_\xi (x) \in B_{\delta_n} (x)$ for each $x \in B_2 (0)$ by definition of $\sigma_n = \sigma_{\frac{\delta_n}{2}}$ (cf.\ Lemma \ref{lem: exp:loc}). Therefore, for each $x\in B_2 (0)$,
	\begin{equation}\label{eq: diamond}
	 \eta \diamond \xi (x) \coloneq \phi_x^{-1} (F_\eta F_\xi (x)) = f_n (x, \xi (x) , \eta (F_\xi (x))) \in B_{\ve_n} (0) \subseteq B_\tau (0)
	\end{equation}
 is defined and yields a smooth map $\eta \diamond \xi \colon B_2 (0) \rightarrow B_{\ve_n} (0) \subseteq B_\tau (0)$. Observe that $\eta, \xi \equiv 0$ implies $\eta \diamond \xi \equiv 0$ by \eqref{eq:fn:ident}. For $(V_{5,n},\kappa_n) \in \fF_5 (K_5)$ and $X \in \hH'$, set $X_{[n]} \coloneq X_{\kappa_n} \circ \kappa_n^{-1}$. Moreover, for $X \in \hH'$ the composition $F_X \coloneq \exp_M \circ X|_{\Omega_{3,K}}$ is defined. Consider $y\in V_{3,n}$ and $X \in \hH'$. By construction $X_{[n]} \in \hH_n'$, whence $X_{[n]} (\kappa_n (y)) \in B_{\mu_n} (0) \subseteq B_\nu (0)$. Since $\set{\kappa_n (y)} \times B_\nu (0) \subseteq T\kappa_n (N_y)$, Lemma \ref{lem: eos:cov} (b) yields for $F_{X_{[n]}}$ as in Lemma \ref{lem: mb:exp} 
    \begin{align*}
    \kappa_n^{-1} F_{X_{[n]}} (\kappa_n (y)) &= \kappa_n^{-1} \exp_n (\kappa_n (y) , X_{[n]} (\kappa_n (y))) = \kappa_n^{-1} \exp_n T\kappa_n \circ X (y) \\
                &= \exp_M T\kappa_n^{-1} T\kappa_n \circ X (y) =  \exp_M \circ X (y) = F_X (y).  
    \end{align*}
Furthermore, a combination of Lemma \ref{lem: eos:cov} (b) and (c) allows us to compute the identity 
\begin{displaymath}
                         T\kappa_n (\exp_M|_{N_y})^{-1} \kappa_n^{-1}|_{\exp_n (T\kappa_n (N_y))} = (\exp_n|_{T\kappa_n (N_y)})^{-1}
\end{displaymath}
for $y \in V_{3,n}$. Set $x \coloneq \kappa_n (y)$ with $y \in V_{2,n}$. Since $\ve_n < \nu$, we conclude $\set{x} \times B_{\ve_n} (0) \subseteq  T\kappa_n (N_y) $. This yields the following identity:
 	\begin{align}
      &(\id_{B_{2} (0)}, X_{[n]} \diamond Y_{[n]}) (x) =  (x, f_n (x, Y_{[n]} (x) , X_{[n]} (F_{Y_{[n]}}(x)))) = (\exp_n|_{\set{x} \times B_{\ve_n} (0)})^{-1} F_{X_{[n]}} F_{Y_{[n]}} (x) \notag \\ 
    =&(\exp_n|_{T\kappa_n (N_y)})^{-1} F_{X_{[n]}} F_{Y_{[n]}} (x) = T\kappa_n (\exp_M|_{N_y})^{-1} \kappa_n^{-1}|_{\exp_n (T\kappa_n (N_y))} F_{X_{[n]}} F_{Y_{[n]}} (x) \notag \\ 
    =& T\kappa_n (\exp_M|_{N_y})^{-1} \kappa_n^{-1} F_{X_{[n]}} F_{Y_{[n]}} (x) = T\kappa_n (\exp_M|_{N_y})^{-1} \exp_M X \kappa_n^{-1} F_{Y_{[n]}} (x)  \notag \\
    =& T\kappa_n (\exp_M|_{N_y})^{-1} \exp_M X \exp_M Y (y) = T\kappa_n (\exp_M|_{N_y})^{-1} \exp_M X (F_Y (y)) .\label{eq: glo:loc}
	\end{align}
 This assignment is defined and smooth on $V_{2,n}$, by \eqref{eq: diamond}. Hence for $X,Y \in \hH'$, we can define $X \diamond Y \colon \Omega_{2,K_5} \rightarrow TM , x \mapsto (\exp_M|_{N_x})^{-1} (\exp_M \circ X \circ \exp_M \circ Y) (x)$, which is an element of $\vect{\Omega_{2,K_5}}$.\glsadd{XdiamondY} The identity \eqref{eq: glo:loc} yields $X \diamond Y \equiv 0$ for $X,Y\equiv 0$. Define for $X, Y \in \hH'$ the map $(X \diamond Y)_{[n]} \coloneq (X\diamond Y)_{\kappa_n|V_{2,n}} \circ \kappa_n^{-1}|_{B_2 (0)}$. Then the above computation \eqref{eq: glo:loc} yields $(X \diamond Y)_{[n]} = X_{[n]} \diamond Y_{[n]}$ on $B_2 (0)$. From \eqref{eq: diamond}, we deduce
	\begin{equation}\label{eq: locdiamest}
		\norm{(X \diamond Y)_{[n]}}_{\overline{B_{\frac{3}{2}} (0)},0} = \norm{X_{[n]} \diamond Y_{[n]}}_{\overline{B_{\frac{3}{2}} (0)},0} < \varepsilon_n < \min \set{\tau , \nu} < R.
	\end{equation} 
  \paragraph{Step 2: A vector field inducing $F_X^{-1}$:} By construction, each $\hH_n'$ for $(V_{5,n},\kappa_n) \in \fF_5 (K_5)$ is contained in a set $\nN_n$ as constructed via Lemma \ref{lem: mb:exp} such that the assertions of Lemma \ref{lem: eos:cov} (c) -(e) hold for $C^\infty (\kappa_n, \RR^d)(\nN_n)$. In particular, we may apply Lemma \ref{lem: vect:emb} with $K = K_5$, the open cover $\fF_5 (K_5)$ and the open sets $(\hH_n')_{(V_{5,n},\kappa_n) \in \fF_5 (K_5)}$: For each chart in $\fF_5 (K_5)$, we obtain an open $C^1$-zero-neighborhood $\hH_n \subseteq C^\infty (\kappa_n, \RR^d) (\hH_n') \subseteq C^\infty (V_{5,n}, \RR^d)$. Then define 
    \begin{displaymath}
     \hH^{\Omega_{5,K_5}}_R \coloneq \bigcap_{(V_{5,n}, \kappa_n) \in \fF_5 (K_5)} (\theta^{\Omega_{5,K_5}}_{\kappa_n})^{-1} (\hH_n) \subseteq \hH'.
    \end{displaymath}
 By Lemma \ref{lem: eos:cov} (e) for each $X \in \hH^{\Omega_{5,K_5}}_R$ the mapping $\exp_M \circ X|_{\Omega_{2, K_5}}$ is a \'{e}tale embedding.
  Consider $X \in \hH^{\Omega_{5,K_5}}_R$ and $(V_{5,n}, \kappa_n) \in \fF_5 (K_5)$. By construction of $\hH_n'$ in Step 1, we deduce with Lemma \ref{lem: mb:exp} (c) that $B_{\frac{5}{4}} (0) \subseteq F_{X_{[n]}} (B_2 (0))$ holds. We already established the identities $F_X (y) =\kappa_n^{-1} F_{X_{[n]}} (\kappa_n (y))$ and $T\kappa_n (\exp_M|_{N_y})^{-1} \kappa_n^{-1}|_{\exp_n (T\kappa_n (N_y))} = (\exp_n|_{T\kappa_n (N_y)})^{-1}$ for $y \in V_{3,n}$ and $X \in \hH^{\Omega_{5,K_5}}_R$. Furthermore, Lemma \ref{lem: mb:exp} (c)-(e) yield a map $X_{[n]}^* \in C^\infty (\Im F_{X_{[n]}} , \RR^d)$ with $F_{X_{[n]}^*} \coloneq \exp_n (\id_{\im F_{X_{[n]}}}, X_{[n]}^*) = F_{X_{[n]}}^{-1}$. This map satisfies $\norm{X_{[n]}^*}_{\overline{B_2 (0)},1} < \rho_n = \min \set{\nu, \tau} < R$. Hence by choice of $\nu$, we deduce $X_{[n]}^* (y) \in T\kappa_n (N_y)$ and thus $F_{X_{[n]}^*} (y) \in \exp_n (T\kappa_n (N_y))$ for each $y \in V_{\frac{5}{4},n}$. Combining these facts we compute for $(V_{5,n}, \kappa_n) \in \fF_5 (K_5)$ and $y \in V_{\frac{5}{4}}$: 
  \begin{align*}
    T\kappa_n^{-1} (\exp_n|_{T\kappa_n (N_y)})^{-1} F_{X_{[n]}^*} (\kappa_n (y)) &= (\exp_M|_{N_y})^{-1} \kappa_n^{-1} (F_{X_{[n]}})^{-1} (\kappa_n (y)) \\
										 &= (\exp_M|_{N_y})^{-1} (\kappa_n^{-1} F_{X_{[n]}} \kappa_n)^{-1} (y)\\
										 &= (\exp_M|_{N_y})^{-1} F_X^{-1} (y) = (\exp_M|_{N_y})^{-1} (\exp_M X|_{\Omega_{2,K_5}})^{-1}(y).
  \end{align*}
 By the computation, we obtain a section of the tangent bundle on $\Omega_{\frac{5}{4}, K_5}$ via 
	\begin{displaymath}
	 X^* \colon \Omega_{\frac{5}{4}, K_5} \rightarrow TM, X^* (y) \coloneq (\exp_M|_{N_y})^{-1} \circ (\exp_M \circ X)^{-1} (y).
	\end{displaymath}
 Let $(V_{5,n}, \kappa_n) \in \fF_5 (K_5)$ and $y \in V_{\frac{5}{4},n}$. Observe that $\exp_n|_{T\kappa_n (N_y)}$ is injective. Furthermore, $F_{X_{[n]}^*} (\kappa_n (y)) = \exp_n (\kappa_n (y), X_{[n]}^* (\kappa_n (y)))$ and $(\kappa_n (y) , X_{[n]}^* (\kappa_n (y))) \in T\kappa_n (N_y)$. These identies imply $(\exp_n|_{T\kappa_n (N_y)})^{-1} F_{X_{[n]}^*} (\kappa_n (y)) = (\kappa_n (y) , X_{[n]}^* (\kappa_n (y))$, whence the local identity above yields\glsadd{Xstar}
 	\begin{equation} \label{eq: inv:loc}
 	 X^* (y) \coloneq (\exp_M|_{N_y})^{-1} \circ F_{X}^{-1} (y)= T\kappa_n^{-1} (\id_{B_2 (0)}, X_{[n]}^*) \kappa_n (y) \text{ for each } y \in V_{\frac{5}{4},n}.
 	\end{equation}
 As $X_{[n]}^*$ is a smooth map by Lemma \ref{lem: mb:exp}, \eqref{eq: inv:loc} shows that $X^*$ is smooth. Hence $X^* \in \vect{\Omega_{\frac{5}{4}, K_5}}$ follows. In addition for each $(V_{5,n},\kappa_n) \in \fF_5 (K_5)$, by choice of $\rho_n$ 
	\begin{equation}\label{eq: locinvest}
	 \norm{X_{[n]}^*}_{\overline{B_2 (0)},1} < \rho_n = \min \set{\nu, \tau } < R. 
	\end{equation}
 Define $\hH_R \coloneq (\res^M_{\Omega_{5,K_5}})^{-1} (\hH^{\Omega_{5,K_5}}_R)$ and observe that the estimates obtained in Step 1 and 2 remain valid for sections in this set.
 \paragraph{Conclusion:} We have constructed $C^1$-neighborhoods of the zero-section 
 \begin{align*}
  \hH^{\Omega_{5,K_5}}_{R} &\coloneq \Gamma^{-1} \left(\prod_{(V_{5,n},\kappa_n)  \in \fF_5 (K_5)} \hH_n\right) \subseteq \vect{\Omega_{5,K_5}}, \\
  \hH_R &\coloneq (\res^{M}_{\Omega_{5,K_5}})^{-1} (\hH^{\Omega_{5,K_5}}_{R}) \subseteq \vect{M}
 \end{align*}
 where $\Gamma \colon \vect{\Omega_{5,K_5}} \rightarrow \prod_{(V_{5,n},\kappa_n)  \in \fF_5 (K_5)}  C^\infty (V_{5,n},\RR^d)$ is the embedding defined in \ref{defn: top:vect} and each $\hH_n \subseteq C^\infty (V_{5,n},\RR^d)$ is an open $C^1$-neighborhood of the zero map.\\
 By construction, $\hH_R$ is contained in the zero-neighborhood $E \cap (\res^M_{\Omega_{5,K}})^{-1} (P)$ chosen in advance. Here $E$ is a neighborhood as in Lemma \ref{lem: vect:emb} and $P \subseteq \vect{\Omega_{5,K}}$ is an open $C^1$-neighborhood of the zero-section. In particular, Lemma \ref{lem: vect:emb} implies that each element of $\hH_R$ satisfies the assertions of Lemma \ref{lem: eos:cov} (d), i.e.:\\ 
 For $(V_{5,n},\kappa_n) \in \fF (K_5)$ and $X \in \hH_R$, we have $X_{\kappa_n}(\overline{V_{1,n}}) \subseteq B_{\nu} (0)$ with $B_2 (0) \times B_{\nu} (0) \subseteq \dom \exp_n$. For a pair $(X,Y) \in \hH_R \times \hH_R$ there are vector fields $X \diamond Y \in \vect{\Omega_{2,K_5}}$ and $X^* \in \vect{\Omega_{\frac{5}{4},K_5}}$, respectively, such that the following identities are satisfied: 
	\begin{align}\label{eq: def1}
	 \exp_M \circ (X \diamond Y)&= \exp_M X \exp_M Y|_{\Omega_2,{K_5}}\\
	 \exp_M \circ X^* &= (\exp_M \circ X|_{\Omega_{2,K_5}})^{-1}|_{\Omega_{\frac{5}{4}, K_5}} .\label{eq: def2}
	\end{align}
 We note that if $X$ and $Y$ are the zero section, then the local formulas \eqref{eq: glo:loc} and \eqref{eq: inv:loc} (with Lemma \ref{lem: mb:exp} (e)) prove that $X \diamond Y$ and $X^*$ are the zero section in $\vect{\Omega_{2, K_5}}$ and $\vect{\Omega_{\frac{5}{4},K_5}}$, respectively. 
\end{Acon}
 The neighborhood $\hH_R$ constructed in this section is used in Section \ref{sect: lgp:str} to obtain symmetric neighborhoods in the space of compactly supported orbisections. The argument in Construction \ref{con: comp:loc} depends only on a finite atlas. Hence the sets constructed are open in $\vect{M}$ with the topology introduced in Definition \ref{defn: top:vect}. Unfortunately, the vector fields $X \diamond Y$ and $X^*$ will thus in general {\bf not} be defined on all of $M$. Because of this, we are not able to prove a statement of the following kind: If $X,Y \in \hH_R$, then $X \diamond Y \in E$ and $X^* \in E$. At the moment, we can only prove the following:

\begin{Acor} \label{cor: comp:tau} \no{cor: comp:tau}
 Consider the setting of Construction \ref{con: comp:loc} and let $\hH_n'$, $(V_{5,n},\kappa_n) \in \fF_5(K_5)$ and $\hH_R$ be as constructed there. For each pair $\eta, \xi \in \hH_n'$, the map $\eta \diamond \xi \colon B_2 (0) \rightarrow B_\tau (0)$ satisfies $\norm{\eta \diamond \xi}_{\overline{B_1 (0)},1} < \tau < R$.
 Hence, by \eqref{eq: glo:loc}, for any pair $(X,Y) \in \hH_R \times \hH_R$ and each chart $(V_{5,n},\kappa_n) \in \fF_{5} (K_5)$, we derive $\norm{(X \diamond Y)_{[n]}}_{\overline{B_1 (0)},1} < \tau$. 
\end{Acor}

In Section \ref{sect: lgp:str}, we consider a setting, which allows $X\diamond Y$ to be extended uniquely to all of $M$. In this case, Corollary \ref{cor: comp:tau} will imply the result mentioned above (cf.\ Proposition \ref{prop: os:glue}). 

 \begin{proof}[Proof of Corollary \ref{cor: comp:tau}]
 By \eqref{eq: diamond}, it suffices to prove that the norm of the derivative is bounded by $\tau$. To do so, we recall the estimates from Step 1 of Construction \ref{con: comp:loc}: Let $x \in \overline{B_1 (0)}, y \in B_2 (0)$ and consider $\xi \in \hH_n'$. Then $F_\xi (x) \in B_2 (0)$ and $\norm{\xi}_{\overline{B_3 (0)},1} < \mu_n $ with $0 <\mu_n < \min \set{\nu, \frac{\sigma_n}{2}, \frac{\tau}{6dt_n}}$. Recall that $\opnorm{d_1 f_n (y_1,y_2,y_3;\cdot)} < \tfrac{\tau}{3}$ holds and $t_n$ is an upper bound for $\opnorm{df_n (y_1,y_2,y_3;\cdot)}$ with $(y_1,y_2,y_3) \in \overline{B_2 (0)} \times \overline{B_{\mu_n}(0)} \times \overline{B_{\mu_n} (0)}$. As $\hH_n' \subseteq \nN_n$ for an open neighborhood $\nN_n$ constructed via Lemma \ref{lem: mb:exp}, we deduce from the proof of the lemma that $\frac{1}{4} \geq \opnorm{dF_\xi (x;\cdot) -\id_{\RR^d}} \geq \opnorm{dF_\xi (x;\cdot)} -1$ for $\norm{x}_\infty < 3$. For each $(x,y) \in \overline{B_1 (0)} \times \overline{B_1 (0)}$ we obtain the estimate $\norm{d\xi (x;y)}_\infty \leq \opnorm{d\xi (x;\cdot)} < \frac{\tau}{6t_n}$. Choose $t_n$ large enough such that $\norm{d\xi (x;y)}_\infty < 2$ on $\overline{B_1 (0)} \times \overline{B_1 (0)}$. Using the rule on partial derivatives and the chain rule with these estimates, we compute for $(x,y) \in \overline{B_1 (0)} \times \overline{B_1 (0)}$: 
	\begin{align*}
	 \norm{d(\eta \diamond \xi) (x;y)}_\infty \stackrel{\eqref{eq: diamond}}{=}& \norm{df_n(x, \xi (x) , \eta (F_\xi (x)), y, d \xi (x,y), d\eta (F_\xi (x), dF_\xi (x,y)))}_\infty  \\
	&\leq \norm{d_1f_n (x, \xi (x), \eta (F_\xi (x)), y)}_\infty + \opnorm{d f_n (x, \xi (x), \eta (F_\xi (x)); \cdot)} \cdot \norm{d\xi (x;y)}_\infty \\
	&\quad \  + \opnorm{df_n (x, \xi (x), \eta (F_\xi (x)); \cdot)} \cdot \opnorm{d\eta (F_\xi (x);\cdot)} \cdot \norm{dF_\xi (x;y)}_\infty \\
	&< \frac{\tau}{3} + \underbrace{\opnorm{df_n (x, \xi (x), \eta (F_\xi (x)); \cdot)}}_{\leq t_n} (\underbrace{\norm{d\xi (x;y)}_\infty}_{\leq d\mu_n} + 2 \underbrace {\opnorm{d\eta (F_\xi (x);\cdot)}}_{\leq d\mu_n})\\
        &\leq \frac{\tau}{3} + \frac{\tau}{6} + \frac{\tau}{3} \leq \tau.
	\end{align*}
 We derive $\norm{\tfrac{\partial}{\partial x_j} (\eta \diamond \xi) (x)}_\infty < \tau$ for $x \in \overline{B_1 (0)}$ and $j \in \set{ 1, 2 ,\ldots, d}$ and thus $\norm{\eta \diamond \xi}_{\overline{B_1 (0)}, 1}  < \tau$.
 \end{proof}

\begin{Alem}\label{lem: comp:sm} \no{lem: comp:sm}
 Consider the open zero neighborhoods $\hH_R^{\Omega_{5,K_5}}$ as in Construction \ref{con: comp:loc}. The maps 
	\begin{align*}
	 &c \colon \hH^{\Omega_{5,K_5}}_R \times \hH^{\Omega_{5,K_5}}_R \rightarrow \vect{\Omega_{2,K_5}},\ (X, Y) \mapsto X \diamond Y \\ 
	 &\iota \colon \hH_R^{\Omega_{5,K_5}} \rightarrow \vect{\Omega_{\frac{5}{4},K_5}} , \ X \mapsto X^*
	\end{align*}
 are smooth.
\end{Alem}

\begin{proof}
 Let $I$ be the finite set indexing $\fF_5 (K_5)$. Following Definition \ref{defn: top:vect} and the definition of $\Omega_{r,K_5}$, the topology on $\vect{\Omega_{r,K_5}}, r \in [1,5]$ is defined via the linear embedding with closed image 
	\begin{displaymath}
	 \Gamma_r \colon \vect{\Omega_{r,K_5}} \rightarrow \prod_{k \in I} C^\infty (V_{r,k}, \RR^d) = \bigoplus_{k \in I} C^\infty (V_{r,k}, \RR^d).
	\end{displaymath}
 Therefore the maps $p_k^r \coloneq \vect{\Omega_{r,K_5}} \rightarrow C^\infty (V_{r,k}, \RR^d) , p_k^r (X) \coloneq X_{\kappa_k}|_{V_{r,k}}, k \in I$ define a patchwork for $\vect{\Omega_{r,K_5}}$ indexed by $I$. Define $$p \colon \vect{\Omega_{5,K_5}} \times \vect{\Omega_{5,K_5}} \rightarrow \bigoplus_{k \in I} C^\infty (V_{5,k} ,\RR^d) \times C^\infty (V_{5,k} ,\RR^d), (X, Y) \mapsto ((p^5_k \times p_k^5) (X, Y))_{k \in I}.$$ 
 Recall that finite products coincide with direct sums in the category of locally convex vector spaces. The universal property of the direct sum therefore assures that the map 
	\begin{align*}
	 L \colon \bigoplus_{k \in I}  C^\infty (V_{5,k} ,\RR^d) \times C^\infty (V_{5,k} ,\RR^d) &\rightarrow \left(\bigoplus_{k \in I}  C^\infty (V_{5,k} ,\RR^d)\right) \times \left(\bigoplus_{k \in I} C^\infty (V_{5,k} ,\RR^d) \right)\\
	 (X_k,Y_k)_{k \in I} &\mapsto \left((X_k)_{k\in I},(Y_k)_{k \in I}\right)
	\end{align*}
 is an isomorphism of locally convex spaces. Furthermore, $L \circ p = \Gamma_5 \times \Gamma_5$ holds. As $\Gamma_5$ is an embedding with closed image, the map $\Gamma_5 \times \Gamma_5$ is a linear embedding with closed image (identifying the domain of $\Gamma_5$ via the embedding with a closed subspace of the codomain of $\Gamma_5$ this follows from \cite[Ch.\ II, \S 2, No.\ 5 Proposition 8 and Corollary 1]{bourbaki1987}). We conclude that $p$ is an embedding with closed image and the family $(p_k^5 \times p_k^5)_{k \in I}$ yields a patchwork for $\vect{\Omega_{5,K_5}} \times \vect{\Omega_{5,K_5}}$.\\
 We claim that the maps $c$ and $\iota$ are patched mappings which are smooth on the patches. If this is true, then the assertion follows from Proposition \ref{prop: pat:loc}. Proceed in two steps and prove the claim first for the map $c$: \\
 Recall from Construction \ref{con: comp:loc} that $\hH_R^{\Omega_{5,K_5}} = \bigcap_{n \in I} (\theta^{\Omega_{5,K_5}}_{\kappa_n})^{-1} (\hH_n)$ holds. Here each of the sets $\hH_n$ is an open neighborhood of the zero-map with $\hH_n \subseteq C^\infty (\kappa_n^{-1}, \RR^d)^{-1} (\hH_n') = C^\infty (\kappa_n, \RR^d) (\hH_n')$ and $\hH_n' \subseteq C^\infty (B_5 (0),\RR^d)$. We define maps 
	\begin{align*}
	 h_n \colon \hH_n' \times \hH_n' &\rightarrow C^\infty (B_2 (0),\RR^d), (\eta, \xi) \mapsto \eta \diamond \xi\\
	 c_n \colon \hH_n \times \hH_n &\rightarrow C^\infty (V_{2,n}, \RR^d),\\ 											
             (X, Y) &\mapsto C^\infty (\kappa_n|_{V_{2,n}} , \RR^d) \circ h_n \circ (C^\infty (\kappa_n^{-1} , \RR^d) \times C^\infty (\kappa_n^{-1} , \RR^d)) (X,Y).
	\end{align*}
 Observe that by Step 1 in Construction \ref{con: comp:loc}, each map $c_n$ maps the zero map $(0,0) \in \hH_n \times \hH_n$ to $0 \in C^\infty (V_{2,n},\RR^d)$. From the definition of $c$ and the identity \eqref{eq: glo:loc}, a trivial computation yields the identity $c_n \circ (p_n^5 \times p_n^5) = p_n^2 \circ c$ for each $n \in I$. Therefore $c$ is a patched mapping whose compatible family is $(c_n)_{n \in I}$. By Proposition  \ref{prop: pat:loc}, the first part of the claim will hold if each $c_n$ is a smooth map. However, $c_n$ will be smooth if and only if $h_n \colon \hH_n' \times \hH_n' \rightarrow C^\infty (B_2 (0) , \RR^d) , (\eta , \xi) \mapsto \eta \diamond \xi$ is smooth, since $C^\infty (\kappa_n^{-1}, \RR^d)$ and $C^\infty (\kappa_n, \RR^d)$ are mutually inverse isomorphisms of topological vector spaces by \cite[Lemma A.1]{hg2004}. Fix $n \in I$ and prove that $h_n$ is a smooth map:\\
 To this end, recall the constants $\ve_n, \delta_n$ obtained in Construction \ref{con: comp:loc}. By Lemma \ref{lem: exp:loc}, we may consider the smooth maps 
	\begin{align*}
	 e_n& \colon B_4 (0) \times B_{\ve_n} (0) \rightarrow \RR^d, (x,y) \mapsto \exp_n (x,y)\\
	 a_n& \colon B_4 (0) \times B_{\delta_n} (0) \rightarrow B_{\ve_n} (0), a_n(x,y) \coloneq b_n(x,x+y).
	\end{align*}
 By \cite[Proposition 4.23 (a)]{hg2004}, these maps induce smooth push-forward maps 
  \begin{align*}
   e_{n*} &\colon \lfloor \overline{B_3 (0)}, B_{\ve_n} (0)\rfloor_\infty \rightarrow C^\infty (B_3 (0) , \RR^d), e_{n*} (\gamma) (x) \coloneq e_n (x,\gamma(x))\\
   a_{n*} &\colon \lfloor\overline{B_2 (0)}, B_{\delta_n} (0) \rfloor_\infty \rightarrow C^\infty (B_2 (0) , \RR^d), a_{n*} (\eta) (x) \coloneq a_n (x,\eta (x)),
  \end{align*}
 where $\lfloor \overline{B_3 (0)}, B_{\ve_n} (0)\rfloor_\infty  \subseteq C^\infty (B_4 (0) , \RR^d)$ and $\lfloor \overline{B_2 (0)}, B_{\delta_n} (0)\rfloor_\infty \subseteq C^\infty (B_{\frac{21}{10}} (0), \RR^d)$ are open sets. Recall from Construction \ref{con: comp:loc} that $\hH_n'$ is a subset of an open set $\nN_n$ which has been constructed by an application of Lemma \ref{lem: mb:exp}. Hence $\eta \in H_n'$ satisfies Lemma \ref{lem: mb:exp} (a), whence $\eta (\overline{B_3 (0)}) \subseteq B_{\ve_n} (0)$ holds. In other words, $\hH_n' \subseteq \lfloor \overline{B_3 (0)}, B_{\ve_n} (0)\rfloor_\infty$ is satisfied (after restricting to $B_4 (0)$, which we suppress in the notation). By definition, $e_{n*} (\eta) = F_\eta$ with $F_\eta$ as defined in Lemma \ref{lem: mb:exp}. Furthermore, applying the estimate \eqref{eq: exp:zero}, we obtain $e_{n*} (\eta) \in \lfloor \overline{B_2 (0)} , B_3 (0)\rfloor_\infty$. By \cite[Lemma 11.4]{hg2004}, there is a smooth composition map
	\begin{displaymath}
	 \Theta \colon C^\infty (B_3 (0) , \RR^d) \times \lfloor \overline{B_2 (0)} , B_3 (0) \rfloor_\infty \rightarrow C^\infty (B_{\frac{21}{20}} (0) , \RR^d), (f , g) \mapsto f\circ g|_{ B_{\frac{21}{20}} (0)}, 
	\end{displaymath}
 where $\lfloor \overline{B_2 (0)}, B_3 (0)\rfloor \subseteq C^\infty (B_3 (0),\RR^d)$. Hence, we conclude that we may compose $\Theta$ and $(e_{n*} \times e_{n*})$ to obtain a smooth map $\Theta \circ (e_{n*} \times e_{n*}) \colon \hH_n' \times \hH_n' \rightarrow C^\infty (B_{\frac{21}{20}} (0) , \RR^d)$. By definition of $\hH_n'$, we derive for $\eta \in H_n'$ the estimate $F_\eta (x) \in B_{\frac{\delta_n}{2}} (x)$ for $x \in B_3 (0)$ (see Lemma \ref{lem: mb:exp} (a)). Thus $\Theta (e_{n*} (\eta) , e_{n*} (\xi))(x) -x \in B_{\delta_n} (0)$ holds for $x \in \overline{B_2 (0)}, \eta , \xi \in \hH_n'$. \\
 Combine the identity \eqref{eq: diamond} with the definition of $f_n$ in Lemma \ref{lem: exp:loc} (c) to deduce the identity  
	\begin{displaymath}
	 h_n (\eta , \xi) = a_{n*} (\Theta (e_{n*} (\eta), e_{n*} (\xi)) - \id_{B_{\frac{21}{20}} (0)}).
	\end{displaymath}
 We conclude that $h_n$ is smooth as composition of smooth maps. Summing up, this proves the first part of the claim.\\
 As a second step, we construct a compatible family for $\iota$. To this end, define maps  
	\begin{align*}
	 i_n &\colon \hH_n' \rightarrow C^\infty (B_{\frac{5}{4}} (0), \RR^d) , \xi \mapsto \xi^*|_{B_{\frac{5}{4}}(0)} \\ 
	\iota_n &\colon \hH_n \rightarrow C^\infty (V_{\frac54,n}, \RR^d) , X \mapsto C^\infty (\kappa_n|_{V_{\frac54,n}} , \RR^d) \circ i_n \circ C^\infty (\kappa_n^{-1} , \RR^d).
	\end{align*}
 From the identity \eqref{eq: inv:loc}, we derive $p_n^{\frac{5}{4}} \iota = \iota_n p^5_n$. Hence $\iota$ is a patched mapping and we have to prove that each $\iota_n$ is smooth. Again $\iota_n$ will be smooth if $i_n$ is smooth.\\
 Recall that $\hH_n' \subseteq \nN_n$ holds for an open set $\nN_n \subseteq C^\infty (B_5 (0), \RR^d)$ with the properties of the set $\nN$ in Lemma \ref{lem: mb:exp}. Hence the map $I_n \colon \nN_n \rightarrow C^\infty (B_2 (0), \RR^d) , \xi \mapsto \xi^*|_{B_2 (0)}$ is smooth by Lemma \ref{lem: mb:exp} (f). Let $\lambda \colon B_{\frac{5}{4}} (0) \hookrightarrow B_{2} (0)$ be the canonical inclusion. The pullback $C^\infty (\lambda, \RR^d)$ is continuous linear, hence smooth. Finally, the identity $i_n = C^\infty (\lambda , \RR^d) \circ I_n|_{\hH_n'}$ assures that $i_n$ is smooth.
\end{proof}
 \cleardoublepage \thispagestyle{empty}
\section{Maps of orbifolds}\label{sect: maps}

In this section, we recall the notion of an orbifold map in local charts which was introduced in \cite{pohl2010} (cf.\ Section \ref{sect: haef} for details on orbifolds). Our exposition follows \cite{pohl2010} and we repeat basic facts for the readers convenience. Orbifold maps in the sense discussed here correspond to maps in a category of groupoids. Our notion of orbifold map developed here is thus equivalent to other types of orbifold maps which are equivalent to maps in the associated groupoid category (cf.\ \cite{chen2006} for the so called Chen-Ruan good map and \cite{alr2007} for the Moerdijk-Pronk strong map, respectively).

\subsection{(Quasi-)Pseudogroups} 

In this section we let $M$ be a smooth manifold.
 
\begin{nota}[Transitions]\label{nota: trans} \no{nota: trans}
 A \ind{transition}{transition on $M$} is a diffeomorphism $f \colon U \rightarrow V$, where $U$, $V$ are open subsets of $M$. Notice that the empty map $\emptyset \rightarrow \emptyset$ is a transition on $M$.\\ The product of two transitions $f \colon U \rightarrow V$, $g \colon U' \rightarrow V'$ is the transition 
		\begin{displaymath}
		 f|^{f(U\cap V')} \circ g|_{g^{-1} (U \cap V')} \colon g^{-1} (U \cap V') \rightarrow f(U\cap V'), x \mapsto f(g(x)).
		\end{displaymath}
 The inverse of $f$ is the inverse of $f$ as a function. If $f \colon U \rightarrow V$ is a map, we denote by $\dom f$\glsadd{dom} the \emph{domain} of $f$ and $\cod f$\glsadd{cod} the \emph{codomain} of $f$. For $x \in \dom f$ we denote by $\germ_x f$ the \ind{}{germ of $f$ at $x$} and the set of all transitions of $M$ by $\aA (M)$.
\end{nota}

\begin{defn}[Pseudogroup] 
 A \ind{}{pseudogroup} on $M$ is a subset $P \subseteq \aA (M)$ which is closed under products and inversion of transitions. We call $P$ a \ind{pseudogroup!full}{full pseudogroup}, if for every open subset $U \subseteq M$ the transition $\id_U$ is contained in $P$. A full pseudogroup is called \ind{pseudogroup!complete}{complete} if it satisfies 
 \begin{compactitem}
  \item[(\emph{Gluing Property})] If $f \in \aA (M)$ and there is an open cover $(U_i)_{i\in I}$ of $\dom f$ such that $f|_{U_i} \in P$ for all $i \in I$, then $f$ is an element of $P$.
 \end{compactitem}
  The pseudogroup $P$ is \ind{pseudogroup!closed under restrictions}{closed under restrictions}, if for any $f\in P$ and open set $U \subseteq \dom f$, the map $f|^{f(U)}_{U} \colon U \rightarrow f(U)$ is in $P$. For example, every full pseudogroup is closed under restrictions.
\end{defn}

\begin{defn}[Quasi-Pseudogroup]\label{defn: qpg} \no{defn: qpg}
 A subset $P$ of $\aA (M)$ is called a \ind{}{quasi-pseudogroup} on $M$ if the following properties are satisfied: 
	\begin{compactenum}
	 \item For each $f \in P$ and $x \in \dom f$, there exist an open set $U$ with $x \in U \subseteq \dom f$ and $g \in P$ together with an open set $V$ such that $f(x) \in V \subseteq \dom g$ and $$(f|_{U})^{-1} = g|_{V}.$$
	 \item If $f,g \in P$ and $x \in f^{-1} (\cod f \cap \dom g)$, then there exists $h \in P$ and an open neighborhood $U \subseteq f^{-1} (\cod f \cap \dom g) \cap \dom h$ of $x$ with $g \circ f|_{U} = h|_{U}$. 
	\end{compactenum}
\end{defn}

Thus inversions and compositions of elements in a quasi-pseudogroup are only required to correspond locally to other elements in the quasi-pseudogroup. For pseudogroups, inverses and composites globally belong to the pseudogroup. Quasi-pseudogroups are designed to work with the germs of their elements. In general, quasi-pseudogroups may be thought of as generators for pseudogroups in the following sense: 

\begin{defn}
 Let $P$ be a pseudogroup on $M$ which satisfies the gluing property and is closed under restrictions. The pseudogroup $P$ is \ind{pseudogroup!generating subset}{generated} by a set $A \subseteq \aA (M)$  if $A \subseteq P$ holds and for each $f \in P$ and $x \in \dom f$ there exists $g \in A$ and an open set $U \subseteq \dom f \cap \dom g$ with $x \in U$ and $f|_{U} = g|_{U}$. Then $P$ is uniquely determined by $A$.\\
 Consider a subset $B$ of $\aA (M)$. If there exists a unique pseudogroup $Q$ on $M$ which satisfies the gluing property, is closed under restrictions and generated by $B$, then we say that $B$ \emph{generates} $Q$.
\end{defn}

\begin{rem}
 \begin{compactenum}
  \item The set $\aA (M)$ is a pseudogroup. Each pseudogroup is a quasi-pseudogroup.
  \item Each quasi-pseudogroup generates a unique pseudogroup which satisfies the gluing property and is closed under restrictions. Vice versa each generating set for such a pseudogroup is necessarily a quasi-pseudogroup.
 \end{compactenum}
\end{rem}

\subsection{Charted orbifold maps}
 In this section, we let $(Q ,\uU )$ and $(Q', \uU')$ be orbifolds. Morphisms of orbifolds will be constructed in several steps, since they arise as equivalence classes of certain objects:

\begin{defn}
 Let $\vV \coloneq \setm{(V_i , G_i , \pi_i)}{i \in I}$ be a representative of $\uU$. We abbreviate the disjoint union of the the chart domains of elements in $\vV$ with 
		\begin{displaymath}
		 V \coloneq \coprod_{i \in I} V_i \quad \text{and define} \quad \pi \colon V \rightarrow Q , x \mapsto \pi_i (x) \text{ for } x \in V_i.
		\end{displaymath}
 Then the subset \glsadd{PSI}
	\begin{displaymath}
	 \PSI (\vV) \coloneq \setm{f \in \aA (V)}{\pi \circ f = \pi|_{\dom f}}
	\end{displaymath}
 of the set of all transitions on $V$ is a complete pseudogroup on $V$ which is closed under restrictions.
\end{defn}

The last definition may be used to associate to each orbifold an \'{e}tale Lie groupoid (as is explained in \cite[2.9 and 2.10]{pohl2010}). Since we are not interested in the correspondence of orbifolds and Lie groupoids, we will not pursue this relation any further. However this relation was invaluable to derive the notion of orbifold map introduced in this section. We refer to \cite{pohl2010} for further details. 

\begin{defn}
 Let $f \colon Q \rightarrow Q'$ be a continuous map. Consider two orbifold charts $(V,G, \pi) \in \uU$ and $(V',G',\pi') \in \uU'$. A smooth map $f_V \colon V \rightarrow V'$ is called \ind{local lift}{local lift of $f$ with respect to $(V,G,\pi)$ and $(V',G',\pi')$} if $\pi' \circ f_V = f \circ \pi$ holds. In this case, $f_V$ is also called a local lift of $f$ at $q$ for each $q \in \pi(V)$.
\end{defn}

\begin{defn}[Representative of an orbifold map] \label{defn: rep:ofdm}\no{defn: rep:ofdm}
 A \ind{orbifold map!representative}{representative of an orbifold map} from an orbifold $(Q ,\uU )$ to an orbifold $(Q', \uU')$ is a tuple 
	\begin{displaymath}
	 \hat{f} \coloneq (f, \set{f_i}_{i \in I} , P,\nu)
	\end{displaymath}
 where \begin{compactitem}
        \item[(R1)] $f \colon Q \rightarrow Q'$ is a continuous map, 
	\item[(R2)] for each $i \in I$, the map $f_i \colon V_i \rightarrow V_i'$ is a local lift of $f$ with respect to orbifold charts $(V_i, G_i, \pi_i) \in \uU$, $(V_i',G_i',\pi_i') \in \uU'$ such that 
			\begin{displaymath} 
			 \bigcup_{i\in I} \pi_i (V_i) = Q
			\end{displaymath} 
	and $(V_i , G_i , \pi_i) \neq (V_j , G_j , \pi_j)$ holds for $i,j \in I, \ i \neq j$,
	\item[(R3)] $P$ is a quasi-pseudogroup which consists of changes of charts of the orbifold atlas 
			\begin{displaymath}
			 \vV \coloneq \setm{(V_i , G_i , \pi_i)}{i\in I}
			\end{displaymath}
	of $(Q,\uU)$ and generates $\PSI (\vV)$,
	\item[(R4)] Set $F \coloneq \coprod_{i \in I} f_i \colon V = \coprod_{i \in I} V_i \rightarrow \coprod_{i \in I} V_i', x \mapsto f_i(x)$ if $x \in V_i$. Choose any orbifold atlas $\vV' \in \uU'$ which contains the set $\set{(V_i',G_i',\pi_i')}_{i \in I}$. Then $\nu \colon P \rightarrow \PSI (\vV')$ is a map which assigns to each $\lambda \in P$ a change of charts  
		\begin{displaymath}
		 \nu (\lambda) \colon (W',H',\chi') \rightarrow (V', G', \vp')
		\end{displaymath}
	between orbifold charts in $\vV'$ such that the following properties are satisfied
		\begin{compactenum}
		 \item[a)] $F \circ \lambda = \nu (\lambda) \circ F|_{\dom \lambda}$ for all $\lambda \in P$,
		 \item[b)] for all $\lambda, \mu \in P$ and all $x \in \dom \lambda \cap \dom \mu$ with $\germ_x \lambda = \germ_x \mu$ we have 
			\begin{displaymath}
			 \germ_{F (x)} \nu (\lambda) = \germ_{F (x)} \nu(\mu),
			\end{displaymath}
		 \item[c)] for all $\lambda , \mu \in P$ and all $x \in \lambda^{-1} (\cod \lambda \cap \dom \mu)$ we have 
			\begin{displaymath}
			 \germ_{F (\lambda (x))} \nu (\mu) \cdot \germ_{F (x)} \nu (\lambda) = \germ_{F (x)} \nu (h)
			\end{displaymath}
		 where $h$ is an element of $P$ such that there is an open set $U$ with
			\begin{displaymath}
			 x \in U \subseteq \lambda^{-1}(\cod \lambda \cap \dom \mu) \cap \dom h
			\end{displaymath}
		 and $\mu \circ \lambda|_{U} = h|_{U}$,
		 \item[d)] for all $\lambda \in P$ and $x \in \dom \lambda$ such that there is an open set $x \in U \subseteq \dom \lambda$ with $\lambda|_{U} = \id_U$ we have $\germ_{F (x)} \nu (\lambda ) = \germ_{F (x)} \id_{U'}$ where $U' \coloneq \coprod_{i \in I} V_i'$.
		\end{compactenum}
       \end{compactitem}
The orbifold atlas $\vV$ is called the \ind{}{domain atlas} of the representative $\hat{f}$, and the set $\setm{(V_i', G_i' , \pi_i')}{i \in I}$ is called the \ind{}{range family} of $\hat{f}$. Note that the range family is not necessarily indexed by $I$. Moreover, the mapping $\nu$ does not depend on the choice of $\vV'$ since it takes its values in $\bigcup_{(i,j) \in I \times I} \CH{V_i'}{V_j'}$ (cf.\ Notation \ref{nota: chch:abr} below). The continuous map $f$ will sometimes be called the \ind{orbifold map!underlying map}{underlying map} of the representative $\hat{f}$. The map $f$ may not be chosen arbitrarily. As \cite[Example 4.5]{pohl2010} shows, it is not even sufficient to require that $f$ be a homeomorphism, to assure that there is a representative $\hat{f}$ with underlying map $f$. 
\end{defn}
The technical condition in (R2) that two orbifold charts in $\vV$ be distinct is required, because in several places $I$ is used as an index set for $\vV$ (cf. property (R3)). \\
In view of Definition \ref{defn: rep:ofdm}, it is useful to have a shorthand for the changes of charts associated to a given orbifold atlas. We fix the following notation.

\begin{nota}\label{nota: chch:abr} \no{nota: chch:abr}
 Let $\vV = \setm{(V_i , G_i, \psi_i)}{i \in I}$ be a representative of $\uU$. Recall the notation for the set of all changes of charts between two orbifold charts (first introduced in Lemma \ref{lem: pr:tofd} (b)): \glsadd{CH}
		\begin{displaymath}
	         \CH{V_i}{V_j} \coloneq \setm{\lambda \colon V_i \supseteq \dom \lambda \rightarrow \cod \lambda \subseteq V_j}{\lambda \text{ is a change of charts}}
		\end{displaymath}
 We define the \ind{orbifold chart!change of charts!set of all changes of charts}{set of all changes of charts} of the atlas $\vV$ via \glsadd{ch}
	\begin{displaymath}
	 \Ch_\vV \coloneq \setm{\lambda \colon V_i \supseteq \dom \lambda \rightarrow \cod \lambda \subseteq V_j}{\lambda \text{ is a change of charts and } i,j \in I} = \bigcup_{(i,j) \in I\times I} \CH{V_i}{V_j}.
	\end{displaymath}
 Observe that $\Ch_\vV$ is a (quasi-)pseudogroup, which generates $\PSI (\vV)$.
\end{nota}

\begin{defn}\label{defn: eq:rofdm} \no{defn: eq:rofdm}
 Let $\hat{f} \coloneq (f, \set{f_i}_{i \in I} , P_1 ,\nu_1)$ and $\hat{g} \coloneq (g, \set{g_i}_{i \in I} , P_2 , \nu_2)$ be two representatives of orbifold maps with the same domain atlas $\vV$ representing the orbifold structure $\uU$ on $Q$ and both range families being contained in the orbifold atlas $\vV'$ of $(Q', \uU')$. Set $F \coloneq \coprod_{i \in I} f_i$. We say that $\hat{f}$ is \emph{equivalent} to $\hat{g}$ if $f=g$, $f_i = g_i$ for all $i \in I$ and 
		\begin{displaymath}
		 \germ_{F (x)} \nu_1 (\lambda_1) = \germ_{F (x)} \nu_2 (\lambda_2)
		\end{displaymath}
 holds for all $\lambda_1 \in P_1, \lambda_2 \in P_2, x \in \dom \lambda_1 \cap \dom \lambda_2$ with $\germ_x \lambda_1 = \germ_x \lambda_2$. This defines an equivalence relation. The equivalence class of $\hat{f}$ will be denoted by  
		\begin{displaymath}
		 (f, \set{f_i}_{i \in I} , [P_1 , \nu_1]).
		\end{displaymath}
 By abuse of notation, we denote by $\hat{f}$ the equivalence class $[\hat{f}]$ of the representative $\hat{f}$, if it is clear that we refer to equivalence classes. The equivalence class of the representative $\hat{f}$ is called \ind{orbifold map}{orbifold map with domain atlas $\vV$ and range atlas $\vV'$}, in short \emph{orbifold map with $(\vV , \vV')$} or, if the specific atlases are not important, a \ind{orbifold map!charted}{charted orbifold map}. Define $\Orb{\vV , \vV'}$ to be \ind{set!of charted orbifold maps}{the set of all orbifold maps with $(\vV , \vV')$}\glsadd{Orb}. To shorten our notation we denote an element $\hat{h} \in \Orb{\vV ,\vV'}$ by $\vV \xra{\hat{h}} \vV'$.
\end{defn}

\begin{rem}\label{rem: psgp} \no{rem: psgp}
 \begin{compactenum}
  \item The results of \cite{pohl2010} apply to the class of second countable orbifolds and the wider class of paracompact orbifolds. We only required orbifolds to be paracompact.  Second countability of all spaces seems to be a standard requirement in the theory of groupoids (cf.\ \cite{follie2003}). However, \cite{cgo2006, msnpc1999} and the survey article by Lerman \cite{ofdstack2008} outline the theory of Lie-groupoids for non second countable manifolds. In particular, the article by Lerman indicates that all desirable properties on the groupoid side are preserved for paracompact orbifolds and manifolds. Hence we require only the weaker condition.
  \item In Definition \ref{defn: rep:ofdm} we used quasi-pseudogroups instead of the pseudogroups $\Ch_\vV$ or $\PSI (\vV)$ since, in general, a quasi-pseudogroup $P$ will be much smaller (sometimes even finite). Observe the following facts, whose proofs we omit here: 
	\begin{compactenum}
	 \item Let $(f, \set{f_i}_{i\in I}, P,\nu)$ be a representative of an orbifold map. Replacing $P$ with a quasi-pseudogroup $P'$ whose elements arise as restrictions of maps in $P$ (if  necessary reducing them to open neighborhoods which are stable with respect to the group action), one may replace $\nu$ with a map $\nu'$ which maps each element in $P'$ to an open embedding in the range family. The pair $(P',\nu')$ may be chosen such that $(f, \set{f_i}_{i\in I}, P,\nu)$ and $(f, \set{f_i}_{i\in I}, P',\nu')$ are in the same equivalence class. 
	 \item Consider a representative of an orbifold map $\hat{f} \colon (Q,\uU) \rightarrow \mM$, where $\mM$ is a connected manifold (without boundary) and range family of the charted map is the atlas $(M,\set{\id_M}, \id_M)$. The map $\nu$ may then be chosen as the map taking each $h \in P$ to $\id_\mM$. 
	 \end{compactenum}
 \end{compactenum}
\end{rem}

\subsection{The identity morphism}

In this section, we construct the identity morphism in the category of reduced orbifolds.
\begin{defn}
 Let $f \colon Q \rightarrow Q'$ be a continuous map between orbifolds $(Q,\uU),(Q',\uU')$. Suppose $f_V$ is a local lift with respect to the orbifold charts $(V,G,\pi) \in \uU$ and $(V',G',\pi') \in \uU'$. Consider embeddings of orbifold charts in $\uU$ and $\uU'$,respectively,
		\begin{displaymath}
		 \lambda \colon (W,K,\chi ) \rightarrow (V,G,\pi) \quad \text{and} \quad \mu \colon (W',K',\chi') \rightarrow (V' ,G', \pi'),
		\end{displaymath}
 such that $f_V (\lambda (W)) \subseteq \mu (W')$ holds. Then the map 
	\begin{displaymath}
	 g \coloneq \mu^{-1} \circ f_V \circ \lambda \colon W \rightarrow W'
	\end{displaymath}
 is a local lift of $f$ with respect to $(W,K,\chi)$ and $(W',K',\chi')$. We say $f_V$ induces the local lift $g$ with respect to $\lambda$ and $\mu$ and call $g$ \ind{orbifold map!induced lift}{induced lift of $f$ with respect to $f_V$, $\lambda$ and $\mu$}.
\end{defn}

\begin{prop}[\hspace{-0.5pt}{\cite[Proposition 5.3]{pohl2010}}] \label{prop: ll:id} \no{prop: ll:id}
 Let $(Q,\uU)$ be an orbifold and $f_V$ be a local lift of $\id_Q$ with respect to $(V,G,\pi), (V',G',\pi') \in \uU$. For each $v \in V$ there exists a restriction $(S, G_S, \pi|_{S})$ of $(V,G,\pi)$ with $v \in S$ and a restriction $(S',G'_{S'},\pi'|_{S'})$ of $(V',G',\pi')$ such that $f_V|_{S}^{S'}$ is diffeomorphism which is a change of charts from $(S,G_S,\pi|_{S})$ to $(S',G'_{S'},\pi'|_{S'})$. In particular, $f_V|_{S}$ induces the identity $\id_{S}$ with respect to the embeddings of orbifold charts $\id_S$ and $(f_V|_{S}^{S'})^{-1}$.
\end{prop}

Proposition \ref{prop: ll:id} shows that every local lift of the identity $\id_Q$ is a local diffeomorphism (but in general it need not be a global diffeomorphism as \cite[Example 5.4]{pohl2010} shows).

\begin{prop}[\hspace{-0.5pt}{\cite[Proposition 5.5]{pohl2010}}]\label{prop: lift:unit} \no{prop: lift:unit}
 Let $(Q,\uU)$ be an orbifold and $\set{f_i}_{i \in I}$ a family of local lifts of $\id_Q$ which satisfies (R2). Then there exists a pair $(P, \nu)$ such that $(\id_Q, \set{f_i}_{i \in I} , P,\nu)$ is a representative of an orbifold map on $(Q ,\uU)$. The pair $(P,\nu)$ is unique up to equivalence of representatives of orbifold maps.
\end{prop}

\begin{prop}[{\hspace{-0.5pt}\cite[Proposition 5.6]{pohl2010}}]\label{prop: eq:ofd} \no{prop: eq:ofd}
 Let $Q$ be a topological space and suppose $\uU$ and $\uU'$ are orbifold structures on $Q$. Consider a charted orbifold map 
	\begin{displaymath}
	 \hat{f} = (\id_Q, \set{f_i}_{i\in I} , [P,\nu ]) 
	\end{displaymath}
 such that the domain atlas $\vV$ is a representative of $\uU$ and the range family $\vV'$, which is an orbifold atlas, is a representative of $\uU'$. If $f_i$ is a local diffeomorphism for each $i \in I$, then $\uU = \uU'$ holds, i.e.\ the orbifolds coincide.
\end{prop}

\begin{defn}\label{defn: idmor} \no{defn: idmor}
 Let $(Q,\uU)$ be an orbifold and $\hat{f} = (f, \set{f_i}_{i\in I} , [P,\nu])$ be a charted orbifold map whose domain atlas is a representative of $\uU$. The representative $\hat{f}$ is called \ind{orbifold map!lift of the identity}{lift of the identity $\ido{}$} if $f = \id_Q$ holds and $f_i$ is a local diffeomorphism for each $i \in I$. We also say that $\hat{f}$ is a representative of $\ido{}$. The set of all lifts of $\ido{}$ is the \ind{orbifold map!identity morphism}{identity morphism} $\ido{}$\glsadd{ido} of $(Q,\uU)$.
\end{defn}


\subsection{Composition of charted orbifold maps}

\begin{con}\label{con: cp:chom} \no{con: cp:chom}
Let $(Q,\uU)$, $(Q',\uU')$ and $(Q'',\uU'')$ be orbifolds, and 
		\begin{displaymath}
		 \vV \coloneq \setm{(V_i , G_i ,\pi_i)}{i \in I}, \quad \vV' \coloneq \setm{(V_j' , G_j' ,\pi_j')}{j \in J}
		\end{displaymath}
 be representatives of $\uU$ and $\uU'$, respectively, where $\vV$ is indexed by $I$ and $\vV'$ by $J$. Furthermore, let $\vV'' \in \uU''$. Consider charted orbifold maps 
	\begin{displaymath}
	 \hat{f} = (f , \set{f_i}_{i \in I} , [P_f,\nu_f]) \in \Orb{\vV ,\vV'}
	\end{displaymath}
 and 
	\begin{displaymath}
		\hat{g} = (g, \set{g_j}_{j \in J} , [P_g , \nu_g]) \in \Orb{\vV' ,\vV''}.
	\end{displaymath}
 Define $\alpha \colon I \rightarrow J$ to be the unique map such that for each $i \in I$, $f_i$ is a local lift of $f$ with respect to $(V_i,G_i,\pi_i)$ and $(V_{\alpha (i)}' , G_{\alpha (i)}' , \pi_{\alpha (i)}')$. We define the composition of $\hat{g}$ and $\hat{f}$: 
	\begin{displaymath}
	 \hat{g} \circ \hat{f} \coloneq \hat{h} = (h, \set{h_i}_{i \in I}, [P_h,\nu_h]) \in \Orb{\vV,\vV''}
	\end{displaymath}
 is given by $h \coloneq g \circ f$ and $h_i \coloneq g_{\alpha (i)} \circ f_i$ for all $i \in I$. To construct a representative $(P_h ,\nu_h)$ of $[P_h , \nu_h]$ fix representatives $(P_f, \nu_f)$ and $(P_g , \nu_g)$ of $[P_f , \nu_f]$ and $[P_g, \nu_g]$, respectively. Consider $\mu \in P_f$ with $\dom \mu \subseteq V_i$, $\cod \mu \subseteq V_j$ for the orbifold charts $(V_i ,G_i , \pi_i)$ and $(V_j , G_j, \pi_j)$ in $\vV$. Property (R4a) assures 
	\begin{displaymath}
	 f_j \circ \mu = \nu_f (\mu) \circ f_i|_{\dom \mu},
	\end{displaymath}
 where $\nu_f (\mu)$ is a change of charts in $\vV'$. For $x \in \dom \mu$ set $y_x \coloneq f_i (x) \in \dom \nu_f (\mu)$. Since $P_g$ generates $\PSI (\vV')$ we may choose $\xi_{\mu, x} \in P_g$ such that there is an open set $y_x \in U_{\mu ,x}' \subseteq \dom \xi_{\mu , x} \cap \dom \nu_f (\mu)$ and the following is satisfied:
		\begin{displaymath}
		 \xi_{\mu , x}|_{U_{\mu ,x}'} = \nu_f (\mu)|_{U_{\mu , x}'}.
		\end{displaymath}
 We may choose an open set $x \in U_{\mu ,x} \subseteq \dom \mu$ such that $f_i(U_{\mu ,x}) \subseteq U_{\mu ,x}'$ holds. By adjusting choices one may achieve that for $\mu_1,\mu_2 \in P_f$ and $x_k \in \dom \mu_k , \ k \in \set{1,2}$ we have 
	\begin{equation}\label{eq: eq:ornot}
	 \mu_1|_{U_{\mu_1 , x_1}} \neq \mu_2|_{U_{\mu_2 ,x_2}} \quad \text{or}\quad \xi_{\mu_1 , x_1} = \xi_{\mu_2 , x_2}.
	\end{equation}
 Define the quasi-pseudogroup
	\begin{displaymath}
	 P_h \coloneq \setm{\mu|_{U_{\mu ,x}}}{\mu \in P_f , \ x \in \dom \mu}
	\end{displaymath}
 and observe that it generates $\PSI (\vV)$ as $P_f$ generates $\PSI (\vV)$. As property \eqref{eq: eq:ornot} holds, we obtain a well defined map 
	\begin{displaymath}
	\nu_h \colon P_h \rightarrow \PSI (\vV''), \nu_h (\mu|_{U_{\mu ,x}}) \coloneq \nu_g (\xi_{\mu , x}).
	\end{displaymath}
 Since $\nu_g$ and $\nu_f$ satisfy the properties (R4a) - (R4d), the same holds for $\nu_h$. Furthermore, the equivalence class of $(P_h , \nu_h)$ does not depend on the choices in the construction of $P_h$ and $\nu_h$. 
\end{con}
%
%
So far, we have only explained the composition of charted orbifold maps in $\Orb{\vV ,\vV'}$ and $\Orb{\vV' , \vV''}$. Obviously we need the composition of maps in $\Orb{\vV , \vV'}$ and maps in $\Orb{\vV'' , \vV'''}$ for arbitrary $\vV'', \vV'''$. The leading idea is to construct a common refinement of the range family and the atlas $\vV''$ together with induced maps, which may then be composed as in Construction \ref{con: cp:chom}. Before we introduce the general construction, we define the notion of induced charted orbifold maps:

\begin{lemdef}[\hspace{-0.5pt}{\cite[Lemma and Definition 5.11]{pohl2010}}]\label{lemdef: ind:ofdm} \no{lemdef: ind:ofdm}
 Let $(Q,\uU)$ and $(Q',\uU')$ be orbifolds. Consider representatives 
	\begin{align*}
	 \vV &= \setm{(V_i,G_i,\pi_i)}{i \in I} \text{ of } \uU  \text{ indexed by } I \\
	 \vV' &= \setm{(V_l',G_l',\pi_l')}{l \in L} \text{ of } \uU'  \text{ indexed by } L, \text{ and a charted map}\\
	 \hat{f} &= (f, \set{f_i}_{i \in I}, [P_f,\nu_f]) \in \Orb{\vV , \vV'} . 
	\end{align*}
Define $\beta \colon I \rightarrow L$ to be the unique map such that for each $i \in I$, $f_i$ is a local lift of $f$ with respect to $(V_i,G_i,\pi_i)$ and $(V_{\beta (i)}' , G_{\beta (i)}' , \pi_{\beta (i)}')$. Suppose there are 
	\begin{compactitem}
	 \item a representative $\wW = \setm{(W_j , H_j , \psi_j)}{j \in J}$ of $\uU$, indexed by $J$,
	 \item a subset $\setm{(W_j' , H_j' , \psi_j')}{j \in J}$ of $\uU'$, indexed by $J$ (not necessarily an orbifold atlas),
	 \item a map $\alpha \colon J \rightarrow I$,
	 \item for each $j \in J$, an embedding of orbifold charts
                                                     \begin{displaymath}
	                                              \lambda_j \colon (W_j , H_j, \psi_j) \rightarrow (V_{\alpha (j)}, G_{\alpha (j)} , \pi_{\alpha (j)}) 
	                                             \end{displaymath}
	       and an embedding of orbifold charts 
                                     \begin{displaymath}
	                              \mu_j \colon (W_j' , H_j', \psi_j') \rightarrow (V_{\beta (\alpha (j))}', G_{\beta (\alpha (j))}' , \pi_{\beta (\alpha (j))}') 
	                             \end{displaymath}
		such that $f_{\alpha (j)} (\lambda_j (W_j) \subseteq \mu_j (W_j')$ holds.
	\end{compactitem}
 For each $j \in J$ we define the smooth map 
                              \begin{displaymath}
                               h_j \coloneq \mu_j^{-1} \circ f_{\alpha (j)} \circ \lambda_j \colon W_j \rightarrow W_j'.
                              \end{displaymath}
 Then the following assertions hold
 \begin{compactenum}
      \item $\ve \coloneq (\id_Q , \set{\lambda_j}_{j \in J} , [P_\ve ,\nu_\ve])$ (with $[P_\ve ,\nu_\ve]$ provided by Proposition \ref{prop: lift:unit}) is a lift of $\id_{(Q,\uU)}$.
      \item The set $\setm{(W_j',H_j', \psi_j')}{j \in J}$ and the family $(\mu_j)_{j \in J}$ may be extended to a representative 
		\begin{displaymath}
		 \wW' = \setm{(W_k', H_k' , \psi_k')}{k \in K}
		\end{displaymath}
	of $\uU'$ and a family $\set{\mu_k}_{k \in K}$ of embeddings of orbifold charts such that 
	\begin{displaymath}
	 \ve' \coloneq (\id_{Q'} , \set{\mu_k}_{k \in K}, [P_{\ve'}, \nu_{\ve'}]) \in \Orb{\wW', \vV'}
	\end{displaymath}
	(with $[P_{\ve'}, \nu_{\ve'}]$ provided by Proposition \ref{prop: lift:unit}) is a lift of the identity $\id_{(Q',\uU')}$.
       \item There is a uniquely determined equivalence class $[P_h ,\nu_h]$ such that 
	\begin{displaymath}
	 \hat{h} \coloneq (f,\set{h_j}_{j \in J} , [P_h , \nu_h]) \in \Orb{\wW , \wW'}
	\end{displaymath}
	and $\hat{f} \circ \ve = \ve' \circ \hat{h}$ holds.
      \end{compactenum}
 We say that the charted orbifold map $\hat{h}$ is induced by $\hat{f}$.
\end{lemdef}

\begin{defn}\label{defn: sim:ofd} \no{defn: sim:ofd}
 Let $(Q,\uU)$ and $(Q', \uU')$ be orbifolds. Further let $\vV_1 , \vV_2$ be representatives of $\uU$ and $\vV_1' , \vV_2'$ be representatives of $\uU'$. Suppose that $\hat{f}_i \in \Orb{\vV_i ,\vV_i'},\ i=1,2$. We call $\hat{f}_1$ and $\hat{f}_2$ \emph{equivalent} $(\hat{f}_1 \sim \hat{f}_2)$ if there are representatives $\wW$ of $\uU$ and $\wW'$ of $\uU'$ together with lifts of the identity $\ve_i \in \Orb{\wW ,\vV_i}$ and $\ve_i' \in \Orb{\wW', \vV_i'},$ respectively (for $i \in \set{1,2}$) and a map $\hat{h} \in \Orb{\wW , \wW'}$ such that the following diagram commutes 
 \begin{displaymath}
  \begin{xy}
  \xymatrix{					
								&	\vV_1 \ar[r]^{\hat{f}_1} 	& \vV_1' 	&				    \\
      \wW \ar[rrr]^{\hat{h}} \ar[ru]_{\ve_1} \ar[rd]^{\ve_2}   	&   				  	&		& \wW' \ar[lu]^{\ve_1'} \ar[ld]^{\ve_2'}  \\
      								& \vV_2 \ar[r]^{\hat{f}_2}		&  \vV_2' 	& 
  }
\end{xy}
\end{displaymath}
\end{defn}

Let $(Q,\uU)$ and $(Q', \uU')$ be orbifolds. The notion of equivalence of charted maps induces an equivalence relation on the set of all charted orbifold maps whose domain atlas is contained in $\uU$ and whose range family is contained in $\uU'$. To prove this fact, in \cite{pohl2010} the following lemmata clarify the relation of induced lifts and induced charted orbifold maps. 

\begin{lem}[{\cite[Lemma 5.13]{pohl2010}}]\label{lem: lft:rel} \no{lem: lft:rel}
 Let $(Q,\uU)$ and $(Q',\uU')$ be orbifolds and  \begin{displaymath}
                                                     \hat{f} \coloneq (f, \set{f_i}_{i \in I}, [P, \nu]) \in \Orb{\vV , \vV'}
                                                    \end{displaymath}
 be a charted orbifold map, where $\vV$ and $\vV'$ are representatives of $\uU$ and $\uU'$, respectively. Assume that there are orbifold charts $(V_\alpha,G_\alpha,\pi_\alpha) \in \vV , \ \alpha =a,b$ and points $x_\alpha \in V_\alpha$ with $\pi_a (x_a) = \pi_b (x_b)$. Then there are arbitrarily small orbifold charts (i.e.\ for each open set $\Omega \subseteq Q$ we may choose charts, which are contained in $\Omega$) $(W,K,\chi) \in \uU$, $(W',K',\chi') \in \uU'$ and embeddings $\lambda_\alpha \colon (W,K,\chi) \rightarrow (V_\alpha , G_\alpha , \pi_\alpha)$, $\mu_\alpha \colon (W',K',\chi') \rightarrow (V_\alpha',G_\alpha' , \pi_\alpha')$ of orbifold charts with $x_\alpha \in \lambda_\alpha (W),\ \alpha= a,b$ such that the induced lift $g$ of $f$ with respect to $f_a$, $\lambda_a$, $\mu_a$ coincides with the one induced by $f_b, \lambda_b, \mu_b$. In other words, we obtain a commutative diagram 
\begin{displaymath}
  \begin{xy}
  \xymatrix{					
								&	V_a \ar[r]^{f_a} 	& V_a' 	&				    \\
      W \ar[rrr]^{g} \ar[ru]_{\lambda_a} \ar[rd]^{\lambda_b}   	&   				  	&		& W' \ar[lu]^{\mu_a} \ar[ld]^{\mu_b}  \\
      								& V_b \ar[r]^{f_b}		&  V_b'	& 
  }
\end{xy}
\end{displaymath}
\end{lem}

\begin{lem}[\hspace{-0.5pt}{\cite[Lemma 5.14]{pohl2010}}]\label{lem: cp:ofdm} \no{lem: cp:ofdm}
 Let $(Q,\uU)$ and $(Q',\uU')$ be orbifolds, $\vV$ a representative of $\uU$, and $\vV'$ one of $\uU'$. Further let $\hat{f} \in \Orb{\vV,\vV'}$. Assume that $\hat{h} \in \Orb{\wW_1, \wW_1'}$ and $\hat{g} \in \Orb{\wW_2 , \wW_2'}$ are both induced by $\hat{f}$.
There are representatives $\wW$ of $\uU$ and $\wW'$ of $\uU'$ together with lifts of the identity $\ve_i \in \Orb{\wW, \wW_i},\ i=1,2$ and $\ve_i' \in \Orb{\wW', \wW_i'} ,\ i=1,2$ such that a charted orbifold map $\hat{k} \in \Orb{\wW, \wW'}$ exists, making the following diagram commutative. 
  
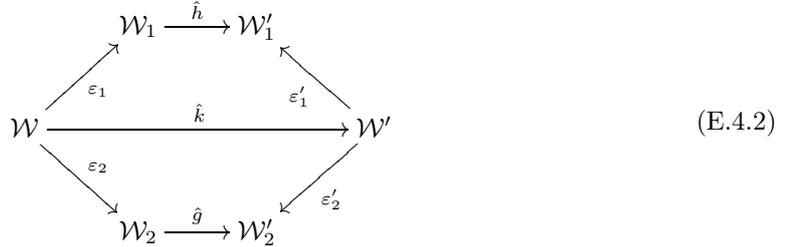
\begin{equation}\begin{gathered}\label{eq: comp:ofm}
  \begin{xy}
  \xymatrix{					
								&	\wW_1 \ar[r]^{\hat{h}} 	& \wW_1' 	&				    \\
      \wW \ar[rrr]^{\hat{k}} \ar[ru]_{\ve_1} \ar[rd]^{\ve_2}   	&   				  	&		& \wW' \ar[lu]^{\ve_1'} \ar[ld]^{\ve_2'}  \\
      								& \wW_2 \ar[r]^{\hat{g}}		&  \wW_2' 	& 
  }
\end{xy}
\end{gathered}
\end{equation}
If the orbifolds are second countable, we may choose $\wW$ and $\wW'$ to be countable.
\end{lem}

\begin{defn} It follows from the last lemma that the relation $\sim$ introduced in Definition \ref{defn: sim:ofd} is indeed an equivalence relation. For details we refer to the exposition in \cite{pohl2010}.\\ Denote the equivalence class of a charted orbifold map $\hat{f}$ with respect to the equivalence relation $\sim$ introduced in Definition \ref{defn: sim:ofd} by $[\hat{f}]$. It will be clear from the context whether $\hat{f}$ is a charted orbifold map and $[\hat{f}]$ denotes its equivalence class, i.e.\ the orbifold morphism, or $\hat{f}$ is a representative of the charted orbifold map and $[\hat{f}]$ is the equivalence class of representatives, which by abuse of notation is also abbreviated as $\hat{f}$.
\end{defn}

\subsection{The orbifold category}\label{sect: ofd_cat}

We have explained how to construct orbifolds and morphisms of orbifolds. Now we introduce the category of orbifolds, which is isomorphic to a full category of certain Lie groupoids (cf. \cite{pohl2010} for details on this topic). 

\begin{defn}\label{defn: orb:cat} \no{defn: orb:cat}
 The category $\ORB$ is defined as follows: The class of objects $\Ob \ORB$ is given by the class of all paracompact Hausdorff orbifolds (as defined in Definition \ref{defn: haef:ofdII}). For two orbifolds $(Q,\uU)$ and $(Q',\uU')$, the morphisms, i.e.\ \ind{orbifold map}{orbifold maps} from $(Q,\uU)$ to $(Q',\uU')$ are the equivalence classes $[\hat{f}]$ of all charted orbifold maps $\hat{f} \in \Orb{\vV,\vV'}$ where $\vV$ is a representative of $\uU$ and $\vV'$ is a representative of $\uU'$, that is, \glsadd{ORB} 
	\begin{displaymath}
	 \ORB ((Q,\uU),(Q',\uU')) \coloneq \setm{\left[\hat{f}\right]}{\hat{f} \in \Orb{\vV , \vV'} , \vV \text{ representative of } \uU, \vV' \text{ representative of } \uU'}.
	\end{displaymath}
 The composition in $\ORB$ is induced by the following construction: Let 
	\begin{displaymath}
	 [\hat{f}] \in \ORB ((Q,\uU),(Q',\uU')) \quad \text{ and }\quad [\hat{g}] \in \ORB ((Q',\uU'),(Q'',\uU'')) 
	\end{displaymath}
 be orbifold maps. Choose representatives $\hat{f} \in \Orb{\vV,\vV'}$ of $[\hat{f}]$ and $\hat{g} \in \Orb{\wW , \wW'}$ of $[\hat{g}]$. Then find representatives $\kK, \kK'$ and $\kK''$ of $\uU, \uU'$ and $\uU''$, respectively, and lifts of the identity $\ve \in \Orb{\kK, \vV}$, $\ve_1' \in \Orb{\kK',\vV'}$, $\ve_2' \in \Orb{\kK', \wW'}$, $\ve'' \in \Orb{\kK'', \wW''}$ together with charted orbifold maps $\hat{h} \in \Orb{\kK, \kK'}$, $\hat{k} \in \Orb{\kK', \kK''}$ such that the diagram 
	
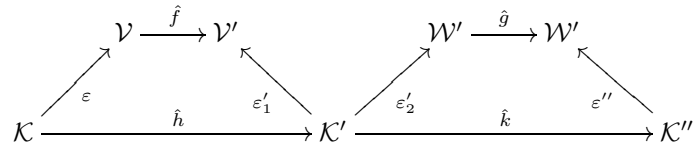
\begin{displaymath}
	 \begin{xy}
  \xymatrix{					
						&	\vV \ar[r]^{\hat{f}} 	& \vV' 	&							& \wW' \ar[r]^{\hat{g}} & \wW' &    \\
      \kK \ar[rrr]^{\hat{h}} \ar[ru]_{\ve}    	&   				&	& \kK' \ar[lu]^{\ve_1'} \ar[ru]_{\ve_2'} \ar[rrr]^{\hat{k}}  & & & \kK'' \ar[lu]^{\ve''}\\
  }
\end{xy}
	\end{displaymath}
 commutes. Define the composition of $[\hat{g}]$ and $[\hat{f}]$ as
  \begin{displaymath}
             [\hat{g}] \circ [\hat{f}] \coloneq [\hat{k} \circ \hat{h}].
  \end{displaymath}
\end{defn}

\begin{prop}[\hspace{-0.5pt}{\cite[Lemma 5.17 and Proposition 5.18]{pohl2010}}] \label{prop: ocat:cp}  \no{prop: ocat:cp}
 It is always possible to compose two orbifold maps in $\ORB ((Q,\uU), (Q',\uU'))$ and $\ORB ((Q',\uU'), (Q'',\uU''))$ and the composition in $\ORB$ is well-defined.
\end{prop}

All equivalence classes of lifts of the identity coincide for a given orbifold $(Q,\uU)$. Hence the \tl identity morphism\tr introduced in Definition \ref{defn: idmor} is the identity morphism of $(Q,\uU)$ in $\ORB$.

\begin{prop}[\hspace{-0.5pt}{\cite[Proposition 5.19]{pohl2010}}]\label{prop: ident:eq}\no{prop: ident:eq}
 Let $(Q,\uU)$ be an orbifold and $\ve$ a lift of $\id_{(Q,\uU)}$. Then the equivalence class $[\ve]$ of $\ve$ consists precisely of all lifts of $\id_{(Q,\uU)}$. Hence the \tl identity morphism\tr $\ido{}$ is the equivalence class $[\ve]$.
\end{prop} \thispagestyle{empty}
\section{Orbifold geodesics: Supplementary Results} \label{App: Haefgeodesics} \no{App: Haefgeodesics} 
In this section, we supply proofs for some of the more technical assertions in Section \ref{Sect: Geod}. 
\begin{Alem}[Lemma \ref{lem: nice:geod}]\label{Alem:ngeod}
 Let $[\hat{c}] \in \ORBM[\iI , (Q,\uU)]$ be an orbifold path and $[a,b] \subseteq \iI$ some compact sub-interval. There exists a charted orbifold map $\hat{g} = (c|_{]x,y[}, \set{g_k}_{1 \leq k \leq N}, [P_g,\nu_g])$ with $x<a<b<y$, $]x,y[ \subseteq \iI$, and $N \in \NN$ such that: 
    \begin{compactitem}
      \item[1.] $[\hat{c}]|_{]x,y[} = [\hat{g}]$,
      \item[2.] $\dom g_k = ]l(k), r(k)[$ for each $1\leq k\leq N$, such that $$x= l(1) < l(2) < r(1) < l(3) < r(2) < \cdots < l(N) < r(N-1) < r(N) = y$$ 
      \item[3.] $P_g = \set{\id_{]l(N),r(N)[}} \cup \setm{\id_{]l(k),r(k)[} , \iota_k^{k+1}, (\iota_k^{k+1})^{-1}}{1\leq k \leq N-1}$, where $\iota_k^{k+1}$ is the canonical inclusion $]l(k),r(k)[ \supseteq ]l(k+1),r(k)[ \hookrightarrow ]l(k+1),r(k+1)[$ .
    \end{compactitem} 
\end{Alem}

\begin{proof}[Proof of Lemma \ref{lem: nice:geod}]\label{proof:LA}
 Consider a representative $\hat{c} =(c,\set{c_i}_{i \in I}, [P_c,\nu_c])$ of $[\hat{c}]$ whose domain atlas is contained in $\aA_\iI$. As $[a,b] \subseteq \iI$ is compact, there is a finite subset $F \subseteq I$ such that $[a,b] \subseteq \bigcup_{i \in F} \dom c_{i}$ and $\dom c_i \cap [a,b] \neq \emptyset$ for all $i \in F$ hold. Set $x \coloneq \inf \bigcup_{i \in F} \dom c_i$ and $y \coloneq \sup \bigcup_{i \in F} \dom c_i$ and consider $\hat{c}|_{]x,y[}$. By construction, for $i \in F$ the set $\dom c_i$ is contained in $]x,y[$. Consider the representative $\hat{\iota}_{]x,y[}$ of the orbifold map $[\hat{\iota}_{]x,y[}]$ whose lifts are given by the family $\set{\id_{\dom c_i}}_{i \in F}$. Following Construction \ref{con: cp:chom} the composition $\hat{h} \coloneq \hat{c}\circ \hat{\iota}_{]x,y[}$ is a representative of $[\hat{c}]|_{]x,y[} \coloneq [\hat{c}] \circ [\hat{\iota}_{]x,y[}]$. By construction, the family of lifts of $\hat{h}$ is $\set{c_i}_{i \in F}$ . As $F$ is finite, we can choose and fix a partition of $]x,y[$ by real numbers $l(k)', r(k)', 1\leq k \leq N \in \NN$ which are ordered as in 2., such that $]l(k)',r(k)'[ \subseteq \dom c_{i_k}$ holds for some $i_k \in F$. Note that each inclusion $\iota_k \colon ]l(k)', r(k)'[ \hookrightarrow \dom c_{i_k}$ is a change of orbifold charts. Apply Lemma \ref{lemdef: ind:ofdm} with respect to the family of pairs $(\iota_k , \id_{\cod c_{i_k}})$, $k \in \set{1, \ldots , N}$ to obtain a representative $\hat{g}'= (c_{]x,y[}, \set{g'_k}_{1\leq k \leq N}, [P_{g'}, \nu_{g'}])$ induced by $\hat{h}$.\\
 Choose $\iota_{k}^{k+1} \in P_{g'}$ with $\dom \iota_k^{k+1} \subseteq ]l(k)',r(k)'[$ and $\cod \iota_k^{k+1} \subseteq ]l(k+1)',r(k+1)'[$. Set $\iota_{k+1}^{k} \coloneq (\iota_k^{k+1})^{-1}$, $l(1) \coloneq x , r(N) \coloneq y$ and 
 \begin{displaymath}
  r(k) \coloneq \sup \dom \iota_k^{k+1}, l(k+1) \coloneq \inf \dom \iota_k^{k+1} \text{ for each } 1 \leq k \leq N-1.
 \end{displaymath}
 By construction $]l(k),r(k)[ \subseteq ]l(k)',r(k)'[$ holds for $1\leq k \leq N$. The numbers $l(k),r(k)$ are ordered as in 2., since the $l(k)',r(k)'$ were ordered in this way. Furthermore, $]x,y[ = \bigcup_{1\leq k\leq N} ]l(k),r(k)[$ is satisfied. With this choice of $\iota_k^{k+1}$, the quasi-pseudogroup $P_g$ as defined in 3. generates the changes of charts for $\setm{]l(k),r(k)[}{1\leq k\leq N}$. Define 
\begin{displaymath}
        \nu_g (\lambda) \coloneq \begin{cases}
                \id_{\cod c_{i_k}} 	  	& \text{ if }  \lambda = \id_{]l(k),r(k)[}) \\
		\nu_{g'} (\iota_{k}^{k+1}) 	& \text{ if } \lambda =\iota_{k}^{k+1}   \\
		\nu_{g'} (\iota_{k}^{k+1})^{-1} & \text{ if } \lambda = (\iota_{k}^{k+1})^{-1} 
	      \end{cases}                                                                                                                                                                                                                                                                                                                                                                                                                                                                                                                                                                                                                                                                                                                                                                                  \end{displaymath}
 to obtain a map $\nu_g \colon P_g \rightarrow \PSI (\aA)$, where $\aA \in \uU$ contains the range family of $\hat{g}'$.\\
 Apply Lemma \ref{lemdef: ind:ofdm} with respect to the pairs $(]l(k),r(k)[ \hookrightarrow ]l(k)',r(k)'[ , \id_{\cod c_{i_k}})$ for $1\leq k\leq N$ to obtain a representative $\hat{g} = (c|_{]x,y[}, \set{g_k}_{1\leq k \leq N}, [P,\nu])$ induced by $\hat{g}'$. Reviewing the construction of $\nu_g$, we see that by construction and property (R4) (d) of Definition \ref{defn: rep:ofdm} the germs of $\nu (\lambda)$ and $\nu_h (\mu)$ must coincide at $g_k(x)$ if $\germ_x \lambda = \germ_x \mu$ holds for $x \in \dom g_k$. Thus $(P_g, \nu_g) \sim (P,\nu)$ follows, whence we may replace the pair $(P,\nu)$ with $(P_g,\nu_g)$. Observe that in each step, we have only applied Lemma \ref{lemdef: ind:ofdm}. Thus $[\hat{g}] = [\hat{c}]|_{]x,y[}$ holds.
\end{proof}

Clearly the definition of the restriction of an orbifold map yields the following Lemma:
 
\begin{Acor}\label{cor: og:ref1} \no{cor: og:ref1}
 If $[\hat{c}] \in \ORBM[\iI,(Q,\uU)]$ is an orbifold geodesic and $[a,b] \subseteq \iI$ compact, then the restriction $[\hat{g}]=[\hat{c}]|_{]x,y[}$ with $x<a<b<y$ constructed in Lemma \ref{lem: nice:geod} is an orbifold geodesic. 
\end{Acor}

\begin{proof}
Simply choose in Lemma \ref{lem: nice:geod} an atlas contained in $\aA_\iI$. 
\end{proof}

\begin{Alem}\label{lem: og:equiv} \no{lem: og:equiv}
 Consider representatives $\hat{c} = (c,\set{c_k}_{k \in A}, [P , \nu]) ,\ \hat{c}'=(c',\set{c_r'}_{r \in B}, [P' , \nu'])$ of orbifold geodesics in $\ORBM[\iI , (Q,\uU)]$, whose domain atlases are contained in $\aA_\iI$. Assume that the lifts satisfy $\cod c_k= U_k$ for $(U_k,G_k,\psi_k) \in \uU$, respectively $\cod c_r' = W_r$ for $(W_r,H_r,\varphi_r) \in \uU$. The following conditions are equivalent:
  \begin{compactenum}
   \item $[\hat{c}] = [\hat{c}']$,
   \item For all $k \in A, r \in B$ and $t \in \dom c_k \cap \dom c_r'$, there is a change of charts $\lambda^{k,r}_t \colon U_k \supseteq \dom \lambda_t^{k,r} \rightarrow W_r$ with $T_t (\lambda^{k,r}_t c_k) (1) = T_tc_r (1)$ (i.e.\ the initial vectors coincide),
   \item for any $t \in \iI$, there is a pair $(k,r) \in A \times B$ and a change of charts $\lambda_t \colon U_k \supseteq \dom \lambda_t \rightarrow W_r$ such that $t \in \dom c_k \cap \dom c_r'$ and $T_t (\lambda_t c_k) (1) = T_tc_r (1)$,
   \item there are representatives $\hat{g} = (c, \set{c_k}_{k \in I}, [P_g , \nu_g])$ of $[\hat{c}]$ and $\hat{g}' = (c, \set{c_k}_{k \in I}, [P_g , \nu_{g}'])$, respectively, of $[\hat{c}']$ whose domain atlases are contained in $\aA_\iI$.\\
   In particular, a geodesic arc in $Q$ is uniquely determined by the initial vector.  
  \end{compactenum}
\end{Alem}

\begin{proof}
 \tl (a) $\Rightarrow $ (b)\tr is a  reformulation of Lemma \ref{lem: osmpath:prop} for orbifold geodesics. \tl (b) $\Rightarrow$ (c)\tr is trivial. 
 To check \tl (c) $\Rightarrow$ (d)\tr , we construct representatives induced by $\hat{c}$ and $\hat{c}'$: The chart domains of the domain atlases of $\hat{c}$ and $\hat{c}'$ are intervals $I_k \coloneq \dom c_k$, $k \in A$, respectively $J_r \coloneq \dom c_r'$, $r \in B$. 
 Pick some $t_0 \in \iI$ together with a pair $(k,r) \in A \times B$ satisfying the hypothesis of (c). There is $\lambda_{t_0} \in \Ch_{U_k,W_r}$ with $T_{t_0} (\lambda_{t_0} c_k) (1) = T_{t_0} c_r' (1)$. Shrinking $\dom \lambda_{t_0}$, we may assume that the set $t_0 \in \dom \lambda_{t_0}$ is $G_k$-stable. Thus it induces an orbifold chart $(\dom \lambda_{t_0}, G_{k,\dom \lambda_{t_0}}, \psi_k|_{\dom \lambda_{t_0}}) \in \uU$. As $c_k$ is a geodesic, we may choose $\ve_{t_0} > 0$ with $c_k ([t_0 - \ve_{t_0} , t_0+\ve_{t_0}]) \subseteq \dom \lambda_{t_0}$ and $[t_0 - \ve_0 , t_0 + \ve_0] \subseteq J_r$. The change of charts $\lambda_{t_0}$ is a Riemannian isometry, since $(Q,\uU,\rho)$ is a Riemannian orbifold. In particular, $\lambda_{t_0}$ maps geodesics of $\dom \lambda_{t_0} \subseteq U_k$ to geodesics of $W_r$. Thus $\lambda_{t_0} \circ c_k \colon ]t_0-\ve_{t_0},t_0+\ve_{t_0}[ \rightarrow W_r$ is a geodesic. Uniqueness of geodesics in Riemannian manifolds implies that $\lambda_{t_0} \circ c_k|_{]t_0 - \ve_{t_0}, t_0 +\ve_{t_0}[} = c_r'|_{]t_0 - \ve_{t_0}, t_0 +\ve_{t_0}[}$, as their derivatives coincide in $t_0$. For the trivial orbifold $\iI$ the set $C_{t_0} \coloneq ]t_0 - \ve_{t_0} , t_0 + \ve_{t_0}[ \subseteq I_k \cap J_r$ induces an orbifold chart via the inclusion of sets. Set $\alpha (t_0) \coloneq k$ and $\beta (t_0) \coloneq r$ and define changes of orbifold charts $\mu_{t_0, \alpha} \colon C_{t_0} \rightarrow I_{\alpha (t_0)}$, $\mu_{t_0, \beta} \colon C_{t_0} \rightarrow J_{\beta (t_0)}$ and $\nu_{t_0,\alpha} \colon \dom \lambda_{t_0} \rightarrow U_{\alpha (t_0)}$ via the inclusion of sets. Furthermore, set $\nu_{t_0, \beta} \coloneq \lambda_{t_0}$. Reviewing the construction, we see that $c_{\alpha (t_0)} \mu_{t_0,\alpha} \subseteq \im \nu_{t_0 , \alpha}$ and $c_{\beta (t_0)}' \mu_{t_0,\beta} \subseteq \im \lambda_{t_0} = \im \nu_{t_0,\beta}$. This implies  
    \begin{equation}\label{eq: ident:geod}
     \nu_{t_0,\alpha}^{-1} c_{\alpha (t_0)} \mu_{t_0,\alpha} = \nu_{t_0,\beta}^{-1} c_{\beta (t_0)}' \mu_{t_0,\beta}.
    \end{equation}
 With respect to the pair $(C_{t_0}, \set{\id_{C_{t_0}}}, C_{t_0}\hookrightarrow \iI)$ and $(\dom \lambda_{t_0}, G_{k, \dom \lambda_{t_0}},\psi_k|_{\dom_{t_0}})$ the lifts of $\hat{c}$ and $\hat{c}'$ coincide.
 The construction did not depend on $t_0$ and may be repeated for each $t \in \iI$. In this way we obtain a (possibly infinite) subset $R \subseteq \iI$ such that $\bigcup_{t \in R} C_{t} = \iI$ and $C_t \neq C_s$ if $t\neq s$. Since these sets cover $\iI$, the construction yields an orbifold atlas $\cC \subseteq \aA_\iI$ for $\iI$. It may happen that the charts $(\dom \lambda_t ,G_{\alpha (t),\dom \lambda_t}, \psi_{\alpha (t)}|_{\dom \lambda_t})$ and $(\dom \lambda_s ,G_{\alpha (s),\dom \lambda_s}, \psi_{\alpha (s)}|_{\dom \lambda_s})$ coincide for $s\neq t$. To satisfy the requirement (R2) in Definition \ref{defn: rep:ofdm}, we redefine the charts: Take $\dom \lambda_s \times \set{s}$ instead of $\dom \lambda_s$ and redefine the group actions, changes of charts etc.\ in the obvious way. Recall that this does not change the equivalence class of $\hat{c}$ and $\hat{c}'$ by virtue of Lemma \ref{lemdef: ind:ofdm}. Without loss of generality we may thus assume $(\dom \lambda_t ,G_{\alpha (t),\dom \lambda_t}, \psi_{\alpha (t)}|_{\dom \lambda_t}) \neq (\dom \lambda_s ,G_{\alpha (s),\dom \lambda_s}, \psi_{\alpha (s)}|_{\dom \lambda_s})$ for $s\neq t$. Using Lemma \ref{lemdef: ind:ofdm}, the charted maps $\hat{c}$ and $\hat{c}'$, induce representatives $\hat{h}$ and $\hat{h}'$ with respect to $\cC$ and an atlas $\wW \in \uU$ which contains $\setm{(\dom \lambda_t ,G_{\alpha (t),\dom \lambda_t}, \psi_{\alpha (t)}|_{\dom \lambda_t}))}{t \in R}$. From \eqref{eq: ident:geod} we deduce that the lifts of $\hat{h}$ and $\hat{h}'$ coincide. Choose a refinement of the domain atlas of $\hat{h}$ as follows: There is a sequence of real numbers in $\iI$
  \begin{displaymath}
   \cdots  < l(-1) < r(-2) < l(0) < r(-1) < l(1) < r(0) < l(2) < r(1) < \cdots 
  \end{displaymath}
  such that $]l(n),r(n)[$ is contained in some chart of the domain atlas of $\hat{h}$ for each $n \in \ZZ$. Apply an argument as in the proof of Lemma \ref{lem: nice:geod} (cf.\ Lemma \ref{Alem:ngeod}) to obtain a cover of $\iI$ by intervals $I_k$ indexed by $\ZZ$, such that the following is satisfied: 
  \begin{compactenum}
   \item[1.] $I_k \cap I_j \neq \emptyset$ if and only if $j \in \set{k-1,k,k+1}$, $k,j \in \ZZ$,
   \item[2.] $\hat{h}$ induces a representative $\hat{g} = (c, \set{g_k}_{k \in \ZZ}, [P_g, \nu_g])$ of $[\hat{c}]$ and $\hat{h}'$ induces a representative $\hat{g}' = (c', \set{g_k'}_{k \in \ZZ}, [P_g', \nu_g'])$ of $[\hat{c}']$ such that $P_g = P_{g'}$ and $P_g = \setm{\id_{]l(k),r(k)[} , \iota_k^{k+1}, (\iota_k^{k+1})^{-1}}{k \in \ZZ}$, where $\iota_k^{k+1}, (\iota_k^{k+1})^{-1}$ are defined as in Lemma \ref{lem: nice:geod}.
   \item[3.] As the lifts of $\hat{h}$ and $\hat{h}'$ coincide, for each $k \in \ZZ$ the lifts $g_k,g_k'$ are given as restriction $g_k =g_k' = h_s|_{]l(k),r(k)[} \colon ]l(k),r(k)[ \rightarrow V_k$, $(V_k,G_k,\psi_k) \in \uU$ of a lift $h_s$ of $\hat{h}$.
  \end{compactenum}
   Shrinking the sets $]l(n),r(n)[, \ n \in \ZZ$, we assume that $g_k (]l(k+1),r(k)[)$ and $g_k(]l(k),r(k-1)[)$ are contained in stable subsets of $\dom \nu_{\hat{g}} (\iota_k^{k+1}) \cap \dom  \nu_{\hat{g}'} (\iota_k^{k+1})$ and $\dom \nu_{\hat{g}} ((\iota_{k-1}^{k})^{-1}) \cap \dom  \nu_{\hat{g}'} ((\iota_{k-1}^{k})^{-1})$, respectively, for each $k \in \ZZ$. Restricting the changes of charts to these stable subsets, by Definition \ref{defn: eq:rofdm} the pairs $(P_g, \nu_{\hat{g}})$ and $(P_g, \nu_{\hat{g}'})$ may be replaced by equivalent pairs such that the maps $\nu_{\hat{g}} (\lambda), \nu_{\hat{g}'} (\lambda )$ are embeddings of orbifold charts with $\dom \nu_{\hat{g}} (\lambda) = \dom \nu_{\hat{g}'} (\lambda )$ for each $\lambda \in P_g$. Unfortunately, $\nu_{\hat{g}}$ and $\nu_{\hat{g}'}$ need not coincide. However, since the lifts coincide we obtain 
    \begin{displaymath}
     \nu_{\hat{g}} (\iota_k^{k+1}) \circ g_k|_{]l(k+1),r(k)[} = g_{k+1} \circ \iota_k^{k+1} = \nu_{\hat{g}'} (\iota_k^{k+1}) \circ g_k|_{]l(k+1),r(k)[}.
    \end{displaymath}
 Hence both geodesic arcs coincide. As $\nu_{\hat{g}} (\iota_k^{k+1})$ and $\nu_{\hat{g}'} (\iota_k^{k+1})$ are embeddings of orbifold charts with the same domain, for each $k \in \ZZ$ there is some $\gamma_{k+1} \in G_{k+1}$ with $\nu_{\hat{g}} (\iota_k^{k+1}) = \gamma_{k+1} . \nu_{\hat{g}'} (\iota_k^{k+1})$. \\
 \tl (d) $\Rightarrow$ (a)\tl \ Consider representatives $\hat{g}$ of $[\hat{c}]$ and $\hat{g}'$ of $[\hat{c}']$ as constructed in Step \tl (c) $\Rightarrow$ (d)\tr . We claim that $[\hat{g}] = [\hat{g}']$ holds. To prove the claim, consider the case that the geodesic arc $\im c$ contains non-singular points. Hence there are $k \in \ZZ$ and $z \in \iI$ such that $c(z) = \psi_k c_k (z)$ is non-singular. For each subset $H_k \subseteq G$ the components of $\bigcap_{g \in H_k} \Sigma_{g}$ are totally geodesic submanifolds of $(V_k,\rho_k)$ by \cite[II.\ Theorem 5.1]{kobayashi1972}). Assume that there is a an open, non-empty set $U$ such that $\im c_k \cap U$ is contained in a component jointly fixed by the elements of some subset $H_k \subseteq G_k$, which contains elements different from the identity $\id_{V_k} \in G_k$. Then the image $\im c_k$ is contained in this component (cf.\ \cite[Proof of Theorem 1.10.15]{klingenberg1995}). This contradicts the choice of $c_k (z)$, whence the non-singular points must be a dense subset of $\im c_k$ with respect to the subspace topology. Changes of charts preserve non-singular points. Hence the same argument may be repeated to prove that the non-singular points must be dense in the image of each $c_k, k\in \ZZ$. In conclusion, we have to consider two cases:
 \paragraph{Case 1:} The geodesic arc of $[\hat{c}]$ (or equivalently the arc of $[\hat{c}']$) contains a non-singular point. The preparatory considerations show that the non-singular points are dense in the image of each lift. Hence $\gamma_{k+1} .\nu_{\hat{g}} (\iota_k^{k+1}) = \nu_{\hat{g}'} (\iota_k^{k+1})$ implies $\gamma_{k+1} = \id_{V_{k+1}}, \forall k \in \ZZ$ as $\im c_{k+1}$ contains non-singular points. We deduce $\nu_{\hat{g}} = \nu_{\hat{g}'}$, whence $\hat{g} = \hat{g}'$ follows. 
 \paragraph{Case 2:} The geodesic arc of $[\hat{c}]$ (or equivalently the arc of $[\hat{c}']$) is contained in the singular locus of $Q$. We construct a representative of $[\hat{c}]$ which coincides with $\hat{g}'$. Apply Lemma \ref{lemdef: ind:ofdm} with suitable changes of charts to $\hat{g}$ and $\hat{g}'$, such that $(V_k,G_k,\psi_k) \neq (V_j,G_j,\psi_j)$ holds if $k \neq j$. Observe that for each choice $(\eta_k)_{k \in \ZZ} \in \prod_{k \in \ZZ} G_k$ the pairs $\set{(\id_{]l(k),r(k)[},\eta_k)}_{k\in \ZZ}$ induce another representative $\hat{h}$ of $[\hat{c}]$ by Lemma \ref{lemdef: ind:ofdm}. Recall from the construction\footnote{Unfortunately, these details are not apparent from the mere statement of Lemma \ref{lemdef: ind:ofdm}. However, the proof of this Lemma in \cite[p.\ 21]{pohl2010} readily entails these facts: Notice that we may choose $P_h = P_{\hat{g}}$, since we applied Lemma \ref{lemdef: ind:ofdm} to the pairs $\set{(\id_{]l(k),r(k)[},\eta_k)}_{k\in \ZZ}$. Here the first embedding of each pair is an identity, whence we need not restrict the elements of $P_{\hat{g}}$ as in \cite[p.\ 21]{pohl2010}. Moreover, the identity \eqref{eq: ident:flip} then follows directly from the proof.} of $\hat{h} = (c, \set{\eta_k \circ c_k}_{k \in \ZZ}, P_h, \nu_h)$ the following details: As $\eta_k \in G_k$ is defined on $V_k$, we may choose $P_h = P_{\hat{g}}$ and $\nu_h$ is uniquely determined by the identity
  \begin{equation}\label{eq: ident:flip}
   \nu_h (\iota_k^{k+1}) = \eta_{k+1}^{-1}\nu_{\hat{g}} (\iota_k^{k+1}) \eta_k^{-1}|_{\eta_k (\dom \nu_{\hat{g}} (\iota_k^{k+1}))}.
  \end{equation}
 We claim that it is possible to inductively (starting from $0$ and consider the cases $\NN_0$ and $\ZZ_0^{-}$ independently) choose the family $(\eta_k)_{k \in \ZZ}$, such that $\eta_k c_k = c_k$ and $\nu_h = \nu_{\hat{g}'}$. Begin with $k=0$. Since $\dom \nu_{\hat{g}} (\iota_{-1}^0) = \dom \nu_{\hat{g}'} (\iota_{-1}^0)$ holds (and these maps are embeddings of orbifold charts by Step \tl (c) $\Rightarrow$ (d)\tr), by Proposition \ref{prop: ch:prop} (d) there is $\gamma_0 \in G_0$ with $\nu_{\hat{g}} (\iota_{-1}^0) = \gamma_0. \nu_{\hat{g}'} (\iota_{-1}^0)$.  The situation is visualized in Figure \ref{fig: lifts}, where we depict the lifts together with the embeddings of orbifold charts. 
 \begin{figure}[t]\centering 
\begin{tikzpicture}[y=0.80pt, x=0.8pt,yscale=-1, inner sep=0pt, outer sep=0pt, >=stealth]
\begin{scope}[shift={(0,-768.89762)}]
  \path[draw=black,line join=miter,line cap=butt,line width=1.226pt]
    (21.7866,968.3939) .. controls (12.8278,972.5175) and (6.0134,981.0205) ..
    (3.9405,990.6624) .. controls (1.8675,1000.3042) and (4.5851,1010.8566) ..
    (11.0576,1018.2977) .. controls (16.0387,1024.0242) and (22.9903,1027.8313) ..
    (30.3178,1029.8091) .. controls (37.6452,1031.7870) and (45.3528,1032.0218) ..
    (52.9169,1031.3976) .. controls (68.0449,1030.1493) and (82.7940,1025.5343) ..
    (97.9623,1024.9515) .. controls (108.8933,1024.5315) and (119.8030,1026.2174)
    .. (130.7420,1026.2160) .. controls (136.2116,1026.2154) and
    (141.7072,1025.7877) .. (147.0184,1024.4814) .. controls (152.3296,1023.1751)
    and (157.4643,1020.9669) .. (161.7997,1017.6323) .. controls
    (167.0348,1013.6056) and (170.9921,1007.9823) .. (173.1902,1001.7542) ..
    controls (175.3883,995.5261) and (175.8426,988.7172) .. (174.6744,982.2167) ..
    controls (172.3379,969.2157) and (163.5067,957.8292) .. (152.1436,951.0940) ..
    controls (139.2681,943.4623) and (123.6340,941.5525) .. (108.7379,943.0117) ..
    controls (93.8419,944.4709) and (79.4567,949.0781) .. (65.2389,953.7555) ..
    controls (50.7204,958.5319) and (36.2358,963.4116) .. (21.7866,968.3939);
  \path[draw=black,line join=miter,line cap=butt,line width=1.320pt]
    (111.5406,794.9883) .. controls (94.8983,801.0506) and (81.5890,815.5654) ..
    (76.9888,832.6697) .. controls (72.3886,849.7740) and (76.6215,869.0068) ..
    (87.9778,882.5992) .. controls (95.4032,891.4868) and (105.5454,897.8921) ..
    (116.4548,901.7794) .. controls (127.3641,905.6667) and (139.0329,907.1222) ..
    (150.6141,907.1060) .. controls (173.7766,907.0735) and (196.5022,901.2892) ..
    (219.4913,898.4598) .. controls (253.1281,894.3199) and (287.5599,896.5703) ..
    (320.3183,905.2572) .. controls (331.2599,908.1587) and (342.1269,911.7901) ..
    (353.4205,912.5591) .. controls (359.0674,912.9436) and (364.7990,912.5952) ..
    (370.2565,911.0950) .. controls (375.7140,909.5948) and (380.8958,906.9108) ..
    (384.9791,902.9914) .. controls (390.5889,897.6067) and (393.8917,890.0767) ..
    (394.8903,882.3652) .. controls (395.8889,874.6537) and (394.6944,866.7680) ..
    (392.2674,859.3805) .. controls (387.4135,844.6055) and (377.8722,831.9201) ..
    (370.1838,818.4016) .. controls (366.6898,812.2580) and (363.5398,805.8725) ..
    (359.2244,800.2752) .. controls (350.5164,788.9807) and (337.4067,781.5411) ..
    (323.5591,778.1298) .. controls (309.7114,774.7184) and (295.1405,775.1151) ..
    (281.0788,777.4949) .. controls (252.9554,782.2544) and (226.6443,794.7147) ..
    (198.6683,800.2752) .. controls (169.8144,806.0102) and (139.4893,804.1701) ..
    (111.5406,794.9883);
  \path[draw=black,line join=miter,line cap=butt,line width=1.186pt]
    (338.7113,940.3390) .. controls (332.5890,938.7551) and (326.0347,938.8777) ..
    (319.9760,940.6896) .. controls (313.9173,942.5015) and (308.3717,945.9974) ..
    (304.1242,950.6824) .. controls (298.0801,957.3491) and (294.7891,966.3029) ..
    (294.4243,975.2943) .. controls (294.0594,984.2856) and (296.5189,993.2835) ..
    (300.8302,1001.1823) .. controls (306.7252,1011.9825) and (316.0233,1020.7619)
    .. (326.7601,1026.7716) .. controls (337.4969,1032.7813) and
    (349.6299,1036.0921) .. (361.8823,1037.2213) .. controls (386.3870,1039.4798)
    and (410.8573,1033.1977) .. (434.2374,1025.5196) .. controls
    (448.0295,1020.9903) and (461.9955,1015.7598) .. (472.9493,1006.2334) ..
    controls (478.4262,1001.4702) and (483.0625,995.6527) .. (485.9677,989.0011)
    .. controls (488.8729,982.3495) and (489.9879,974.8407) .. (488.5884,967.7185)
    .. controls (487.5390,962.3779) and (485.1007,957.3424) .. (481.7332,953.0666)
    .. controls (478.3656,948.7907) and (474.0837,945.2663) .. (469.3465,942.5862)
    .. controls (459.8722,937.2261) and (448.7592,935.3000) .. (437.8742,935.4079)
    .. controls (416.1043,935.6237) and (395.0442,943.5529) .. (373.2983,944.5980)
    .. controls (361.6392,945.1584) and (349.8862,943.7111) ..
    (338.7113,940.3390);
  \path[draw=black,line join=miter,line cap=butt,line width=1.107pt]
    (107.3790,849.7895) .. controls (115.6300,851.2499) and (123.6403,852.3907) ..
    (131.2827,852.7952) .. controls (146.5675,853.6043) and (160.3804,851.4683) ..
    (171.7024,843.0539) .. controls (216.9906,809.3961) and (227.2559,798.0042) ..
    (275.5633,807.8427) .. controls (287.6402,810.3023) and (295.6788,813.7328) ..
    (302.2834,817.7296) .. controls (308.8879,821.7265) and (314.0583,826.2897) ..
    (320.3986,831.0148) .. controls (326.7390,835.7398) and (334.2493,840.6266) ..
    (345.5336,845.2708) .. controls (351.1757,847.5928) and (357.7614,849.8542) ..
    (365.6161,852.0043);
  \path[draw=black,line join=miter,line cap=butt,line width=1.163pt]
    (18.1445,982.5877) .. controls (25.4127,989.1217) and (33.0663,993.6551) ..
    (40.9815,996.6606) .. controls (48.8966,999.6660) and (57.0734,1001.1435) ..
    (65.3877,1001.5657) .. controls (82.0164,1002.4100) and (99.1958,999.0329) ..
    (115.9346,995.2152) .. controls (124.3041,993.3064) and (132.5634,991.2875) ..
    (140.5887,989.6310);
  \path[draw=black,line join=miter,line cap=butt,line width=0.800pt]
    (331.1655,970.6328) .. controls (335.4071,974.3190) and (340.4313,978.1313) ..
    (347.9803,981.7544) .. controls (355.5293,985.3774) and (365.6031,988.8110) ..
    (379.9437,991.7398) .. controls (387.1140,993.2041) and (394.0696,993.6964) ..
    (400.7523,993.4369) .. controls (407.4350,993.1773) and (413.8447,992.1658) ..
    (419.9231,990.6226) .. controls (426.0014,989.0793) and (431.7484,987.0043) ..
    (437.1056,984.6176);
  \path[draw=black,dash pattern=on 2.23pt off 1.12pt,line join=miter,line
    cap=butt,miter limit=4.00,line width=1.115pt] (133.7753,971.9051) --
    (129.1792,974.1132) .. controls (126.8570,973.5504) and (124.3155,973.9575) ..
    (122.2852,975.2173) .. controls (120.9236,976.0622) and (119.8097,977.2551) ..
    (118.8382,978.5295) .. controls (116.8885,981.0874) and (115.4275,984.0457) ..
    (113.0932,986.2581) .. controls (112.3212,986.9898) and (111.4583,987.6347) ..
    (110.7952,988.4663) .. controls (109.9281,989.5537) and (109.4431,990.9398) ..
    (109.4431,992.3306) .. controls (109.4431,993.7213) and (109.9281,995.1075) ..
    (110.7952,996.1948) -- (115.3912,1002.8193) .. controls (116.4992,1003.6170)
    and (117.6504,1004.3545) .. (118.8382,1005.0275) .. controls
    (120.3181,1005.8660) and (121.8548,1006.6043) .. (123.4342,1007.2357) --
    (126.8812,1009.4437) .. controls (131.0932,1010.5518) and (135.3062,1011.6560)
    .. (139.5203,1012.7561) .. controls (142.0449,1013.4152) and
    (144.6019,1014.0772) .. (147.2111,1014.0931) .. controls (149.8203,1014.1091)
    and (152.5296,1013.4103) .. (154.4573,1011.6519) .. controls
    (155.2486,1010.9301) and (155.8831,1010.0566) .. (156.5814,1009.2444) ..
    controls (157.2797,1008.4323) and (158.0702,1007.6609) .. (159.0533,1007.2357)
    .. controls (159.6088,1006.9955) and (160.2067,1006.8741) ..
    (160.7959,1006.7360) .. controls (161.3851,1006.5979) and (161.9786,1006.4383)
    .. (162.5003,1006.1316) .. controls (163.0282,1005.8213) and
    (163.4651,1005.3681) .. (163.7933,1004.8512) .. controls (164.1216,1004.3343)
    and (164.3439,1003.7549) .. (164.4891,1003.1600) .. controls
    (164.7793,1001.9704) and (164.7643,1000.7312) .. (164.7983,999.5071) ..
    controls (164.8724,996.8396) and (165.1867,994.1726) .. (165.0269,991.5088) ..
    controls (164.8672,988.8451) and (164.1849,986.1196) .. (162.5003,984.0500) ..
    controls (160.8378,982.0075) and (158.2397,980.7592) .. (155.6063,980.7377) --
    (149.8613,977.4255) -- (133.7752,971.9051);
  \path[draw=black,dash pattern=on 2.89pt off 1.45pt,line join=miter,line
    cap=butt,miter limit=4.00,line width=1.446pt] (153.2752,835.7435) .. controls
    (144.4008,830.5768) and (134.5751,827.0500) .. (124.4407,825.3937) .. controls
    (120.1415,824.6911) and (115.6893,824.3291) .. (111.4651,825.3937) .. controls
    (107.4009,826.4180) and (103.7003,828.7636) .. (100.9399,831.9175) .. controls
    (98.1795,835.0715) and (96.3550,839.0124) .. (95.6061,843.1363) .. controls
    (95.0749,846.0615) and (95.0749,849.0823) .. (95.6061,852.0075) --
    (106.8104,865.5679) -- (112.9068,871.4135) -- (124.4407,881.5784) --
    (134.5327,888.9712) -- (147.5083,893.4068) .. controls (149.9096,893.2999) and
    (152.3156,893.2999) .. (154.7169,893.4068) .. controls (159.3348,893.6124) and
    (163.9330,894.2127) .. (168.5550,894.2731) .. controls (173.1771,894.3334) and
    (177.8953,893.8273) .. (182.1097,891.9283) .. controls (184.9828,890.6336) and
    (187.6080,888.6609) .. (189.3183,886.0141) .. controls (190.9090,883.5525) and
    (191.6512,880.5537) .. (191.3924,877.6343) .. controls (191.1336,874.7150) and
    (189.8755,871.8935) .. (187.8766,869.7501) -- (183.5514,863.8359) .. controls
    (179.4883,854.9164) and (172.9261,847.1512) .. (164.8090,841.6577) .. controls
    (161.2221,839.2301) and (157.3398,837.2394) .. (153.2752,835.7435);
  \path[draw=black,dash pattern=on 2.24pt off 1.12pt,line join=miter,line
    cap=butt,miter limit=4.00,line width=1.122pt] (432.2055,950.1728) .. controls
    (424.6979,951.6034) and (417.5999,955.2775) .. (412.2545,960.7400) .. controls
    (406.9092,966.2025) and (403.3703,973.4508) .. (402.5700,981.0515) .. controls
    (401.6166,990.1056) and (404.6524,999.5082) .. (410.7205,1006.2952) ..
    controls (416.7886,1013.0823) and (425.7941,1017.1477) .. (434.8980,1017.2098)
    .. controls (444.0020,1017.2718) and (453.0621,1013.3298) ..
    (459.2222,1006.6262) .. controls (465.3824,999.9226) and (468.5463,990.5623)
    .. (467.7165,981.4961) -- (464.4180,958.1415) .. controls (455.4975,951.0258)
    and (443.4147,948.0368) .. (432.2055,950.1728);
  \path[draw=black,dash pattern=on 1.88pt off 0.94pt,line join=miter,line
    cap=butt,miter limit=4.00,line width=0.940pt] (355.3858,831.6022) .. controls
    (352.8090,829.1484) and (349.4941,827.5094) .. (346.0060,826.8063) .. controls
    (342.5178,826.1033) and (338.8658,826.3217) .. (335.4426,827.2926) .. controls
    (328.5961,829.2343) and (322.8341,834.1103) .. (318.7896,839.9659) .. controls
    (317.2215,842.2363) and (315.8533,844.7495) .. (315.5606,847.4932) .. controls
    (315.3612,849.3618) and (315.6776,851.2731) .. (316.3956,853.0097) .. controls
    (317.1137,854.7462) and (318.2267,856.3093) .. (319.5785,857.6147) .. controls
    (322.2819,860.2255) and (325.8769,861.7668) .. (329.5532,862.5479) .. controls
    (332.7602,863.2292) and (336.0733,863.3729) .. (339.2404,864.2206) .. controls
    (344.6803,865.6766) and (349.3398,869.0994) .. (354.3094,871.7480) .. controls
    (358.1146,873.7759) and (362.2023,875.3816) .. (366.4835,875.8952) .. controls
    (370.7647,876.4088) and (375.2655,875.7684) .. (378.9729,873.5667) .. controls
    (380.8266,872.4659) and (382.4624,870.9847) .. (383.6804,869.2058) .. controls
    (384.8984,867.4269) and (385.6917,865.3486) .. (385.8954,863.2023) .. controls
    (386.0990,861.0560) and (385.7038,858.8461) .. (384.7036,856.9363) .. controls
    (383.7033,855.0264) and (382.0925,853.4298) .. (380.1420,852.5114) .. controls
    (379.5907,852.2519) and (379.0157,852.0462) .. (378.4594,851.7976) .. controls
    (377.9031,851.5489) and (377.3597,851.2530) .. (376.9129,850.8387) .. controls
    (376.4539,850.4130) and (376.1107,849.8744) .. (375.8489,849.3057) .. controls
    (375.5870,848.7371) and (375.4029,848.1363) .. (375.2263,847.5357) .. controls
    (374.8731,846.3344) and (374.5217,845.0782) .. (373.6839,844.1477) .. controls
    (373.2330,843.6470) and (372.6594,843.2662) .. (372.0456,842.9882) .. controls
    (371.4318,842.7102) and (370.7771,842.5319) .. (370.1166,842.3988) .. controls
    (368.7955,842.1327) and (367.4352,842.0419) .. (366.1494,841.6386) .. controls
    (363.7773,840.8948) and (361.8516,839.1499) .. (360.2263,837.2689) .. controls
    (358.6010,835.3879) and (357.1861,833.3165) .. (355.3858,831.6022);
  \path[draw=black,dash pattern=on 2.15pt off 1.07pt,line join=miter,line
    cap=butt,miter limit=4.00,line width=1.074pt] (338.2419,942.7908) .. controls
    (326.7005,946.1210) and (316.0882,952.6121) .. (307.8666,961.3700) .. controls
    (303.4110,966.1163) and (299.5687,971.6965) .. (298.1976,978.0605) .. controls
    (297.5121,981.2425) and (297.4685,984.5834) .. (298.2642,987.7396) .. controls
    (299.0599,990.8959) and (300.7149,993.8586) .. (303.1205,996.0513) .. controls
    (306.6592,999.2770) and (311.5632,1000.6444) .. (316.3457,1000.8788) ..
    controls (321.1281,1001.1132) and (325.8924,1000.3243) .. (330.6481,999.7671)
    .. controls (336.0627,999.1327) and (341.5359,998.7923) .. (346.8324,997.5009)
    .. controls (352.1288,996.2095) and (357.3315,993.8696) .. (361.0233,989.8582)
    .. controls (364.9522,985.5892) and (366.8626,979.6563) .. (366.6846,973.8573)
    .. controls (366.5065,968.0582) and (364.3403,962.4142) .. (361.0233,957.6542)
    .. controls (355.7097,950.0290) and (347.3616,944.5823) ..
    (338.2419,942.7908);
  
  \path[fill=black] (77.952759,995.66925) node[above right] (text3859) {$c_{-1}$};
  \path[fill=black] (233.85826,818.50391) node[above right] (text3863) {$c_{0}$};
  \path[fill=black] (370,1005) node[above right] (text3867) {$c_{1}$};
  \path[draw=black,dash pattern=on 2.23pt off 1.12pt,line join=miter,line
    cap=butt,miter limit=4.00,line width=1.115pt] (28.1267,972.5526) .. controls
    (25.8045,971.9898) and (23.2630,972.3968) .. (21.2327,973.6567) .. controls
    (19.8710,974.5016) and (18.7571,975.6945) .. (17.7856,976.9689) .. controls
    (15.8359,979.5267) and (14.3750,982.4851) .. (12.0406,984.6975) .. controls
    (11.2687,985.4291) and (10.4057,986.0741) .. (9.7426,986.9057) .. controls
    (8.8755,987.9930) and (8.3905,989.3791) .. (8.3905,990.7699) .. controls
    (8.3905,992.1607) and (8.8755,993.5468) .. (9.7426,994.6342) --
    (14.3386,1001.2587) .. controls (15.4466,1002.0564) and (16.5978,1002.7939) ..
    (17.7856,1003.4669) .. controls (19.2656,1004.3054) and (20.8022,1005.0437) ..
    (22.3817,1005.6751) -- (25.8287,1007.8831) .. controls (30.2163,1008.0433) and
    (34.5655,1009.1831) .. (38.4677,1011.1955) .. controls (42.1902,1013.1152) and
    (45.5629,1015.8366) .. (49.6063,1016.9291) .. controls (52.4502,1017.6975) and
    (55.5722,1017.5839) .. (58.2468,1016.3492) .. controls (59.5842,1015.7318) and
    (60.7989,1014.8413) .. (61.7623,1013.7271) .. controls (62.7258,1012.6130) and
    (63.4347,1011.2745) .. (63.7795,1009.8425) .. controls (64.3350,1007.5355) and
    (63.9395,1005.1235) .. (63.7795,1002.7559) .. controls (63.6714,1001.1555) and
    (63.6726,999.5489) .. (63.7457,997.9465) .. controls (63.8675,995.2776) and
    (64.1885,992.6080) .. (64.0198,989.9417) .. controls (63.8511,987.2754) and
    (63.1443,984.5531) .. (61.4477,982.4893) .. controls (59.7766,980.4564) and
    (57.1851,979.2114) .. (54.5537,979.1771) -- (28.1267,972.5526);
  \path[->,draw=black,dash pattern=on 2.39pt off 1.20pt,line join=miter,line
    cap=butt,miter limit=4.00,line width=1.196pt] (154,988.7131) -- (154,865);
  \path[->,draw=black,dash pattern=on 2.39pt off 1.20pt,line join=miter,line
    cap=butt,miter limit=4.00,line width=1.366pt] (352.6917,860) --
    (353.0243,970.2640);
  \path[->,draw=black,dash pattern=on 2.39pt off 1.20pt,line join=miter,line
    cap=butt,miter limit=4.00,line width=1.196pt] (129,979.0761) -- (129,860);
  \path[->,draw=black,dash pattern=on 2.39pt off 1.20pt,line join=miter,line
    cap=butt,miter limit=4.00,line width=1.366pt] (323.4950,852.1590) --
    (323.8277,963.7482);
  \path[fill=black] (162,938.97638) node[above right] (text4032)
    {$\nu_{\hat{g}} (\iota_{-1}^0)$};
  \path[fill=black] (160,955) node[above right] (text4032)
    {$= \gamma_0 . \nu_{\hat{g}'} (\iota_{-1}^0)$};
  \path[fill=black] (358,931.88977) node[above right] (text4032-8)
    {$\nu_{\hat{g}} (\iota_0^1)$};
  \path[fill=black] (78,924.80316) node[above right] (text4032-5)
    {$\nu_{\hat{g}'} (\iota_{-1}^0)$};
  \path[fill=black] (275,924.80316) node[above right] (text4032-9)
    {$\nu_{\hat{g}'} (\iota_0^1)$};
\end{scope}
\end{tikzpicture}
\caption[Orbifold geodesics and pairs of changes of charts]{Lifts of orbifold geodesics in the singular locus related by pairs of embeddings.}\label{fig: lifts}
\end{figure}
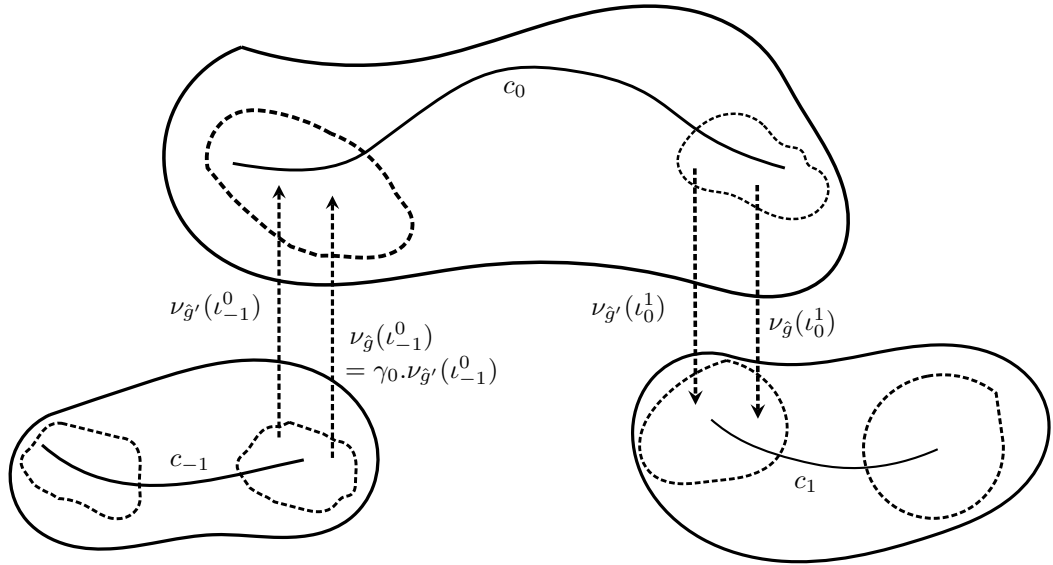

\noindent
The isometry $\gamma_0$ fixes the geodesic $c_0$ pointwise on the set $\im c_0 \cap \cod \nu_{\hat{g}'} (\iota_{-1}^0)$ since
  \begin{equation}\label{eq: commute}
   \gamma_0 c_0|_{]l(0),r(-1)[} = \gamma_0 \nu_{\hat{g}'} (\iota_{-1}^0) c_{-1}|_{]l(0),r(-1)[} = \nu_{\hat{g}} (\iota_{-1}^0) c_{-1}|_{]l(0),r(-1)[} = c_0|_{]l(0),r(-1)[}.
  \end{equation}
Hence $\gamma_0.c_0 = c_0$ follows. Set $\eta_0 \coloneq \gamma_0^{-1}$ and $\eta_{-1} \coloneq \id_{V_{-1}}$ to obtain $\eta_0 .c_0 = c_0$ and $\eta_{-1}.c_{-1}=c_{-1}$. Furthermore, \eqref{eq: ident:flip} yields $\nu_h (\iota_{-1}^0) = \eta_0 \nu_{\hat{g}} (\iota_{-1}^0) \id_{V_{-1}} = \gamma_0^{-1} \nu_{\hat{g}} (\iota_{-1}^0) = \nu_{\hat{g}'} (\iota_{-1}^0)$. Proceed by induction on $k\geq 1$: Consider $k\geq 1$ such that for $0\leq l <k$ elements $\eta_l \in G_l$ have been chosen with 
  \begin{displaymath}
   \eta_l . c_l = c_l \text{ and } \nu_h ( \iota_{l-1}^l) = \eta_{l}. \nu_{\hat{g}} (\iota_{l-1}^l) \eta_{l-1}^{-1}|_{\dom \nu_{\hat{g}} (\iota_{l-1}^l)} = \nu_{\hat{g}'} (\iota_{l-1}^l).
  \end{displaymath}
 We have to choose $\eta_k$ with $\eta_k.c_k =c_k$ and $ \nu_h (\iota_{k-1}^{k}) = \eta_{k}\nu_{\hat{g}} (\iota_{k-1}^{k}) \eta_{k-1}^{-1}|_{\dom \nu_{\hat{g}} (\iota_{k-1}^k)}$. Argue as in the case $k=0$: Since the embeddings of orbifold charts share the same domain, there is $\gamma_k \in G_k$ with $\gamma_k. \nu_{\hat{g}} (\iota_{k-1}^k) = \nu_{\hat{g}'} (\iota_{k-1}^k)$. A computation as \eqref{eq: commute} shows that $\gamma_k$ fixes $\im c_k$ pointwise. Since $\dom \nu_{\hat{g}} (\iota_{k-1}^k)$ is $G_{k-1}$-stable and $\eta_{k-1}$ fixes $\im c_{k-1}$ pointwise, $\eta_{k-1} (\dom \nu_{\hat{g}} (\iota_{k-1}^k)) = \dom \nu_{\hat{g}} (\iota_{k-1}^k)$ follows. Thus we consider the embedding of orbifold charts $\lambda \coloneq \nu_{\hat{g}} (\iota_{k-1}^k) \eta_{k-1}^{-1}|_{\dom \nu_{\hat{g}} (\iota_{k-1}^k)}$. Since $\dom \lambda = \dom \gamma_k \nu_{\hat{g}'} (\iota_{k-1}^k)$, Proposition \ref{prop: ch:prop} (d) yields a unique $h_k \in G_k$ with $\lambda = h_k .\gamma_k . \nu_{\hat{g}'} (\iota_{k-1}^k)$. Define $\eta_k$ via the formula $\eta_k \coloneq (h_k \cdot \gamma_k)^{-1} \in G_k$. We compute the following identities:
  \begin{align*}
   \nu_h (\iota_{k-1}^{k}) &= \eta_k. \nu_{\hat{g}} (\iota_{k-1}^k) \eta_{k-1}^{-1}|_{\dom \nu_{\hat{g}} (\iota_{k-1}^k)} = \eta_k \lambda = \eta_k. \eta_k^{-1} .\nu_{\hat{g}'} (\iota_{k-1}^k) =\nu_{\hat{g}'} (\iota_{k-1}^k) \\
   \eta_k. c_k|_{]l(k),r(k-1)[} &= \eta_k .\nu_{\hat{g}'} (\iota_{k-1}^k) \circ c_{k-1}|_{]l(k),r(k-1)[} = \nu_{\hat{g}} (\iota_{k-1}^k) \eta_{k-1}^{-1}. c_{k-1}|_{]l(k),r(k-1)[} \\ 
                               &=  \nu_{\hat{g}} (\iota_{k-1}^k) \circ c_{k-1}|_{]l(k),r(k-1)[} = c_k|_{]l(k),r(k-1)[}.
  \end{align*}
 Thus the isometry $\eta_k$ fixes the geodesic $c_k$ pointwise on $\im c_k \cap \cod \nu_{\hat{g}'} (\iota_{k-1}^k)$, whence $\eta_k$ fixes all of $\im c_k$ pointwise. We may thus inductively choose elements in $G_k , k \geq 1$, with the required properties. Observe that by (R4) (c) and (d) of Definition \ref{defn: rep:ofdm}, $\nu_{\hat{g}} (\iota_{k}^{k-1})|_{\im \nu_{\hat{g}} (\iota_{k-1}^k)} = \nu_{\hat{g}} (\iota_{k-1}^k)^{-1}$. Instead of choosing $\eta_k$ for $k < 0$ such that $\eta_{k+1} \nu_{\hat{g}} (\iota_{k}^{k+1}) \eta_{k}^{-1}|_{\dom \nu_{\hat{g}} (\iota_{k}^{k+1})} = \nu_{\hat{g}'} (\iota_{k}^{k+1})$, it suffices to choose $\eta_k$ with $\eta_{k} \nu_{\hat{g}} (\iota_{k+1}^{k}) \eta_{k+1}^{-1}|_{\dom \nu_{\hat{g}} (\iota_{k+1}^{k})} = \nu_{\hat{g}'} (\iota_{k+1}^{k})$. If we require that $\eta_k$ fixes $c_k$ pointwise, then an argument as in the case $k \geq 1$ allows us to inductively choose $\eta_k$ for $k < -1$ with the desired properties. Summing up, there is a family $(\eta_k)_{k \in \ZZ}$ such that $\hat{h} = \hat{g}'$ holds, where $\hat{h}$ was constructed via Lemma \ref{lemdef: ind:ofdm} with respect to the pairs $\set{(\id_{]l(k),r(k)[},\eta_k)}_{k\in \ZZ}$. By Lemma \ref{lemdef: ind:ofdm} $\hat{g} \sim \hat{h} = \hat{g}'$. Hence in both cases $[\hat{c}] =[\hat{g}] = [\hat{g}']= [\hat{c}']$ follows from Definition \ref{defn: sim:ofd}.  
\end{proof}

The next lemma is a restatement of Lemma \ref{lem: join:ofg} together with a detailed proof. We shall demonstrate that two orbifold geodesics whose initial vectors coincide in some point induce a well defined join, i.e.\ an orbifold geodesic defined on the union of their respective domains.

\begin{Alem}[Lemma \ref{lem: join:ofg}]\label{GLUE!}
 Consider an orbifold geodesic $[\hat{c}] \in \ORBM[\iI ,(Q,\uU)]$ together with an orbifold geodesic $[\hat{c}'] \in \ORBM[\iI', (Q, \uU)]$ such that for some $x_0 \in \iI \cap \iI'$ their initial vectors coincide. There is an unique orbifold geodesic $[\hat{c} \vee \hat{c}'] \in\ORBM[\iI \cup \iI' , (Q,\uU)]$ such that $[\hat{c} \vee \hat{c}']|_{\iI'} = [\hat{c}']$ and $[\hat{c} \vee \hat{c}']|_{\iI} = [\hat{c}]$.
\end{Alem}

 \begin{proof}[Proof of Lemma \ref{lem: join:ofg}]
  As a first step, we construct an orbifold geodesic on $\iI \cup \iI'$, with the same initial vector at $x_0$: If $\iI \subseteq \iI'$ holds, we set $[\hat{c} \vee \hat{c}'] \coloneq [\hat{c}]$. If $\iI' \subseteq \iI$ holds set $[\hat{c} \vee \hat{c}'] \coloneq [\hat{c}]'$. For these cases, the assertion follows from Proposition \ref{prop: unq:sdat} (b). Interchanging the roles of $[\hat{c}]$ and $[\hat{c}']$ if necessary, it suffices to consider the case $\iI = ]a,b[$ and $\iI' = ]x,y[$ with $a < x < b < y$.\\
  Fix $t_0 \in ]x,b[$ with $t_0 > x_0$. We construct an orbifold geodesic by gluing several pieces: Choose representatives $\hat{c} = (c, \set{c_k}_{k \in A}, [P_{\hat{c}}, \nu_{\hat{c}}])$ of $[\hat{c}]$ and $\hat{c}' = (c', \set{c_r'}_{r \in B}, [P_{\hat{c}'}, \nu_{\hat{c}'}])$ of $[\hat{c}']$ such that the lifts are defined on charts, which are contained in $\aA_\iI$ and $\aA_{\iI'}$, respectively. Since the initial vectors of $[\hat{c}]$ and $[\hat{c}']$ at $x_0$ coincide, they coincide at each point in $\iI \cap \iI' = ]x,b[$ by Proposition \ref{prop: unq:sdat}. By a combination of Lemma \ref{lem: og:equiv} (d) and Lemma \ref{lemdef: ind:ofdm} we may thus assume that there are $k_{t_0} \in A, r_{t_0} \in B$ with $t_0 \in \dom c_{k_{t_0}} = \dom c_{r_{t_0}}' \subseteq ]x,b[$, such that $c_{r_{t_0}}' = c_{k_{t_0}}$ holds. Proposition \ref{prop: unq:sdat} implies that 
  \begin{displaymath} 
   c \vee c' \colon ]a,y[ \rightarrow Q , t \mapsto \begin{cases}
							c(t) & t \in ]a,b[ \\
							c'(t) & t \in ]x,y[																	    \end{cases}
  \end{displaymath}
  is a continuous map. Restricting the lifts (cf.\ proof of Lemma \ref{Alem:ngeod}), we obtain representatives $\hat{c}|_{]a,t_0[}$ induced by $\hat{c}$ and $\hat{c}'|_{]t_0,y[}$ induced by $\hat{c}'$: \\
  The lifts of these mappings are precicsely the restrictions of lifts $c_k, c_r'$ such that the intersections $\dom c_k \cap ]a,t_0[$ and $\dom c_r' \cap ]t_0,y[$ are non-empty. As these intersections may coincide, we choose new index sets $R,S$ for these atlases. Since the domain atlases of $\hat{c}$ and $\hat{c}'$ are contained in $\aA_\iI$ an $\aA_{\iI'}$ , respectively, the domain atlas of $\hat{c}\hat{c}|_{]a,t_0[}$ is contained in $\aA_{]a,t_0[}$ and the domain atlas of $\hat{c}'|_{]t_0,y[}$ is contained in $\aA_{]t_0,y[}$. By construction, $\hat{c}|_{]a,t_0[} = (c|_{]a,t_0[}, \set{g_k}_{k \in R}, [P_{]a,t_0[}, \nu_{]a,t_0[}])$ is obtained by restriction of all data to the open set $]a,t_0[$, i.e: There is a map $\alpha \colon R \rightarrow A$ such that the lifts satisfy $g_k = c_{\alpha (k)}|_{\dom c_{\alpha (k)} \cap ]a,t_0[}$. Each element in $P_{]a,t_0[}$ is constructed as the restriction of an element in $P_{\hat{c}}$ to an open subset of its domain and $\nu_{]a,t_0[}(\mu|_{\dom \mu \cap ]a,t_0[}) \coloneq \nu_{\hat{c}}(\mu)$. As $U_{t_0} \coloneq \dom c_{k_{t_0}}\cap ]a,t_0[ \neq \emptyset$ holds, this chart is contained in the domain atlas $\wW_{]a,t_0[}$ of $\hat{c}|_{]a,t_0[}$. Let $i \colon U_{t_0} \rightarrow \dom c_{k_{t_0}}$ be the inclusion of sets. Define change of charts as follows: For $\lambda \in P_{]a,t_0[}$ and $(W,G,\psi) \in \wW_{]a,t_0[} \in \aA_{]a,t_0[}$, 
      \begin{displaymath}
       \lambda_{t_0} \coloneq \begin{cases}
                                 \lambda \circ (i|^{\im i \cap i(\dom \lambda)})^{-1} & \text{if } \lambda \in \CH{U_{t_0}}{W}\\
				 i \circ \lambda &  \text{if } \lambda \in \CH{W}{U_{t_0}}\\
				 i \circ \lambda \circ (i|^{\im i\cap i(\dom \lambda)})^{-1} & \text{if } \lambda \in \CH{U_{t_0}}{U_{t_0}}.
                                \end{cases}
      \end{displaymath}
  Each of these changes of charts is well defined and $\lambda \neq \mu$ implies $\lambda_{t_0} \neq \mu_{t_0}$. Thus we may define $\nu_{t_0} (\lambda_{t_{0}}) \coloneq \nu_{]a,t_0[} (\lambda)$. Furthermore, set $\nu_{t_0} (\id_{\dom c_{k_{t_0}}}) \coloneq \id_{V_{k_{t_0}}}$, $\nu_{t_0} (i) \coloneq \id_{V_{k_{t_0}}}$ and $\nu_{t_0} (i^{-1}) \coloneq \id_{V_{k_{t_0}}}$. We obtain a set of changes of charts 
    \begin{displaymath}
     C_{t_0} \coloneq \setm{\lambda_{t_{0}}}{\lambda \in \CH{U_{t_0}}{W} \cup \CH{W}{U_{t_0}} \cup \CH{U_{t_0}}{U_{t_0}}, \ W \in \wW_{]a,t_0[}} \cup \set{\id_{\dom c_{k_{t_0}}}, i, i^{-1}}.
    \end{displaymath}
  Since $P_{]a,t_0[}$ is a quasi-pseudogroup, the construction implies that $C \coloneq C_{t_0} \sqcup P_{]a,t_0[}$ is a quasi-pseudogroup which generates $\PSI (\wW_{]a,x[} \cup \set{(\dom c_{k_{t_0}^x}, \set{\id_{\dom c_{k_{t_0}^x}}}, \dom c_{k_{t_0}^x} \hookrightarrow ]a,\sup \dom c_{k_{t_0}^x}[)})$. Our previous observations imply that for an atlas $\bB \in \uU$ containing the codomains of the lifts $\set{g_k}_{k \in R}$, the map 
    \begin{displaymath} 
     \nu_C \colon C \rightarrow \PSI (\bB), \lambda \mapsto \begin{cases}
			\nu_{t_0} (\lambda)	& \text{ if } \lambda \in C_{t_0} \\
			\nu_{]a,t_0[} (\lambda)	& \text{ if } \lambda \in P_{]a,t_0[}																							    \end{cases}
    \end{displaymath}
  is well defined. Consider $\hat{c}_{a,t_0} \coloneq (c|_{]a,\sup \dom c_{k_{t_0}}[}, \set{\dom g_{k}}_{k \in R} \cup \set{c_{k_{t_0}}}, C,\nu_C)$. The map $\nu_{]a,t_0[}$ satisfies property (R4) of Definition \ref{defn: rep:ofdm}. Together with the definition of $\lambda_{t_0}$ and $\nu_C$, this implies that $\nu_C$ satisfies the property (R4). Hence $\hat{c}_{a,t_0}$ is a representative of an orbifold map such that each lift is a geodesic defined on a chart in $\aA_{]a,\sup \dom c_{k_{t_0}}[}$. In other words, $[\hat{c}_{a,t_0}]$ is an orbifold geodesic whose initial vector at any point in its domain coincides with the corresponding one of $[\hat{c}]$.\\ 
  Note that in the domain atlas of $\hat{c}_{a,t_0}$, only $(\dom c_{k_{t_0}}, \set{\id_{\dom c_{k_{t_0}}}}, \dom c_{k_{t_0}} \hookrightarrow ]a,\sup \dom c_{k_{t_0}}[)$ intersects $[t_0,b[$. We may thus interpret this chart as an \tl adhesive joint\tr .\\
  Repeat the construction for $\hat{c}'$: We obtain $\hat{c}_{t_0,y} \coloneq (c'|_{]\inf \dom c_{r_{t_0}'},y[}, \set{h_k}_{k \in S} \cup \set{c_{r_{t_0}}'}, D, \nu_D)$. Again only the chart with domain $\dom c_{r_{t_0}}' = \dom c_{k_{t_0}}$ in its domain atlas intersects $]a,t_0]$.\\
  We will glue the geodesics $\hat{c}_{a,t_0}, \hat{c}_{t_0,y}$ at their \tl adhesive joints\tr to obtain a geodesic on $]a,y[$: With the exception of $\id_{\dom c_{k_{t_0}}} = \id_{\dom c_{r_{t_0}}}$, the quasi-pseudogroups $C$ and $D$ contain only changes of charts, whose domains are contained in $]a,t_0[$ (for $C$) respectively in $]t_0,y[$ (for $D$). In particular, $C\cap D = \set{\id_{\dom c_{k_{t_0}}}}$ holds, whence we obtain a disjoint union:
    \begin{displaymath}
     C \cup D = \set{\id_{\dom c_{k_{t_0}}}} \sqcup C \setminus \set{\id_{\dom c_{k_{t_0}}}} \sqcup D \setminus \set{\id_{\dom c_{k_{t_0}}}} .                                                                                                                                                                                                                                                                                                                                                                                                                                                                                                                                                                                                                                                                                                                                               \end{displaymath}
 Consider $\lambda, \mu \in C \cup D$. If $\lambda \in C \setminus D$ and $\mu \in D$ such that the composition is defined on some open subset of their domains, then $\mu = \id_{\dom c_{k_{t_0}}} \in C$. Vice versa, an analogous condition holds for elements in $D \setminus C$. Thus any pair in $(C\setminus D) \times (D \setminus C)$ may not be composed on any open subset of their respective domains. As both sets $C,D$ are quasi-pseudogroups, $P^\star \coloneq C \cup D$ is a quasi-pseudogroup which generates the changes of charts of the atlas whose domains are given by $\setm{\dom h_s}{s \in S} \cup \setm{g_k}{k \in R} \cup \set{\dom c'_{r_{t_0}}}$. Define 
  \begin{displaymath}
   \nu^\star (\lambda) \coloneq \begin{cases}
                        \nu_D (\lambda) & \text{ if } \lambda \in D \\
                        \nu_C (\lambda) & \text{ if } \lambda \in C .
                       \end{cases}
  \end{displaymath}
 As $\nu_C (\id_{c_{k_{t_0}}}) = \id_{V_{k_{t_0}}} = \id_{V_{r_{t_0}}} = \nu_D (\id_{c_{k_{t_0}}})$ holds, the map $\nu^\star$ is well defined. Since $\nu_C$ and $\nu_D$ satisfy condition (R4), the same holds for $\nu^\star$ with respect to the lifts $\setm{h_s}{s\in S} \cup \set{c_{k_{t_0}}} \cup \setm{g_k}{k \in R}$. Hence $\hat{c}^\star \coloneq (c \vee c', \setm{h_s}{s\in S} \cup \set{c_{k_{t_0}}} \cup \setm{g_k}{k \in R}, P^\star, \nu^\star)$ is a representative of an orbifold geodesic on $]a,y[$.\\
 Observe that the initial vector of $\hat{c}^\star$ at $x_0$ coincides by construction with the initial vector of $[\hat{c}]$ at $x_0$. As the initial vector of $[\hat{c}]$ coincides with the one of $[\hat{c}']$ in $x_0$, $[\hat{c} \vee \hat{c}'] \coloneq [\hat{c}^\star]$ satisfies the first assertion by Proposition \ref{prop: unq:sdat}.
\end{proof}
\end{appendix}

\bibliographystyle{babplain-lf}
\thispagestyle{empty}
\phantomsection
\addcontentsline{toc}{section}{References}
\bibliography{Orbifoldmain}


\newpage \phantomsection
\addcontentsline{toc}{section}{List of Figures}
\listoffigures

\newpage 
\phantomsection
\addcontentsline{toc}{section}{Index}
\printindex
\end{document}